%% file: main.tex
\documentclass[hidelinks, 10pt]{amsbook}
\usepackage{amsthm,amssymb,amscd,amsmath,amsfonts,amssymb,amscd,stmaryrd,hyperref,mathrsfs,yhmath,mathtools, amsrefs}                                                                                                                    
\usepackage[all]{xy}
\usepackage{hyperref}
\usepackage{enumerate}
\usepackage{physics}
\usepackage{ragged2e}
\usepackage{xfrac, faktor}
\usepackage{lipsum}  
\usepackage{tikz-cd}
\usepackage{enumitem}
\usepackage{imakeidx}
\usepackage{subfiles}
\usepackage{bbm}
\usepackage{quiver}
\usepackage{mdframed}
\usepackage[margin = 1.2in]{geometry}
\usepackage{makecell}
\usepackage{chngcntr}
\counterwithin{table}{chapter}
\usepackage{booktabs}
\usepackage{algorithm}
\usepackage{algpseudocode}
\usepackage{changepage}

\usepackage{amsmath,amsfonts,amscd,amssymb,verbatim,amsthm}

\algrenewcommand\algorithmicrequire{\textbf{Input:}}
\algrenewcommand\algorithmicensure{\textbf{Output:}}

\newcommand{\HB}{\mathrm{B}}
\newcommand{\HC}{\mathrm{C}}
\newcommand{\HH}{\mathrm{H}}
\newcommand{\HZ}{\mathrm{Z}}
\newcommand{\HO}{\mathrm{O}}

\DeclareMathOperator{\Hom}{Hom}
\DeclareMathOperator{\Ext}{Ext}

\DeclareMathOperator{\Mod}{\mathbf{Mod}}

\DeclareMathOperator{\Gal}{Gal}

\DeclareMathOperator{\charfield}{char}

\DeclareMathOperator{\Vect}{\mathbf{Vect}}

\DeclareMathOperator{\SD}{SD}
\DeclareMathOperator{\Rep}{\mathbf{Rep}}
\DeclareMathOperator{\End}{End}

\DeclareMathOperator{\Aff}{Aff}
\DeclareMathOperator{\Euc}{Euc}
\DeclareMathOperator{\Out}{Out}
\DeclareMathOperator{\Inn}{Inn}

\DeclareMathOperator{\Char}{Char}

\newcommand{\bLG}{\bL_{\sigma} \hspace{-0.1em}[G]}

\DeclareMathOperator{\fd}{fd}
\DeclareMathOperator{\semisimple}{ss}
\DeclareMathOperator{\id}{id}

\DeclareMathOperator{\Ind}{Ind}
\DeclareMathOperator{\Cnd}{Coind}
\DeclareMathOperator{\Inf}{Inf}
\DeclareMathOperator{\Irr}{Irr}
\DeclareMathOperator{\UIrr}{UIrr}
\DeclareMathOperator{\URep}{URep}
\DeclareMathOperator{\Prof}{Prof}
\DeclareMathOperator{\GL}{GL}
\DeclareMathOperator{\PGL}{PGL}
\DeclareMathOperator{\Br}{Br}
\DeclareMathOperator{\BrW}{BW}

\DeclareMathOperator{\Fun}{Fun}

\DeclareMathOperator{\Ad}{Ad}
\DeclareMathOperator{\opp}{op}

\DeclareMathOperator{\UR}{U_R}
\DeclareMathOperator{\UU}{U}

\DeclareMathOperator{\Aut}{Aut}
\DeclareMathOperator{\Stab}{Stab}

\DeclareMathOperator{\GradDiv}{GradDiv}

\DeclareMathOperator{\Ab}{\textbf{Ab}}

\newcommand{\twopartdef}[4]
{
	\left\{
		\begin{array}{ll}
			#1 & \mbox{if } #2 \\
			#3 & \mbox{if } #4
		\end{array}
	\right.
}

\makeatletter
\newcommand{\colim@}[2]{%
  \vtop{\m@th\ialign{##\cr
    \hfil$#1\operator@font colim$\hfil\cr
    \noalign{\nointerlineskip\kern1.5\ex@}#2\cr
    \noalign{\nointerlineskip\kern-\ex@}\cr}}%
}
\newcommand{\colim}{%
  \mathop{\mathpalette\colim@{\rightarrowfill@\textstyle}}\nmlimits@
}
\makeatother

\newcommand{\bC}{{\mathbb C}}

\newcommand{\bF}{{\mathbb F}}

\newcommand{\bH}{{\mathbb H}}

\newcommand{\bK}{{\mathbb K}}
\newcommand{\bL}{{\mathbb L}}
\newcommand{\bN}{{\mathbb N}}

\newcommand{\bQ}{{\mathbb Q}}
\newcommand{\bR}{{\mathbb R}}
\newcommand{\bZ}{{\mathbb Z}}

\newcommand{\bbL}{\, {\mathbb L}}
\newcommand{\bbK}{\, {\mathbb K}}

\newcommand{\sA}{{\mathcal A}}
\newcommand{\sB}{{\mathcal B}}
\newcommand{\sC}{{\mathcal C}}

\newcommand{\sF}{{\mathcal F}}
\newcommand{\sG}{{\mathcal G}}
\newcommand{\sH}{{\mathcal H}}
\newcommand{\sI}{{\mathcal I}}

\newcommand{\sL}{{\mathcal L}}
\newcommand{\sM}{{\mathcal M}}
\newcommand{\sN}{{\mathcal N}}

\newcommand{\sQ}{{\mathcal Q}}
\newcommand{\sR}{{\mathcal R}}

\newcommand{\sT}{{\mathcal T}}

\newcommand{\sW}{{\mathcal W}}

\newcommand{\fC}{{\mathfrak C}}
\newcommand{\fD}{{\mathfrak D}}

\newcommand{\fm}{{\mathfrak m}}

\newcommand{\fR}{{\mathfrak R}}
\newcommand{\fS}{{\mathfrak S}}
\newcommand{\fT}{{\mathfrak T}}
\newcommand{\fL}{{\mathfrak L}}

\newtheorem{thm}{Theorem}
\numberwithin{thm}{section}
\newtheorem{lemma}[thm]{Lemma} 
\newtheorem{prop}[thm]{Proposition} 
\newtheorem{cor}[thm]{Corollary} 
 
\newtheorem{prob}[thm]{Problem} 
\newtheorem{just}[thm]{Justification} 
\newtheorem*{theoremA}{Theorem A}

\newtheorem*{theoremB}{Theorem B}
\newtheorem*{theoremC}{Theorem C}

\newtheorem*{theorem*}{Theorem}
\newtheorem*{question1*}{Question 1}
\newtheorem*{question2*}{Question 2}

\theoremstyle{definition}
\newtheorem{defn}[thm]{Definition}
\newtheorem{eg}[thm]{Example}
\newtheorem{remark}[thm]{Remark}

\newtheorem*{definition*}{Definition}

\title{Dyson's Tenfold Way, Revisited}
\author{Dmitriy Rumynin}
\address{Department of Mathematics, University of Warwick, 
Coventry, CV4 7AL, UK}
\email{D.Rumynin@warwick.ac.uk}
\author{James Taylor}
\address{School of Mathematics, University of Bristol, Bristol, UK}
\email{james.h.a.taylor@bristol.ac.uk}

\begin{document}
\begin{abstract}
  In this work, we study semilinear representations of groups acting on degree two Galois extensions. Our results include a classification of irreducible semilinear representations, which generalises Dyson's Tenfold Way. Using this, we completely describe the semilinear representation theory of certain groups, such as symmetric groups, and of finite groups over cyclotomic fields. 

We also prove many results concerning semilinear representations over $\bC$, which includes an induction theorem which simultaneously generalises Brauer's Induction Theorem and Serre's Induction Theorem. We are also interested in applications of our classification theorem to crystallography, and include a detailed study of the semilinear representation theory of magnetic wallpaper groups.
\end{abstract}
\maketitle

\setcounter{tocdepth}{1}
\tableofcontents

\include{section-introduction}

\include{chapter-realstructures}
\include{chapter-algebras}
\include{chapter-twelveways}

\include{chapter-tenways}
\include{chapter-physics}

\appendix
\include{appendix-tablesfor24types}
\include{appendix-classandcohomC2modules}
\include{appendix-gradedalgebras}

\include{appendix-classicalresults}
\include{appendix-wallpaper}
\include{appendix-repswallpaper}

\backmatter
\bibliography{biblio}{}
\bibliographystyle{plain}

\end{document}

%% file: section-introduction.tex
\chapter*{Introduction}

Antilinear representations play an important role in algebra, physics and topology \cite{AtiSeg, AtiKtheory, NeebOlaf, DYSON, BradDav, FreMoo, TopIns, Bootstrap}. Besides physics, they also appear in chemistry and crystallography \cite{Bro, Ceulemans, JensenSaue}. They were first introduced by Wigner \cite{WIG} under the name of co-representations to model time-reversal symmetry in quantum mechanics. 

Algebraically, for $G$ a group with an index two subgroup $H$, an antilinear representation of $G$ is an action of $G$ on a complex vector space such elements of $H$ act linearly, whilst elements of $G\setminus H$ act antilinearly. The key result in the theory is {\em Dyson's Tenfold Way}, which classifies irreducible antilinear representations into 10 fundamentally different classes. The goal of this text is to generalise Dyson's Tenfold Way to the following setting:
\begin{itemize}
\item the group $G$ is any (discrete) group, 
\item instead of complex vector spaces, we work with any pair $(\bL, \tau)$, where $\bL$ is a field and $\tau$ is an order two automorphism of $\bL$. We study actions of $G$ on finite-dimensional vector spaces over $\bL$, such that elements of $H$ act linearly and elements of $G \setminus H$ act $\tau$-linearly. We write
  \[ \Irr_{\bbL}^{\hspace{0.1em} \star}(G)\]
  for the set of (isomorphism classes of) such irreducible representations of $G$ over $\bL$.
\end{itemize}

The main outcome of our work is that we show these representations fall into different {\em types}, that capture key algebraic properties of a representation. For $\bL = \bC$, we recover Dyson's Tenfold Way. 

As we will use the term ``type'' in a few different ways, let us summarise, without giving precise definitions, the various notations of a type we will use. For $V$ a representation, we have:
\begin{itemize}
\item {\em The type} $\fT(V)$. This is the main overriding type that keeps track of all the algebraic information about $V$.
As a result of our classification, there are 24 types; the reader can find them in Table~\ref{tab:24types}.
By $\fT(V)$ we understand either the number from 1 to 24, or the actual information in the corresponding row, some sort of a ``barcode''.
\item {\em The D-type} $\fD(V)$, or {\em Dyson type}. 
This is an alternative name for the type $\fT(V)$ which we use in certain situations, such as when $\bL$ is algebraically closed, when only 12 types are possible.
$\fD(V)$ is the Roman numeral in the leftmost column of Table~\ref{tab:24types}.
It is called the Dyson type because when $\bL = \bC$, only 10 of the 12 are possible, and $\fD(V)$ is simply the label of Dyson's Tenfold Way.
\item {\em The S-type} $\fS(V)$, or {\em semilinear type}.
This is the information in the column $\bL[H] \subset \bL_{\star \hspace{-0.1em}}[G]$ of Table~\ref{tab:24types}.
There are 3 possible S-types: \textbf{A}, \textbf{B} or \textbf{C}.
\item {\em The C-type} $\fC(V)$, or {\em classical type}.
This is the information in the column $\bK[H] \subset \bL[H]$ of Table~\ref{tab:24types}.
There are 3 possible C-types: \textbf{A}, \textbf{B} or \textbf{C}.
\item {\em The L-type} $\fL(V)$, or {\em linear type}.
  This is the information in the four columns
of Table~\ref{tab:24types}
that do not contain  $\bL_{\star \hspace{-0.1em}}[G]$. By inspection,  
there are 18 possible L-types.
\end{itemize}  
For instance, if $\fT(V)=5$, then we examine row 5 of Table \ref{tab:24types} and conclude that
\[ \fD(V) = \mbox{II}, \ \ \fS(V) = \textbf{C}, \ \ \fC(V) = \textbf{A}, \ \ \fL(V) =( \textbf{B},\textbf{A},\textbf{1},\textbf{3(b)}).\]
Now we outline the main ideas behind these notions of types.

\section*{D-Type}

Let us first recall, through our lens of types, the theory of $\bR$-linear and antilinear representations, and Dyson's Tenfold Way. The types that play a role here are the C-type, the S-type, and the D-type.

\subsection*{C-Type and $\bR$-Linear Representations}

Let $H$ be a (discrete) group.
Let us examine the relationship between the sets
$\Irr_{\bR}(H)$ and $\Irr_{\bC}(H)$
of finite-dimensional irreducible $\bR$-linear and $\bC$-linear representations.
Explicitly, for any $V \in \Irr_{\bC}(H)$, there is a unique $X \in \Irr_{\bR}(H)$ with $\Hom_{\bR[H]}(X, V|_{\bR}) \neq 0$.
For this $X$, the endomorphism ring $\End_{\bR[H]}(X)$ is a finite-dimensional division algebra over $\bR$.
By Frobenius' Theorem there are three possibilities:

\begin{itemize}
    \item[(\textbf{A})] $\End_{\bR[H]}(X) \cong \bR$, $X_{\bC} \cong V$, and $V|_{\bR} \cong X \oplus X$.
    \item[(\textbf{B})] $\End_{\bR[H]}(X) \cong \bC$, $X_{\bC} \cong V \oplus \overline{V}$, and $V|_{\bR} \cong X$.
    \item[(\textbf{C})] $\End_{\bR[H]}(X) \cong \bH$, $X_{\bC} = V \oplus V$, and $V|_{\bR} \cong X$.
\end{itemize}

In our language of types, still yet to be defined, the label (\textbf{A}), (\textbf{B}) or (\textbf{C}) that occurs is exactly the C-type $\fC(X)$ of $X$, or equally the C-type $\fC(V)$ of $V$. 

Observe that there are two different algebraic ways to determine the C-type. The first is via the endomorphism algebra $\End_{\bR[H]}(X)$. The second is through the interplay between induction and restriction functors. This interplay is the main perspective behind our theory of types.

\subsection*{The S-Type and Antilinear Representations}
Now let $G$ be a group which acts on $\bC$ non-trivially, in the sense that there is a surjective group homomorphism $\sigma \colon G \twoheadrightarrow \Gal(\bC / \bR)$.

Similarly to $\bR$-linear representations, we can relate irreducible antilinear representations $\Irr_{\bC}^{\sigma}(G)$ of $G$ to $\Irr_{\bC}(H)$, where now $H$ denotes the index two normal subgroup $H \coloneqq \ker(\sigma)$ of $G$.
Consider the induction and restriction functors
\[
(-)|_H \colon \Rep_{\bC}^{\sigma}(G) \rightarrow \Rep_{\bC}(H), \qquad \Ind_H^G(-) \colon \Rep_{\bC}(H) \rightarrow \Rep_{\bC}^{\sigma}(G),
\]
and an action of $C_2 = \langle \epsilon \rangle$ on $\Irr_{\bC}(H)$, which is defined by fixing some $g \in G \setminus H$, and for $V \in \Irr_{\bC}(H)$ setting $\epsilon * V$ to be $\overline{V}$ as a $\bC$ vector space, with action of $h \in H$ by $h * v \coloneqq (g^{-1} hg)(v)$.

For any $V \in \Irr_{\bC}(H)$, there is a unique $W \in \Irr_{\bC}^{\sigma}(G)$ with $\Hom_{\bC_{\sigma}\hspace{-0.1em}[G]}(W, \Ind_H^G(V)) \neq 0$. For this $W$, the endomorphism algebra $\End_{\bC_{\sigma}\hspace{-0.1em}[G]}(W)$ is a finite-dimensional division algebra over $\bR$, and Frobenius' Theorem again results in three possibilities:
\begin{itemize}
    \item[(\textbf{A})] $\End_{\bC_{\sigma}\hspace{-0.1em}[G]}(W) \cong \bR$, $W|_H \cong V$, and $\Ind_H^G(V) \cong W \oplus W$.
    \item[(\textbf{B})] $\End_{\bC_{\sigma}\hspace{-0.1em}[G]}(W) \cong \bC$, $W|_H \cong V \oplus \epsilon * V$, and $\Ind_H^G(V) \cong W$.
    \item[(\textbf{C})] $\End_{\bC_{\sigma}\hspace{-0.1em}[G]}(W) \cong \bH$, $W|_H \cong V \oplus V$, and $\Ind_H^G(V) \cong W$.
\end{itemize}
In our language, this label is exactly the S-type $\fS(W)$ of $W$, or equally the S-type $\fS(V)$ of $V$.

\subsection*{The D-Type and Dyson's Tenfold Way}
Given a pair $(G,\sigma)$ as before, let us consider both interplays above simultaneously
\[
\Irr^{\sigma}_{\bC}(G) \longleftrightarrow \Irr_{\bC}(H)\longleftrightarrow \Irr_{\bR}(H), \quad 
W \leftrightarrow V \leftrightarrow X.
\]
In this situation, it is natural to consider how $W$ and $X$ relate to one another. For example, to ask how the S-type for the pair $(V,W)$ is related to the C-type for the pair $(V,X)$.

To answer this, another linear representation category $\Rep_{\bR}(G)$ enters the scene.
There are induction and restriction functors now between the four categories: 
\begin{equation} \label{2categories}
  \begin{tikzcd}[sep=2.25em]
	& {U \in \Rep_{\bR}(G)} & \\
	{X\in \Rep_{\bR}(H)} && {\Rep_{\bC}^{\sigma}(G) \ni W} \\
	& {V\in \Rep_{\bC}(H)}
	\arrow[from=1-2, to=2-1]
	\arrow[from=1-2, to=2-3]
	\arrow[from=2-1, to=1-2]
	\arrow[from=2-1, to=3-2]
	\arrow[from=2-3, to=1-2]
	\arrow[from=2-3, to=3-2]
	\arrow[from=3-2, to=2-1]
	\arrow[from=3-2, to=2-3]
\end{tikzcd}\end{equation}

Let $U \in \Irr_{\bR}(G)$ such that $\Hom_{\bR[G]}(U, W|_{\bR}) \neq 0$, or equivalently that $\Hom_{\bR[G]}(U, \Ind_{H}^G(X)) \neq 0$. This $U$ is unique up to replacement with $U \otimes \pi$, the tensor product of $U$ with the sign representation $\pi$ of $G$ over $\bR$ (with kernel $H$). For either choice, we have that $\End_{\bR[G]}(U) = \End_{\bR[G]}(U \otimes \pi)$. 

Now pick any $V \in \Irr_{\bC}(H)$, then choose $X, W$ and $U$ as above. We mark where these modules belong on the diagram~\eqref{2categories} for the convenience of the reader. The S-type and the C-type of V can be recovered from the the triple 
\[
(\End_{\bR[H]}(X), \End_{\bR[G]}(U), \End_{\bC_{\sigma}\hspace{-0.1em}[G]}(W)) \in \{\bR, \bC, \bH\}^3 \, ,
\]
which is independent of the choices of $X$, $U$ and $W$. The D-type $\fD(V)$ of $V$, still yet to be defined, is an invariant which records how induction and restriction works between $V$, $X$, $W$ and $U$ in diagram~\eqref{2categories}. It is sufficiently rich to determine the endomorphism ring triple above.

The fundamental result of Dyson, phrased in our language, is the following.
\begin{theorem*}[Dyson's Tenfold Way \cite{DYSON}]
Suppose $G$ is a (discrete) group.
Given $V \in \Irr_{\bC}(H)$, there are ten possibilities for $\fD(V)$, listed with corresponding endomorphism algebra triples:
\[
\renewcommand{\arraystretch}{1.1}
\begin{array}{@{}c@{\;}l@{\qquad}
                c@{\;}l@{\qquad}
                c@{\;}l@{\qquad}
                c@{\;}l@{\qquad}
                c@{\;}l@{}}
\mathrm{I}    & (\bR,\bR,\bR) &
\mathrm{III}  & (\bR,\bR,\bC) &
\mathrm{V}    & (\bC,\bR,\bR) &
\mathrm{VII}  & (\bC,\bC,\bC) &
\mathrm{IX}   & (\bH,\bC,\bR) \\

\mathrm{II}   & (\bR,\bC,\bH) &
\mathrm{IV}   & (\bC,\bC,\bC) &
\mathrm{VI}   & (\bC,\bH,\bH) &
\mathrm{VIII} & (\bH,\bH,\bH) &
\mathrm{X}    & (\bH,\bH,\bC)
\end{array}
\]
\end{theorem*}

We thus see that the endomorphism ring triple is almost enough to determine the D-type, but not quite. The unique repeat is in D-types IV and VII. In Dyson's original notation these are called $\sC \sC 2$ and $\sC \sC 1$ respectively - he uses the pair $(\End_{\bC_{\sigma}\hspace{-0.1em}[G]}(W), \End_{\bR[H]}(X))$ as his label for the type, so II is called $\sQ \sR$ and so on (with letters denoting $\sR$eal, $\sC$omplex and $\sQ$uaternionic).

Conceptually, the main reason that there are 10 possible D-types is that the Brauer-Wall monoid $\BrW(\bR) \sqcup \BrW(\bC)$ has order 10 (see \S \ref{sect:recoverDyson10}). For us to obtain finiteness in our classification of types over general fields, it turns out that the induction-restriction perspective is more profitable, as for general fields the corresponding Brauer-Wall monoids will typically be infinite.

\section*{Semilinear Representations for Degree Two Extensions}
Let us move on from the complex numbers to more general fields.
Suppose now that $G$ is a group, and $\bL$ is a field, upon which $G$ acts through a group homomorphism $\sigma \colon G \rightarrow \Aut(\bL)$.
In this context, a \emph{semilinear representation} of $G$ over $\bL$ is a finite-dimensional $\bL$-vector space $W$, with a homomorphism $\rho \colon G \rightarrow \Aut_{\bK}(W)$, such that each $\rho(g)$ is $\sigma(g)$-linear.

If $\bL / \bK$ is finite, where we write $\bK \coloneqq \bL^G$, then automatically $\bL / \bK$ is Galois, and the map $\sigma \colon G \rightarrow \Gal(\bL / \bK)$ is surjective. Under this assumption, the category of such semilinear representations $\Rep_{\bbL}^{\sigma}(G)$ can be described in terms of $\Rep_{\bbL}(H)$, where $H \coloneqq \ker(\sigma \colon G \rightarrow \Gal(\bL/\bK))$ \cite{TAYSchur}.

\subsection*{S-Type for Degree Two Extensions} 
The focus of the present work is on understanding $\Rep_{\bbL}^{\sigma}(G)$ in the smallest non-trivial case, when $\bL / \bK$ has degree two. From now until the end of introduction, this is an assumption.
Of course, this includes the case when $\bL / \bK$ is $\bC / \bR$.

As for antilinear representations over $\bC$, there is an analogously defined action of $C_2 = \langle \epsilon \rangle$ on $\Irr_{\bL}(H)$, and semilinear induction and restriction functors
\[
(-)|_H \colon \Rep_{\bbL}^{\sigma}(G) \rightarrow \Rep_{\bbL}(H), \qquad \Ind_H^G(-) \colon \Rep_{\bbL}(H) \rightarrow \Rep_{\bbL}^{\sigma}(G).
\]
For any $V \in \Irr_{\bC}^{\sigma}(G)$, there is still a unique $W \in \Irr_{\bbL}^{\sigma}(G)$ with $\Hom_{\bbL_{\sigma}\hspace{-0.1em}[G]}(W, \Ind_H^G(V)) \neq 0$.
The crucial difference is that instead of $\End_{\bbL_{\sigma}\hspace{-0.1em}[G]}(W)$ the key quantity is the ratio:
\[
d(V,W) \coloneqq \dim_{\bbK}\End_{\bbL_{\sigma}\hspace{-0.1em}[G]}(W) / \dim_{\bbL}\End_{\bbL[H]}(V).
\]
With this notation, there are three possibilities, which we call the S-type $\fS(V) = \fS(W)$:
\begin{itemize}
    \item[(\textbf{A})] $d(V,W) = 1$, $W|_H \cong V$, and $\Ind_H^G(V) \cong W \oplus W$.
    \item[(\textbf{B})] $d(V,W) = 2$, $W|_H \cong V \oplus \epsilon * V$, and $\Ind_H^G(V) \cong W$.
    \item[(\textbf{C})] $d(V,W) = 4$, $W|_H \cong V \oplus V$, and $\Ind_H^G(V) \cong W$.
\end{itemize}

\subsection*{Classification of Semilinear Representations}
The reason for our interest in the classification into three types \textbf{A}, \textbf{B} and \textbf{C}, is that it provides a method to construct all irreducible semilinear representations of $G$ over $\bL$ from $\Irr_{\bbL}(H)$, provided that one already knows the S-type of each $V \in \Irr_{\bbL}(H)$.
Concretely, pick $V \in \Irr_{\bbL}(H)$.

\vspace{0.2em}

\noindent \hspace{0.7em} $\bullet$ If $\fS(V)=$ \textbf{A}, then $V$ extends uniquely to an irreducible semilinear representation of $G$ over $\bL$.

\noindent \hspace{0.7em} $\bullet$ If $\fS(V)\in\{\textbf{B},\textbf{C}\}$, then $\Ind_H^G(V)$ is an irreducible semilinear representation of $G$ over $\bL$.

\vspace{0.2em}

\noindent Furthermore, all $W \in \Irr_{\bbL}^{\sigma}(G)$ arise this way.

The question is thus how to determine, given $V \in \Irr_{\bbL}(H)$, whether $V$ has type \textbf{A}, \textbf{B} or \textbf{C}.
We can immediately see that the $C_2$-action allows us to distinguish S-type \textbf{B}:
\[
\text{\centering{$\fS(V)=\textbf{B}$ if and only if $V \not\cong \epsilon * V$.}}
\]
The following question is therefore of fundamental interest:

\begin{question1*}
For $V \in \Irr_{\bbL}(H)$ with $V \cong \epsilon * V$, how can one tell if $V$ has S-type \emph{\textbf{A}} or \emph{\textbf{C}}?
\end{question1*}

\section*{The Type and L-Type of $V$}
How is Question 1 answered when $\bL / \bK$ is $\bC / \bR$?
This is where Dyson's Tenfold Way enters the picture.
If $V \in \Irr_{\bbL}(H)$ has $V \cong \epsilon * V$, then $\End_{\bC_{\sigma}\hspace{-0.1em}[G]}(W)$ is either $\bR$ or $\bH$. Examining the possibilities in Dyson's Tenfold Way above, we see that 
\[
\fS(V) = \twopartdef{\text{\textbf{A}}}{(\End_{\bR[H]}(X), \End_{\bR[G]}(U)) \in \{(\bR, \bR), (\bC, \bR), (\bH, \bC)\}}{\text{\textbf{C}}}{(\End_{\bR[H]}(X), \End_{\bR[G]}(U)) \in \{(\bR, \bC), (\bC, \bH), (\bH, \bH)\}}.
\]
Therefore Question 1 is answered by examining the pair $(\End_{\bR[H]}(X), \End_{\bR[G]}(U))$, a purely linear representation-theoretic quantity, which is computable using classical (linear) methods.

For general fields, we want to understand if something similar is possible in general. Is there some linear data we can associate to $V$, which determines the S-type of $V$?

Besides the four categories on the diagram~\eqref{2categories},
we can bring one more category into play, arriving at {\em the diamond}
\begin{equation}\label{eq:maincats}
\begin{tikzcd}
	& {\Rep_{\bbK}(G)} \\
	{\Rep_{\bbL}(G)} & {\Rep_{\bbK}(H)} & {\Rep_{\bbL}^{\hspace{0.1em}\sigma}(G)} \\
	& {\Rep_{\bbL}(H)}
	\arrow["3"', tail reversed, from=1-2, to=2-1]
	\arrow["4", tail reversed, from=1-2, to=2-3]
	\arrow["4", tail reversed, from=2-2, to=1-2]
	\arrow["3"', tail reversed, from=2-2, to=3-2]
	\arrow["4", tail reversed, from=3-2, to=2-1]
	\arrow["3"', tail reversed, from=3-2, to=2-3]
\end{tikzcd}
\end{equation}
of categories, connected by canonical induction and restriction functors. 

Let us from now on assume that $\charfield(\bK) \neq 2$.
We already know that there are three possibilities for induction and restriction between $\Rep_{\bbL}(H)$ and $\Rep_{\bbL}^{\hspace{0.1em}\sigma}(G)$, recorded in the S-type. It turns out that this finiteness extends to the other pairs: there are either three or four possibilities for induction and restriction between each category in adjacent pair of categories in diagram~\eqref{eq:maincats}.

We label each extension with the exact number. When there are three possibilities, these are called \textbf{A}, \textbf{B} and \textbf{C}, and when there are four these are called $\mathbf{1}, \mathbf{2}$, \textbf{3(a)} and \textbf{3(b)}.

\begin{definition*}
The \emph{type} of $V$ is the tuple
\[
\fT(V) \in \{\mathbf{A},\mathbf{B},\mathbf{C}\}^3 \times \{\text{\textbf{1}}, \text{\textbf{2}}, \text{\textbf{3(a)}}, \text{\textbf{3(b)}}\}^3
\]
which records induction and restriction in diagram~\eqref{eq:maincats}, read from left to right.
\end{definition*}

The S-type is the third entry of $\fT(V)$, recording induction and restriction between $\Rep_{\bbL}(H)$ and $\Rep_{\bbL}^{\hspace{0.1em}\sigma}(G)$. The C-type is the second entry of $\fT(V)$, recording induction and restriction between $\Rep_{\bbL}(H)$ and $\Rep_{\bbK}(H)$. 

\begin{definition*}
The \emph{linear type}, or L-type, of $V$ is the tuple
\[
\fL(V) \in \{\mathbf{A},\mathbf{B},\mathbf{C}\}^2 \times \{\text{\textbf{1}}, \text{\textbf{2}}, \text{\textbf{3(a)}}, \text{\textbf{3(b)}}\}^2
\]
obtained from $\fT(V)$ by forgetting the entries which involve $\Rep_{\bbL}^{\hspace{0.1em}\sigma}(G)$.
\end{definition*}

By definition, the L-type consists of exactly the part of $\fT(V)$ which involves only classical (linear) representation theory.
Using $\fT(V)$ and $\fL(V)$, we can make Question 1 more precise.

\begin{question2*}
    For $V \in \Irr_{\bbL}(H)$, can $\fT(V)$ be recovered from $\fL(V)$?
\end{question2*}

\section*{Classification of Types}

Our first main result is a generalisation of Dyson's Tenfold Way to arbitrary degree two extensions, which yields a complete answer to Question 2.

\begin{theoremA}
For $V \in \Irr_{\bL}(H)$, there are 24 possibilities for $\fT(V)$, listed in Table \ref{tab:24types}.
\end{theoremA}

Examining the 24 types listed in Table \ref{tab:24types}, we see that there are 18 possible values for $\fL(V)$, and therefore that $\fL(V)$ does not determine $\fT(V)$ in general. This resolves Question~2.

In many cases, however, Theorem A does allow us to deduce $\fT(V)$ from $\fL(V)$. Examining the 24 semilinear types listed in Table \ref{tab:24types} in more detail, we see that
\begin{itemize}
    \item in 12 out of the 24 types, the L-type $\fL(V)$ \emph{does} determine $\fT(V)$, however, and
    \item in the remaining 12 types, the L-type $\fL(V)$ determines two possibilities for $\fT(V)$.
\end{itemize}

\section*{Applications of the Classification of Types}
The significance of Theorem A is that it provides a means by which we can prove results for semilinear representations using results about linear representations. Let us now describe these, together with other applications of Theorem A.

\subsection*{Symmetric Groups}
The possible types $\fT(V)$
are often significantly restricted when one considers specific groups.
For example, suppose that $\bL / \bK$ has characteristic $0$
and the group-subgroup pair is $A_n \leq S_n$, an alternating/symmetric group with $n \geq 3$.

Given $V \in \Irr_{\bL}(A_n)$, using Theorem A one can show that the possible types $\fT(V)$ are
only 1, 2, 3 and 12 of Tables~\ref{tab:24types5} and \ref{tab:24types6} (Corollary \ref{cor:symgroups}).
In particular, the L-type $\fL(V)$ determines $\fT(V)$.

More importantly,  the only possible S-types $\fS(V)$ are \textbf{A} or \textbf{B}.
Hence, the semilinear representation theory of $S_n$ over $\bL$ is completely determined by $\Irr_{\bL}(A_n)$ and its action by $C_2$.
In particular, everything is completely determined by ordinary character theory. For $V \in \Irr_{\bL}(A_n)$:

\vspace{0.2em}

\noindent \hspace{0.7em} $\bullet$ If $V \cong \epsilon * V$, then $V$ extends uniquely to an irreducible semilinear representation of $S_n$ over $\bL$.

\noindent \hspace{0.7em} $\bullet$ If $V \not\cong \epsilon * V$, then $\Irr_H^G(V)$ is an irreducible semilinear representation of $S_n$ over $\bL$.

\noindent \hspace{0.7em} $\bullet$ All irreducible semilinear representations of $S_n$ over $\bL$ arise this way.

\subsection*{Cyclotomic Fields}
The possible types $\fT(V)$ are also often significantly restricted when one concentrates on a specific extension $\bL / \bK$.
For example, if $\bL$ is a finite field, then its Brauer group is trivial.
It follows that $\fC(V)$ and $\fS(V)$ are never \textbf{C} \cite[Cor.\ 9.4]{TAYSchur}.
Hence, by Theorem~A only 10 types $\fT(V)$ are possible, as seen by consulting Table \ref{tab:24types}.

Another important extension is $\bC / \bR$.
Here, by Dyson's Tenfold Way, only 10 types are possible. In Section~\ref{sect:recoverDyson10} we give a short proof of Dyson's Tenfold Way using Theorem~A, and show how in light of Theorem~A it follows directly from the classification of $C_2$-graded division algebras over $\bR$.

Another interesting class of extensions are the cyclotomic extensions.
Associated to the pair $(G, \sigma)$, we can consider the extension $\bQ(\zeta_e)_+ \leq \bQ(\zeta_e)$, where $e$ denotes the exponent of $H = \ker(\sigma)$ (assumed large enough so that this is a degree two extension). Our second main result is the following.

\begin{theoremB}
There are $15$ possibilities for the type $\fT(V)$ over $\bQ(\zeta_e)_+ \leq \bQ(\zeta_e)$, listed in Table~\ref{tab:15typesQe}.
\end{theoremB}

Examining in more detail the possible types, this allows us to construct all irreducible semilinear representations of $G$ over $\bQ(\zeta_e)_+ \leq \bQ(\zeta_e)$ using only the complex character table of $H$ and its action by $C_2$, together with the complex character table of $G$ and its action by $\Gal(\bQ(\zeta_{2e}) / \bQ(\zeta_e))$.

\section*{Other Results}
\subsection*{Antilinear Induction Theorem}

The crux of the proof of Theorem B relies on the following induction theorem for antilinear representations (semilinear representations for the extension $\bC / \bR$).

\begin{theoremC}
Suppose that $G$ is finite and $W \in \Irr_{\bC}^{\sigma}(G)$. Then there are subgroups $H_i \leq H$ and $G_j \leq G$ with each $G_j \not\subset H$, one-dimensional representations $V_i \in \Irr_{\bC}(H_i)$ and $W_j \in \Irr_{\bC}^{\sigma}(G_i)$, and $m_i, n_j \in \bZ$, such that
\[
    [V] = \sum_i m_i \cdot [\Ind_{H_i}^G(V_j)] + \sum_j n_j \cdot [\Ind_{G_j}^G(W_j)]
\]
in $K_0(\Rep_{\bC}^{\sigma}(G))$.
\end{theoremC}

We note that Theorem C requires considering induction from \emph{all} subgroups of $G$, not just the non-trivially graded subgroups $G_j$, and therefore blends the $\bC$-linear and antilinear theories.

Along with Brauer's induction Theorem, which concerns $\bC$-linear representations, Serre proved an induction theorem for $\bR$-linear representations, which states that any irreducible $\bR$-linear representation of a finite group is induced from representations of subgroups, which themselves are inflated from irreducible representations of cyclic or dihedral groups \cite[Lem.\ 73.19]{CR2}.

Theorem C follows from a slightly stronger statement we prove, Theorem \ref{thm:antilinearinductionthm}, which also applies to groups which also act trivially on $\bC$, and can be viewed as a simultaneous generalisation of both Brauer's Induction Theorem and Serre's Induction Theorem (see Section \ref{sect:recoveryofclassical}).

This generalisation furthermore explains the appearance of dihedral groups in Serre's Induction Theorem. Indeed, viewing each dihedral group $D_{2n}$ with it's natural action on $\bC$ with kernel $C_n$, dihedral groups are characterised as exactly those finite groups $G$ which admit a faithful $1$-dimensional antilinear representation (Lemma \ref{lem:1dimfftypeRisdihedral}).

From Theorem C, we also obtain the consequence that all antilinear representations admit models over $\bQ(\zeta_e)$ (Theorem \ref{thm:antilinearcyclotomicfinal}). Similarly, this can be viewed as a generalisation of the fact that any $V \in \Irr_{\bC}(H)$ admits a model over $\bQ(\zeta_e)$ and that any $X \in \Irr_{\bR}(H)$ admits a model over $\bQ(\zeta_e)_+$ \cite{PAS}.

\subsection*{Classification of Real Structures}

For a fixed extension $\bL / \bK$, to give a surjective homomorphism $\sigma \colon G \rightarrow \Gal(\bL / \bK)$ is the same as to give an subgroup $H \leq G$ of index two.

One can instead shift perspective from $G$ and begin with a group $H$, and consider \emph{Real structures} on $H$, which are pairs $(G,\iota)$ where $G$ is a group and $\iota \colon H \hookrightarrow G$ is an injective group homomorphism with $\iota(H)$ index two in $G$. For example, in this way one can view the $\bR$-linear representation theory of $H$ as just one among the possible antilinear representation theories $\Rep_{\bC}^{\sigma}(G)$.
One simply considers $G = H \times C_2$,  the trivial Real structure.

In Chapter \ref{chapter:realstructures}, we classify Real structures on a group $H$, first in a explicit and elementary fashion, and them cohomologically. We relate these approaches and provide algorithms for enumerating all Real structures on a given group $H$.

\subsection*{Classification of $C_2$-Modules}
The cohomological approach to the classification of Real structures involves the cohomology for the group $C_2$.
With this in mind, we prove a classification theorem for finitely-generated $\bZ[C_2]$-modules. We show that this category satisfies the Krull-Schmidt Theorem, meaning every object is a direct sum of indecomposable objects. We then describe the indecomposable objects in terms of certain finite linear data over the field $\bF_2$, which is algorithmically computable.

\subsection*{Antilinear Representations of Magnetic Groups}
Antilinear representations are of interest in Crystallography \cite{Bro, Wills}.
In Chapter~\ref{chap:physics} we consider \emph{magnetic groups}, which are Real structures on crystallographic groups.
We consider the 2-dimensional case, of Real structures on wallpaper groups.

We rederive the classification of 80 $2$-dimensional magnetic groups from our classification of Real structures and provide an in-depth analysis of their antilinear representation theory. In particular, we extend the theory of the Frobenius-Schur indicator to this setting, show that Dyson's Tenfold Way applies in this setting, and study the types of their irreducible representations. We calculate all relevant invariants of these groups and provide full algorithms in magma to justify the correctness of the tabulated results.

\section*{Loose Ends}
There are many aspects of our theory of types that deserve further explanation and additional exploration.
We touch upon some of them but not in sufficient depth.

\subsection*{Brauer-Wall Groups}
While we discuss them briefly in Appendix~\ref{appendix:gradedalgebras}, they deserve much more thorough investigation.
Consider a field $\bK$ of characteristic zero, its finite extension $\bF\geq \bK$ and an element of the Brauer-Wall group $A\in\BrW(\bK)$.
Can one describe all $C_2$-graded groups $G$ and quadratic extensions $\bL\geq \bK$ that ``realise'' $A$, in a sense of existence of
a surjective $C_2$-graded homomorphism $\bL_{\star \hspace{-0.1em}}[G] \twoheadrightarrow M_n (A)$?
What possible types $\fT (V)$ can the corresponding antilinear representation have? These questions appear intriguing.

\subsection*{Topological Groups}
The classical Dyson's Tenfold Way works equally well in the context of unitary representations of compact groups.
This obviously covers the case of finite groups.

Our work extends the Tenfold Way to finite-dimensional (discrete) representations of infinite groups. The next obvious step to try locally compact groups, and continuous representations.

Furthermore, we feel that the Tenfold Way is quite a general phenomenon that may work in several, more abstract contexts: more general topological groups, groups in tensor categories, topological vectors spaces, ind-vector spaces, etc.

\subsection*{Higher Dimensional Crystallographic Groups}
Chapter \ref{chap:physics} and Appendix~\ref{appendix_wallpaper} focus on the wallpaper groups, 2-dimensional crystallographic groups. Yet almost all of its provisions work in higher dimensions, and more generally for any virtually free abelian group. It would be interesting to see this theory applied in the 3-dimensional case, to classify antilinear representations of 3-dimensional magnetic groups.

\subsection*{Other Galois Groups}
Is the group $C_2$ special? What about more general Galois groups?
It feels like the next step is to consider groups $G$ acting on fields $\bL$ such that
the image of $G$ in the Galois group is cyclic of prime order. We believe that much of our theory can be done in this generality.

A much harder proposition is the question of more general Galois groups.
For instance, it would be interesting to investigate the case when the image of $G$ in the Galois group is a finite simple group.

\subsection*{Topology and Physics}
The authors are neither physicists nor topologists. Nevertheless, the theory developed here appears to have substantive consequences in both areas. We hope that our work will spark interest in specialists in those fields, and that these directions be explored further.

\section*{Acknowledgements and Author Statements}

\subsection*{Acknowledgements} The second author was supported by the departmental grant Progetto Sviluppo Dipartimentale - UNIPD PSDIP23O88 at the University of Padova, and the Heilbronn Institute of Mathematical Research at the University of Bristol.

\subsection*{AI Usage Statement} The key ideas and main results of this paper were all established without the use of LLMs, whilst free versions of LLMs available online were used for certain tasks, such as writing the associated magma code or searching the literature. The paper was entirely written without the use of LLMs.

%% file: chapter-realstructures.tex
\chapter{Real Structures on a Group}\label{chapter:realstructures}

The key object of our study are pairs $H\leq G$, where $G$ is a group and $H$ is an index 2 subgroup. In this chapter we consider such pairs from the point of view of $H$.

A \emph{Real structure} on $H$ is an injective group homomorphism $\phi \colon H \rightarrow G$,
where $\phi(H) \leq {G}$ is an index two subgroup.
An \emph{isomorphism of Real structures} 
    \[
        \psi \colon (\phi_1 \colon H \rightarrow {G}_1) \rightarrow (\phi_2 \colon H \rightarrow {G}_2)
    \]
    is a group isomorphism $\psi \colon {G}_1 \rightarrow {G}_2$ which restricts to an inner automorphism of $H$.
    We denote it by $(\phi_1 \colon H \rightarrow {G}_1) \simeq (\phi_2 \colon H \rightarrow {G}_2)$
    or simply ${G}_1 \simeq {G}_2$, when the maps are clear from the context. 
Our goal in this chapter is to classify Real structures on $H$ up to an isomorphism.    

We consider this definition of isomorphism, because inner automorphisms do not ``change'' representations of $H$, which is an important consideration for us. Instead of asking $\psi|_H$ to be an inner automorphism, however, a one might demand only that $\psi (H) = H$. We describe how the results from this chapter can also be used to understand this weaker notion of isomorphism in Section \ref{sect:weakequiv}.

\section{Elementary Approach}\label{sect:elementarydesc}
Let $H$ be a group. For $\alpha \in \Aut(H)$ and $h \in H$, consider a disjoint union of two copies of $H$
\[
G(\alpha, h) \coloneqq H \coprod Hw,
\]
where the second is written as $Hw$, using a formal symbol $w$. 
Define a binary operation $\cdot$ by
\[
\begin{array}{c|cc}
   \cdot & g_2 & g_2 w \\ \hline
   g_1   & g_1 g_2 & g_1 g_2 w \\
   g_1 w & g_1 \alpha(g_2) w & g_1 \alpha(g_2) h
\end{array}
\]
for all $g_1, g_2 \in H$. In particular, $w^2 = h$ and $w g = \alpha(g) w$ for any $g \in H$. Recall the adjoint map
\[
\Ad : H \rightarrow \Aut (H), \ \ \Ad (g): x \mapsto gxg^{-1} \, . 
\]
\begin{lemma}\label{lem:elementary1}
	The following statements hold for any group $H$. 
	\begin{itemize}
        \item $G(\alpha, h)$ defines a group if and only if
          \begin{equation} \label{Conditions_alpha_h}
            \Ad(h) = \alpha^2 \ \mbox{ and } \ \alpha(h) = h.
            \end{equation}
	\item If  $G(\alpha, h)$ is a group, $H \leq G(\alpha, h)$ is a Real structure on $H$.
	\item Suppose $H \leq {G}$ is a Real structure. Choose  $u \in {G} \setminus H$. Set
          $h \coloneqq u^2$ and $\alpha(h) \coloneqq \Ad(u)$. Then there is a canonical isomorphism
          \[ {G} \xrightarrow{\cong} G(\alpha, h), \ g\mapsto g, \ gu \mapsto gw\]
          of Real structures on $H$, which is the identity on $H$.
\end{itemize}
\end{lemma}

\begin{proof}
	Suppose first that $G(\alpha, h)$ defines a group. Then for any $g \in H$, by associativity
    \[
        (w \cdot w) \cdot g = w \cdot (w \cdot g)
    \]
    so $hg = \alpha^2(g)h$ and, hence, $\Ad(h) = \alpha^2$. Furthermore, $w$ has a left inverse $h^{-1} w$ and a right inverse $\alpha(h^{-1}) w$, so these must coincide and $\alpha(h) = h$.
    
    Conversely, suppose that $\Ad(h) = \alpha^2$ and $\alpha(h) = h$. The identity $e \in H$ provides an identity of $G(\alpha, h)$ because $\alpha(e) = e$. For associativity, suppose that $g_1 w^a, g_2 w^b, g_3 w^c \in G(\alpha, h)$. We want to show that
	\[
		g_1 w^a \cdot (g_2 w^b \cdot g_3 w^c) = (g_1 w^a \cdot g_2 w^b) \cdot g_3 w^c.
	\] 
	If $a = 0$ this is easy verified, so we may assume that $a = 1$. If $b = 0$, the the above identity becomes
	\[
	g_1 w \cdot (g_2 g_3 w^c) = (g_1 \alpha(g_2) w) \cdot g_3 w^c.
	\]
	which follows because $\alpha(g_2 g_3) = \alpha(g_2) \alpha(g_3)$. So we may also assume that $b = 1$. When $c= 0$ the identity becomes
	\[
	g_1 w \cdot (g_2 \alpha(g_3) w) = (g_1 \alpha(g_2) h) \cdot g_3
	\]
	or equivalently, simplifying both sides,
	\[
	g_1 \alpha(g_2) \alpha^2(g_3) h = g_1 \alpha(g_2) h g_3
	\]
	which follows from the fact that $\alpha^2(g_3) h = \Ad(h)(g_3) h = h g_3 h^{-1} h = h g_3$.
	
	When $c = 1$, the identity is
		\[
	g_1 w \cdot (g_2 \alpha(g_3) h) = (g_1 \alpha(g_2) h) \cdot g_3 w
	\]
	or equivalently, simplifying both sides,
	\[
	g_1 \alpha(g_2) \alpha^2(g_3) \alpha(h) w = g_1 \alpha(g_2) h g_3 w,
	\]
	which follows from the fact that $\alpha(h) = h$ and the same reasoning as above.
	
	For inverses, any element of $H$ has a two-sided inverse, so consider $g_1 w \in Hw$. An element $g_2 w \in Hw$ satisfies $(g_1 w)(g_2 w) = e$
        if and only if $g_1 \alpha(g_2) = h^{-1}$, so set $g_2 = \alpha^{-1}(g_1^{-1} h^{-1})$. Then $g_2 w$ defines a right inverse to $g_1 w$, which is also a left inverse because $G(\alpha, h)$ is a semigroup.
	
	The remaining statements are straightforward to verify.
\end{proof}
We can reinterpret this result in the following way. Pick any $\alpha \in \Aut(H)$ and $h \in H$ and
consider the following group, a quotient of the free product 
\[
G' (\alpha, h) \coloneqq \frac{H\ast \bZ}{N} = \frac{H \ast \langle w \rangle}{\langle \langle w^2=h, \  \ \forall g\in H  \ wg = \alpha (g) w \rangle \rangle} \; .
\]
This group comes with a natural map $H\rightarrow G' (\alpha, h), \ g\mapsto g$. 
\begin{cor}\label{cor:elementary} The following statements are equivalent.
\begin{enumerate}  
\item The natural map $H\rightarrow G' (\alpha, h)$ is a Real structure on $H$.
\item The natural map $H\rightarrow G' (\alpha, h)$ is injective.
\item The conditions~\eqref{Conditions_alpha_h} hold.
\end{enumerate} 
If these are satisfied, then the natural map $G'(\alpha, h) \rightarrow G(\alpha, h)$ is an isomorphism of Real structures.
\end{cor}

\begin{proof}
  Clearly (1) implies (2). To show that (2) implies (3), suppose that $H\rightarrow G' (\alpha, h)$ is injective. Then for any $g \in G$, $hgh^{-1} = w^2 g w^{-2} = \alpha^2(g)$ in $G'(\alpha, h)$, and thus by injectivity $\alpha^2 = \Ad(h)$. Similarly, $w$ commutes with $h$ in $G'(\alpha, h)$, and thus $hw = wh = \alpha(h)w$, so $h = \alpha(h)$.
  
  Finally, to see that (3) implies (1), note that if the conditions ~\eqref{Conditions_alpha_h} are satisfied then $G(\alpha, h)$ is a group by Lemma \ref{lem:elementary1}, and so from the presentation of $G'(\alpha,h)$ there is a surjective homomorphism $G'(\alpha,h) \rightarrow G(\alpha,h)$ which commutes with $\iota \colon H \rightarrow G'(\alpha,h)$ and $H \hookrightarrow G(\alpha,h)$. From the surjectivity, $[G'(\alpha,h) : \iota(H)] \geq [G(\alpha,h): H] = 2$, whereas from the presentation of $G'(\alpha,h)$, $[G'(\alpha,h) : \iota(H)] \leq 2$, so $[G'(\alpha,h) : \iota(H)] = 2$, which forces the surjection $G'(\alpha,h) \rightarrow G(\alpha,h)$ to be an isomorphism.
\end{proof}

Let us write $\fR(H)$ for the set of Real structures on $H$ up to isomorphism.
Lemma~\ref{lem:elementary1} yields the following description of this set. Set
\[
  \widetilde{\fR}(H) \coloneqq \{ (\alpha,h)\in \Aut(H)\times H \,\mid\, \mbox{equations }\eqref{Conditions_alpha_h}\mbox{ hold}\}.
\]
The Lemma~\ref{lem:elementary1} tells us that the map
\[
\widetilde{\fR}(H) \rightarrow \fR(H), \qquad (\alpha, h) \mapsto G(\alpha, h)
\]
is well-defined and surjective. It remains only to pinpoint when the Real structures $G(\alpha_1, h_1)$ and $G(\alpha_2, h_2)$ are isomorphic. 
We denote their corresponding elements $w$ by $w_1$ and $w_2$.

\begin{lemma}\label{lem:modifiedelementary2}
  Suppose that $G(\alpha_1, h_1)$ and $G(\alpha_2, h_2)$ define groups, and $a, b \in H$. Then the formula
  \[
  \phi_{a,b}(g) = a^{-1} g a, \qquad \phi_{a,b}(g w_1) = (a^{-1} g a b) w_2,
  \]
  for $g \in G$ defines an isomorphism $\phi_{a,b} \colon G(\alpha_1, h_1) \rightarrow G(\alpha_2, h_2)$ of Real structures if and only if
  \begin{equation}\label{copyH2action}
  a b  \alpha_2(b) h_2 a^{-1} = h_1
  \ \mbox{ and } \ 
    \Ad(a b  \alpha_2(a)^{-1}) \alpha_2 = \alpha_1 \, ,
\end{equation}
and any isomorphism from $G(\alpha_1, h_1)$ to $G(\alpha_2, h_2)$ is of this form.
\end{lemma}
\begin{proof}
  Using Corollary \ref{cor:elementary}, such a formula defines a well-defined group homomorphism if and only if the defining relations of $G'(\alpha_1, h_1)$ are preserved by $\phi_{a,b}$. A routine calculation in $G(\alpha_2, h_2)$ yields that these are equivalent to the relations (\ref{copyH2action}) in $H$. It is clear that once $\phi_{a,b}$ is well-defined, it is an isomorphism.

  Finally, consider a isomorphism of Real structures $\phi: G(\alpha_1, h_1) \to G(\alpha_2, h_2)$.
  On $H$, $\phi$ is inner, which determines $a$, and $b$ is determined by $\phi (w_1) = bw_2$.
  Inevitably, $\phi = \phi_{a,b}$.
\end{proof}

Consider the semidirect $H \ltimes H$, where $H$ acts on itself via $\Ad \colon H \rightarrow \Aut(H)$.

\begin{cor}\label{cor:HHaction}
Suppose that $G(\alpha_1, h_1)$ and $G(\alpha_2, h_2)$ define groups. Then $G(\alpha_1, h_1) \cong G(\alpha_2, h_2)$ as Real structures if and only if there exist $a,b \in H$ which satisfy the relations (\ref{copyH2action}).

In particular, the group $H \ltimes H$ acts on the set $\widetilde{\fR}(H)$ on the left by the formula 
  \begin{equation}\label{copyHHaction_explicit}
     (a,b) : (\alpha , h ) \mapsto
    (\Ad(a b  \alpha(a)^{-1}) \alpha,  a b  \alpha(b) h a^{-1}) 
  \end{equation}
  so that the set of Real structures on $H$ can be described as a quotient
  $\fR(H) = (H \ltimes H) \backslash \widetilde{\fR}(H)$. 
\end{cor}

We can give a simpler description of $\fR(H)$ involving only one copy of $H$.

\begin{lemma}\label{lem:elementary3}
  Suppose that $G(\alpha_1, h_1)$ and $G(\alpha_2, h_2)$ define groups. Then $G(\alpha_1, h_1) \cong G(\alpha_2, h_2)$ as Real structures if and only if there is some $s \in H$ with
  \begin{equation}\label{copyH2action2}
  s  \alpha_2(s) h_2  = h_1
  \ \mbox{ and } \ 
    \Ad(s) \alpha_2 = \alpha_1 \, .
\end{equation} 
\end{lemma}
\begin{proof}
  Given such an $s$, we can define an isomorphism $G(\alpha_1, h_1) \rightarrow G(\alpha_2, h_2)$ using Lemma \ref{lem:modifiedelementary2} taking $a = 1$ and $b = s$. Conversely, given any isomorphism $\phi \colon G(\alpha_1, h_1) \xrightarrow{\sim} G(\alpha_2, h_2)$, if $\phi|_H = \Ad(t)$, then we obtain another isomorphism
  \[
  \psi \coloneqq \phi \circ \Ad(t^{-1}) \colon G(\alpha_1, h_1) \xrightarrow{\sim} G(\alpha_1, h_1) \xrightarrow{\sim} G(\alpha_2, h_2)
  \]
  which is the identity on $H$.
  Suppose $\psi(w_1) = s w_2$.  Then $\psi = \phi_{1,s}$ by Lemma \ref{lem:modifiedelementary2}, and therefore $s$ satisfies the required relations, again by Lemma \ref{lem:modifiedelementary2}.
\end{proof}

\begin{thm}\label{cor:Haction}
The group $H$ acts on the set $\widetilde{\fR}(H)$ on the left by the formula 
  \begin{equation}\label{Haction_explicit}
     s : (\alpha , h ) \mapsto
    (\Ad(s) \alpha,  s  \alpha(s) h) 
  \end{equation}
  so that the set of Real structures on $H$ can be described as a quotient
  $\fR(H) = H \backslash \widetilde{\fR}(H)$. 
\end{thm}

\begin{remark}
Corollary \ref{cor:HHaction} and Corollary \ref{cor:Haction} tell us that the $H \ltimes H$ orbits on $\tilde{\fR}(H)$ are the same as those of $H$. It is therefore natural to ask: given $(\alpha_1, h_1)$ and $(\alpha_2, h_2)$ in the same $H \ltimes H$-orbit, what is an explicit element of $H$ which maps $(\alpha_1, h_1)$ and $(\alpha_2, h_2)$?

Concretely, we have that
\begin{equation}\label{copycopyH2action}
  a b  \alpha_2(b) h_2 a^{-1} = h_1
  \ \mbox{ and } \ 
    \Ad(a b  \alpha_2(a)^{-1}) \alpha_2 = \alpha_1 \, ,
\end{equation}
hold, and we want to find $s \in H$ with
\begin{equation}\label{copycopyH2action2}
  s  \alpha_2(s) h_2  = h_1
  \ \mbox{ and } \ 
    \Ad(s) \alpha_2 = \alpha_1 \, .
\end{equation} 

It turns out that $s \coloneqq a b  \alpha_2(a)^{-1}$ does the trick. The second relation is clearly satisfied, and one can verify the first:
\begin{align*}
         s  \alpha_2(s) h_2  = a b  \alpha_2(a)^{-1}\alpha_2(a b  \alpha_2(a)^{-1}) h_2  & = 
         a b  \alpha_2(b)  \alpha^2_2(a^{-1}) h_2  = \\
         = a b  \alpha_2(b)  h_2 (a^{-1}) h_2^{-1} h_2 &
         = a b  \alpha_2(b)  h_2 a^{-1} = h_1.
  \end{align*}  
  \end{remark}

\begin{eg}
We can also see that the orbits of the trivial pair $(1,1)$ (corresponding to the trivial Real structure $G(1,1) = H \times C_2$) under both $H$ and $H \ltimes H$ are the same as follows. The map
\[
\mu \colon H \rightarrow \widetilde{\fR}(H), \quad x \mapsto ( \Ad (x),x^2)\, ,
\]
is $H \ltimes H$-equivariant, when $H$ is given the action of $H \ltimes H$ by $(a,b) \cdot x \coloneqq abxa^{-1}$. Indeed, 
\[
\mu((a,b) \cdot x) = \mu(abxa^{-1}) = (abxa^{-1}, a (bx)^2 a^{-1}),
\]
which is equal to
\[
(a, b) \cdot \mu(x) = (a,b) \cdot (\Ad (x),x^2) = (\Ad(ab\Ad(x)(a)^{-1})\Ad(x), ab \Ad(x)(b)x^2a^{-1}).
\]

The action of $H \ltimes H$ on $H$ is transitive, and therefore the image of $H$ under $\mu$ is exactly the $H \ltimes H$-orbit of $(1,1) = \mu(1)$. Therefore, we see that the $H$-orbit of $(1,1)$, which is, by definition, the image of $\mu$, agrees with the $H \ltimes H$ orbit of $(1,1)$.
\end{eg}

\section{Cohomological Approach}
Let $H$ be a group. In Section \ref{sect:elementarydesc} we classified Real structures on $H$ using an elementary approach. In this section we consider an alternative point of view and describe how Real structures on $H$ are classified via cohomology and the Eilenberg-MacLane Theorem \cite{EML}. As a reference we follow \cite[\S 1]{DJR}.

We first note that a Real structure on $H$ is the same as an extension of $H$ by $C_2$.
Moreover, $(\phi_1 \colon H \rightarrow {G}_1) \simeq (\phi_2 \colon H \rightarrow {G}_2)$
if and only if the corresponding extensions are \emph{inner virtually equivalent},
in the sense that there is an isomorphism $\psi \colon {G}_1 \rightarrow {G}_1$ such that the diagram
\[\begin{tikzcd}
	1 & H & {{G}_1} & {C_2} & 1 \\
	1 & H & {{G}_2} & {C_2} & 1
	\arrow[from=1-1, to=1-2]
	\arrow["{\phi_1}", from=1-2, to=1-3]
	\arrow["\psi|_H"', from=1-2, to=2-2]
	\arrow[from=1-3, to=1-4]
	\arrow["\psi"', from=1-3, to=2-3]
	\arrow[from=1-4, to=1-5]
	\arrow["\id"', from=1-4, to=2-4]
	\arrow[from=2-1, to=2-2]
	\arrow["{\phi_2}", from=2-2, to=2-3]
	\arrow[from=2-3, to=2-4]
	\arrow[from=2-4, to=2-5]
\end{tikzcd}\]
commutes and $\psi|_H$ is an inner automorphism. We call these extensions {\em equivalent} if additionally $\psi|_H=\id$. Observe that, as in the proof of Lemma~\ref{lem:elementary3}, two extensions are equivalent if and only if they are inner virtually equivalent.

Let $\Out (H)$ denote the group of outer automorphisms of $H$, the quotient of $\Aut(H)$ by the subgroup of inner automorphisms.

\begin{defn}
  The \emph{coupling} $\theta$ of a Real structure $\phi \colon H \rightarrow {G}$ is the outer automorphism $\theta$ of $H$
  induced by conjugation by any $g\in G\setminus H$.
\end{defn}

As $\theta$ is viewed as an outer automorphism, this is independent of the choice of odd element $g \in G \setminus H$. Furthermore, isomorphic Real structures have the same coupling.
Note that, because $g^2\in H$ for any $g\in G$, we have $\theta^2 = 1$ in $\Out(G)$. In particular, the data of the coupling $\theta$ is equivalent to the data of a group homomorphism $\Phi_{\theta} \colon C_2 \rightarrow \Out(G)$.

\begin{defn}
We write $\fR_{\theta}(H)$ for the set of Real structures on $H$ with a fixed coupling $\theta$.
\end{defn}
By definition, there is a disjoint union decomposition
\begin{equation} \label{RealStr:RTheta}
\fR(H) = \bigsqcup_{\theta \in \Out(H), \ \theta^2 =1} \fR_{\theta}(H) \, . 
\end{equation}

\begin{defn}
  For $\theta \in \Out(H)$, let $Z_{\theta} \coloneqq Z(H)$, with $\bZ[C_2]$-module structure $\epsilon * z = \alpha(z)$ for any lift $\alpha \in \Aut(H)$ of $\theta$.
\end{defn}

Note that this is independent of the choice of $\alpha$, and for any such $\alpha$, $\alpha^2$ is inner, so $\alpha^2$ is trivial on $Z$ and the $\bZ[C_2]$-module structure is well-defined. 

From $\theta \in \Out(H)$, we therefore have the $C_2$-module $Z_{\theta}$ and may consider the cohomology groups $\HH^2(C_2, Z_{\theta})$ and $\HH^3(C_2, Z_{\theta})$. In order to describe $\fR_{\theta}$ in terms of these groups, we need a particular element $[\theta] \in \HH^3(C_2, Z_{\theta})$ associated to $\theta$, defined as follows.

We follow \cite[\S 7]{EML}. Choose a representative $\alpha \in \Aut(H)$ for $\theta$. Define the normalised function
\[ \tilde{\alpha} \colon C_2 \rightarrow \Aut(H), \ \ \tilde{\alpha}(1) = 1, \ \ \tilde{\alpha}(\epsilon) = \alpha \, .\]
As $\theta^2=1 \in \Out(H)$, $\alpha^2$ is inner, so we may choose some $t \in H$ such that $\alpha^2 = \Ad(t)$. This extends to another normalised function
\[ t: C_2 \times C_2 \rightarrow H, \ \ t(\epsilon, \epsilon) = t, \mbox{ and } t(x,y) = 1 \mbox{ if } 1\in \{x,y\}.\]
Then one obtains a normalised cocycle $f_\theta \in \HZ^3 (C_2 , Z_{\theta})$ defined by
\[
    f_\theta(x,y,z) \coloneqq \tilde{\alpha}(x)(t(y,z))t(x,yz)t(xy,z)^{-1}t(x,y)^{-1}.
\]
\begin{defn}
We set $[\theta] \coloneqq [f_\theta]\in \HH^3(C_2, Z_{\theta})$. 
\end{defn}
From Proposition \ref{prop:C2cohomology} there is a canonical identification between $\HH^3(C_2, Z_{\theta})$ defined in terms of normalised cocycles and the group $Z_{\theta}^1 / (\epsilon - 1) Z_{\theta}$, defined by sending a normalised cocycle $f$ to $f(\epsilon, \epsilon, \epsilon)$. In our case, the associated element is 
\begin{equation}
z_\theta \coloneqq  f_\theta(\epsilon , \epsilon , \epsilon ) =  \alpha (t)t^{-1} \in Z_{\theta}^1 / (\epsilon - 1) Z_{\theta}.
\end{equation}

With this element $[\theta]$, we can now state the Eilenberg-MacLane Theorem.
\begin{thm}[{\cite[Prop.\ 1.4]{DJR}}]\label{thm:EMLv1}
    Let $\theta \in \Out(H)$ such that $\theta^2=1$.
    \begin{itemize}
        \item If $[\theta] \neq 0$ in $\HH^3(C_2, Z_{\theta})$, then $\fR_{\theta}(H) = \emptyset$.
        \item If $[\theta] = 0$ in $\HH^3(C_2, Z_{\theta})$, then $\fR_{\theta}(H)$ has a free transitive group action of $\HH^2(C_2, Z_{\theta})$.
    \end{itemize}
\end{thm}

\section{Unification}\label{sect:unification} 
In this section we reconcile the elementary and cohomological approaches to describing $\fR(H)$.

\subsection{Existence}

We first consider the question of whether $\fR_{\theta}(H) \neq \emptyset$ from each point of view.

From Lemma~\ref{lem:elementary1} we have that $\fR_{\theta}(H) \neq \emptyset$ if and only if there exist $\alpha \in \Aut(H)$ and $h \in H$ with $\overline{\alpha} = \theta$, $\Ad(h) = \alpha^2$ and $\alpha(h) = h$. 

We can show directly that this is consistent with Theorem~\ref{thm:EMLv1}. 

\begin{lemma}\label{lem:deg3cohom}
    Suppose that $\theta \in \Out(H)$. Then $[f_{\theta}] = 0$ in $\HH^3(Z_{\theta})$ if and only if there are $\alpha \in \Aut(H)$ and $h \in H$ with $\overline{\alpha} = \theta$, $\Ad(h) = \alpha^2$ and $\alpha(h) = h$.
\end{lemma}

\begin{proof}
  By construction, $[f_{\theta}] = 0$ in $\HH^3(C_2, Z_{\theta})$ if and only if $z_{\theta} = 1$ in $Z_{\theta}^1 / (\epsilon - 1) Z_{\theta}$.

  If so, let $\alpha \in \Aut(H)$ and $t \in H$ be the elements involved in the construction of $z_{\theta}$, so $\overline{\alpha} = \theta$, $\Ad(t) = \alpha^2$, and $z_{\theta} = \alpha(t)t^{-1}$. As $z_{\theta} = 1$, there is some $z \in Z_{\theta}$ with $\alpha(t)t^{-1} = \alpha(z) z^{-1}$, and so we may set $h \coloneqq tz^{-1}$ and compute
    \[
        \alpha(h)h^{-1} = \alpha(t)\alpha(z)^{-1}t^{-1} z = 1 \ 
        \mbox{ and } \ \Ad(h) = \Ad(t) = \alpha^2.
    \]

    Conversely, given $\alpha \in \Aut(H)$ and $h \in H$ with $\overline{\alpha} = \theta$, $\Ad(h) = \alpha^2$ and $\alpha(h) = h$, we may use these elements to construct the element $z_{\theta}$, the value of $z_{\theta} \in Z_{\theta}^1 / (\epsilon - 1) Z_{\theta}$ being independent of this choice. For this pair, $z_{\theta} = \alpha(h)h^{-1} = 1$ on the nose, and so in particular $z_{\theta} = 1$ in $Z_{\theta}^1 / (\epsilon - 1) Z_{\theta}$.
\end{proof}

\subsection{Group Action}\label{sect:freetransgroupaction}

The second part of Theorem \ref{thm:EMLv1} says that when $\fR_{\theta}(H)$ is non-empty, there is a free and transitive group action of $\HH^2(C_2, Z_{\theta})$ on $\fR_{\theta}(H)$. Let us now also show this from the elementary point of view, and the description of $\HH^2(C_2, Z_{\theta})$ as $Z^{C_2}_{\theta} / (\epsilon + 1)Z_{\theta}$ of Proposition \ref{prop:C2cohomology}.

Suppose that $\theta \in \Out(H)$ with $\fR_{\theta}(H) \neq \emptyset$, and pick an element of $\fR_{\theta}(H)$, which by Lemma \ref{lem:elementary1} we may choose to represented by some $G(\alpha, h)$ with $\alpha \in \Aut (H)$ and $h \in H$ such that $\alpha^2 = \Ad(h)$ and $\alpha(h) = h$.

We can then define an action of $Z_{\theta}^{C_2} = \{z \in Z_{\theta} \mid \alpha(z) = z\}$ on $\fR_{\theta}(H)$ as follows.

\begin{defn}
  For $z \in Z_{\theta}^{C_2}$, let $z*G(\alpha, h)\coloneqq G(\alpha, zh)$.
\end{defn}

This is well-defined because $\Ad(zh) = \Ad(h)$ and $\alpha(zh) = \alpha(z) \alpha(h) = zh$. This is furthermore independent of the choice of representative: if $G(\alpha_1, h_1) \simeq G(\alpha_2, h_2)$, then by Lemma~\ref{lem:elementary3} there is some $s \in H$ such that $\Ad(s) \circ \alpha_2 = \alpha_1$ and $h_1 = s \alpha_2(s) h_2$, and so $h_1 z = s \alpha_2(s) h_2 z$. Hence, $G(\alpha_1, h_1 z) \simeq G(\alpha_2, h_2 z)$ by Lemma~\ref{lem:elementary3}.

\begin{lemma}\label{lem:freetransaction}
    The construction above defines a transitive group action of $Z_{\theta}^{C_2}$ on $\fR_{\theta}(H)$ with all stabilisers equal to $\{\alpha(z)z \mid z \in Z_{\theta}\}$. In particular, the induced action of $\HH^2(C_2, Z_{\theta})$ on $\fR_{\theta}(H)$ is transitive and free.
\end{lemma}

\begin{proof}
  It is clear that this defines a group action. For transitivity, consider $G(\alpha_1, h_1),G(\alpha_2, h_2)\in \fR_{\theta}(G)$. As $G(\alpha_1, h_1)$ and $G(\alpha_2, h_2)$ have the same coupling, $\Ad(t) \circ \alpha_2 = \alpha_1$ for some $t \in H$, and therefore by Lemma~\ref{lem:elementary3},
    \[
    G(\alpha_2, h_2) \;\simeq\; G(\Ad(t) \circ \alpha_2, t \alpha_2(t)h_2) \;\simeq\; G(\alpha_1, t \alpha_2(t) h_2) \, .
    \]
    In particular, we may assume that $\alpha_1$ and $\alpha_2$ are both equal to some common $\alpha \in \Aut(H)$. 

    We have that $\Ad(h_1h_2^{-1}) = \Ad(h_1) \Ad(h_2)^{-1} = 1$, hence $h_2 = h_1 z$ for some $z \in Z$. Furthermore, $h_2 = \alpha(h_2) = \alpha(h_1) \alpha(z) = h_2 \alpha(z)$, so $\alpha(z) = z$, hence $z \in Z_{\theta}^{C_2}$ and $z * G(\alpha, h_1) = G(\alpha, h_2)$ in $\fR_{\theta}(H)$.

    To compute the stabilisers, suppose that $z \in Z_{\theta}^{C_2}$ and $G(\alpha, h) \in \fR_{\theta}(H)$. Then $z *  G(\alpha, h) \simeq G(\alpha, h)$ if and only if $G(\alpha, zh) \simeq G(\alpha, h)$, which holds if and only if there is some $s \in Z$ with $zh = s \alpha(s) h$, or equivalently, if and only if $z \in \{\alpha(z)z \mid z \in Z_{\theta}\}$.
\end{proof}

\section{Examples}

In this section we give examples of some families of groups for which we can use the theory developed above to describe all Real structures.

\subsection{Abelian Groups}\label{sect:realstrabgroups}

In this section we assume that $H$ is abelian. Then $\Out(H) = \Aut(H)$, and for any $\alpha \in \Aut(H)[2]$, we always have that $\fR_{\alpha}(H) \neq \emptyset$, as one can always explicitly take $h = 1$ in Lemma \ref{lem:deg3cohom}. Furthermore, $Z(H) = H$, and so there is a canonical identification
\[
  \frac{\{h \in H \mid \alpha(h) = h\}}{\{\alpha(z)z \mid z \in H\}} \xrightarrow{\sim} \fR_{\alpha}(H), \qquad h \mapsto G(\alpha, h).
\]

\begin{eg}
Suppose that $H = \bZ$ is infinite cyclic. Then $\Aut(\bZ)[2] = \{\pm 1\}$. If $\alpha$ is the identity, then 
\[
\frac{\{h \in \bZ \mid \alpha(h) = h\}}{\{\alpha(z) + z \mid z \in \bZ\}} = \bZ / 2 \bZ.
\]
The trivial element corresponds to $\bZ \leq \bZ \times \bZ / 2 \bZ$, and the non-trivial element to $\bZ \leq \frac{1}{2}\bZ$.

If instead $\alpha(1) = -1$, then $\{z \in \bZ \mid -z = z\} = 0$, so $|\fR_{\alpha}(H)| = 1$ and there is a unique Real structure on $\bZ$ which induces this automorphism on $\bZ$. It is the infinite dihedral group $\bZ \leq D_{\infty}$.
\end{eg}

Now let us consider the case of finite cyclic groups $C_n$, for $n \geq 1$. In what follows we use the canonical identification
\[
(\bZ / n \bZ)^{\times} \xrightarrow{\sim} \Aut(C_{n}), \qquad t \mapsto [\alpha_t \colon x \mapsto x^t].
\]
We first treat the case of cyclic groups of prime power order.

\begin{eg}\label{eg:RealStructuresOddPrimePower}
Suppose that $H = C_{p^n} = \langle x \rangle$, where $n \geq 1$ and $p$ is an odd prime. From the isomorphism
\[
(\bZ / p^n \bZ)^{\times} \cong \bZ / (p-1) \bZ \times \bZ / p^{n-1} \bZ
\]
and the fact that $p$ is odd, we have that $|\Aut(C_{p^n})[2]| = 2$, and therefore $\Aut(C_{p^n})[2] = \{\alpha_1, \alpha_{-1}\}$. 

If $\alpha = \alpha_1$, then 
\[
\frac{\{h \in C_{p^n} \mid \alpha_1(h) = h\}}{\{\alpha_1(z)z \mid z \in C_{p^n}\}} = C_{p^n} / C_{p^n}^2 = 0,
\]
and therefore $|\fR_{\alpha_1}(C_{p^n})| = 1$, with unique Real structure $C_{p^n} \leq C_{p^n} \times C_2$ (which is isomorphic to the Real structure $C_{p^n} \leq C_{2 p^n}$). If instead $\alpha = \alpha_{-1}$, then
\[
\{h \in C_{p^n} \mid \alpha_{-1}(h) = z\} = \{h \in C_{p^n} \mid h^2 = 1\} = 0,
\]
and therefore $|\fR_{\alpha_{-1}}(C_{p^n})| = 1$, consisting of $C_{p^n} \leq D_{2p^n}$.
\end{eg}

There are more Real structures, however, when the prime $p$ is $2$.

\begin{eg}\label{eg:RealStructuresEvenPrimePower}
  Suppose that $H = C_{2^n} = \langle x \rangle$, where $n \geq 1$. In general, when $\alpha = \alpha_1$,
  \[
  \frac{\{h \in C_{2^n} \mid \alpha_1(h) = h\}}{\{\alpha_1(z)z \mid z \in C_{2^n}\}} = C_{2^n} / C_{2^n}^2,
  \]
  and therefore there are two Real structures, namely $G(\alpha_1, 1)$ and $G(\alpha_1, x)$, which are given by $C_{2^n} \leq C_{2^n} \times C_2$ and $C_{2^n} \leq C_{2^{n+1}}$ respectively.
  
  To determine $\Aut(C_{2^n})[2]$, there are now three cases, depending on $n$.
  
  $\bullet$ When $n = 1$, $\Aut(C_2) = 1$, and $\fR_{\alpha_1}(C_2)$ consists of the two Real structures $C_{2} \leq C_{2}\times C_2$ and $C_2 \leq C_4$ by the above.

  $\bullet$ When $n = 2$, $\Aut(C_4) = \{\pm 1\}$, and $\fR_{\alpha_1}(C_2)$ consists of the two Real structures $C_{4} \leq C_{4}\times C_2$ and $C_4 \leq C_8$ by the above. For $\fR_{\alpha_{-1}}(C_2)$, we have that
  \[
  \frac{\{h \in C_4 \mid \alpha_{-1}(h) = h\}}{\{\alpha_1(z)z \mid z \in C_4\}} = C_4^2 / 1,
  \]
  and so $|\fR_{\alpha_{-1}}(C_2)| = 2$, given by $C_4 \leq D_8$ and $C_4 \leq Q_8$.

  $\bullet$ When $n \geq 3$,
  \[
  (\bZ / 2^n \bZ)^{\times} \cong \bZ / 2 \bZ \times \bZ / 2^{n-2} \bZ,
  \]
  and so $|\Aut(C_{2^n})[2]| = 4$. The elements $\{\pm 1, 2^{n-1} \pm 1\}$ all square to $1$ in $(\bZ / 2^n \bZ)^{\times}$, and therefore $\Aut(C_{2^n})[2] = \{\alpha_t \mid t \in \{\pm 1, 2^{n-1} \pm 1\}\}$. 

  When $t = \pm 1$, 
  \begin{align*}
    \sR_{\alpha_1}(C_{2^n}) &= \{C_{2^n} \leq C_{2^n} \times C_2, \: C_{2^n} \leq C_{2^{n+1}}\}, \\
    \sR_{\alpha_1}(C_{2^n}) &= \{C_{2^n} \leq D_{2^{n+1}}, \: C_{2^n} \leq Q_{2^{n+1}}\},
  \end{align*}
  similarly to above, where $Q_{2^{n+1}}$ is the generalised quaternion group of order $2^{n+1}$. 
  
  When $t = 2^{n-1} + 1$, we have that 
  \[
  \frac{\{h \in C_{2^n} \mid \alpha_t(h) = h\}}{\{\alpha_t(z)z \mid z \in C_{2^n}\}} = C_{2^n}^2 / C_{2^n}^2 = 0
  \]
  and therefore $|\fR_{\alpha_t}(C_{2^n})| = 1$, given by 
\[
C_{2^n} = \langle x \rangle \leq \langle x , w \mid x^{2^n}, w^2, wxw^{-1} = x^{2^{n-1} + 1} \rangle \eqqcolon M_{2^{n+1}},
\]
the \emph{modular maximal-cyclic group} or \emph{modular group} of order $2^{n+1}$.

When $t = 2^{n-1} - 1$, we have that
\[
\frac{\{h \in C_{2^n} \mid \alpha_t(h) = h\}}{\{\alpha_t(z)z \mid z \in C_{2^n}\}} = C_{2^n}^{2^{n-1}} / C_{2^n}^{2^{n-1}} = 0
\]
and therefore $|\fR_{\alpha_t}(C_{2^n})| = 1$, given by 
\[
C_{2^n} = \langle x \rangle \leq \langle x , w \mid x^{2^n}, w^2, wxw^{-1} = x^{2^{n-1} - 1} \rangle \eqqcolon \SD_{2^{n+1}},
\]
the \emph{semidihedral group} of order $2^{n+1}$.
\end{eg}

We can use the prime power order cases above to describe the general case. We write
\[
\phi \colon \bZ / n\bZ \xrightarrow{\sim} \bZ / 2^r \bZ \times \bZ / p_1^{n_1} \bZ \times \cdots \times \bZ / p_k^{n_k} \bZ,
\]
for the canonical isomorphism induced by projection, where $n = 2^r p_1^{n_1} \cdots p_k^{n_k}$.

\begin{cor}
Suppose that $n = 2^r p_1^{n_1} \cdots p_k^{n_k}$ a prime factorisation, with $r \geq 0$, $k \geq 0$, and $n_i \geq 1$ for all $i$, and $C_n = \langle x \rangle$. Then for $t \in (\bZ / n\bZ)^{\times}[2]$, $\sR_{\alpha_t}(C_n)$ is as in one of the following cases.
\begin{enumerate}
  \item[\emph{($r = 0$)}] $|\Aut(C_n)[2]| = 2^k$, $|\fR_{\alpha_t}(C_n)| = 1$, and $\fR_{\alpha_t}(C_n)$ is described by
  \[
    \langle x , w \mid x^n, w^2, wxw^{-1} = x^t \rangle.
  \]
  \item[\emph{($r = 1$)}] $|\Aut(C_n)[2]| = 2^k$, $|\fR_{\alpha_t}(C_n)| = 2$, and $\fR_{\alpha_t}(C_n)$ is described by
  \[
    \langle x , w \mid x^n, w^2, wxw^{-1} = x^t \rangle, \qquad \langle x , w \mid x^n, w^2 = x^{p_1^{n_1} \cdots p_k^{n_k}}, wxw^{-1} = x^t \rangle.
  \]
  \item[\emph{($r = 2$)}] $|\Aut(C_n)[2]| = 2^{k+1}$, $|\fR_{\alpha_t}(C_n)| = 2$, and there are two cases:
  \begin{itemize}
    \item If $t \equiv 1 \!\! \mod 4$, $\fR_{\alpha_t}(C_n)$ is described by
    \[
    \langle x , w \mid x^n, w^2, wxw^{-1} = x^t \rangle, \qquad \langle x , w \mid x^n, w^2 = x^u, wxw^{-1} = x^t \rangle,
    \]
   where $\phi(u) = (1,0,0, ... , 0)$.
  \item If $t \equiv -1 \!\! \mod 4$, $\fR_{\alpha_t}(C_n)$ is described by
  \[
    \langle x , w \mid x^n, w^2, wxw^{-1} = x^t \rangle, \qquad \langle x , w \mid x^n, w^2 = x^u, wxw^{-1} = x^t \rangle,
  \]
  where $\phi(u) =  (2,0,0, ... , 0)$.
  \end{itemize}
  \item[\emph{($r \geq 3$)}] $|\Aut(C_n)[2]| = 2^{k+2}$, and there are three cases:
  \begin{itemize}
  \item If $t \equiv 1 \!\! \mod 2^r$, $|\fR_{\alpha_t}(C_n)| = 2$, and $\fR_{\alpha_t}(C_n)$ is described by
    \[
    \langle x , w \mid x^n, w^2, wxw^{-1} = x^t \rangle, \qquad \langle x , w \mid x^n, w^2 = x, wxw^{-1} = x^t \rangle,
  \]
  where $\phi(u) = (1,0,0, ... , 0)$.
  \item If $t \equiv -1 \!\! \mod 2^r$, $|\fR_{\alpha_t}(C_n)| = 2$, and $\fR_{\alpha_t}(C_n)$ is described by
    \[
    \langle x , w \mid x^n, w^2, wxw^{-1} = x^t \rangle, \qquad \langle x , w \mid x^n, w^2 = x^{2^{r-1}}, wxw^{-1} = x^t \rangle,
  \]
  where $\phi(u) = (2^{r-1},0,0, ... , 0)$.
  \item If $t \equiv 2^{r-1} \pm 1 \!\! \mod 2^r$, $|\fR_{\alpha_t}(C_n)| = 1$, and $\fR_{\alpha_t}(C_n)$ is described by
    \[
    \langle x , w \mid x^n, w^2, wxw^{-1} = x^t \rangle.
  \]
  \end{itemize}
\end{enumerate}
In particular, 
\[
|\fR(C_n)| = \twopartdef{2^{k+r}}{r \leq 2}{2^{k+2} + 2^{k+1}}{r \geq 3}.
\]
\end{cor}

\begin{proof}
  From the isomorphism
  \[
  C_n \xrightarrow{\sim} C_{2^r} \times C_{p_1^{n_1}} \cdots C_{p_k^{n_k}} 
  \]
  there is a square of isomorphisms
\[\begin{tikzcd}
	{(\bZ / n\bZ)^{\times}} & {(\bZ / 2^r \bZ)^{\times} \times (\bZ / p_1^{n_1} \bZ)^{\times} \times \cdots \times (\bZ / p_k^{n_k} \bZ)^{\times}} \\
	{\Aut(C_n)} & {\Aut(C_{2^n}) \times \Aut(C_{p_1^{n_1}}) \times \cdots \times \Aut(C_{p_k^{n_k}}) }
	\arrow["\sim", from=1-1, to=1-2]
	\arrow["\sim"', from=1-1, to=2-1]
	\arrow["\sim"', from=1-2, to=2-2]
	\arrow["\sim"', from=2-2, to=2-1]
\end{tikzcd}\] 
and therefore (in the notation of Section \ref{sect:cohomC2modules}) an isomorphism
\[
\frac{(C_n)^0}{(\epsilon + 1)C_n} \xrightarrow{\sim} \frac{(C_{2^r})^0}{(\epsilon + 1)C_{2^r}} \times \frac{(C_{p_1^{n_1}})^0}{(\epsilon + 1)C_{p_1^{n_1}}} \times \cdots \times \frac{(C_{p_k^{n_k}})^0}{(\epsilon + 1)C_{p_k^{n_k}}}.
\]
Then the classification follows from Examples \ref{eg:RealStructuresOddPrimePower} and \ref{eg:RealStructuresEvenPrimePower}, which describe the right-hand side. 
\end{proof}

\subsection{Groups of Odd Order}

\begin{lemma}\label{lem:oddorderdegree2automorphisms}
  Suppose that $H$ has odd order. Then the natural map
  \[
    \Aut(H)[2] \rightarrow \Out(H)[2]
  \]
  is an isomorphism.
\end{lemma}

\begin{proof}
  For injectivity, if $\alpha \in \Aut(H)[2]$ is inner, then $\Ad(g) = \alpha$ for some $g \in H$ and therefore $\alpha$ has odd order. We also have that $\alpha^2 = 1$, hence $\alpha = 1$. For surjectivity, suppose that $\theta^2 = 1$ in $\Out(H)$, and that $\alpha \in \Aut(H)$ is any lift of $\theta$. Then $\alpha^2 = \Ad(g)$ for some $g \in H$, and hence $\alpha^{2n} = 1$ for some odd $n \geq 1$. Then $\alpha^n$ maps to $\theta^n$, and $\theta^n = \theta$ as $n$ is odd, and $(\alpha^n)^2 = 1$, so $\alpha^n \in \Aut(H)[2]$.
\end{proof}

\begin{cor}\label{cor:oddorderRealstructures}
  Suppose that $H$ has odd order. Then
  \[
    \Aut(H)[2] \rightarrow \fR(H), \qquad \alpha \mapsto G(\alpha,1)
  \]
  is a bijection.
\end{cor}

\begin{proof}
  From Lemma \ref{lem:oddorderdegree2automorphisms} we have that the Real structures are partitioned by $\alpha \in \Aut(H)[2]$ rather than $\theta \in \Out(H)[2]$. As $\alpha^2 = 1$, we may always take $h = 1$ and form the semidirect product $G(\alpha, 1)$. By Lemma \ref{cohom_div_by2}, $\HH^2(C_2, Z_{\theta}) = 0$, hence this is the unique Real structure inducing this outer automorphism.
\end{proof}

In particular, we see that all Real structures
\[
1 \rightarrow H \rightarrow G \rightarrow C_2 \rightarrow 1
\]
are split when $H$ has odd order, and the Real structures are exactly $G(\alpha, 1)$ for $\Aut(H)[2]$.

We can upgrade this to the case when $H = H_1 \times H_2$, with $H_2$ of odd order and $\gcd(|H_1|,|H_2|) = 1$. In this case, for $\theta \in \Out(H)[2]$, we write $(\theta_1, \alpha_2)$ for the image of $\theta$ under the isomorphism
\[
\Out(H)[2] \xrightarrow{\sim} \Out(H_1)[2] \times \Out(H_2)[2] \xrightarrow{\sim} \Out(H_1)[2] \times \Aut(H_2)[2]
\]
induced from Lemma \ref{lem:oddorderdegree2automorphisms}.
\begin{cor}\label{cor:Realstructuresoddorderfactor}
  Suppose that $H = H_1 \times H_2$, with $H_2$ of odd order and $\gcd(|H_1|,|H_2|) = 1$. Then for any $\theta \in \Out(H)[2]$, $\fR_{\theta}(H) \neq \emptyset$ if and only if $\fR_{\theta_1}(H_1) \neq \emptyset$, in which case the natural map
  \begin{equation}\label{eqn:oddorderfactorisom}
  \fR_{\theta_1}(H_1) \rightarrow \fR_{\theta}(H), \qquad [G(\alpha_1, h)] \mapsto [G(\alpha_1 \times \alpha_2, (h,1))]
  \end{equation}
  is an isomorphism. In particular, this induces an isomorphism
  \[
    \sR(H_1) \times \Aut(H_2)[2] \xrightarrow{\sim} \sR(H).
  \]
\end{cor}

\begin{remark}
  The surjectivity of the isomorphism (\ref{eqn:oddorderfactorisom}) shows that for any Real structure $G$ on $H$, there is an element $w \in G \setminus H$ with $w^2 \in H_1$. In particular, this allows one to define an explicit inverse to (\ref{eqn:oddorderfactorisom}), by sending $G$ to the subgroup of $G$ generated by $H_1$ and any such $w \in G \setminus H$, which as $w^2 \in H_1$ will be a Real structure on $H_1$.
\end{remark}

\begin{proof}
  First, fix a lift $\alpha_1 \in \Aut(H_1)$ of $\theta_1 \in \Out(H_1)[2]$, and note that $\fR_{\theta_1}(H_1) \neq \emptyset$ if and only if there is some $h \in H_1$ with $\alpha_1(h) = h$ and $\alpha_1^2 = \Ad(h)$ by Lemma \ref{lem:elementary1}. This is equivalent to the existence of a pair $(h,h_2) \in H_1 \times H_2$ with $(\alpha_1(h), \alpha_2(h_2))= (h, h_2)$ and $(\alpha_1 \times \alpha_2)^2 = \Ad(h, h_2)$, i.e.\ that $\fR_{\theta}(H) \neq \emptyset$, as we can always take $h_2 = 1$ because $\alpha_2^2 = 1$.

  Is it straightforward to show that (\ref{eqn:oddorderfactorisom}) is well-defined an injective using Lemma \ref{lem:elementary1} and Lemma \ref{lem:elementary3}. For the surjectivity, any Real structure in $\fR_{\theta}(H)$ has the form $G(\alpha_1 \times \alpha_2, (h,h_2))$ for some $(h, h_2) \in H_1 \times H_2$ and lift $\alpha_1 \in \Aut(H_1)$ of $\theta_1 \in \Out(H_1)[2]$, by Lemma \ref{lem:elementary1} and Lemma \ref{lem:elementary3}. Now $\alpha_2(h_2) = h_2$ and $\alpha_2^2 = \Ad(h_2)$, and so $G(\alpha_2, h_2)$ is a Real structure on $H_2$. By Corollary \ref{cor:oddorderRealstructures}, $G(\alpha_2, h_2) \cong G(\alpha_2, 1)$, and therefore there is some $z \in Z(H_2)$ with $h_2 = z \alpha_2(z)$. In particular, taking $s = (1,z^{-1})$ in Lemma \ref{lem:elementary3}, we see that
\[
G(\alpha_1 \times \alpha_2, (h,h_2)) \cong G(\alpha_1 \times \alpha_2, (h,1))
\]
which shows that (\ref{eqn:oddorderfactorisom}) is surjective.
\end{proof}

\subsection{Centreless Groups} Immediately from Theorem \ref{thm:EMLv1}, we have the following.

\begin{cor}\label{cor:Realstructuresnocenter}
Suppose that $Z(H) = 1$. Then the coupling map
\[
    \fR(H) \rightarrow \Out(H)[2]
\]
is an isomorphism.
\end{cor}

For example, this applies to all non-abelian simple groups. Let us now examine what this means for the symmetric and alternating groups.

\begin{eg}
Suppose that $H=S_n$, for $n \geq 3$, so $Z(S_n) = 1$. When $n \neq 6$, $\Out(S_n) = 1$, and therefore the only Real structure is the trivial Real structure $S_n \leq S_n \times C_2$. When $n = 6$, however, there is an exotic outer automorphism making $|\Out(S_6)| = 2$, and so a second Real structure. This is explicitly given by the inclusion $S_6 \hookrightarrow \Aut(S_6)$.
\end{eg}

\begin{eg}
Suppose that $H=A_n$, for $n \geq 4$, so $Z(A_n) = 1$. When $n \neq 6$, $\Out(A_n)$ has order two, and therefore there are two Real structures on $A_6$ given by $A_n \leq A_n \times C_2$ and $A_n \leq S_n$. When $n = 6$, however, $\Out(A_6) \cong C_2 \times C_2$, and there are four Real structures. In addition to the previous two, there are two other groups which contain $A_6$ as an index two subgroup, namely $\text{PGL}_2(\bF_9)$ and the Mathieu group $M_{10}$.
\end{eg}

\begin{remark}
  When $Z(H) = 1$, from Corollary \ref{cor:Realstructuresnocenter} we know that given $\theta \in \Out(H)[2]$ there is a unique Real structure with coupling $\theta$. If $\alpha \in \Aut(H)$ is a lift of $\theta$, then by Lemma \ref{lem:elementary1} this can be represented by $G(\alpha, h)$ for some $h \in H$. From Lemma \ref{lem:elementary3}, if $G(\alpha, h_1) \cong G(\alpha, h_2)$ for some $h_1, h_2 \in H$, then there is some $s \in H$ with $\Ad(s) = 1$ and $s\alpha(s)h_1 = h_2$, and thus $s = 1$ and $h_1 = h_2$, because $Z(H) = 1$. This shows that for any $\alpha \in \Aut(H)$ which lifts $\theta$, there is a unique $h \in H$ such that $G(\alpha,h)$ forms a group. From the construction of $[\theta] \in \HH^3(C_2, Z_{\theta})$ and the proof of Lemma \ref{lem:deg3cohom}, we can see that this $h$ is the unique $h \in H$ such that $\Ad(h) = \alpha^2$, using again that $Z(H) = 1$.

  In particular, this gives an explicit inverse to the isomorphism of Corollary \ref{cor:Realstructuresnocenter}.
\end{remark}

\section{Weak-Equivalence}\label{sect:weakequiv}

Let $H$ be a group. In this section, we describe how to use the classification of Real structures on $H$ up to isomorphism to classify Real structures on $H$ up to weak equivalence. Here a weak equivalence from $(G_1, H)$ to $(G_2,H)$ is an group isomorphism $\psi \colon G_1 \rightarrow G_2$ for which $\psi(H) = H$. An isomorphism has the additional requirement that $\psi|_H \in \Aut(H)$ is an inner automorphism. We write $G_1 \cong G_2$ (resp.\ $G_1 \simeq G_2$) to denote that there is an isomorphism (resp.\ a weak equivalence) from $G_1$ to $G_2$.

Recall that we write 
\[
\widetilde{\fR}(H) \coloneqq \{ (\alpha,h)\in \Aut(H)\times H \,\mid\, \mbox{equations }\eqref{Conditions_alpha_h}\mbox{ hold}\},
\]
and view $\fR(H)$ as a quotient of this set. The following defines an action of $\Out(H)$ on $\fR(H)$, and shows that the orbits exactly the weak equivalence classes of Real structures on $H$.

\begin{prop}
  Let $H$ be a group. Then:
  \begin{enumerate}
    \item There is an action of $\Aut(H)$ on $\widetilde{\fR}(H)$ defined by setting
    \[
    \beta * (\alpha, h) \coloneqq (\beta  \alpha  \beta^{-1}, \beta(h))
    \]
    for $\beta \in \Aut(H)$ and $(\alpha, h) \in \widetilde{\fR}(H)$. 

  \item If $G(\alpha_1, h_1) \cong G(\alpha_2, h_2)$, then
  \[
  G(\beta  \alpha_1  \beta^{-1} , \beta(h_1)) \cong G(\beta  \alpha_2  \beta^{-1}, \beta(h_2)),
  \]
  and therefore this action induces an action of $\Aut(H)$ on $\fR(H)$.
  
  \item Two equivalence classes $[(\alpha_1,h_1)], [(\alpha_2, h_2)] \in \fR(H)$ are in the same orbit of $\Aut(H)$ if and only if $G(\alpha_1, h_1) \simeq G(\alpha_2, h_2)$.
  
  \item The induced action of $\Aut(H)$ on $\fR(H)$ is trivial on $\Inn(H)$, and thus induces an action of $\Out(H)$ on $\fR(H)$, for which the orbits correspond to weak equivalence classes of Real structures on $H$.
  \end{enumerate}
\end{prop}

\begin{proof}
  For point $(1)$, note that this is a well-defined action: if $\alpha^2 = \Ad(h)$ and $\alpha(h) = h$, then
  \[
  (\beta  \alpha  \beta^{-1})^2 = \beta  \alpha^2  \beta^{-1} = \beta  \Ad(h)  \beta^{-1} = \Ad(\beta(h))
  \]
  and $(\beta  \alpha  \beta^{-1})(\beta(h)) = \beta(h)$, so $\beta * (\alpha, h) \in \widetilde{\fR}(H)$.
  
  For point (2), if $G(\alpha_1, h_1) \cong G(\alpha_2, h_2)$ then there is some $s \in H$ with $s \alpha_2(s) h_2 = h_1$ and $\Ad(s) \alpha_2 = \alpha_1$ by Lemma \ref{lem:elementary3}. Then, we see that 
  \[
  \beta(s) (\beta \alpha_2 \beta^{-1})(\beta(s)) \beta(h_2) = \beta(h_1), \: \text{ and } \: \Ad(\beta(s)) (\beta\alpha_2\beta^{-1}) = \beta\alpha_1\beta^{-1}
  \]
  and therefore $G(\beta  \alpha_1  \beta^{-1} , \beta(h_1)) \cong G(\beta  \alpha_2  \beta^{-1}, \beta(h_2))$, again by Lemma \ref{lem:elementary3}.

  For point (3), suppose that $(\alpha, h) \in \widetilde{\fR}(H)$, and $\beta \in \Aut(H)$. Letting $w$ and $w'$ denote the formal symbols defining $G(\alpha,h)$ and $G(\beta \alpha \beta^{-1}, \beta(h))$ respectively, we define a map
  \[
  \psi \colon G(\alpha,h) \rightarrow G(\beta \alpha \beta^{-1}, \beta(h)) \: \text{ by } \: \psi(h) = \beta(h), \psi(hw) = \beta(h) w',
  \]
  which is clearly bijective by definition. Furthermore, we may check that $\psi(w)^2 = (w')^2 = \beta(h) = \psi(w^2)$, and
  \[
  \psi(wx) = \psi(\alpha(x) w) = \beta(\alpha(x))w' = w' \beta(x) = \psi(w) \psi(x)
  \]
  and therefore this is a well-defined group homomorphism by Corollary \ref{cor:elementary}.

  Conversely, suppose that there is an weak equivalence $\phi \colon G(\alpha_1,h_1) \rightarrow G(\alpha_2, h_2)$. We set $\beta \coloneqq \phi|_{H} \in \Aut(H)$, and claim that
  \[
  G(\alpha_2, h_2) \cong \beta * G(\alpha_1,h_1) = G(\beta \alpha_1 \beta^{-1},\beta(h_1)).
  \]
  By Lemma \ref{lem:elementary3}, this is equivalent to showing that there is some $s \in H$ with
  \[
    s \alpha_2(s) h_2 = \beta(h_1) \text{ and } \Ad(s) \alpha_2 = \beta \alpha_1 \beta^{-1}. 
  \]
  For this, set $s \coloneqq \phi(w_1)w_2^{-1} \in H$, so $\phi(w_1) = s w_2$. Then we may check that
  \[
    \beta(h_1) = \phi(w_1^2) = \phi(w_1)^2 = (s w_2)^2 = s \alpha_2(s)h_2
  \] 
  and that for any $x \in H$,
  \[
  \beta(\alpha_1(\beta^{-1}(x))) s w_2 = \phi(\alpha_1(\beta^{-1}(x)) w_1) = \phi(w_1 \beta^{-1}(x)) = \phi(w_1)\phi(\beta^{-1}(x)) = s w_2 x = s \alpha_2(x) w_2
  \]
  and therefore $\Ad(s) \alpha_2 = \beta \alpha_1 \beta^{-1}$.

  Finally, for point (4), note that for $\beta \in \Aut(H)$, the weak equivalence $\psi \colon G(\alpha,h) \rightarrow G(\beta \alpha \beta^{-1}, \beta(h))$ defined above is actually an isomorphism of Real structures if $\beta \in \Inn(H)$, as $\psi|_H = \beta$. Therefore the action of $\Aut(H)$ on $\Out(H)$ is trivial on $\Inn(H)$.
\end{proof}

\begin{eg}
The smallest example where $\Out(H) \backslash \fR(H)$ is a proper quotient of $\fR(H)$, so there are Real structures which are weak equivalent but not isomorphic, is when $H = C_2 \times C_2 = \langle x \rangle \times \langle y \rangle$. 

The outer automorphism group is
\[
\Out(C_2 \times C_2) = \Aut(\{x,y,xy\}) = S_3,
\]
and therefore the order two elements are $1, (12), (13), (23)$. As in Section \ref{sect:realstrabgroups}, because $C_2 \times C_2$ is abelian we have that for any $\alpha \in \Out(C_2 \times C_2)[2]$ the map 
\[
\frac{\{z \in C_2 \times C_2 \mid \alpha(z) = z\}}{\{z\alpha(z) \mid z \in C_2 \times C_2\}} \rightarrow \fR_{\alpha}(H), \qquad z \mapsto G(\alpha, z)
\]
is an isomorphism. In particular, if $\alpha \neq 1$, then $|\fR_{\alpha}(H)| = 1$. The action of $\Out(C_2 \times C_2)$ on $\Out(C_2 \times C_2)[2]$ by conjugation has two orbits, $\{1\}$ and $\{(12), (13), (23)\}$, and therefore all three Real structures with non-trivial coupling are weak equivalent.

When $\alpha = 1$ instead, $|\fR_{1}(H)| = 4$, and $\fR_{1}(H)$ is explicitly given by 
\[
G(1, 1), \: G(1,x), \: G(1,y), \: G(1,xy).
\]
The first is isomorphic as a group to $C_2^3$, and the remaining three to $C_2 \times C_4$. The stabiliser of $\alpha$ is now $S_3$, and the induced action of $S_3$ on $\fR_{1}(H)$ fixes $G(1,1)$ and identifies all remaining elements.

We therefore see that $H = C_2 \times C_2$ provides examples of distinct Real structures with different couplings which are weak equivalent, and distinct Real structures with the same coupling that are weak equivalent. Overall, $|\fR(H)| = 7$, whereas $|\Out(H) \backslash \fR(H)| = 3$. 
\end{eg}

We now consider how this action of $\Out(H)$ interacts with the decomposition \eqref{RealStr:RTheta} of $\fR(H)$ by coupling, which we use to compute $\Out(H) \backslash \fR(H)$ in practice.

\begin{lemma} \label{lemma:weak_eq}
  For $\eta \in \Out(H)$, $\eta(\fR_{\theta}(H)) = \fR_{\eta \theta \eta^{-1}}(H)$. Furthermore, for any representatives
  \[
  \{\theta_t \mid t \in \Out(H) \backslash \Out(H)[2]\} \subset \Out(H)[2]
  \]
  of $\Out(H) \backslash \Out(H)[2]$, where $\Out(H)$ acts on $\Out(H)[2]$ by conjugation, the natural map
  \[
  \bigsqcup_{t \in \Out(H) \backslash \Out(H)[2]} \Out_{\theta_t}(H) \backslash \fR_{\theta_t}(H) \longrightarrow \Out(H) \backslash \fR(H)
  \]
  is an isomorphism, where $\Out_{\theta_t}(H) \coloneqq C_{\Out(H)}(\theta_t)$ denotes the centraliser of $\theta_t$ in $\Out(H)$.
\end{lemma}

\begin{proof}
  For any $\theta, \eta \in \Out(H)$, any element of $\fR_{\theta}(H)$ is represented by some $G(\alpha, h)$ with $[\alpha] = \theta$ in $\Out(H)$, and therefore $\eta * G(\alpha, h)$ has coupling $\eta [\alpha] \eta^{-1} = \eta \theta \eta^{-1}$.

  The second claim follows from the following general fact. Suppose that $X$ is a left $G$-set, and $X = \sqcup_{i \in I} X_i$ is a disjoint union decomposition which is permuted by $G$, meaning that for any $g \in G$ and $i \in I$, $g(X_i) = X_{g(i)}$ for some $g(i) \in I$. Then this defines an action of $G$ on $I$, and for any choice of representatives $\{i_t \mid t \in G \backslash I\} \subset I$ of $G \backslash I$, the natural map
  \[
  \bigsqcup_{t \in G \setminus I} \Stab_{G}(i_t) \backslash X_{i_t} \rightarrow G \backslash X
  \]
  is a bijection.
\end{proof}

Given $\theta \in \Out(H)[2]$, if $\fR_{\theta}(H)$ is non-empty, then we can describe $\Out_{\theta}(H) \backslash \fR_{\theta}(H)$ explicitly as follows. Fix a lift $\alpha \in \Aut(H)$ of $\theta$, and let $h \in H$ be such that $G(\alpha,h)$ forms a group, or equivalently that $\Ad(h) = \alpha^2$ and $\alpha(h) = h$. For example, this may always be done by taking $t \in H$ with $\Ad(t) = \alpha^2$, finding $z \in Z(H)$ with $\alpha(t)t^{-1} = \alpha(z) z^{-1}$ and setting $h \coloneqq tz^{-1}$, as in Lemma \ref{lem:deg3cohom}.

From Lemma \ref{lem:freetransaction}, we see that the map
\[
\Phi_{\alpha, h} \colon Q_{\theta} \coloneqq Z_{\theta}^{C_2} / (\epsilon + 1)Z_{\theta} \rightarrow \fR_{\theta}(H), \qquad \Phi_{\alpha, h}(z) \coloneqq G(\alpha, z h)
\]
is a bijection. The following lemma describes an action of $\Out_{\theta}(H)$ on $Q_{\theta}$ which makes this bijection equivariant with respect to the action described above of $\Out_{\theta}(H)$ on $\fR_{\theta}(H)$.

\begin{lemma}\label{lem:torsoraction}
Let $\alpha\in\Aut(H)$ and $h\in H$ satisfy $\alpha^2=\Ad(h)$ and $\alpha(h)=h$. Let $\beta\in\Aut(H)$ satisfy
$[\alpha,\beta]=1$ in $\Out(H)$ and choose $s\in H$ such that $\beta\alpha\beta^{-1}=\Ad(s)\alpha$. Set
\[
z_\beta\coloneqq
(s\alpha(s))^{-1}\beta(h)h^{-1}.
\]
Then $z_\beta\in Z_\theta^{C_2}$, its class in $Q_{\theta}$ is independent of the choice of $s$, and
\[
[\beta]\ast[z]\coloneqq[\beta(z)z_\beta]
\]
defines a well-defined action of $\Out_\theta(H)$ on $Q_{\theta}$ making $\Phi_{\alpha,h}$ $\Out_\theta(H)$-equivariant.
\end{lemma}

\begin{proof}
  It is possible to check this directly by hand, but we instead proceed as follows, and first compute the induced action of $\Out_\theta(H)$ on $Q_{\theta}$ transported via $\Phi_{\alpha,h}$. For $[z] \in Q_{\theta}$, and $[\beta] \in \Out_{\theta}(H)$, 
  \[
  [\beta] * \Phi_{\alpha, h}([z]) = [\beta] * G(\alpha, zh) = G(\beta \alpha \beta^{-1}, \beta(z) \beta(h)),
  \]
  and therefore $[\beta] * [z] \in Q_{\theta}$ is the unique element of $Q_{\theta}$ such that there is some $s \in H$ with $\beta \alpha \beta^{-1} = \Ad(s)$ and $s \alpha(s) ([\beta] * [z]) h = \beta(z) \beta(h)$, by Lemma \ref{lem:elementary3}. Note that if $s' \in H$ is any element with $\beta \alpha \beta^{-1} = \Ad(s)$, then $s' = sy$ for some $y \in Z(H)$, and therefore
  \[
  (s' \alpha(s'))^{-1}\beta(z) \beta(h)h^{-1} = (y \alpha(y))^{-1} (s \alpha(s))^{-1}\beta(z) \beta(h)h^{-1} = (s \alpha(s))^{-1}\beta(z) \beta(h)h^{-1}
  \]
  in $Q_{\theta}$, and therefore 
  \[
  [\beta] * [z] = \beta(z) (s \alpha(s))^{-1} \beta(h) h^{-1} = \beta(z) z_{\beta}
  \]
  in $Q_{\theta}$ for \emph{any} $s \in H$ with $\beta \alpha \beta^{-1} = \Ad(s)$. As $[\beta] * [z] \in Q_{\theta}$, this tells us that $\beta(z) z_{\beta}$ and thus $z_{\beta}$ are in $Z_{\theta}^{C_2}$, $[z_{\beta}] \in Q_{\theta}$ is independent of $s$, and defining $[\beta]\ast[z]\coloneqq[\beta(z)z_\beta]$ gives exactly the transported action of $\Out_\theta(H)$ on $Q_{\theta}$ under $\Phi_{\alpha,h}$, as required.
\end{proof}

\subsection{Algorithm} In \cite{GITHUB}, we give an algorithm in Magma which, for a finite group $H$, returns all possible Real structures on $H$ up to isomorphism, and up to weak equivalence.

Below we describe this algorithm, which is based in the results of Section \ref{sect:unification} and Section \ref{sect:weakequiv}. This applies to any group, finite or not, for which each subroutine of the algorithm can be implemented.

\begin{algorithm}
\caption{Compute $\fR_\theta(H)$ and $\Out_{\theta}(H) \backslash \fR_{\theta}(H)$ for an order two element $\theta \in \mathrm{Out}(H)$}
\begin{algorithmic}[1]

\Require A group $H$, an order two element $\theta \in \mathrm{Out}(H)$, an element $\alpha \in \mathrm{Aut}(H)$ lifting $\theta$
\Ensure A list of elements of $H$, $\mathfrak{h}$, such that $\{G(\alpha, h) : h \in \mathfrak{h}\}$ is a set of representatives for $\fR_\theta(H)$, and a subset $\mathfrak{h}_0$ such that $\{G(\alpha, h) : h \in \mathfrak{h}_0\}$ is a set of representatives for $\Out_{\theta}(H) \backslash \fR_{\theta}(H)$

    \State Find $t \in H$ such that $\alpha^2 = \mathrm{Ad}(t)$
    \State Set $\lambda := \alpha(t) t^{-1} \in H$
    \If{there exists $z \in Z(H)$ such that $\lambda = \alpha(z) z^{-1}$}
        \State Fix such a $z$ and set $h := tz^{-1}$
        \State Compute the subgroups of $Z(H)$:
        \[
            Z^\alpha := \{ z \in Z(H) : \alpha(z) = z \}
        \]
        \[
            B^\alpha := \{ \alpha(z) \cdot z : z \in Z(H) \}
        \]
        \State Set $Q \coloneqq Z^\alpha / B^\alpha, \pi \colon Z^{\alpha} \twoheadrightarrow Q$ the quotient map, $\sigma \colon Q \rightarrow Z^{\alpha}$ a section
        \State Let $\mathfrak{B}$ be a set of coset representatives in $\Aut(H)$ of $\eta \in \Out_{\theta}(H)$
        \State For each $\beta \in \mathfrak{B}$, find $s \in H$ with $\beta \alpha \beta^{-1} = \Ad(s) \alpha$, and set
        \[
            z_{\beta} \coloneqq (s \alpha(s))^{-1}\beta(h) h^{-1}
        \]
        \State Set $\mathrm{orbitreps} \leftarrow []$, $\mathrm{tosort} \leftarrow Q$
        \While{$\mathrm{tosort} \neq \emptyset$}
          \State Pick $z \in \mathrm{tosort}$
          \State $\mathrm{orbit} \coloneqq \{\beta(z)z_{\beta} \mid \beta \in \mathfrak{B}\} \subset Q$
          \State Add $z$ to $\mathrm{orbitreps}$
          \State Remove $\mathrm{orbit}$ from $\mathrm{tosort}$
        \EndWhile

        \State $\mathfrak{h} \leftarrow [\sigma(z)h : z \in Q]$
        \State $\mathfrak{h}_0 \leftarrow [\sigma(z)h : z \in \mathrm{orbitreps}]$
    \Else
        \State $\mathfrak{h} \leftarrow [\,]$
        \State $\mathfrak{h}_0 \leftarrow [\,]$
    \EndIf

\Return $\mathfrak{h}, \mathfrak{h}_0$

\end{algorithmic}
\end{algorithm}

\begin{algorithm}
\caption{Compute $\fR(H)$ and $\Out(H) \backslash \fR(H)$}
\begin{algorithmic}[1]

\Require A group $H$
\Ensure A list of automorphisms $\mathfrak{a}$ of $H$, which represent $\Out(H)[2]$ in $\Aut(H)$, and a subset $\mathfrak{a}_0 \subset \mathfrak{a}$ which whose reductions in $\Out(H)[2]$ represent $\Out(H) \backslash \Out(H)[2]$. For each $\alpha \in \mathfrak{a}$, a list of elements $\mathfrak{h}_{\alpha}$ of $H$, such that $\{G(\alpha, h) : h \in \mathfrak{h}_{\alpha}\}$ is a set of representatives for $\fR_{[\alpha]}(H)$. For each $\alpha \in \mathfrak{a}_0$, a list of elements $\mathfrak{h}_{\alpha, 0}$ of $H$, such that $\{G(\alpha, h) : h \in \mathfrak{h}_{\alpha, 0}\}$ is a set of representatives for $\Out_{[\alpha]}(H) \backslash \fR_{[\alpha]}(H)$.

In particular,
\[
\{G(\alpha, h) : \alpha \in \mathfrak{a}, h \in \mathfrak{h}_{\alpha}\}
\]
is a set of representatives of $\fR(H)$, and 
\[
\{G(\alpha, h) : \alpha \in \mathfrak{a}_0, h \in \mathfrak{h}_{\alpha,0}\}
\]
is a set of representatives for $\Out(H) \backslash \fR(H)$.

\State Compute $\Out(H)[2]$
\State Pick a section $s \colon \Out(H) \rightarrow \Aut(H)$
\State Pick $\{\theta_t \colon t \in \Out(H) \backslash \Out(H)[2]\}$, representatives for $\Out(H) \backslash \Out(H)[2]$
\State $\mathfrak{a} \leftarrow [s(\theta) \mid \theta \in \Out(H)[2]]$
\State $\mathfrak{a}_0 \leftarrow [s(\theta_t) \mid t \in \Out(H) \backslash \Out(H)[2]]$

\State For each $\alpha \in \mathfrak{a}$, $\mathfrak{h}_{\alpha}$ is the $\mathfrak{h}$ provided by Algorithm 1
\State For each $\alpha \in \mathfrak{a}_0$, $\mathfrak{h}_{\alpha, 0}$ is the $\mathfrak{h}_0$ provided by Algorithm 1

\Return $\mathfrak{a}, \mathfrak{a}_0, [\mathfrak{h}_{\alpha} \mid \alpha \in \mathfrak{a}], [\mathfrak{h}_{\alpha, 0} \mid \alpha \in \mathfrak{a}_0]$

\end{algorithmic}
\end{algorithm}

%% file: chapter-algebras.tex
\chapter{Semilinear Representations and Associated Algebras}\label{chapter:algebras}

Suppose that $G$ is a group, $\bL$ is a field, and $\sigma \colon G \rightarrow \Aut(\bL)$ is group homomorphism from $G$ to the group of field automorphisms of $L$. For $g \in G$ and $\alpha \in \bL$, we often use the shorthand ${}^{g}\alpha\coloneqq \sigma (g) (\alpha)$.
Let $\bK \coloneqq \bL^{\!G}$ be the fixed-point subfield.

The central objects we are interested in are semilinear representations of $G$. In this chapter we review their basic theory.

\section{Semilinear Representations}

\begin{defn}
    A \emph{semilinear representation} of $G$ is a finite dimensional $\bL$-vector spaces $V$ with a group homomorphism $\rho \colon G \rightarrow \Aut_{\bK}(V)$ such that
    \[
    \rho(g)(\alpha \cdot v) = {}^{g}\hspace{-0.1em}\alpha \cdot \rho(g)(v) \ \mbox{  for all  } \ g\in G, \alpha\in\bL , v\in V.
    \]
    A \emph{morphism of semilinear representations} from $V$ to $W$ is an $\bL$-linear map $\phi \colon V \rightarrow W$ such that 
    \[
        \rho_W(g) \circ \phi = \phi \circ \rho_V(g) \ \mbox{  for all  } \ g\in G.
    \]
    We denote the category of such representations by $\Rep_{\bbL}^{\hspace{0.1em} \sigma}(G)$.
\end{defn}

\begin{eg}
    When $\sigma$ is trivial, so $\sigma(g)$ is the identity of $\bL$ for all $g \in G$, $\Rep_{\bbL}^{\hspace{0.1em} \sigma}(G)$ is simply the category of finite-dimensional linear representations of $G$ over $\bL$. In this case we write omit $\sigma$ from the notation and write $\Rep_{\bbL}(G)$ for $\Rep_{\bbL}^{\hspace{0.1em} \sigma}(G)$.
\end{eg}

\begin{defn}
An $\bL$-subspace $W$ of a semilinear representation $V$ is a \emph{subrepresentation} if
\[
        \rho(g)(w) \in W \ \mbox{  for all  } \ g\in G, w \in W.
\]
We say that $V$ is \emph{irreducible} if $V$ is non-zero and $0$ and $V$ are the only subrepresentations of $V$.
We write $\Irr_{\bbL}^{\hspace{0.1em} \sigma}(G)$ for the set of isomorphism classes of irreducible objects of $\Rep_{\bbL}^{\hspace{0.1em} \sigma}(G)$.
\end{defn}

We can perform certain constructions with semilinear representations, which generalise those of linear representations.

\begin{eg}\label{eg:tensorandhoms}
    Suppose that $V, W \in \Rep_{\bbL}^{\hspace{0.1em} \sigma}(G)$. Then
    \begin{itemize}
        \item $V \oplus W$ is a semilinear representation of $G$, where $g \in G$ acts by
    \[
        (v, w) \mapsto (\rho_V(g)(v), \rho_W(g)(w)).
    \]
        \item $V \otimes_{\bbL} W$ is a semilinear representation of $G$, where $g \in G$ acts by
    \[
       v \otimes w \mapsto \rho_V(g)(v) \otimes \rho_W(g)(w).
    \]
    \item $\Hom_{\bbL}(V, W)$ is a semilinear representation of $G$, where $g \in G$ acts by
    \[
        f \mapsto \rho_W(g) \circ f \circ \rho_V(g^{-1}).
    \]
    \end{itemize}
\end{eg}

\begin{eg}\label{eg:trivialrep}
We always have the trivial representation $\bL$, with its natural action of $G$. This is always irreducible and has endomorphism ring $\bK =\bL^{\!G}$.
\end{eg}

\begin{eg}\label{eg:permreps}
Suppose that $G$ acts on a finite set $S$. Then we can form the permutation representation $\bL[S]$, which has $\bL$-basis the elements of $S$, where each $g \in G$ acts via the $\sigma(g)$-linear extension of the permutation action on this basis.
\end{eg}

\begin{eg}\label{eg:inflation}
  Suppose that $f \colon G \twoheadrightarrow Q$ is a surjective morphism of groups, with $\ker(f) \subset \ker(\sigma)$. This kernel condition defines $\sigma$ also on $Q$. There is a pullback or inflation functor
\[
f^* (-) = \Inf_Q^G(-) \colon \Rep_{\bbL}^{\hspace{0.1em} \sigma}(Q) \rightarrow \Rep_{\bbL}^{\hspace{0.1em} \sigma}(G)
\]
defined by letting $G$ act on any $W \in \Rep_{\bbL}^{\hspace{0.1em} \sigma}(Q)$ by the composition of $f$ with $\rho_W \colon Q \rightarrow \Aut_{\bbK}(W)$. 

In particular, this applies to the quotient map $q \colon G \twoheadrightarrow G/N$ for any normal subgroup $N$ of $G$ which is contained in $\ker(\sigma)$. In general, if $V$ is a non-zero semilinear representation of $G$, then $N \coloneqq \ker(\rho_V)$ automatically acts trivially on $\bL$, so $N \subset \ker(\sigma)$, and we can view $V$ as $q^* W$, where $W$ is $V$ with its natural action of $G/N$.
\end{eg}

\section{Matrix Form of Semilinear Representations}\label{sect:matrixform}

We can describe the category $\Rep_{\bbL}^{\hspace{0.1em} \sigma}(G)$ more concretely in terms of matrices. We write $\Vect_{\bbL}$ for the category of finite-dimensional $\bL$-vector spaces.
\subsection{Semilinear Algebra}
For any $\gamma \in \Aut(\bL)$ and $V \in \Vect_{\bbL}$, we write ${}^{\gamma}V$ for $V$ as an abelian group but with $\bL$-vector space structure $\lambda * v = \gamma^{-1}(\lambda)v$. This defines a functor
\[
    {}^{\gamma}(-) \colon \Vect_{\bbL} \rightarrow \Vect_{\bbL},
\] 
which is the identity on morphisms and satisfies ${}^{\gamma}({}^{\mu}V) = {}^{\gamma \mu}V$.

Let $V, W \in \Vect_{\bbL}$ and fix bases $\sA$ and $\sB$ for $V$ and $W$ respectively, which induce an $\bL$-vector space identification
\[
    \End_{\bbL}(V,W) \rightarrow M_{m \times n}(\bL), \qquad [T \colon V \rightarrow W] \mapsto {}_{\sB}[T]_{\sA}.
\]
For any $\gamma \in \Aut(\bL)$, $\sA$ and $\sB$ also form bases for ${}^{\gamma}V$ and ${}^{\gamma}W$ respectively, and there is similarly an $\bL$-vector space identification
\[
    \End_{\bbL}({}^{\gamma}V,{}^{\gamma}W) \rightarrow M_{m \times n}(\bL), \qquad [T \colon V \rightarrow W] \mapsto {}_{\sB}[T]_{\sA},
\]
In fact, there is an equality of sets
\[
    \End_{\bbL}(V,W) = \End_{\bbL}({}^{\gamma}V,{}^{\gamma}W),
\]
and under this equality, the two identifications above are related by the commutativity of the diagram
\[\begin{tikzcd}
	{\End_{\bbL}(V,W)} & {M_{m \times n}(\bL)} \\
	{\End_{\bbL}({}^{\gamma}V,{}^{\gamma}W)} & {M_{m \times n}(\bL)}
	\arrow["\sim", from=1-1, to=1-2]
	\arrow["{=}"', from=1-1, to=2-1]
	\arrow["{{}^{\gamma}(-)}", from=1-2, to=2-2]
	\arrow["\sim", from=2-1, to=2-2]
\end{tikzcd}\]
where ${}^{\gamma}(-) \colon M_{m \times n}(\bL) \rightarrow M_{m \times n}(\bL)$ is the $\bK$-linear automorphism defined by applying $\gamma$ to each entry.

Concretely, if $A \in M_{m \times n}(\bL)$ represents an $\bL$-linear map $T \colon V \rightarrow W$, then ${}^{\gamma \!}A$ represents the $\bL$-linear map ${}^{\gamma}T \colon {}^{\gamma}V \rightarrow {}^{\gamma}W$ (which is set-theoretically simply the map $T$).

Suppose now that $V, W \in \Vect_{\bbL}$, and $S \colon V \rightarrow W$ is a \emph{$\gamma$-semilinear map}, meaning that
\[
S(\alpha \cdot v) = {}^{\gamma \hspace{-0.1em}}\alpha \cdot S(v) \ \mbox{  for all  } \ \alpha \in \bL, v \in V,
\]
or equivalently, that the map $S \colon {}^{\gamma} V \rightarrow W$ is $\bL$-linear. Let $B$ be the matrix associated to the $\bL$-linear map $S \colon {}^{\gamma} V \rightarrow W$ with respect to the bases $\sA$ and $\sB$, and write $[-]$ for the identifications $V \xrightarrow{\sim} \bL^n$, $W \xrightarrow{\sim} \bL^m$ induced by these choices of basis. Then
\[
[S(v)] = B \cdot {}^{\gamma}[v],
\]
where ${}^{\gamma}[v]$ is $\gamma$, applied to all the entries of $[v]$.

\subsection{Semilinear Matrix Representations}
Suppose now that $V \in \Rep_{\bbL}^{\hspace{0.1em} \sigma}(G)$. Then each $\rho(g) \colon V \rightarrow V$ is $\sigma_g$-semilinear and $\rho(g_1 g_2) = \rho(g_1) \circ \rho(g_2)$, so $V$ determines $\bL$-linear maps
\[
(\rho(g) \colon {}^{\sigma(g)}V \rightarrow V)_{g \in G},
\]
which satisfy $\rho(g_1 g_2) =  \rho(g_1) \circ {}^{\sigma(g_1)}\rho(g_2)$. Conversely, given $V \in \Vect_{\bbL}$, any such a family of $\bL$-linear maps defines a semilinear representation structure on $V$.

Let $n \coloneqq \dim_{\bbL}(V)$, and fix a basis of the $\bL$-vector space $V$. For any $g \in G$, write $A_g \in M_{n \times n}(\bL)$ for the matrix of $\rho(g) \colon {}^{\sigma(g)} V \rightarrow V$. The compatibility relations above in matrix form become
\begin{equation}\label{eqn:matrixrelation}
    A_{g_1g_2} = A_{g_1} \cdot {}^{\sigma(g_1) \!} A_{g_2} \ \mbox{  for all  } \ g_1,g_2 \in G.
\end{equation}
In particular, as $I_n = A_{g g^{-1}} = A_{g} \cdot {}^{\sigma(g)} A_{g^{-1}}$, each $A_g \in \GL_n(\bL)$. The following is direct to verify.

\begin{prop}\label{prop:matrixdesc}
Suppose that $V \in \Vect_{\bbL}$ has dimension $n \geq 1$, and fix a basis of $V$.

Then taking the matrix of a semilinear map defines a canonical bijection between semilinear representation structures $\rho$ on $V$ and families $(A_g)_{g \in G}$ in $\GL_n(\bL)$ satisfying (\ref{eqn:matrixrelation}). Under this bijection,
\begin{equation}\label{eqn:matrixevalrelation}
[\rho(g)v] = A_g \cdot {}^{\sigma(g)}[v] \ \mbox{  for all  } \ g \in G, v \in V.
\end{equation}
In particular, an $\bL$-subspace $W$ of $V$ is a subobject of $V$ in $\Rep_{\bbL}^{\hspace{0.1em} \sigma}(G)$ exactly when
\[
        A_g \cdot {}^{\sigma(g)}[w] \in W \ \mbox{  for all  } \ g\in G, w \in W.
\]

Suppose that $W \in \Vect_{\bbL}$ has dimension $m \geq 1$, and fix a basis of $W$. Let $V$ has semilinear representation structure corresponding to $(A_g)_{g \in G}$ in $\GL_n(\bL)$ and let $W$ have semilinear representation structure corresponding to $(B_g)_{g \in G}$ in $\GL_m(\bL)$.

Then an $\bL$-linear morphism $V \rightarrow W$ corresponding to $M \in M_{n \times m}(\bL)$ is a morphism in $\Rep_{\bbL}^{\hspace{0.1em} \sigma}(G)$ if and only if $M$ satisfies
\begin{equation}\label{eqn:matrixhomrelation}
B_g \cdot {}^{\sigma(g)} M = M \cdot A_g \ \mbox{  for all  } \ g \in G.
\end{equation}
\end{prop}

\begin{eg}
Suppose that $V \in \Rep_{\bbL}^{\hspace{0.1em} \sigma}(G)$ has dimension $n \geq 1$, and $\sA$ and $\sB$ are $\bL$-bases of $V$. Suppose that $(A_g)_{g \in G}$ and $(B_g)_{g \in G}$ represent $V$ with respect to $\sA$ and $\sB$ respectively, and that $P \in \GL_n(\bL)$ is the change of basis matrix from $\sA$ to $\sB$, meaning that $P \circ [-]_{\sA} = [-]_{\sB}$. Then
\[
B_g = P \cdot A_g \cdot {}^{\sigma(g)}(P^{-1}).
\]
\end{eg}

\begin{eg}
We can give matrix descriptions of the constructions of Examples \ref{eg:tensorandhoms}, \ref{eg:trivialrep} and \ref{eg:permreps}. The trivial representation has matrix representation by $(I_1)_{g \in G}$, for $I_1$ the identity of $M_1(\bL)$, and if $V$ and $W$ in $\Rep_{\bbL}^{\hspace{0.1em} \sigma}(G)$ have matrix representations $(A_g)_{g \in G}$ and $(B_g)_{g \in G}$ respectively, then 
\begin{itemize}
    \item $V \oplus W$ is represented by the block sum $(A_g \oplus B_g)_{g \in G}$,
    \item $V \otimes_{\bbL} W$ is represented by the Kronecker product $(A_g \otimes B_g)_{g \in G}$,
    \item $\Hom_{\bbL}(V,\bL)$ is represented by $((A_g^{-1})^T)_{g \in G}$.
\end{itemize}

If $G$ acts on a finite set $S$, then $\bL[S]$ has matrix representation by the same matrices as those which represent the linear representation $\bK[S]$, where each $A_g$ is a permutation matrix.
\end{eg}

\section{Semilinear Representations and Skew Group Rings}

Semilinear representations can also be viewed as modules over skew group rings. By $\bL[G]$ we denote the group algebra of $G$.

\begin{defn}
    The skew group ring $\bLG$ is $\bL[G]$ as an abelian group, but with multiplication $*$, where
    \begin{equation} \label{LG_product}
        (\alpha_1 \cdot g_1) * (\alpha_2 \cdot g_2) \coloneqq \alpha_1 ({}^{g_1}\alpha_2) \cdot g_1g_2 \ \mbox{  for all  } \ \alpha_1, \alpha_2 \in \bL, \ g_1,g_2 \in G,
    \end{equation}
    extended to all of $\bL[G]$ by $\bZ$-linearly.
\end{defn}

An alternative interpretation of $\bLG$ is that of the set of functions $\varphi : G \rightarrow \bL$ with finite support, where the multiplication formula~\eqref{LG_product} turns into a convolution
    \begin{equation} \label{LG_product2}
      (\varphi \ast \psi) (g) = \sum_{x\in G}  \varphi (x) \; {}^{x}\psi(x^{-1}g) ='' \int_G \varphi (x) \; {}^{x}\psi(x^{-1}g) d\mu_G(x)''.
     \end{equation}
In this section we discuss the properties of $\bLG$ as well as its completions.

\subsection{Categorical properties}
Note that $\bL$ is a subalgebra in $\bL_{\sigma} \hspace{-0.1em}[G]$, but it is not, in general, central.
The central subalgebra is only $\bK$. We therefore naturally view $\bL_{\sigma} \hspace{-0.1em}[G]$ as a $\bK$-algebra.
The following is well-known, and easy to verify.

\begin{prop}
  Suppose that $V \in \Vect_{\bbL}$. Then there is a canonical bijection between left $\bL_{\sigma} \hspace{-0.1em}[G]$-module structures on $V$ and semilinear structures on $V$, under which a morphism in $\Vect_{\bbL}$ is a morphism in $\Rep_{\bbL}^{\hspace{0.1em} \sigma}(G)$ if and only
  if it is a morphism in $\Mod_{\bL_{\sigma} \hspace{-0.1em}[G]}$.
\end{prop}

In particular, we see that $\Rep_{\bbL}^{\hspace{0.1em} \sigma}(G)$ is canonically identified with the full subcategory of $\Mod_{\bL_{\sigma} \hspace{-0.1em}[G]}$ where the underling $\bL$-vector space is finite-dimensional.

\begin{remark}
We note that $\bL$ is still a subalgebra of $\bL[G]$, which is not central whenever the action of $G$ on $\bL$ is non-trivial. Therefore there are two $\bL$-vector space structures on $\bLG$, given by left and right multiplications, and the precise algebraic structure of the skew group ring $\bLG$ is that of a Hopf $\bL|\bK$-algebra \cite{Nich85, Rum00}.
\end{remark}

Let us summarise the tensor properties of semilinear representations in the following, easy to prove, proposition.

\begin{prop} The operations in Example~\ref{eg:tensorandhoms} can be also applied to the category
$\Mod_{\bL_{\sigma} \hspace{-0.1em}[G]}$.
  Under these operations, the following properties hold.
  \begin{enumerate}
  \item $\bL_{\sigma} \hspace{-0.1em}[G]$ is a symmetric closed monoidal $\bK$-linear category.
  \item Its identity object is $\bL$, with the action given by $\sigma$, cf.\ Example~\ref{eg:trivialrep}. 
  \item $\Rep_{\bbL}^{\hspace{0.1em} \sigma}(G)$ is a symmetric  monoidal rigid $\bK$-linear category.
  \item The forgetful functor $\Rep_{\bbL}^{\hspace{0.1em} \sigma}(G)\rightarrow \Vect_{\bbL}$ is tensor $\bK$-linear and faithful, making
    $\Rep_{\bbL}^{\hspace{0.1em} \sigma}(G)$ a Tannakian category \cite{DelMil}.
  \item  The Tannakian category $\Rep_{\bbL}^{\hspace{0.1em} \sigma}(G)$ is neutral if and only if the action $\sigma$ is trivial 
  \end{enumerate}
\end{prop}

We also note here that Maschke's Theorem extends to semilinear representations \cite[Lem.\ 1.3]{KUN}.

\begin{prop}\label{prop:Maschke}
    Suppose $G$ is finite. Then $\Rep_{\bbL}^{\rtimes}(G)$ is semisimple if and only if $|H| \in \bL^{\times}$.
\end{prop}

\subsection{Associated linear representation categories}

There are linear representation categories that interact with $\Rep_{\bbL}^{\hspace{0.1em} \sigma}(G)$. 

\begin{defn}\label{def:H}
We write $H \coloneqq \ker(\sigma \colon G \rightarrow \Aut_{\bbK}(\bL))$.
\end{defn}

With $H$ defined as above, there are $\bK$-algebra inclusions
\begin{equation}\label{eqn:diamondrings}\begin{tikzcd}
	& {\bK[G]} \\
	{\bL[G]} & {\bK[H]} & {\bL_{\sigma \hspace{-0.1em}}[G]} \\
	& {\bL[H]}
	\arrow[hook, from=1-2, to=2-1]
	\arrow[hook, from=1-2, to=2-3]
	\arrow[hook, from=2-2, to=1-2]
	\arrow[hook, from=2-2, to=3-2]
	\arrow[hook, from=3-2, to=2-1]
	\arrow[hook, from=3-2, to=2-3]
\end{tikzcd}
\end{equation}
which induce induction and restriction functors between categories
\begin{equation}\label{eqn:diamond}
\begin{tikzcd}
	& {\Rep_{\bbK}(G)} \\
	{\Rep_{\bbL}(G)} & {\Rep_{\bbK}(H)} & {\Rep_{\bbL}^{\hspace{0.1em}\sigma}(G)} \\
	& {\Rep_{\bbL}(H)}
	\arrow[tail reversed, from=1-2, to=2-1]
	\arrow[tail reversed, from=1-2, to=2-3]
	\arrow[tail reversed, from=2-2, to=1-2]
	\arrow[tail reversed, from=2-2, to=3-2]
	\arrow[tail reversed, from=3-2, to=2-1]
	\arrow[tail reversed, from=3-2, to=2-3]
\end{tikzcd}
\end{equation}

For example, in terms of matrices, the induction functor
\[
    \bL_{\sigma \hspace{-0.1em}}[G] \otimes_{\bK[G]} - \colon \Rep_{\bbK}(G) \rightarrow \Rep_{\bbL}^{\hspace{0.1em}\sigma}(G)
\]
simply takes a linear matrix representation $\rho \colon G \rightarrow \GL_n(\bK)$ of $G$ and view this as a semilinear representation under the inclusion $\GL_n(\bK) \hookrightarrow \GL_n(\bL)$. Indeed, the matrices relations (\ref{eqn:matrixrelation}) are satisfied, as all matrices are fixed by $\Aut(\bL / \bK)$.

We will use these linear representation categories to understand $\Rep_{\bbL}^{\hspace{0.1em}\sigma}(G)$.   

\subsection{Bohr completion} \label{Bohr_completion}

In this section, we let $A$ denote a $\bK$-algebra and consider a certain completion of $A$ that allows one to consider simple $A$-modules in terms of idempotents. We will later take $A$ to be one of the algebras of Diagram (\ref{eqn:diamondrings}).

Consider the forgetful functor
\[
\sF : \Mod^{\fd,\semisimple}_A \to \Vect_{\bK}
\]
from the category of finite-dimensional semisimple $A$-modules.
\begin{defn}
  The {\em Bohr completion} of $A$ is $\widehat{A} \coloneqq \End (\sF)$.
\end{defn}

Let $\Irr(A)$ denote the set of isomorphism classes of finite-dimensional simple $A$-modules. For $V \in \Mod^{\fd, \semisimple}_A$, set $E(V) \coloneqq \End_A V$, and let $I(V)$ denote the annihilator of $V$ in $A$.

\begin{lemma} \label{Bohr2}
  Let $\sR$ be a set of representatives of $\Irr(A)$. Then the natural maps
\[ \widehat{A} \xlongrightarrow{\sim}  \prod_{V\in \sR} \End_{E(V)}(V) \xlongleftarrow{\sim} \prod_{V\in \sR} A/I(V)\]
are isomorphisms.
\end{lemma}

\begin{proof}
  The right-hand map is an isomorphism: the natural map $A \rightarrow \End_{E(V)}(V)$ has kernel $I(V)$, and is surjective when $V$ is irreducible by, for example, the Jacobson Density Theorem.

  For the left-hand map, note that an element of $\widehat{A}$ is a collection $(x_W)$ of $x_W \in \End_{E(W)}(W)$ for each $W \in \Mod^{\fd,\semisimple}_A$, such that $x_{W'} \circ f = f \circ x_{W}$ for any $A$-linear $f \colon W \rightarrow W'$.
  
  If $x_V = 0$ for any $V \in \sR$, then for any $W \in \Mod^{\fd,\semisimple}_A$, $W = W_1 \oplus \cdots \oplus W_n$ for some irreducible submodules $W_i$. Choosing isomorphisms $f_i \colon V_i \xrightarrow{\sim} W_i$ for some $V_i \in \sR$,
  \[
  0 = (f_i \circ \iota_i) \circ x_{V_i} = x_W \circ (f_i \circ \iota_i),
  \]
  where $\iota_i \colon V_i \rightarrow W$ is the inclusion. Then $x_W = 0$, the image of $\oplus_i (f_i \circ \iota_i)$ being all of $W$.
  
  Suppose conversely that $(y_V)_{V \in \sR}$ is any collection of $y_V \in \End_{E(V)}(V)$, subject to no compatibility requirements. To define an extension to a compatible family $(x_W)$ as above, we first set, for $V_1, \ldots , V_n \in \sR$,
  \[
  x_{\underline{V}} \coloneqq x_{V_1} \oplus \cdots \oplus x_{V_n}, \, \text{ where } \, \underline{V} \coloneqq V_1 \oplus \cdots \oplus V_n.
  \]

  Then for any $W \in \Mod^{\fd,\semisimple}_A$ and isomorphism $\phi \colon W \xrightarrow{\sim} \underline{V}$ to such a $\underline{V}$, set
  \[
    x_W^{\phi} \coloneqq \phi^{-1} \circ x_{\underline{V}} \circ \phi.
  \]

 All that remains to show is reduced to the following property: for any morphism $f \colon W \rightarrow X$ in $\Mod^{\fd,\semisimple}_A$, and isomorphisms $\phi \colon W \xrightarrow{\sim} \underline{V}$, $\psi \colon X \xrightarrow{\sim} \underline{U}$, to direct sums of objects of $\sR$, then $f \circ x_X^{\psi} = x_W^{\phi} \circ f$.
  
  Indeed, this shown, then given any $W \in \Mod^{\fd,\semisimple}_A$, we can find some such direct sum $\underline{V}$, and isomorphism $\phi \colon W \xrightarrow{\sim} \underline{V}$, and may set $x_W \coloneqq x_W^{\phi}$. This is independent of the choice of $\underline{V}$ and $\phi$, taking $X = W$ and $f = \id_W$ in the property above, which also tells us that the resulting family $(x_W)_W$ is a well-defined element of $\widehat{A}$.

  To establish the claimed property, note that in the diagram 
\[\begin{tikzcd}
	{\underline{V}} & W & X & {\underline{U}} \\
	{\underline{V}} & W & X & {\underline{U}}
	\arrow["{x_{\underline{V}}}", from=1-1, to=2-1]
	\arrow["\phi"', from=1-2, to=1-1]
	\arrow["f", from=1-2, to=1-3]
	\arrow["{x_{W}^{\phi}}", from=1-2, to=2-2]
	\arrow["\psi", from=1-3, to=1-4]
	\arrow["{x_{X}^{\psi}}", from=1-3, to=2-3]
	\arrow["{x_{\underline{U}}}", from=1-4, to=2-4]
	\arrow["\phi"', from=2-2, to=2-1]
	\arrow["f", from=2-2, to=2-3]
	\arrow["\psi", from=2-3, to=2-4]
\end{tikzcd}\]
  the commutativity of the central square is equivalent to the commutativity of the outer square. Therefore, it is sufficient to show that for any $A$-linear $g \colon \underline{V} \rightarrow \underline{U}$, $g \circ x_{\underline{V}} = x_{\underline{U}} \circ g$. To see this, label write $\underline{U} = U_1 \oplus \dots \oplus U_m$ for some $U_j \in \sR$. Letting $\iota_i \colon V_i \hookrightarrow \underline{V}$ and $\pi_j \colon \underline{U} \rightarrow U_j$ denote the canonical inclusions and projections, it is sufficient to show that for any $i = 1, ... n$ and $j = 1, ... m$,
  \[
  \pi_j \circ g \circ x_{\underline{V}} \circ \iota_i = \pi_j \circ x_{\underline{U}} \circ g \circ \iota_i.
  \]
  We have that $x_{\underline{V}} \circ \iota_i = x_{V_i} \circ \iota_i$, and $\pi_j \circ x_{\underline{U}} = x_{W_j}\circ \pi_j$, so this is equivalent to the claim that
  \[
  (\pi_j \circ g \circ \iota_i) \circ x_{V_i} = x_{W_j} \circ (\pi_j \circ g \circ \iota_i).
  \]
  If $V_i \neq W_j$, then $V_i \not\cong W_j$ and both maps are zero. If instead $V_i = W_j$, then $\pi_j \circ g \circ \iota_i \in E(V_i)$, and the claim then follows from the fact that $x_{V_i} = x_{W_j} \in \End_{E(V_i)}(V_i)$.
\end{proof}

\begin{remark}
  Lemma~\ref{Bohr2} yields an alternative description of the canonical homomorphism $\iota_A \colon A\to\widehat{A}$, as the map
  \[
  A \rightarrow \prod_{V\in \sR} A/I(V), \qquad a \mapsto (a+I(V))_{V \in \sR}.
  \]
\end{remark}

\begin{eg}
Suppose that $A = \bK[T]$ and $\bK$ is algebraically closed. Then for the theory of the Jordan normal form, there is a canonical identification of $\iota_A \colon A \rightarrow \widehat{A}$ with the product of evaluation homomorphisms
\[
\prod_{\lambda \in \bK} \text{ev}_{\lambda} \colon \bK[T] \hookrightarrow \prod_{\lambda \in \bK} \bK.
\]
\end{eg}

\begin{eg}
The map $\iota_A \colon A \rightarrow \widehat{A}$ doesn't always have to be injective, however. If $A$ is a finite dimensional $\bK$-algebra, then there are finitely many simple modules, and each simple module is finite dimensional, hence $\iota_A \colon A \rightarrow \widehat{A}$ has kernel $J(A)$, the Jacobson radical of $A$, and $\iota_A \colon A \rightarrow \widehat{A}$ is canonically identified with the natural map $A \rightarrow A / J(A)$ to the semisimplification $A / J(A)$ of $A$.
\end{eg}

Not all $\bK$-algebra homomorphisms $\phi : A \rightarrow B$ extend to Bohr completions. This is true, however, in the following setting. Suppose that the restriction functor
\[
(-)|_A \colon \Mod_B \to \Mod_A
\]
preserves finite-dimensional semisimple modules. Then there is a canonical $\bK$-algebra homomorphism $\widehat{\phi} : \widehat{A} \rightarrow \widehat{B}$, defined by
\[
\widehat{\phi}((x_W)_W) = (y_U)_U, \qquad y_U \coloneqq x_{U|_A}.
\]
This is compatible with the canonical homomorphisms $\iota_A \colon A\to\widehat{A}$ and $\iota_B \colon B\to\widehat{B}$ and $\phi \colon A \rightarrow B$.

We can characterise when $\widehat{\phi} : \widehat{A} \rightarrow \widehat{B}$ is injective as follows.

\begin{lemma} \label{hom_cond}
  Suppose that the functor $(-)|_A \colon \Mod_B \to \Mod_A$ preserves semisimplicity and finite-dimensionality. Suppose also that the induction functor 
  \[ \Mod_A \to \Mod_B, \quad M \mapsto B\otimes_{A}M \]
  preserves semisimplicity and finite-dimensionality. Then $\widehat{\phi} : \widehat{A} \rightarrow \widehat{B}$ is injective if and only if $B \otimes_A V \neq 0$ for any $V \in \Irr(A)$.
\end{lemma}
\begin{proof}
Let $\sR$ be a set of representatives of $\Irr(A)$. Suppose first that the induction functor takes finite-dimensional simple $A$-modules to non-zero (semisimple) $B$-modules, and that $\widehat{\phi}((x_W)_W) = 0$, meaning that $x_{U|_A} = 0$ for any $U \in \Mod^{\fd,\semisimple}_B$. For any $V \in \sR$, there is a canonical $A$-linear map $\iota \colon V \rightarrow B \otimes_A V$ defined by $\iota(v) = 1 \otimes v$. By assumption $B \otimes_A V \neq 0$, and therefore $\iota(V) \neq 0$, hence $\iota$ is injective, by the simplicity of $V$. From the fact that $(x_W)_W \in \widehat{A}$ defines a natural transformation,
\[
\iota \circ x_V = x_{(B \otimes_A V)|_A} \circ \iota = 0,
\]
and therefore $x_V = 0$ from the injectivity of $\iota$. Then $(x_W)_W = 0$, by Lemma \ref{Bohr2}.

Suppose conversely that $\widehat{\phi} : \widehat{A} \rightarrow \widehat{B}$ is injective, and that, for a contradiction, $B \otimes_A V_0 = 0$ for some $V_0 \in \sR$. Consider the element $(x_W)_W \in \widehat{A}$ mapping to $(\delta_{V,V_0})_{V \in \sR}$ under the isomorphism of Lemma \ref{Bohr2}, where $\delta_{V,V} = \id_V$ and $\delta_{V,V_0} = 0$ if $V_0 \neq V$. Then for any $U \in \Mod^{\fd,\semisimple}_B$, 
\[
\Hom_A(V, U|_A) \cong \Hom_B(B \otimes_A V, U) = 0,
\]
and thus for any isomorphism $U|_A \cong V_1 \oplus \cdots \oplus V_n$ with $V_i \in \sR$, each $V_i \neq V$. In particular, as each $x_{V_i} = 0$, then $y_U = x_{U|_A} = 0$, and so, as $U$ was arbitrary, $\widehat{\phi}((x_W)_W) = 0$, whilst $(x_W)_W \neq 0$.
\end{proof}

What is the advantage of $\widehat{A}$ over $A$? It allows us to use idempotents when dealing with simple modules. Let $\sI_A$ denote the set of central primitive idempotents of $\widehat{A}$, and recall that $\Irr(A)$ denote the set of isomorphism classes of irreducible $A$-modules.

Given a simple $A$-module $V$, let $\sR$ be any set of representatives of $\Irr(A)$ containing $V$, and let $e_V \in \widehat{A}$ denote the element of $\widehat{A}$ mapping to $(\delta_{V,V_0})_{V \in \sR}$ under the isomorphism of Lemma \ref{Bohr2}, where $\delta_{V,V} = \id_V$ and $\delta_{V,V_0} = 0$ if $V_0 \neq V$. Concretely, $e_V$ is constructed as in the proof of Lemma \ref{Bohr2}, and acts as the identity on any simple $A$-module $W$ isomorphic to $V$, and as zero on all other simple $A$-modules. In particular, $e_V$ is independent of the choice of representatives $\sR$ and depends only on the isomorphism class of $V$.

\begin{lemma}\label{lem:idemirrcorresp}
  The map $V \mapsto e_V$ defines a bijection $\Irr(A) \rightarrow \sI_A$.
\end{lemma}

\begin{proof}
  The map $V \mapsto e_V$ is injective by construction, and surjective from the description of $\widehat{A}$ of Lemma \ref{Bohr2} as a product of matrix rings over division rings.
\end{proof}

\begin{remark}\label{rem:bohrcompletionexample}
  Given $V \in \Irr(A)$, one can naturally view $V$ as a $\widehat{A}$-module via pullback along the projection $\widehat{A} \twoheadrightarrow \End_{E(V)} (V)$, and in this map we obtain an inclusion
  \[
  \Mod_A^{\fd, \semisimple} / \cong \, \lhook\joinrel\longrightarrow \, \Mod_{\widehat{A}}^{\fd,\semisimple} / \cong.
  \]

  We note that there can be finite-dimensional simple $\widehat{A}$-modules which are not in the image of this map. To construct an example, suppose that $\End_{E(V)}(V) \cong \bF_p$ for any $V \in \Irr(A)$ for some fixed prime $p$, so that any choice $\sR$ of representatives $\sR$ of $\Irr(A)$ induces an isomorphism
  \[
  \widehat{A} \xrightarrow{\sim} \prod_{i \in \sR} \bF_p
  \]
  by Lemma \ref{Bohr2}. Being a commutative ring, the the simple modules of $\widehat{A}$ correspond to the maximal ideals of $\widehat{A}$, which, as $\widehat{A}$ is an product of fields, are naturally identified with ultrafilters on $I$. Under the correspondence the principal ultrafilters correspond the simple modules $V \in \sR$, and the existence of free (i.e.\ non-principal) ultrafilters shows that there exist simple $\widehat{A}$-modules outside the image of the above map whenever $\Irr(A)$ is infinite. To see that these are all finite-dimensional, note that if $\fm$ is any maximal ideal of $\widehat{A}$, then we can view the left $\widehat{A}$-module $M = \widehat{A} / \fm$ as a field extension of $\bF_p$. However, for any surjection 
  \[
  \prod_{i \in \sR} \bF_p \twoheadrightarrow k
  \]
  to a field $k$, $p \cdot k = 0$ and $x^p = x$ for any $x \in k$, hence $k = \bF_p$, and so $M$ is $1$-dimensional over $\bF_p$. 

  To see that such an example is actually possible, consider the group ring $A = \bF_p[G]$, where $G = C_2^{\bN}$ and $p$ is an odd prime. Then finite-dimensional irreducible representations of $G$ over $\bF_p$ are all $1$-dimensional, and correspond to elements of $\{\pm 1\}^{\bN}$, which follows from the fact that an any family of commuting diagonalisable operators on a finite-dimensional vector space can be simultaneously diagonalised. In particular $\Irr(A)$ is infinite, and $\End_{E(V)}(V) = \bF_p$ for all $V \in \Irr(A)$.
\end{remark}

We are interested in $\widehat{A}$ for the following reason. Let us assume that all the maps on Diagram~\eqref{eqn:diamondrings} extend to Bohr completions. Then the completed version, Diagram~\eqref{eqn:diamondrings_complete} below, exhibits five algebras whose central primitive idempotents correspond to the simple finite-dimensional modules of the algebras on Diagram~\eqref{eqn:diamondrings} itself:
\begin{equation}\label{eqn:diamondrings_complete}
 \begin{tikzcd}
	& \widehat{\bK[G]} \\
	\widehat{\bL[G]} & \widehat{\bK[H]} & \widehat{\bL_{\sigma \hspace{-0.1em}}[G]} \\
	& \widehat{\bL[H]}
	\arrow[ from=1-2, to=2-1]
	\arrow[ from=1-2, to=2-3]
	\arrow[ from=2-2, to=1-2]
	\arrow[ from=2-2, to=3-2]
	\arrow[ from=3-2, to=2-1]
	\arrow[ from=3-2, to=2-3]
 \end{tikzcd}
\end{equation}

In Section \ref{sect:blockinterpretation} we will consider central idempotents of $\widehat{\bK[H]}$ which remain central in each larger ring, and use the resulting diagram of algebras, an \emph{S-block}, to study how semilinear representations of $G$ over $\bL$ interact with the linear representation theory of $G$ and $H$ over $\bK$ and $\bL$.

We will use the following lemma, which we advise be skipped and referred back to when needed.

\begin{lemma}\label{lem:gradedzeropartcompletion}
    Suppose that $f \colon A \rightarrow B$ is a $\bK$-algebra morphism, for which the functors
\[\begin{tikzcd}
	{\Mod_A} & {\Mod_B}
	\arrow["{B \otimes_A -}", curve={height=-12pt}, from=1-1, to=1-2]
	\arrow["{(-)|_A}", curve={height=-12pt}, from=1-2, to=1-1]
\end{tikzcd}\]
    preserve semisimplicity and finite-dimensionality. Suppose that $\phi \colon B \rightarrow B$ is an order two $\bK$-linear automorphism, $\charfield(\bK) \neq 2$, and $f(A) = B^{\phi = 1}$. Then
    \[
    \widehat{f}(\widehat{A}) = (\widehat{B})^{\widehat{\phi} = 1}.
    \]
\end{lemma}

\begin{proof}
    For any $x \in \widehat{A}$ we have that 
    \[
    \widehat{\phi} \widehat{f}(x) = (x_{\phi_* f_* N})_N = (x_{f_* N})_N = \widehat{f}(x),
    \]
    as $f = \phi \circ f$, and thus $\widehat{f}(\widehat{A}) \subset \widehat{B}^{\widehat{\phi} = 1}$. Conversely, suppose that $z \in \widehat{B}$ with $\widehat{\phi}(z) = z$. Because $\phi$ has order two and $\charfield(\bK) \neq 2$, we may define an $A$-linear projection operator
    \[
    \varphi \colon B \rightarrow A, \qquad \varphi(b) = \frac{b + \phi(b)}{2}.
    \]
    We define $x \in \widehat{A}$ by
    \[
        x_M \coloneqq \varphi_M \circ z_{B \otimes_A M} \circ \iota_M,
    \]
    where 
    \[
    \iota_M \colon M \rightarrow B \otimes_A M, \qquad \iota_M(m) = 1 \otimes m
    \]
    is the $A$-linear inclusion and 
    \[
    \varphi_M \colon B \otimes_M \rightarrow M, \qquad \varphi_M(b \otimes m) \coloneqq \varphi(b)m
    \]
    is the $A$-linear map defined using the projection $\varphi$. It is straightforward to see that $x$, defined in this way, gives a well-defined element of $\widehat{A}$, and so it remains to show that $\widehat{f}(x) = z$. This is equivalent to showing that for any $N \in \Mod_B^{\fd,\semisimple}$,
    \[
    x_{f_* N} = \widehat{f}(x)_N = z_N,
    \]
    or equivalently, that
    \[
    \varphi_{f_* N} \circ z_{B \otimes_A N} \circ \iota_{f_* N} = z_N.
    \]
    There is a $B$-linear map
    \[
    \pi_N \colon B \otimes_A N \rightarrow N, \qquad \pi_N(b \otimes n) = bn,
    \]
    and therefore 
    \[
    \pi_N \circ z_{B \otimes_A N} = z_N \circ \pi_N,
    \]
    and pre-composing with $\iota_{f_* N}$ we have that
    \[
    \pi_N \circ z_{B \otimes_A N} \circ \iota_{f_* N} = z_N
    \]
    hence it suffices to show that
    \[
    \varphi_{f_* N} \circ z_{B \otimes_A N} \circ \iota_{f_* N} =  \pi_N \circ z_{B \otimes_A N} \circ \iota_{f_* N}.
    \]
    From the fact that $\widehat{\phi}(z) = z$ means that $z_{\phi_* (B \otimes_A N)} = z_{B \otimes_A N}$. We have a $B$-linear map
    \[
    \phi \otimes 1 \colon B \otimes_A N \rightarrow \phi_*(B \otimes_A N),
    \]
    and therefore (using that $z \in \widehat{B}$),
    \[
    (\phi \otimes 1) \circ z_{\phi_* (B \otimes_A N)} = z_{B \otimes_A N} \circ (\phi \otimes 1)
    \]
    and thus
    \[
    (\phi \otimes 1) \circ z_{B \otimes_A N} = z_{B \otimes_A N} \circ (\phi \otimes 1).
    \]
    By definition, 
    \[
    \varphi_{f_* N} = \frac{1}{2}\left( \pi_N \circ (\phi \otimes 1) + \pi_N \right),
    \]
    and therefore we can compute
    \begin{align*}
       \varphi_{f_* N} \circ z_{B \otimes_A N} \circ \iota_{f_* N} 
       &= \frac{1}{2}\left( \pi_N \circ (\phi \otimes 1) \circ z_{B \otimes_A N} \circ \iota_{f_* N} + \pi_N \circ z_{B \otimes_A N} \circ \iota_{f_* N} \right), \\
       &= \frac{1}{2}\left( \pi_N \circ z_{B \otimes_A N} \circ (\phi \otimes 1) \circ \iota_{f_* N} + \pi_N \circ z_{B \otimes_A N} \circ \iota_{f_* N} \right), \\
    &= \pi_N \circ z_{B \otimes_A N} \circ \iota_{f_* N},
    \end{align*}
    using that $(\phi \otimes 1) \circ \iota_{f_* N} = \iota_{f_* N}$. This shows that $x_{f_* N} = z_N$, and therefore $\widehat{f}(x) = z$, as required.
\end{proof}

\subsection{Bohr Completion for Group Ring Extensions}\label{sect:SblockStructure}
       How common are homomorphisms satisfying the conditions of Lemma~\ref{hom_cond}?
      Our goal in this section is to show that the algebra extensions of Diagram~\eqref{eqn:diamondrings} all satisfy these assumptions, subject to mild additional conditions.
       
   \begin{lemma} \label{Sblock_1}
     Suppose $\bK\leq \bL$ is a finite field extension, $A$ is a $\bK$-algebra. Let $B=\bL \otimes_{\bbK} A$ be the extension of scalars.
     Then the restriction functor the restriction functor  $\Mod_B \to \Mod_A$ preserves finite-dimensional semisimple modules.

If, furthermore, the extension $\bK\leq \bL$ is separable, then the induction functor 
  $\Mod_A \to \Mod_B$
  takes finite-dimensional simple $A$-modules to non-zero finite-dimensional semisimple modules.
\end{lemma}
\begin{proof} 
  Let $V$ be a finite dimension simple $B$-module. The algebra of operators on $V$, coming from $B$ is $B/I(V)$. As in the proof of Lemma \ref{Bohr2}, $B / I(V)$ is isomorphic to a matrix algebra $M_n(D)$, where $D$ is the finite-dimensional division $\bL$-algebra $D \coloneqq \End_B(V)$. Let $A^\prime$ be the image of $A$ in $B/I(V)$. Observe that $A^\prime$ is semisimple: any nilpotent ideal $J$ of $A^\prime$ will give a nilpotent ideal $\bL \otimes_{\bbK} J$ in $B/I(V)$, but $B/I(V)$ is simple, so $J=0$.
  Thus, $V$ is a semisimple $A$-module.

  For the second statement, consider a finite-dimensional simple $A$-module $M$.
  Similarly, $A/I(M)$ is isomorphic to a matrix algebra $M_n(D)$ for some finite-dimensional division $\bK$-algebra $D$. Under this identification, the induced module $\bL \otimes_{\bK} M$ is identified with $(\bL \otimes_{\bK} D)^n$ and its algebra $B'$ of operators coming from $B$ is identified with $M_n (\bL \otimes_{\bK} D)$.
  Separability implies that $B'$ is semisimple and, hence, the induced $B$-module is also semisimple.
\end{proof}
In light of Lemma~\ref{hom_cond}, we have the following corollary.
\begin{cor}
If $\bK\leq \bL$ is a finite separable field extension, then the canonical $\bK$-algebra homomorphisms
	$\widehat{\bK[H]} \to \widehat{\bL[H]}$ and $\widehat{\bK[G]} \to \widehat{\bL[G]}$ are well-defined and injective.
\end{cor}  

Let $\Gamma$ be a finite group. The second batch of algebras in Diagram~\eqref{eqn:diamondrings} are $\Gamma$-crossed products, for $\Gamma = G/H$.
Recall that if $A$ is a $\bK$-algebra, then a crossed product $A\sharp_{\alpha,\gamma}\Gamma$ is an algebra structure on the group algebra $(A[\Gamma],+)$ given by
\[ ag \cdot bh = a \gamma(g)(b) \alpha(g,h) gh
\ \mbox{ for all } \ a,b\in A, \ g,h \in \Gamma, \]
where $\alpha : \Gamma \times \Gamma \to A$ is a 2-cocycle, $\gamma : \Gamma \to \Aut (A)$ is a homomorphism. This has to satisfy natural conditions ensuring that $A\sharp_{\alpha,\gamma}\Gamma$ is an associative $\bK$-algebra with identity $1_A 1_\Gamma$.

   \begin{lemma} \label{Sblock_2}
     Suppose that $B$ is a crossed product $A\sharp_{\alpha,\gamma}\Gamma$ for some $\bK$-algebra $A$ and finite group $\Gamma$.
      Then the restriction functor  $\Mod_B \to \Mod_A$ preserves finite-dimensional semisimple modules.

If, furthermore, $\Gamma$ is finite and $|\Gamma| \in \bL^{\times}$, then the induction functor 
  $\Mod_A \to \Mod_B$ takes finite-dimensional simple $A$-modules to non-zero finite-dimensional semisimple modules.
\end{lemma}
   \begin{proof} Observe that the induction of an $A$-module $M$ can be written as
     \[
     B \otimes_A M = \bigoplus_{g\in \Gamma} g \otimes_{\bK} M = \bigoplus_{g\in \Gamma} g \ast M,
     \]
     where the action of $A$ on $g \ast M$ is given by
     \[ a \cdot (g \ast m) = agm = g \gamma(g)^{-1} (a) m = g \ast (\gamma(g)^{-1} (a) \cdot m) . \]
     So it is just the module twisted by the automorphism $\gamma (g)^{-1}$. In particular, if $M$ is semisimple, then $B\otimes_A M$
     is semisimple as an $A$-module.

     Let $V$ be a finite-dimensional simple $B$-module. There exists a simple $A$-submodule $M$ inside $V$. By simplicity of $V$,
     the induced homomorphism $B\otimes_A M \to V$ is surjective. Thus, $V$ is a quotient of a semisimple $A$-module.
     So it is semisimple as an $A$-module.

     Before proving the second statement, observe that our condition on $\Gamma$ implies that the extension $A \leq B$ is separable in the sense of \cite[Def.\ 2]{HS}. Indeed, its separating idempotent is
     \[ e = \sum_i x_i\otimes_A y_i = \frac{1}{|\Gamma|} \sum_{g \in \Gamma}  \alpha(g,g^{-1})^{-1}g \otimes_A g^{-1} . \]
     This means that any extension of $B$-modules, split as an extension of $A$-modules, is also split as an extension of $B$-modules \cite[Prop.\ 2.6, Thm.\ 1.1]{HS}. Explicitly, given an $A$-splitting $\Phi$, then $\sum_i x_i \Phi y_i$ is the required $B$-splitting (which in our case is just Maschke averaging over the group $\Gamma$).

     Now consider a finite-dimensional simple $A$-module $M$. $B\otimes_A M$ is non-zero by the description above.
     Pick any simple $B$-submodule $V$ inside the induced module $B\otimes_A M$.
  We observed that $B\otimes_A M$ is $A$-semisimple. Thus, the extension $V \hookrightarrow B \otimes_A M$ is $A$-split.
By separability, it is also $B$-split. Hence, $B\otimes_A M$ is $B$-semisimple. 
\end{proof}

\begin{cor}\label{cor:secondbatch}
 If $\bL / \bK$ is finite and $[\bL : \bK] \in \bL^{\times}$, 
 then the canonical $\bK$-algebra homomorphisms
 $\widehat{\bK[H]} \to \widehat{\bK[G]}$,
  $\widehat{\bL[H]} \to \widehat{\bL[G]}$,
 and $\widehat{\bL[H]} \to \widehat{\bL_{\sigma}[G]}$ are well-defined and injective.
\end{cor}  

\begin{proof}
  This follows from Lemma \ref{hom_cond} and Lemma \ref{Sblock_2}, noting that by the assumption that $\bL / \bK$ is finite, $\bL / \bK$ is automatically Galois and the induced map $G / H \rightarrow \Gal(\bL / \bK)$ is an isomorphism \cite[Rem.\ 3.5]{TAYSchur}, hence $\Gamma \coloneqq G/H$ satisfies the assumptions of Lemma \ref{Sblock_2}.
\end{proof}

\begin{remark}
  Actually, for the injectivity of $\widehat{\bL[H]} \to \widehat{\bL_{\sigma}[G]}$, you do not need that $[\bL : \bK] \in \bL^{\times}$, only that $\bL / \bK$ is finite. This will follow from Theorem \ref{thm:facts} below.
\end{remark}

The sixth homomorphism on the Diagram~\eqref{eqn:diamondrings} does not immediately fall into an easily identifiable construction but
we can deal with it by using the result for the other homomorphisms.
\begin{thm} \label{SBlock_3}
Suppose that $\bL / \bK$ is finite and $[\bL : \bK] \in \bL^{\times}$.
Then all homomorphisms in Diagram~\eqref{eqn:diamondrings} satisfy the conditions of Lemma~\ref{hom_cond}, and thus 
all homomorphisms in Diagram~\eqref{eqn:diamondrings_complete} exist and are injective. 
\end{thm}

\begin{proof}
  After the last two lemmas, we know this is true for five of the homomorphisms of Diagram \eqref{eqn:diamondrings},
  and it remains to cobble something together for the sixth homomorphism  ${\bK[G]} \to {\bL_{\sigma}[G]}$. It is clear that induction and restriction between these $\bK$-algebras preserves finite-dimensionality.

  Suppose first that $W$ is a finite-dimensional simple $\bL_{\sigma}[G]$-module. By the commutativity of Diagram \eqref{eqn:diamondrings} and the result for the other homomorphisms of Diagram \eqref{eqn:diamondrings}, $(W|_{\bbK[G]})_{\bbK[H]}$ is semisimple, hence $W|_{\bbK[G]}$ is semisimple being the quotient of the semisimple module
  \[
  \bK[G] \otimes_{\bbK[H]} (W|_{\bbK[G]})_{\bbK[H]} \twoheadrightarrow W|_{\bbK[G]}
  \]
  under the natural surjection.

  If instead $U$ is a finite-dimensional simple $\bK[G]$-module, $\bL_{\sigma}[G] \otimes_{\bK[G]} U$ is non-zero (it can be identified with $\bL \otimes_{\bK} U$ with the diagonal action of $G$ \cite[\S 4.1]{TAYSchur}), and we again have a canonical $\bK[G]$-module surjection
  \[
      \bK[G] \otimes_{\bbK[H]} U|_{\bbK[H]} \twoheadrightarrow U
  \]
  which induces an $\bL_{\sigma}[G]$-module surjection
  \[
  \bL_{\sigma}[G] \otimes_{\bbK[H]} U|_{\bbK[H]} \equiv \bL_{\sigma}[G] \otimes_{\bK[G]} (\bK[G] \otimes_{\bbK[H]} U|_{\bbK[H]}) \twoheadrightarrow \bL_{\sigma}[G] \otimes_{\bK[G]} U
  \]
  from the right-exactness of the tensor product. The $\bL_{\sigma}[G]$-module $\bL_{\sigma}[G] \otimes_{\bbK[H]} U|_{\bbK[H]}$ is semisimple by the commutativity of Diagram \eqref{eqn:diamondrings} and the result for the other homomorphisms of Diagram \eqref{eqn:diamondrings}, and therefore the quotient $\bL_{\sigma}[G] \otimes_{\bK[G]} U$ is too. 
\end{proof}

\section{Semilinear Representations for Finite Galois Extensions}\label{sect:finiteGaloisextensions}

From now on in this chapter we specialise and assume that $\bL / \bK$ is a finite extension. We note that this automatically implies that $\bL / \bK$ is Galois, and that the induced map $\sigma \colon G \rightarrow \Gamma \coloneqq \Gal(\bL / \bK)$ is surjective \cite[Rem.\ 3.5]{TAYSchur}. For example, the assumption that $\bL / \bK$ is finite is automatic if $G$ is finite or $\bL$ is a finite extension of its prime subfield. In this section we recall some results from \cite{TAYSchur} and describe $\Irr_{\bbL}^{\hspace{0.1em} \sigma}(G)$ in terms of $\Irr_{\bbL}(H)$, where $H = \ker(\sigma \colon G \rightarrow \Gal(\bL / \bK))$. 

There is a natural action of $\Gamma$ on $\Irr_{\bbL}(H)$ defined as follows.

\begin{defn}
For $g \in G$ and $V \in \Irr_{\bbL}(H)$ we let $g * V$ be the $\bL[H]$-module which is $V$ as a $\bK$-vector space, but with $\bL[H]$-module structure defined by
\[
\lambda h \cdot v = {}^g\lambda \cdot (g^{-1} h g)(v)
\]
for $\lambda \in \bL$, $h \in H$ and $v \in V$.
\end{defn}
This defines an left action of $G$ on $\Irr_{\bbL}(H)$ which is trivial on $H$, and therefore induces an action of the quotient $G/H$ and thus of $\Gamma$.

\begin{defn}
For $V \in \Irr_{\bbL}(H)$, we let $\Gamma_V$ denote the stabiliser of $V$ in $\Gamma$.
\end{defn}
There are induction and restriction functors between $\Rep_{\bbL}^{\hspace{0.1em} \sigma}(G)$ and $\Rep_{\bbL}(H)$ where
\[
(-)|_H \colon \Rep_{\bbL}^{\hspace{0.1em} \sigma}(G) \rightarrow \Rep_{\bbL}(H)
\]
is just restriction along the inclusion $\bL[H] \rightarrow \bL_{\sigma} \hspace{-0.1em}[G]$ and
\[
\Ind_H^G(-) \coloneqq \bL_{\sigma} \hspace{-0.1em}[G] \otimes_{\bL[H]} - \colon \Rep_{\bbL}(H) \rightarrow \Rep_{\bbL}^{\hspace{0.1em} \sigma}(G).
\]
The following theorem describes how $\Irr_{\bbL}^{\hspace{0.1em} \sigma}(G)$ can be understood in terms of $\Irr_{\bbL}(H)$.
\begin{thm}\label{thm:facts}
Suppose that $\bL/ \bK$ is finite. If $W_1, W_2 \in \Rep_{\bbL}^{\hspace{0.1em} \sigma}(G)$, then:
\begin{itemize}
    \item If $W_1|_H \cong W_2|_H$ then $W_1 \cong W_2$.
    \item The natural map
    \[
    \bL \otimes_{\bbK} \Hom_{\bL_{\sigma} \hspace{-0.1em}[G]}(W_1, W_2) \rightarrow \Hom_{\bL[H]}(W_1|_H, W_2|_H)
    \]
    is an isomorphism.
    \item In particular, if $W \in \Rep_{\bbL}^{\hspace{0.1em} \sigma}(G)$ then the natural map
    \[
        \bL \otimes_{\bbK} \End_{\bL_{\sigma} \hspace{-0.1em}[G]}(W) \rightarrow \End_{\bL[H]}(W|_H)
    \]
    is a $\bK$-algebra isomorphism.
\end{itemize}

Let $W \in \Irr_{\bbL}^{\hspace{0.1em} \sigma}(G)$ and $V \in \Irr_{\bbL}(H)$, where $V$ is an irreducible submodule of $W|_H$, or equivalently that $W$ is an irreducible submodule of $\Ind_H^G(V)$. Then:
\begin{itemize}
    \item $W|_H$ and $\Ind_H^G(V)$ are described by
    \[
        W|_H \cong \bigoplus_{\gamma \in \Gamma / \Gamma_V} (\gamma * V)^{m(W)}, \qquad \Ind_H^G(V) \cong W^{\oplus |\Gamma_V| / m(W)}
    \]
    for some integer $m(W) \geq 1$.
    \item In particular, $W|_H$ and $\Ind_H^G(V)$ are semisimple, and $m(W) \mid |\Gamma_V|$.
    \item $\Ind_H^G(V) \cong \Ind_H^G(\gamma * V)$ for any $\gamma \in \Gamma$.
    \item Any $V \in \Irr_{\bbL}(H)$ is contained in $W|_H$ for a unique $W \in \Irr_{\bbL}^{\hspace{0.1em} \sigma}(G)$. This defines a bijective correspondence
    \[
        \Irr_{\bbL}^{\hspace{0.1em} \sigma}(G) \leftrightarrow \Irr_{\bbL}(H) / \hspace{0.1em}\Gamma.
    \]
\end{itemize}
\end{thm}
\begin{proof}
    This is Theorem A and Theorem B from \cite{TAYSchur}.
\end{proof}

\begin{remark}
  Theorem A and B of \cite{TAYSchur} is stated more generally for indecomposable modules, but we restrict ourselves to the study of irreducible modules here.
\end{remark}

\subsection{Semilinear Schur indices}\label{sect:semilinearschurindex}

Theorem \ref{thm:facts} gives an almost complete description of $\Irr_{\bbL}^{\hspace{0.1em} \sigma}(G)$ purely in terms of $\Irr_{\bbL}(H)$ and its action of $\Gamma$ as follows. 

\begin{defn}
For $V \in \Irr_{\bbL}(H)$, the \emph{semilinear Schur index} of $V$ is $m_{\bbK}^{\bbL}(V) \coloneqq m(W)$, where $W \in \Irr_{\bbL}^{\hspace{0.1em} \sigma}(G)$ is the unique irreducible semilinear $\bL$-representation of $G$ corresponding to $V$.
\end{defn}

Then to construct all objects of $\Irr_{\bbL}^{\hspace{0.1em} \sigma}(G)$, first compute the orbits of the $\Gamma$-action on $\Irr_{\bbL}(H)$. For each such orbit $\sA$, $m_{\bbK}^{\bbL}(V)$ copies of the direct sum
\[
\bigoplus_{V \in \sA} V
\]
extends uniquely to an object of $\Irr_{\bbL}^{\hspace{0.1em} \sigma}(G)$, and all objects of $\Irr_{\bbL}^{\hspace{0.1em} \sigma}(G)$ arise this way.

This description is not quite complete, because it still rests on a determination of the semilinear Schur indices $m_{\bbK}^{\bbL}(V)$ for $V \in \Irr_{\bbL}(H)$, which, like their classical counterparts (of which these are a generalisation \cite[\S 10]{TAYSchur}), appear very hard to compute in general.

It is natural to ask how much information about the numbers $m_{\bbK}^{\bbL}(V)$ is contained in the other associated linear representation theory categories of diagram (\ref{eqn:diamond}). We shall give a satisfactory answer to this question in Chapter \ref{chapter:twelvefoldway}, where we consider the simplest non-trivial case, that $\bL / \bK$ has degree 2.

\subsection{Associated linear representation categories}

The assumption that $\bL / \bK$ is finite allows us to say more about the extensions of diagram (\ref{eqn:diamondrings}). Recall that an extension $R \rightarrow S$ of $\bK$-algebras is called \emph{Frobenius} if the induction and coinduction functors $\Mod_R \rightarrow \Mod_S$ are naturally isomorphic \cite[Def.\ 1.1]{KAD}. The following is Proposition 4.1 of \cite{TAYSchur}.

\begin{prop}\label{prop:frobenius}
If $\bL / \bK$ is finite then all extensions of diagram (\ref{eqn:diamondrings}) are Frobenius.
\end{prop}

\subsection{Induction from subgroups}

For the remained of this chapter we continue with our assumption that $\bL/\bK$ is a finite extension and record some basic properties concerning induction from subgroups. We recommend this be skipped and only referred back to when necessary.

Recall that $\Gamma \coloneqq \Gal(\bL/\bK)$, and $\sigma \colon G \twoheadrightarrow \Gamma$ denotes the action map with kernel $H$. For a subgroup $S$, we can consider the subalgebra $\bL_{\sigma} \hspace{-0.1em}[S]$ of $\bbL_{\sigma} \hspace{-0.1em}[G]$.

\begin{defn}
    For a finite-index subgroup $S \leq G$ we set
    \[
        \Ind_S^G(-) \coloneqq \bL_{\sigma} \hspace{-0.1em}[G] \otimes_{\bbL_{\sigma} \hspace{-0.1em}[S]} (-) \colon \Rep_{\bbL}^{\hspace{0.1em} \sigma}(S) \rightarrow \Rep_{\bbL}^{\hspace{0.1em} \sigma}(G),
    \]
    and for an arbitrary subgroup $S \leq G$ we set
    \[
        (-)|_S \coloneqq \Res_{\bbL_{\sigma} \hspace{-0.1em}[S]}^{\bbL_{\sigma} \hspace{-0.1em}[G]}(-) \colon \Rep_{\bbL}^{\hspace{0.1em} \sigma}(G) \rightarrow \Rep_{\bbL}^{\hspace{0.1em} \sigma}(S).
    \]
\end{defn}

\begin{lemma}\label{lem:inductioncompatibility}
    Let $S$ be a finite-index subgroup of $G$, and $T$ be the kernel of the induced map $S \rightarrow \Gamma$. Then the diagrams
\[\begin{tikzcd}
	{\Rep_{\bbL}^{\hspace{0.1em} \sigma}(S)} & {\Rep_{\bbL}^{\hspace{0.1em} \sigma}(G)} && {\Rep_{\bbL}^{\hspace{0.1em} \sigma}(S)} & {\Rep_{\bbL}^{\hspace{0.1em} \sigma}(G)} \\
	{\Rep_{\bbL}(T)} & {\Rep_{\bbL}(H)} && {\Rep_{\bbL}(T)} & {\Rep_{\bbL}(H)}
	\arrow["{\Ind_{S}^T}", from=1-1, to=1-2]
	\arrow["{(-)|_T}"', from=1-4, to=2-4]
	\arrow["{(-)|_S}"', from=1-5, to=1-4]
	\arrow["{(-)|_H}", from=1-5, to=2-5]
	\arrow["{\Ind_{T}^S}", from=2-1, to=1-1]
	\arrow["{\Ind_{T}^H}"', from=2-1, to=2-2]
	\arrow["{\Ind_{H}^G}"', from=2-2, to=1-2]
	\arrow["{(-)|_T}", from=2-5, to=2-4]
\end{tikzcd}\]
    commute up to natural isomorphism. If additionally $S$ surjects onto $\Gamma$, then the diagrams
\[\begin{tikzcd}
	{\Rep_{\bbL}^{\hspace{0.1em} \sigma}(S)} & {\Rep_{\bbL}^{\hspace{0.1em} \sigma}(G)} && {\Rep_{\bbL}^{\hspace{0.1em} \sigma}(S)} & {\Rep_{\bbL}^{\hspace{0.1em} \sigma}(G)} \\
	{\Rep_{\bbL}(T)} & {\Rep_{\bbL}(H)} && {\Rep_{\bbL}(T)} & {\Rep_{\bbL}(H)}
	\arrow["{\Ind_{S}^T}", from=1-1, to=1-2]
	\arrow["{(-)|_T}", from=1-1, to=2-1]
	\arrow["{(-)|_H}", from=1-2, to=2-2]
	\arrow["{(-)|_S}", from=1-5, to=1-4]
	\arrow["{\Ind_{T}^H}"', from=2-1, to=2-2]
	\arrow["{\Ind_{T}^{S}}", from=2-4, to=1-4]
	\arrow["{\Ind_H^G}"', from=2-5, to=1-5]
	\arrow["{(-)|_T}", from=2-5, to=2-4]
\end{tikzcd}\]
    commute up to natural isomorphism.
\end{lemma}

\begin{proof}
    The commutativity of the first two diagrams follows directly from the transitivity of induction and restriction. Suppose now that the induced map $\sigma \colon S \rightarrow \Gamma$ is surjective, which means that $[G:H] = [S:T]$ and $[G:S] = [H:T]$. The first isomorphism is given for $W \in \Rep_{\bbL}^{\hspace{0.1em} \sigma}(S)$ by
    \[
        \phi \colon \bL[H] \otimes_{\bbL[T]} (W|_T) \rightarrow \left( \bL_{\sigma} \hspace{-0.1em}[G] \otimes_{\bbL_{\sigma} \hspace{-0.1em}[S]} W \right)|_H, \qquad \phi(h \otimes w) \coloneqq h \otimes w.
    \]  
    Because $[G:S] = [H:T]$, both sides have $\bL$-dimension $[G:S] \dim_{\bbL}(W)$, and so it is sufficient to show that $\phi$ is surjective. Indeed, for any $g \in G$, $g = h s$ for $h \in H$ and $s$ in $S$, because the coset representatives of $T$ in $S$ are coset representatives for $H$ in $G$, and so for any $w \in W$,
    \[
        \phi(h \otimes s(w)) = h \otimes s(w) = g \otimes w.
    \]
    The second isomorphism is similarly given for $V \in \Rep_{\bbL}(H)$ by
    \[
        \psi \colon \bL_{\sigma} \hspace{-0.1em}[S] \otimes_{\bbL[T]} (V|_T) \rightarrow \left( \bL_{\sigma} \hspace{-0.1em}[G] \otimes_{\bbL[H]} V \right)|_S, \qquad \phi(s \otimes v) \coloneqq s \otimes v,
    \]  
    which can similarly be seen to be an isomorphism using that any $g \in G$ has the form $g = s h$ for $s \in S$ and $h \in H$ and that both sides have the same $\bL$-dimension as $[G:H] = [S:T]$.
\end{proof}

\begin{lemma}\label{lem:projectionformula}
    Suppose that $S \leq G$ is a finite-index subgroup. Then there is a natural isomorphism
    \[
        \Ind_S^G(W \otimes_{\bbL} V|_S) \xrightarrow{\sim} \Ind_S^G(W) \otimes_{\bbL} V
    \]
    in $V \in \Rep_{\bbL}^{\hspace{0.1em} \sigma}(G)$ and $W \in \Rep_{\bbL}^{\hspace{0.1em} \sigma}(S)$.
\end{lemma}

\begin{proof}
    The natural map of $\bL_{\sigma} \hspace{-0.1em}[S]$-modules
    \[
        W \otimes_{\bbL} V|_S \rightarrow (\bL_{\sigma} \hspace{-0.1em}[G] \otimes_{ \bbL_{\sigma} \hspace{-0.1em}[S]} W) \otimes_{\bbL} V,
    \]
    extends to
    \[
    \bL_{\sigma} \hspace{-0.1em}[G] \otimes_{ \bbL_{\sigma} \hspace{-0.1em}[S]}(W \otimes_{\bbL} V|_S) \rightarrow (\bL_{\sigma} \hspace{-0.1em}[G] \otimes_{ \bbL_{\sigma} \hspace{-0.1em}[S]} W) \otimes_{\bbL} V
    \]
    which is surjective, hence an isomorphism, both sides having the same $\bK$-dimension.
\end{proof}

\subsubsection{Comparison with Classical Induction}

We also record here that this induction is compatible with induction in the split case, under the equivalence of \cite[\S 6]{TAYSchur}.

\begin{lemma}\label{lem:compatiblewithsplitinduction}
    Suppose that $S$ is a subgroup of $H$. Then the diagrams
\[\begin{tikzcd}
	{\Rep_{\bbK}(H)} & {\Rep_{\bbL}^{\hspace{0.1em} \sigma}(H \times \Gamma)} && {\Rep_{\bbK}(H)} & {\Rep_{\bbL}^{\hspace{0.1em} \sigma}(H \times \Gamma)} \\
	{\Rep_{\bbK}(S)} & {\Rep_{\bbL}^{\hspace{0.1em} \sigma}(S \times \Gamma)} && {\Rep_{\bbK}(S)} & {\Rep_{\bbL}^{\hspace{0.1em} \sigma}(S \times \Gamma)}
	\arrow["\sim", from=1-1, to=1-2]
	\arrow["\sim", from=1-4, to=1-5]
	\arrow["{(-)|_S}"', from=1-4, to=2-4]
	\arrow["{(-)|_{S \times \Gamma}}", from=1-5, to=2-5]
	\arrow["{\Ind_{S}^H}", from=2-1, to=1-1]
	\arrow["\sim", from=2-1, to=2-2]
	\arrow["{\Ind_{S \times \Gamma}^{H \times \Gamma}}"', from=2-2, to=1-2]
	\arrow["\sim", from=2-4, to=2-5]
\end{tikzcd}\]
commute up to natural isomorphism, where the horizontal arrows are the equivalences of \cite[\S 6]{TAYSchur}.
\end{lemma}

\begin{proof}
    For $U \in \Rep_{\bbK}(S)$, the first isomorphism is given by
    \[
        \bL \otimes_{\bbK} (\bK[H] \otimes_{\bbK[S]} U) \rightarrow \bbL_{\sigma} \hspace{-0.1em}[H \times \Gamma] \otimes_{\bbL_{\sigma} \hspace{-0.1em}[S \times \Gamma]} (\bL \otimes_{\bbK} U), \qquad \lambda \otimes h \otimes u \mapsto h \otimes \lambda \otimes u
    \]
    for $\lambda \in \bL$, $h \in H$ and $u \in U$. This is clearly surjective, hence an isomorphism as both sides have the same $\bL$-dimension. The second isomorphism is clear.
\end{proof}

%% file: chapter-twelveways.tex
\chapter{The 24 Types of Semilinear Representation}\label{chapter:twelvefoldway}

In this chapter we conduct an in-depth study of the irreducible objects of $\Rep_{\bbL}^{\hspace{0.1em} \sigma}(G)$ in the case that the extension $\bL / \bK$ has degree $2$. We will associate to each $W \in \Irr_{\bbL}^{\hspace{0.1em} \sigma}(G)$ a \emph{type}, which characterises the behaviour of $W$ under the induction and restriction functors of diagram (\ref{eqn:diamond}), and in many cases allows the semilinear Schur index $m(W)$ of $W$ (cf.\ Theorem \ref{thm:facts}) to be determined by purely linear representation theoretic data. We make the following standing assumption.

\vspace{0.5em}
\begin{center}
\emph{We assume in this chapter that $\bL / \bK$ has degree $2$.}
\end{center}
\vspace{0.5em}

In particular, this assumption implies that $H \leq G$ is a Real structure on $H$ (cf.\ Chapter \ref{chapter:realstructures}). We will additionally assume that $\bK$ has characteristic two from Section \ref{sect:diamondanalysis} onwards.
 
 \subsection{Notation}\label{sect:chap3notation} Because $\bL / \bK$ is quadratic, $\sigma \colon G \twoheadrightarrow \Gamma = \Gal(\bL / \bK)$ is determined by the requirement that $H$ act trivially. We thus suppress $\sigma$ from the notation and write:
 \begin{itemize}
    \item $\Rep_{\bbL}^{\hspace{0.1em} \star}(G) \coloneqq \Rep_{\bbL}^{\hspace{0.1em} \sigma}(G)$.
    \item $\Irr_{\bbL}^{\hspace{0.1em} \star}(G) \coloneqq \Irr_{\bbL}^{\hspace{0.1em} \sigma}(G)$.
    \item $\bL_{\star \hspace{-0.1em}}[G] \coloneqq \bL_{\sigma} \hspace{-0.1em}[G]$.
    \item $\overline{(-)} \colon \bL \rightarrow \bL$ is the unique non-trivial element of $\Gamma$.
    \item For $V \in \Vect_{\bbL}$, $\overline{V} \in \Vect_{\bbL}$ is $V$ but with $\bL$-vector space structure by
    \[
    \alpha * v \coloneqq \overline{\alpha} v \ \mbox{  for all  } \ \alpha\in\bL , v\in V.
    \]
 \end{itemize}

 \subsection{Action of $C_2$ on $\Irr_{\bbL}(H)$} \label{sect:C2actionIrrH}
 Recall that we write $C_2 = \langle \epsilon \rangle$ for the group with two elements. There is an action of $C_2$ on $\Irr_{\bbL} (H)$,
 defined on $V \in \Irr_{\bbL}(H)$ by setting $\epsilon * V \coloneqq {}^{\omega}\overline{V}$, which is $\overline{V}$
 as a $\bL$-vector space but with action of $H$ defined by picking an element $\omega \in G \setminus H$, and setting 
 \[
 h * v \coloneqq (\omega^{-1} h \omega)(v)
 \]
 for any $h \in H$ and $v \in V$. This gives ${}^{\omega}\overline{V}$ the structure of a $\bL[H]$-module,
 which is independent (up to isomorphism) of the choice of $\omega$. Hence it defines an action of $C_2$ on the set $\Irr_{\bbL}(H)$.

 \subsection{The three S-types of irreducible representation}

 Under our assumption that $\bL / \bK$ has degree $2$, Theorem \ref{thm:facts} specialises to the following.

 \begin{thm}\label{thm:threeoptions}
    Suppose that $W \in \Irr_{\bbL}^{\hspace{0.1em} \star}(G)$ and $V \in \Irr_{\bbL}(H)$ are such that $V$ is an irreducible submodule of $W|_H$, or equivalently that $W$ is an irreducible submodule of $\Ind_H^G(V)$. Then $V$ and $W$ are in either case \emph{A}, \emph{B} or \emph{C} of Table \ref{tab:thewaysRHSsquare}.
 \end{thm}
 \begin{table}[h]
\caption{Induction and Restriction Possibilities Between $\Irr_{\bbL}(H)$ and $\Irr_{\bbL}^{\hspace{0.1em} \star}(G)$}
\label{tab:thewaysRHSsquare}
\begin{center}
\def\arraystretch{1.1}
\begin{tabular}{ c c c c }
    \toprule 
Case & \textbf{A} & \textbf{B} & \textbf{C} \\
\midrule
 $W|_H$ & $V$ & $V \oplus {}^{\omega}\overline{V}$ & $V \oplus V$ \\ 
 $\Ind_H^G(V)$ & $W \oplus W$ & $W$ & $W$\\ 
 $V \cong  {}^{\omega}\overline{V}$? & Yes & No & Yes \\
  $m(W)$ & 1 & 1 & 2 \\
  $\dim_{\bL} \End_{\bL[H]}(V) : \dim_{\bK} \End_{\bL_{\star \hspace{-0.1em}}[G]}(W)$ & $1:1$ & $1:2$ & $1:4$ \\ 
  \bottomrule
\end{tabular}
\end{center} 
\end{table}
 \begin{proof}
    Because $\bL / \bK$ is finite we may apply Theorem \ref{thm:facts}, noting that the action of $C_2$-defined above agrees with that of $\Gamma$ in Theorem \ref{thm:facts}. The $C_2$-orbit of $V$ has size either $1$ or $2$. If $V \not\cong  {}^{\omega}\overline{V}$, then $\Gamma_V = 1$, and so by Theorem \ref{thm:facts} $m(W) = 1$ and
    \[
        W|_H \cong V \oplus  {}^{\omega}\overline{V}, \qquad \Ind_{H}^G(V) \cong W.
    \]
    Furthermore,
    \begin{align*}
    \dim_{\bK} \End_{\bL_{\star \hspace{-0.1em}}[G]}(W) &= \dim_{\bL} \End_{\bL[H]}(V \oplus {}^{\omega}\overline{V}), \\
    &= \dim_{\bL} \End_{\bL[H]}(V) +  \dim_{\bL} \End_{\bL[H]}({}^{\omega}\overline{V}), \\
    &= 2 \dim_{\bL} \End_{\bL[H]}(V).
    \end{align*}
    If instead $\Gamma_V = \Gamma$, then as $m(W) | |\Gamma_V| = 2$, there are two options: $m(W) = 1$ or $m(W) = 2$. Similarly to the above, these correspond to the two remaining cases (A) and (C).
 \end{proof}

In particular, we see that Case (B) is distinguished from Cases (A) and (C) from the action of $C_2$ on $\Irr_{\bbL}(H)$. However this action doesn't differentiate between cases (A) and (C).

The main aim of this chapter is to improve on this and obtain a finer classification of $\Irr_{\bbL}^{\hspace{0.1em} \star}(G)$, by considering the other linear representation categories of Diagram \eqref{eqn:diamond} beyond $\Rep_{\bbL}(H)$.

\section{Index Two Induction and Restriction}\label{sect:fourfoldway}

We first consider general induction and restriction behaviours between algebras equipped with certain endofunctors on their categories of representations, which we will apply to the extensions of Diagram \eqref{eqn:diamondrings}. This will give us a way to establish versions of Theorem \ref{thm:threeoptions} for all other extensions of diagram (\ref{eqn:diamondrings}) under our standing assumption that $\bL / \bK$ has degree $2$.

Suppose in this section that $\bK$ is any field, $A$ and $B$ are $\bK$-algebras, and $f \colon A \rightarrow B$ is a $\bK$-algebra homomorphism. The restriction functor $\Res \colon \Mod_B \rightarrow \Mod_A$
has left and right adjoint functors: induction and coinduction 
\[
\Ind \coloneqq B\otimes_{A}(-), \ \Cnd \coloneqq \Hom_A(B,-) \colon \Mod_A \rightarrow \Mod_B \, .
\]
We assume existence of covariant $\bK$-linear dualities
\begin{align*}
    (-)^* &\colon \Mod_A \rightarrow \Mod_A, \\
    \overline{(-)} \hspace{0.4em} &\colon \Mod_B \rightarrow \Mod_B ,
\end{align*}
meaning covariant functors whose composition with themselves are naturally isomorphic to the identity functor, which are related in the sense that there are isomorphisms
\begin{align*}
    \Res(\Cnd(M)) &\cong M \oplus M^*, \\
    \Cnd(\Res(N)) &\cong N \oplus \overline{N}, \\
    \Cnd(M) &\cong \Cnd(M^*), \\
    \Res(N) &\cong \Res(\overline{N}).
\end{align*}
for all $M \in \Mod_A$, $N \in \Mod_B$.
If  $M_1, M_2 \in \Rep (A)$ and $N_1, N_2 \in \Rep (B)$, where $\Rep (A)$ denotes the full subcategory of modules finite-dimensional over $\bK$, and similarly for $B$, we set
\begin{align*}
    (M_1, M_2)_A &\coloneqq \dim_{\bbK} \Hom_A(M_1, M_2), \\
    (N_1, N_2)_B &\coloneqq \dim_{\bbK} \Hom_B(N_1, N_2).
\end{align*}
These distribute over direct sums in each variable and, because $\Res$ is left adjoint to $\Cnd$, satisfy the identity.
\begin{equation} \label{ResCnd}
    (\Res(N), M)_A = (N, \Cnd(M))_B.
\end{equation}
In particular, using~\eqref{ResCnd} and the properties above, we arrive at the  identities
\begin{align}
    (\Res(N), \Res(N))_A &= (N,N)_B + (N, \overline{N})_B, \label{eqn:innerprod1} \\
    (\Cnd(M), \Cnd(M))_B &= (M,M)_A + (M, M^*)_A. \label{eqn:innerprod2}
\end{align}

\begin{prop}\label{prop:algindressummary}
Suppose that $A$ and $B$ are as described at the start of this section, $M \in \Rep(A)$ and $N \in \Rep(B)$, finite-dimensional over $\bK$, and that $0 \neq (\Res(N),M)_A$, or equivalently, $0 \neq (N, \Cnd(M))_B$. Then there are four mutually exclusive possibilities, described in Table \ref{tab:indres}.

In particular, $\Res(N)$ and $\Cnd(M)$ are semisimple, the factors of $\Cnd(M)$ are a subset of $\{N, \overline{N}\}$, and the factors of $\Res(N)$ are a subset of $\{M, M^*\}$.
\end{prop}

\begin{table}
\caption{Induction and Restriction Possibilities Between $A$ and $B$}
\label{tab:indres}
\begin{center}
\def\arraystretch{1.1}
\begin{tabular}{ c c c c c } 
    \toprule
Case & \textbf{1} & \textbf{2} & \textbf{3(a)} & \textbf{3(b)} \\
\midrule
 $\Cnd(M)$ & $N \oplus \overline{N}$ & $N$ & $N \oplus N$ & $N$\\ 
 $\Res(N)$ & $M$ & $M \oplus M^*$ & $M$ &  $M \oplus M$ \\ 
 $M \cong M^*$? & Yes & No & Yes & Yes \\ 
 $N \cong \overline{N}$? & No & Yes & Yes & Yes \\ 
  $(M,M)_A : (N,N)_B$ & $1:1$ & $1:1$ & $2:1$ & $1:2$ \\ 
  \bottomrule
\end{tabular}
\end{center} 
\end{table}

\begin{proof}
Suppose that $M \in \Rep(A)$ and $N \in \Rep(B)$, so $M$ and $N$ are finite dimensional over $\bK$, and that $0 \neq (\Res(N),M)_A$, or equivalently, $0 \neq (N, \Cnd(M))_B$. We have the following cases.

\vspace{0.5em}
\textbf{Case 1:} $N \not\cong \overline{N}$. Then because $\Res(N) = \Res(\overline{N})$,
\begin{equation} \label{NoNCndM}
N \oplus \overline{N} \hookrightarrow \Cnd(M),
\end{equation}
and so because $\Res$ is faithful,
\[
\Res(N) \oplus \Res(N) \cong \Res(N \oplus \overline{N}) \hookrightarrow \Res(\Cnd(M)) \cong M \oplus M^*.
\]
Therefore $\Res(N) \cong M$ and $M \cong M^*$. Furthermore, the map in~\eqref{NoNCndM} must have been an isomorphism, and so $\Cnd(M) \cong N \oplus \overline{N}$. Comparing dimensions of endomorphism rings,
\begin{align*}
    2(M,M)_A &= (M,M)_A + (M,M^*)_A, \\
    &= (\Cnd(M), \Cnd(M))_B, \\
    &= (N \oplus \overline{N}, \Cnd(M))_B, \\
    &= (N,\Cnd(M))_B + (\overline{N}, \Cnd(M))_B, \\
    &= 2(N,\Cnd(M))_B, \\
    &= 2 (N,N)_B,
\end{align*}
hence $(M,M)_A = (N,N)_B$. Note that we have used the equalities
\[
(\overline{N}, \Cnd(M))_B = (\Res(\overline{N}), M)_A = (\Res(N), M)_A = (N, \Cnd(M)).
\]

\vspace{0.5em}
\textbf{Case 2:} $M \not\cong M^*$. We have that $N \hookrightarrow \Cnd(M)$, and therefore because $\Res$ is faithful,
\[
\Res(N) \hookrightarrow \Res(\Cnd(M)) \cong M \oplus M^*,
\]
and $\Res(N)$ is semisimple. Because $\Cnd(M) = \Cnd(M^*)$, $(\Res(N), M) = (\Res(N),M^*)$, hence
\[
M \oplus M^* \hookrightarrow \Res(N)
\]
because $\Res(N)$ is semisimple, and thus $\Res(N) \cong M \oplus M^*$. In particular,
\[
N \oplus \overline{N} \cong \Cnd(\Res(N)) \cong \Cnd(M) \oplus \Cnd(M^*),
\]
and $\Cnd(M) \cong N$, $N \cong \overline{N}$. Similarly to the previous case, $(M,M)_A = (N,N)_B$.

\vspace{0.5em}
\textbf{Case 3:} $N \cong \overline{N}$ and $M \cong M^*$. We have that $N \hookrightarrow \Cnd(M)$, and so, $\Res$ being faithful,
\[
\Res(N) \hookrightarrow \Res(\Cnd(M)) \cong M \oplus M.
\]
Therefore, there are two options.

\vspace{0.5em}
\textbf{Case 3(a):} $\Res(N) = M$. Then
\[
\Cnd(M) = \Cnd(\Res(N)) = N \oplus N.
\]
In this case, $(M,M)_A = 2 (N,N)_B$, as
\begin{align*}
    2(N,N)_B &= (N,N)_B + (N,\overline{N})_B, \\
    &= (\Res(N), \Res(N))_A, \\
    &= (M,M)_A.
\end{align*}

\vspace{0.5em}
\textbf{Case 3(b):} $\Res(N) = M \oplus M$. Then
\[
\Cnd(M) \oplus \Cnd(M) = \Cnd(\Res(N)) = N \oplus N,
\]
so $\Cnd(M) = N$. Similarly, $(N,N)_B = 2 (M,M)_A$.
\end{proof}

\section{Examples From the Diamond}\label{sect:examplesfromdiamond} Suppose from now we are in the context described at the start of this chapter, so $\bL / \bK$ is a Galois extension of degree two, and $H \leq G$ is a Real structure on $H$. We also fix an element $\omega \in G \setminus H$.

In this section we concentrate on understanding better the structure of the \emph{diamond}, by which we mean diagram (\ref{eqn:diamond}). More precisely, in this section we show that when $\charfield(\bK)\neq 2$, all six extensions of the diamond fit into the formalism developed in the previous section.

All six extensions of the diamond are Frobenius by Proposition \ref{prop:frobenius}, and so in what follows we consider $\Ind = B \otimes_A (-)$, rather than coinduction $\Cnd$ which we used above.

\subsection{The Extension $\bL[H] \hookrightarrow \bL_{\star \hspace{-0.1em}}[G]$}

Here the duality on $\Mod_{\bL_{\star \hspace{-0.1em}}[G]}$ is the identity functor, and the duality on $\Mod_{\bL[H]}$ sends $V \mapsto {}^{\omega}\overline{V}$. That this satisfies the required identities follows from \cite[Cor.\ 8.5]{TAYSchur}, \cite[Prop.\ 8.17]{TAYSchur} and Theorem \ref{thm:facts}. Because the duality on $\Mod_{\bL_{\star \hspace{-0.1em}}[G]}$ is trivial Case 1 doesn't occur. The remaining three cases recover the three cases of Table \ref{tab:thewaysRHSsquare}: Case 2 is Case (B), Case 3(a) is Case (A) and Case 3(b) is Case (C).

\subsection{The Extensions $\bK[H] \hookrightarrow \bL[H]$ and $\bK[G] \hookrightarrow \bL[G]$}\label{sect:classicalindextwofieldextoptions} Let us describe everything for $\bK[H] \hookrightarrow \bL[H]$, the case for $\bK[G] \hookrightarrow \bL[G]$ obtained by replacing $H$ by $G$ everywhere. This can be viewed as a special case of the previous example in light of \cite[\S10]{TAYSchur}.

For $\bK[H] \hookrightarrow \bL[H]$, the duality on $\Mod_{\bK[H]}$ is the identity functor, and on $\Mod_{\bL[H]}$ sends $V \mapsto \overline{V}$. Now it is the duality on $\Mod_{\bK[H]}$ that is trivial, so Case 2 doesn't occur.

We also have three cases arising from Theorem \ref{thm:facts}: if $X \in \Irr_{\bbK}(H)$ and $V \in \Irr_{\bbL}(H)$ with $X$ a submodule of $V|_{\bK[G]}$, or equivalently $V$ a submodule of $\bL \otimes_{\bbK} X$, then the three options from Theorem \ref{thm:facts} are summarised in Table \ref{tab:thewaysLHSsquare}. The three options from Proposition \ref{prop:algindressummary} recover these three cases: Case 1 is Case (B), Case 3(a) is Case (C) and Case 3(b) is Case (A).

\begin{table}
\caption{Induction and Restriction Possibilities Between $\bK[H]$ and $\bL[H]$}
\label{tab:thewaysLHSsquare}
\begin{center}
\def\arraystretch{1.3}
\begin{tabular}{ c c c c } 
    \toprule
Case & \textbf{A} & \textbf{B} & \textbf{C} \\
\midrule
 $\bL \otimes_{\bbK} X$ & $V$ & $V \oplus \overline{V}$ & $V \oplus V$ \\ 
 $V|_{\bK[H]}$ & $X \oplus X$ & $X$ & $X$\\ 
 $V \cong \overline{V}$? & Yes & No & Yes \\ 
  $m_c(V)$ & 1 & 1 & 2 \\
  $\dim_{\bbL} \End_{\bL[H]}(V) : \dim_{\bbK} \End_{\bK[H]}(X)$ & $1:1$ & $1:2$ & $1:4$ \\ 
  \bottomrule
\end{tabular}
\end{center} 
\end{table}

\subsection{The Extensions $\bK[H] \hookrightarrow \bK[G]$ and $\bL[H] \hookrightarrow \bL[G]$}\label{sect:indextwosubgroup}
Let us describe everything for $\bK[H] \hookrightarrow \bK[G]$, the case for $\bL[H] \hookrightarrow \bL[G]$ obtained by replacing $\bK$ by $\bL$ everywhere.

Write $\pi$ for the $1$-dimensional sign representation of $G$ over $H$, which takes value $1$ on $H$ and $-1$ on $G \setminus H$. The duality on $\Mod_{\bK[G]}$ is 
\[
(-) \otimes \pi \colon \Mod_{\bK[G]} \rightarrow \Mod_{\bK[G]}.
\]
It is direct to see that $\Res(V) = \Res(V \otimes \pi)$. The duality on $\Mod_{\bK[H]}$ is
\[
{}^\omega(-) \colon \Mod_{\bK[H]} \rightarrow \Mod_{\bK[H]},
\]
where ${}^{\omega} W$ is $W$ as a $\bK$-vector space and action of $h \in H$ by the action of $\omega^{-1} h \omega$ on $W$. There is an isomorphism of $\bK[G]$-modules
\[
    \bK[G] \otimes_{\bK[H]} {}^{\omega} W \xrightarrow{\sim} \bK[G] \otimes_{\bK[H]} W, \qquad g \otimes w \mapsto \omega^{-1} g \omega \otimes w
\]
and so $\Ind(W) = \Ind({}^{\omega} W)$. The following shows that the other properties also hold when $\charfield(\bK) \neq 2$.

\begin{prop}\label{prop:resindKHG}
Suppose that $W \in \Mod_{\bK[H]}$, $U \in \Mod_{\bK[G]}$, and $\charfield(K) \neq 2$. Then 
\begin{align*}
    \Res(\Ind(W)) &\cong W \oplus {}^{\omega} W, \\
    \Ind(\Res(U)) &\cong U \oplus (U \otimes \pi).
\end{align*}
\end{prop}

\begin{proof}
    The first isomorphism is given by the $K[H]$-module identifications
    \[  
        \bK[G] \otimes_{\bK[H]} W \equiv (\bK[H] \oplus \bK[H] \omega) \otimes_{\bK[H]} W \equiv W \oplus {}^{\omega} W.
    \]
    The second is given by the $K[G]$-module identification
    \[
        \bK[G] \otimes_{\bK[H]} U \xrightarrow{\sim} \bK[G/H] \otimes_{\bbK} U, \qquad g \otimes u \mapsto g \otimes g(u)
    \]
    and noting that the representation $\bK[G/H] = \mathbbm{1} \oplus \pi$ as $\bK[G]$-modules because $\charfield(\bK) \neq 2$.
\end{proof}

\begin{remark}
    Proposition \ref{prop:resindKHG} is false when $\charfield(\bK) \neq 2$: if it were true, then by Proposition \ref{prop:algindressummary} $\Ind(W)$ would always be semisimple for any simple  $W \in \Mod_{\bK[H]}$. However, in characteristic $2$ this isn't true when one takes $W$ to be the trivial representation.
\end{remark}

\subsection{The Extension $\bK[G] \hookrightarrow \bL_{\star \hspace{-0.1em}}[G]$} As in the previous section, we have the duality
\[
(-) \otimes \pi \colon \Mod_{\bK[G]} \rightarrow \Mod_{\bK[G]}.
\]
We have that the $\bK[G]$-module $\pi$ maps under $\Ind$ to the trivial representation of $\bL_{\star \hspace{-0.1em}}[G]$, by \cite[Thm.\ A]{TAYSchur} and the fact that $\Ind(\pi)$ restricts to the trivial representation as a $\bL[H]$-module. Therefore because $\Ind = \bL \otimes_{\bbK} -$ is a tensor functor (see \cite[\S4.1]{TAYSchur}), $\Ind(U) = \Ind(U \otimes \pi)$.

For $\Mod_{\bL_{\star \hspace{-0.1em}}[G]}$, we have a $\bK$-algebra automorphism $\overline{\cdot} \colon \bL_{\star \hspace{-0.1em}}[G] \rightarrow \bL_{\star \hspace{-0.1em}}[G]$ sending $\lambda g \mapsto \overline{\lambda} g$ which is well-defined because $\Gamma$ is abelian. We can therefore precompose any $\bL_{\star \hspace{-0.1em}}[G]$-module to with $\overline{\cdot}$, and to obtain the duality
\[
\overline{(-)} \colon \Mod_{\bL_{\star \hspace{-0.1em}}[G]} \rightarrow \Mod_{\bL_{\star \hspace{-0.1em}}[G]}.
\]
As the automorphism $\overline{\cdot}$ of $\bL_{\star \hspace{-0.1em}}[G]$ fixes $\bK[G]$, we have that $\Res(X) = \Res(\overline{X})$ for any $X \in \Mod_{\bL_{\star \hspace{-0.1em}}[G]}$. We can again check that the other properties hold when $\charfield(\bK) \neq 2$.

\begin{prop}\label{prop:resindKGLG}
Suppose that $U \in \Mod_{\bK[G]}$, $X \in \Mod_{\bL_{\star \hspace{-0.1em}}[G]}$, and $\charfield(\bK) \neq 2$. Then 
\begin{align*}
    \Res(\Ind(U)) &\cong U \oplus (U \otimes \pi), \\
    \Ind(\Res(X)) &\cong X \oplus \overline{X}.
\end{align*}
\end{prop}

\begin{proof}
    For the first isomorphism, $\bL \cong \bK[\Gamma]$ as a $\bK[G]$-module because $\bL/\bK$ is Galois, and $\bK[\Gamma] = \mathbbm{1} \oplus \pi$ because $\charfield(\bK) \neq 2$. Therefore
    \[
        \Res(\Ind(U)) = \bL \otimes_{\bbK} U = \bK[\Gamma] \otimes_{\bbK} U = U \oplus (U \otimes \pi).
    \]
    The second isomorphism is given by
    \[
        \bL \otimes_{\bbK} X \rightarrow X \oplus \overline{X}, \qquad \lambda \otimes x \mapsto (\lambda x, \overline{\lambda} x),
    \]
    which is easily verified to be $\bL_{\star \hspace{-0.1em}}[G]$-linear. To show it is bijective, by dimension considerations it suffices to show surjectivity. To this end, let $e = \sum_i \alpha_i \otimes \beta_i \in \bL \otimes_{\bbK} \bL$ be the unique element which maps to $(1,0)$ under the Galois isomorphism $\bL \otimes_{\bbK} \bL \rightarrow \bL \oplus \bL$, $a \otimes b \mapsto (ab,\overline{a} b)$. Then for any $x \in X$, 
    \[
        \sum_i \alpha_i \otimes \beta_i x \mapsto \left(\sum_i \alpha_i \beta_i x, \sum_i \overline{\alpha_i} \beta_i x \right) = (x,0).
    \]
    Similarly considering $(0,1)$ we can obtain any element $(0,x)$ and hence the map is surjective.
\end{proof}

\section{Analysis of the Diamond}\label{sect:diamondanalysis}

We continue with the notation from above, so $\bL / \bK$ is a Galois extension of degree two, $H \leq G$ is a Real structure on $H$, and $\omega$ is a fixed element of $G \setminus H$. Furthermore:

\vspace{0.5em}
\begin{center}
\emph{From now on in this chapter we additionally assume that $\charfield(\bK) \neq 2$.}
\end{center}
\vspace{0.0em}

\subsection{Group Actions on Sets of Irreducible Representations}

Let us start with the following version of the diamond~\eqref{eqn:diamond}, where we restrict to finite-dimensional irreducible representations:
\begin{equation}\label{eqn:irrdiamond}
\begin{tikzcd}
	& {  \Irr_{\bbK}(G) \, \rotatebox[origin=c]{90}{$\circlearrowleft$} \, C_2} & \\
  {{\Irr_{\bbL}(G)} \, \rotatebox[origin=c]{90}{$\circlearrowleft$} \, K_4}
  &
  {\Irr_{\bbK}(H) \, \rotatebox[origin=c]{90}{$\circlearrowleft$} \, C_2}
    &
    { {\Irr_{\bbL}^{\hspace{0.1em} \star}(G)}\, \rotatebox[origin=c]{90}{$\circlearrowleft$} \, C_2} \\
     & {\Irr_{\bbL}(H) \, \rotatebox[origin=c]{90}{$\circlearrowleft$} \, K_4}
	\arrow[tail reversed, from=1-2, to=2-1]
	\arrow[tail reversed, from=1-2, to=2-3]
	\arrow[tail reversed, from=2-2, to=1-2]
	\arrow[tail reversed, from=2-2, to=3-2]
	\arrow[tail reversed, from=3-2, to=2-1]
	\arrow[tail reversed, from=3-2, to=2-3]
\end{tikzcd}
\end{equation}

To describe the actions of the groups $K_4 = \{1,a,b,c\}$ and $C_2 = \{1,\epsilon\}$, we describe how these groups act on the sets of (isomorphism classes of) all representations, which restrict to the actions on the sets of (isomorphism classes of) irreducible representations above.

First, pick $Z = (Z,\rho) \in \Rep_{\bbL}(G)$. We can obtain three new representations:
\begin{itemize}
\item the Galois twist $\overline{Z}$ (as described in Section \ref{sect:chap3notation}), with the same action of $G$,
\item the twist ${}^{\omega}{Z} = (Z,{}^{\omega}\rho)$ where ${}^{\omega}\rho (x) = \rho (\omega^{-1} x \omega)$ for all $x\in G$. 
\item the tensor product $Z \otimes \pi \coloneqq Z \otimes_{\bbL} \pi$ with the sign character $\pi : G \to \{ \pm 1\}$. 
\end{itemize}  
The twist by $\omega$ is an inner twist, hence $\rho(\omega) \colon Z \xrightarrow{\sim} {}^{\omega}{Z}$ is an isomorphism and this operation gives nothing new. The remaining two operations commute, and we thus have an action of $K_4$:
\[
    a * Z = \overline{Z}, \qquad b*Z = Z \otimes \pi, \qquad c * Z = \overline{Z} \otimes \pi.
    \]

Consider now $V \in \Rep_{\bbL}(H)$. Here $\pi$ is not available, so the group $K_4$ acts by (cf.\ Section~\ref{sect:C2actionIrrH})
\[
    a * V = \overline{V}, \qquad b * V = {}^{\omega} V, \qquad c * V = \epsilon * V = {}^{\omega}\overline{V},
\]

As is true of any transitive $K_4$-set, the structure of the $K_4$-set $\{V, \overline{V}, {}^{\omega} V, {}^{\omega} \overline{V}\}$ is uniquely determined by $\Stab_{K_4}(V)$. Similarly, the $K_4$-set $\{Z, \overline{Z}, Z \otimes \pi, \overline{Z} \otimes \pi\}$ is determined by $\Stab_{K_4}(Z)$.

Now pick $U  \in \Rep_{\bbK}(G)$. Here the Galois twist is no longer available, while the twist by $\omega$ similarly produces an isomorphic representation. Hence, we have a $C_2$-action
\[
    \epsilon  * U =  U \otimes \pi \coloneqq U \otimes_{\bbK} \pi.
    \]

    Now consider $X  \in \Rep_{\bbK}(H)$. Here the Galois twist and the sign character are off limits.
All that remains is the twist by $\omega$, which yields a $C_2$-action
\[
    \epsilon  * X =  {}^{\omega}X.
    \]

    Finally, let us examine $W = (W,\rho) \in {\Irr_{\bbL}^{\hspace{0.1em} \star}(G)}$. All three operations appear available.
    However, we first note that $W \cong W \otimes_{\bbL} \pi$, because the semilinear sign representation $(\bbL, \pi)$ is isomorphic to the trivial semilinear representation $(\bbL, 1)$.
    Indeed, pick any non-zero ``imaginary'' element, that is, an $\alpha\in \bbL$ with $\overline{\alpha} = - \alpha \neq 0$. Then multiplication by $\alpha$ is an isomorphism $(\bbL, 1)\to (\bbL, \pi)$.

    What remains are the twist ${}^{\omega} W$ by $\omega$, and the conjugate representation $\overline{W}$, which defines a semilinear representation with the same action of $G$ because $\Gal(\bL / \bK)$ is abelian: for any $w \in W$, $\alpha \in \bL$, and $g \in G$,
    \[
    \rho(g)(\alpha * v) = \rho(g)(\overline{\alpha}v) = \sigma(g)(\overline{\alpha}) \rho(g)(v) = \overline{\sigma(g)(\alpha)} \rho(g)(v) = \sigma(g)(\alpha) * \rho(g)(v).
    \]
    
    Although both of ${}^{\omega} W$ and $\overline{W}$ differ from $W$ in general, they turn out to be the same: $\rho (\omega)$ can be viewed as an isomorphism $\overline{W} \to {}^{\omega}W$. Hence, we have only a $C_2$-action
\[
    \epsilon  * W =  {}^{\omega}W = \overline{W}.
    \]

\subsection{Propagation}\label{sect:propagation}

Let $V \in \Irr_{\bbL}(H)$. From this we can ``propagate out'' along the diamond~\eqref{eqn:irrdiamond} and define: \vspace{0.5em}
\begin{itemize}
    \item $W \in \Irr_{\bbL}^{\hspace{0.1em} \star}(G)$, an irreducible submodule of $\Ind_{\bL[H]}^{\bL_{\star \hspace{-0.1em}}[G]}(V)$, \vspace{0.5em}
    \item $X \in \Irr_{\bbK}(H)$, an irreducible submodule of $V|_{\bK[H]}$, \vspace{0.5em}
    \item $Z \in \Irr_{\bbL}(G)$, an irreducible submodule of $\Ind_{\bL[H]}^{\bL[G]}(V)$, \vspace{0.5em}
    \item $U \in \Irr_{\bbK}(G)$, an irreducible submodule of $Z|_{\bK[G]}$.
\end{itemize}

Using these modules, we can restrict the diamond~\eqref{eqn:diamond} to orbits under the corresponding groups.
The meaning of the labelling of the arrows will become clear in the next section.
\begin{equation}\label{eqn:labelleddiamond}
\begin{tikzcd}
	& {  \{U,U\otimes \pi \} \, \rotatebox[origin=c]{90}{$\circlearrowleft$} \, C_2} & \\
  {\{Z,\overline{Z},Z\otimes \pi ,\overline{Z}\otimes \pi  \} \, \rotatebox[origin=c]{90}{$\circlearrowleft$} \, K_4}
  &
  { \{X, {}^{\omega}X\} \, \rotatebox[origin=c]{90}{$\circlearrowleft$} \, C_2}
    &
    {\{W, \overline{W} \} \, \rotatebox[origin=c]{90}{$\circlearrowleft$} \, C_2} \\
    & {\{V,\overline{V},{}^{\omega}V,{}^{\omega}\overline{V} \} \, \rotatebox[origin=c]{90}{$\circlearrowleft$} \, K_4}
	\arrow["3"', tail reversed, from=1-2, to=2-1]
	\arrow["4", tail reversed, from=1-2, to=2-3]
	\arrow["4", tail reversed, from=2-2, to=1-2]
	\arrow["3"', tail reversed, from=2-2, to=3-2]
	\arrow["4", tail reversed, from=3-2, to=2-1]
	\arrow["3"', tail reversed, from=3-2, to=2-3]
\end{tikzcd}
\end{equation}

We have used four of the double arrows in the diamond~\eqref{eqn:irrdiamond} to define these modules. Equally, we could have used the 
remaining two, as clear from the next fact.
\begin{lemma}
    $U$ is a submodule of $\Ind_{\bK[H]}^{\bK[G]}(X)$ and $W$ is a submodule of $\Ind_{\bK[G]}^{\bL_{\star \hspace{-0.1em}}[G]}(U)$.
\end{lemma}

\begin{proof}
    For the first claim, first note that
    \[
    \left(Z, \Ind_{\bK[H]}^{\bL[G]}(X)\right) \neq 0,
    \]
    as $Z$ is a submodule of $\Ind_{\bL[H]}^{\bL[G]}(V)$ and $V \subset \Ind_{\bK[H]}^{\bL[H]}(X)$ by Frobenius reciprocity. Then 
    \[
        \left(Z, \Ind_{\bK[G]}^{\bL[G]}\left(\Ind_{\bK[H]}^{\bK[G]}(X)\right) \right) \neq 0
    \]
    so again by Frobenius reciprocity
    \[
    \left(Z|_{\bK[G]}, \Ind_{\bK[H]}^{\bK[G]}(X)\right) \neq 0.
    \]
    As $Z|_{\bK[G]}$ is a direct sum of copies of $U$ the first claim follows.

    For the second claim, we have that
    \[
    \left(W, \Ind_{\bK[H]}^{\bL_{\star \hspace{-0.1em}}[G]}(X)\right) \neq 0,
    \]
    as $V \subset \Ind_{\bK[H]}^{\bL[H]}(X)$ by Frobenius reciprocity and $W \subset \Ind_{\bL[H]}^{\bL_{\star \hspace{-0.1em}}[G]}(V)$. We have that
    \[
        \left(W, \Ind_{\bK[G]}^{\bL_{\star \hspace{-0.1em}}[G]}\left(\Ind_{\bK[H]}^{\bK[G]}(X)\right) \right) \neq 0,
    \]
    that $\Ind_{\bK[H]}^{\bK[G]}(X)$ is a direct sum of $U$ and $U \otimes \pi$ in equal multiplicity, and that
    \[
    \Ind_{\bK[G]}^{\bL_{\star \hspace{-0.1em}}[G]}(U) = \Ind_{\bK[G]}^{\bL_{\star \hspace{-0.1em}}[G]}(U \otimes \pi),
    \]
    hence
    \[
    \left(W, \Ind_{\bK[G]}^{\bL_{\star \hspace{-0.1em}}[G]}(U) \right) \neq 0,
    \]
    which establishes the claim.
\end{proof}

\subsection{The type and edge labels}
Each pair of modules adjacent in the diamond~\eqref{eqn:labelleddiamond} falls into one of the four cases 1, 2, 3(a) and 3(b) of Section \ref{sect:fourfoldway}. For three of these pairs there are only three possibilities: we label those instead by A, B and C as in Tables \ref{tab:thewaysRHSsquare} and \ref{tab:thewaysLHSsquare}. The double edges of the diamond~\eqref{eqn:labelleddiamond} are marked with the number of possibilities for each edge
and
we can record the induction and restriction behaviour of this family of modules as an element of
\[
T \coloneqq \{\text{A,\:B,\:C}\}^3 \times \{\text{1,\:2,\:3(a),\:3(b)}\}^3
\]
corresponding to the ordered extensions
\[
(\bK[G] \subset \bL[G], \: \bK[H] \subset \bL[H], \: \bL[H] \subset \bL_{\star \hspace{-0.1em}}[G]) \times (\bL[H] \subset \bL[G], \: \bK[H] \subset \bK[G], \: \bK[G] \subset \bL_{\star \hspace{-0.1em}}[G]).
\]
The ordering here is chosen to match the ordering obtained by reading the diamond from left to right.

\begin{defn}
The \emph{type} of $V$, $\sT(V)$, is the associated element of $T$.
\end{defn}

\begin{lemma}\label{lem:typeindependence}
    The type of $V$ is independent of the choices of the irreducible representations $X$, $Z$, $W$, $U$, and depends only on the $K_4$-orbit of $V$.
\end{lemma}

\begin{proof}
    From the explicit definition of each action, on any orbit of Diagram \eqref{eqn:labelleddiamond} there is an identification of endomorphism rings. Furthermore, for any extension of Diagram \eqref{eqn:labelleddiamond}, irreducible representations in the same orbit will have corresponding irreducible representations in the same orbit. Therefore, for any irreducible representation, and any extension of Diagram \eqref{eqn:labelleddiamond}, the endomorphism ring ratios of Proposition \ref{prop:algindressummary} will be constant on orbits, as will the properties of being fixed by the $C_2$-action of Proposition \ref{prop:algindressummary} (because $K_4$ and $C_2$ are abelian), and therefore the type for each extension will also be constant on orbits.    
\end{proof}

\begin{remark}
    We have chosen to start with $V$, but we may equally define the type of any other irreducible representation from any of the other four categories of Diagram \eqref{eqn:diamond}, which similarly will be independent of the choices involved and depend only on its orbit. We give an alternative point of view on the type in terms of the notion of an \emph{S-block} in Section \ref{sect:blockinterpretation} below.
\end{remark}

A priori, many possible types can occur, as the set $T$ has order $3^3 \cdot 4^3 = 1728$. The main result of this section is the following classification, which shows that only 24 possibilities can occur.

\begin{thm}\label{thm:24ways}
The possible types for $V$ that can occur are those 24 types listed in Table \ref{tab:24types}.
\end{thm}

In Section \ref{sect:realise24ways} we consider the question of realisation of the 24 types of Table \ref{tab:24types}. In the rest of this section we prove Theorem \ref{thm:24ways}. First we make some useful observations.

\begin{lemma}\label{lem:UnotUpi}
    If $U \not\cong U \otimes \pi$, then the type of $V$ is one of
    \emph{\begin{itemize}
        \item \textbf{(A, A, A, 1, 1, 2)},
        \item \textbf{(B, B, B, 1, 1, 2)},
        \item \textbf{(C, C, C, 1, 1, 2)}.
    \end{itemize}}
\end{lemma}
\begin{proof}
    If $U \not\cong U \otimes \pi$, then $\bK[H] \subset \bK[G]$ is in Case 1, and $\bK[G] \subset \bL_{\star \hspace{-0.1em}}[G]$ is in Case 2. As $\Res(Z \otimes \pi) = \Res(Z) \otimes \pi$, we have that $Z \not\cong Z \otimes \pi$, so $\bL[H] \subset \bL[G]$ is in Case 1. From this, we have that $(W,W) = (U,U) = (X,X)$ and $(Z,Z) = (V,V)$. The three remaining cases for each of the remaining extensions are determined by the ratio $(V,V) : (X,X)$ on endomorphism rings, and a choice for one extension determines the other two. The three possibilities are listed above.
\end{proof}

\begin{lemma}\label{lem:WnotomegaW}
    If $X \not\cong {}^{\omega}X$, then the type of $V$ is one of
    \emph{\begin{itemize}
        \item \textbf{(A, A, B, 2, 2, 3(b))},
        \item \textbf{(B, B, B, 2, 2, 1)},
        \item \textbf{(C, C, B, 2, 2, 3(a))}.
    \end{itemize}}
\end{lemma}
\begin{proof}
    If $X \not\cong {}^{\omega}X$, then $\bK[H] \subset \bK[G]$ is in Case 2, and $(U,U) = (X,X)$. As $X \not\cong {}^{\omega}X$, $\{V, \overline{V}\} \neq \{{}^{\omega}V, {}^{\omega}\overline{V}\}$, so $V \not\cong {}^{\omega} \overline{V}$, meaning that $\bL[H] \subset \bL_{\star \hspace{-0.1em}}[G]$ is in Case B and $V \not\cong {}^{\omega} \overline{V}$, meaning that $\bL[H] \subset \bL[G]$ is in Case 2. So the type is (?, ?, B, 2, 2, ?), and $(W,W) = (V,V) = (Z,Z)$. We now consider the possible options for $\bK[H] \subset \bK[G]$, and show these determine the rest.
    \begin{itemize}
        \item In Case A, $2(X,X) = (V,V)$, so $2(U,U) = (W,W)$, and $\bK[G] \subset \bL_{\star \hspace{-0.1em}}[G]$ is in Case 3(b). Furthermore, $2(U,U) = (Z,Z)$, so $\bK[G] \subset \bL[G]$ is in Case A.
        \item In Case B, $(X,X) = (V,V)$, so $(U,U) = (W,W)$, and as $U \cong U \otimes \pi$, $\bK[G] \subset \bL_{\star \hspace{-0.1em}}[G]$ is in Case 1. Furthermore, $(U,U) = (Z,Z)$, so $\bK[G] \subset \bL[G]$ is in Case B.
        \item In Case C, $(X,X) = 2(V,V)$, so $(U,U) = 2(W,W)$, and $\bK[G] \subset \bL_{\star \hspace{-0.1em}}[G]$ is in Case 3(a). Furthermore, $(U,U) = 2(Z,Z)$, so $\bK[G] \subset \bL[G]$ is in Case C. \qedhere
    \end{itemize}
\end{proof}

\subsection{Proof of Theorem \ref{thm:24ways}} We now being our analysis of the possible types of $V$ in earnest. It is essentially a case analysis, indexed first by the extension $\bK[H] \subset \bL[H]$ and then by $\bK[H] \subset \bK[G]$.

\vspace{0.5em}
\textcolor{gray}{\hrule}

\vspace{0.9em}
\begin{center}
    \large{\textbf{Case A} (?, A, ?, ?, ?, ?)}
\end{center}
\[
 V \cong \overline{V}, \ \Ind(X) = V, \ \Res(V) = X \oplus X, \ (X,X) = \dim_{\bbL} \End_{\bL[H]}(V).
\]
\vspace{-0.3em}
\textcolor{gray}{\hrule}
\vspace{0.5em}

\vspace{0.5em}
\noindent \textbf{Case A.1} (?, A, ?, ?, 1, ?)

From Lemma \ref{lem:UnotUpi}, the type of $V$ must be \textbf{(A, A, A, 1, 1, 2)}.

\vspace{0.5em}
\noindent \textbf{Case A.2} (?, A, ?, ?, 2, ?)

From Lemma \ref{lem:WnotomegaW}, the type of $V$ must be \textbf{(A, A, B, 2, 2, 3(b))}.

\vspace{0.5em}
\noindent \textbf{Case A.3(a)} (?, A, ?, ?, 3(a), ?)

In this case we have that
\[
U \cong U \otimes \pi, \ X \cong {}^{\omega} X,\ \Ind(X) = U \oplus U, \ \Res(U) = X, \ (X,X) = 2(U,U).
\]
Because $X \cong {}^{\omega} X$ and $\Ind(X) \cong V$, $V \cong {}^{\omega} V$, and therefore, as $V \cong \overline{V}$, $\Stab_{K_4}(V) = K_4$. In particular, $\bL[H] \subset \bL_{\star \hspace{-0.1em}}[G]$ is either in Case A or Case C. In Case C, we would have $(W,W) = 4 \dim_{\bbL} \End_{\bL[H]}(V)$, and thus would have $(W,W) = 4 (X,X) = 8(U,U)$, which cannot occur. Therefore, we must be in Case A, with type (?, A, A, ?, 3(a), ?) and $(W,W) = \dim_{\bbL} \End_{\bL[H]}(V)$. In particular,
\[
(W,W) = \dim_{\bbL} \End_{\bL[H]}(V) = (X,X) = 2(U,U),
\]
and so $\bK[G] \subset \bL_{\star \hspace{-0.1em}}[G]$ is in Case 3(b) and the type is (?, A, A, ?, 3(a), 3(b)).

There are then three options for the extension $\bK[G] \subset \bL[G]$.
\begin{itemize}
\item In Case A, type (A, A, A, ?, 3(a), 3(b)), $(U,U) = \dim_{\bbL}\End_{\bL[G]}(Z)$. Therefore
    \[
    \dim_{\bbL}\End_{\bL[H]}(V) = (X,X) = 2(U,U) = 2 \dim_{\bbL}\End_{\bL[G]}(Z),
    \]
    hence $\bL[H] \subset \bL[G]$ is in Case 3(a). The type is \textbf{(A, A, A, 3(a), 3(a), 3(b))}.
\item In Case B, $(U,U) = 2\dim_{\bbL}\End_{\bL[G]}(Z)$, hence
\[
\dim_{\bbL}\End_{\bL[H]}(V) = (X,X) = 2(U,U) = 4 \dim_{\bbL}\End_{\bL[G]}(Z),
\]
which cannot occur.
\item In Case C, $(U,U) = 4\dim_{\bbL}\End_{\bL[G]}(Z)$, hence
\[
\dim_{\bbL}\End_{\bL[H]}(V) = (X,X) = 2(U,U) = 8 \dim_{\bbL}\End_{\bL[G]}(Z),
\]
which cannot occur.
\end{itemize}

\vspace{0.5em}
\noindent \textbf{Case A.3(b)} (?, A, ?, ?, 3(b), ?)

In this case we have that 
\[
U \cong U \otimes \pi, \ X \cong {}^{\omega} X, \ \Ind(X) = U, \ \Res(U) = X \oplus X, \ 2(X,X) = (U,U).
\]
Because $X \cong {}^{\omega} X$ and $\Ind(X) \cong V$, $V \cong {}^{\omega} V$, and therefore, as $V \cong \overline{V}$, $\Stab_{K_4}(V) = K_4$. In particular, $\bL[H] \subset \bL_{\star \hspace{-0.1em}}[G]$ is either in Case A or Case C. Let us first observe that being A or C determines the case for $\bK[G] \subset \bL_{\star \hspace{-0.1em}}[G]$: 

\begin{itemize}
    \item In Case A, $(W,W) = \dim_{\bbL} \End_{\bL[H]}(V)$, so $(U,U) = 2(X,X) = 2 (W,W)$ and $\bK[G] \subset \bL_{\star \hspace{-0.1em}}[G]$ is in Case 3(a). The type is (?, A, A, ?, 3(b), 3(a)).
    \item In Case C, $(W,W) = 4 \dim_{\bbL} \End_{\bL[H]}(V)$, and thus $(W,W) = 4 (X,X) = 2(U,U)$, so $\bK[G] \subset \bL_{\star \hspace{-0.1em}}[G]$ is in Case 3(b). The type is (?, A, C, ?, 3(b), 3(b)).
\end{itemize}

Note also that in either case $(V,V) = 2 \dim_{\bbL} \End_{\bL[H]}(V) = 2(X,X) = (U,U)$. We can use this to see that the three cases for $\bK[G] \subset \bL[G]$ determine those for $\bL[H] \subset \bL[G]$: Case 2 cannot occur for $\bL[H] \subset \bL[G]$, as $V \cong {}^{\omega}V$, and the three remaining cases are determined by $(V,V):(Z,Z)$, which is determined by the case for $\bK[G] \subset \bL[G]$. The three options are (B, A, ?, 1, 3(b), ?), (C, A, ?, 3(a), 3(b), ?), and (A, A, ?, 3(b), 3(b), ?).

Pairing these two options for the extension $\bL[H] \subset \bL_{\star \hspace{-0.1em}}[G]$ with the three options for $\bK[G] \subset \bL[G]$, we obtain six options:
\begin{itemize}
     \item \textbf{(B, A, A, 1, 3(b), 3(a))},
     \item \textbf{(B, A, C, 1, 3(b), 3(b))},
     \item \textbf{(C, A, A, 3(a), 3(b), 3(a))},
     \item \textbf{(C, A, C, 3(a), 3(b), 3(b))},
     \item \textbf{(A, A, A, 3(b), 3(b), 3(a))},
     \item \textbf{(A, A, C, 3(b), 3(b), 3(b))},
\end{itemize}

\vspace{0.5em}
\textcolor{gray}{\hrule}

\vspace{0.9em}
\begin{center}
    \large{\textbf{Case B} (?, B, ?, ?, ?, ?)}
\end{center}
\[
 V \not\cong \overline{V}, \ \Ind(X) = V \oplus \overline{V}, \ \Res(V) = X, \ (X,X) = 2\dim_{\bbL} \End_{\bL[H]}(V).
\]
\vspace{-0.3em}
\textcolor{gray}{\hrule}
\vspace{0.5em}

\vspace{0.5em}
\noindent \textbf{Case B.1} (?, B, ?, ?, 1, ?)

From Lemma \ref{lem:UnotUpi}, the type of $V$ must be \textbf{(B, B, B, 1, 1, 2)}.

\vspace{0.5em}
\noindent \textbf{Case B.2} (?, B, ?, ?, 2, ?)

From Lemma \ref{lem:WnotomegaW}, the type of $V$ must be \textbf{(B, B, B, 2, 2, 1)}.

\vspace{0.5em}
\noindent \textbf{Case B.3(a)} (?, B, ?, ?, 3(a), ?)

In this case we have that 
\[
U \cong U \otimes \pi, \ X \cong {}^{\omega} X, \ \Ind(X) = U \oplus U, \ \Res(U) = X, \ (X,X) = 2(U,U).
\]
We cannot have that $\bL[H] \subset \bL_{\star \hspace{-0.1em}}[G]$ is in Case C, otherwise $4(U,U) = 2(X,X) = 2(V,V) = (W,W)$, which cannot happen for the extension $\bK[G] \subset \bL_{\star \hspace{-0.1em}}[G]$. Therefore we have two options for $\bL[H] \subset \bL_{\star \hspace{-0.1em}}[G]$:
\begin{itemize}
    \item In Case A, $(V,V) = 2(W,W)$, so as $(V,V) = 2(U,U)$, $(W,W) = (U,U)$, and $\bK[G] \subset \bL_{\star \hspace{-0.1em}}[G]$ is in Case 1, with type (?, B, A, ?, 3(a), 1). As $V \not\cong \overline{V}$ and $V \cong {}^{\omega} \overline{V}$, $\Stab_{K_4}(V) = \langle c \rangle$ and $V \not\cong {}^{\omega} V$ and $\bL[H] \subset \bL[G]$ is in Case 2, with $(Z,Z) = (V,V)$. Then $(Z,Z) = 2(U,U)$, so $\bK[G]\subset \bL[G]$ is in Case A, and the type is \textbf{(A, B, A, 2, 3(a), 1)}.
    \item In Case B, $(V,V) = (W,W)$ so as $(V,V) = 2(U,U)$, $(W,W) = 2(U,U)$, and $\bK[G] \subset \bL_{\star \hspace{-0.1em}}[G]$ is in Case 3(b), with type (?, B, B, ?, 3(a), 3(b)). If $Z \cong \overline{Z}$, then $\{V, {}^{\omega} V \} = \{\overline{V}, {}^{\omega}\overline{V}\}$, which is a contradiction, as $V \not\cong \overline{V}$ and $V \not\cong {}^{\omega}\overline{V}$. Therefore $Z \not\cong \overline{Z}$, and $\bK[G] \subset \bL[G]$ is in Case B. This forces $(V,V) = 2(Z,Z)$, so $\bL[H] \subset \bL[G]$ is in Case 3(a). The type is \textbf{(B, B, B, 3(a), 3(a), 3(b))}.
\end{itemize}

\vspace{0.5em}
\noindent \textbf{Case B.3(b)} (?, B, ?, ?, 3(b), ?)

In this we have that
\[
U \cong U \otimes \pi, \ X \cong {}^{\omega} X, \ \Ind(X) = U, \ \Res(U) = X \oplus X, \ 2(X,X) = (U,U).
\]
Suppose first that $W \cong \overline{W}$. Then $\{V, {}^{\omega} \overline{V} \} = \{\overline{V}, {}^{\omega}V\}$, so, as $V \not\cong \overline{V}$, $V \cong {}^{\omega} V$, and thus $\Stab_{K_4}(V) = \langle b \rangle$, so $V \not\cong {}^{\omega} \overline{V}$ and $\bL[H] \subset \bL_{\star \hspace{-0.1em}}[G]$ is in Case B. As $(U,U) = 2(W,W)$, $\bK[G] \subset \bL_{\star \hspace{-0.1em}}[G]$ is in Case 3(a), with type (?, B, B, ?, 3(b), 3(a)). If $Z \cong \overline{Z}$, then $\{V, {}^{\omega} V \} = \{\overline{V}, {}^{\omega}\overline{V}\}$, which is a contradiction, as $V \not\cong \overline{V}$ and $V \not\cong {}^{\omega}\overline{V}$. Therefore $Z \not\cong \overline{Z}$, and $\bK[G] \subset \bL[G]$ is in Case B. Therefore $(Z,Z) = 2(V,V)$, and $\bL[H] \subset \bL[G]$ is in Case 3(b). The type is \textbf{(B, B, B, 3(b), 3(b), 3(a))}.

If instead $W \not\cong \overline{W}$, then $\bK[G] \subset \bL_{\star \hspace{-0.1em}}[G]$ is in Case 1, so $(U,U) = (W,W)$. This means that $(W,W) = 2(V,V)$, and $\bL[H] \subset \bL_{\star \hspace{-0.1em}}[G]$ is in Case C. The type is (?, B, C, ?, 3(b), 1). From this, $\Stab_{K_4}(V) = \langle c \rangle$, hence $V \not\cong {}^{\omega} V$, and $\bL[H] \subset \bL[G]$ is in Case 2, and $2(Z,Z) = (U,U)$. Therefore $\bK[G] \subset \bL[G]$ is in Case C, and the type is \textbf{(C, B, C, 2, 3(b), 1)}.

\vspace{0.5em}
\textcolor{gray}{\hrule}

\vspace{0.9em}
\begin{center}
    \large{\textbf{Case C} (?, C, ?, ?, ?, ?)}
\end{center}
\[
 V \cong \overline{V}, \ \Ind(X) = V \oplus V, \ \Res(V) = X, \ (X,X) = 4\dim_{\bbL} \End_{\bL[H]}(V).
\]
\vspace{-0.3em}
\textcolor{gray}{\hrule}
\vspace{0.5em}

\vspace{0.5em}
\noindent \textbf{Case C.1} (?, C, ?, ?, 1, ?)

From Lemma \ref{lem:UnotUpi}, the type of $V$ must be \textbf{(C, C, C, 1, 1, 2)}.

\vspace{0.5em}
\noindent \textbf{Case C.2} (?, C, ?, ?, 2, ?)

From Lemma \ref{lem:WnotomegaW}, the type of $V$ must be \textbf{(C, C, B, 2, 2, 3(a))}.

\vspace{0.5em}
\noindent \textbf{Case C.3(a)} (?, C, ?, ?, 3(a), ?)

In this case we have that
\[
U \cong U \otimes \pi, \ X \cong {}^{\omega} X, \ \Ind(X) = U \oplus U, \ \Res(U) = X, \ (X,X) = 2(U,U).
\]
Because $X \cong {}^{\omega} X$, $V \cong {}^{\omega} V$, and therefore, as $V \cong \overline{V}$, $\Stab_{K_4}(V) = K_4$. In particular, $\bL[H] \subset \bL_{\star \hspace{-0.1em}}[G]$ is either in Case A or Case C. Let us first observe that this determines the case for $\bK[G] \subset \bL_{\star \hspace{-0.1em}}[G]$: 

$\bullet$ In Case A, $2(W,W) = (V,V)$, so as $(U,U) = (V,V)$, $(U,U) = 2 (W,W)$ and $\bK[G] \subset \bL_{\star \hspace{-0.1em}}[G]$ is in Case 3(a). The type is (?, C, A, ?, 3(a), 3(a)).

$\bullet$ In Case C, $(W,W) = 2 (V,V)$, so as $(U,U) = (V,V)$, $2(U,U) = (W,W)$ and $\bK[G] \subset \bL_{\star \hspace{-0.1em}}[G]$ is in Case 3(b). The type is (?, C, C, ?, 3(a), 3(b)).

As $(V,V) = (U,U)$, we see also that the three cases for $\bK[G] \subset \bL[G]$ determine those for $\bL[H] \subset \bL[G]$: Case 2 cannot occur for $\bL[H] \subset \bL[G]$, as $V \cong {}^{\omega}V$, and the three remaining cases are determined by $(V,V):(Z,Z)$, which is determined by the case for $\bK[G] \subset \bL[G]$. The three options are (B, C, ?, 1, 3(a), ?), (C, C, ?, 3(a), 3(a), ?), and (A, C, ?, 3(b), 3(a), ?).

Pairing these two options for the extension $\bL[H] \subset \bL_{\star \hspace{-0.1em}}[G]$ with the three options for $\bK[G] \subset \bL[G]$, we obtain six options:
\begin{itemize}
     \item \textbf{(B, C, A, 1, 3(a), 3(a))},
     \item \textbf{(B, C, C, 1, 3(a), 3(b))},
     \item \textbf{(C, C, A, 3(a), 3(a), 3(a))},
     \item \textbf{(C, C, C, 3(a), 3(a), 3(b))},
     \item \textbf{(A, C, A, 3(b), 3(a), 3(a))},
     \item \textbf{(A, C, C, 3(b), 3(a), 3(b))}.
\end{itemize}

\vspace{0.5em}
\noindent \textbf{Case C.3(b)} (?, C, ?, ?, 3(b), ?)
In this case we have that
\[
U \cong U \otimes \pi, \ X \cong {}^{\omega} X, \ \Ind(X) = U, \ \Res(U) = X \oplus X, \ 2(X,X) = (U,U).
\]
As the ratio $(U,U) : (X,X) : (V,V)$ along the central axis of the diamond is is $4:2:1$, and no pair has a ratio greater than $2$, we must have that $(Z,Z) = 2(V,V)$, $(U,U) = 2(Z,Z)$, $(W,W) = 2 (V,V)$, and $(U,U) = 2(W,W)$. Therefore the type is \textbf{(C, C, C, 3(b), 3(b), 3(a))}. \qedsymbol

\section{Behaviour of Each Type}

In Section \ref{sect:realise24ways} below, we consider which of the 24 types of Table~\ref{tab:24types} actually occur. Before that, in this section we record the behaviour of each type in Tables \ref{tab:24types2}, \ref{tab:24types4}, \ref{tab:24types5}, \ref{tab:24types5.5} and \ref{tab:24types6}. We will use this to identify the types in the examples of Section \ref{sect:realise24ways} below.

\subsection{Justification of Tables \ref{tab:24types2}, \ref{tab:24types4}, \ref{tab:24types5}, \ref{tab:24types5.5} and \ref{tab:24types6}}

Table \ref{tab:24types2}, Table \ref{tab:24types4}, and Table \ref{tab:24types6} each follow directly from Theorem \ref{thm:threeoptions}, Proposition \ref{prop:algindressummary} and Table \ref{tab:thewaysLHSsquare}. 

These also give the dimension ratios appearing in Table \ref{tab:24types5.5}. When all entries of a row of Table \ref{tab:24types5.5} are divisible by two (resp. four), this extra divisibility property of the dimensions follows from the fact that the dimensions of the endomorphism rings of Table \ref{tab:24types6} force the dimension of the representation of minimal dimension (which is always $V$) of a row to be a multiple of two (resp. four). The rows with minimal dimension equal to two are (2), (6), (7), (13), (16), (17), (18), (19), (22), (23), (24) and the rows with minimal dimension equal to four are (20) and (21).

Finally, Table \ref{tab:24types5} also follows from Theorem \ref{thm:threeoptions}, Proposition \ref{prop:algindressummary} and Table \ref{tab:thewaysLHSsquare}, the only part requiring further justification being the computation of $\Stab_{K_4}(Z)$ and $\Stab_{K_4}(V)$.

\subsubsection{Computation of $\Stab_{K_4}(Z)$}

To compute $\Stab_{K_4}(Z)$, note that if at least one of $a * Z \cong Z$ or $b * Z \cong Z$ is true, then $\Stab_{K_4}(Z)$ is determined by $a * Z$ and $b * Z$, which can be deduced from the type. The remaining cases to consider are (4), (5), (10), (18), (19), where $a * Z \not\cong Z$ or $b * Z \not\cong Z$, and so $\Stab_{K_4}(Z)$ could be either $1$ or $\langle c \rangle$. This can be deduced from whether or not $V \cong \overline{V}$, as follows.

If $V \cong \overline{V}$, then
\[
\Ind_{\bL[H]}^{\bL[G]}(V) \cong \Ind_{\bL[H]}^{\bL[G]}(\overline{V}),
\]
hence
\[
Z \oplus (Z \otimes \pi) \cong \overline{Z} \oplus (\overline{Z} \otimes \pi),
\]
so, as $Z \not\cong \overline{Z}$, $Z \cong \overline{Z} \otimes \pi$ and the stabiliser is $\langle c \rangle$.

Conversely, if $\overline{Z} \cong Z \otimes \pi$, then
\[
\Res_{\bL[H]}^{\bL[G]}(\overline{Z}) \cong \Res_{\bL[H]}^{\bL[G]}(Z \otimes \pi) = \Res(Z),
\]
so $V \cong \overline{V}$. 

Therefore in cases (4), (5), (18) and (19) the stabiliser is $\langle c \rangle$, and $\Stab_{K_4}(Z) = 1$ in case (10).

\subsubsection{Computation of $\Stab_{K_4}(V)$}

Whether or not $\overline{V} = a * V \cong V$, ${}^{\omega} V = b * V \cong V$, or $\epsilon * V = c * V \cong V$ can be read off from the type of $\bK[G] \subset \bL[G]$, $\bL[H] \subset \bL[G]$ and $\bL[H] \subset \bL_{\star \hspace{-0.1em}}[G]$ respectively.

\subsection{When $\bL$ Splits $H$ and $G$}\label{sect:LsplitsHG}

Certain types of the 24 of Theorem \ref{thm:24ways} appearing in Table \ref{tab:24types} can be ruled out depending on whether certain modules are absolutely irreducible. Indeed, from Table \ref{tab:24types6} we can see that
\begin{itemize}
    \item If $\End_{\bK[H]}(X) = \bK$, then there cannot be 3(a) in the 5th column,
    \item If $\End_{\bK[G]}(U) = \bK$, then there cannot be 3(b) in the 5th column, or 3(a) in the 6th column,
    \item If $\End_{\bL[H]}(V) = \bL$, then there cannot be 3(a) in the 4th column.
    \item If $\End_{\bL[G]}(Z) = \bL$, then there cannot be 3(b) in the 4th column.
\end{itemize}

For example, the number of possible types is at most
\begin{itemize}
    \item 12, when $\bK$ splits $H$,
    \item 10, when $\bK$ splits $G$,
    \item 5, when $\bK$ splits $H$ and $G$,
    \item 18, when $\bL$ splits $H$,
    \item 18, when $\bL$ splits $G$,
    \item 12, when $\bL$ splits $H$ and $G$.
\end{itemize}

Of particular interest is the final case when $\bL$ is a splitting field for $H$ and $G$. As described above, in this case 12 of the 24 types of Theorem \ref{thm:24ways} are ruled out, and we choose to give a alternative simpler name the remaining types.

\begin{defn}\label{def:twelvetypes}
For $V \in \Irr_{\bbL}(H)$, with $\End_{\bbL[H]}(V) = \bL$ and $\End_{\bbL[G]}(Z) = \bL$, the D-type of $V$, $\fD(V)$, is the roman numeral associated to the type as below. \vspace{0.5em}

\begin{minipage}[t]{0.48\textwidth}  
\begin{tabular}{c l c}
I & \textbf{(A, A, A, 1, 1, 2)} & (1) \\
II & \textbf{(B, A, C, 1, 3(b), 3(b))} & (5) \\
II.V & \textbf{(B, A, A, 1, 3(b), 3(a))} & (4) \\
III & \textbf{(A, A, B, 2, 2, 3(b))} & (2) \\
IV & \textbf{(B, B, B, 1, 1, 2)} & (10) \\
V & \textbf{(A, B, A, 2, 3(a), 1)} & (12) \\
\end{tabular}
\end{minipage}%
\hfill
\begin{minipage}[t]{0.48\textwidth}  
\begin{tabular}{c l c}
VI & \textbf{(C, B, C, 2, 3(b), 1)} & (15) \\
VII & \textbf{(B, B, B, 2, 2, 1)} & (11) \\
VIII & \textbf{(C, C, C, 1, 1, 2)} & (16) \\
IX & \textbf{(B, C, A, 1, 3(a), 3(a))} & (18)\\
IX.V & \textbf{(B, C, C, 1, 3(a), 3(b))} & (19)\\
X & \textbf{(C, C, B, 2, 2, 3(a))} & (17) \\
\end{tabular}
\end{minipage}
\end{defn}

We will see in Section \ref{sect:recoverDyson10} of Chapter \ref{chapter:tenfoldway} how these are consistent with the labelling of Dyson's tenfold way, with the two new cases that occur here labelled by II.V and IX.V. The motivation behind the naming of the new cases II.V and IX.V is that these D-types have the same L-type (see Section \ref{sect:forgetsemilineardata} below) as II and IX respectively. More precisely, under the map $\Psi$ of Section \ref{sect:forgetsemilineardata} below, $\Psi(\text{II}) = \Psi(\text{II.V})$ and $\Psi(\text{IX}) = \Psi(\text{IX.V})$. Note however that $\Psi$ \emph{is} injective when restricted to the other eight types.

There are other nice facts that happen when $\bL$ splits $H$ and $G$, such as a more explicit description of the endomorphism rings involved in cases A, B and C.

\begin{prop}\label{prop:endringwhensplit}
    Suppose that $V, W$ are as in Theorem \ref{thm:threeoptions}, and $\End_{\bL[H]}(V) = \bL$. Then the S-type of $V$ is characterised by the following conditions:
    \begin{enumerate}
        \item[\emph{(A)}] $\End_{\bL_{\star \hspace{-0.1em}}[G]}(W) = \bK$,
        \item[\emph{(B)}] $\End_{\bL_{\star \hspace{-0.1em}}[G]}(W) = \bL$,
        \item[\emph{(C)}] $\End_{\bL_{\star \hspace{-0.1em}}[G]}(W) = D$, where $D$ is a central division algebra over $\bK$ of degree $2$.
    \end{enumerate}
\end{prop}

\begin{proof}
    This follows directly from \cite[Prop.\ 8.28]{TAYSchur}.
\end{proof}

\begin{cor}
    When $\bL$ splits $H$ and $G$, the endomorphism rings for each type are described in Table \ref{tab:12types}.
\end{cor}

\begin{proof}
We know the first two columns must be $\bL$ from the fact that $\bL$ splits $H$ and $G$. We therefore know the dimensions of the remaining entries of each row by Table \ref{tab:24types6}. If the dimension is $2$, then the endomorphism ring is $\bL$ by Proposition \ref{prop:endringwhensplit}. Otherwise, from Proposition \ref{prop:endringwhensplit} the endomorphism is a central division algebra $D$ of dimension $4$ over $\bK$. A priori, these endomorphism rings are different, however in all but D-type IX.V the division algebra is the same.

Indeed, for any pair of adjacent columns among the final three, at least one of induction or restriction between the corresponding algebras preserves irreducibility, by Table \ref{tab:indres}, and therefore by functoriality there is a map from one endomorphism ring to the other. This has to be injective, because each endomorphism ring is a central simple algebra, and therefore an isomorphism by comparing the dimension.
\end{proof}

\begin{remark}
    The $\bK$-algebras $D_1$ and $D_2$ in D-type IX.V of Table \ref{tab:12types} are non-isomorphic in general. See Remark \ref{rem:differentdivalg} of Section \ref{sect:realiseIXV}.
\end{remark}

\section{Forgetting the Semilinear Data}\label{sect:forgetsemilineardata}

Given $V \in \Irr_{\bbL}(H)$, we can assign to $V$ the type of $V$, which we may consider as a map
\[
    \fT \colon  \Irr_{\bbL}(H) \rightarrow T,
\]
which by Theorem \ref{thm:24ways} lands in the subset $T_0 \subset T$ of the 24 possible types. Setting
\[
S \coloneqq \{\text{A,\:B,\:C}\}^2 \times \{\text{1,\:2,\:3(a),\:3(b)}\}^2,
\]
we also have a forgetful map
\[
\Psi \colon T \rightarrow S
\]
defined by forgetting the 3rd and 6th entry of an element of $T$. By construction this map remembers exactly the linear representation theoretic information about the type of $V$, which forms the left-hand side of the diamond (\ref{eqn:labelleddiamond}).

One of our original motivations for considering the type was to understand if the semilinear Schur index can always be determined the linear part of the type: the part obtainable without reference to semilinear representations. The following shows that this is not possible in general.

\begin{cor}
The map $\Psi|_{\fT(\Irr_{\bbL}(H))} \colon \fT(\Irr_{\bbL}(H)) \rightarrow S$ is not injective in general.
\end{cor}

In other words: it is impossible in general to recover the type of an irreducible $V \in \Irr_{\bbL}(H)$ from linear representation-theoretic data alone.

\begin{proof}
    Let $\bL / \bK$ be the extension $\bQ(i) / \bQ$, and let $G = \text{SmallGroup}(32,9)$, which has a unique subgroup $H$ isomorphic to $D_8 \times C_2$. Considering $G$ with the presentation
    \[
G := \langle a,b,c \mid 
a^4 = b^2 = c^4 = 1,\;
bab = c a c^{-1} = a^{-1},\;
c b c^{-1} = ab
\rangle
\]
the subgroup $H$ is given by $H = \langle a, b, c^2 \rangle$. The types of the two $2$-dimensional irreducible representations of $H$ over $\bQ(i)$ are (8) and (9), which have the same linear type. To see this, one can note that the field $\bQ(i)$ is a splitting field for $H$, as $H$ has exponent $4$, and therefore one can use Theorem \ref{thm:15waysforcyclotomic} to determine the type by computing Frobenius-Schur indicators. Alternatively, one can use the algorithm to compute the linear type found at \cite{GITHUB}.
\end{proof}

\subsection{The Linear Type}\label{sect:lineartype}

We include a repository \cite{GITHUB} containing magma code for computing the \emph{linear type} of $V \in \Irr_{\bbL}(H)$, or in other words, an algorithm for the function
\[
   \fL \coloneqq \colon \Psi \circ \fT \colon \Irr_{\bbL}(H) \rightarrow T_0 \rightarrow S.
\]
We note that the fibers of the map
\[
    \Psi|_{T_0} \colon T_0 \rightarrow S
\]
which have size at least two are:
\[
    \{(4), (5)\}, \quad \{(6), (7)\}, \quad \{(8), (9)\}, \quad \{(18), (19)\}, \quad \{(20), (21)\}, \quad \{(22), (23)\}.
\]
In particular, we have the following. 
\begin{cor}\label{cor:canrecoverSLtype}
Let $V \in \Irr_{\bbL}(H)$. If the linear type $\fL(V)$ satisfies
\[
\fL(V) \in \{(1), (2), (3), (10), (11), (12), (13), (14), (15), (16), (17), (24)\},
\]
then the type $\fT(V)$ can be uniquely recovered from $\fL(V)$.
\end{cor}

\section{Realisation of The 24 Types}\label{sect:realise24ways}

In this section we consider realisability of the 24 types of Theorem \ref{thm:24ways} and Table \ref{tab:24types}. We give explicit realisations of all but one type, leaving the realisation of type (21) as an open question. We give all our realisations over number fields, and in the twelve cases of Section \ref{sect:LsplitsHG} we further take $V \in \Irr_{\bbL}(H)$ to be such that $\End_{\bL[H]}(V) = \bL$ and $\End_{\bL[G]}(Z) = \bL$.

In light of Section \ref{sect:lineartype}, to determine the type of a representation, we typically first compute the linear type, which can either be done by hand or using the algorithm of Section \ref{sect:lineartype}. This either determines the type, or reduces the type to two possibilities.

These two possibilities are distinguished by the S-type of $V$, which to determine one has to typically either give an explicit extension of $V$ to a semilinear representation of $G$, or show by hand that no such extension is possible.

\subsection{Realisation of Type (1) / I}

Let $H \leq G$ be any Real structure, and let $\bL / \bK$ be any degree two extension with $\charfield(\bK) \neq 2$ (as assumed throughout this section). Then the trivial representation of $H$ over $\bL$, has linear type \textbf{(A, A, 1, 1)}, and therefore type \textbf{(A, A, A, 1, 1, 2)}.

\subsection{Realisation of Type (2) / III}

Let $\bL / \bK$ be the extension $\bQ(i) / \bQ$, and $G = D_8 = \langle a,b \mid a^4 = b^2 = e, bab = a^{-1} \rangle$. Let $H$ be the subgroup $\{e, b, a^2, a^2 b\} \cong K_4$, and let $V$ be the $1$-dimensional representation of $H$ over $\bQ(i)$ with character $\chi$, where $\chi(b) = -1$ and $\chi(a^2) = -1$. Then the linear type of $V$ is \textbf{(A, A, 2, 2)}, and so the type of $V$ is \textbf{(A, A, B, 2, 2, 3(b))}.

\subsection{Realisation of Type (3)}

Let $\bL / \bK$ be the extension $\bQ(i) / \bQ$, and let $H = C_3$ inside $G = D_6$. Let $V$ be the non-trivial irreducible representation of $H$ over $\bQ(i)$, which is $2$-dimensional. The linear type of $V$ is \textbf{(A, A, 3(a), 3(a))}, and so the type of $V$ is \textbf{(A, A, A, 3(a), 3(a), 3(b))}.

\subsection{Realisation of Type (4) / II.V}

Let $\bL = \bQ(\zeta)$, where $\zeta$ is a primitive $8$th root of $1$. This is a degree $4$ extension of $\bQ$, and we take $\bK$ to be the intermediate field $\bQ(i)$. As $\bL = \bK(\zeta)$, with $\zeta^2 = i$, the non-trivial automorphism $\alpha$ of $\bL/\bK$ sends $\zeta \mapsto - \zeta$. Let $H = C_4 = \langle x \rangle$ and $G = C_8 = \langle y \rangle$, with $y^2 = x$. Let $W$ be the $1$-dimensional $\bK[H]$-module where $x$ acts by multiplication on $\bK$ by $i$, and let $V$ be the base change to $\bL$.

We claim that $V$ is extendable to a $\bL \rtimes G$-module. Define an $\alpha$-linear action of $y$ by
\[
    \rho_y(v) = (\zeta i) \cdot \alpha(v),
\]
for $v \in V$. Then we can compute that
\begin{align*}
\rho_y(\rho_y(v)) &= \rho_y((\zeta i)\alpha(v)), \\
&= (\zeta i) \alpha(\zeta i) \cdot \alpha^2(v), \\
&= \zeta (-\zeta) i^2 \cdot v, \\
&= \zeta^2 \cdot v, \\
&= i \cdot v,
\end{align*}
and therefore, as this is exactly the action of $x = y^2$, $V$ extends via this action of $y$ to an $\bL \rtimes G$-module.

To see that this example is of type \textbf{(B, A, A, 1, 3(b), 3(a))}, from the fact that $V$ is induced from a $\bK[H]$-module and extends to a $\bL \rtimes G$-module, it must be of type \textbf{(?, A, A, ?, ?, ?)}. The only other type with these specified entries is \textbf{(A, A, A, 1, 1, 2)}. In our case, we also have that $W$ is not the restriction of any irreducible $\bK[G]$-module, and therefore the extension $\bK[H] \subset \bK[G]$ must be in type 2 or 3(b), and so this example must have type \textbf{(B, A, A, 1, 3(b), 3(a))}.

\subsection{Realisation of Type (5) / II}

Let $\bL / \bK$ be the extension $\bQ(i) / \bQ$, $G = C_4$, and $H = C_2$. Then the linear type of the non-trivial $1$-dimensional $L$-representation $V$ is \textbf{(B, A, 1, 3(b))}. Therefore, to show that the type of $V$ is \textbf{(B, A, C, 1, 3(b), 3(b))}, it is sufficient to show that $V$ does not extend to a semilinear representation of $G$ over $L$.

Suppose there was such an extension, call it $\rho$. Writing $C_4 = \langle y \rangle$, $C_2 = \langle x \rangle$ with $y^2 = x$, let $a \in \bQ(i)$ be such that $y$ acts by matrix $a$. Then, from the relation $y^2 = x$, we have that
\[
a \overline{a} = -1
\]
which is impossible, as $a \overline{a} \in \bR_{\geq 0}$.

Alternatively, if such an extension of $V$ to a semilinear representation existed, it would also exist for $\bC / \bR$, which is impossible as the Frobenius-Schur indicator of $V$ is $-1$ (cf.\ Chapter \ref{chapter:tenfoldway}).

\subsection{Realisation of Type (6)}

Let $\bL / \bK$ be the extension $\bQ(\zeta_{21}) / \bQ(\zeta_7)$, let $G$ be the quaternion group $Q_8 = \langle a,x \mid a^4, x^2 = a^2, x^{-1}ax = a^{-1} \rangle$, and let $H = \langle a \rangle$. Let $V$ be the unique $2$-dimensional irreducible representation of $H$ over $\bQ(\zeta_{21})$, which has linear type \textbf{(C, A, 3(a), 3(b))}. To see that the type of $V$ is \textbf{(C, A, A, 3(a), 3(b), 3(a))} it is sufficient to show the following. 

\begin{lemma}
    There is an extension of $V$ to a semilinear representation of $G$ over $\bQ(\zeta_{21})$.
\end{lemma}

\begin{proof}
If $\mu \in \bQ(\zeta_7)$ and $a \in \bQ(\zeta_{21})$ with $N_{\bQ(\zeta_{21}) / \bQ(\zeta_{7})}(a)(1+\mu^2) = -1$, then the matrices
    \[
A = \begin{pmatrix}
    0 & -1 \\ 1 & 0
\end{pmatrix}, \qquad X = a \begin{pmatrix}
    1 & \mu \\ \mu & -1
\end{pmatrix}
\]
define such an extension of $V$, as $A^4 = I_2$, $X \cdot {}^{\sigma \!}X = - I_2 = A^2$, and $A \cdot X = X \cdot {}^{\sigma \!}A^{-1}$, where $\sigma$ denotes the non-trivial element of $\Gal(\bQ(\zeta_{21}) / \bQ(\zeta_{7}))$. Such a pair of elements $(\mu, a)$ can indeed be found: for example, taking $\mu = 2 \sqrt{-7} = 2(\zeta_7 + \zeta_7^2 - \zeta_7^3 + \zeta_7^4 - \zeta_7^5 - \zeta_7^6)$, and $a = \frac{1}{9}(\zeta_{21}^7 + 2)$.
\end{proof}

\subsection{Realisation of Type (7)}

Let $\bL / \bK$ be the extension $\bQ(\omega) / \bQ$, where $\omega$ is a primitive third root of $1$, and let $G$ be the quaternion group $Q_8 = \langle a,x \mid a^4, x^2 = a^2, x^{-1}ax = a^{-1} \rangle$. Taking $H = \langle a \rangle$, let $V$ be the unique $2$-dimensional irreducible representation of $H$ over $\bQ(\omega)$, which has linear type \textbf{(C, A, 3(a), 3(b))}. The following shows that the type of $V$ is \textbf{(C, A, C, 3(a), 3(b), 3(b))}. 

\begin{lemma}
    There is no extension of $V$ to a semilinear representation of $G$ over $\bQ(\omega)$.
\end{lemma}

\begin{proof}
A matrix form of $V$ over $\bQ(\omega)$ is given by 
\[
\rho(a) = \begin{pmatrix}
    0 & -1 \\ 1 & 0
\end{pmatrix},
\]
and an extension of $\rho$ to a semilinear representation of $G$ over $\bQ(\omega)$ is equivalent to the choice of an $X \in \GL_2(\bQ(\omega))$ with $A X = X \overline{A}^{-1}$ and $X \cdot \overline{X} = A^2 = -I_2$. The first condition implies that
\[
X = \begin{pmatrix}
    a & b \\ b & -a
\end{pmatrix}
\]
for some $a , b \in \bQ(\omega)$. The second condition is equivalent to $\overline{X} = - X^{-1}$, which mean that $a = \lambda \overline{a}$, $b = \lambda \overline{b}$, where $\lambda = \det(X) = - (a^2 + b^2)$. If $a \neq 0$, then $b/a = \lambda \overline{b} / \lambda \overline{a} = \overline{b/a}$, and similarly if $b \neq 0$ then $a/b \in \bQ$, so in either case there are $r,s \in \bR$ and $z \in S^1$ with $a = rz$ and $b = sz$. We have that $\lambda = -(a^2 + b^2) = - z^2 (r^2 + s^2)$, and from the equation $a = \lambda \overline{a}$ and the fact that $\overline{z} = z^{-1}$,
\[
r z = \lambda \overline{a} = - z^2 (r^2 + s^2) r \overline{z} = - r z (r^2 + s^2),
\]
so $r = 0$, as $r^2 + s^2 > 0$ because $(r,s) \in \bR^2 \setminus \mathbf{0}$. Similarly, this also forces $s = 0$, a contradiction.
\end{proof}

\subsection{Realisation of Type (8)}

Let $\bL / \bK$ be the extension $\bQ(i) / \bQ$, and $G = D_{16} = \langle a,b \mid a^8 = b^2 = e, bab = a^{-1} \rangle$. Let $H$ be the subgroup $\langle a^2, b \rangle \cong D_4$, and let $V$ be the unique $2$-dimensional representation of $H$ over $\bQ(i)$. Then the linear type of $V$ is \textbf{(A, A, 3(b), 3(b))}. To prove that the type of $V$ is \textbf{(A, A, A, 3(b), 3(b), 3(a))} it is sufficient to show that $V$ extends to a semilinear representation of $G$ over $\bQ(i)$. To this end, note that $V$ has a matrix representation $\rho$ with
\[
\rho(a^2) = \begin{pmatrix}
    i & 0 \\ 0 & -i
\end{pmatrix}, \qquad \rho(b) = \begin{pmatrix}
    0 & 1 \\ 1 & 0
\end{pmatrix}
\]
and so we may define
\[
A \coloneqq \begin{pmatrix}
    0 & -1 \\ i & 0
\end{pmatrix}
\]
which is chosen to satisfy the relations $A \overline{A} = \rho(a^2)$ and $A \overline{X} \overline{A} = X$. These are exactly the relations to ensure that $\rho(a) \coloneqq A$ defines an extension of $V$ to a semilinear representation of $G$ over $\bQ(i)$.

\subsection{Realisation of Type (9)} Let $G = C_{12}$, $H = C_6$, and $\bL / \bK$ be the extension $\bQ(\omega) / \bQ$ where $\omega$ is a primitive third root of $1$. Let $V$ be the unique order two character of $H$ over $\bQ(\omega)$. The linear type of $V$ is \textbf{(A, A, 3(b), 3(b))}, and to prove that the type of $V$ is \textbf{(A, A, C, 3(b), 3(b), 3(b))} it is sufficient to show that $V$ does not extend to a semilinear representation of $G$ over $\bQ(\omega)$.

Similarly to the realisation of type (5), such an extension is equivalent to the existence of an $a \in \bQ(\omega)$ with $a \overline{a} = -1$, which is impossible as the norm of an element of $\bQ(\omega)$ lies in $\bQ_{\geq 0}$, $\bQ(\omega) / \bQ$ being imaginary quadratic.

\subsection{Realisation of Type (10) / IV}

Let $H = C_3$, $G = C_6$, and $\bL/\bK$ be the extension $\bQ(\omega) / \bQ$, where $\omega$ is a primitive third root of $1$. Let $V$ be the $1$-dimensional representation of $H$ over $\bQ(\omega)$ where the generator acts by multiplication by $\omega$. Then the linear type of $V$ is \textbf{(B, B, 1, 1)}, and so the type of $V$ is \textbf{(B, B, B, 1, 1, 2)}.

\subsection{Realisation of Type (11) / VII}

Let $\bL / \bK$ be the extension $\bQ(i) / \bQ$, and let $G$ denote the modular group of order 16, identified as SmallGroup(16,6), which has presentation given by
\[
    G = \langle a, x \mid a^8 = x^2 = e, xax = a^5 \rangle.
\]  
This has a unique subgroup $H$ isomorphic to $C_4 \times C_2$, explicitly defined as $H = \langle a^2, x \rangle$. Let $V$ be the $1$-dimensional representation of $H$ over $\bQ(i)$ with character $\chi$ such that $\chi(a^2) = i$ and $\chi(x) = 1$. Then it is direct to check that the three characters $\chi, {}^a \chi, \overline{\chi}$ are distinct, and therefore $\Stab_{K_4}(V) = 1$, so from Table \ref{tab:24types5} the type of $V$ must be \textbf{(B, B, B, 2, 2, 1)}.

\subsection{Realisation of Type (12) / V}

Let $H = C_3$, $G = D_6$, and $\bL/\bK$ be the extension $\bQ(\omega) / \bQ$, where $\omega$ is a primitive third root of $1$. Let $V$ be the $1$-dimensional representation of $H$ over $\bQ(\omega)$ where the generator acts by multiplication by $\omega$. Then the linear type of $V$ is \textbf{(A, B, 2, 3(a))}, and so the type of $V$ is \textbf{(A, B, A, 2, 3(a), 1)}.

\subsection{Realisation of Type (13)}

Let $\bL / \bK$ be the extension $\bQ(i) / \bQ$, and let $G$ by the modular group of order 16 with presentation
\[
G = \langle a, x \mid a^8, x^2, xax^{-1} = a^5 \rangle.
\]
We take $H$ to be the cyclic subgroup $H = \langle a \rangle$, and $V$ to be either of the irreducible $2$-dimensional representations of $H$ over $\bQ(i)$. The linear type of $V$ is \textbf{(B, B, 3(a), 3(a))}, and so the type of $V$ is \textbf{(B, B, B, 3(a), 3(a), 3(b))}.

\subsection{Realisation of Type (14)}

Let $\bL / \bK$ be the extension $\bQ(i) / \bQ$, $G = C_8$ and $H = C_4$. Take $V$ to be either of the faithful characters of $H$ over $\bQ(i)$. The linear type of $V$ is \textbf{(B, B, 3(b), 3(b))}, and so the type of $V$ is \textbf{(B, B, B, 3(b), 3(b), 3(a))}.

\subsection{Realisation of Type (15) / VI}

Let $\bL / \bK$ be the extension $\bQ(i) / \bQ$, $G = Q_8$, and $H = \langle j \rangle \leq Q_8$. Take $V$ to be $1$-dimensional representation of $H$ with $\chi(j) = i$. The linear type of $V$ is \textbf{(C, B, 2, 3(b))}, and so the type of $V$ is \textbf{(C, B, C, 2, 3(b), 1)}.

\subsection{Realisation of Type (16) / VIII}

Let $\bL / \bK$ be the extension $\bQ(i) / \bQ$, $H = Q_8 \leq G = Q_8 \times C_2$. Take $V$ to be the unique $2$-dimensional representation of $H$. The linear type of $V$ is \textbf{(C, C, 1, 1)}, and so the type of $V$ is \textbf{(C, C, C, 1, 1, 2)}.

\subsection{Realisation of Type (17) / X}

Let $\bL / \bK$ be $\bQ(i) / \bQ$, and let $G$ be SmallGroup(32,8) in the database of small groups. The group $G$ has presentation
\[
G \;=\;
\left\langle 
  a, b, c \;\middle|\;
  a^{8} = b^{4} = c^{2} = 1,\;
  a^{4} = b^{2},\;
  bc = cb,\;
  a b a^{-1} = b^{3} c,\;
  a c a^{-1} = b^{2} c
\right\rangle
\]
and the field $\bQ(i)$ is a splitting field for $G$. We take $H$ to be the subgroup
\[
H = \langle a^2 , b, c \rangle
\]
which is isomorphic to $Q_8 \times C_2$, where, for $i$ and $j$ the standard generators of $Q_8$, $a^2 \leftrightarrow i$ and $b \leftrightarrow j$.

Let $V$ be the inflation of the unique $2$-dimensional representation of $Q_8$ along the projection $H \twoheadrightarrow Q_8$, tensored with the sign character of the projection $H \twoheadrightarrow C_2$. In particular, we have that $\chi(c) = -2$, $\chi(b^2) = -2$, and $\chi(b^2 c) = 2$. This character is $\bQ$-valued, but not induced from a representation over $\bQ$, and therefore the classical Schur index $m_c(V) = 2$. We can also check that ${}^{\omega} V \not\cong \overline{V}$, so $\Stab_{K_4}(V) = \langle a \rangle$, as $\chi(c) = -2$, but 
\[
{}^{a}\chi(c) = \chi(aca^{-1}) = \chi(b^2 c) = 2.
\] 
Consulting Table \ref{tab:24types5}, these two properties force the type of $V$ to be \textbf{(C, C, B, 2, 2, 3(a))}.

\subsection{Realisation of Type (18) / IX}

Let $\bL / \bK$ be the extension $\bQ(i) / \bQ$, and let $G$ be SmallGroup(16,13), also know as the Pauli group, with presentation
\[
G = \langle a, x, y \mid a^4, y^4, x^2, a^2 = y^2, xax^{-1} = a^{-1}, [a,y], [x,y] \rangle.
\]
There is a unique copy of $Q_8$ as the subgroup $H \coloneqq \langle a, xy \rangle$, where $i \mapsto a, j \mapsto xy$. Let $(V, \rho)$ be the unique $2$-dimensional irreducible representation of $H$, with matrix form
\[
\rho(a) = \begin{pmatrix}
    i & 0 \\ 0 & -i
\end{pmatrix}, \qquad \rho(xy) = \begin{pmatrix}
    0 & -1 \\ 1 & 0
\end{pmatrix}.
\]
This has linear type \textbf{(B, C, 1, 3(a))}. From Table \ref{tab:24types}, to show this has type \textbf{(B, C, A, 1, 3(a), 3(a))} it is sufficient to show that $\rho$ extends to a semilinear representation of $G$ over $\bQ(i)$. Set
\[
A \coloneqq \rho(a), \qquad X \coloneqq \begin{pmatrix}
    i & 0 \\ 0 & i
\end{pmatrix}, \qquad Y \coloneqq \begin{pmatrix}
    0 & -i \\ i & 0
\end{pmatrix}
\]
Then letting the generators $a, x, y$ of $G$ act by the matrices $A,X,Y$ respectively, one can check that this defines a semilinear representation which extends $\rho$, by explicitly verifying the relations
\[
\rho(a) = A, \ \rho(xy) = X \overline{Y} , \ Y \overline{Y} = A^2, \ X \overline{X} = I_2, \ X \overline{A} = A^{-1} X, \ A Y = \overline{Y} A, \ X \overline{Y} = Y \overline{X}.
\]

\subsection{Realisation of Type (19) / IX.V}\label{sect:realiseIXV} Let $\bL / \bK$ be the extension $\bQ(i) / \bQ$, and let $G$ be $\text{SmallGroup}(24,5)$, with presentation
\[
G = \langle w, x, y, z \;|\;
w^2, \;
x^2 = y, \;
y^2, \;
z^3, \;
[w,x], \;
[w,y], \;
[x,y], \;
[y,z], \;
[x,z], \;
w z w^{-1} = z^{-1}
\rangle.
\]
Let $H = \langle yz, wx \rangle$, which is isomorphic to
\[
\text{Dic}_{12} \coloneqq \langle a,b \mid a^6, b^2 = a^3,  b^{-1} a b = a^{-1} \rangle
\]
where $a$ corresponds to $yz$ and $b$ corresponds to $wx$. Let $V$ be the $2$-dimensional representation of $H$ over $\bQ(i)$ defined by 
\[
\rho(yz) \coloneqq \begin{pmatrix}
    0 & -1 \\ 1 & 1
\end{pmatrix}
\qquad
\rho(wx) \coloneqq \begin{pmatrix}
    0 & i \\ i & 0
\end{pmatrix}.
\]

We could compute the linear type of $V$ using the algorithm (cf.\ Section \ref{sect:lineartype}), but we instead give a more involved example than the previous types to show how this can be computed by hand.

\begin{lemma}
    The linear type of $V$ is \emph{\textbf{(B, C, 1, 3(a))}}.
\end{lemma}

\begin{proof}
The complex character table of $H$ is
\[
\begin{array}{c|c|rrrrrr}
 & m(\psi_i) & (1) & (a^3) & (a^2) & (b) & (ab) & (a) \\
\hline
\text{Size}  &   & 1 & 1 & 2 & 3 & 3 & 2 \\
\text{Order} &   & 1 & 2 & 3 & 4 & 4 & 6 \\
\hline
\psi_{1} & 1 & 1 & 1 & 1 & 1 & 1 & 1 \\
\psi_{2} & 1 & 1 & 1 & 1 & -1 & -1 & 1 \\
\psi_{3} & 1 & 1 & -1 & 1 & i & -i & -1 \\
\psi_{4} & 1 & 1 & -1 & 1 & -i & i & -1 \\
\psi_{5} & 2 & 2 & -2 & -1 & 0 & 0 & 1 \\
\psi_{6} & 1 & 2 & 2 & -1 & 0 & 0 & -1 \\
\end{array}
\]
where $m(\psi)$ is the Schur index of $\psi$ over $\bQ$. This is the same as the character table of $H$ over $\bQ(i)$: for $\psi_i$ with $i \neq 5$ this follows from the fact that $\psi_i$ is $\bQ(i)$-valued, and for $\psi_5$ this is because $\psi_5$ is the character of the representation $\rho$, and so admits a model over $\bQ(i)$. We know that the type of $\psi_5$ has C-type (the type of the extension $\bQ[H] \leq \bQ(i)[H]$) $\mathbf{C}$, and therefore that the type of $V$ is one of the nine types in the range (16) to (24) of Table \ref{tab:24types}. Furthermore, because $\psi_5$ is the unique representation with Schur index $2$, we see that conjugation by the odd element $w$ fixes $\psi_5$, and therefore by Table \ref{tab:indres} (which applies by Section \ref{sect:indextwosubgroup}) $V$ is the restriction of exactly two irreducible characters $\chi$, $\chi \otimes \pi$ of $G$ over $\bQ(i)$ with $\chi \neq \chi \otimes \pi$, where $\pi$ is the sign character of $G$ over $H$. Setting $Z \in \Irr_{\bQ(i)}(G)$ to be the representation corresponding to $\chi$, this immediately tells us that $b \not\in \Stab_{K_4}(Z)$, and therefore that the type of $V$ is one of (16), (18) or (19).

To see that the type of $V$ cannot be (16), consider the complex character table of $G$
\[
\begin{array}{c|rrrrrrrrrrrr}
\text{Class} & 1 & 2 & 3 & 4 & 5 & 6 & 7 & 8 & 9 & 10 & 11 & 12 \\
\hline
\text{Size}  & 1 & 1 & 3 & 3 & 2 & 1 & 1 & 3 & 3 & 2 & 2 & 2 \\
\text{Order} & 1 & 2 & 2 & 2 & 3 & 4 & 4 & 4 & 4 & 6 & 12 & 12 \\
\hline
\chi_1  & 1 & 1 & 1 & 1 & 1 & 1 & 1 & 1 & 1 & 1 & 1 & 1 \\
\chi_2  & 1 & 1 & -1 & -1 & 1 & 1 & 1 & -1 & -1 & 1 & 1 & 1 \\
\chi_3  & 1 & 1 & 1 & 1 & 1 & -1 & -1 & -1 & -1 & 1 & -1 & -1 \\
\chi_4  & 1 & 1 & -1 & -1 & 1 & -1 & -1 & 1 & 1 & 1 & -1 & -1 \\
\chi_5  & 1 & -1 & 1 & -1 & 1 & i & -i & i & -i & -1 & -i & i \\
\chi_6  & 1 & -1 & -1 & 1 & 1 & i & -i & -i & i & -1 & -i & i \\
\chi_7  & 1 & -1 & 1 & -1 & 1 & -i & i & -i & i & -1 & i & -i \\
\chi_8  & 1 & -1 & -1 & 1 & 1 & -i & i & i & -i & -1 & i & -i \\
\chi_9  & 2 & 2 & 0 & 0 & -1 & 2 & 2 & 0 & 0 & -1 & -1 & -1 \\
\chi_{10} & 2 & 2 & 0 & 0 & -1 & -2 & -2 & 0 & 0 & -1 & 1 & 1 \\
\chi_{11} & 2 & -2 & 0 & 0 & -1 & -2i & 2i & 0 & 0 & 1 & -i & i \\
\chi_{12} & 2 & -2 & 0 & 0 & -1 & 2i & -2i & 0 & 0 & 1 & i & -i \\
\end{array}
\]
The element $y$ is central in $G$, and order two, and therefore the 2nd column is the conjugacy class of $y$. We therefore see that the type characters $\chi$ and $\chi \otimes \pi$ which restrict to $\psi_5$ are $\chi_{11}$ and $\chi_{12}$, and so $\Stab_{K_4}(Z) = \langle c \rangle$ because $\overline{\chi_{11}} = \chi_{12}$. This tells us that the type of $V$ is one of (18) or (19), both of which have linear type \textbf{(B, C, 1, 3(a))}.
\end{proof}

Finally, to show that the type of $V$ is indeed \textbf{(B, C, C, 1, 3(a), 3(a))}, it remains to see that $V$ does not extend to a semilinear representation of $G$ over $\bQ(i)$.

\begin{lemma}\label{lem:noextension}
    There is no extension of $V$ to a semilinear representation of $G$ over $\bQ(i)$.
\end{lemma}
\begin{proof}
Suppose there exists such an extension. We have that $(yz)^3 = y$ and $(yz)^4 = z$, so $\rho(y) = -I_2$ and $\rho(z) = -\rho(yz)$. Furthermore, both generators $w$ and $x$ are odd (in $G \setminus H$), and therefore such an extension is equivalent to the data of matrices $W,X$ over $\bQ(i)$ which satisfy the appropriate relations coming from the relations in $G$ and $\rho(wx) = W \overline{X}$. In particular, the relations $w^2 = 1$, $x^2 = y$ and $xz = zx$ mean that $W$ and $X$ satisfy the relations
\begin{align*}
W \overline{W} &= I_2, \\ 
X \overline{X} &= \rho(y), \\
X \overline{\rho(z)} &= \rho(z) X.
\end{align*}
Writing
\[
X = \begin{pmatrix} a & b \\ c & d \end{pmatrix},
\]
then as $\rho(z)$ is $\bQ$-valued
\[
X \rho(z) = \rho(z) X,
\]
and therefore
\[
\begin{pmatrix} a & b \\ c & d \end{pmatrix} \begin{pmatrix} 0 & 1 \\ -1 & -1 \end{pmatrix} 
=
\begin{pmatrix} 0 & 1 \\ -1 & -1 \end{pmatrix}\begin{pmatrix} a & b \\ c & d \end{pmatrix}
\]
or equivalently
\[
\begin{pmatrix} -b & a-b \\ -d & c-d \end{pmatrix} = \begin{pmatrix} c & d \\ -a-c & -b-d \end{pmatrix}.
\]
In particular, $c = -b$ and $d = a-b$, so
\[
X = \begin{pmatrix} a & b \\ -b & a-b \end{pmatrix}.
\]
Because $b = wx$ and $w$ commutes with $x$, $x$ commutes with $b$ and therefore
\[
X \overline{\rho(b)} = \rho(b) X
\]
so
\[
\begin{pmatrix} a & b \\ -b & a-b \end{pmatrix} \begin{pmatrix} 0 & -i \\ -i & 0 \end{pmatrix}
=
\begin{pmatrix} 0 & i \\ i & 0 \end{pmatrix} \begin{pmatrix} a & b \\ -b & a-b \end{pmatrix}
\]
or equivalently
\[
\begin{pmatrix} -b & a-b \\ a & b \end{pmatrix} = \begin{pmatrix} -b & -a \\ b-a & b \end{pmatrix}.
\]
In particular $b -a =a$ so $b = 2a$ and
\[
X = \begin{pmatrix} a & 2a \\ -2a & -a \end{pmatrix}.
\]
Finally, we have that
\[
X \overline{X} = \rho(y) = -I_2,
\]
so $X^{-1} = - \overline{X}$, hence
\[
- \det(X) \begin{pmatrix} \overline{a} & \overline{2a} \\ -\overline{2a} & -\overline{a} \end{pmatrix} = \begin{pmatrix} -a & -2a \\ 2a & a \end{pmatrix}.
\]
Concentrating on any fixed entry,
\[
a = \overline{a} \det(X) = 3 a^2 \overline{a},
\]
hence
\[
\frac{1}{3} = a \overline{a}.
\]
However, this is impossible: the $3$-adic valuation of the norm of an element of $\bQ(i)$ is always even, because $3$ is inert in $\bQ(i)$.
\end{proof}

\begin{remark}
If $\bL$ was bigger then this representation over $\bK$ would extend. For example, if $\bL / \bK$ was $\bC / \bR$ we could take $a = 1 / \sqrt{3}$ and defining $X$ as above would give the required extension.
\end{remark}

\begin{remark}\label{rem:differentdivalg}
    Let us identify the division algebras $D_1 \coloneqq \End_{\bQ[H]}(X)$ and $D_2 \coloneqq \End_{\bQ(i) \rtimes G}(W)$ in $\Br(\bQ)$. Using the \texttt{SchurIndices} function in magma, one can check that the non-split places of $D_1$ are $3$ and $\infty$. For $D_2$, even though we do not know this algebra explicitly (we haven't even constructed explicitly the semilinear representation $W$), we can still identify $D_2$ by considering the splitting behaviour at places of $\bQ$.
    
    For a place $v$ of $\bQ$, $D_{2,v}$ is split whenever $\bQ(i)_v \coloneqq \bQ_v \otimes_{\bQ} \bQ(i)$ is not a field \cite[Cor.\ 11.4]{TAYSchur}. Suppose instead that $\bQ(i)_v$ is a field, which happens exactly when $v = \infty$ or $v = p$ is a finite prime $p \geq 2$ with $p = 2$ or $p \equiv 3 \!\! \mod 4$. Note that in this case $\bQ(i)_v = \bQ_v(i)$, and
    \[
    V_v \coloneqq \bQ(i)_v \otimes_{\bQ(i)} V \in \Rep_{\bQ(i)_v}(H)
    \]
    is irreducible, as $V$ is absolutely irreducible over $\bQ(i)$, and therefore $V_v$ has a Schur index $m_{\bQ_v}^{\bQ_v(i)}(V_v) \in \{1,2\}$. This Schur index is related to the splitting of $D_{2,v}$ by \cite[Cor.\ 11.4]{TAYSchur}, which tells us that $D_{2,v}$ is split if and only if $m_{\bQ_v}^{\bQ_v(i)}(V_v) = 1$. From the proof of Lemma \ref{lem:noextension}, one sees that $m_{\bQ_v}^{\bQ_v(i)}(V_v) = 2$ (equiv.\ $V_v$ does not extend to a semilinear representation of $G$ over $\bQ(i)_v$) if and only if $1/3$ cannot be written as the sum of two squares in $\bQ_v$. Over $\bR$, $1/3$ is the sum of squares, so $D_{2,v}$ is split over $\bR$ (and so in particular $D_1 \not\cong D_2$). 

    For the other places, if $p \geq 2$ is a finite prime then $1/3$ is the sum of two squares in $\bQ_p$ if and only if $3p^{2e}$ is the sum of the squares of two elements of $\bZ_p$ for some $e \geq 1$, which can be seen by clearing denominators. For primes $p \equiv 3 \!\! \mod 4$, a non-zero $p$-adic integer $x \in \bZ_p$ is a sum of two squares of $\bZ_p$ if and only if $\nu_p(x)$ is even, and a non-zero $2$-adic integer $x = 2^r u \in \bZ_2$ with $u \in \bZ_2^{\times}$ is a sum of two squares in $\bZ_2$ if and only if $u \equiv 1 \!\! \mod 4$ \cite{KCON}. In particular, taking $x \coloneqq 3p^{2e}$, we have that $3p^{2e}$ is not a sum of two elements of $\bZ_p$ if and only if $p = 2$ or $p = 3$.

    In summary, we can conclude that places over $\bQ$ where $D_2$ splits are the finite places $2$ and $3$. 
\end{remark}

\subsection{Realisation of Type (20)}

Let $\bL / \bK$ be the extension $\bQ(i) / \bQ$ and $G = \text{SmallGroup}(32,44)$, with presentation
\[
G = \langle x, y, z \mid x^8,\; y^2,\; z^2,\; yxy^{-1} = x^3,\; zxz^{-1} = x^5,\; zyz^{-1} = x^4 y\rangle.
\]
We take $H \coloneqq \langle x, yz \rangle$, which is isomorphic to $Q_{16}$, and let $V$ be the unique $4$-dimensional irreducible representation of $H$ over $\bQ(i)$, explicitly given by
\[
\rho(x) = \begin{pmatrix}
    0 & i & 0 & 0 \\
1 & 0 & 0 & 0 \\
0 & 0 & 0 & 1 \\
0 & 0 & -i & 0
\end{pmatrix}, \qquad \rho(yz) = \begin{pmatrix}
0 & 0 & -1 & 0 \\
0 & 0 & 0 & -1 \\
1 & 0 & 0 & 0 \\
0 & 1 & 0 & 0
\end{pmatrix}.
\]
The linear type of $V$ is \textbf{(C, C, 3(a), 3(a))}, and therefore to show that the type of $V$ is \textbf{(C, C, A, 3(a), 3(a), 3(a))} it is sufficient to show that $V$ extends to a semilinear representation of $G$ over $\bQ(i)$. Set
\[
Y = \begin{pmatrix}
-1 & 0 & 0 & 0 \\
0 & -i & 0 & 0 \\
0 & 0 & 1 & 0 \\
0 & 0 & 0 & i
\end{pmatrix}, \qquad Z = \begin{pmatrix}
0 & 0 & 1 & 0 \\
0 & 0 & 0 & i \\
1 & 0 & 0 & 0 \\
0 & i & 0 & 0
\end{pmatrix}.
\]
Then letting the generators $y$ and $z$ of $G$ act by the matrices $Y$ and $Z$ respectively, one can check that this defines a semilinear representation which extends $\rho$, by explicitly verifying the relations
\[
\rho(yz) = Y \overline{Z} , \quad Y \overline{Y} = I_2, \quad Z \overline{Z} = I_2, \quad Y \overline{\rho(x)} = \rho(x)^3 Y, \quad Z \overline{\rho(x)} = \rho(x)^5 Z, \quad Z \overline{Y} = \rho(x)^4 Y \overline{Z}. 
\]

\subsection{Analysis of Type (21)}

We do not know if type (21) can be realised. This has the same linear type as type (20), and we have checked by hand that all representations $V \in \Irr_{\bbL}(H)$ for $\bL / \bK$ either $\bQ(i) / \bQ$ or $\bQ(\omega) / \bQ$ and $|H| \leq 50$ with the correct linear type have type (20), not type (21). Therefore we leave the following as an open question.

\begin{prob}\label{prob:case21}
Is type (21) realisable, for some finite extension $\bL / \bK$ and finite group $G$?
\end{prob}

It is unclear if this is true even without assuming that $G$ is finite. Let us give a brief analysis of the behaviour of this type, to explain the constraints on the linear representations that it imposes, in the case when $V$ has minimal possible dimension $4$ over $\bL$.

In this case, $D_U \coloneqq \End(U)$ is a $4$-dimensional central division $\bK$-algebra. Let $D_X \coloneqq \End(X)$, and $D_W \coloneqq \End(W)$, which both have dimension $8$ over $\bK$. We have that $U|_H = X$, and $\bL_{\star \hspace{-0.1em}}[G] \otimes_{\bL[H]} U = W$, and so there are maps $D_U \rightarrow D_X$ and $D_U \rightarrow D_W$. We therefore have isomorphisms
\[
Z_X \otimes_{\bbK} D_U \xrightarrow{\sim} D_X, \qquad Z_W \otimes_{\bbK} D_U \xrightarrow{\sim} D_W,
\]
because these maps are injective (from central simple algebras) and the image has dimension $8$.

We have that $\End(V)$ is a degree four field extension of $\bK$. We also have that $X_L = V \oplus V$, and $W|_H = V \oplus V$, hence
\[
\bL \otimes_{\bbK} D_X \xrightarrow{\sim} M_2(\End(V)), \qquad \bL \otimes_{\bbK} D_W \xrightarrow{\sim} M_2(\End(V)),
\]
and taking the centre
\[
\bL \otimes_{\bbK} Z_X \xrightarrow{\sim} \End(V), \qquad \bL \otimes_{\bbK} Z_W \xrightarrow{\sim} \End(V).
\]

Consider the extension $\bL_X \coloneqq \bL \otimes_{\bbK} Z_X / Z_X$. Then
\[
\End(X_{Z_X}) \cong Z_X \otimes_{\bbK} D_X \cong Z_X \otimes_{\bbK} Z_X \otimes_{\bbK} D_U \cong (Z_X \times Z_X) \otimes_{\bbK} D_U = D_X \times D_X.
\]
Therefore $X_{Z_X} = X_1 \oplus X_2$, with $X_1 \not\cong X_2$ and $\End(X_1) \cong D_X \cong \End(X_2)$ as $Z_X$-algebras. Furthermore, as $X_1 \not\cong X_2$, then $V_1 \coloneqq X_{1, \bbL_X} \not\cong V_2 \coloneqq X_{2, \bbL_X}$. In particular,
\[
V_{\bbL_X} \cong V_1 \oplus V_2
\]
with $V_1 \not\cong V_2$. Now $Z$ is absolutely irreducible, and $V_1$ and $V_2$ both correspond under induction-restriction to $Z_{\bbL_X}$, and $Z \otimes \pi \cong Z$, hence $V_i \not\cong V_i^{\omega}$, and thus $V_1 \cong V_2^{\omega}$, as $V^{\omega} = V$. Furthermore $V_i \cong \overline{V_i}$, applying the non-trivial automorphism of $\bL_X / Z_X$ as each is defined over $Z_X$, and therefore $V_1 \not\cong \overline{V_1}^{\omega} \cong V_2$. Both $V_1$ and $V_2$ are absolutely irreducible, and therefore there is an irreducible semilinear representation $W_0 / \bL_X$ with $W_0|_H = V_1 \oplus V_2$ and $\End(W_0) = \bL_X$.

Now $W|_H = V \oplus V$, hence
\[
W_{\bbL_X}|_H = (W|_H)_{\bbL_X} = V_1 \oplus V_2 \oplus V_1 \oplus V_2.
\]
In particular, $W_{\bbL_X}|_H = (W_0 \oplus W_0)|_H$, hence $W_{\bbL_X} = W_0 \oplus W_0$, and
\[
Z_X \otimes_{\bbK} \End(W) \cong \End(W_{\bbL_X}) \cong M_2(\End(W_0)) \cong M_2(\bL_X).
\]
For example, taking the centre, we have that 
\[
Z_X \otimes_{\bbK} Z_W \xrightarrow{\sim} \bL_X.
\]
In particular, we see that $Z_X$ and $Z_W$ are linearly disjoint. They are two of the three index two subfields of $\bL_X = \End(V)$, the other being $\bL$. For example, this means that $Z_X \neq Z_W$.

\subsection{Realisation of Type (22)}

Let $\bL / \bK$ be the extension $\bQ(i) / \bQ$, and let $G = \text{SD}_{16}$, the semidihedral group of order 16, with presentation
\[
G = \langle a,x \mid a^8, x^2, xax^{-1} = a^3\rangle.
\]
We take $H \coloneqq \langle a^2, ax \rangle$, which is isomorphic to $Q_8$, and let $V$ be the unique $2$-dimensional irreducible representation of $H$ over $\bQ(i)$. The linear type of $V$ is \textbf{(A, C, 3(b), 3(a))}, and therefore to show that the type of $V$ is \textbf{(A, C, A, 3(b), 3(a), 3(a))} it is sufficient to show that $V$ extends to a semilinear representation of $G$ over $\bQ(i)$. The matrix form of $V$ is 
\[
\rho(a^2) = \begin{pmatrix}
    i & 0 \\ 0 & -i
\end{pmatrix}, \qquad \rho(ax) = \begin{pmatrix}
    0 & -1 \\ 1 & 0
\end{pmatrix}.
\]
Set
\[
A \coloneqq \begin{pmatrix}
    0 & -1 \\ i & 0
\end{pmatrix}, \qquad X \coloneqq \begin{pmatrix}
    i & 0 \\ 0 & 1
\end{pmatrix}.
\]
Then letting the generators $a$ and $x$ of $G$ act by the matrices $A$ and $X$ respectively, one can check that this defines a semilinear representation which extends $\rho$, by explicitly verifying the relations
\[
\rho(a^2) = A \overline{A} , \quad \rho(ax) = A \overline{X}, \quad X \overline{X} = I_2, \quad X \overline{A} X = A \overline{A} A.
\]

\subsection{Realisation of Type (23)}

Let $\bL / \bK$ be the extension $\bQ(\omega) / \bQ$, where $\omega$ is a primitive third root of $1$, and let $G = \text{SmallGroup}(16,13)$, with presentation
\[
G = \langle x,y,z \mid x^4, z^2, y^2 = z^2, [x,y], [x,z], zyz^{-1} = y^{-1}\rangle.
\]
This has a unique subgroup isomorphic to $Q_8$, given by $H = \langle y, xz\rangle$, where $y$ and $xz$ correspond to the generators $i$, $j$ of $Q_8$. Let $V$ be the $2$-dimensional representation irreducible representation of $Q_8$ over $\bQ(\omega)$. The linear type of $V$ is \textbf{(A, C, 3(b), 3(a))}, and therefore to show that the type of $V$ is \textbf{(A, C, C, 3(b), 3(a), 3(b))} it is sufficient to show the following.

\begin{lemma}\label{lem:noextension23}
    There is no extension of $V$ to a semilinear representation of $G$ over $\bQ(\omega)$.
\end{lemma}

\begin{proof}
This is similar to the proof of Lemma \ref{lem:noextension}. A matrix form of $V$ over $\bQ(\omega)$ is given by 
\[
\rho(y) = \begin{pmatrix}
    \omega & \omega^2 \\ \omega^2 & -\omega
\end{pmatrix}, \qquad
\rho(xy) = \begin{pmatrix}
    0 & -1 \\ 1 & 0
\end{pmatrix}
\]
and an extension of $\rho$ to a semilinear representation of $G$ over $\bQ(\omega)$ is equivalent to the data of a $X,Z \in \GL_2(\bQ(\omega))$ such that the matrix relations corresponding to the relations of $G$ are satisfied and $X\overline{Z} = \rho(xy)$. From the fact that $zy = y^{-1}z$, we have that $Z \overline{\rho(y)} = \rho(y^{-1}) Z$, from which one obtains that $Z$ has the form 
\[
Z = \begin{pmatrix}
    \omega b + \omega^2 c & b \\ c & -\omega^2 b - \omega c
\end{pmatrix}
\]
for some $b,c \in \bQ(\omega)$. From $z^2 = 1$, $\overline{Z} Z = I_2$, and from this one can compute that $b = -\det(Z)\overline{b}$ and $c = - \det(Z) \overline{c}$. From $x^2 = b^2$, $X \overline{X} = \rho(b^2) = - I_2$, and $X = \rho(xy)Z$ from $X\overline{Z} = \rho(xy)$, so $\rho(xy) Z \rho(xy) \overline{Z} = -I_2$. From this one can again compute that $b = \det(Z) \overline{c}$ and $c = \det(Z) \overline{b}$. Combining, we have that $c = -b$, and computing the determinant that $\det(Z) = -2 b^2$. Putting these back into $b = -\det(Z) \overline{b}$, $b = 2b^2\overline{b}$, so $1/2 = b \overline{b}$. However this is impossible: the $2$-adic valuation of a norm of an element of $\bQ(\omega)$ is always even, because the prime $2$ is inert in $\bQ(\omega)$.
\end{proof}

\subsection{Realisation of Type (24)}

Let $\bL / \bK$ be the extension $\bQ(i) / \bQ$, and let $G = Q_{16}$, the generalised quaternion group of order 16, with presentation
\[
G = \langle a,b,c \mid a^4 = b^2 = c^2 = abc \rangle.
\]
We take $H \coloneqq \langle a^2, b \rangle$, which is isomorphic to $Q_8$, and let $V$ be the unique $2$-dimensional irreducible representation of $H$ over $\bQ(i)$. The linear type of $V$ is \textbf{(C, C, 3(b), 3(b))}, and so the type of $V$ is \textbf{(C, C, C, 3(b), 3(b), 3(a))}.

\section{Interpretation of the 24 Types in Terms of Blocks}\label{sect:blockinterpretation}

We can interpret the 24 types of Theorem \ref{thm:24ways} in terms of extensions of matrix algebras coming from blocks of Diagram \eqref{eqn:diamondrings_complete}. From Theorem \ref{SBlock_3}, and $\bK$ has characteristic different from $2$, we know that the inclusions of Diagram \eqref{eqn:diamondrings} induce maps of Bohr completions, and the each resulting morphism of Diagram \eqref{eqn:diamondrings_complete} is injective.

Suppose that $e \in \widehat{\bK[H]}$ is a central primitive idempotent. There is no guarantee that the image of $e$ is central in any of the larger algebras of Diagram \eqref{eqn:diamondrings_complete}, so we consider
\[
\zeta_{\omega} \colon \widehat{\bK[H]} \rightarrow \widehat{\bK[H]},
\]
the $\bK$-algebra automorphism defined by completing the $\bK$-algebra automorphism
\[
\omega^{-1} - \omega \colon \bK[H] \rightarrow \bK[H]
\]
of $\bK[H]$, and set
\[
f \coloneqq \twopartdef{e}{e = \zeta_{\omega}(e)}{e + \zeta_{\omega}(e)}{e \neq \zeta_{\omega}(e)}.
\] 
Concretely, if $e$ corresponds to the simple module $X \in \Irr_{\bbK}(H)$ under the correspondence of Lemma \ref{lem:idemirrcorresp}, then $\zeta_{\omega}(e)$ corresponds to ${}^{\omega} X \in \Irr_{\bbK}(H)$.

For $B$ any of the $\bK$-algebras of Diagram \ref{eqn:diamondrings}, we write $\widehat{\iota}_B \colon \widehat{\bK[H]} \hookrightarrow \widehat{B}$ for the inclusion map.

\begin{defn}\label{def:ideminbigger}
    For any central primitive idempotent $e \in \widehat{\bK[H]}$, the \emph{semilinear block}, or \emph{S-block}, corresponding to $e$ is the diagram of $\bK$-algebra inclusions
    \[
     \begin{tikzcd}
	& \widehat{\iota}_{\bbK[G]}(f) \cdot \widehat{\bK[G]} \\
	\widehat{\iota}_{\bbL[G]}(f) \cdot \widehat{\bL[G]} & f \cdot \widehat{\bK[H]} & \widehat{\iota}_{\bL_{\star \hspace{-0.1em}}[G]}(f) \cdot \widehat{\bL_{\star \hspace{-0.1em}}[G]} \\
	& \widehat{\iota}_{\bbL[H]}(f)\cdot \widehat{\bL[H]}
	\arrow[hook, from=1-2, to=2-1]
	\arrow[hook, from=1-2, to=2-3]
	\arrow[hook, from=2-2, to=1-2]
	\arrow[hook, from=2-2, to=3-2]
	\arrow[hook, from=3-2, to=2-1]
	\arrow[hook, from=3-2, to=2-3]
\end{tikzcd}
\]
\end{defn}

We now interpret the 24 types of Theorem \ref{thm:24ways} in terms of the S-block. In particular, this shows that the image of $e$ in each larger algebra is a central idempotent, and that each algebra of the S-block is a matrix algebra.

\begin{prop}\label{prop:descofSblock}
    Suppose that $V \in \Irr_{\bbL}(H)$, and that $U$, $W$, $X$ and $Z$ correspond to $V$ as in Section \ref{sect:propagation}. Let $e \in \widehat{\bK[H]}$ be the central primitive idempotent corresponding to $X$, and define $f \in \widehat{\bK[H]}$ as above. Then:
    \begin{itemize}
        \item $f = \sum_{T \in C_2 \cdot X} e_T$, and for $n_X \coloneqq \dim_{\bbK}(X) / \dim_{\bbK}(\End_{\bbK[H]}(X))$,
        \[
            f \cdot \widehat{\bK[H]} \cong \prod_{T \in C_2 \cdot X} M_{n_X}(\End_{\bbK[H]}(X)).
        \]
        \item $\widehat{\iota}_{\bbK[G]}(f) = \sum_{T \in C_2 \cdot U} e_T$, and for $n_U \coloneqq \dim_{\bbK}(U) / \dim_{\bbK}(\End_{\bbK[G]}(U))$,
        \[
            \widehat{\iota}_{\bbK[G]}(f) \cdot \widehat{\bK[G]} \cong \prod_{T \in C_2 \cdot U} M_{n_U}(\End_{\bbK[G]}(U)).
        \]
        \item $\widehat{\iota}_{\bbL[H]}(f) = \sum_{T \in K_4 \cdot V} e_T$, and for $n_V \coloneqq \dim_{\bbL}(V) / \dim_{\bbL}(\End_{\bbL[H]}(V))$,
        \[
            \widehat{\iota}_{\bbL[H]}(f) \cdot \widehat{\bL[H]} \cong \prod_{T \in K_4 \cdot V} M_{n_V}(\End_{\bbL[H]}(V)).
        \]
        \item $\widehat{\iota}_{\bbL_{\star \hspace{-0.1em}}[G]}(f) = \sum_{T \in C_2 \cdot W} e_T$, and for $n_W \coloneqq \dim_{\bbK}(W) / \dim_{\bbK}(\End_{\bbL_{\star \hspace{-0.1em}}[G]}(W))$,
        \[
            \widehat{\iota}_{\bbL_{\star \hspace{-0.1em}}[G]}(f) \cdot \widehat{\bL_{\star \hspace{-0.1em}}[G]} \cong \prod_{T \in C_2 \cdot W} M_{n_W}(\End_{\bbL_{\star \hspace{-0.1em}}[G]}(W)).
        \]
        \item $\widehat{\iota}_{\bbL[G]}(f) = \sum_{T \in K_4 \cdot Z} e_T$, and for $n_Z \coloneqq \dim_{\bbL}(Z) / \dim_{\bbL}(\End_{\bbL[H]}(V))$,
        \[
            \widehat{\iota}_{\bbL[G]}(f) \cdot \widehat{\bL[G]} \cong \prod_{T \in K_4 \cdot Z} M_{n_Z}(\End_{\bbL[G]}(Z)).
        \]
    \end{itemize}
\end{prop}

\begin{proof}
We give the argument for $\bK[G]$, the remaining points being deduced in exactly the same manner. Note that $\widehat{\iota}_{\bbK[G]}(f)$ acts as the identity on simple objects of $\Rep_{\bbK}(G)$ isomorphic to $U$ or $U \otimes \pi$, because these restrict in $\Rep_{\bbK}(H)$ to a direct sum of $X$ and ${}^{\omega}X$, and as zero on all other simple representations of $\bK$ over $G$, as these contain no copies of $X$ or ${}^{\omega}X$ in their restriction to $\Rep_{\bbK}(H)$. This exactly means that $\widehat{\iota}_{\bbK[G]}(f)$ is the sum of $e_{U}$ and $e_{U \otimes \pi}$, the central primitive idempotents corresponding to $X$ and ${}^{\omega} X$ respectively under the correspondence of Lemma \ref{lem:idemirrcorresp}. The second statement follows immediately, from the definition of the idempotents and the description of the Bohr completion as the product of matrix rings over division rings of Lemma \ref{Bohr2}, noting that the endomorphism rings of both $U$ and $U \otimes \pi$ are the same. 
\end{proof}

\begin{eg}
Suppose that $V$ is the trivial representation of $H$. Then endomorphism rings of the corresponding irreducible representations are listed in Table \ref{tab:12types}, the orbits in Table \ref{tab:24types5}, and the dimensions in Table \ref{tab:24types5.5}, and thus by Proposition \ref{prop:descofSblock} the S-block is
\[
     \begin{tikzcd}
	& \bK \times \bK &  \\
	\bL \times \bL & \bK & M_2(\bK) \\
	& \bL &
	\arrow[hook, from=1-2, to=2-1]
	\arrow[hook, from=1-2, to=2-3]
	\arrow[hook, from=2-2, to=1-2]
	\arrow[hook, from=2-2, to=3-2]
	\arrow[hook, from=3-2, to=2-1]
	\arrow[hook, from=3-2, to=2-3]
\end{tikzcd}
\]
\end{eg}

\begin{eg}
For a more interesting example, suppose that $V$ is absolutely irreducible of type VI $ = (15)$, of dimension $n$. As above, consulting Tables \ref{tab:24types5}, \ref{tab:24types5.5} and \ref{tab:12types}, the S-block for $V$ is
\[
     \qquad \qquad \begin{tikzcd}
	& M_n(D) \\
	M_{2n}(\bL) & M_n(\bL) & M_n(D) \times M_n(D) \\
	& M_n(\bL) \times M_n(\bL)
	\arrow[hook, from=1-2, to=2-1]
	\arrow[hook, from=1-2, to=2-3]
	\arrow[hook, from=2-2, to=1-2]
	\arrow[hook, from=2-2, to=3-2]
	\arrow[hook, from=3-2, to=2-1]
	\arrow[hook, from=3-2, to=2-3]
\end{tikzcd}
\]
where $D$ is a $4$-dimensional central division algebra over $\bK$.
\end{eg}
 
\section{Associated $C_2$-Graded Algebras}\label{sect:BrauerWall}

In this section we associate to any semilinear representation $W \in \Irr_{\bbL}^{\hspace{0.1em} \star}(G)$ a $C_2$-graded division algebra $D(W)$, and a simple $C_2$-graded algebra $A(W)$, and explore how these put more restrictions on the structure of the S-block, or equivalently the type of $W$.

For basic notions concerning $C_2$-graded algebras, we direct the reader to Appendix \ref{appendix:gradedalgebras}, which contains the main definitions and certain results which we will make use of below.

We first consider $D(W)$.

\begin{defn}
If $W \in \Irr_{\bbL}^{\hspace{0.1em} \star}(G)$, we set
\[
D(W) \coloneqq \End_{\bbK[G]}(W),
\]
with $C_2$-grading
\begin{align*}
D(W)_0 &= \{f \in \End_{\bbK[G]}(W) \mid f(\lambda w) = \lambda f(w) \text{ for any $w \in W$, $\lambda \in \bL$} \}, \\
D(W)_1 &= \{f \in \End_{\bbK[G]}(W) \mid f(\lambda w) = \overline{\lambda} f(w) \text{ for any $w \in W$, $\lambda \in \bL$}\}.
\end{align*}
\end{defn}

To see that this is a well-defined $C_2$-grading, meaning that $D(W) = D(W)_0 \oplus D(W)_1$, take $\alpha \in \bL \setminus \bK$ with $\overline{\alpha} = - \alpha$, which exists from the eigenspace decomposition for $\overline{(-)} \colon \bL \rightarrow \bL$,
\[
\bL = \{\alpha \in \bL \mid \overline{\alpha} = \alpha\} \oplus \{\alpha \in \bL \mid \overline{\alpha} = -\alpha\},
\]
which uses that $\overline{(-)} \colon \bL \rightarrow \bL$ has order two and $\charfield(\bK) \neq 2$. We can define the automorphism
\[
t \coloneqq t_{\alpha} \colon \End_{\bbK[G]}(W) \rightarrow \End_{\bbK[G]}(W), \qquad t_{\alpha}(f)(w) = \frac{1}{\alpha} f(\alpha w),
\]
which has order two because $\alpha^2 \in \bK$, and is independent of the choice of $\alpha$.

Then, just as with $\bL$, because $\charfield(\bK) \neq 2$ we can decompose $D(W)$ into $\pm 1$ eigenspaces for $t_{\alpha}$, which are exactly $D(W)_0$ and $D(W)_1$. Explicitly, given $f \in D(W)$,
\[
f = f_0 + f_1, \qquad f_0 \coloneqq \frac{\alpha f(v) + f(\alpha v)}{2 \alpha} \in D(W)_0, \qquad f_1 \coloneqq \frac{\alpha f(v) - f(\alpha v)}{2 \alpha} \in D(W)_1.
\]

To see this more canonically, note that there is a $\bL_{\star \hspace{-0.1em}}[G]$-linear isomorphism
\[
\bL \otimes_{\bbK} W \xrightarrow{\sim} W \oplus \overline{W}, \qquad \lambda \otimes w \mapsto (\lambda w, \overline{\lambda}w)
\]
as in the proof of Proposition \ref{prop:resindKGLG} (which incidently doesn't require that $\charfield(\bK) \neq 2$), and thus we can compose the isomorphisms
\begin{align*}
\Hom_{\bbL_{\star \hspace{-0.1em}}[G]}(W,W) \oplus \Hom_{\bbL_{\star \hspace{-0.1em}}[G]}(\overline{W},W) \xrightarrow{\sim} \Hom_{\bbL_{\star \hspace{-0.1em}}[G]}(\bL \otimes_{\bbK} W , W) \xrightarrow{\sim} \Hom_{\bbK[G]}(W,W),
\end{align*}
where
\[
f_1 \oplus f_2 \mapsto [w \mapsto f_1(w) + f_2(w)],
\]
under which $\Hom_{\bbL_{\star \hspace{-0.1em}}[G]}(W,W)$ has image $D(W)_0$ and $\Hom_{\bbL_{\star \hspace{-0.1em}}[G]}(\overline{W},W)$ has image $D(W)_1$. This perspective allows us to immediately deduce the following basic properties of $D(W)$.

\begin{prop}\label{prop:basicpropDW}
    For $W \in \Irr_{\bbL}^{\hspace{0.1em} \star}(G)$, $D(W)$ satisfies the following properties:
    \begin{enumerate}
        \item $D(W)$ is a $C_2$-graded division algebra,
        \item $D(W)$ is ungraded \emph{(}meaning $D(W)_1 = 0$\emph{)} if and only if $W \not\cong \overline{W}$,
        \item $D(W)$ is simple as a $\bK$-algebra if and only if $U \cong U \otimes \pi$,
        \item The graded centre of $D(W)$ is contained in $D(W)_0$, so
        \[
            \widehat{Z}(D(W)) = D(W)_0 \cap Z(D(W)) \leq Z(D(W)_0).
        \] 
    \end{enumerate}
\end{prop}

\begin{proof}
    The first two properties follow immediately from the above discussion and Schur's lemma, noting that the inverse of any non-zero element of $x \in \Hom_{\bbL_{\star \hspace{-0.1em}}[G]}(\overline{W},W)$ can be viewed as an element $y \in \Hom_{\bbL_{\star \hspace{-0.1em}}[G]}(\overline{W},W)$ for which $xy = yx = 1$. For the third point note that $W|_{\bbK[G]} \cong U^{\oplus n}$ for $n \in \{1,2\}$ if $U \cong U \otimes \pi$ and $W|_{\bbK[G]} \cong U \oplus (U \otimes \pi)$ if $U \not\cong U \otimes \pi$ by Proposition \ref{prop:algindressummary}. Then as a $\bK$-algebra, $D(W) = \End_{\bbK[G]}(W)$ is $M_n(\End_{\bbK[G]}(U))$ in the first case and $\End_{\bbK[G]}(U) \times \End_{\bbK[G]}(U \otimes \pi)$ in the latter. The final point follows from Lemma \ref{lem:centreofC2graddiv}.
\end{proof}

Now we want to describe $\End_{\bR[H]}(W)$ in terms of $D(W)$. For this, we introduce the following.

\begin{defn}
For $A$ a $C_2$-graded algebra over $\bK$ and $\alpha \in \bL \setminus \bK$, we set
\[
A_{\alpha} \coloneqq A \, \hat{\otimes} \, \bK[\alpha],
\]
where $\bK[\alpha] = \bL$ as a $\bK$-algebra, with $C_2$-grading $\bK[\alpha]_0 = \bK$, $\bK[\alpha]_1 = \bK \cdot \alpha$.
\end{defn}

To establish some basic properties of $A_{\alpha}$, recall that for any $C_2$-graded algebra $A$ one can consider the \emph{main automorphism} of $A$, $t \colon A \rightarrow A$, which acts as the identity on $A_0$ and multiplication by $-1$ on $A_1$ \cite[Def.\ IV.3.7]{LAM}. When $A = D(W)$, this agrees with $t \colon D(W) \rightarrow D(W)$ as defined above.

\begin{lemma}\label{lem:Aalphaproperties}
    Suppose that $A$ is a $C_2$-graded algebra over $\bK$ and $\alpha \in \bL \setminus \bK$. Then:
    \begin{itemize}
        \item For a formal symbol $\underline{\alpha}$, the $\bK$-linear identification 
        \[
            \phi \colon A \oplus A \cdot \underline{\alpha} \xrightarrow{\sim} A_{\alpha}, \qquad \phi(a + b \cdot \underline{\alpha}) = a \otimes 1 + b \otimes \alpha
        \] 
        induces on $A \oplus A \cdot \underline{\alpha}$ multiplication determined by
\[
\underline{\alpha}^2 = \alpha^2, \qquad \underline{\alpha} x = t(x) \underline{\alpha}
\]
for $x \in A$, and so explicitly given by, for $x,y, u ,v \in A$,
        \[
        (x + y \cdot \underline{\alpha}) * (u + v \cdot \underline{\alpha}) = (xu + \alpha^2 y t(v)) + (xv + y t(u)) \cdot \underline{\alpha}.
        \]
        The induced $C_2$-grading is 
        \[
        (A_{\alpha})_0 = A_0 \oplus A_1 \cdot \underline{\alpha}, \qquad \text{and} \qquad (A_{\alpha})_1 = A_1 \oplus A_0 \cdot \underline{\alpha}.
        \]
        \item If $A$ is $C_2$-graded simple, then $A_{\alpha}$ is $C_2$-graded simple.
        \item If $\beta \in \bL \setminus \bK$ with $\langle \alpha \rangle_{\bbK} = \langle \beta \rangle_{\bbK}$, then $A_{\alpha} \cong B_{\alpha}$ as $C_2$-graded algebras.
        \item The inclusion $A \hookrightarrow A_{\alpha}$ induces identifications $\widehat{Z}(A) = \widehat{Z}(A_{\alpha})$ and
        \[
            Z(A_{\alpha}) = \widehat{Z}(A)_0 \oplus S \cdot \underline{\alpha}, \text{ where } S = \{x \in A_0 \mid x t(a) = ax \text{ for all } a \in A\}.
        \]
        When $A = D$ is a $C_2$-graded division algebra, then $\bF \coloneqq\widehat{Z}(D)$ is a field and this simplifies to
        \[
        Z(D_{\alpha}) = \twopartdef{\bF \oplus \bF u \cdot \underline{\alpha}}{\text{there exists $u \in D_0 \setminus 0$ with $u^{-1}xu = t(x)$ for any $x \in D$,}}{\bF}{\text{otherwise.}}
        \]
        The first case occurs if and only if $D$ is simple as an $\bF$-algebra and $Z(D) = \bF$, and the second case occurs if and only if $D_{\alpha}$ is simple as an $\bF$-algebra and has centre $\bF$.
    \end{itemize}
\end{lemma}

\begin{proof}
    The first point follows by direct computation. The second follows from the fact that the grading on $\bK[\alpha]$ only depends on $\langle \alpha \rangle_{\bbK}$. If $A$ is $C_2$-graded simple, then $A_{\alpha}$ is $C_2$-graded simple, and $A \hookrightarrow A_{\alpha}$ induces an identification $\widehat{Z}(A) = \widehat{Z}(A_{\alpha})$, both of which follow from \cite[Thm.\ IV.2.3]{LAM}, as $\bK[\alpha]$ is $C_2$-graded simple with graded centre $\bK$. It remains only to compute $Z(A_{\alpha})$, which, being graded, is determined by $Z(A_{\alpha})_0$ and $Z(A_{\alpha})_1$. We have that $Z(A_{\alpha})_0 = \widehat{Z}(A_{\alpha})_0 = \widehat{Z}(A)_0$.

    For $x \in A$, we have that $x \underline{\alpha} \in Z(A_{\alpha})_1$ if and only if for any $a \in A$, $x \underline{\alpha}$ commutes with $a$ and $a \underline{\alpha}$, or equivalently $x t(a) = ax$ and $x t(a) = a t(x)$ respectively. Taking $a = 1$, $t(x) = x$, and thus
    \[
    \{ x \in A \mid  x \underline{\alpha} \in Z(A_{\alpha})_1 \} = \{x \in A_0 \mid x t(a) = ax \text{ for all } a \in A\}.
    \]

    When $A = D$ is a $C_2$-graded division algebra, $\widehat{Z}(D) = \widehat{Z}(D)_0$ by Lemma \ref{lem:centreofC2graddiv}, and for $Z(A_{\alpha})_1$, if there is such a non-zero element $u$, then right-multiplication by $u$ defines an isomorphism
    \[
    Z(D) \cap D_0 = \widehat{Z}(D) \xrightarrow{\sim} \{x \in D_0 \mid x t(y) = yx\},
    \]
    from which we yield the required description. The final statement is \cite[Thm.\ IV.3.4]{LAM}.
\end{proof}

The following is the main result of this section, which allows us to relate $D(W)$ with $\End_{\bK[H]}(W)$. We consider $\End_{\bK[H]}(W)$ with $C_2$-grading induced by the order two automorphism
\[
t \coloneqq t_{\alpha,g} \colon \End_{\bbK[H]}(W) \rightarrow \End_{\bbK[H]}(W), \qquad t(f) = \rho(g)^{-1} \circ m_{\alpha^{-1}} \circ f \circ m_{\alpha} \circ \rho(g)
\]
for any $g \in G \setminus H$ and $\alpha \in \bL \setminus \bK$ with $\overline{\alpha} = - \alpha$, where $m_{\alpha} \colon W \rightarrow W$ denotes left-multiplication by $\alpha$. This is independent of the choice of $g$ and $\alpha$.

\begin{prop}\label{prop:recoverKHendo}
    Let $W \in \Irr_{\bbL}^{\hspace{0.1em} \star}(G)$ and $\alpha \in \bL \setminus \bK$ with $\overline{\alpha} = - \alpha$. Then the natural map
    \[
    \phi_{\alpha} \colon D(W)_{\alpha} \rightarrow \End_{\bbK[H]}(W), \qquad \phi_{\alpha}(f \otimes \lambda) = f \circ m_{\lambda},
    \]
    where $m_{\lambda}$ denotes left-multiplication by $\lambda \in \bK[\alpha]$, is an isomorphism of $C_2$-graded $\bK$-algebras.
\end{prop}

\begin{proof}
    With respect to the $C_2$-grading on $\End_{\bbK[H]}(W)$, both $D(W)$ and $\bK[\alpha]$ are $C_2$-graded subalgebras of $\End_{\bbK[H]}(W)$, where we view $\bK[\alpha] \hookrightarrow \End_{\bbK[H]}(W)$ via $\lambda \mapsto m_{\lambda}$. By definition, 
    \[
        \widehat{C}_{\End_{\bbK[H]}(W)}(\bK[\alpha]) = \{ f_0 + f_1 \in \End_{\bbK[H]}(W) \mid m_{\alpha^{-1}} \circ f_0 \circ m_{\alpha} = f_0, m_{\alpha^{-1}} \circ f_1 \circ m_{\alpha} = -f_1\},
    \]
    which, as $t(f_0) = f_0$ and $t(f_1) = -f_1$, are exactly the condition that both $f_0$ and $f_1$ commute with $\rho(g)$. As $f_0,f_1 \in \End_{\bbK[H]}(W)$ and thus commute with all $\rho(h)$ for $h \in H$, we thus obtain that 
    \[
    \widehat{C}_{\End_{\bbK[H]}(W)}(\bK[\alpha]) = \End_{\bbK[G]}(W) = D(W),
    \]
    from which the result follows by Proposition \ref{prop:graddoublecent}.
\end{proof}

The key point is that the isomorphism type of $D(W)$ as a $C_2$-graded division algebra determines the centre of the endomorphism algebra of $X$, as described in the following.

\begin{cor}\label{cor:Dalphaisom}
If $W \in \Irr_{\bbL}^{\hspace{0.1em} \star}(G)$, and $D \cong D(W)$ as $C_2$-graded $\bK$-algebras, then for any $\alpha \in \bL \setminus \bK$ with $\overline{\alpha} = - \alpha$,
\[
    D_{\alpha} \cong \End_{\bK[H]}(W)
\]
as $\bK$-algebras. In particular, $D_{\alpha}$ is simple as a $\bK$-algebra if and only if $X \cong {}^{\omega}X$, and
\[
Z(D_{\alpha}) \cong Z(\End_{\bbK[H]}(W)) \cong \twopartdef{Z(\End_{\bbK[H]}(X))}{X \cong {}^{\omega}X}{Z(\End_{\bbK[H]}(X)) \times Z(\End_{\bbK[H]}(X))}{X \not\cong {}^{\omega}X}
\]
as $\bK$-algebras.
\end{cor}

\begin{remark}
We note that, in practice, $Z(D_{\alpha})$ can be described using Lemma \ref{lem:Aalphaproperties}.
\end{remark}

\begin{proof}
    After Proposition \ref{prop:recoverKHendo}, all that remains to show is the final claim regarding $X$. We have that $W|_{k[H]} \cong X^{\oplus n}$ for some $n \geq 1$ if $X \cong {}^{\omega}X$ and $W|_{k[H]} \cong X^{\oplus n} \oplus ({}^{\omega}X)^{\oplus n}$ for some $n \geq 1$ if $X \not\cong {}^{\omega}X$ (which can be seen considering the restriction possibilities for simple modules from $\Rep_{\bbL}^{\hspace{0.1em} \star}(G)$ to $\Rep_{\bbK}(G)$, and $\Rep_{\bbK}(G)$ to $\Rep_{\bbK}(H)$). The claim then follows taking endomorphism rings.
\end{proof}

\begin{eg}
The extra restrictions on the type that Corollary \ref{cor:Dalphaisom} do not seem enough to rule out Case 21 and resolve Problem \ref{prob:case21}. However, for a given extension $\bL / \bK$, Corollary \ref{cor:Dalphaisom} does place significant restrictions on the type: for example, see Section \ref{sect:recoverDyson10} for the case $\bC / \bR$.
\end{eg}

Now we turn to $A(W)$, which we define in terms of blocks. Fix $\alpha \in \bL \setminus \bK$ with $\overline{\alpha} = - \alpha$. Conjugation by $\alpha$ defines an order-two automorphism of $\bK[G]$
\[
\phi \colon \bK[G] \rightarrow \bK[G], \qquad \phi(x) = \alpha x \alpha^{-1},
\]
which is independent of the choice of $\alpha$ and induces an order-two automorphism 
\[
\hat{\phi} \colon \widehat{\bK[G]} \rightarrow \widehat{\bK[G]}, \qquad \hat{\phi}((x_M)_M) = (x_{\phi_* N})_N.
\]
Concretely, $\phi(h) = h$ for each $h \in H$ and $\phi(g) = - g$ for each $g \in G \setminus H$, and thus $\phi_* N = N \otimes \pi$ is nothing but the twist of $N$ by the sign representation $\pi$ of $G$. Because $\phi$ and $\widehat{\phi}$ have order two and $\charfield(\bK) \neq 2$, the $\pm 1$ eigenspaces of $\phi$ and $\hat{\phi}$ define a $C_2$-gradings on $\bK[G]$ and $\widehat{\bK[G]}$ respectively.

It is clear that 
\[
\bK[G]_0 = \bK[H], \qquad \bK[G]_0 = \bK[G \setminus H].
\]
On the other hand, we need to use Lemma \ref{lem:gradedzeropartcompletion} to show the corresponding result for $\widehat{\bK[G]}$:

\begin{cor}\label{cor:degzeropartofcompletion}
    The inclusion
    \[
    \widehat{\iota}_{\bbK[G]} \colon \widehat{\bK[H]} \hookrightarrow \widehat{\bK[G]}
    \]
    induces an identification
    \[
    \widehat{\iota}_{\bbK[G]}(\widehat{\bK[H]}) = \widehat{\bK[G]}_0.
    \]
\end{cor}

For any central primitive idempotent $e \in \widehat{\bK[H]}$, we can consider the idempotent $f$ of $\widehat{\bK[H]}$ from Definition \ref{def:ideminbigger}. By definition, the idempotent $\iota_{\bbK[G]}(f)$ is central in $\widehat{\bK[G]}$, and therefore lies in $\widehat{\bK[G]}_0$, and thus the following is a well-defined $C_2$-grading.

\begin{defn}
    For $W \in \Irr_{\bbL}^{\hspace{0.1em} \star}(G)$, $A(W) \coloneqq \widehat{\iota}_{\bbK[G]}(f) \cdot \widehat{\bK[G]}$, with $C_2$-grading
    \[
        A(W)_0 = \widehat{\iota}_{\bbK[G]}(f) \cdot \widehat{\bK[G]}_0, \qquad A(W)_1 = \widehat{\iota}_{\bbK[G]}(f) \cdot \widehat{\bK[G]}_1.
    \]
\end{defn}

\begin{lemma}\label{lem:basicpropertiesAW}
    For $W \in \Irr_{\bbL}^{\hspace{0.1em} \star}(G)$, $A(W)$ is simple and non-trivially graded as a $C_2$-graded algebra over $\bK$, and
    \[
        A(W)_0 = \widehat{\iota}_{\bbK[G]}(f \cdot \widehat{\bK[H]}).
    \]
    Furthermore, $A(W)$ \emph{(}resp.\ $A(W)_0$\emph{)} is simple as a $\bK$-algebra if and only if $U \cong U \otimes \pi$ \emph{(}resp.\ $X \cong {}^{\omega}X$\emph{)}.
\end{lemma}

\begin{proof}
    The description of $A(W)_0$ follows directly from Corollary \ref{cor:degzeropartofcompletion}. To see that $A(W)$ is non-trivially graded, note that the image of any $g \in G \setminus H$ under the $C_2$-graded $\bK$-algebra homomorphism 
    \[
    \bK[G] \rightarrow \widehat{\bK[G]} \twoheadrightarrow \widehat{\iota}_{\bbK[G]}(f) \cdot \widehat{\bK[G]}
    \]
    is a unit, as $g$ is a unit of $\bK[G]$, and lies in $(\widehat{\iota}_{\bbK[G]}(f) \cdot \widehat{\bK[G]})_1$. 

    For the simplicity of $A(W)$, note that using Proposition \ref{prop:descofSblock} there are $\bK$-algebras isomorphisms
        \[
            A(W)_0 \cong \prod_{T \in C_2 \cdot X} M_{n_X}(\End_{\bbK[H]}(X)), \qquad A(W) \cong \prod_{T \in C_2 \cdot U} M_{n_U}(\End_{\bbK[G]}(U))
        \]
        for some $n_X, n_U \geq 0$. If $C_2 \cdot U = \{U\}$, then $A(W)$ is simple even as a $\bK$-algebra, so suppose instead that $U \not\cong U \otimes \pi$. From the four cases of Proposition \ref{prop:algindressummary} applied to the extension $\bK[H] \rightarrow \bK[G]$, the fact that $U \not\cong U \otimes \pi$ forces $X \cong {}^{\omega}X$, and thus $A(W)_0$ is simple as a $\bK$-algebra. Then $A(W)$ is $C_2$-graded simple by Lemma \ref{lem:invoddeltgradsimple}, using the invertible odd element of $A(W)$ described above.
\end{proof}

We can now describe, for $W \in \Irr_{\bbL}^{\hspace{0.1em} \star}(G)$, how $A(W)$ and $D(W)$ are related. To do this, we realise both $A(W)$ and $D(W)$ as $C_2$-graded subalgebras of $\End_{\bbK}(W)$.

Here, for any $W \in \Vect_{\bbL}$, we view $\End_{\bbK}(W)$ as a $C_2$-graded algebra with grading by linear and antilinear maps, as in the definition of $D(W)$. We write $\bK[\epsilon]$ for $\bK[\epsilon] / (\epsilon^2 - 1)$, with $C_2$-grading $\bK[\epsilon]_0 = \bK$, $\bK[\epsilon]_1 = \bK \cdot \epsilon$, and in the following view $M_{n}(\bK[\epsilon])$ with $C_2$-grading
\[
M_{n}(\bK[\epsilon])_0 = M_{n}(\bK), \qquad M_{n}(\bK[\epsilon])_1 = M_{n}(\bK \cdot \epsilon).
\]

We first describe the class of $\End_{\bbK}(W)$ in $\BrW(\bK)$.

\begin{lemma}
    If $W \in \Vect_{\bbL}$, then for any $\bL$-basis $e_1, ... ,e_n$ of $W$ the natural map
    \[
        \varphi \colon M_{n}(\bK[\epsilon]) \rightarrow \End_{\bbK}(W), \qquad \varphi((a_{ij}) + (b_{ij})\epsilon)\left(\sum_{k} \lambda_k e_k \right) \coloneqq \sum_{j,k} (a_{kj}\lambda_j + b_{kj} \overline{\lambda_j}) e_k
    \]
    is an isomorphism of $C_2$-graded $\bK$-algebras. In particular, $[\bK[\epsilon]] = [\End_{\bbK}(W)]$ in $\BrW(\bK)$.
\end{lemma}

\begin{proof}
    Because $\bL / \bK$ is Galois, the natural map $\bK[\Gal(\bL / \bK)] \rightarrow \End_{\bbK}(\bL)$ is an isomorphism. Then the map $\varphi$ is the composition of the canonical isomorphisms
    \[
    M_n(\bK[\epsilon]) \equiv M_n(\bK[\Gal(\bL / \bK)]) \xrightarrow{\sim} M_n(\End_{\bbK}(\bL)) \xrightarrow{\sim} \End_{\bbK}(\bL^n) \xrightarrow{\sim} \End_{\bbK}(W),
    \]
    and is easily seen to preserve the $C_2$-grading.
\end{proof}

The $C_2$-graded division $\bK$-algebra $D(W)$ is naturally a $C_2$-graded $\bK$-subalgebra of $\End_{\bbK}(W)$. To realise $A(W)$ as such a subalgebra requires a little more work. There is a natural $\bK$-algebra homomorphism
\[
 \widehat{\bK[G]} \rightarrow \End_{\bK}(W), \qquad (x_M)_M \mapsto x_{W|_{\bbK[G]}},
\]
which, as $U$ and $U \otimes \pi$ are the only constituents of $W|_{\bbK[G]}$, induces an injection of $\bK$-algebras
\[
\pi_W \colon A(W) = \widehat{\iota}_{\bbK[G]}(f) \cdot \widehat{\bK[G]} \hookrightarrow \End_{\bbK}(W).
\]

\begin{lemma}
    The inclusion
    \[
        \pi_W \colon A(W) \hookrightarrow \End_{\bbK}(W)
    \]
    preserves the $C_2$-grading, and thus is a homomorphism of $C_2$-graded $\bK$-algebras.
\end{lemma}

\begin{proof}
    There is similarly a natural $\bK$-algebra homomorphism
\[
 \widehat{\bK[H]} \rightarrow \End_{\bK}(W), \qquad (y_N)_N \mapsto y_{W|_{\bbK[H]}},
\]
which lands inside $\End_{\bbK}(W)_0 = \End_{\bbL}(W)$, because $m_i$ is a $\bK[H]$-linear endomorphism of $W|_{\bbK[H]}$, and so commutes with $y_{W|_{\bbK[H]}}$. This similarly induces an injection
\[
\eta_W \colon f \cdot \widehat{\bK[H]} \hookrightarrow \End_{\bbK}(W)_0
\]
because $X$ and ${}^{\omega} X$ are the only constituents of $W|_{\bbK[H]}$, which fits into the commutative diagram
\[\begin{tikzcd}
	{f \cdot\widehat{\bK[H]}} & \\
	{\hat{\iota}_{\bbK[G]}(f) \cdot \widehat{\bK[G]}} & {\End_{\bbK}(W)}
	\arrow["{\hat{\iota}_{\bbK[G]}}"', hook, from=1-1, to=2-1]
	\arrow["{\eta_W}", hook, from=1-1, to=2-2]
	\arrow["{\pi_W}", hook, from=2-1, to=2-2]
\end{tikzcd}\]
In particular, using Corollary \ref{cor:degzeropartofcompletion},
\[
\pi_W(A(W)_0) = \pi_W(\hat{\iota}_{\bbK[G]}(f \cdot\widehat{\bK[H]})) = \eta_W(f \cdot\widehat{\bK[H]}) \subset \End_{\bbK}(W)_0.
\]
For $A(W)_1$, the action map $\bK[G] \rightarrow \End_{\bbK}(W)$ is $C_2$-graded and is part of the commutative diagram
\[\begin{tikzcd}
	{\bK[G]} & \\
	{\widehat{\bK[G]}} & {\End_{\bbK}(W)}
	\arrow[hook, from=1-1, to=2-1]
	\arrow[from=1-1, to=2-2]
	\arrow[from=2-1, to=2-2]
\end{tikzcd}\]
Thus any $g \in G \setminus H$ defines an element of $A(W)^{\times} \cap A(W)_1$ for which $\pi_W(g) \in \End_{\bbK}(W)_1$, so
\[
\pi_W(A(W)_1) = \pi_W(g \cdot A(W)_0) = \pi_W(g) \cdot \pi_W(A(W)_0) \subset \pi_W(g) \cdot \End_{\bbK}(W)_0 = \End_{\bbK}(W)_1
\]
as required.
\end{proof}

We want to relate $D(W)$ to the graded centraliser of the image of $A(W)$ in $\End_{\bbK}(W)$. 

\begin{defn}
For $W \in \Irr_{\bbL}^{\hspace{0.1em} \star}(G)$ we set
\[
\tilde{D}(W) \coloneqq \widehat{C}_{\End_{\bbK}(W)}(\pi_W(A(W))).
\]
\end{defn}

The $C_2$-graded $\bK$-subalgebras $D(W)$ and $\tilde{D}(W)$ are not equal:
\[
D(W)_0 = \tilde{D}(W)_0, \qquad D(W)_1 \cap \tilde{D}(W)_1 = 0.
\]
However, for any $\alpha \in \bL \setminus \bK$ with $\overline{\alpha} = - \alpha$, multiplication by $\alpha$ defines a $\bK$-linear isomorphism $D(W)_1 \xrightarrow{\sim} \tilde{D}(W)_1$. The below shows that if additionally $\alpha^2 = -1$, this extends to a $C_2$-graded $\bK$-algebra isomorphism $\psi_W \colon D(W) \rightarrow \tilde{D}(W)$, and thus allows us to relate $A(W)$ and $D(W)$.

\begin{prop}\label{prop:relateAWDW}
Suppose that $\bL = \bK(i)$ where $i^2 = -1$, and $W \in \Irr_{\bbL}^{\hspace{0.1em} \star}(G)$. Then the map
\[
   \psi_W \colon D(W) \rightarrow \tilde{D}(W), \qquad f_0 + f_1 \mapsto f_0 + if_1
\]
is an isomorphism of $C_2$-graded $\bK$-algebras.

Consequently, $\widehat{Z}(A(W)) = \bK$ if and only if $\widehat{Z}(D(W)) = \bK$, in which case the natural map
\[
\pi_W \otimes \psi_W \colon A(W) \, \hat{\otimes} \, D(W) \rightarrow \End_{\bbK}(W)
\]
is an isomorphism of $C_2$-graded $\bK$-algebras and thus an equality
\[
[A(W)] \cdot [D(W)] = [\bK[\epsilon]]
\]
in $\BrW(\bK)$.
\end{prop}

\begin{proof}
    The image $\pi_W(A(W))$ is the same as the image of the action map $\bK[G] \rightarrow \End_{\bbK}(W)$, both being equal to $\End_{\End_{\bbK[G]}(W)}(W)$. In particular, $\tilde{D}(W)_0 = \End_{\bbK[G]}(W)_0 = D(W)_0$ and 
    \[
        \tilde{D}(W)_1 = \{f \in \End_{\bbK}(W) \mid  f \circ m_i = - m_i \circ f, \rho(g) \circ f = \sigma(g) f \circ \rho(g) \text{ for all } g\in G \}
    \]
    where $\sigma \colon G \twoheadrightarrow \{\pm 1\}$ is the structure map of $G$. It is clear that
    \[
    D(W)_1 \rightarrow \tilde{D}(W)_1, \qquad f_1 \mapsto if_1
    \]
    is a $\bK$-linear isomorphism. The condition that $i^2 = -1$ is used to show that $\psi_W$ is multiplicative:
    \begin{align*}
        \psi_W(x_0 + x_1)\psi_W(y_0 + y_1) &= (x_0 + i x_1)(y_0 + i y_1), \\
        &= x_0y_0 + (ix_1)(i y_1) + x_0 (iy_1) + ix_1y_0, \\
        &= x_0y_0 - i^2 x_1y_1 + i x_0y_1 + ix_1y_0, \\
        &= (x_0y_0 + x_1y_1) + i(x_0y_1 + x_1 y_0), \\
        &= \psi_W((x_0 + x_1)(y_0 + y_1)),
    \end{align*}
    so $\psi_W$ is an isomorphism of $C_2$-graded $\bK$-algebras. Proposition \ref{prop:graddoublecent} shows the remaining claims.
\end{proof}

\section{Examples}

\subsection{The Split Case}

Suppose that $G = H \times C_2$, and $V \in \Irr_{\bL}(H)$. Then $U \not\cong U \otimes \pi$, and consulting Table \ref{tab:24types} and using Theorem \ref{thm:24ways}, we thus have the following.

\begin{cor}
    Suppose that $G = H \times C_2$, and $V \in \Irr_{\bL}(H)$. Then $\fT(V) \in \{(1), (10), (16)\}$.
\end{cor}

These correspond to the semilinear type of $V$ being A, B or C respectively, as described in Theorem \ref{thm:threeoptions}. Note also that the type of $V$ for the extension $\bK[H] \subset \bL[H]$ is the same as the type of $V$ for the extension $\bL[H] \subset \bL_{\star \hspace{-0.1em}}[G]$, which is consistent with the equivalence $\Rep_{\bbL}^{\hspace{0.1em} \star}(G) \cong \Rep_{\bK}(H)$ of \cite[\S 6]{TAYSchur}.

\subsection{Symmetric Groups}

Suppose that $G = S_n$ and $H = A_n$ for some $n \geq 3$ (to ensure that $H$ is index two in $G$). Suppose that $\charfield(\bK) = 0$ or $\charfield(\bK) > n$, and let $V \in \Irr_{\bbL}(A_n)$. The key property that restricts the type of $V$ is that every linear representation of $S_n$ over $\bK$ is absolutely irreducible. This forces $V$ to have type A for the extension $\bK[S_n] \leq \bL[S_n]$ and $\dim_{\bbL} \End_{\bbL[S_n]}(Z) = 1$. Consulting Tables \ref{tab:24types} and \ref{tab:24types6}, and using Theorem \ref{thm:24ways}, we thus immediately have the following.

\begin{cor}\label{cor:symgroups}
Suppose that $n \geq 3$, and $V \in \Irr_{\bbL}(A_n)$. Then $\fT(V) \in \{(1), (2), (3), (12)\}$.

In particular, the semilinear Schur index $m_{\bbK}^{\bbL}(V)$ of $V$ is $1$.
\end{cor}

In particular, as described in Section \ref{sect:semilinearschurindex}, this allows on to completely determine the semilinear representation theory of $S_n$ over $\bL$ from the character table of $A_n$ over $\bL$ together with its action of $C_2$ (cf.\ Section \ref{sect:C2actionIrrH}).

\section{Further Questions}\label{sect:24typesfurtherquestions}

We end this chapter with some further questions, in addition to Problem \ref{prob:case21}, on detailed aspects of the structure of the 24 types.

\subsection{Realisation for fixed extension}

In Section \ref{sect:realise24ways} we considered the realisability of each type for some quadratic extension and finite Real structure. The first question is to obtain a finer understanding once the field extension is fixed.

\begin{prob}
Suppose that $\bL / \bK$ is a given quadratic field extension. Which of the 24 types occur for some finite group $H$, Real structure $G$ on $H$, and irreducible representation $V \in \Irr_{\bbL}(H)$?
\end{prob}

For example, if $\bL / \bK$ is a number field, it is natural to ask which arithmetic properties of $\bL / \bK$ determine which types occur. The following shows that when $\bK$ is totally real, not all of the twelve types of Section \ref{sect:LsplitsHG} are realised.

\begin{eg}
Suppose that $\bL$ is a CM number field, $\bK$ is the index two totally real subfield, and $V \in \Irr_{\bbL}(H)$ is absolutely irreducible. Then the D-type of $V$ cannot be II.V.

Indeed, if the linear type of $V$ is (B, A, 1, 3(b)), then by Table \ref{tab:24types6}, $Z$ is also absolutely irreducible, and the pair $(Z_{\bC}, V_{\bC})$ correspond over $\bC$. We have that $\Stab_{K_4}(Z_{\bC}) = \Stab_{K_4}(Z) = \langle c \rangle$ and $\Stab_{K_4}(V_{\bC}) = \Stab_{K_4}(V) = K_4$, and therefore the linear type of $V_{\bC}$ is also (B, A, 1, 3(b)) from Table \ref{tab:24types4}. This forces $V_{\bC}$ to have type II, as we will see in Theorem \ref{thm:dyson10}, and therefore $V_{\bC}$ does not extend to a semilinear representation of $G$ over $\bC$. This means that $V$ also doesn't extend to a semilinear representation of $G$: any such extension $W$ would define an semilinear extension $W_{\bC}$ of $V_{\bC}$. This forces $m_{\bK}^{\bL}(V) = 1$, and therefore $V$ to have D-type II, not D-type II.V.
\end{eg}

\begin{eg}
Suppose that $H$ has exponent $e$ and $G$ has exponent $f$. Then the possible types for the cyclotomic extensions $\bQ(\zeta_e)_+ \leq \bQ(\zeta_e)$ and $\bQ(\zeta_f)_+ \leq \bQ(\zeta_f)$ are described by Theorem \ref{thm:15waysforcyclotomic}, which also provides an computable method to determine the type.
\end{eg}

\subsection{Realisation of Brauer-Wall classes}

It is well-known that not all classes of the Brauer group of a field $\bK$ are realised as endomorphism rings of irreducible representations of finite groups over $\bK$, such classes which are realised in the way defining the Schur subgroup of $\Br(\bK)$, which is typically not equal to all of $\Br(\bK)$.

On the other hand, in \cite[\S 16]{TAYSchur} it is shown that any class $[A] \in \Br(\bK)$ is realised as $[\End_{\bbL_{\sigma}[G]}(V)]$ for some finite Galois extension $\bL / \bK$, group $G$, surjective homomorphism $\sigma \colon G \rightarrow \Gal(\bL / \bK)$ and $V \in \Irr_{\bL}^{\sigma}(G)$. The group $G$ can furthermore be taken to be finite if $[A]$ has order coprime to the characteristic of $\bK$. 

When $\bL / \bK$ has degree two, we saw in Section \ref{sect:BrauerWall} a way to associate to any S-block an element
\[
[D(W)] \in \BrW(\bK)
\]
whenever $\widehat{Z}(D(W)) = \bK$ (for example, this happens whenever $Z(\End_{\bbL_{\star \hspace{-0.1em}}[G]}(W)) = \bK$). When $\bL / \bK$ is $\bC / \bR$, $\BrW(\bR) \cong \bZ / 8 \bZ$ (see Proposition \ref{prop:gradeddivalgoverR} below), and this association is exhaustive.

It is natural to ask if the same holds for any quadratic $\bL / \bK$, and if not to determine the image.

\begin{prob}
Suppose that $\bL / \bK$ is a given quadratic field extension. Which classes of $\BrW(\bK)$ are realised for S-blocks of finite Real structures $H \leq G$ over $\bL / \bK$?
\end{prob}

\subsection{Prime degree classification}

Let $p$ be a prime number. Much of the theory of Section \ref{sect:fourfoldway} and this entire chapter goes through for prime degree extensions $\bL / \bK$, provided that $\bK$ contains a primitive $p$th root of $1$. 

\begin{prob}
Generalise the classification into 24 types to arbitrary prime degree $\bL / \bK$.
\end{prob}

More generally, one could try to obtain a similar classification into types of semilinear representations for $\Gal(\bL/ \bK)$ a fixed group, the case here being for the group $C_2$.

%% file: chapter-tenways.tex
\chapter{Antilinear Representations: The Case of $\bC / \bR$}\label{chapter:tenfoldway}

In this chapter we continue with the setup of Chapter \ref{chapter:twelvefoldway}, and specialise to the case when $\bL / \bK$ is the extension $\bC / \bR$. In this case, we refer to semilinear representations as \emph{antilinear} representations.

\begin{itemize}
    \item 
In Section \ref{sect:antilinearoverview} we recall the key properties that antilinear representations satisfy as opposed to general semilinear representations.

\item In Section \ref{sect:recoverDyson10} we compare Dyson's Tenfold Way and our classification into types.

\item In Section \ref{sect:antilinearvsrealreps} we consider results for $\bR$-linear representations and use Dyson's Tenfold Way to understand how these extend to results for general antilinear representations.

\item In Section \ref{sect:realisablecyclo}, we prove an induction theorem for antilinear representations and use this to show that antilinear representations of finite groups admit models over cyclotomic fields.

\item In Section \ref{sect:unitaryreps} we extend the theory of the Frobenius-Schur indicator to unitary representations of certain infinite groups, including compact groups and discrete amenable groups.

\item In Section \ref{sect:deponrealstr} we study for a fixed group $H$ how the D-Type of a representation of $H$ varies as we alter the Real structure on $H$.
\end{itemize}

\section{Overview of Antilinear Representations}\label{sect:antilinearoverview}

The two main facts which restrict the structure of antilinear representations as compared to general semilinear representations for degree two extensions are that

\begin{itemize}
    \item $\bC$ is algebraically closed, so any $V \in \Irr_{\bC}(H)$ is absolutely irreducible,
    \item Any finite-dimensional division algebra over $\bR$ is either $\bR$, $\bC$ or $\bH$ (Frobenius' Theorem).
\end{itemize}

From the absolute irreducibility, the possible types that can occur are a subset of those twelve types of Definition \ref{def:twelvetypes}, and so for antilinear representations we use the D-type notation instead of the type. It actually happens that D-types II.V and IX.V cannot occur (and the remaining ten D-types do), but this fact, known as Dyson's Tenfold Way, is not immediate, and relies of the classification of $C_2$-graded division algebras over $\bR$ (see Section \ref{sect:recoverDyson10} below).

From the classification of finite-dimensional division algebras over $\bR$, the three possible S-types for an antilinear representation $W$ are characterised by the endomorphism ring of $W$ being $\bR$, $\bC$ or $\bH$, and in particular the division algebras appearing in Table \ref{tab:24types6} are either $\bR$, $\bC$ or $\bH$.

The other main feature of the case when $\bL / \bK$ is $\bC / \bR$ is the existence of the Frobenius-Schur indicator. If $H$ is finite, and $V \in \Irr_{\bC}(H)$ has character $\chi$, we define
\[
\sF^{\sharp}(V) \coloneqq \sF^{\sharp}(\chi) \coloneqq \sum_{g \in G \setminus H} \chi(g^2).
\]

This invariant allows one to compute the S-type (cf.\ Theorem \ref{thm:threeoptions}) of any $V \in \Irr_{\bC}(H)$.

\begin{prop}{\cite[Thm.\ 4.2]{RUMTAY}}\label{prop:FSfinitegroups}
Suppose that $H$ is finite and $V \in \Irr_{\bC}(H)$. Then:
\begin{itemize}
    \item $\fS(V) = \mathbf{A}$ if and only if $\sF^{\sharp}(V) = 1$,
    \item $\fS(V) = \mathbf{B}$ if and only if $\sF^{\sharp}(V) = 0$,
    \item $\fS(V) = \mathbf{C}$ if and only if $\sF^{\sharp}(V) = -1$.
\end{itemize}
\end{prop}

As a consequence of the induction theorem, we prove in Section \ref{sect:realisablecyclo}, we obtain that Proposition \ref{prop:FSfinitegroups} also applies when $H$ is finite to certain cyclotomic extensions (Theorem \ref{thm:15waysforcyclotomic}, Remark \ref{rem:FSforcyclo}).

We also extend Proposition \ref{prop:FSfinitegroups} to certain infinite groups in Section \ref{sect:unitaryreps}.

\section{The Tenfold Way for $\bC / \bR$}\label{sect:recoverDyson10}

In this section we give a simple proof of the tenfold way for $\bC / \bR$ by showing that for $\bC / \bR$ cases II.V and IX.V of Definition \ref{def:twelvetypes} cannot occur. The key input is the classification of $C_2$-graded division algebras over $\bR$, which we now recall.

By Lemma \ref{lem:centreofC2graddiv}, the graded centre of any $C_2$-graded division algebra over $\bR$ is either $\bR$ or $\bC$, and so any $C_2$-graded division algebra can be viewed as an element of $\BrW(\bR)$ or $\BrW(\bC)$ (see Definition \ref{def:BW}). As with the ordinary Brauer group, any class of $\BrW(\bF)$ has a unique representative by a central graded division algebra over $\bF$ by Corollary \ref{cor:divalgrep}, and therefore the set of $C_2$-graded division algebras over $\bR$ is canonically identified with the disjoint union
\[
\BrW(\bR) \sqcup \BrW(\bC).
\]

From the exact sequence \eqref{eqn:exactseqbrauer} of Appendix \ref{appendix:gradedalgebras}, we have that $\BrW(\bC) \cong \bZ / 2 \bZ$, and that $\BrW(\bR)$ fits into an exact sequence
\begin{equation}\label{eqn:exactBrWR}
    0 \rightarrow \bZ / 2 \bZ \rightarrow \BrW(\bR) \rightarrow \bZ / 4 \bZ \rightarrow 0,
\end{equation}
as $\Br(\bR) \cong \bZ / 2 \bZ$ and $Q(\bR)$ is cyclic with generator $(1,1)$.

To produce non-trivially graded elements of $\BrW(\bR)$ and $\BrW(\bC)$, we consider the following.

\begin{defn}
For an finite-dimensional division algebra $A$ over $\bR$, we set
\[
A[\epsilon] \coloneqq A[\epsilon] / (\epsilon^2 - 1), \qquad A[\delta] \coloneqq A[\delta] / (\delta^2 + 1),
\]
where $\epsilon$ and $\delta$ are variables which commute with all elements of $A$. We view these as a $C_2$-graded $\bR$-algebras where $A$ has degree $0$ and $\epsilon$ and $\delta$ have degree $1$.

When $A = \bC$, we also write 
\[
\bC_{\hspace{-0.1em}\rtimes \!}[\epsilon], \qquad \bC_{\hspace{-0.1em}\rtimes \!}[\delta] 
\]
for the same abelian groups, but now viewed as $\bR$-algebras via the skew multiplication:
\[
\epsilon z = \overline{z} \epsilon, \: \epsilon^2 = 1, \qquad \qquad \delta z = \overline{z} \delta, \: \delta^2 = -1.
\]
We view $\bC_{\hspace{-0.1em}\rtimes \!}[\epsilon]$ and $\bC_{\hspace{-0.1em}\rtimes \!}[\delta]$ as $C_2$-graded $\bR$-algebras where $\bC$ has degree $0$ and $\epsilon$ and $\delta$ have degree $1$.
\end{defn}

We note that $\widehat{Z}(A[\epsilon]) = Z(A) = \widehat{Z}(A[\delta])$, and $\widehat{Z}(\bC_{\hspace{-0.1em}\rtimes \!}[\epsilon]) = \bR = \widehat{Z}(\bC_{\hspace{-0.1em}\rtimes \!}[\delta])$.

We can use this to state the classification of $C_2$-graded division algebras over $\bR$. Motivated by Corollary \ref{cor:Dalphaisom}, for such a $D$ we also compute the algebra $D_i$, which we use in Theorem \ref{thm:dyson10} below.

\begin{prop}\label{prop:gradeddivalgoverR}
The $C_2$-graded division algebras over $\bR$ up to isomorphism are listed in the first column of Table \ref{tab:graddivoverR}. The group $\BrW(\bR)$ is generated by $\bR[\epsilon]$, and for $1 \leq k \leq 8$, the $C_2$-graded division algebra over $\bR$ in the $k$th row is equivalent in $\BrW(\bR)$ to the $k$th power of $\bR[\epsilon]$.

For $D$ such an algebra, the graded centre $\widehat{Z}(D)$, $j(D) \in Q(\widehat{Z}(D))$, ungraded part $D_0$, isomorphism class as semisimple $\bR$-algebra $D_\text{ss}$, and associated $\bR$-algebra $D_i$, are as in Table \ref{tab:graddivoverR}. 
\end{prop}

\begin{table}
\caption{D-Types and $C_2$-graded division algebras over $\bR$}
\label{tab:graddivoverR}
\begin{center}
\def\arraystretch{1.1}
\begin{tabular}{ c | c c c c c | c} 
    \toprule
$D$ & $\widehat{Z}(D)$ & $j(D)$ & $D_0$ & $D_{\text{ss}}$ & $D_{i}$ & $\fD(W)$ \\
\midrule
 $(\bR)$ & $\bR$ & $(0,1)$ & $\bR$ & $\bR$ & $\bC$ & V \\ 
 $\bR[\epsilon]$ & $\bR$ & $(1,1)$ & $\bR$ & $\bR \times \bR$ & $M_2(\bR)$ & I \\ 
 $\bC_{\hspace{-0.1em}\rtimes \!}[\epsilon]$ & $\bR$ & $(0,-1)$ & $\bC$  & $M_2(\bR)$ & $M_2(\bR) \times M_2(\bR)$ & III \\ 
 $\bH[\delta]$ & $\bR$ & $(1,-1)$ & $\bH$  & $M_2(\bC)$ & $M_4(\bR)$ & II \\ 
 $(\bH)$ & $\bR$ & $(0,1)$ & $\bH$  & $\bH$ & $M_2(\bC)$ & VI \\ 
 $\bH[\epsilon]$ & $\bR$ & $(1,1)$ & $\bH$  & $\bH \times \bH$ & $M_2(\bH)$ & VIII \\ 
 $\bC_{\hspace{-0.1em}\rtimes \!}[\delta]$ & $\bR$ & $(0,-1)$ & $\bC$  & $\bH$ & $\bH \times \bH$ & X \\ 
 $\bR[\delta]$ & $\bR$ & $(1,-1)$ & $\bR$  & $\bC$ & $\bH$ & IX \\ 
 $(\bC)$ & $\bC$ & $(0,1)$ & $\bC$  & $\bC$ & $\bC \times \bC$ & VII \\ 
 $\bC[\epsilon]$ & $\bC$ & $(1,1)$ & $\bC$  & $\bC \times \bC$ & $M_2(\bC)$ & IV \\
  \bottomrule
\end{tabular}
\end{center} 
\end{table}

\begin{proof}
    It is clear that all the listed $C_2$-graded algebras $D$ are $C_2$-graded division algebras, and the values of $\widehat{Z}(D)$, $j(D)$, and $D_0$ are direct to compute. It is also clear that no pair are isomorphic as $C_2$-graded algebras, by considering the invariant $j(D)$ and the $\bR$-dimension. Conversely, given any $C_2$-graded division algebra $D$ over $\bR$, this defines an element $[D] \in \BrW(\bR) \sqcup \BrW(\bC)$ as above. These have orders $|\BrW(\bR)| = 8$ and $|\BrW(\bC)| = 2$, and therefore $[D] = [D']$ for some $D'$ on our list, because they all define distinct classes by the injectivity of Corollary \ref{cor:divalgrep} and therefore must exhaust $\BrW(\bR) \sqcup \BrW(\bC)$. But then $D \cong D'$, again by the injectivity of Corollary \ref{cor:divalgrep}.
    
    To understand the group structure of $\BrW(\bR)$ in terms of these representatives, note that there is a $C_2$-graded isomorphism
    \[
    \bR[\epsilon] \, \hat{\otimes} \, \bR[\epsilon] \xrightarrow{\sim} \bC_{\hspace{-0.1em}\rtimes \!}[\epsilon],
    \]
    where $1 \, \hat{\otimes} \, 1$, $1 \, \hat{\otimes} \, \epsilon$, $\epsilon \, \hat{\otimes} \, 1$ and $\epsilon \, \hat{\otimes}\, \epsilon$ map to $1$, $\epsilon$, $i \epsilon$ and $i$ respectively. The graded opposite algebra (the inverse in the Brauer-Wall group) of $\bC_{\hspace{-0.1em}\rtimes \!}[\epsilon]$ is $\bC_{\hspace{-0.1em}\rtimes \!}[\delta]$, which follows from the definition, and therefore $\bC_{\hspace{-0.1em}\rtimes \!}[\epsilon]$ is not self-inverse, hence has order neither $1$ or $2$, so $\bR[\epsilon]$ must have order $8$. This implies that $4 \cdot [\bR[\epsilon]] = [(\bH)]$, and the description of the remaining rows as powers of $\bR[\epsilon]$ follows from noting that $\bR[\delta]$ is the graded opposite algebra of $\bR[\epsilon]$ (so $[\bR[\delta]]$ is the inverse of $[\bR[\epsilon]]$), and that
    \[
        \bH[\epsilon] \cong (\bH) \, \hat{\otimes} \, \bR[\epsilon], \qquad \bH[\delta] \cong (\bH) \, \hat{\otimes} \, \bR[\delta].
    \]

    It remains to compute $D_{\text{ss}}$ and $D_i$. For the ungraded rows, it is clear that $(\bR)_{i} \cong \bC$ and
    \[
    (\bC)_i \cong \bC[x] / (x^2 + 1) \cong \bC \times \bC,
    \]
    and for $(\bH)$ there is an $\bR$-algebra isomorphism
\[
(\bH)_{i} \xrightarrow{\sim} M_2(\bC),
\]
where
\[
\mathbf{i} \mapsto \begin{pmatrix}
    0 & -1 \\ 1 & 0
\end{pmatrix}, \quad \mathbf{j} \mapsto \begin{pmatrix}
    0 & i \\ i & 0
\end{pmatrix}, \quad \mathbf{k} \mapsto \begin{pmatrix}
    -i & 0 \\ 0 & i
\end{pmatrix}, \quad \underline{i} \mapsto \begin{pmatrix}
    i & 0 \\ 0 & i
\end{pmatrix}.
\]

For rows $2, 6$ and $10$, suppose that $A \in \{\bR, \bC, \bH\}$. Then as $Z(A)$-algebras,
\[
A[\epsilon] \xrightarrow{\sim} A \times A, \qquad x + y \epsilon \mapsto (x+y, x-y).
\]
Furthermore,
\[
A[\epsilon]_i = A[\epsilon] \oplus A[\epsilon] \cdot \underline{i} \xrightarrow{\sim} M_2(A),
\]
where
\[
(x + y \epsilon, (u + v \epsilon) \cdot \underline{i}) \mapsto \begin{pmatrix}
    x + y & -u - v \\ u - v & x - y
\end{pmatrix}.
\]

All that remains are rows $3$, $4$, $7$ and $8$. For $\bC_{\hspace{-0.1em}\rtimes \!}[\epsilon]$, there is an isomorphism
\[
\bC_{\hspace{-0.1em}\rtimes \!}[\epsilon] \xrightarrow{\sim} M_2(\bR), \qquad (a+ib + (c + id) \epsilon) \mapsto \begin{pmatrix}
    x + y & -u - v \\ u - v & x - y
\end{pmatrix},
\]
which comes from noting that $\bC_{\hspace{-0.1em}\rtimes \!}[\epsilon]$ is nothing but $\bC_{\star \hspace{-0.1em}}[C_2]$ for the trivial real structure $1 \leq C_2$: the isomorphism is the action map $\bC_{\star \hspace{-0.1em}}[C_2] \xrightarrow{\sim} \End_{\bR}(\bC)$ on the trivial module $\bC$, and taking the basis $\{1, i\}$. To see that $\bC_{\hspace{-0.1em}\rtimes \!}[\epsilon]_i \cong M_2(\bR) \times M_2(\bR)$, we first note that $\bC_{\hspace{-0.1em}\rtimes \!}[\epsilon]_i$ is a semisimple $\bR$-algebra, which follows by \cite[Def.\ IV.3.5, Thm.\ IV.3.6]{LAM} from the fact that it is $C_2$-graded simple (Lemma \ref{lem:Aalphaproperties}). From Lemma \ref{lem:Aalphaproperties} we also have that
\[
Z(\bC_{\hspace{-0.1em}\rtimes \!}[\epsilon]_i) = \bR \oplus \bR \cdot i \cdot \underline{i} \cong \bR[t] / (t^2 - 1) \cong \bR \times \bR
\]
so $\bC_{\hspace{-0.1em}\rtimes \!}[\epsilon]_i \cong A \times A'$ for some central simple $\bR$-algebras $A, A'$. The composition of the inclusion $M_2(\bR) \hookrightarrow \bC_{\hspace{-0.1em}\rtimes \!}[\epsilon]_i$ followed by the projection to either factor is injective (by the simplicity of $M_2(\bR)$), hence $A \cong M_2(\bR) \cong A'$ by dimension considerations.

For $\bC_{\hspace{-0.1em}\rtimes \!}[\delta]$, there is a non-zero $\bR$-algebra homomorphism
\[
\bH \xrightarrow{\sim} \bC_{\hspace{-0.1em}\rtimes \!}[\delta], \qquad \mathbf{i} \mapsto i, \: \mathbf{j} \mapsto \delta, \: \mathbf{k} \mapsto i \delta,
\]
which is injective by the simplicity of $\bH$, hence an isomorphism. Similarly to above, $\bC_{\hspace{-0.1em}\rtimes \!}[\delta]_i$ is a semisimple $\bR$-algebra with centre $\bR \times \bR$ and an $\bR$-subalgebra $\bH \hookrightarrow \bC_{\hspace{-0.1em}\rtimes \!}[\delta]_i$ hence $\bC_{\hspace{-0.1em}\rtimes \!}[\delta]_i \cong \bH \times \bH$.

For $\bR[\delta]$, it is clear that $\bR[\delta] \cong \bC$ as $\bR$-algebras, and the isomorphism $\bR[\delta]_i \cong \bH$, follows because $\bR[\delta]_i$ is nothing but $\bC_{\hspace{-0.1em}\rtimes \!}[\delta]$.

Finally, for $\bH[\delta]$, $\bH[\delta] \cong (\bH)_i$, where $\delta$ corresponds to $\underline{i}$, and therefore $\bH[\delta] \cong M_2(\bC)$ from the case of $(\bH)$ above. To see that $\bH[\delta]_i \cong M_4(\bR)$, note that the centraliser of $\bH$ in $\bH[\delta]_i$ is 
\[
\bR[\delta]_i = \bR[\delta] \oplus \bR[\delta] \cdot \underline{i} \leq \bH[\delta]_i.
\]
We have that $\bR[\delta]_i \cong \bH$ from the case of $\bR[\delta]$ above, and therefore there is an $\bR$-algebra homomorphism
\[
\bH \otimes_{\bR} \bH \rightarrow \bH[\delta]_i,
\]
Because $\bH$ has order two in $\Br(\bR)$, $\bH \otimes_{\bR} \bH \cong M_4(\bR)$, and therefore the $\bR$-algebra homomorphism
\[
M_4(\bR) \xrightarrow{\sim} \bH \otimes_{\bR} \bH \rightarrow \bH[\delta]_i
\]
is an isomorphism, as $M_4(\bR)$ is simple and $\bH[\delta]_i$ has the same $\bR$-dimension.
\end{proof}

We can use this to give a short proof of Dyson's Tenfold Way.

\begin{thm}[Dyson's Tenfold Way]\label{thm:dyson10}
Let $W \in \Irr_{\bC}^{\hspace{0.1em} \star}(G)$. Then the pair $(D(W),\fD(W))$ is as in one of the rows of Table \ref{tab:graddivoverR}, for $D \coloneqq D(W)$. In particular, $W$ cannot have D-type II.V or IX.V, and $\fD(W)$ uniquely determines and is determined by the isomorphism type of $D(W)$ as a $C_2$-graded $\bR$-algebra.
\end{thm}

\begin{proof}
    We have that as $\bR$-algebras, $D(W)_0 = \End_{\bC_{\star \hspace{-0.1em}}[G]}(W)$, $D(W)_{\text{ss}} = \End_{\bR[G]}(W)$, and $D(W)_i \cong \End_{\bR[H]}(W)$ by Corollary \ref{cor:Dalphaisom}, and therefore the triples 
    \[
    (D(W)_0, D(W)_{\text{ss}}, D(W)_i) \:\: \text{ and } \: \: (\End_{\bC_{\star \hspace{-0.1em}}[G]}(W), \End_{\bR[G]}(U), \End_{\bR[H]}(X))
    \]
    have Brauer-equivalent simple factors in each component. If $W$ has D-type II.V, this triple would be $(\bR,\bC,\bR)$, using Table \ref{tab:12types}, which occurs for no line of Table \ref{tab:graddivoverR}, and therefore this D-type cannot occur. Similarly, if $W$ has D-type IX.V, this tripe would be $(\bH, \bC, \bH)$, which does not appear in Table \ref{tab:graddivoverR}.

    For the remaining 10 D-types, we can similarly read off the isomorphism class of $D(W)$ as a $C_2$-graded division algebra. For example, if $W$ has D-type II, then by Table \ref{tab:12types} the triple is $(\bH, \bC, \bR)$, which shows that $D(W) \cong \bH[\delta]$ as $C_2$-graded algebras. This determines $D(W)$ in all D-types apart from when $W$ has D-type IV or VII, which both have triple $(\bC, \bC, \bC)$, and therefore there are two options for $D(W)$ as either $(\bC)$ or $\bC[\epsilon]$. These can be distinguished by noting that $D(W)$ has $D(W)_1 = 0$ if and only if $W \not\cong \overline{W}$ (Proposition \ref{prop:basicpropDW}), which occurs in case VII and not in case IV (using Table \ref{tab:24types5}), and so $D(W) \cong \bC[\epsilon]$ in case IV and $D(W) \cong (\bC)$ in D-type VII.
\end{proof}

\begin{remark}
    All of the ten D-types of Table 4.1 are realised by irreducible antilinear representations of finite groups \cite[Thm.\ 3.6]{RUMTAY}, and therefore any $C_2$-graded division algebra over $\bR$ is realised as $D(W)$ for some irreducible antilinear representation $W$ of a finite group by Theorem \ref{thm:dyson10}.
\end{remark}

We can use Theorem \ref{thm:dyson10} to alternatively characterise the 10 D-types in terms of $A(W)$.

\begin{cor}\label{cor:DysonforAW}
If $W \in \Irr_{\bC}^{\hspace{0.1em} \star}(G)$, then $\widehat{Z}(A(W)) = \widehat{Z}(D(W))$. When $\widehat{Z}(D(W)) = \bR$, 
\[
        [A(W)] \cdot [D(W)] = [\bR[\epsilon]]
\]
in $\BrW(\bR)$. When $\widehat{Z}(D(W)) = \bC$, $[A(W)] = [D(W)]$ in $\BrW(\bC)$.
\end{cor}

\begin{proof}
    By Proposition \ref{prop:relateAWDW}, $\widehat{Z}(A(W)) = \bR$ if and only if $\widehat{Z}(D(W)) = \bR$, in which case
    \[
        [A(W)] \cdot [D(W)] = [\bR[\epsilon]]
    \]
    in $\BrW(\bR)$. That $\widehat{Z}(D(W)) = \bR$ occurs in all D-types but IV and VII, and so it only remains to understand $A(W)$ in D-types IV and VII, for which $\widehat{Z}(A(W)) \neq \bR$.

    When $W$ has D-type IV, $X \cong {}^{\omega} X$ and $U \not\cong U \otimes \pi$ from Table \ref{tab:24types5}, and thus $A(W)_0$ is simple with $Z(A(W)_0) \cong \bC$ and $Z(A(W)) \cong \bC \times \bC$ by Proposition \ref{prop:descofSblock} and Lemma \ref{lem:basicpropertiesAW}. In particular, $\widehat{Z}(A(W)) \leq Z(A(W)_0) \cong \bC$ by Lemma \ref{lem:centreofC2graddiv}, and so $\widehat{Z}(A(W)) \cong \bC$. Because $Z(A(W)) > \widehat{Z}(A(W))$, $A(W)$ has \emph{odd type} (see \cite[Def.\ IV.3.5]{LAM}), and thus $[A(W)] = [\bC[\epsilon]]$ in $\BrW(\bC)$, as $[\bC[\epsilon]]$ is the unique central $C_2$-graded division algebra over $\bC$ of odd type, and being of odd type depends only on the class in $\BrW(\bC)$ \cite[Thm.\ IV.4.3]{LAM}.

    Similarly, when $W$ has D-type VII, $X \not\cong {}^{\omega} X$ and $U \cong U \otimes \pi$ from Table \ref{tab:24types5}, so $A(W)$ is simple with $Z(A(W)) = \bC$ by Proposition \ref{prop:descofSblock} and Lemma \ref{lem:basicpropertiesAW}. Now $\widehat{Z}(A(W)) \leq Z(A(W)) \cong \bC$ by Lemma \ref{lem:centreofC2graddiv}, so $\widehat{Z}(A(W)) \cong \bC$, but now $Z(A(W)) = \widehat{Z}(A(W))$, so $A(W)$ has \emph{even type}, and thus $[A(W)] = [(\bC)]$ in $\BrW(\bC)$.
\end{proof}

\begin{remark}
    The paper \cite{GeiMoo} also computes the class of $A(W)$ and $D(W)$ in the appropriate Brauer-Wall group for each of the ten D-types. What we call $A(W)$ is denote $\sA$ there, and their algebra $\sW$ is our $\widetilde{D}(W)$, which by Proposition \ref{prop:relateAWDW} is isomorphic to $D(W)$.
    
    Whilst it allows us to compute the classes of $A(W)$ and $D(W)$, for us the main interest in Theorem \ref{thm:dyson10} is the fact that it gives a very simple proof of Dyson's Tenfold Way using only the classification of $C_2$-graded division algebras over $\bR$ and Corollary \ref{cor:Dalphaisom}, which holds for general types.
    
    As was observed in \cite[\S 5]{GeiMoo}, one can define a ``double'' Frobenius-Schur indicator,
    \[
    \sF'(W) \coloneqq \frac{\sF(W) + i \tilde{\sF}(W)}{|\sF(W) + i \tilde{\sF}(W)|} \in \{e^{\frac{2 \pi i n}{8}}\}_{n \in \bZ / 8 \bZ},
    \]
    which has $\sF'(W) = 0$ when $W$ has D-type IV or VII, and otherwise determines the D-type of $W$. Now,
    \[
        \sF'(W) = e^{\frac{2\pi i n}{8}}, \: \text{ where $n \in \bZ / 8 \bZ$ is defined by } [D(W)] = n \cdot [\bR[\epsilon]]
    \]
    by Theorem \ref{thm:dyson10} and therefore there is an identification of abelian groups of $\{e^{\frac{2 \pi i n}{8}}\}_{n \in \bZ / 8 \bZ}$ with $\BrW(\bR)$. Considering instead $A(W)$, we obtain a second identification of these groups from the relation
    \[
    \sF'(W) = e^{\frac{2\pi i n}{8}}, \: \text{ where $n \in \bZ / 8 \bZ$ is defined by } [A(W)] = (1-n) \cdot [\bR[\epsilon]]
    \]
    by Corollary \ref{cor:DysonforAW}. The two identifications are related by a reflection in the line $\langle 1 + i \rangle_{\bR}$.
\end{remark}

\section{Comparison between Antilinear Representations and Representations over $\bR$}\label{sect:antilinearvsrealreps}

In this section we review some well-known properties of real-representations and observe to what extent these generalise to the antilinear setting. Throughout $H$ is a finite group.

\subsection{Dimension of complex representations of C-type $\bH$}

\begin{prop}\label{prop:evendim}
    Suppose that $V \in \Irr_{\bC}(H)$ has C-type $\bH$. Then $\dim_{\bC}(V)$ is even.
\end{prop}

We note that, unlike the other results in this section, this doesn't actually require that $H$ is finite.

\begin{proof}
    If $V$ corresponds to $W \in \Irr_{\bR}(H)$, then $\dim_{\bR}(W) = 2 \dim_{\bC}(V)$, and $4 \mid \dim_{\bR}(W)$ because $W$ is a module over $\bH$, hence $\dim_{\bC}(V)$ is even.
\end{proof}

This is not true for antilinear representations, which are now complex vector spaces.

\begin{eg}\label{eg:dim1typeH}
    Let $C_2$ have Real structure $C_2 \leq C_4$, and let $V$ be the non-trivial $1$-dimensional representation of $C_2$. Then, writing $C_2 = \langle x \rangle \leq \langle y \rangle = C_4$ with $y^2 = x$, we many compute
    \[
        \sF(V) = \frac{1}{2} \left( \chi_V(y^2) + \chi_V((y^3)^2) \right) = \frac{1}{2} \left( \chi_V(x) + \chi(x) \right) = -1
    \]
    hence $V$ has S-type $\bH$, but dimension $1$.
\end{eg}

\subsection{Abelian groups have no complex representations of C-type $\bH$}

\begin{cor}
Suppose that $H$ is abelian. Then any $V \in \Irr_{\bC}(H)$ has C-type $\bR$ or $\bC$.
\end{cor}

\begin{proof}
    This follows from Proposition~\ref{prop:evendim} or can be seen directly noting that $\bR[H]$ is commutative, so there can be no copy of $\bH$ in its Artin-Wedderburn decomposition.
\end{proof}

This is not true for antilinear representations: from Example \ref{eg:dim1typeH}, $C_2 \leq C_4$ provides an example of Real structure $H \leq G$ where both $H$ and $G$ are abelian, but there is a representation of S-type $\bH$.

\subsection{Burnside's Theorem}

As a consequence of Proposition \ref{prop:evendim} we have the following.

\begin{cor}\label{cor:weakburnside}
Suppose that $|H|$ is odd. Then any $V \in \Irr_{\bC}(H)$ has C-type $\bR$ or $\bC$.
\end{cor}

\begin{proof}
    If there was $V \in \Irr_{\bC}(H)$ of C-type $\bH$, then $V$ has even dimension by Proposition \ref{prop:evendim}, which is impossible as $\dim_{\bC}(V) \mid |H|$ because $V$ is absolutely irreducible.
\end{proof}

Burnside proved the following stronger statement.

\begin{prop}\label{prop:burnside}
    Suppose that $|H|$ is odd. Then any $V \in \Irr_{\bC}(H)$ which is non-trivial has C-type $\bC$.
\end{prop}

\begin{proof}
    As $H$ has odd order, $h \mapsto h^2$ is injective, hence bijective on $H$, so
    \[
        \sF(V) = \sum_{h \in H} \chi(h^2) = \sum_{h \in H}\chi(h) = \dim_{\bC} V^H
    \]
    which is $0$ apart from the case when $V$ is the trivial representation. Alternatively, one can use \cite[Cor.\ 5.19]{RUMTAY} for the Real structure $H \leq H \times C_2$: the odd order condition means that there is only one odd element which squares to $1$, and hence, given that there are no irreducible characters of C-type $\bH$ by Proposition \ref{prop:evendim}, there also are no non-trivial characters of C-type $\bR$.
\end{proof}

The following example shows that Proposition \ref{prop:burnside} is not true for antilinear representations.

\begin{eg}
Consider $C_3 \leq D_6$, and let $V$ be any non-trivial $1$-dimensional representation of $C_3$ over $\bC$. Writing $C_3 = \langle x \rangle \leq \langle x, y \mid x^3, y^2, yxy^{-1} = x^{-1} \rangle$, we can compute that
\[
    \sF(V) = \frac{1}{3} \left( \chi_V(y^2) + \chi_V((xy)^2) + \chi_V((x^2y)^2) \right)
    = \frac{1}{3} \left( \chi_V(1) + \chi_V(1) + \chi_V(1) \right) = 1,
\]
and so the S-type of any $V \in \Irr_{\bC}(C_3)$ for $C_3 \leq D_6$ is $\bR$.
\end{eg}

\subsection{Order of a group with a quaternionic representation}

One can go even further than Corollary \ref{cor:weakburnside}, which says that if there is some $V \in \Irr_{\bC}(H)$ of C-type $\bH$, then $|H|$ is even.

\begin{prop}\label{prop:strongweakburnside}
    Suppose that there is some $V \in \Irr_{\bC}(H)$ of C-type $\bH$. Then $4 \mid |H|$.
\end{prop}

\begin{proof}
    If instead $|H| = 4k+2$, then any Sylow $2$-subgroup of $H$ is cyclic of order $2$, and the claim follows from \cite[Thm.\ 10.9]{ISA}, which states that if $H$ has a Sylow $2$-subgroup which is elementary abelian then $H$ has no irreducible characters of C-type $\bH$. 
\end{proof}

This Proposition \ref{prop:strongweakburnside} is not true for antilinear representations, as shown by Example \ref{eg:dim1typeH}.

\subsection{Antilinear representations of groups of odd order}

We can however use Proposition \ref{prop:strongweakburnside} to show that Corollary \ref{cor:weakburnside} extends to antilinear representations.

\begin{prop}\label{prop:typesforoddorder}
    Suppose that $H$ has odd order. Then every non-trivial $V \in \Irr_{\bC}(H)$ has D-type IV, V or VII. In particular, $V$ has S-type $\bR$ or $\bC$.
\end{prop}

Note that the proof uses the strong form Proposition \ref{prop:strongweakburnside} of Corollary \ref{cor:weakburnside} from the classical case, and the tenfold way.

\begin{proof}
Because the C-type of $V$ is $\bC$, the D-type of $V$ is one of IV, V, VI or VII. In D-type VI, the $U \in \Irr_{\bR}(G)$ corresponding to $V$ has C-type $\bH$, which can't happen by Proposition \ref{prop:strongweakburnside} because the order of the Real structure $G$ is not a multiple of $4$.
\end{proof}

\section{Realisability over Cyclotomic Fields}\label{sect:realisablecyclo}

We continue in this section with the assumption that $H$ is finite. Let $e \coloneqq \exp(H)$, the exponent of $H$, when $\exp(H) > 2$ and $e = 4$ when $\exp(H) \leq 2$, and set $E \coloneqq \bQ(\zeta_e)$ to be the cyclotomic field generated by a primitive $e$th root of unity. The modification when $\exp(H) \leq 2$ is to ensure that $e > 2$, and therefore that the maximal totally real subfield
\[
E_{+} \coloneqq \bR \cap E = \bQ(\zeta_e + \zeta_e^{-1})
\]
has index two in $E$, so that we can consider semilinear $E$-representations for a Real structure $H \leq G$.

The main result of \cite{PAS} is that irreducible $\bR$-representations always have models over $E_+$. Note that when $\exp(H) = 2$, $E_+ = \bQ(i) \cap \bR = \bQ(\zeta_2) \cap \bR = \bQ$, and so our $E_+$ agrees with that of \cite{PAS}.

\begin{thm}\label{thm:classicalcyclotomic}
    Suppose that $W \in \Irr_{\bR}(H)$. Then there is $\tilde{W} \in \Irr_{E_+}(H)$ with $\tilde{W}_{\bR} \cong W$.
\end{thm}

Concretely, this means that any matrix representation of $W$ over $\bR$ can always be conjugated over $\bR$ to matrices which all have values in $E_{+}$. We give the proof for the readers convenience, adapted from that of \cite{PAS} into a form which will be easier for us to generalise to antilinear representations.

\begin{proof}
Let $V \in \Irr_{\bC}(H)$ correspond to $W$, and let $\tilde{V} \in \Irr_{E}(H)$ be a model of $V$ over $E$, so $\tilde{V}_{\bC} \cong V$ and $\chi_{\tilde{V}} = \chi_V$, which always exists because $E$ is a splitting field for $H$ \cite[\S12.3]{SERRE}.

If $V$ has C-type $\bC$ or $\bH$, then $\chi_W = \chi_V + \overline{\chi_V}$. Setting $\tilde{W} \coloneqq \tilde{V}|_{E_+}$, we have that $\chi_{\tilde{W}} = \chi_V + \overline{\chi_V}$, and thus $\tilde{W}_{\bR} \cong W$ by comparing characters.

We may therefore assume that $V$ has C-type $\bR$. There are two options for $\overline{m}_{E_+}^E(\tilde{V}) \in \{1,2\}$. If $\overline{m}_{E_+}^E(\tilde{V}) = 1$, then there is some $\tilde{W} \in \Irr_{E_+}(H)$ with $\tilde{W}_{E} \cong V$, so $\chi_{\tilde{W}} = \chi_V = \chi_W$, hence $\tilde{W}_{\bR} \cong W$ and we're done.

Suppose instead, for a contradiction, that $\overline{m}_{E_+}^E(\tilde{V}) = 2$. The key fact we use is that this means that $2 \mid \langle \chi_W , \chi_U \rangle$ for any $U \in \Rep_{E_+}(H)$ by \cite[Cor.\ 10.3(4)]{TAYSchur} and the fact that $\chi_W = \chi_V = \chi_{\tilde{V}}$. The second input is Serre's induction theorem \cite[Lem.\ 73.19]{CR2}, which says that we can write
    \[
        \chi_W = \sum_{i=1}^r a_i \cdot \Ind_{H_i}^H \psi_i
    \]
    where $r \geq 1$, $a_i \in \bZ$, each $H_i$ is a subgroup of $H$, and $\psi_i$ is a character of $W_i \in \Irr_{\bR}(H_i)$ such that $W_i$ is inflated from a cyclic or dihedral group. Because $W$ is irreducible of C-type $\bR$,
    \[
        1 = \langle \chi_W, \chi_W \rangle = \sum_{i = 1}^r a_i \langle \chi_W, \Ind_{H_i}^H \psi_i \rangle.
    \]
    
    Now each $W_i$ admits a model over $E_+$, as the irreducible $\bR$-representations of cyclic and dihedral groups are defined over $E_+$, and therefore $\Ind_{H_i}^H \psi_i$ is also the character of an object of $\Rep_{E_+}(H)$. This means that each $\langle \chi_W, \Ind_{H_i}^H \psi_i \rangle$ is even, and so $1$ is a sum of even numbers, a contradiction.
\end{proof}

The goal of this section is to prove the analogous statement for antilinear representations.

\begin{thm}\label{thm:antilinearcyclotomicfinal}
    Suppose $W \in \Irr_{\bC}^{\hspace{0.1em} \star}(G)$. Then there is $\tilde{W} \in \Irr_{E}^{\hspace{0.1em} \star}(H)$ with $\tilde{W}_{\bC} \cong W$.
\end{thm}

Concretely, this is equivalent to the claim that any antilinear matrix representation of $G$ with values in $\bC$ can always be (antilinear) conjugated over $\bC$ to matrices which all have values in $E$. In the case of split extensions this recovers Theorem \ref{thm:classicalcyclotomic}.

Let us note two cases when Theorem \ref{thm:antilinearcyclotomicfinal} is immediate. The first is when $W$ is $1$-dimensional.

\begin{prop}
Suppose that $V \in \Irr_E(H)$ is $1$-dimensional, and extends to a antilinear representation of $G$ over $\bC$. Then $V$ extends to a antilinear representation of $G$ over $E$.
\end{prop}

\begin{proof}
    We use the matrix form of semilinear representations from Section \ref{sect:matrixform} of Chapter \ref{chapter:algebras}. In this framework, a $1$-dimensional antilinear representation over $\bC$ can be described by a group homomorphism $\rho \colon G \rightarrow \GL_1^{\star}(\bC)$, where $\GL_1^{\star}(\bC) = \GL_1^1(\bC) \sqcup \GL_1^{-1}(\bC)$, where elements are written as $(a)_1$ and $(a)_{-1}$ respectively for $a \in \bC^{\times}$. Multiplication is given by
    \[
    (a)_{1} \cdot (b)_{\epsilon} = (ab)_{\epsilon} \: \text{ and } \: (a)_{-1} \cdot (b)_{\epsilon} = (a \overline{b})_{-\epsilon}.
    \]
    
    If $\rho(g) = (r_g)_{\epsilon(g)}$, then the equivalent matrix representations are $\rho_x$ for some $x \in \bC^{\times}$, where
    \[
        \rho_x(g) = \twopartdef{(r_g)_1}{\epsilon = 1}{(\overline{x} x^{-1} r_g)_{-1}}{\epsilon = -1}.
    \]

    The map $\bC^{\times} \rightarrow \bC^{\times}$, $x \mapsto \overline{x} x^{-1}$ has image $S^1$, and so the equivalent representations are parametrised by $y \in S^1$, where
    \[
    \rho^y(g) = \twopartdef{(r_g)_1}{\epsilon(g) = 1}{(y r_g)_{-1}}{\epsilon(g) = -1}.
    \]

    Let $\chi \colon H \rightarrow E^{\times}$ be the character of $V$, and suppose that $\rho \colon G \rightarrow \GL_1^{\star}(\bC)$ extends $\chi$, meaning that $r_g = \chi(g)$ for all $g \in H$. To show that $\chi$ extends to an antilinear representation over $E$ is to ask that one can find $y \in S^1$ with $\rho^y$ valued in $\GL_1^{\star}(E)$, meaning that all $y r_g \in E^{\times}$.

    For any $g_1, g_2 \in G \setminus H$, $g_1 g_2 \in H$, so $(r_{g_1} \overline{r_{g_2}})_1 = \rho(g_1g_2) = (r_{g_1 g_2})_1$. We know that $r_{g_1 g_2}$ is an $e$th root of $1$, and therefore (as $e$th roots are closed under multiplication and complex conjugation) $r_{g_1}$ is an $e$th root of $1$ if and only if $r_{g_2}$ is an $e$th root of $1$.
    
    In particular, we may fix any $\omega \in G \setminus H$, and set $y = r_\omega^{-1}$. Then $r^y(\omega) = (1)_{-1}$, and $\rho^y$ is valued in $\GL_1^{\star}(E)$. To see that this is well-defined, meaning that $y \in S^1$, note that because $g^2 \in H$, $|r_g| = r_g \overline{r_g}$ is an $e$th root of $1$. We also have that $r_g \overline{r_g} \in \bR_{> 0}$, and therefore $|r_g| = 1$.
\end{proof}

We second case when the conclusion of Theorem \ref{thm:antilinearcyclotomicfinal} is immediate is when $W$ has S-type $\bC$ or $\bH$.

\begin{lemma}\label{lem:onlycareabouttypeR}
    Suppose $W \in \Irr_{\bC}^{\hspace{0.1em} \star}(G)$ has S-type $\bC$ or $\bH$. Then there is $\tilde{W} \in \Irr_{E}^{\hspace{0.1em} \star}(H)$ with $\tilde{W}_{\bC} \cong W$. 
\end{lemma}

\begin{proof}
    If $V \in \Irr_{\bC}(H)$ corresponds to $W$ and has S-type $\bC$ or $\bH$, then $\chi_W = \chi_V + {}^{\omega}\overline{\chi_V}$. Setting $\tilde{W} \coloneqq \Ind_H^G(\tilde{V})$, where $\tilde{V}$ is a model of $V$ over $E$, we have that 
    \[
    \chi_{\tilde{W}} = \chi_{\tilde{V}} + {}^{\omega}\overline{\chi_{\tilde{V}}} = \chi_V + {}^{\omega}\overline{\chi_V}
    \]
    hence $\tilde{W}_{\bC} \cong W$ as they have the same character.
\end{proof}

Therefore, to prove Theorem \ref{thm:antilinearcyclotomicfinal}, one must only consider the case when $W$ has D-type I, V, or IX. To understand the potential obstructions in these cases, let $V \in \Irr_{\bC}(G)$ correspond to $W$, and let $\tilde{V}$ be a model of $V$ over $E$. If $\tilde{W} \in \Irr_{E}^{\hspace{0.1em} \star}(G)$ corresponds to $\tilde{V}$, then we either have that $\tilde{W}_{\bC} \cong W$ (when $\tilde{W}$ has S-type $\bR$) or $\tilde{W}_{\bC} \cong W \oplus W$ (when $\tilde{W}$ has S-type $\bH$). To prove Theorem \ref{thm:antilinearcyclotomicfinal} is to prove that this second case cannot occur.

The classification into 24 ways alone is not enough to deduce this. For example, a priori, the pair $(V,W)$ could have D-type I, but $(\tilde{V}, \tilde{W})$ have type (9). Similarly, $(V,W)$ could have D-type V, and $(\tilde{V}, \tilde{W})$ have type (15), or $(V,W)$ could have D-type IX, whilst $(\tilde{V}, \tilde{W})$ have type (19).

Instead, in order to prove Theorem \ref{thm:antilinearcyclotomicfinal}, we will follow the general strategy of the classical proofs for cyclotomic realisability over $\bC$ and $\bR$, and deduce this as a consequence of a semilinear induction theorem for antilinear representations.

\subsection{Fundamental Representations}\label{sect:fundreps}

There are two natural operations one can consider on antilinear representations. The first is inflation, and the second is induction from a subgroup.

In this section we consider families of Real groups which (in a precise sense) generate all semilinear representations of all Real groups using these two operations. Throughout this section, $\bL / \bK$ is a characteristic $0$ degree two Galois extension, which we will later take to be $\bC / \bR$.

It turns out that it is most natural to work with the category of $C_2$-graded groups, which we consider as pairs $(G, \mu_G \colon G \rightarrow C_2)$, where a morphism $(G, \mu_G) \rightarrow (H,\mu_H)$ is a group homomorphism $f \colon G \rightarrow H$ such that $\mu_H \circ f = \mu_G$. The category of Real groups is recovered as the full subcategory which has objects those $(G,\mu_G)$ such that $\mu_G \colon G \rightarrow C_2$ is surjective.

A $C_2$-graded subgroup (resp.\ quotient) is simply a morphism of $C_2$-graded groups which is injective (resp.\ surjective). For a subgroup $H$ of a $C_2$-graded group $G$, $H$ has a canonical $C_2$-graded group structure, which is the unique structure making the inclusion a morphism of $C_2$-graded groups. For a quotient $f \colon G \twoheadrightarrow Q$ of $C_2$-graded group, $Q$ has at most one structure as a $C_2$-graded group making $f \colon G \twoheadrightarrow Q$ a morphism of $C_2$-graded groups, which exists if and only if $\ker(f) \subset \ker(\mu_G)$.

For a $C_2$-graded group $G$, we write $\bL_{\star \hspace{-0.1em}}[G]$ for the skew group algebra of $G$, which has our standard meaning when $\mu_G \colon C \rightarrow C_2 \equiv \Gal(\bL / \bK)$ is surjective, but which we also use to denote the ordinary group algebra $\bL[G]$ when $\mu_G$ is trivial. Similarly, we extend the notation $\Rep_{\bbL}^{\hspace{0.1em} \star}(G)$ and $\Irr_{\bbL}^{\hspace{0.1em} \star}(G)$ to the case where $\mu_G$ is potentially trivial. With this definition, $\bL_{\star \hspace{-0.1em}}[-]$ can be viewed as a functor from $C_2$-graded groups to $\bK$-algebras. In particular, for a $C_2$-graded subgroup $H \hookrightarrow G$ we can define $C_2$-graded induction by
\[
\Ind_S^G(-) \coloneqq \bL_{\star \hspace{-0.1em}}[G] \otimes_{\bbL_{\star \hspace{-0.1em}}[H]} -,
\]
and for a $C_2$-graded quotient $f \colon G \twoheadrightarrow Q$, we can define $C_2$-graded inflation $\Inf_Q^G(-)$ as restriction along the natural map $\bL_{\star \hspace{-0.1em}}[G] \rightarrow \bL_{\star \hspace{-0.1em}}[Q]$.

\begin{defn}
    A $C_2$-graded subquotient of $(G, \mu_G)$ is a $C_2$-graded quotient of a $C_2$-graded subgroup of $G$.
\end{defn}

The following shows that the $C_2$-graded subquotient relation is transitive and establishes a property of representations inflated and induced from $C_2$-graded subquotients, which we will use later.

\begin{prop}\label{prop:subquotofsubqout}
    Suppose that we have a diagram of $C_2$-graded groups:
\[\begin{tikzcd}
	G \\
	H & Q \\
	& S & P
	\arrow[hook, from=2-1, to=1-1]
	\arrow["f", two heads, from=2-1, to=2-2]
	\arrow[hook, from=3-2, to=2-2]
	\arrow[two heads, from=3-2, to=3-3]
\end{tikzcd}\]
Then viewing $f^{-1}(S)$ as a $C_2$-graded group using the inclusion $f^{-1}(S) \hookrightarrow G$, the map $f^{-1}(S) \twoheadrightarrow P$ is a $C_2$-graded quotient, and there is a natural isomorphism
\[
\Inf_Q^H(\Ind_S^Q(\Inf_P^S(V))) \cong \Ind^H_{f^{-1}(S)}(\Inf_P^{f^{-1}(S)}(V))
\]
for any $V \in \Rep_{\bbL}^{\hspace{0.1em} \star}(P)$.
\end{prop}

\begin{proof}
    The first claim follows from the commutativity of the diagram
\[\begin{tikzcd}
	G && {C_2} \\
	H & Q \\
	{f^{-1}(S)} & S & P
	\arrow["{\mu_G}", curve={height=-6pt}, from=1-1, to=1-3]
	\arrow[hook, from=2-1, to=1-1]
	\arrow["{\mu_H}", curve={height=-6pt}, from=2-1, to=1-3]
	\arrow["f", two heads, from=2-1, to=2-2]
	\arrow["{\mu_Q}", from=2-2, to=1-3]
	\arrow[hook, from=3-1, to=2-1]
	\arrow["f", two heads, from=3-1, to=3-2]
	\arrow["{\mu_S}"', curve={height=6pt}, from=3-2, to=1-3]
	\arrow[hook, from=3-2, to=2-2]
	\arrow[two heads, from=3-2, to=3-3]
	\arrow["{\mu_P}"', curve={height=6pt}, from=3-3, to=1-3]
\end{tikzcd}\]

For the second, given $V \in \Rep_{\bbL}^{\hspace{0.1em} \star}(P)$, we define
\[
\phi \colon \bL_{\star \hspace{-0.1em}}[H] \otimes_{\bbL_{\star \hspace{-0.1em}}[f^{-1}(S)]} \Inf_P^{f^{-1}(S)}(V) \rightarrow \Inf_Q^H(\bL_{\star \hspace{-0.1em}}[Q] \otimes_{\bbL_{\star \hspace{-0.1em}}[S]} \Inf^S_P(V)), \qquad \phi(h \otimes v) \coloneqq f(h) \otimes v,
\]
which is well-defined and surjective. As both sides have the same $\bK$-dimension $\dim_{\bbK}(V) \cdot [H:f^{-1}(S)] = \dim_{\bbK}(V) \cdot [Q:S]$, $\phi$ is an isomorphism.
\end{proof}
\begin{defn}
We write $\Char^{\hspace{0.1em} \star}(G,\bL)$ for the collection of characters of representations of the $C_2$-graded group $G$ over $\bL$.
\end{defn}

\begin{defn}
For a list $\sL$ of isomorphism classes of $C_2$-graded groups, we say that $\sL$ is \emph{generating} if $\sL$ \emph{satisfies condition} $(G)$: for any $C_2$-graded group $G$ and $\chi \in \Char^{\hspace{0.1em} \star}(G, \bL)$ one can write 
\[
\chi = \sum_{i} n_i \cdot \Ind_{G_i}^G(\Inf_{Q_i}^{G_i}(\psi_i)),
\]
where $n_i \in \bZ$, $G_i$ are $C_2$-graded subgroups of $G$, $G_i \twoheadrightarrow Q_i$ are $C_2$-graded quotients, $\psi_i \in \Char^{\hspace{0.1em} \star}(Q_i, \bL)$, and each $Q_i \in \sL$.

We say that such a list $\sL$ is \emph{minimal} if after removing any element of $\sL$, $\sL$ no longer satisfies $(G)$.
\end{defn}

\begin{defn}\label{def:orderC2gradedgroup}
We define the \emph{order} $o(G)$ of a $C_2$-graded group $G$ by 
\[
o(G) \coloneqq 2 |\ker(\mu_G)| + |\mu_G(G)| - 2 \in \bZ_{\geq 1}.
\]
\end{defn}

Concretely, $o(G) = 2|G| - 1$ if $G$ is ungraded, and $o(G) = |G|$ if $G$ is non-trivially graded. With this definition, the order of a $C_2$-graded subgroup (resp.\ $C_2$-graded quotient) of $G$ is always at most the order of $G$, with equality if and only if the inclusion map (resp.\ the quotient map) is an isomorphism. For example, when $G$ is non-trivially graded the proper inclusion $\ker(\mu_G) \hookrightarrow G$ is between groups of order $|G| - 1$ and $|G|$ respectively.

Given a list of isomorphism classes of $C_2$-graded groups $\sL$ and $n \geq 1$, we set $\sL_{\leq n}$ for the subset of $C_2$-graded groups of order $\leq n$, and $\sL_n$ for the subset of $C_2$-graded groups of order exactly $n$.

\begin{prop}\label{prop:minimallistexists}
    A minimal list $\sG$ satisfying \emph{$(G)$} exists, and can be described inductively:
    \[
    \sG_1 = \{(C_1, C_1 \hookrightarrow C_2)\},
    \]
    and for any $n \geq 1$, $\sG_{n+1}$ is defined as the set of $C_2$-graded groups of order $n+1$ for which there is some $\chi \in \Char^{\hspace{0.1em} \star}(G, \bL)$ such that you cannot write
    \[
    \chi = \sum_{i} n_i \cdot \Ind_{G_i}^G(\Inf_{Q_i}^{G_i}(\psi_i)),
    \]
    where $n_i \in \bZ$, $G_i$ are $C_2$-graded subgroups of $G$, $G_i \twoheadrightarrow Q_i$ are $C_2$-graded quotients, $\psi_i \in \Char^{\hspace{0.1em} \star}(Q_i, \bL)$, and each $Q_i \in \sG_{\leq n}$. 
\end{prop}

\begin{proof}
    It is clear that $\sG$, as defined inductively above, satisfies $(G)$. To see that it is minimal, let $G \in \sG_{n+1}$, and suppose for a contradiction that $\sG \setminus \{G\}$ satisfied $(G)$. In particular, for any character $\chi \in \Char^{\hspace{0.1em} \star}(G, \bL)$ one can write
\[
\chi = \sum_{i} n_i \cdot \Ind_{G_i}^G(\Inf_{Q_i}^{G_i}(\psi_i)),
\]
where $n_i \in \bZ$, $G_i$ are $C_2$-graded subgroups of $G$, $G_i \twoheadrightarrow Q_i$ are $C_2$-graded quotients, $\psi_i \in \Char^{\hspace{0.1em} \star}(Q_i, \bL)$, and each $Q_i \in \sG \setminus \{G\}$. Any proper $C_2$-graded subquotient of $G$ has order $\leq n$, and thus each $Q_i \in \sG_{\leq n}$. However this contradictions the fact that $G \in \sG_{n+1}$.
\end{proof}

\begin{prop}\label{prop:minimallistexistsuniquely}
    A minimal list $\sG$ satisfying \emph{$(G)$} exists and is unique. The list $\sG$ can be explicitly described as the set of $C_2$-graded groups $G$ for which there is some $\chi \in \Char^{\hspace{0.1em} \star}(G, \bL)$ such that you cannot write
    \[
    \chi = \sum_{i} n_i \cdot \Ind_{G_i}^G(\Inf_{Q_i}^{G_i}(\psi_i)),
    \]
    for $n_i \in \bZ$, $G_i$ are $C_2$-graded subgroups of $G$, $G_i \twoheadrightarrow Q_i$ are $C_2$-graded quotients, $\psi_i \in \Char^{\hspace{0.1em} \star}(Q_i, \bL)$, and each $Q_i$ is a proper subquotient $G$.

    Furthermore, any list $\sL$ of Real groups satisfying \emph{$(G)$} contains $\sG$.
\end{prop}

\begin{remark}
    The inductive definition of Proposition \ref{prop:minimallistexistsuniquely} is computationally useful if one wants to describe $\sG$ by an inductive search, as one does not need to check condition $(G)$ over all subquotients, but only over those already in your list. The explicit definition of $\sG$ is theoretically useful, as it allows you to determine if certain infinite families are in $\sG$, without a complete knowledge of $\sG$.
\end{remark}

\begin{proof}
    Let $\sH$ denote the list of $C_2$-graded groups defined via the explicit description above. It is clear that for any list $\sL$ of $C_2$-graded groups satisfying $(G)$, $\sH \subset \sL$. We will show that $\sH$ satisfies $(G)$. From this it is immediate that any minimal list $\sL$ satisfying $(G)$ has to be equal to $\sH$, so in particular $\sG = \sH$, and that any list of $C_2$-graded groups satisfying $(G)$ contains $\sG$.

    We prove that $\sH$ satisfies $(G)$ by induction on the order of the $C_2$-graded group $G$ in condition $(G)$. If $o(G) = 1$, then by definition $G \in \sH$, and so $(G)$ is satisfied. Suppose now that $n \geq 1$ and $\sH$ satisfies $(G)$ for all groups of order $\leq n$, and take $G$ a $C_2$-graded group with $o(G) = n + 1$. If $G \in \sH$ then $(G)$ is trivially satisfied, so suppose that $G \not\in \sH$, meaning that for all $\chi \in \Char^{\hspace{0.1em} \star}(G, \bL)$
\[
\chi = \sum_{i} n_i \cdot \Ind_{G_i}^G(\Inf_{Q_i}^{G_i}(\psi_i)),
\]
where $n_i \in \bZ$, $G_i$ are $C_2$-graded subgroups of $G$, $G_i \twoheadrightarrow Q_i$ are $C_2$-graded quotients, $\psi_i \in \Char^{\hspace{0.1em} \star}(Q_i, \bL)$, and each $Q_i$ is a proper subquotient of $G$. For each $i$, because $Q_i$ is a proper subquotient, $o(Q_i) \leq n$, and therefore by our induction hypothesis we can write 
\[
\psi_i = \sum_{j} m_j \Ind_{S_j}^{Q_i}(\Inf_{P_j}^{S_j}(\mu_j)),
\]
where each $P_j \in \sH$. Using Proposition \ref{prop:subquotofsubqout}, if $f \colon H_i \twoheadrightarrow Q_i$ is the quotient map, then $f^{-1}(S_j)$ is a $C_2$-graded subgroup of $G$, $P_j$ is a $C_2$-graded quotient of $f^{-1}(S_j)$, and
\[
\Inf_{Q_i}^{H_i}(\Ind_{S_j}^{Q_i}(\Inf_{P_j}^{S_j}(\mu_j))) \cong \Ind^{H_i}_{f^{-1}(S_j)}(\Inf_{P_j}^{f^{-1}(S_j)}(\mu_j)).
\]
In particular, 
\begin{align*}
\Ind_{G_i}^G(\Inf_{Q_i}^{G_i}(\psi_i)) 
&= \sum_{j} m_j \Ind_{G_i}^G(\Inf_{Q_i}^{G_i}(\Ind_{S_j}^{Q_i}(\Inf_{P_j}^{S_j}(\mu_j)))), \\
&= \sum_{j} m_j \Ind_{G_i}^G(\Ind^{G_i}_{f^{-1}(S_j)}(\Inf_{P_j}^{f^{-1}(S_j)}(\mu_j))), \\
&= \sum_{j} m_j \Ind^{G}_{f^{-1}(S_j)}(\Inf_{P_j}^{f^{-1}(S_j)}(\mu_j)).
\end{align*}
Replacing each $\Ind_{G_i}^G(\Inf_{Q_i}^{G_i}(\psi_i))$ in this way, we can write  
\[
\chi = \sum_{i} a_i \cdot \Ind_{G_i'}^G(\Inf_{Q_i'}^{G_i'}(\psi_i')),
\]
where $a_i \in \bZ$, $G_i'$ are $C_2$-graded subgroups of $G$, $G_i' \twoheadrightarrow Q_i'$ are $C_2$-graded quotients, $\psi_i' \in \Char^{\hspace{0.1em} \star}(Q_i', \bL)$, and each $Q_i' \in \sH_{\leq n}$. In particular, $\sH$ satisfies $(G)$ for $C_2$-graded groups of order $n+1$, and we are done by induction.
\end{proof}

Therefore, we see that for any characteristic $0$ degree two Galois extension $\bL / \bK$, there is a unique minimal list $\sG_{\bL / \bK}$ satisfying $(G)$. We refer to this $\sG_{\bL / \bK}$ as the list of \emph{fundamental representations}. For any given $\bL / \bK$, the main question is whether or not $\sG_{\bL / \bK}$ admits a reasonable description.

\subsection{Dihedral Representations}

The most important groups for our description of $\sG_{\bC / \bR}$ will be dihedral representations, that is, representations of the Real structure $C_n \leq D_{2n}$. Here we consider $C_n = \langle x \rangle \subset D_{2n} = \langle x ,y \mid x^n, y^2, yxy^{-1} = x^{-1} \rangle$, which makes sense whenever $n \geq 1$ (and so for us $D_2 \coloneqq C_2$ and $D_4 \coloneqq C_2 \times C_2$).

Let us first see that the antilinear representations of $C_n \leq D_{2n}$ are always defined over $\bQ(\zeta_n)$.

\begin{lemma}\label{lem:dihedralrealisable}
    Let $n \geq 1$. Then every irreducible representation of $C_n$ over $\bQ(\zeta_n)$ extends to an antilinear representation of the Real structure $C_n \leq D_{2n}$ over $\bQ(\zeta_n)$.
\end{lemma}

\begin{proof}
    Writing $D_{2n} = \langle x, y \mid x^n, y^2, yxy^{-1} = a \rangle$, if $x$ acts on $F$ by multiplication by the $n$th root of unity $\zeta$, then one can let $y$ act on $F$ by complex conjugation.
\end{proof}

Conversely, the conceptual reason why dihedral representations are important for the antilinear representation theory of all Real groups is the following.

\begin{lemma}\label{lem:1dimfftypeRisdihedral}
    Suppose that $H \leq G$ is a Real group with a $1$-dimensional faithful antilinear representation. Then $(H \leq G) \cong (C_n \leq D_{2n})$.
\end{lemma}

In the proof, we fix a basis an view antilinear representations as morphisms of Real groups $G \rightarrow \GL_n^{*}(\bC)$ (as in \cite[\S 2.2]{RUMTAY}), where
\[
\GL_n^*(\bC) = \{A_1 \mid A \in \GL_n(\bC)\} \sqcup \{A_{-1} \mid A \in \GL_n(\bC)\}
\]
with multiplication 
\[
A_{\epsilon}B_{\delta} = (A \cdot {}^{\epsilon}B)_{\epsilon\delta} = \twopartdef{(A B)_{\epsilon\delta}}{\epsilon = 1}{(A \overline{B})_{\epsilon\delta}}{\epsilon = -1.}
\]
\begin{proof}
    Let $\rho \colon G \rightarrow \GL_1^{*}(\bC)$ be a $1$-dimensional faithful antilinear representation of $G$. Then $\rho|_H \colon H \rightarrow \GL_1^1(\bC)$ is faithful, which means that $H$ is cyclic, so $H = \langle x \rangle$ for some $x \in H$. Let $y \in G$ be any odd element, and write $\rho(x) = [a]_1 \in \GL_1(\bC)$ and $\rho(y) = [b]_{-1} \in \GL_1^{-1}(\bC)$. 
    
    If $H$ has order $n$, then $a^n = 1$ in $\bC^{\times}$. Because $yxy^{-1} = x^i$ for some $i \in \bZ / n \bZ$, $yx = x^{i}y$, hence
    \[
    [b a^{-1}]_{-1} = [b \overline{a}]_{-1} = [b]_{-1} [a]_1 = \rho(yx) = \rho(x^{i}y) = [a^i]_1[b]_{-1} = [a^i b]_{-1},
    \]
    and so $a^{i} = a^{-1}$. In particular, $\rho(x^i) = \rho(x^{-1})$, hence by faithfulness $x^i = x^{-1}$ and so $i = -1$.
    
    We also have that $y^2 \in H$, hence $y^2 = x^j$ for some $j \in \bZ / n \bZ$, and
    \[
    [b \overline{b}]_1 = [b]_{-1}[b]_{-1} = \rho(y^2) = \rho(x^j) = [a^j]_1.
    \]
    so $b \overline{b} = a^j$. Now $a^j$ is a root of $1$, and $b \overline{b} \in \bR_{> 0}$, hence $\rho(y^2) = [1]_1$, and thus by faithfulness $y^2 = 1$.
    
    In summary, we have that $G = \langle x, y \rangle$, $G$ has order $2n$, and $G$ satisfies the relations $x^n = 1$, $y^2 = 1$, $yxy^{-1} = x^{-1}$, hence $G \cong D_{2n}$ as Real groups.
\end{proof}

\subsection{A Semilinear Induction Theorem}

The aim of this section is to describe the list $\sG_{\bC / \bR}$ of fundamental representations, in the special case when $\bL / \bK$ is $\bC / \bR$. Because induction and inflation both preserve cyclotomic realisability, this will allow us to prove cyclotomic realisability for all antilinear representations as an immediate consequence, in Section \ref{sect:proofofcyclorealise} below.

\begin{thm}\label{thm:fundamentalrepsexplicit}
    The list $\sG_{\bC / \bR} = \sG_{\bC / \bR}^{\, 1} \sqcup \sG_{\bC / \bR}^{\,-1}$ consists of the ungraded groups
    \[
       \sG_{\bC / \bR}^{\, 1} =  \{C_n \mid n \geq 1, n \neq 2 \!\!\! \mod 4\}
    \]
    and the Real groups
    \[
       \sG_{\bC / \bR}^{\, -1} = \{(C_n \leq D_{2n}) \mid n \geq 1, n \neq 2 \!\!\! \mod 4 \}.
    \]
\end{thm}

In order to prove Theorem \ref{thm:fundamentalrepsexplicit}, we need to choose a class of Real groups which is big enough to generate all antilinear representations of all Real groups under induction, but small enough that all representations are inflated from dihedral representations.

\begin{defn}\label{def:realelementary}
    A \emph{Real-elementary} group is a Real group of the form $C \rtimes P_0 \leq C \rtimes P$, where
    \begin{itemize}
        \item $C$ is cyclic of odd order,
        \item $P_0 \leq P$ is a Real $2$-group,
        \item $P_0$ acts on $C$ trivially, and $P \setminus P_0$ acts on $C$ by $(-)^{-1}$.
    \end{itemize}
\end{defn}

\begin{prop}\label{prop:realelementarygrpsaremonomial}
    Suppose that $C \rtimes P$ is a Real-elementary group. Then every irreducible antilinear representation of $C \rtimes P$ is induced from a $1$-dimensional representation of a subgroup.
\end{prop}

\begin{proof}
    Let $W$ be an irreducible antilinear representation of $C \rtimes P$. If $W$ has S-type $\bC$ or $\bH$, then $W \cong \Ind_{C \times P_0}^{C \rtimes P}(V)$, where $V \in \Irr_{\bC}(C \times P_0)$ corresponds to $W$, from which the claim follows as $C \times P_0$ is supersolvable, hence monomial by \cite[\S 8.5 Thm.\ 16]{SERRE}.

    Suppose instead that $W$ has S-type $\bR$, so $V \coloneqq W|_{C \times P_0}$ is irreducible. We adapt the proof of \cite[\S 8.5 Thm.\ 16]{SERRE}, and first show the claim when $C$ is trivial. We proceed by induction on the order of $P_0$. If $P_0$ is trivial then $P_0 \leq P$ is the Real group $C_1 \leq C_2$ and there is nothing to prove. Suppose instead that $P_0$ is non-trivial, so by induction we may assume that $W$, or equivalently $V$, is faithful.
    
    If $P_0$ is abelian, then because $V$ is irreducible, $W$ has dimension $1$, and we're done. If instead $P_0$ is not abelian, then $1 \neq P_0 / Z(P_0)$, and as $Z(P_0)$ is characteristic, $Z(P_0) \leq P$ is normal, so we can view $P_0 / Z(P_0) \leq P / Z(P_0)$. Because $P / Z(P_0)$ is a $2$-group, the normal subgroup $P_0 / Z(P_0)$ intersects $Z(P / Z(P_0))$ non-trivially, and so taking any element of order $2$ we can find a cyclic order $2$ normal subgroup $Q \triangleleft P / Z(P_0)$ which is contained in $P_0 / Z(P_0)$. Let $A$ denote the preimage of $Q$ in $P_0$. Because $A / Z(P_0)$ is cyclic, $A$ is abelian, and because $Q$ is normal in $P / Z(P_0)$, $A$ is normal in $P$.

    In summary, we have $A \leq P_0$, which is abelian, normal in $P$, and properly contains $Z(P_0)$. Consider $W|_A$. The representation $W|_{P_0}$ is faithful, so $W|_A$ cannot be isotypic: if it were then each $a \in A$ would have $\rho_V(a)$ a scalar matrix, and so commute with each $\rho_V(x)$ for $x \in P_0$, and so we would have $A \subset Z(P_0)$ by faithfulness, a contradiction. Therefore, the canonical decomposition $W|_A = V_1 \oplus \cdots \oplus V_k$ into isotypic factors necessarily has $k > 1$.
    
    Because $A$ is normal in $P$, there is a group action of $P$ on the set of isomorphism classes of representations of $A$: for $U \in \Rep_{\bC}(A)$ and $z \in P$, $z * U$ is $U$ as a $\bR$-vector space, with action of $a \in A$ by $z^{-1} a z$ and complex structure $\lambda * u \coloneqq \lambda u$ if $z \in P_0$ and $\lambda * u = \overline{\lambda} u$ is $z \in P \setminus P_0$. The action of each $z \in P$ defines an isomorphism $\rho_V(z) \colon z * W|_A \xrightarrow{\sim} W|_A$, which induces a permutation action of $P$ on the set $\{V_1, ... , V_k\}$. Let $S \coloneqq \Stab_P(V_1)$, noting that $S < P$ because $k > 1$ and the action is transitive ($W$ being irreducible). The natural map
    \[
        \bC_{\star \hspace{-0.1em}}[G] \otimes_{\bC_{\star \hspace{-0.1em}}[S]} V_1 \rightarrow W
    \]
    is surjective, the action being transitive, and therefore an isomorphism, both sides having the same dimension. In particular, $W \cong \Ind_S^G(V_1)$, and we are done by induction.

    This proves the claim when $C = 1$. In general, if $W$ has S-type $\bR$, then $W|_{C \times P_0}$ is irreducible, and because $C$ is abelian $C$ acts on $W|_{C \rtimes P_0}$ through some character $\chi \colon C \rightarrow \bC^{\times}$. In particular, $W|_{P_0}$ is still irreducible, and also of S-type $\bR$ (admitting the extension $W|_P$ to an antilinear representation). We can therefore write $W|_{P} = \Ind_S^P(\chi)$, for some Real subgroup $S \leq P$ and $1$-dimensional $U \in \Irr_{\bC}^{\hspace{-0.1em} \star }(S)$.
    
    Define an antilinear representation $\chi \otimes U$ of $C \rtimes S$, which is $U$ as an antilinear $S$-representation, with action of $C$ by scalar multiplication by $\chi$. The is well-defined: for $y \in C$ and $z \in S$, then if $z \in S_0$, $zy = y z$ in $C \rtimes S$ and indeed
    \[
        \rho_{\chi \otimes U}(z) \rho_{\chi \otimes U}(y)u
        = \rho_U(z)(\chi(y)(u)) = \chi(y) \rho_U(z)(u) =  \rho_{\chi \otimes U}(y) \rho_{\chi \otimes U}(z)u.
    \] 
    If on the other hand $z \in S \setminus S_0$, then $zy = y^{-1}z$ in $C \rtimes S$, and indeed
    \[
    \rho_{\chi \otimes U}(z) \rho_{\chi \otimes U}(y)u
        = \rho_U(z)(\chi(y)(u)) = \overline{\chi(y)} \rho_U(z)(u) =  \rho_{\chi \otimes U}(y^{-1}) \rho_{\chi \otimes U}(z)u.
    \]
    Then the fact that $W \cong \Ind_{C \rtimes S}^{C \rtimes P}(\chi \otimes U)$ follows from comparing characters.
\end{proof}

\begin{remark}
    Proposition \ref{prop:realelementarygrpsaremonomial} says that Real-elementary groups are \emph{Real-monomial}. Classically, supersolvable groups are monomial, and any monomial group is solvable, so monomial groups are somewhere between these two notions. In the antilinear setting, from computations it seems that a Real group $H \leq G$ is Real-monomial whenever $H$ is supersolvable, but we were unable to prove this.
\end{remark}

The second main result we need is an antilinear induction theorem for Real-elementary subgroups.

\begin{thm}\label{thm:antilinearinductionthm}
For any $C_2$-graded group $G$ and $\chi \in \Char^{\hspace{0.1em} \star}(G, \bC)$, there are $C_2$-graded subgroups $G_i$ of $G$, $n_i \in \bZ$, and $\chi_i \in \Char^{\hspace{0.1em} \star}(G_i, \bL)$ such that each $G_i$ is either ungraded or Real-elementary and
\[
    \chi = \sum_{i} n_i \cdot \Ind_{G_i}^G(\chi_i).
\]
\end{thm}

\subsection{Proof Of Theorem \ref{thm:antilinearinductionthm}}

In this section we prove Theorem \ref{thm:antilinearinductionthm}. The claim is trivial if $G$ is ungraded, so suppose that $G$ is non-trivially graded. Let $\sR^{\hspace{0.1em} \star}_{\bZ}(G)$ denote the Grothendieck group of antilinear representations of $G$. Let $I_1$ denote the image of the $\bZ$-linear map
\[
    \Ind_{\sI_1} \coloneqq \bigoplus_{H \leq \sI_1} \Ind_H^G \colon \bigoplus_{H \in \sI_1} \sR_{\bZ}(H) \rightarrow \sR^{\hspace{0.1em} \star}_{\bZ}(G)
\]
where $\sI_1$ is the collection of ungraded subgroups of $G$, and let $I_2$ denote the image of the $\bZ$-linear map
\[
    \Ind_{\sI_2} \coloneqq \bigoplus_{S \in \sI_2} \Ind_S^G \colon \bigoplus_{S \in \sI_2} \sR_{\bZ}^{\hspace{0.1em} \star}(S) \rightarrow \sR^{\hspace{0.1em} \star}_{\bZ}(G),
\]
where $\sI_2$ is the collection of Real-elementary subgroups of $G$. From the correspondence between irreducible representations of $G_0$ and antilinear representations of $G$ (Theorem \ref{thm:threeoptions}), we have that $I_1$ has finite index in $\sR^{\hspace{0.1em} \star}_{\bZ}(G)$ equal to $2^k$, where $k$ is the number of irreducible antilinear representations of $G$ of S-type $\bR$. Therefore, it is sufficient to show that $I \coloneqq I_1 + I_2$ has odd index in $\sR^{\hspace{0.1em} \star}_{\bZ}(G)$, for then $I$ must have both odd and even index, and thus $I = \sR^{\hspace{0.1em} \star}_{\bZ}(G)$. Because $I$ of $\sR^{\hspace{0.1em} \star}_{\bZ}(G)$ is an ideal by the projection formula, Lemma \ref{lem:projectionformula}, we are reduced to the following.

\begin{prop}\label{prop:imageoddindex}
If $|G| = 2^r m$ with $m$ odd, then $m \in I$.
\end{prop}

To prove Proposition \ref{prop:imageoddindex} we adapt the proof of \cite[\S 12.7]{SERRE} to our setting. Set $n \coloneqq |G_0|$, so $n = 2^{r-1} m$, and let $A$ denote the subring of $\overline{\bZ}$ generated by the $n$th roots of $1$. There is an induction map
\[
A \otimes_{\bZ} \Ind_{\sI} \colon \bigoplus_{S \in \sI} A \otimes_{\bZ} \sR_{\bZ}^{\hspace{0.1em} \star}(S) \rightarrow A \otimes_{\bZ} \sR^{\hspace{0.1em} \star}_{\bZ}(G),
\]
and just as in \cite[\S 10.2, Lem.\ 5]{SERRE} we have the following, from which we see that to prove Proposition \ref{prop:imageoddindex} it is sufficient to show that $m \in A \otimes_{\bZ} I$.

\begin{lemma}\label{lem:cyclotomicintegers}
   The image of $A \otimes_{\bZ} \Ind_{\sI}$ is $A \otimes_{\bZ} I$, and $(A \otimes I) \cap \sR^{\hspace{0.1em} \star}_{\bZ}(G) = I$.
\end{lemma}

Let $G_0 /\!/ G_0$ denote the set of conjugacy classes of $G_0$, and let $G_0 /\!/ G$ denote the set of Real conjugacy classes \cite[\S 1.3]{RUMTAY} of $G_0$ with respect to the Real structure $G$. Because antilinear representations of $G$ are Real class function and are determined by their restriction to $G_0$ the injective ring homomorphism
\[
\Phi_0 \colon A \otimes_{\bZ} \sR_{\bZ}(G_0) \hookrightarrow \Fun(G_0 /\!/ G_0, A),
\]
to the ring of functions from $G_0 /\!/ G_0$ to $A$ restricts to an injective ring homomorphism
\[
\Phi \colon A \otimes_{\bZ} \sR^{\hspace{0.1em} \star}_{\bZ}(G) \hookrightarrow \Fun(G_0 /\!/ G, A).
\]
to the ring of functions from $G_0 /\!/ G$ to $A$. The first ingredient in the proof of Proposition \ref{prop:imageoddindex} is the following sufficient criterion for such a function be in $A \otimes_{\bZ} I$.

\begin{lemma}\label{lem:intvaluedconstantfunctions}
    Suppose that $\varphi \in \Fun(G_0 /\!/ G, A)$ has values in $|G| \cdot \bZ$. Then $\varphi \in \Phi(A \otimes_{\bZ} I_1)$. 
\end{lemma}

In the proof we use the following action of $C_2 = \langle \epsilon \rangle$ on $\Fun(G_0 /\!/ G, A)$. Fix $\omega \in G \setminus G_0$, and set 
\[
(\epsilon * \varphi)(x) \coloneqq \varphi(\omega x^{-1} \omega^{-1}).
\]
With this definition, from \cite[Prop.\ 2.8]{RUMTAY} we have that for any $\varphi_0 \in A \otimes_{\bZ} \sR^{\hspace{0.1em} \star}_{\bZ}(G)$,
\[
\Phi(\Ind_{G_0}^G(\varphi_0)) = \Phi_0(\varphi_0) + \epsilon * \Phi_{0}(\varphi_0).
\]
\begin{proof}
    Let $\varphi_0 \coloneqq \varphi_0 / 2$, which is a Real class function with values in $n \bZ$. By \cite[\S10.3, Lem.\ 6]{SERRE}, $\varphi_0 = \Phi_0(x_0)$ for some $x_0 \in A \otimes_{\bZ} \sR_{\bZ}(G_0)$. Set $x \coloneqq \Ind_{G_0}^G(x_0) \in A \otimes_{\bZ} I_1$. We have that
    \begin{align*}
        \Phi(x) = \Phi(\Ind_{G_0}^G(x_0)) =  \Phi_0(x_0) + \epsilon * \Phi_{0}(x_0) = \varphi_0 + \epsilon * \varphi_0 = \varphi,
    \end{align*}
    because $\varphi_0$ is a Real class function, and thus fixed by the $C_2$ action.
\end{proof}

We now use $I_2$ and Real-elementary subgroups, which we haven't needed until now.

\begin{defn}
For any $x \in G_0$ of odd order, set
\[
Z(x) \coloneqq \{g \in G_0 \mid gxg^{-1} = x\} \sqcup \{z \in G \setminus G_0 \mid zxz^{-1} = x^{-1}\},
\]
fix a $2$-Sylow subgroup $P$ of $Z(x) / \langle x \rangle$, and let $R(x)$ be the preimage of $P$ in $Z(x)$.
\end{defn}

Given $x \in G_0$ of odd order, either $R(x)$ is ungraded (contained in $G_0$) or $R(x)$ is Real-elementary (Definition \ref{def:realelementary}), because $R(x)$ is the semi-direct product of $\langle x \rangle$ and any $2$-Sylow subgroup of $R(x)$.

\begin{lemma}\label{lem:localspecialelt}
    Suppose that $x \in G_0$ has odd order. Then there is $\varphi_x \in A \otimes_{\bZ} \sR^{\hspace{0.1em} \star}_{\bZ}(R(x))$ such that
    \[
    \varphi \coloneqq \Ind_{R(x)}^G(\varphi_x) \in A \otimes_{\bZ} I
    \]
    satisfies the following:
    \begin{enumerate}
        \item $\Phi(\varphi)$ is $\bZ$-valued,
        \item $\Phi(\varphi)(x) \not\equiv 0 \!\! \mod 2$,
        \item $\Phi(\varphi)(y) = 0$ for any odd order $y \in G_0$ which is not Real conjugate to $x$.
    \end{enumerate}
\end{lemma}

\begin{proof}
    We have that $R(x)_0$ is a $2$-elementary subgroup of $G_0$, and thus by \cite[\S 10.3, Lem.\ 8]{SERRE} there is $\varphi_{x,0} \in A \otimes_{\bZ} \sR_{\bZ}(R(x)_0)$ such that, setting $\varphi_0 \coloneqq \Ind_{R(x)_0}^{G_0}(\varphi_{x,0})$, 
    \begin{enumerate}
        \item $\Phi(\varphi_0)$ is $\bZ$-valued,
        \item $\Phi(\varphi_0)(x) \not\equiv 0 \!\! \mod 2$,
        \item $\Phi(\varphi_0)(y) = 0$ for any odd order $y \in G_0$ which is not conjugate to $x$.
    \end{enumerate}
    
    There are now two cases. Suppose first that $R(x)$ is ungraded, so $R(x) \subset G_0$ and $R(x) = R(x)_0$. Set $\varphi_x \coloneqq \varphi_{x,0}$ and $\varphi \coloneqq \Ind_{R(x)}^G(\varphi_x) = \Ind_{G_0}^{G}(\varphi_0)$. We have that
    \[
        \Phi(\varphi) = \Phi_0(\varphi_0) + \epsilon * \Phi_{0}(\varphi_0),
    \]
    and so $\Phi(\varphi)$ is $\bZ$-valued. To verify the second condition, suppose that $y \in G_0$ has odd order. Then
    \[
        \Phi(\varphi)(y) = \Phi_0(\varphi_0)(y) + \Phi_{0}(\varphi_0)(\omega y^{-1} \omega^{-1}).
    \]
    When $y = x$, $\omega x^{-1} \omega^{-1}$ is not $G_0$-conjugate to $x$ because $Z(x) \subset G_0$, so $\Phi_{0}(\varphi_0)(\omega x^{-1} \omega^{-1}) \equiv 0 \!\! \mod 2$ and $\Phi(\varphi)(x) \equiv 0 \!\! \mod 2$. When $y$ is not Real conjugate to $x$, then both $y$ and $\omega y^{-1} \omega^{-1}$ are not $G_0$-conjugate to $x$, hence $\Phi(\varphi)(y) = 0$.

    Suppose instead that $R(x)$ is a Real subgroup of $G$, and choose a $2$-Sylow subgroup $P$ of $R(x)$ so that we can view $R(x) = C \rtimes P$ with $R(x)_0 = C \times P_0$, where $P_0$ is an index two subgroup of $P$. The element $\varphi_{x,0} \in A \otimes_{\bZ} \sR_{\bZ}(R(x)_0)$ is constructed as
    \[
        \varphi_{x,0} = \sum_{\chi} \chi(x^{-1}) \otimes (\chi \otimes \chi_1)
    \]
    as $\chi$ runs over the irreducible characters of $C$, and $\chi_1$ is the trivial character of $P_0$. Because $\chi_1$ extends to an antilinear representation of $P$, the following the final part of the proof of Proposition \ref{prop:realelementarygrpsaremonomial}, each $\chi \otimes \chi_1$ is the character of an object of $\sR^{\hspace{0.1em} \star}_{\bZ}(R(x))$ and we can view $\varphi_{x,0}$ as the restriction of $\varphi_x \in A \otimes_{\bZ} \sR^{\hspace{0.1em} \star}_{\bZ}(R(x))$. Setting $\varphi \coloneqq \Ind_{R(x)}^G(\varphi_x)$ we have by Lemma \ref{lem:inductioncompatibility} that $\Phi(\varphi) = \Phi_0(\varphi_0)$, and therefore $\Phi(\varphi)$ is $\bZ$-valued, $\Phi(\varphi)(x) \not\equiv 0 \!\! \mod 2$, and $\Phi(\varphi)(y) = 0$ for any odd order $y \in G_0$ which is not Real conjugate to $x$, as the fact that $Z(x) \not\subset G_0$ means that the Real conjugacy of $x$ is the same as the $G_0$-conjugacy class of $x$.
\end{proof}

\begin{cor}\label{cor:particularelt}
There exists $\varphi \in A \otimes_{\bZ} I$ such that
\begin{enumerate}
    \item $\Phi(\varphi)$ is $\bZ$-valued,
    \item $\Phi(\varphi)(x) \not\equiv 0 \!\! \mod 2$ for any $x \in G_0$.
\end{enumerate}
\end{cor}

\begin{proof}
    Let $x_1, ... , x_k$ denote representations for the Real conjugacy classes of $G_0$ which consist of elements of odd order. Let $\varphi_i \in A \otimes_{\bZ} I$ be the element provided by Lemma \ref{lem:localspecialelt}, and set $\varphi \coloneqq \sum_i \varphi_i$. Clearly $\Phi(\varphi)$ is $\bZ$-valued, this being true of each $\Phi(\varphi_i)$. For the second condition, let $x \in G_0$, and write $x'$ for the $2$-regular part of $x$ as in \cite[\S 10.1]{SERRE}. From \cite[\S 10.3, Lem.\ 7]{SERRE}, $
    \Phi(\varphi_i)(x) \equiv \Phi(\varphi_i)(x') \!\! \mod 2$ for any $i = 1, ... , k$. Furthermore, $x'$ is Real conjugate among the $x_1, ... , x_k$ to exactly one element $x_i$, and therefore
    \[
        \Phi(\varphi)(x) \equiv \sum_j \Phi(\varphi_j)(x') \equiv  \Phi(\varphi_i)(x_i) \not\equiv 0 \!\! \mod 2. \qedhere
    \]
\end{proof}

We are now in a position to prove Proposition \ref{prop:imageoddindex}.

\begin{proof}[Proof of Proposition \ref{prop:imageoddindex}]
    Let $\varphi \in A \otimes_{\bZ} I_2$ be as in Corollary \ref{cor:particularelt}, and set $N \coloneqq 2^r$, the order of $(\bZ / 2^{r+1} \bZ)^{\times}$. Consider the function $\phi \coloneqq m \Phi(\varphi^N - 1) \in \Fun(G_0 / \! / G, A)$. Because $\Phi(\varphi)$ is $\bZ$-valued, the same is true for $\phi$, and because $\Phi(\varphi)(x) \not\equiv 0 \!\! \mod 2$ for any $x \in G_0$, $\Phi(\varphi)(x)^N \equiv 1 \!\! \mod 2$, hence $\phi$ has values in $|G| \cdot \bZ$. By Lemma \ref{lem:intvaluedconstantfunctions} there is $\psi \in A \otimes_{\bZ} I_1$ with $\Phi(\psi) = \phi$, and so by the injectivity of $\Phi$, $\psi = m\varphi^N - m$. We have that $m \varphi^N \in A \otimes_{\bZ} I$ because $A \otimes_{\bZ} I$ is an ideal of $A \otimes_{\bZ} \sR^{\hspace{0.1em} \star}_{\bZ}(G)$ and $\varphi \in A \otimes_{\bZ} I$ by Corollary \ref{cor:particularelt}. It then follows that $m = \psi - m \varphi^N \in A \otimes_{\bZ} I$, from which the claim follows from Lemma \ref{lem:cyclotomicintegers}.
\end{proof}

\subsection{Proof of Theorem \ref{thm:fundamentalrepsexplicit}}

We are now in a position to prove Theorem \ref{thm:fundamentalrepsexplicit}. We first show that $\sG_{\bC / \bR} \subset \sH$, where $\sH \coloneqq \sH^{1} \sqcup \sH^{-1}$, for $\sH^1$ the ungraded list
\[
\sH^{1} = \{C_n \mid n \geq 1\},
\]
and $\sH^{-1}$ the non-trivially graded list
\[
\sH^{-1} = \{C_n \leq D_{2n} \mid n \geq 1\}.
\]
We then compute exactly $\sG_{\bC / \bR} \subset \sH$ to be the sublist as described by Theorem \ref{thm:fundamentalrepsexplicit}.

We show that $\sG_{\bC / \bR} \subset \sH$ by induction on the order (see Definition \ref{def:orderC2gradedgroup}) of the $C_2$-graded group $G$. When the $o(G) = 1$, the only option for $G$ is that $G = C_1$ with the trivial grading, and this lies in both $\sG_{\bC / \bR}$ and $\sH$. Suppose that $\sG_{\bC / \bR, \leq n} \subset \sH_{\leq n}$ for some $n \geq 1$, and take $G \in \sG_{\bC / \bR, n+1}$. By Proposition \ref{prop:minimallistexistsuniquely}, there is some $\chi \in \Char^{\hspace{0.1em} \star}(G, \bC)$ such that you cannot write
    \[
    \chi = \sum_{i} n_i \cdot \Ind_{G_i}^G(\Inf_{Q_i}^{G_i}(\psi_i)),
    \]
where $n_i \in \bZ$, $G_i$ are $C_2$-graded subgroups of $G$, $G_i \twoheadrightarrow Q_i$ are $C_2$-graded quotients, $\psi_i \in \Char^{\hspace{0.1em} \star}(Q_i, \bC)$, and each $o(Q_i) < o(G)$.

Suppose that $n$ is even, so $G$ is ungraded. By Brauer's induction theorem we can write $\chi$ as a linear combination of characters induced from one-dimensional representations of ungraded subgroups, each of which, being one-dimensional, is inflated from a representation of a cyclic group. In particular, if $G$ were not cyclic, then it would fail to satisfy the above property, and so we must have that $G \in \sH^{1}$.

Suppose instead that $n$ is odd, so $G$ is non-trivially graded. By Theorem \ref{thm:antilinearinductionthm}, we immediate see that $G$ is Real-elementary. Decomposing $\chi$ into irreducible components, we may also assume that $\chi$ is the character of an irreducible antilinear representation of $G$. Then by Proposition \ref{prop:realelementarygrpsaremonomial}, $\chi$ is one-dimensional. As $G$ satisfies the above property, $\chi$ must be faithful, and therefore $G \cong (C_{m} \leq D_{2m})$ by Lemma \ref{lem:1dimfftypeRisdihedral}, where $m \coloneqq (n+1)/2$. In particular, $G \in \sH^{-1}$, and thus by induction $\sG_{\bC / \bR} \subset \sH$.

Now we consider exactly which elements of $\sH$ are in $\sG_{\bC / \bR}$. If $G \in \sH^1$, so $G = C_n$ for some $n \geq 1$, then all subquotients of $G$ are also ungraded, and thus $G \in \sG_{\bC / \bR}$ if and only if $n \equiv 2 \!\! \mod 4$ by Theorem \ref{thm:wheniscyclicnotneeded}. Suppose instead $G \in \sH^{-1}$, so $G = D_{2n}$ for some $n \geq 1$. From Lemma \ref{lem:dihedralrealisable} the characters of antilinear representations of $G$ are the same as the characters of the ungraded subgroup $C_n$. The Real structure $C_n \leq D_{2n}$ has the property that every ungraded subgroup $C_m \leq C_n$ lies in the Real subgroup $(C_m \leq D_{2m}) \leq (C_n \leq D_{2n})$, and any Real quotient of $(C_n \leq D_{2n})$ is of the form $(C_n \leq D_{2n}) \twoheadrightarrow (C_m \leq D_{2m})$. In particular, by Lemma \ref{lem:inductioncompatibility} the $\bZ$-submodule of $\Char^{\hspace{0.1em} \star}(D_{2n}, \bC)$ generated by characters of proper $C_2$-graded subquotients under inflation and induction is the same as the $\bZ$-submodule of $\Char(C_{n}, \bC)$ generated by proper subquotients of $C_n$ under inflation and induction. We therefore have that the answer is the same as the ungraded $C_n$ case, and so $(C_n \leq D_{2n}) \in \sG_{\bC / \bR}$ if and only if $n \equiv 2 \!\! \mod 4$ by Theorem \ref{thm:wheniscyclicnotneeded}.

\subsection{Proof of Theorem \ref{thm:antilinearcyclotomicfinal}}\label{sect:proofofcyclorealise}

From Theorem \ref{thm:fundamentalrepsexplicit} and Lemma \ref{lem:dihedralrealisable} there are $X_1, ... , X_n \in \Rep_{E}^{\hspace{0.1em} \star}(G)$, with characters $\psi_1, ... ,\psi_n$ such that
\[
\chi_W = \sum_{i = 1}^n a_{i} \psi_i
\]
for some $a_i \in \bZ$. If $W$ has S-type $\bC$ or $\bH$ we are done by Lemma \ref{lem:onlycareabouttypeR}. Otherwise, $W$ has S-type $\bR$, and therefore $\langle \chi_W, \chi_W \rangle = 1$, so
\[
1 = \langle \chi_W, \chi_W \rangle = \sum_{i = 1}^n a_{i} \langle \chi_W, \psi_i \rangle.
\]

Let $V \in \Irr_{\bC}(H)$ correspond to $W$, and let $\tilde{V} \in \Irr_E(H)$ be a model of $V$ over $E$, using \cite[\S 12.3]{SERRE}. Let $\tilde{W} \in \Irr_{E}^{\hspace{0.1em} \star}(G)$ correspond to $\tilde{V}$, and note that $\chi_W = \chi_V = \chi_{\tilde{V}}$. We have that $m_{E_+}^E(\tilde{V}) \in \{1,2\}$, and $m_{E_+}^E(\tilde{V}) \mid \langle \chi_{\tilde{V}}, \chi_U \rangle$ for any $U \in \Rep_{E}^{\hspace{0.1em} \star}(G)$ by \cite[Prop.\ 9.3(4)]{TAYSchur}. In particular, if $m_{E_+}^E(\tilde{V}) = 2$, then each term $\langle \chi_W, \psi_i \rangle$ is divisible by two, which is impossible as the terms sum to $1$. Therefore, $m_{E_+}^E(\tilde{V}) = 1$, so $\chi_{\tilde{W}} = \chi_{\tilde{V}}$, and $\tilde{W}_{\bC} \cong W$, as their characters are the same. \qedsymbol

\subsection{Recovery of Serre's Induction Theorem over $\bR$}\label{sect:recoveryofclassical}

Let $G$ be a $C_2$-graded group with trivial grading. Then every $C_2$-graded subgroup and $C_2$-graded quotient also has the trivial grading, and we see that Theorem \ref{thm:fundamentalrepsexplicit} recovers Brauer's Induction Theorem \cite[\S 10.5]{SERRE}, that any complex representation of a finite group is a $\bZ$-linear combination of representations induced from one-dimensional representations of its subgroups.

We now show how Theorem \ref{thm:fundamentalrepsexplicit} also recovers Serre's Induction Theorem \cite[Lem.\ 73.19]{CR2}, for $\bR$-linear representations.

\begin{cor}\label{cor:recoverserre}
    Let $H$ be a finite group, and $\chi \in \Char(H, \bR)$. Then
    \[
    \chi = \sum_{i} n_i \Ind_{H_i}^H(\psi_i),
    \]
    where $n_i \in \bZ$, $H_i \leq H$, and $\psi_i$ is inflated from an $\bR$-representation of a cyclic or dihedral group. 
\end{cor}

Before giving the proof, we need the following result regarding induction from subgroups of $H \times C_2$ which are not the product of subgroups of $H$ and $C_2$.

Recall that by Goursat's Lemma, subgroups of $H \times C_2$ come in three flavours: 
\begin{enumerate}
    \item $K \times 1$, for $K \leq H$,
    \item $K \times C_2$, for $K \leq H$,
    \item $S(K,u) \coloneqq \{(k, u(k)) \mid k \in K\}$, where $K \leq H$ and $u \colon K \twoheadrightarrow C_2$.
\end{enumerate}

It is these subgroups of the third class that do most of the work when Theorem \ref{thm:antilinearinductionthm} is applied to split Real structures. For these subgroups, note that there is an canonical isomorphism $K \xrightarrow{\sim} S(K,u)$, and using this we have a restriction functor
\[
    (-)|_{\bR[K]} \colon \Rep_{\bC}^{\hspace{0.1em} \star}(S(K,u)) \rightarrow \Rep_{\bR}(K).
\]

\begin{lemma}\label{lem:twistedsubgroupinductionsplit}
    The diagram
\[\begin{tikzcd}[column sep=large]
	{\Rep_{\bC}^{\hspace{0.1em} \star}(S(K,u)) } & {\Rep_{\bC}^{\hspace{0.1em} \star}(H \times C_2)} \\
	{\Rep_{\bR}(K)} & {\Rep_{\bR}(H)}
	\arrow["{\Ind_{S(K,u)}^{H \times C_2}}", from=1-1, to=1-2]
	\arrow["{(-)|_{\bR[K]}}", from=1-1, to=2-1]
	\arrow["\sim", from=1-2, to=2-2]
	\arrow["{\Ind_K^H}"', from=2-1, to=2-2]
\end{tikzcd}\]
    commutes up to isomorphism.
\end{lemma}

\begin{remark}
    This isomorphism can probably be made natural, but the above is sufficient for our purposes.
\end{remark}

\begin{proof}
    Let $V \in \Rep_{\bC}^{\hspace{0.1em} \star}(S(K,u))$ and let $\chi$ denote the character of $V$. If $\psi$ denote the character of $V|_{\bR[K]}$, then by \cite[\S 2.6]{RUMTAY} $\psi(k) = \chi(k) + \overline{\chi}(k)$ when $k \in K_0$ and $\psi(k) = 0$ when $k \in K \setminus K_0$. We thus compute that $\Ind_K^H(\psi)$ has character with
    \begin{align*}
    |K| \cdot \Ind_K^H(\psi)(h) &= \sum_{\substack{x \in H \\ h \in xKx^{-1}}} \psi(x^{-1} h x), \\
    &= \sum_{\substack{x \in H \\ h \in xK_0x^{-1}}} \chi(x^{-1} h x) + \overline{\chi(x^{-1} h x)}, \\
    &= \sum_{\substack{x \in H \\ h \in xK_0x^{-1}}} \chi(x^{-1} h x) + \chi((x\omega)^{-1} h (x \omega)), \\
    &= 2 \cdot \sum_{\substack{x \in H \\ h \in xK_0x^{-1}}} \chi(x^{-1} h x),
    \end{align*}
    and therefore
    \[
    \Ind_K^H(\psi)(h) = \frac{1}{|K_0|} \sum_{\substack{x \in H \\ h \in xK_0x^{-1}}} \chi(x^{-1} h x),
    \]
    which is exactly the character of $\Ind_{S(K,u)}^{H \times C_2}(V)$ by Lemma \ref{lem:inductioncompatibility}. Then we see that
    \[
    \Ind_K^H(V|_{\bR[G]}) \cong \Ind_{S(K,u)}^{H \times C_2}(V)^{C_2},
    \]
    because they have the same character by the above and \cite[Lem.\ 6.2]{TAYSchur}. 
\end{proof}

\begin{proof}[Proof of Corollary \ref{cor:recoverserre}]
    Let $\chi \in \Char(H, \bR)$, and let $U \in \Rep_{\bR}(H)$ be the corresponding representation. We have that $\bC \otimes_{\bR} U \in \Rep_{\bC}^{\hspace{0.1em} \star}(H \times C_2)$, which also has character $\chi$ by \cite[Lem.\ 6.2]{TAYSchur}. From Theorem \ref{thm:fundamentalrepsexplicit} there are
    \begin{enumerate}
        \item ungraded subgroups $H_i \leq H$ and Real subgroups $S_j \leq H \times C_2$,
        \item $C_2$-graded quotients $H_i \twoheadrightarrow C_{n_i}$ and $S_j \twoheadrightarrow D_{2n_j}$,
        \item $\chi_i \in \Char(C_{n_i},\bC)$ and $\chi_j \in \Char^{\hspace{0.1em} \star}(D_{2n_j},\bC)$
    \end{enumerate}
    such that
    \[
        \chi = \sum_i m_i \Ind_{H_i \times 1}^{H \times C_2}(\chi_i) + \sum_j m_j \Ind_{S_j}^{H \times C_2}(\chi_j)
    \]
    for some $m_i, m_j \in \bZ$. For any $H_i$, there is a commutative diagram
\[\begin{tikzcd}
	{\Rep_{\bC}(C_{n_i})} & {\Rep_{\bC}(H_i)} && {\Rep_{\bC}^{\hspace{0.1em} \star}(H \times C_2)} \\
	&& {\Rep_{\bC}(H)} && {\Rep_{\bC}(H)} \\
	{\Rep_{\bR}(C_{n_i})} & {\Rep_{\bR}(H_i)} && {\Rep_{\bR}(H)}
	\arrow["{\Inf_{C_n}^{H_i}}", from=1-1, to=1-2]
	\arrow["{(-)|_{\bR[C_{n_i}]}}"', from=1-1, to=3-1]
	\arrow["{\Ind_{H_i \times 1}^{H \times C_2}}", from=1-2, to=1-4]
	\arrow["{\Ind_{H_i}^{H}}"', from=1-2, to=2-3]
	\arrow["{(-)|_{\bR[H_i]}}"', from=1-2, to=3-2]
	\arrow["{(-)|_H}", from=1-4, to=2-5]
	\arrow["\sim"', from=1-4, to=3-4]
	\arrow["{\Ind_{H \times 1}^{H \times C_2}}"', from=2-3, to=1-4]
	\arrow["{(-)|_{\bR}}"', from=2-3, to=3-4]
	\arrow["{\Inf_{C_{n_i}}^{H_i}}"', from=3-1, to=3-2]
	\arrow["{\Ind_K^H}"', from=3-2, to=3-4]
	\arrow["{(-)^{C_2}}"', draw=none, from=3-4, to=1-4]
	\arrow["{(-)_{\bC}}"', from=3-4, to=2-5]
\end{tikzcd}\]
The top triangle commutes by the transitivity of induction, and the two triangles on the right by \cite[Lem.\ 6.2]{TAYSchur} and \cite[Lem.\ 6.3]{TAYSchur}. The left-hand square and bottom-left triangle commute because restriction from $\bC$ to $\bR$ commutes with induction and inflation. In particular, each character $\Ind_{H_i \times 1}^{H \times C_2}(\chi_i)$, obtained by going along the top from $\Rep_{\bC}(C_{n_i})$ to $\Rep_{\bC}(H)$, is the same as the character obtained by going around the other outside route, and therefore in the image of of $\Rep_{\bR}(C_{n_i})$.

For the remaining subgroups, note that if $S \leq H \times C_2$ is a Real subgroup which surjects onto $D_{2n}$ for some $n \geq 3$, then $S = S(K,u)$ for some $K \leq H$ and $u \colon K \twoheadrightarrow C_2$. Indeed, the only other option is that $S = K \times C_2$, and such a $q \colon K \times C_2 \twoheadrightarrow D_{2n}$ being $C_2$-graded would induce an isomorphism $K / \ker(q) \times C_2 \cong D_{2n}$, which is obviously impossible. We cannot have $n = 2$ by Theorem \ref{thm:fundamentalrepsexplicit}. If $n = 1$, then $S = 1 \times C_2$, and by Lemma \ref{lem:inductioncompatibility} any induced representation is a sum of copies of $\bC[H]$, which corresponds to the representation $\bR[H]$ of $H$ over $\bR$. However, this comes from inducing the trivial representation of the trivial subgroup of $H$ over $\bR$, and so is permitted. We may therefore restrict to $S_i = S(K_i,u_i)$, for which there is a commutative diagram,
\[\begin{tikzcd}[column sep=large]
	{\Rep_{\bC}^{\hspace{0.1em} \star}(D_{2n_j}) } & {\Rep_{\bC}^{\hspace{0.1em} \star}(S(K_i,u_i)) } & {\Rep_{\bC}^{\hspace{0.1em} \star}(H \times C_2)} \\
	{\Rep_{\bR}(D_{2n_j}) } & {\Rep_{\bR}(K_i)} & {\Rep_{\bR}(H)} & {\Rep_{\bC}(H)}
	\arrow["{\Inf_{D_{2 n_j}}^{S(K_i, u_i)}}", from=1-1, to=1-2]
	\arrow["{(-)|_{\bR[D_{2n_j}]}}"', from=1-1, to=2-1]
	\arrow["{\Ind_{S(K_i,u_i)}^{H \times C_2}}", from=1-2, to=1-3]
	\arrow["{(-)|_{\bR[K_i]}}", from=1-2, to=2-2]
	\arrow["\sim", from=1-3, to=2-3]
	\arrow["{(-)|_H}", from=1-3, to=2-4]
	\arrow["{\Inf_{D_{2n_j}}^{K_i}}"', from=2-1, to=2-2]
	\arrow["{\Ind_{K_i}^H}"', from=2-2, to=2-3]
	\arrow["{(-)_{\bC}}"', from=2-3, to=2-4]
\end{tikzcd}\]
the left-hand square commuting by definition, the right-hand square by Lemma \ref{lem:twistedsubgroupinductionsplit}, and the right-hand triangle by \cite[Lem.\ 6.2]{TAYSchur}.
\end{proof}

\subsection{Classification of Types Over Cyclotomic Fields}

We can use Theorem \ref{thm:antilinearcyclotomicfinal} to give a complete description of the possible types of representations over cyclotomic fields and to describe the type in terms of easily computable invariants.

In addition to the cyclotomic field $E = \bQ(\zeta_e)$ and its index two maximal totally real subfield $E_+$, we set $f \coloneqq \exp(G)$ when $\exp(H) > 2$ and $f \coloneqq 4$ when $\exp(H) \leq 2$, and similarly define $F \coloneqq \bQ(\zeta_f)$ and $F_+ \coloneqq F \cap \bR = \bQ(\zeta_f + \zeta_f^{-1})$. Because $H$ has index two in $G$, either $e = f$ or $2e = f$, and so either $(F,F_+) = (E,E_+)$ or $F$ and $F_+$ are degree two extensions of $E$ and $E_+$ respectively.

\begin{thm}\label{thm:15waysforcyclotomic}
Let $V \in \Irr_{E}(H)$ and let $Z \in \Irr_{\bC}(G)$ correspond to $V_{\bC}$. Then the types of $V$ for the extensions $E_+ \leq E$, $F_+ \leq F$ and $\bR \leq \bC$ are in one of the cases of Table \ref{tab:15typesQe}. \emph{
\begin{table}[h!] 
\centering
\bgroup
\caption{The Possible Types of $V$ over $\bQ(\zeta_e)$ and $\bQ(\zeta_f)$}
\label{tab:15typesQe}
\def\arraystretch{1.1}
\begin{tabular}{ c  c  c  c }
    \toprule
    $\bR \leq \bC$ & $\bQ(\chi_Z) \subset E$ & $E_+ \leq E$ & $F_+ \leq F$ \\
    \midrule 
 I & Yes & (1) & (1) \\
 I & No & (8) & (1) \\ \midrule
 II & Yes & (5) & (5) \\
 II & No & (9) & (5) \\ \midrule
 III & Yes & (2) & (2) \\ \midrule
 IV & Yes & (10) & (10) \\
 IV & No & (14) & (10) \\ \midrule
 V & Yes & (12) & (12) \\ \midrule 
 VI & Yes & (15) & (15) \\ \midrule 
 VII & Yes & (11) & (11) \\ \midrule
 VIII & Yes & (16) & (16) \\
 VIII & No & (24) & (16) \\ \midrule 
 IX & Yes & (18) & (18) \\
 IX & No & (22) & (18) \\ \midrule 
 X & Yes & (17) & (17) \\ \bottomrule 
\end{tabular}
\egroup
\end{table}
}
\end{thm}

Here $\bQ(\chi_Z)$ is the extension of $\bQ$ generated by the character values of $\chi_Z$.

\begin{remark}\label{rem:FSforcyclo}
    The type of $V_{\bC}$ is determined by three computable quantities: the classical and antilinear Frobenius-Schur indicators, and whether or not $\chi_V \cong {}^{\omega}\chi_V$. From the Theorem, the type of $V_{F}$ and $V_{\bC}$ are the same (just using different labels), and therefore the type of $V_{F}$ is also determined by these quantities. We also see that the type of $V$ over $E$ is determined by these quantities and whether or not $\bQ(\chi_Z) \subset E$. These are all computationally implementable, and we magma code for this at \cite{GITHUB}.
\end{remark}

\begin{remark}
    All cases of Theorem \ref{thm:15waysforcyclotomic} do actually occur! Indeed, the examples of Section \ref{sect:realise24ways} for these 15 cases are over $E_+ \leq E$.
\end{remark}

\begin{proof}
    The type of $V$ for the extension $\bL[H] \leq \bL_{\star \hspace{-0.1em}}[G]$ is the same for all three choices of $\bK \leq \bL$ by Theorem \ref{thm:antilinearcyclotomicfinal}. Similarly, the type of $V$ for the extension $\bK[H] \leq \bL[H]$ is the same for all three choices of $\bK \leq \bL$ by Theorem \ref{thm:classicalcyclotomic} (which itself follows from Theorem \ref{thm:antilinearcyclotomicfinal}, see Section \ref{sect:recoveryofclassical}). In particular, the type of $V_{\bC}$ is uniquely determined by the type of $V$ over $E$, as it is uniquely by the type of $V_{\bC}$ for the extensions $\bR[H] \leq \bC[H]$ and $\bC[H] \leq \bC_{\star \hspace{-0.1em}}[G]$ and whether or not $\overline{\chi}_V \cong {}^{\omega}\chi_V$. Furthermore, the types of $V_{F}$ for the extensions $F_+[G] \leq F[G]$ and $F[H] \leq F[G]$ are the same as the types of $V_{\bC}$ for the extension $\bR[G] \leq \bC[G]$ and $\bR[H] \leq \bC[G]$ respectively, by Theorem \ref{thm:classicalcyclotomic} and the fact that $F$ is a splitting field for $H$ and $G$. Therefore, the type of $V_{F}$ is equal to the type of $V_{\bC}$, any row of Table \ref{tab:24types} being determined by its first four columns.

    Note also that whether or not $\bQ(\chi_Z) \not\subset E$ is determined by the type of $V$ over $E$, as this is equivalent to $\End_{E[G]}(Z_0)$ being $2$-dimensional over $E$, where $Z_0 \in \Irr_{E}(H)$ corresponds to $V$. Indeed, if $e=f$, then neither can happen, and if $f = 2e$ then both are equivalent to being in Case B of Table \ref{tab:thewaysLHSsquare} for the extension $E[G] \leq F[G]$.
    
    Therefore, it suffices to show that the only types that can occur for $V$ over $E$ are those 15 listed in Table \ref{tab:15typesQe} out of the possible 24 of Theorem \ref{thm:24ways}. Cases (3), (6), (7), (13), (20) and (21) are ruled out immediately using Table \ref{tab:24types6} because $E$ is a splitting field for $H$, and so we are left with types (4), (19) and (23).

    If $V$ had type (4), then $V_{F}$ would have type $(1)$. If $Z_e \in \Irr_{E}(G)$ corresponds to $V$, then as $V$ has type (4), $Z_e$ is absolutely irreducible by Table \ref{tab:24types6}, and so $Z_{e, F}$ corresponds to $V_{F}$. Then $\Stab_{K_4}(Z_{e, F}) = \Stab_{K_4}(Z_e) = \langle c \rangle$ by Table \ref{tab:24types5} as $V$ has type (4), which contradicts the fact that $\Stab_{K_4}(Z_{e, F}) = \langle a \rangle$ because $V_{F}$ has type $(1)$, again by Table \ref{tab:24types5}.

    If $V$ had type (19), then $V_{F}$ would have type $(16)$. Similarly, if $Z_e \in \Irr_{E}(G)$ corresponds to $V$, then as $V$ has type (4), $Z_e$ is absolutely irreducible by Table \ref{tab:24types6}, and therefore $Z_{e,F}$ corresponds to $V_{F}$. We again obtain a contradiction: from Table \ref{tab:24types5} this would force $\Stab_{K_4}(Z_{e,F})$ to be simultaneously $\langle a \rangle$ and $\langle c \rangle$.

    Finally, if $V$ had type (23), then $V_{F}$ would have type $(16)$. If $U_e \in \Irr_{E_+}(G)$ corresponds to $V$, then as $V$ has type (23) $U_{e, E}$ is irreducible, and $\End_{E[G]}(U_{e, E})$ is $2$-dimensional over $E$. In particular, as $F$ splits $G$ then from Table \ref{tab:thewaysLHSsquare} of Section \ref{sect:classicalindextwofieldextoptions}
    \[
    (U_{e, E})_{F} = U_{e, F} = Z_1 \oplus Z_2
    \]
    for some non-isomorphic $Z_1, Z_2 \in \Irr_{F}(G)$. On the other hand, if $U_f \in \Irr_{F_+}(G)$ corresponds to $V_{F}$, then as $V_{F}$ has type $(16)$, $U_{f, F} = Z_3 \oplus Z_3$ for some $Z_3 \in \Irr_{F}(G)$ by Table \ref{tab:24types2}. However the irreducible factors of $U_{e, F_+}$ correspond to $V_{F}$, and therefore $U_f \hookrightarrow U_{e, F_+}$. In particular, 
    \[
        Z_3 \oplus Z_3 = U_{f, F} \hookrightarrow U_{e, F} = Z_1 \oplus Z_2,
    \]
    so $Z_1 \cong Z_3 \cong Z_2$, which contradicts the fact that $Z_1$ and $Z_2$ are non-isomorphic.
\end{proof}

\subsection{Applications to Classical Schur Indices}

We can give some application of Theorem \ref{thm:15waysforcyclotomic} to classical Schur indices for groups which admit an index two subgroup. In what follows, whenever $f > e$, let $\sigma$ denote the non-trivial element of $\Gal(F / E)$, which we also use to denote its restriction to the non-trivial element of $\Gal(F_+ / E_+)$.

\begin{prop}\label{prop:indextwosubgroupclassicalCextension}
    Suppose that $G$ is a finite group with an index two subgroup $H$ of exponent~$e$. For $Z \in \Irr_{\bC}(G)$, let $Z_f \in \Irr_{F}(G)$ and $Z_e \in \Irr_{E}(G)$ correspond to $Z$. Then $Z$ is in one of the two cases of Table \ref{tab:2optionsZQef}
    \emph{
\begin{table}[h!] 
\centering
\bgroup
\caption{The Possible Models for $Z$ over $\bQ(\zeta_e)$ and $\bQ(\zeta_f)$}
\label{tab:2optionsZQef}
\def\arraystretch{1.1}
\begin{tabular}{ c c c}
    \toprule 
    $\bQ(\chi_Z) \subset E$ & $Z_{f, \bC}$ & $Z_{e, F}$ \\ \midrule
 Yes & $Z$ & $Z_f$ \\
 No & $Z$ & $Z_f \oplus Z_f^{\sigma}$ \\ \bottomrule
\end{tabular}
\egroup
\end{table}
}
In particular $\overline{m}_{E}(Z) = 1$: $Z$ admits a model over $E(\chi_Z)$.
\end{prop}

\begin{proof}
    First, because $F$ splits $G$ \cite[\S 12.3]{SERRE}, $Z_{f, \bC} \cong Z$. If $E = F$ then we are in the first case, so suppose not, meaning that $2e = f$ and the extension $F / E$ has degree $2$.
    
    From Section \ref{sect:classicalindextwofieldextoptions}, there are three cases possible for the extension $E[G] \leq F[G]$, listed in Table \ref{tab:thewaysLHSsquare}. Case A corresponds to the first line of Table \ref{tab:2optionsZQef}, and Case B to the second, and therefore the claim is that Case C cannot occur. Indeed, in Case C, $\End_{E[G]}(Z_e)$ has dimension $4$ over $E$, from Table \ref{tab:thewaysLHSsquare} and the fact that $\End_{F[G]}(Z_f)$ has dimension $1$ over $F$. However this is impossible: $Z_e$ is also in one of the four cases of Proposition \ref{prop:algindressummary} for the extension $E[H] \subset E[G]$ (cf.\ Section \ref{sect:indextwosubgroup}), and hence $\End_{E[G]}(Z_f)$ has dimension $1$ or $2$, because $E$ splits $H$.
\end{proof}

The proof of Proposition \ref{prop:indextwosubgroupclassicalCextension} is relatively straightforward using Proposition \ref{prop:algindressummary}, and doesn't use Theorem \ref{thm:antilinearcyclotomicfinal} and its consequences. It could instead be proven using the classification of Theorem \ref{thm:15waysforcyclotomic} as we do now in the following, where the use of Theorem \ref{thm:15waysforcyclotomic} is now essential.

\begin{prop}\label{prop:indextwosubgroupclassicalRextension}
    Suppose that $G$ is a finite group with an index two subgroup $H$ of exponent~$e$. For $U \in \Irr_{\bR}(G)$, let $U_f \in \Irr_{F_+}(G)$ and $U_e \in \Irr_{E_+}(G)$ correspond to $U$. Then $U$ is in one of the two cases of Table \ref{tab:2optionsUQef}
    \emph{
\begin{table}[h!] 
\centering
\bgroup
\caption{The Possible Models for $U$ over $\bQ(\zeta_e)_+$ and $\bQ(\zeta_f)_+$}
\label{tab:2optionsUQef}
\def\arraystretch{1.1}
\begin{tabular}{c c c }
    \toprule
    $\bQ(\chi_U) \subset E_+$ & $U_{f, \bR}$ & $U_{e, F_+}$ \\ \midrule 
 Yes & $U$ & $U_f$ \\
 No & $U$ & $U_f \oplus U_f^{\sigma}$ \\ \bottomrule
\end{tabular}
\egroup
\end{table}
}
In particular $\overline{m}_{E_+}(U) = 1$: $U$ admits a model over $E_+(\chi_U)$.
\end{prop}

\begin{proof}
    By Theorem \ref{thm:classicalcyclotomic} we have that $U_{f, \bC} \cong U$. If $E = F$ then we are in the first case, so suppose not, meaning that $2e = f$ and the extension $F / E$ has degree $2$. From Section \ref{sect:classicalindextwofieldextoptions}, there are three possibilities for the extension $E_+[G] \leq F_+[G]$. Case A corresponds to the first line of Table \ref{tab:2optionsUQef}, and Case B to the second, and therefore the claim is that Case C cannot occur.
    
    To see this, let $V \in \Irr_{E}(H)$ correspond to $U_e$, and so $V_{F}$ corresponds to $U_f$. One can check using Table \ref{tab:24types6} that in each of the 15 cases of Table \ref{tab:15typesQe} the ratio
    \[
    \dim_{F} \End_{F[G]}(U_f) : \dim_{E} \End_{E[G]}(U_e)
    \]
    is either $1 : 1$ or $1 : 2$, and therefore that Case C, when the ratio would be $1 : 4$, never occurs.
\end{proof}

\begin{remark}
    Propositions \ref{prop:indextwosubgroupclassicalCextension} and \ref{prop:indextwosubgroupclassicalRextension} are not true without the existence of the index two subgroup. For example, if one instead assumes that $G$ has even order and takes $m \coloneqq \text{Exponent}(G) / 2$, then the conclusion of Theorem \ref{prop:indextwosubgroupclassicalRextension} is false. An example is given by SmallGroup(48,41). However, it seems that the statement remains is also true if one instead uses the assumption that $G$ has even order and takes $m \coloneqq |G|/2$ (we have checked this computationally up to $G$ of order 200).
\end{remark}

\section{Unitary Representations}\label{sect:unitaryreps}

In this section we no longer assume that $G$ is finite, but instead restrict our attention to unitary (antilinear) representations. This gives us several additional tools that make the theory reminiscent of that of a finite group.

\begin{defn}
  Let $(V,\tau)$ be a finite-dimensional hermitian vector space (over $\bC$). We write $\UU(V)$ for the unitary group of $V$,
  and $\UU^\ast (V)$ for the set of antilinear transformations $f: V\to V$, ``flipping'' the hermitian form, i.e.,
\[ \tau (f(x), f(y)) = \tau (y,x) = \overline{\tau (x,y)} \, . \]
  We write $\UR (V)\coloneqq \UU(V) \cup \UU^\ast (V)$, which we view as a $C_2$-graded group (cf.\ Section \ref{sect:fundreps}) under composition.

  An \emph{(antilinear) unitary representation} of a $C_2$-graded group $(G,\epsilon_G)$ is a homomorphism of $C_2$-graded groups $G\to \UR (V)$, where $V$ is a finite-dimensional hermitian vector space.
\end{defn}

Observe that an antilinear map will not preserve a hermitian form in general. Indeed, writing vectors as columns in an orthonormal basis, the form becomes $\tau (x,y) = x^T \overline{y}$, while an antilinear transformation is given by $x \mapsto A \overline{x}$. Then the flipping/preserving conditions become
\[ \tau (A\overline{x},A \overline{y}) = \overline{x}^T A^T \overline{A} y  \stackrel{?}{=}
\begin{cases}
	\tau(x,y) =  x^T \overline{y} & \ \mbox{preserving},\\
	\overline{\tau(x,y)} =  \overline{x}^T {y} & \ \mbox{flipping}.
\end{cases}\]
for all $x,y$. Clearly, the preserving condition leads to a contradiction by choosing various $x$ and $y$, while the flipping condition is
just $A^T \overline{A}=I$ which, for instance, all real orthogonal matrices satisfy.

Note that if $G$ is finite, all antilinear representations over $\bC/\bR$ are unitary by Weyl's Unitary Trick:
pick any positive definite hermitian form $\tau$ on $V$, then integrate it over $G$ to get a new form
\[ \widetilde{\tau} (x,y) = \sum_{g\in H} \tau (gx,gy) + \sum_{g\in G\setminus H} \tau (gy,gx)\, .\]

\subsection{Bohr compactification}
The Bohr compactification of a group $G$ is the universal homomorphism $G\to \widehat{G}$ to a compact group.
In the literature it is often used for topological groups but here we are working only with discrete groups. 

The best way to think of the Bohr compactification $\widehat{G}$ is that it is the group one gets from Tannaka-Krein Duality, applied to the category of unitary (linear) representations of $G$. Let us adopt a more down-to-earth definition.

\begin{defn}\label{def:bohrcompactification}
The Bohr compactification $\widehat{G}$ is the closure of the image of $G$ in $\prod_V \UU(V)$ where the product is taken over isomorphism classed over irreducible unitary (linear) representations of $G$.
\end{defn}

From the definition, we have a continuous map $\alpha : G \to \widehat{G}$. Since the Bohr compactification is not really our focus, we will just list its standard properties. In what follows, we let $\beta \colon G \rightarrow \Prof(G)$ denote the profinite completion of $G$.
We write $\URep(\widehat{G})$ for the set of isomorphism classes of finite-dimensional continuous unitary representations of $\widehat{G}$, and $\UIrr(\widehat{G})$ for the subset of irreducible ones. We write the same for $G$, noting that the continuity requirement is vacuous and can be dropped.

\begin{prop}\label{prop:bohrcompactificationproperties} The following properties hold for a discrete group $G$.
  \begin{enumerate}
  \item $\alpha (G)$ is dense in $\widehat{G}$.
  \item If $G$ is abelian, then $\widehat{G}$ is abelian.
  \item If $\beta: G \to K$ is a homomorphism to a compact Hausdorff group $K$, then there is a unique continuous homomorphism $f: \widehat{G}\to K$ such that $\beta = f \circ \alpha$.
  \item In particular, the association $G \mapsto \widehat{G}$ is functorial, and is left adjoint to the forgetful functor from compact groups to discrete groups.
  \item Sending $(V, \rho) \mapsto (V, \rho \circ \alpha)$ defines bijections
  \[
  \URep(\widehat{G}) \xlongrightarrow{\sim} \URep(G), \qquad \UIrr(\widehat{G}) \xlongrightarrow{\sim} \UIrr(G).
  \]
  \item The natural map $\widehat{G} \rightarrow \Prof(G)$ induced by $\beta \colon G \rightarrow \Prof(G)$ induces an isomorphism 
  \[
  \pi_0(\widehat{G}) = \widehat{G} / \widehat{G}_0 \xrightarrow{\sim} \Prof(G),
  \]
  as topological groups, where $G_0$ is the connected component of the identity in $\widehat{G}$.
  \item If $H$ is a finite-index normal subgroup of $G$, then the induced map $\widehat{H} \rightarrow \widehat{G}$ is injective, the image $\widehat{H}$ is a finite-index normal subgroup of $\widehat{G}$, and $\alpha_G \colon G \rightarrow \widehat{G}$ induces an isomorphism
  \[
  G/H \cong \widehat{G}/ \widehat{H}.
  \]
  \end{enumerate}  
\end{prop}  

\begin{proof}
    The first two points follow from the definition, the second because when $G$ is abelian each factor $U(V)$ is isomorphic to $S^1$. Points (3) to (5) follow from \cite[Prop.\ 5]{BEK1}. The sixth point is \cite[Prop.\ 7]{BEK1}, and the final point follows from \cite[Prop.\ 9]{BEK1} and \cite[Prop.\ 10]{BEK1}.
\end{proof}
 
Note that the map $\alpha : G \to \widehat{G}$ is not necessarily injective. For example, if $G$ is a non-linear simple group, then $\widehat{G}$ is the trivial group. If $\alpha$ is injective, then the topology on $G$ induced from $\widehat{G}$ is discrete if and only if $G$ is finite because $\widehat{G}$ is compact Hausdorff and $G$ is dense in $\widehat{G}$. In particular, the Bohr compactification is not a topological compactification in the usual sense, and perhaps ``Bohr completion'' would be a more appropriate term here as well.

We also note the similarity between the Bohr compactification of a topological group $G$ and the Bohr completion of a $\bK$-algebra $A$ as defined in Section \ref{Bohr_completion} (for example, compare Definition \ref{def:bohrcompactification} with Lemma \ref{Bohr2}). Certain properties that hold for one may not translate obviously to the other: for example the correspondence of simple modules of Proposition \ref{prop:bohrcompactificationproperties} fails for Bohr completions of algebras (see Remark \ref{rem:bohrcompletionexample}), but one instead has Lemma \ref{lem:idemirrcorresp}.

\subsection{Frobenius-Schur indicators for compact groups} \label{section:FS_UIrr} 

If $G$ is a finite $C_2$-graded group, we have three Frobenius-Schur indicators in play \cite{RUMTAY, RUMTAY3}. The indicators involve summations over $G$ or $H$, the index 2 subgroup. For infinite groups, they work in the same way but only for unitary representations.

First, let $G$ be a compact Hausdorff topological group, and $H$ an index two (normal) subgroup of $G$. Let $\mu_G$, $\mu_H$ denote Haar measures on $G$ and $H$ normalised in such a way that $\mu_G(G) = 1$ and $\mu_H(H) = 1$. Define functions 
\begin{equation} \label{FS_formulas}
\sF,{\sF^{\sharp}} : \UIrr (H) \rightarrow \bC , \quad
\sF (V) =
\int_{H} \chi_V (x^2) \mu_G(dx), \quad
    {\sF^{\sharp}} (V) =
    \int_{G\setminus H} \chi_V (x^2) \mu_G(dx),
    \end{equation}

We note that for $V \in\UIrr (H)$, $V$ is finite-dimensional, so we may consider $V$ as an irreducible representation of $H$ viewed as an abstract group, and associate to $V$ a type $\fT(V)$, a D-type $\fD(V)$, an S-type $\fS(V)$ and a C-type $\fC(V)$.

The following is a straightforward generalisation to compact Hausdorff groups of the relationship between the Frobenius-Schur indicator and the classical and semilinear types.

\begin{prop}\label{prop:FSexttocompact}
    Suppose that $V \in \UIrr (H)$. Then 
    \[
        \sF(V) = \left\{
		\begin{array}{cl}
			1 & \mbox{if } \: \fC(V) = \textbf{\emph{A}} \\
            0 & \mbox{if } \: \fC(V) = \textbf{\emph{B}} \\
			-1 & \mbox{if } \: \fC(V) = \textbf{\emph{C}}
		\end{array}
        \right. ,
        \qquad 
        \sF^{\sharp}(V) = \left\{
		\begin{array}{cl}
			1 & \mbox{if } \: \fS(V) = \textbf{\emph{A}} \\
            0 & \mbox{if } \: \fS(V) = \textbf{\emph{B}} \\
			-1 & \mbox{if } \: \fS(V) = \textbf{\emph{C}}
		\end{array}
        \right. .
    \]
\end{prop}

\begin{proof}
    For $\sF(V)$, this is well-known, for example, see \cite[Prop.\ 6.8]{DieckBrocker}. For $\sF^{\sharp}(V)$, this follows exactly as in the proof of \cite[Thm.\ 4.2]{RUMTAY}, by computing $\sF^{\sharp}(V)$ for each possible D-type.
\end{proof}

\subsection{Frobenius-Schur indicators for discrete groups}

Now suppose that $G$ is a discrete, potentially infinite, group, with an index two subgroup $H$.

\begin{lemma}
    Let $V \in \UIrr(H)$, and let $\widehat{V} \in \UIrr(\widehat{H})$ be the corresponding irreducible unitary representation of $\widehat{H}$. Then $\fC(V) = \fC(\widehat{V})$ and $\fS(V) = \fS(\widehat{V})$.
\end{lemma}

\begin{proof}
    To ease notation, we write $(-)^{\wedge}$ instead $\widehat{(-)}$. We first show that $\fC(V) = \fC(V^{\wedge})$. Let $X^{\wedge} \in \Irr_{\bR}(\widehat{H})$ denote the irreducible $\bR$-linear representation corresponding to $V^{\wedge}$, which is continuous because $X^{\wedge}$ is a subrepresentation of the restriction to $\bR$ of $V^{\wedge} \oplus \overline{V^{\wedge}}$. We can pull back along $\alpha$ to obtain an $\bR$-linear representation $X$ of $H$, which is irreducible by continuity, because $H$ is dense in $\widehat{H}$. By construction, $X$ is a constituent of $V|_{\bR}$, as $X^{\wedge}$ is a constituent of $V^{\wedge}|_{\bR}$. Again by continuity, the restriction map
    \[
    \End_{\bR \widehat{H}}(X^{\wedge}) \rightarrow \End_{\bR H}(X)
    \]
    is an isomorphism, and thus $\fC(V) = \fC(V^{\wedge})$. A similar argument shows that $\fS(V) = \fS(V^{\wedge})$.
\end{proof}

Therefore, we can compute the D-type of $V$ using integration on $\widehat{H}$ and $\widehat{G}$, by Proposition \ref{prop:FSexttocompact}.
We note that this is not clear how to do this for a general group, where a useful description of the Bohr compactification is often not available.

\subsection{Frobenius-Schur indicators for amenable groups}
Countable amenable groups are interesting because
the integrals~\ref{FS_formulas} can be computed intrinsically, without a reference to the Bohr compactification.
As amenable groups are not our focus, we take the following definition, it being the most convenient (see \cite[Prop.\ 4.7.5]{CSCoornaert}). 

\begin{defn}
A discrete countable group $G$ is amenable if it admits a F{\o}lner sequence $(F_n)_n$. The latter consists of a sequence of non-empty finite subsets $F_n\subseteq G$ for $n \geq 1$ such that
\[ \forall g\in G \quad \lim_{n \rightarrow \infty} \frac{|gF_n \oplus F_n|}{|F_n|} = 0\, ,\]
where $A \oplus B$ is the symmetric difference of the subsets.
\end{defn}

Let $\sM(G) \subset \ell^{\infty}(G,\mathbb{R})^\ast$ denote the space of means on $G$, consisting of linear functionals $\mu \in \ell^{\infty}(G,\mathbb{R})^\ast$ which are positive ($\mu(f) \geq 0$ whenever $f \geq 0$) and satisfy $\mu(1)=1$.
It follows from the Banach–Alaoglu Theorem that $\sM(G)$ is compact in the weak-$\ast$ topology \cite[Thm.\ 4.2.1]{CSCoornaert}.
The F{\o}lner sequence gives a sequence of averaging functionals
\[
  \mu_n \in \sM(G), \qquad
  \mu_n(f) = \frac{1}{|F_n|}\sum_{x\in F_n} f(x),
\]
which has a cluster point $\mu_G \in \sM(G)$ in the weak-$\ast$ topology by the compactness of $\sM(G)$.
We note that $\mu_G$ is not unique in general. For $\mu \in \ell^{\infty}(G,\mathbb{R})^*$ and $f \in \ell^{\infty}(G,\mathbb{C})$, we define 
\[
\mu(f) \coloneqq \mu(\Re(f)) + i \mu(\Im(f))
\]
and denote
\[
\int_G f(x) d \mu(x) \coloneqq \mu(f).
\]

Recall that $\alpha\colon G\to\widehat G$ denotes the Bohr compactification of $G$,
and $\mu_{\widehat G}$ denotes the normalized Haar measure on $\widehat G$.
We are going to use the following lemma only for matrix elements.
We suspect it is known to specialists but we could not a reference, so we provide a full proof.

\begin{lemma}\label{lem:folnerintegral}
Let $G$ be a discrete countable amenable group with F{\o}lner sequence
$(F_n)_n$, $(\mu_n)_n$ -- the corresponding sequence of means,
$\mu_G$ -- any weak-$\ast$ cluster point of $(\mu_n)_n$. Then for every
$f\in C(\widehat G,\mathbb{C})$, we have
\[
  \int_{\widehat G} f(z)\,d\mu_{\widehat G}(z)
  \;=\; \int_G (f\circ\alpha)(x)\,d\mu_G(x)
  \;=\; \lim_{n\to\infty}\mu_n(f\circ\alpha).
\]
\end{lemma}

\begin{proof}
  Intuitively, it suffices to prove the statement for matrix elements because they are dense (Stone-Weierstrass Theorem), and it is true
  for matrix elements because averaging defines an orthogonal projection $L^2 (\widehat{G},\bC) \to L^2 (\widehat{G},\bC)^G = \bC$.
  Now let us give full details.

Pick  $(V,\pi) \in \UIrr(\widehat{G})$, so that $(V,\rho \coloneqq \pi\alpha) \in \UIrr({G})$. Set
\[
  A_n=\frac{1}{|F_n|}\sum_{x\in F_n}\rho(x),
  \qquad
  P=\int_{\widehat G}\pi(z)\,d\mu_{\widehat G}(z)
  \qquad\in\mathrm{End}_{\bC}(V) \, .
  \]
We claim that the sequences of matrices $A_n$ converges to the matrix $P$.

To prove the claim, examine $P$. It is an endomorphism of $V$ as a $\widehat{G}$-module. By Schur's Lemma $P=\alpha I$.
Moreover, $P(v)$ is $\widehat{G}$-invariant: $\pi (g) (P(v)) = P(v)$. Hence, $P = 0$ if $V$ is non-trivial, and $P=I$ if $V$ is trivial.

To examine $A_n$, we use the operator norm on $\mathrm{End}_{\bC}(V)$.
Each $\rho(x)$ is unitary, so $\| \rho(x) \| =1$ and $\|A_n\|\leq 1$.
Thus, all $A_n$ belong to the (compact) unit ball and the sequence $(A_n)$ has convergent subsequences.
To prove the claim, it remains to show that $P$ is the limit of every convergent subsequence.

Since for each
$g\in G$ we have $\rho(g)\sum_{x\in F_n}\rho(x)=\sum_{y\in gF_n}\rho(y)$, we can estimate 
\[
  \|\rho(g)A_n-A_n\|
  =\frac{1}{|F_n|}\Big\|\sum_{y\in gF_n}\rho(y)-\sum_{x\in F_n}\rho(x)\Big\|
  \leq\frac{|gF_n\oplus F_n|}{|F_n|}\xrightarrow[n\to\infty]{}0 .
  \]
  Thus, the limit of any convergent subsequence is a linear operator $V \to V^G = V^{\widetilde{G}}$.
  The claim follows immediately for a non-trivial $V$. For the trivial $V$ we need an additional observation that $A_n = I$.

  At this stage, we have proved that
\[
  \int_{\widehat G} f(z)\,d\mu_{\widehat G}(z)
  \;=\; \lim_{n\to\infty}\mu_n(f\circ\alpha), 
\]
  if $f$ is a matrix coefficient of a finite-dimensional unitary representation, i.e., 
  \[ f (x)=\langle\pi(x)v,w\rangle, \ \mbox{ for some } v,w\in (V,\pi) \in \URep(\widehat{G}) \, , \]
  since we can extend the statement of the proposition from matrix coefficients of irreducible representations to their finite linear combinations.
  Let us extend this to all continuous functions.
  
Let $\mathcal T\subset C(\widehat G,\mathbb{C})$ be the linear span of the
matrix coefficients of all finite-dimensional continuous unitary
representations of $\widehat G$.
By the combination of Peter-Weyl and Stone-Weierstrass theorems,
$\mathcal T$ is dense in $C(\widehat G,\mathbb{C})$ in the sup norm.
Now consider
$f\in C(\widehat G,\mathbb{C})$.
Given  $\varepsilon>0$, choose $h\in\mathcal T$
such that  $\|f-h\|_\infty<\varepsilon / 3$.
By the triangle inequality, 
\[
  \Big|\mu_n(f\circ\alpha)-\int_{\widehat{G}} f\,d\mu_{\widehat G}\Big|
  \leq
  \Big|\mu_n((f-h)\circ\alpha)\Big|
  +\Big|\mu_n(h\circ\alpha)-\int_{\widehat{G}} h\,d\mu_{\widehat G}\Big|
  +\Big|\int_{\widehat{G}}(h-f)\,d\mu_{\widehat G}\Big|. 
  \]
Since
$\|(f-h)\circ\alpha\|_\infty\leq\|f-h\|_\infty < \varepsilon/3$  
and $\int\int_{\widehat{G}} d \mu_{\widehat G} =1$, we can estimate the first and third terms
\[
  \Big|\mu_n(f\circ\alpha)-\int_{\widehat{G}} f\,d\mu_{\widehat G}\Big|
  < \frac{2}{3} \varepsilon
   +\Big|\mu_n(h\circ\alpha)-\int_{\widehat{G}} h\,d\mu_{\widehat G}\Big|.
  \]
  Now choose $N$ such that for $n>N$ we have $|\mu_n(h\circ\alpha)-\int_{\widehat{G}} h\,d\mu_{\widehat G}| < \varepsilon/3$.
  Then for $n>N$,
  \[|\mu_n(f\circ\alpha)-\int_{\widehat{G}} f\,d\mu_{\widehat G}| < \varepsilon \, . \]

  It remains to examine the middle integral $\int_G (f\circ\alpha)(x)\,d\mu_G(x)$. Obviously, our functions belong to
$\ell^\infty(G,\mathbb{R})$ and for any
$h\in\ell^\infty(G,\mathbb{R})$ such that $\mu_n(h)\to c$ for some number $c$, we have
$\mu_G(h)=c$.
\end{proof}

If $G$ is an abstract group, and $H$ is an index two subgroup, then $G$ is amenable if and only if $H$ is amenable. Given an odd element $\omega \in G \setminus H$ and a F{\o}lner sequence $F_n$ of $H$, then $(F_n \cup \omega F_n)_n$ is a F{\o}lner sequence for $G$.
 
\begin{prop}\label{prop:FSamenable}
    Suppose that $H$ is a discrete countable amenable group, $G$ is a Real structure on $H$, $F_n$ is a F{\o}lner sequence for $H$ and $\omega \in G \setminus H$. Then for any $V \in \UIrr(H)$ the quantities
    \[
    \sF_{\Sigma}(V) \coloneqq \lim_{n \rightarrow \infty} \left( \frac{1}{|F_n|} \sum_{x \in F_n} \chi_V(x^2)\right),
    \qquad
    \sF_{\Sigma}^{\sharp}(V) \coloneqq \lim_{n \rightarrow \infty} \left( \frac{1}{|F_n|} \sum_{x \in F_n} \chi_V((\omega x)^2)\right)
    \]
    converge and satisfy
    \[
        \sF_{\Sigma}(V) = \left\{
		\begin{array}{cl}
			1 & \mbox{if } \: \fC(V) = \textbf{\emph{A}} \\
            0 & \mbox{if } \: \fC(V) = \textbf{\emph{B}} \\
			-1 & \mbox{if } \: \fC(V) = \textbf{\emph{C}}
		\end{array}
        \right. ,
        \qquad 
        \sF_{\Sigma}^{\sharp}(V) = \left\{
		\begin{array}{cl}
			1 & \mbox{if } \: \fS(V) = \textbf{\emph{A}} \\
            0 & \mbox{if } \: \fS(V) = \textbf{\emph{B}} \\
			-1 & \mbox{if } \: \fS(V) = \textbf{\emph{C}}
		\end{array}
        \right. .
    \]
\end{prop}

\begin{proof}
    Let $\widehat{V} \in \UIrr(\widehat{G})$ be the unique extension of $V$ from Proposition \ref{prop:bohrcompactificationproperties}, and let $f \in C(\widehat{G}, \bC)$ defined by $f(x) \coloneqq f(x) \coloneqq \chi_{\widehat{V}}(x^2)$, so that 
    \[
    \sF(\widehat{V}) = \int_{\widehat{H}} f(x) d \mu_{\widehat{H}}(x), \qquad \sF^{\sharp}(\widehat{V}) = \int_{\widehat{G} \setminus \widehat{H}} f(x) d \mu_{\widehat{G}}(x).
    \]
    By Proposition \ref{prop:FSexttocompact} it suffices to show that $\sF(\widehat{V}) = \sF_{\Sigma}(V)$ and $\sF^{\sharp}(\widehat{V}) = \sF^{\sharp}_{\Sigma}(V)$. By definition,
    \[
    \sF_{\Sigma}(V) = \lim_{n\to\infty}\mu_{H,n}(f\circ\alpha_H),
    \]
    and therefore that $\sF(\widehat{V}) = \sF_{\Sigma}(V)$ follows directly from Lemma \ref{lem:folnerintegral}. For the second equality, note that by Proposition \ref{prop:bohrcompactificationproperties}(7), the induced map $\widehat{H} \rightarrow \widehat{G}$ is injective and we can view $\widehat{H}$ as an index two normal subgroup of $\widehat{G}$. Then we again use Lemma \ref{lem:folnerintegral} for both $H$ and $G$ to compute that
    \begin{align*}
        \sF^{\sharp}_{\Sigma}(V)
        &= \lim_{n\to\infty}\mu_{G,n}(f\circ\alpha_G) - \lim_{n\to\infty}\mu_{H,n}(f\circ\alpha_H), \\
        &= \int_{\widehat{G}} f(x) d \mu_{\widehat{G}}(x) - \int_{\widehat{H}} f(x) d \mu_{\widehat{H}}(x), \\
        &= \int_{\widehat{G} \setminus \widehat{H}} f(x) d \mu_{\widehat{G}}(x), \\
        &= \sF^{\sharp}(\widehat{V})
    \end{align*}
    as required.
\end{proof}

\section{Dependence of the D-Type on the Real Structure}\label{sect:deponrealstr}

In this section we consider how, for a group $H$, the D-type of some $V \in \Irr_{\bC}(H)$ depends on a given Real structure $G$.

\begin{lemma}\label{lem:indepK4action}
    The action of $K_4$ on $\Irr_{\bC}(H)$ depends only on the coupling of $G$.
\end{lemma}

\begin{proof}
    This is clear, as this action is defined in terms of conjugation by an element of $G \setminus H$.
\end{proof}

\begin{defn}
    For any $V \in \Irr_{\bC}(H)$, we write $\mathfrak{D}(V; G)$ for the D-type of $V$ with respect to the real structure $G$. When $H$ is finite, we set
    \[
    \sF(V ; G) = \frac{1}{|H|} \sum_{g \in G \setminus H} \chi_V(g^2).
    \]
\end{defn}

\begin{cor}
    Let $V \in \Irr_{\bC}(H)$. Then the properties that
    \begin{itemize}
        \item $\mathfrak{D}(V; G) = \text{\emph{III}}$,
        \item $\mathfrak{D}(V; G) = \text{\emph{IV}}$,
        \item $\mathfrak{D}(V; G) = \text{\emph{VII}}$,
        \item $\mathfrak{D}(V; G) = \text{\emph{X}}$,
    \end{itemize}
    depend only on the coupling of the Real structure $G$.
\end{cor}

\begin{proof}
    By Lemma \ref{lem:indepK4action}, $\Stab_{K_4}(V)$ depends only on the coupling. Furthermore, the C-type of $V$ being $\bR$, $\bC$ or $\bH$ depends only on $G$, and not on the Real structure at all. The D-types which are determined by these two pieces of data are those listed above, by Table 2 of \cite{RUMTAY3}.
\end{proof}

Recall that by Lemma \ref{lem:elementary1} any Real structure on $H$ is of the form $G = G(\alpha, h)$, where $\alpha \in \Aut(H)$ and $h \in H$ with $\alpha^2 = \Ad(h)$ and $\alpha(h) = h$. Let us fix such a model. Then a set of representatives of all Real structures on $G$ with the same coupling $\theta \coloneqq \overline{\alpha}$ as $G(\alpha, h)$ is $G(\alpha, zh)$, where $z$ varies $\HH^2(Z_{\theta})$ which we view as the quotient $HH^2(Z_{\theta}) \coloneqq \{z \in Z(H) \mid \alpha(z) = z\} / \{z\alpha(z) \mid z \in Z(H)\}$.

\begin{prop}\label{prop:variationofFS}
    Suppose that $H$ is finite and $V \in \Irr_{\bC}(H)$ with character $\chi$. Then
    \[
        \sF(V; G(\alpha, zh)) = \frac{\chi(z)}{\dim(V)} \cdot \sF(V; G(\alpha, h)).
    \]
    In particular, if $\mathfrak{D}(V; G(\alpha, h))$ is not equal to \emph{III, IV, VII} or \emph{X}, $\chi(z) \in \{\pm \dim(V) \}$. In this case, if $\chi(z) = \dim(V)$, then 
    \[
    \mathfrak{D}(V; G(\alpha, h)) = \mathfrak{D}(V; G(\alpha, zh)),
    \]
    and if $\chi(z) = -\dim(V)$, then $\{\mathfrak{D}(V; G(\alpha, h)), \mathfrak{D}(V; G(\alpha, zh))\}$ is one of $\{\text{\emph{I}}, \text{\emph{II}}\}$, $\{\text{\emph{V}}, \text{\emph{VI}}\}$, $\{\text{\emph{VIII}}, \text{\emph{IX}}\}$.
\end{prop}

\begin{proof}
    We can compute
    \begin{align*}
         |G| \cdot \sF(V; G(\alpha, h)) &= \sum_{g \in G(\alpha, zh) \setminus G} \chi(g^2), \\
         &= \sum_{g \in G} \chi((gz)^2), \\
         &= \sum_{g \in G} \chi(g \alpha(g) h).
    \end{align*}
    Replacing $h$ by $zh$, we then have that
    \begin{align*}
        |G| \cdot \sF(V; G(\alpha, zh)) &= \sum_{g \in G} \chi(g \alpha(g) h z), \\
        &= \frac{\chi(z)}{\dim(V)} \sum_{g \in G} \chi(g \alpha(g) h),
    \end{align*}
    using the fact that $\dim(V) \chi(gy) = \chi(g) \chi(y)$ for any central element $y$ of $G$, because $V$ is irreducible and so $y$ acts by a scalar. The remaining points then follow from Table 2 of \cite{RUMTAY3}.
\end{proof}

\begin{remark}
    Suppose that $\mathfrak{D}(V; G(\alpha, h))$ is not equal to III, IV, VII or X, and let $z \in Z(H)$. Because $[G(\alpha, h)] = [G(\alpha, z\alpha(z)h)]$, $\sF(V; G(\alpha, zh)) = \sF(V; G(\alpha, h))$ and therefore by the above we must have that $\chi(z \alpha(z)) = \dim(V)$, or equivalently that $\chi(z) \chi(\alpha(z)) = \dim(V)^2$. One can also see directly that this is true, because outside of D-types III, IV, VII or X, $\overline{\chi} = \chi(\alpha(-))$ because $c \in K_4$ stabilises $V$, hence $\chi(z) \chi(\alpha(z)) = |\chi(z)|^2 = \dim(V)^2$. 

    We also know that $\chi(z) = \pm \dim(V)$, by Proposition \ref{prop:variationofFS}, and therefore $\chi(z) = \chi(\alpha(z))$.
\end{remark}

This allows us to say identify immediately many $V \in \Irr_{\bC}(H)$ that will have S-type $\bC$.

\begin{cor}
    Suppose that $H$ is finite, $V \in \Irr_{\bC}(H)$, and $G$ is a Real structure on $H$ with coupling $\theta$. Then then if $\chi(z) \not\in \bR$ for some $z \in \HH^2(Z_{\theta})$, or $\chi(z) \neq \chi(\alpha(z))$ for some $z \in Z(H)$ and lift $\alpha$ of $\theta$, then $\sF(V, \widehat{G}) = 0$.
\end{cor}

%% file: chapter-physics.tex
\chapter{Wallpaper Groups: A User's Guide}\label{chap:physics}
Corepresentations are particularly important for the crystallographic groups, also known as the wallpaper groups.
In this chapter we present a guide into their linear and antilinear representations
in a style that makes it easier to read for the users of these groups: crystallographers, chemists, mathematicians, physicists.
We won't prove theorems in this chapter. Instead, we state facts, enumerated using uppercase Roman numerals and
present {\em justifications} why some of these facts hold.

\section{Crystallographic Groups} \label{section:algebra_cryst}
Consider the groups of affine transformations and euclidean isometries
\[ \Aff_n \coloneqq \GL_n (\bR)\ltimes \bR^n \, , \quad \Euc_n \coloneqq \HO_n (\bR)\ltimes \bR^n \]
of the euclidean space $(\bR^n,\phi)$. Note that the subgroups $\bR^n$ are canonical and normal: they consist of all the translations.
The subgroups $\GL_n (\bR)$  and $\HO_n (\bR)$ depend on the choice of an origin:
they consist of all the transformations preserving the origin.

A subgroup $G \leq \Euc_n$ is {\em crystallographic}, if it is {\em uniform} (the quotient $\bR^n/G$ is compact)
and {\em discrete} (as a subset of $\Euc_n$). When $n=2$, we refer to $G$ as {\em a wallpaper group}.

The fundamental facts about $G$ are the three {\em Bieberbach Theorems}:
\begin{enumerate}[label=\textnormal{(\Roman*)}]
\item {\em The translation group} $T\coloneqq G \cap \bR^n$ is a lattice, i.e., $T=\bZ^n = \sum_{i=1}^n \bZ e_i$, spanned over $\bZ$ by a basis $e_1, \ldots, e_n$ of $\bR^n$.
{\em The point group} $Q\coloneqq G/T$ is finite.  
\item If $G,G'$ are two crystallographic group of the same dimension $n$ and $\psi: G \to G'$ is an isomorphism of groups,
  then $\psi$ is a conjugation in $\Aff_n$, i.e., there exists $x\in\Aff_n$ such that $\psi (g) = xgx^{-1}$ for all $g\in G$.
\item Up to conjugation in $\Aff_n$, there exist only finitely many crystallographic groups of each dimension $n$. 
\end{enumerate}  

To work with the algebraic structure of $G$,
we choose a representative $\breve{q}\in G$
for each element $q\in Q$, 
The finite set of all these representatives $\breve{Q} \coloneqq \{ \breve{q}\mid q\in Q\}$ is called {\em the transversal}.
The point group $Q$ naturally acts on the lattice $T$ via conjugation: $q(t) \coloneqq \breve{q}t\breve{q}^{-1}$.

\begin{just}
  The action axiom is $q_1(q_2 (t)) = q_1q_2 (t)$. Its verification goes through because $T$ is abelian.
  Naturally in this context means that the action does not depend on the choice of origin and the choice of transversal.
  Both are true because the action can be written
as $q(t)(x) = q(t(x)) - q(x)$ for any $x\in \bR^n$.
\end{just}  

It is convenient to write all the matrices in the basis $e_1, \ldots, e_n$.
Let $\Phi = \begin{pmatrix}\phi (e_i,e_j)\end{pmatrix}$ be the Gram  matrix of this basis.
Then $Q$ is a subgroup of $\GL_n (\bZ)\cap \HO_n (\bR, \Phi) = \{ x\in \GL_n (\bZ) \mid x^{T}\Phi x=\Phi\}$.
The transversal element of $q\in Q$ can be written as $\breve{q} = [q|a]$, which is
a transformation $\breve{q} (x) = q(x) + a$.

Each point $x\in \bR^n$ gives also {\em the isotropy group} $G_x \coloneqq \{ g\in G \mid g(x) = x\}$, often called {\em the site-symmetry group}.
Each $G_x$ is naturally a subgroup of $Q$ under the map $G_x \to G \twoheadrightarrow Q$. Every isotropy group gives a semidirect product subgroup
$G_x \ltimes T$ inside $G$. The next statement is a standard algebra fact.
\begin{enumerate}[resume, label=\textnormal{(\Roman*)}]
\item The following five statements about $G$ are equivalent:
  \begin{itemize}
\item   $G$ is a semidirect product $Q\ltimes T$.
\item  The following exact sequence splits:
  \begin{equation} \label{sequence:TGQ}
    1\to T\to G \xrightarrow{\varpi} Q \to 1 \, .
    \end{equation}
\item  There exists a point $x\in \bR^n$ with $G_x=Q$.
\item One can choose a transversal $\breve{Q}$, which is a subgroup. 
\item With a suitable choice of the origin, one can choose a transversal $\breve{Q}$, where $\breve{q} = [q|0]$ for each $q\in Q$. 
\end{itemize}
  \end{enumerate}

Let us choose an origin from now on.
Then every element $g\in G$ can be written as $g =[q|a]$ for unique $q\in Q$ and $a\in \bR^n$. It is worth writing the formulas for multiplication and inverse in $G$ and the conjugation action on $T$ in these terms:
\begin{equation} \label{mult_wallpaper}
  [q_1|a_1][q_2|a_2] =[q_1q_2|q_1(a_2)+a_1]\, , \quad 
  [q|a]^{-1}= [q^{-1}|-q^{-1}(a)]\, , \quad  [q|a]t[q|a]^{-1} = q(t) \, .
\end{equation}  
\begin{just}
  The first formula is verified by applying both sides to $x\in\bR^n$.
The second and third formulas are deduced from the first one.
\end{just}  

We need a description of the centre $Z=Z(G) \coloneqq \{ g\in G \mid \, \forall h\in G \; gh=hg\}$.
Let us restrict the exact sequence~\eqref{sequence:TGQ} to the centre:
\begin{equation} \label{sequence:TGQ_centre}
  1 \to Z \cap T \to Z \xrightarrow{\varpi} Z(Q) \, .
\end{equation}  
\begin{enumerate}[resume, label=\textnormal{(\Roman*)}]
\item In the exact sequence~\eqref{sequence:TGQ_centre}, the kernel $Z\cap T$ is equal to $T^Q = \{ t\in T \mid \, \forall q\in Q \; q \cdot t = t \}$
  and the image of $\varpi$ is trivial.  Therefore, $Z(G) =T^Q$.
\end{enumerate}
\begin{just} The exact sequence is the restriction of the map $\varpi: G \to Q$ to $Z$. Write $t\in T$ as $[I|t]\in T$ and use~\eqref{mult_wallpaper}:
  \[ [q|x][I|t] = [I|t][q|x] \ \Leftrightarrow \  [q| q(t)+x]= [q| t+x] \ \Leftrightarrow \ q(t) = t \, .\]
  This proves both facts. For $t\in T$ to belong to $Z$, we need $q(t)=t$ for all $q\in Q$. And for
$[q|x]\in Z$ we need $q(t)=t$ for all $t\in T$, forcing $q=I$.  
\end{just}

Another algebraic piece of the puzzle is {\em the original cocycle}, 
    \[ \beta_{o} \in \HZ^1 (Q,\bR^n/T), \quad  [ q | \beta_o (q) ] = \varpi^{-1} (q) = G \cap [ q | \bR^n ]\, . \]
    To understand it better, choose an origin $o$ and write the transversal elements explicitly as $\breve{q} = [q| a_q]$.
    Then $\beta_o (q) = a_q+T$, so that $\beta_o$ essentially describes the fibres of the map $\varpi : G \to Q$.
    It depends on the origin but it becomes canonical, once an origin is chosen. Hence the name.

    Notice that the original cocycle maps to
    the 2-cohomology responsible for the original extension~\eqref{sequence:TGQ} under the maps
    \[ \HZ^1 (Q,\bR^n/T) \twoheadrightarrow \HH^1 (Q,\bR^n/T) \xrightarrow{\cong} \HH^2 (Q,T) \, . \]

    Finally consider {\em the original lattice} $T_o$, generated by $T$ and all $a_q$, $q\in Q$.
    It is always a free abelian group but it may or may not be a discrete subgroup $\bR^n$. It depends on the origin but does not depend on the choice of transversal.
  \begin{enumerate}[resume, label=\textnormal{(\Roman*)}]
  \item The original lattice is a $Q$-submodule of $\bR^n$ containing $T$. One can choose an origin so that
    $T_o \subseteq \frac{1}{|Q|} T$.
  \end{enumerate}
  \begin{just}
    By definition, $T_o$ is an abelian group and $T\subseteq T_o \subset \bR^n$. Pick $s,t\in Q$. By~\eqref{mult_wallpaper},
    $q(a_s)+a_q \in T_o$ for all $q,s\in Q$. Thus, $q(a_s) \in T_o$.

    The second statement follows immediately from the transfer principle in group cohomology but we give an elementary recipe.
    Start with any origin $o$, then move it to $o+ \frac{1}{|Q|} \sum_{q\in Q} a_q$. The reader can easily check that the new origin satisfies
    the requirements.
\end{just}

Notice that $T_o$ is not a $\GL_n(\bZ)$-module, in general. But $\widetilde{T}\coloneqq \frac1{|Q|}T$ is a $\GL_n(\bZ)$-module.
  Let us call an origin {\em rational} if $T_o \subseteq \widetilde{T}$. With a rational origin, we can always realise the group $G$ as a finite index subgroup in a split group
  \[ G \leq Q \ltimes T_o \leq Q \ltimes \widetilde{T} \, . \]

\section{Automorphism Groups}

Let us fix a rational origin $o$ throughout this section.

The second Bieberbach Theorem yields a description of the automorphism groups $\Aut (G)$ and $\Out (G)$ via
the normaliser $N_{\Aff_n} (G) = \{ x \in \Aff_n \mid \forall g\in G \ xgx^{-1}\in G \}$ and
the centraliser $C_{\Aff_n} (G) = \{ x \in \Aff_n \mid \forall g\in G \ xgx^{-1}=g \}$:
\begin{equation} \label{formula:AUT1}
  \Aut (G)  = N_{\Aff_n} (G) /C_{\Aff_n} (G), \quad
  \Out (G)  = N_{\Aff_n} (G) /G \cdot C_{\Aff_n} (G).
\end{equation}  

We would like another, more efficient description.
Note that $\varphi (T) =T$ for any automorphism $\varphi \in \Aut (G)$ as follows from the second Bieberbach Theorem.
Consequently, any automorphism $\varphi \in \Aut (G)$ has three components
\begin{equation} \label{formula:AUT0}
  \varphi [q|x] = [ \varphi_Q (q) | \varphi_T (x) + \varphi_m (q,x)] \, , 
\end{equation}
where $\varphi_Q\in \Aut (Q)$ and the remaining components
$\varphi_T$ and $\varphi_m$ may be represented by an element $\Aut (T_o) = \GL_n (\bZ)$ and a map $Q \times T_o \to T_o$.
Let us make it precise. 
  \begin{enumerate}[resume, label=\textnormal{(\Roman*)}]
  \item For every $\varphi \in \Aut (G)$ there exist unique $A\in N_{\GL_n (\bZ)}(Q)$ and $\alpha \in \HZ^1 (Q,T_o)$
  such that   
\begin{equation} \label{formula:AUT}
  \varphi [q|x] = [ AqA^{-1} | A(x + \alpha (q))] \, .
\end{equation}
  \end{enumerate}
  \begin{just}
    By the second Bieberbach Theorem, $\varphi$ is a conjugation by $[A|a]\in \Aff_n$.  Using~\eqref{mult_wallpaper}, we can write
    \[
    \varphi [q|x] =
            [ A|a] [q|x][A|a]^{-1}  =
            [ Aq | Ax +a][A^{-1}| - A^{-1}a]  =
                       [ AqA^{-1} |Ax +a - AqA^{-1}a] \, .
    \]
    Since $\phi [I|t] = [I|At]$ for all $t\in T$, we conclude that $A\in \GL_n (\bZ)$.
    Since $AqA^{-1}\in Q$ for all $q\in Q$, $A$ must actually belong to the normaliser $N_{\GL_n (\bZ)}(Q)$.
    Finally, the function 
    \[ \alpha : Q \to \bR^n, \quad \alpha (q) \coloneqq A^{-1}a - qA^{-1}a\]
    takes values in $T_o$. It is a 1-cocycle, by a direct calculation, i.e.,  $d \alpha \equiv 0_T$ where
    \[d \alpha (q,s) \coloneqq q\cdot \alpha (s) - \alpha (qs)+ \alpha (q) \, . \]
    Finally, suppose $(A_1,\alpha_1)$ and $(A_2,\alpha_2)$ represent the same automorphism.
    Then
    \[
    (A,\alpha) \coloneqq (A_1,\alpha_1)^{-1}\cdot(A_2,\alpha_2) = (A_1^{-1}A_2, A_1^{-1}\alpha_2 - A_1^{-1}\alpha_1)^{-1} \in N_{\GL_n (\bZ)}(Q) \ltimes \HZ^1 (Q,T_o)
    \] 
    represents the trivial automorphism $\varphi$. Applying~\eqref{formula:AUT}, $A(x+\alpha(I)) = x+\alpha(I)$ for all $x\in T$.
    Thus, $A=I$.
    Applying~\eqref{formula:AUT} again, $a_q+\alpha(q) =a_q$ for all $q\in Q$. Hence, $\alpha =0$. 
  \end{just}
  
  So far  we have realised $\Aut (G)$ as a subgroup of the semidirect product
  $N_{\GL_n (\bZ)}(Q) \ltimes \HZ^1 (Q,T_o)$. To describe it more precisely,
  we need little more algebraic yoga. Notice that the original cocycle $\beta_o$ has values in $T_o/T$.
  We can think that $\beta_o \in \HZ^1 (Q, T_o/T)$.
  We want to push $\alpha \in \HZ^1 (Q, T_o)$ forward to the same group:
  let $\overline{\alpha} \coloneqq (T_o \to T_o/T)_\ast (\alpha) \in \HZ^1 (Q, T_o/T)$.
        Finally, we need to pinpoint the following right $N_{\GL_n (\bZ)}(Q)$-action on $\HZ^1 (Q, \widetilde{T}/T)$
    \[ \gamma^A : s \mapsto A^{-1} \cdot \gamma (AsA^{-1}) \, . \]
Under this action, $\HZ^1 (Q, T_o/T)$ is a $Q$-submodule of $\HZ^1 (Q, \widetilde{T}/T)$.
Observe that the induced $Q$-action on $\HH^1 (Q, \widetilde{T}/T)$ is trivial.
Therefore,  $N_{\GL_n (\bZ)}(Q, o)$,
the stabiliser of the original cohomology class $[\beta_o]\in \HH^1 (Q,\widetilde{T}/T)$ inside the normaliser $N_{\GL_n (\bZ)}(Q)$,
contains $Q$.
\begin{enumerate}[resume, label=\textnormal{(\Roman*)}]
\item
  A pair $(A,\alpha)\in N_{\GL_n (\bZ)}(Q)\times \HZ^1 (Q,T_o)$ defines an automorphism $\varphi = (A,\alpha)$ by the formula~\eqref{formula:AUT}
  if and only if
  \[A \in N_{\GL_n (\bZ)}(Q, o) \ \mbox{ and } \ \beta_o^A = \beta_o + \overline{\alpha} \, . \]
   Moreover,
    the full automorphism group fits into a short exact sequence 
    \begin{equation} \label{Physics_AutSequence}
      1 \longrightarrow \HZ^1 (Q,T) \longrightarrow \Aut (G) \longrightarrow   N_{\GL_n (\bZ)}(Q,o) \longrightarrow 1\, .
    \end{equation} 
\end{enumerate}

\begin{just}
  Pick $(A,\alpha)\in N_{\GL_n (\bZ)}(Q)\times \HZ^1 (Q,T_o)$.
  A straightforward algebraic check tells that  the formula~\eqref{formula:AUT} defines an injective homomorphism $G\to \Aff_n$.
  Thus, we just need to control its image. Using the transversal, observe that
  $(A,\alpha)$ being an automorphism is equivalent to the condition
  \[ \forall q\in Q \quad (A,\alpha) [q| a_q] = [ AqA^{-1} | A(a_q + \alpha (q)) ] \in G \, . \]
  Clearly, this condition is equivalent to another condition with the original cocycle
  \[ \forall q\in Q \quad Aa_q +A \alpha(q)\in  a_{AqA^{-1}} + T  = \beta_o (AqA^{-1}) \, . \]
  Modulo $T$, this is equivalent to
  \[a_q + \alpha(q)+ T = A^{-1}\beta_o (AqA^{-1}) \quad \mbox{or} \quad
\beta_o + \overline{\alpha} = \beta_o^A \, ,\]
as stated. Now observe that $\overline{\alpha}$ is a 1-coboundary. Hence, 
$[\beta_o]  = [\beta_o]^A$ and $A\in N_{\GL_n (\bZ)}(Q,o)$. It follows that
\[
f: (A,\alpha) \mapsto A, \quad \Aut (G) \longrightarrow   N_{\GL_n (\bZ)}(Q,o) 
\]
is a surjective group homomorphism. Its kernel consists of
of pairs $(I,\alpha)$. But $(I,\alpha)$ defines an automorphism of $G$ if and only if $\overline{\alpha}=0$. In other words,
$\alpha \in \HZ^1 (Q,T)$. Therefore,  $\HZ^1 (Q,T)$ is the kernel of $f$.
\end{just}

As observed above, $Q$ is a subgroup of $N_{\GL_n (\bZ)}(Q, o)$.
It must be normal since $Q$ is normal in $N_{\GL_n (\bZ)}(Q)$.
The group of inner automorphisms $\Inn (G)$ is the image of $G$ inside $\Aut (G)$.
It can be described in terms of the exact sequence~\eqref{Physics_AutSequence} as
\begin{equation} \label{Physics_InnSequence}
  1 \longrightarrow \HB^1 (Q,T) \longrightarrow \Inn (G) \longrightarrow  Q \longrightarrow 1\, .
\end{equation} 
This yields a short exact sequence that describes $\Out (G)$ 
    \begin{equation} \label{Physics_OutSequence}
      1 \to \HH^1 (Q,T) \to \Out (G) \to N(Q, o) \to 1\, , \ \mbox{ where } N(Q, o) \coloneqq N_{\GL_n (\bZ)}(Q, o)/Q \, .
      \end{equation}
    If the sequence~\eqref{sequence:TGQ} splits, then three simplifications take place:
    $\beta =0$, $N_{\GL_n (\bZ)}(Q, o) = N_{\GL_n (\bZ)}(Q)$ and the sequence~\eqref{Physics_OutSequence} also splits.
We tabulate all the necessary data in the 2-dimensional case in Appendix~\ref{appendix_wallpaper}.

\section{Magnetic Wallpaper Groups} \label{sect:are_there}

A \emph{magnetic group} is subgroup $M \subset \Euc_n \times C_2$, such that the image $\pi_1(M)$ under projection onto the first factor $\Euc_n$ is crystallographic. Two magnetic groups are considered \emph{equivalent} if conjugation by an element of $\Aff_n$ 
takes one to the other.

Now magnetic groups come in three classes, corresponding to the three possible classes of subgroup of a direct product $\Euc_n \times C_2$:
\begin{itemize}
  \item Ordinary crystallographic groups, $M = G \times 1 < \Euc_n \times C_2$,
  \item ``Grey'' groups, $M = G \times C_2 < \Euc_n \times C_2$, where $G$ is crystallographic, and
  \item ``Black-and-white'' or ``two-colour'' groups
    $M = \{(g, \phi(g)) \mid g \in G\}$
  for some surjective homomorphism $\phi \colon G \twoheadrightarrow C_2$.
\end{itemize}

Observe that the data which determines a black-and-white group is a pair $(G, \phi \colon G \twoheadrightarrow C_2)$ where $G$ is crystallographic, which is the same as the data of the pair $(G,\ker(\phi))$, i.e. a crystallographic group and an index two subgroup. Why are they called black-and-white? If $G$ is $2$-dimensional, a wallpaper group, then $G$ is the symmetry group of some wallpaper pattern. The map $\phi$ corresponds to a colouring of this pattern in black and white so that $g$ preserves the colouring if $\phi(g) = 1$ and swaps the colours if $\phi(g) \neq 1$, see \cite{Schat}.

We can also give a more algebraic point of view, which makes them easier to classify. Let us define a \emph{Magnetic group} (with a capital $M$) to be a $C_2$-graded group (meaning a group $G$ equipped with a homomorphism $\varphi \colon G \rightarrow C_2$) such that the ungraded subgroup $H = \ker(\varphi)$ is crystallographic. Two Magnetic groups are the same (are \emph{equivalent}) if they are isomorphic as $C_2$-graded groups: there is some group isomorphism $\psi \colon G_1 \rightarrow G_2$ such that $\psi(H_1) = H_2$.

From a magnetic grey or black-and-white group $M \subset \Euc_n \times C_2$, one obtains a Magnetic group $(M, \pi_2 \colon M \rightarrow C_2)$ via projection to the second factor $C_2$. The key property is that every Magnetic group is magnetic.

\begin{enumerate}[label=\textnormal{(\Roman*)}]
\item
  Suppose $H$ is a crystallographic group and $H \leq G$ is a Real structure. Then $G$ is either grey, i.e.,
  $G \cong H \times C_2$, or black-and-white, i.e., $G$ is also crystallographic (in the way compatible to the crystallographic structure of $H$) of the same dimension as $H$.
\end{enumerate}

Note that $G$ cannot be both grey and black-and-white: $H \times C_2$, for $H$ crystallographic, has torsion elements in $Z(H \times C_2)$, which is impossible for crystallographic groups, cf. Fact~\ref{section:algebra_cryst}(V).

\begin{just}
  Let $T\subset H$ be the translation lattice and $Q$ the point group of $H$.
Set $V\coloneqq (\bR^n)^Q$. Observe first that
\[
C_{\Aff_n}(H) \;=\; C_{\Euc_n}(H) \;=\; \{\,[I|v] \mid v\in V\,\}.
\]
Indeed, if $[A|b]$ commutes with every $[I|t]$, $t\in T$, then $At=t$ for all $t\in T$; since $T$ spans $\bR^n$, $A=I$.
Then $[I|b]$ commutes with $[B|c]\in H$ if and only if $Bb=b$. In particular, $C_{\Aff_n}(H)$ is an $\bR$-vector space, so for any $C_2$-action $\HH^1 (C_2, C_{\Aff_n}(H))=0$ -- this is the heart of the proof below.

Pick $g\in G\setminus H$ and consider the automorphism $\gamma_g\in\Aut(H)$, $\gamma_g(h)=ghg^{-1}$.
By the second Bieberbach theorem, $\gamma_g=\gamma_x$ for some $x\in\Aff_n$.
This $x$ is unique up to right multiplication by $C_{\Aff_n}(H)$: one may replace it by $x[I|v]$ for any $v\in V$.
Let $A$ be the linear part of $x$. Since $x$ normalises $H$, conjugation by $x$ preserves $C_{\Aff_n}(H)$, so $AV=V$, and
\[
[I|v]\,x \;=\; x\,[I|A^{-1}v] \qquad (v\in V).
\]
We hope to choose $v$ so that
\[
\psi : G\to\Aff_n,\qquad \psi(h)=h,\quad \psi(hg)=h\,x[I|v]\qquad (h\in H)
\]
is a homomorphism. Note that $G$ is generated by $H$ and $g$ subject to the following relations:
\[ h_1h_2 = h_1h_2 \in H, \ \ ghg^{-1}=\gamma_g(h) \ \mbox{ and } \ g^2=g^2\in H \, . \]
Since $\gamma_{x[I|v]}=\gamma_x=\gamma_g$ on $H$,
$\psi$ is a well-defined homomorphism if and only if
\[
g^2 \;=\; (x[I|v])^2 \;=\; x^2\,[I|A^{-1}v+v]\,.
\]
Now $g^2$ and $x^2$ induce the same automorphism of $H$, so $x^{-2}g^2\in C_{\Aff_n}(H)$, i.e.\ $g^2=x^2[I|w]$ for some $w\in V$.
Since $g^2\in H$ commutes with $g$ and $\gamma_x=\gamma_g$, we have $xg^2=\gamma_x(g^2)\,x=\gamma_g(g^2)\,x=g^2x$, hence
\[
g^2[I|-w]\,x \;=\; x^2\cdot x \;=\; x\cdot x^2 \;=\; x\,g^2[I|-w] \;=\; g^2\,x\,[I|-w]\,,
\]
so $[I|-w]x=x[I|-w]$, i.e.\ $Aw=w$.
Let $v\coloneqq\tfrac12 w$. Then $A^{-1}v=v$ and
\[
(x[I|v])^2 \;=\; x^2\,[I|A^{-1}v+v] \;=\; x^2\,[I|w] \;=\; g^2 ,
\]
so $\psi$ is a well-defined homomorphism, injective on $H$.

Let $K\coloneqq\ker\psi$. Since $K\cap H=1$ and $[G:H]=2$, either $K=1$ or $K=\langle k\rangle$ has order $2$.

If $|K|=2$, then $k\notin H$ and, being a normal subgroup of order $2$, $\langle k\rangle$ is central; hence
$G=H\times\langle k\rangle\cong H\times C_2$ with the $C_2$-factor in $G\setminus H$, i.e.\ $G$ is grey.

If $K=1$, then $\Gamma\coloneqq\psi(G)\subset\Aff_n$ contains $T$ as a normal subgroup of finite index.
Its group $Q_\Gamma\subset\GL_n(\bR)$ of linear parts is finite, since $T$ lies in the kernel of the linear-part map.
Averaging, $\langle u,v\rangle'\coloneqq\sum_{B\in Q_\Gamma}\langle Bu,Bv\rangle$ is a $Q_\Gamma$-invariant inner product, and so there is some $S\in\GL_n(\bR)$ with $\langle u,v\rangle'=\langle Su,Sv\rangle$. Then $SQ_\Gamma S^{-1}\subset\mathrm{O}(n)$, and since $[S|0][B|b][S|0]^{-1}=[SBS^{-1}|Sb]$, the group $\Gamma' \coloneqq [S|0]\,\Gamma\,[S|0]^{-1}$ lies in $\Euc_n$. As the crystallographic subgroup $H' = [S|0]\,\Gamma\,[S|0]^{-1}$ has index two in $\Gamma'$, we see that $G$ is crystallographic of the same dimension as $H$.
\end{just} 

The following is easy to see using (I) above and Bieberbach (II).

\begin{enumerate}[resume, label=\textnormal{(\Roman*)}]
\item
  Up to equivalence, every Magnetic group $(G, \varphi \colon G \rightarrow C_2)$ is of the form $(M, \pi_2)$ for some magnetic group $M$.
  If two magnetic groups $M, M' \subset \Euc_n \times C_2$ are equivalent, they are of the same class.
  Two ordinary magnetic groups $M, M' \subset \Euc_n \times 1$ are equivalent
  if and only if they are isomorphic.
Two grey or black-and-white magnetic groups $M, M' \subset \Euc_n \times C_2$ are equivalent
if and only if there exists an isomorphism $\psi: M \to M'$ that agrees with the second projection, that is,
$\pi_2' \psi = \pi_2$.
\end{enumerate}

The main consequence is our notation of a Magnetic group contains is essentially the same as that of a magnetic group, and this makes the classification of magnetic groups more of an algebraic problem, for example as one can now use the theory from Chapter \ref{chapter:realstructures}.

We are interested in magnetic groups in dimension $2$ -- how many are there? It was shown by Heesch \cite{Heesch} in 1929 that there are 80. See also \cite{Schat} for a mathematical treatment and \cite{PaZa, Souv, Zamo} for generalisations. Similar counts are known in dimensions 3 and 4 \cite{Souv} but we won't address them here.

In dimension 2 the count goes like this:
  \begin{itemize}
  \item 17 ordinary wallpaper groups $G$,
  \item 17 grey groups $G\times C_2$ for each wallpaper group $G$,
  \item 46 black-and-white groups $H\leq G$, where $G$ is a wallpaper group, $H$ an index two subgroup.
  \end{itemize}

We are interested in the corepresentations of magnetic groups, and as we need to work with these groups explicitly, we review this classification in Appendix \ref{appendix_wallpaper}. In fact, we derive this classification in two ways: by counting the index two subgroups, and by the more abstract classification of Real structures on wallpaper groups, using the results of Chapter \ref{chapter:realstructures} and the Magnetic groups point of view.

The interesting case are the black-and-white groups. We tabulate these in Table \ref{table:wallpaper3}.

\section{Classification of Irreducible Corepresentations}\label{sect:cliffordtheory}
Consider a crystallographic group $H$ and a Real structure $H\leq G$.
In this section we give a classification of irreducible corepresentations of $G$ via
a version of Clifford Theory. 

Let us do a quick recap of its classical version, where $K$ is a group and $N$ is its normal subgroup.
There is a left action of $K$ and $K/N$ on the set ${\Irr_{\bC}(N)}$ by conjugation:
\begin{equation} \label{Cliff_1a}
  {}^{k}(V,\rho) \coloneqq (V, n \mapsto \rho (k^{-1}nk)) \ \mbox{ where } \ n\in N, \, k\in K.
  \end{equation} 

Clifford Theory can be viewed as a correspondence between certain subsets of ${\Irr_{\bC}(K)}$ and the $K$-orbits of ${\Irr_{\bC}(N)}$ along the following lines.

  \begin{enumerate}[label=\textnormal{(\Roman*)}]
  \item Given $V\in {\Irr_{\bC}(K)}$, its restriction to $N$ is a sum of irreducible representations over the orbit of a single $W\in {\Irr_{\bC}(N)}$:
\begin{equation} \label{Cliff_1b}
  V|_N \cong \left(\oplus_{U \in {}^{K}W} U\right)^{\oplus m}\, . 
\end{equation} 
\item Given $(W,\rho) \in {\Irr_{\bC}(N)}$, consider the stabiliser of its isomorphism class $K_W$. Note that $K_W$ contains $N$:
$\rho(n)$ is an isomorphism ${}^n W \to W$. Now for each $k\in K_W \setminus N$, pick an isomorphism of $N$-modules $\rho(k) : {}^k W \to W$.
The newly extended map $\rho : K_W \to \GL (W)$ is a projective representation, in general, defining
a cocycle and the cohomology class:
\begin{equation} \label{Cliff_1c}
d\rho \in \HZ^2 (K_W , \bC^\times), \ d \rho (t,s) = \rho (t) \rho (s) \rho (ts)^{-1}, \quad
c_W = [ d \rho]^{-1} \in \HH^2 (K_W , \bC^\times) \, .
\end{equation}  
\item The cohomology class $c_W$ does not depend on the choice of $\rho$, only on $W$.
It represents an obstruction to turning $W$ into a representation of $K_W$.  
In particular, its restriction to $N$ is trivial.
This may or may not lift\footnote{There is no naive exact sequence connecting cohomology of $N$, $K_W$ and $K_W/N$. The relations is governed by the Lyndon-Hochschild-Serre spectral sequence. This yields the seven-term inflation-restriction sequence, where one can see the obstructions to existence and uniqueness of $\overline{c_W}$.}
to a class $\overline{c_W} \in  \HH^2 (K_W/N , \bC^\times)$.
\item Suppose such a lift $\overline{c_W}$ exists (which in general will not be unique).
  The key set of representations is ${\Irr_{W}(K_W/N)}$, the set of irreducible finite-dimensional $\overline{c_W}$-projective representations of $K_W/N$. For any $X\in {\Irr_{W}(K_W/N)}$, the tensor product $W\otimes X$
  is a linear representation of $K_W$ -- the cocycles of the two tensor multiplicands cancel each other out.
\item Finally given $W$ and $\overline{c_W}$ as above, there is a one-to-one correspondence between
  ${\Irr_{W}(K_W/N)}$ and those $V\in {\Irr_{\bC}(K)}$ that contain $W$ in their restriction to $N$, defined by
\begin{equation} \label{Cliff_1d}
  X \mapsto \Ind_{K_W}^K W \otimes X, \, \qquad V \mapsto \Hom_N(W,V), \,
\end{equation}  
where $\Hom_N(W,V)$ is a $c_W$-projective representation of $K_W/N$ via $k * f \coloneqq \rho_V(k) \circ f \circ \rho(k)^{-1}$.
\end{enumerate}

Let us now readjust everything to Real structures and corepresentations. We have a crystallographic group $H$ sitting inside a Real structure $G$ and $T$ is the lattice of $H$. The pair $(K,N)$ above will be replaced with the pair $(G,T)$. Note that the translation subgroup $T$ remains normal in $G$ because $G$ acts on $H$ by automorphisms and each automorphism preserves translation by the second Bieberbach Theorem.

We fix a transversal $\breve{S}$ of the finite group $S \coloneqq G/T$, a section of the map $G \twoheadrightarrow S$, such that $\breve{1_S} = 1_G$. Because $T$ is ungraded, the group $S$ inherits the grading $\pi : T \to C_2$ from $G$.
Since $T$ is abelian, all representations ${\Irr_{\bC}(T)}$ are one dimensional, i.e.,
they are characters, elements of $T^\vee \coloneqq \hom (T,\bC^\times)$.
The group $S$ acts on $T$ by conjugation $t^s = \breve{s}^{-1} t \breve{s}$,
which translates into {\em the left Real action} on $T^\vee$: for $s\in S$ and $\chi\in T^\vee$,
\begin{equation} \label{Cliff_2a}
{}^{s}\chi (t)  \coloneqq {}^{\pi (s)} \chi (\breve{s}^{-1}t\breve{s}) = \begin{cases}
    \overline{\chi (\breve{s}^{-1}t\breve{s})}  , & \text{if $s$ odd,}\\
    \chi (\breve{s}^{-1}t\breve{s}) , & \text{if $s$ even.}\\
    \end{cases}
  \end{equation}  
This is compatible with the natural action of $S$ on $\Irr_{\bC}(T)$ as considered in Chapter \ref{chapter:twelvefoldway}, i.e.
\[
  s * \bC_{\chi} = \bC_{{}^{s} \chi},
\]
where $\bC_{\chi} \in \Irr_{\bC}(T)$ is $\bC$ with action of $T$ by multiplication by $\chi$.

Let us now describe Clifford Theory in this context, which allows one to describe ${\Irr_{\bC}^{\hspace{0.1em} \star}(G)}$ in terms of the $S$-orbits on $T^\vee$ and projective antilinear representations of quotients of $S$.

  \begin{enumerate}[resume, label=\textnormal{(\Roman*)}]
  \item Given $V\in {\Irr_{\bC}^{\hspace{0.1em} \star}(G)}$, its restriction to $T$ is a sum of characters over a single orbit, i.e.,
\begin{equation} \label{Cliff_2b}
  V|_T \cong (\oplus_{\chi \in S \cdot \mu} \bC_\chi) ^{\oplus m}\, .
\end{equation}
\end{enumerate}

\begin{just}
  Elements of $T$ give a commuting set of linear operators on the finite dimensional $\bC$-vector space $V$.
  Hence, they have a common eigenvector, which spans $\bC_\chi$ for some $\chi\in T^\vee$.
  Then $V_\chi$, the sum of all $T$-submodules of $V$ isomorphic to $\bC_\chi$, is not zero. Pick $s\in S$ and observe the $T$-action on $\breve{s}V_\chi$:
\begin{equation} \label{Cliff_3a}
  t\breve{s} v = \breve{s}(\breve{s}^{-1}t\breve{s})v = \breve{s} \chi (\breve{s}^{-1}t\breve{s}) v =
  \begin{cases}
    \overline{\chi (\breve{s}^{-1}t\breve{s})}   \breve{s}v , & \text{if $s$ odd,}\\
    \chi (\breve{s}^{-1}t\breve{s})          \breve{s}v , & \text{if $s$ even.}\\
  \end{cases}
\end{equation}
This justifies~\eqref{Cliff_2a}. Furthermore, $\sum_{s\in S} \breve{s}V_\chi$ is a non-zero  $H$-submodule. Since $V$ is irreducible,
it must be $V$: this justifies~\eqref{Cliff_2b}.
\end{just}

\begin{enumerate}[resume, label=\textnormal{(\Roman*)}]
\item Given a character $\chi\in T^\vee$, we need its stabilisers $S_\chi$ and $G_\chi$ under the action~\eqref{Cliff_2a}.
Similarly to~\eqref{Cliff_1c}, we have a cochain and the cocycle, which, in our case, can be always reduced to the smaller group $S_\chi$ and found by an explicit formula. To do that, we consider the cocycle defining the product on the transversal, $\gamma \in \HZ^2 (S, T)$, which is defined by the formula
\begin{equation} \label{Cliff_2c_0}
\gamma (h,k) \cdot \breve{g}  = \breve{h} \cdot \breve{k}
\end{equation}
for $h,k\in S$, where $g \coloneqq hk$. Then $d\chi$ defined by
\begin{equation} \label{Cliff_2c}
d\chi \in \HZ^2 (S_\chi , \bC^\times), \  
d \chi (h,k) = \chi (\gamma (h,k))^{-1}, \quad
c_\chi \coloneqq [ d \chi]^{-1} \in \HH^2 (S_\chi , \bC^\times)\, ,
\end{equation}  
is well-defined, where $\bC^{\times}$ is \emph{not the trivial $S_{\chi}$-module} but rather has the action where elements of $(S_{\chi})_0$ act as the identity and all other elements of $S_{\chi}$ act by complex conjugation. Of course, if $S_{\chi}$ happens to be trivially graded, then $\bC^{\times}$ will be the trivial $S_{\chi}$-module.
\item In analogy with the classical picture, taking $W = \bC_{\chi}$, so $G_{\chi} = K_W$, the image of $c_\chi$ under the inflation map $\HH^2(S_\chi, \bC^{\times}) \rightarrow \HH^2(G_{\chi}, \bC^{\times})$ is the the class $c_W$ defined by
\begin{equation}\label{Cliff_1c_new}
d\rho \in \HZ^2 (G_{\chi} , \bC^\times), \ d \rho (t,s) \rho (ts) = \rho (t) \rho (s), \quad
c_W = [ d \rho]^{-1} \in \HH^2 (G_{\chi} , \bC^\times) \, 
\end{equation}  
where for each $g \in G_{\chi}$ we have chosen an isomorphism $\rho(g) \colon g * W \xrightarrow{\sim} W$, and similarly $\bC^{\times}$ has this potentially non-trivial action of $G_{\chi}$.
\item The class $c_\chi$ does not depend on the choice of transversal, only on the character $\chi$.
  If $G_\chi$ is a semidirect product of $S_\chi \ltimes T$,
    then the class $c_\chi \in \HH^2 (S_\chi , \bC^\times)$ is automatically trivial. 
\end{enumerate}

\begin{just}    
  To see that~\eqref{Cliff_2c} is well-defined, note that because $\chi$ is fixed by $S_{\chi}$, $\chi \colon T \rightarrow \bC^{\times}$ can be considered as a homomorphism of $S_{\chi}$-modules with respect to conjugation on $T$ and our given action on $\bC^{\times}$. Thus $d \chi$ is just the image of $\gamma$ under of the map $ \chi \circ (-) \colon \HC^2(S_{\chi},T) \rightarrow \HC^2(S_{\chi},\bC^{\times})$.

  To see that the image of $c_{\chi}$ is $c_W$, observe that the key simplification compared to~\eqref{Cliff_1c}: the representations of $T$ are one-dimensional, so we can choose $\rho (g): g * W \to W$ to be essentially the identity maps. More precisely, for $t \in T$ and $s \in S_{\chi}$ we define
  \begin{equation} \label{Cliff_2g}
    \rho (t \breve{s}) (a) \coloneqq \chi(t)^{-1} \pi(s)(a) = \chi(t)^{-1} {}^{s}a
 =
  \begin{cases}
    \chi(t)^{-1} a , & \text{if $s$ even,}\\
    \chi(t)^{-1} \overline{a} , & \text{if $s$ odd.}\\
  \end{cases}
  \end{equation}
Then for $t_1 \breve{s_1}, t_2 \breve{s_2} \in G_{\chi}$ we can compute
\[
\rho(t_1 \breve{s_1}) \rho(t_2 \breve{s_2}) = \chi(t_1)^{-1} \pi(s_1) \circ \chi(t_2)^{-1} \pi(s_2) = \chi(t_1)^{-1} {}^{s_1}(\chi(t_2)^{-1}) \pi(s_1s_2),
\]
and
\[
\rho(t_1 \breve{s_1} t_2 \breve{s_2}) = \rho(t_1 \breve{s_1}t_2\breve{s_1}^{-1} \breve{s_1} \breve{s_2}) = \chi(t_1)^{-1} \chi( \breve{s_1}t_2\breve{s_1}^{-1})^{-1} \chi(\gamma(s_1,s_2))^{-1} \pi( s_1 s_2).
\]
Because $\chi$ is fixed by $s_2$, ${}^{s_1}(\chi(t_2)) = \chi( \breve{s_1}t_2\breve{s_1}^{-1})$, and thus
\[
d \rho(t_1 \breve{s_1}, t_2 \breve{s_2}) = \chi(\gamma(s_1,s_2))^{-1},
\]
which is exactly the pullback of $d \chi$. Independence of the choice of transversal follows from the independence of the choice of transversal of $\gamma \in \HZ^2 (S, T)$. Similarly, if $S_\chi$ is a semidirect product, the transversal is the group $S_\chi$ itself, forcing $\gamma$, and thus $d\chi$ and $c_\chi$, to be trivial.
\end{just}    
  
\begin{enumerate}[resume, label=\textnormal{(\Roman*)}]  
\item The key set of representations is ${\Irr_{\chi}^{\hspace{0.1em} \star}(S_\chi)}$, the set of irreducible finite dimensional antilinear $c_{\chi}$-projective representations of $S_\chi$. These are vector spaces $(V,\rho)$ equipped with $\pi(g)$-linear operator $\rho(g)$ for each $g\in S_\chi$ such that
\begin{equation} \label{Cliff_2d}
\chi(\gamma (h,k)) \rho (hk) = \rho(h) \rho(k), \ \mbox{ for all } \ h,k \in S_\chi \, .
\end{equation}
If $\chi(\gamma (h,k))$ is trivial, then these are just antilinear representations of $S_\chi$. In general, these can be viewed as antilinear representations of another finite group $\widetilde{S_\chi}$, a central extension of $S_\chi$ (see Section \ref{sect:projantireps}), and so we can understand ${\Irr_{\chi}^{\hspace{0.1em} \star}(S_\chi)}$ using facts about antilinear representations of finite groups.
\item Given $\chi \in T^\vee$ and $(W,\sigma)\in {\Irr_{\chi}^{\hspace{0.1em} \star}(S_\chi)}$,
the tensor product $(\bC_\chi \otimes W,\mu)$ over $\bC$ is an antilinear representation of $G_\chi$ by
\begin{equation} \label{Cliff_2e}
\mu (t \breve{s}) ( 1 \otimes w) = \chi(t) \otimes \sigma(s)(v), \mbox{ where } t\in T, s\in S_\chi, w\in W\, .
\end{equation}
\end{enumerate}

\begin{just}
Each $\mu(t \breve{s})$ is $\pi(s)$-linear, so all that needs to be checked is that $\mu(t_1 \breve{s_1}) \circ \mu(t_2 \breve{s_2}) = \mu(t_1 \breve{s_2} t_2 \breve{s_2})$, for $t_1,t_2 \in T$ and $s_1,s_2 \in S_{\chi}$. We can compute
\begin{align*}
\mu(t_1 \breve{s_1} t_2 \breve{s_2})(1 \otimes w) &= \mu(t_1 \breve{s_1} t_2 \breve{s_1}^{-1} \gamma(s_1,s_2) \breve{s_1s_2})(1 \otimes w), \\
&= \chi(t_1) \chi(\breve{s_1} t_2 \breve{s_1}^{-1}) \chi(\gamma(s_1,s_2)) \otimes \sigma(s_1s_2)(w),
\end{align*}
and 
\begin{align*}
\mu(t_1 \breve{s_1})(\mu(t_2 \breve{s_2})(1 \otimes w)) &= \mu(t_1 \breve{s_1})(\chi(t_2) \otimes \sigma(s_2)(w)), \\
&= \chi(t_1) \otimes \sigma(s_1)(\chi(t_2) \sigma(s_2)(w)), \\
&= \chi(t_2)\: {}^{s_1}(\chi(t_2)) \chi(\gamma(s_1,s_2)) \otimes \sigma(s_1 s_2)(w),\\
&= \chi(t_2)\: {}^{s_1}(\chi(t_2)) \chi(\gamma(s_1,s_2)) \otimes \sigma(s_1 s_2)(w)
\end{align*}
using that ${}^{s_1} \chi = \chi$, so ${}^{s_1}(\chi(t_2)) = \chi(\breve{s_1} t_2 \breve{s_1}^{-1})$.
\end{just}

\begin{enumerate}[resume, label=\textnormal{(\Roman*)}]  
\item
  These are all the irreducible representations.
  There is a one-to-one correspondence
\begin{equation} \label{Cliff_2f}
  \bigcup_{\chi \in S \backslash T^{\vee}}  {\Irr_{\chi}^{\hspace{0.1em} \star}(S_\chi)} \stackrel{1:1}{\longleftrightarrow}
{\Irr_{\bC}^{\hspace{0.1em} \star}(G)}, \quad
(\chi, W) \mapsto \Ind_{G_\chi}^G (\bC_\chi \otimes W , \mu)  
\end{equation}
where the union is over representatives of $S$-orbits and $(\bC_\chi \otimes W, \mu)$  is defined in~\eqref{Cliff_2e}. 
    \end{enumerate}

\begin{just}
Let $V \in \Rep_{\bC}^{\hspace{0.1em} \star}(G)$ be finite-dimensional with $V = \sum_{\chi' \in S \cdot \chi} V_{\chi'}$. By \eqref{Cliff_3a}, $\breve{r} V_{\chi} = V_{{}^{r}\chi}$ for $r \in S$, so $V = \bigoplus_{r \in S/S_{\chi}} \breve{r} V_{\chi}$, the sum being direct because eigenvectors for different characters are linearly independent. In particular, $V_{\chi}$ is stable under $G_{\chi}$, and $V \cong \Ind_{G_{\chi}}^{G} V_{\chi}$. On $W \coloneqq V_{\chi}$ the group $T$ acts by $\chi$, and setting $\sigma(s) \coloneqq \mu(\breve{s})|_{W}$ for $s \in S_{\chi}$ gives
\[
\sigma(h)\sigma(k) = \mu(\gamma(h,k)) \mu(\breve{hk}) = \chi(\gamma(h,k)) \sigma(hk),
\]
so $(W, \sigma)$ is a $c_{\chi}$-projective antilinear representation of $S_{\chi}$, and $V_{\chi} = \bC_{\chi} \otimes W$ as in \eqref{Cliff_2e}. Conversely, for any $(W, \sigma)$ the $\chi$-eigenspace of $\Ind_{G_{\chi}}^{G}(\bC_{\chi} \otimes W)$ is $\bC_{\chi} \otimes W$, as the other summands have characters ${}^{r}\chi \neq \chi$.

A $G$-stable subspace $V' \subseteq V$ is $T$-stable, hence the sum of its eigenspaces, so $V' = \bigoplus_{r} \breve{r} V'_{\chi}$ with $V'_{\chi} \subseteq W$ stable under $S_{\chi}$. Conversely an $S_{\chi}$-stable subspace $W' \subseteq W$ gives the $G$-stable subspace $\bigoplus_{r} \breve{r} W'$. Hence $V$ is irreducible if and only if $W$ is.

It follows that the map \eqref{Cliff_2f} is well-defined and that it is surjective, as every $V \in \Irr_{\bC}^{\hspace{0.1em} \star}(G)$ satisfies the hypothesis for some $\chi$ by \eqref{Cliff_2b}. It is injective, as the orbit of $\chi$ is determined by $V$ as the set of characters with $V_{\chi} \neq 0$ and $W \cong V_{\chi}$.
\end{just}
  
\begin{enumerate}[resume, label=\textnormal{(\Roman*)}]  
\item The representation $V\in {\Irr_{\bC}^{\hspace{0.1em} \star}(G)}$ is unitary if and only if the corresponding character
  $\chi\in T^\vee$ is unitary. This means that $\chi$ takes values on the unit circle $\UU(1) \leq \bC^\times$. 
    \end{enumerate} 

\begin{just}
Suppose $(V,\mu)$ is unitary. Then $\mu(t)$ is a unitary operator for each $t \in T$, so its eigenvalues lie in $\UU(1)$. This proves the only if part.

Conversely, suppose that $\chi$ is unitary, and let $V \cong \Ind_{G_{\chi}}^{G}(\bC_{\chi} \otimes W)$. As $|\chi(\gamma(h,k))| = 1$, Weyl's Unitary Trick applies to the projective representation $W$: for any positive definite hermitian form $\tau$ on $W$,
\[
\widetilde{\tau}(x,y) \coloneqq \sum_{s \in (S_{\chi})_0} \tau(\sigma(s)x, \sigma(s)y) + \sum_{s \in S_{\chi} \setminus (S_{\chi})_0} \tau(\sigma(s)y, \sigma(s)x)
\]
is preserved by $\sigma(s)$ for $s$ even and flipped for $s$ odd. As $\chi$ is unitary, the same holds for the action \eqref{Cliff_2e} of $G_{\chi}$ on $V_{\chi} = \bC_{\chi} \otimes W$. Now $V = \bigoplus_{r \in S/S_{\chi}} \breve{r} V_{\chi}$, and we define a form on $V$ by declaring the summands orthogonal and setting
\[
\langle \breve{r}u, \breve{r}v \rangle \coloneqq \begin{cases} \widetilde{\tau}(u,v), & r \text{ even}, \\ \widetilde{\tau}(v,u), & r \text{ odd}, \end{cases} \qquad u,v \in V_{\chi}.
\]
For $g \in G$ we have $g \breve{r} = \breve{r}' h$ for some $h \in G_{\chi}$, where $\pi(g) + \pi(r) = \pi(r') + \pi(h)$, and so this form is preserved by $g$ if $g$ is even and flipped if $g$ is odd.
\end{just}

\section{The Extended Stabiliser and the Dyson Type}\label{sect:extendedstabiliser}

Let $G$, $H$ and $T$ be as above, so $H$ is an index two normal subgroup of $G$ and $T$ is a finite-index normal subgroup of $G$ contained in $H$. To any character $\chi \in T^{\vee}$ and $W \in \Irr_{\chi}^{\hspace{0.1em} \star}(S_\chi)$, we saw in the last section how to obtain an irreducible antilinear representation of $G$
\[V(\chi, W) \coloneqq \Ind_{G_\chi}^G (\bC_\chi \otimes W)\, .\]

The aim of this section is to answer the question: how does the D-type of $V(\chi, W)$ depend on $(\chi,W)$? For example, if $\chi_1, \chi_2 \in T^{\vee}$ with $S_{\chi_1} = S_{\chi_2}$ and $c_{\chi_1} = c_{\chi_2} \in \HH^2(S_{\chi_1}, (\bC^{\times}, \pi))$, then for any set-map $b \colon S_{\chi_1} \rightarrow \bC^{\times}$ realising this, so $d \chi_1 = d\chi_2 db$, then there is a bijection
\[
\Omega_b \colon \Irr_{\chi_1}^{\hspace{0.1em} \star}(S_{\chi_1}) \xrightarrow{\sim} \Irr_{\chi_2}^{\hspace{0.1em} \star}(S_{\chi_2}), \qquad (W, \sigma) \mapsto (W, b\sigma), \: \text{ where }  (b\sigma)(s) \coloneqq b(s) \sigma(s).
\]
What is the relationship between the D-types of $V(\chi_1, W)$ and $V(\chi_2, \Omega_b(W))$ for $W \in \Irr_{\chi_1}^{\hspace{0.1em} \star}(S_{\chi_1})$?

It turns out that the S-types will always be the same, but  the D-types can be different. In other words, $(S_{\chi}, c_{\chi})$ does  not contain enough information to determine the D-type.

%

The remedy is to consider {\em the extended group} $S \times C_2 = S \times \langle \epsilon \rangle$. This acts on $T^{\vee}$ by
\[
\big( (s, \epsilon^k) \cdot \chi \big)(t) \coloneqq \begin{cases} \chi(\breve{s}^{-1} t \breve{s}), & k = 0, \\ \overline{\chi(\breve{s}^{-1} t \breve{s})}, & k = 1. \end{cases}
\]

{\em The extended stabiliser} $S'_{\chi} \leq S \times C_2$ of $\chi \in T^{\vee}$ is the stabiliser of $\chi$ for this action. It has two homomorphisms to $C_2$
\[
\pi(s, \epsilon^k) \coloneqq \pi(s) \qquad \mbox{and} \qquad \kappa(s, \epsilon^k) \coloneqq \epsilon^k.
\]
We let $S'_{\chi}$ act on $\bC^{\times}$ through $\kappa$: elements $(s,1)$ act trivially and elements $(s,c)$ act by complex conjugation. We write $\bC^{\times}_{\kappa}$ for this $S'_{\chi}$- module and define
\[
c'_{\chi} \coloneqq [d \chi']^{-1} \in \HH^2(S'_{\chi}, \bC^{\times}_{\kappa}), \qquad d\chi'\big( (s,\epsilon^j), (t, \epsilon^k) \big) \coloneqq \chi(\gamma(s,t))^{-1}.
\]

The action through $\kappa$ is the unique action of $S'_{\chi}$ on $\bC^{\times}$ by the identity and complex conjugation for which $\chi \colon T \rightarrow \bC^{\times}$ is $S'_{\chi}$-equivariant for all $\chi$. Hence, $d\chi'$ is a $2$-cocycle since it is  the image of $\gamma$ under $\chi^{-1}$. The subgroup of $S'_{\chi}$ on which $\pi = \kappa$ is $\{ (s, \epsilon^{\pi(s)}) \mid s \in S_{\chi} \}$, which we identify with $S_{\chi}$. It has index $1$ or $2$, and $d\chi'$ and $c'_{\chi}$ restrict to $d\chi$ and $c_{\chi}$ respectively on it.
We are ready to state the key fundamental fact.


To state it, consider $\chi_1, \chi_2 \in T^{\vee}$ with $S_{\chi_1}' = S_{\chi_2}'$ and $c_{\chi_1}' = c_{\chi_2}' \in \HH^2(S_{\chi_1}', \bC^{\times}_{\kappa})$.
Choose any function $b \colon S_{\chi_1}' \rightarrow \bC^{\times}$ such that $d \chi_1' = d \chi_2' db$. Since $S_{\chi_1}' = S_{\chi_2}'$, we know that $S_{\chi_1} = S_{\chi_2}$.
Consider the induced bijection
\[
\Omega_b \colon \Irr_{\chi_1}^{\hspace{0.1em} \star}(S_{\chi_1}) \xrightarrow{\sim} \Irr_{\chi_2}^{\hspace{0.1em} \star}(S_{\chi_2}), \qquad (W, \sigma) \mapsto (W, b\sigma), \: \text{ where }  (b\sigma)(s) \coloneqq b(s) \sigma(s) \, .
\]
\begin{enumerate}[label=\textnormal{(\Roman*)}]
\item  
In the notation above, the induced bijection preserves D-type. In other words, $\fD(V(\chi_1, W)) = \fD(V(\chi_2, \Omega_b(W)))$ for any $W \in \Irr_{\chi_1}^{\hspace{0.1em} \star}(S_{\chi_1})$.
\end{enumerate}
\begin{just}
  The justification takes the rest of the chapter. By Theorem~\ref{thm:dyson10}, it is sufficient to show that the associated $C_2$-graded division algebras
  \[D(\chi_1, W) \coloneqq D(V(\chi_1, W)) \quad \mbox{ and } \quad D(\chi_2, \Omega_b(W)) \coloneqq D(V(\chi_2, \Omega_b(W)))\]
  are isomorphic as $C_2$-graded algebras. We can produce an explicit isomorphism.
To do that, we describe, for a pair $(\chi, W)$, the $C_2$-graded algebra $D(\chi, W) \coloneqq D(V(\chi, W))$ purely in terms of $W$. By definition, $D(V(\chi, W)) = \End_{\bR G}(V(\chi, W))$, graded by $\bC$-linear and $\bC$-antilinear maps.

Let $\rho(t \breve{s}) \coloneqq \chi(t) \sigma(s)$ denote the action of $t \breve{s} \in G_{\chi}$ on $U \coloneqq \bC_{\chi} \otimes W \subseteq V(\chi, W)$. For $y \in S'_{\chi}$ and $\lambda \in \bC$ we write
\[{}^{y}\lambda
\coloneqq
\begin{cases}
  \lambda & \mbox{ if } \kappa(y) = 1,\\
  \overline{\lambda} & \mbox{ if } \kappa(y) = \epsilon.
\end{cases}\]
We write $\breve{y}$ for the lift of the $S$-component of $y$. If $S'_{\chi} \neq S_{\chi}$, we fix $x \in S'_{\chi} \setminus S_{\chi}$. Then $\breve{x}$ normalises $G_{\chi}$ and $\breve{x}^2 \in G_{\chi}$.
Using this element $x$, we define the even and odd components of a $C_2$-graded algebra by
\[
D_0(U) \coloneqq \End_{\bC_{\star \hspace{-0.1em}}[G_{\chi}]}(U) = \{ \varphi_0 \colon U \rightarrow U \mid \varphi_0 \text{ is $\bC$-linear, and } \varphi_0 \, \rho(g) = \rho(g) \, \varphi_0 \text{ for all } g \in G_{\chi} \},
\]
and, if $S'_{\chi} \neq S_{\chi}$, 
\[
D_1(U) \coloneqq \{ \varphi \colon U \rightarrow U \mid \varphi(\lambda u) = {}^{x}\lambda \, \varphi(u), \text{ for all $\lambda \in \bC, u \in U$, } \varphi \, \rho(g) = \rho(\breve{x} g \breve{x}^{-1}) \, \varphi \text{ for all } g \in G_{\chi} \}.
\]
If $S'_{\chi} = S_{\chi}$, we instead set $D_1(U) \coloneqq 0$.

We need the following lemma, which we prove later.

\hrulefill
\begin{lemma}\label{lem:modelonW}
  Suppose that $\chi \in T^{\vee}$, $W \in \Irr_{\chi}^{\hspace{0.1em} \star}(S_{\chi})$, and set $V \coloneqq V(\chi, W)$. Then the map 
  \[
  D(V)_0 \rightarrow D_0(U), \qquad \Phi \mapsto \Phi|_U
  \]
  is a bijection. If $S'_{\chi} = S_{\chi}$, then $D(V)_1 = 0$. If not, then the map
  \[
  D(V)_1 \rightarrow D_1(U), \qquad \Phi \mapsto \Phi \breve{x}|_U
  \]
  is a bijection. Writing $D(U) \coloneqq D_0(U) \oplus D_1(U)$, the product on $D(V)$ transports to the product 
\[
\varphi_0 * \varphi \coloneqq \varphi_0 \circ \varphi, \qquad \varphi * \varphi_0 \coloneqq \varphi \circ \varphi_0, \qquad \varphi * \varphi' \coloneqq \varphi \circ \varphi' \circ \rho(\breve{x}^2)^{-1},
\]
on $D(U)$, where $\varphi_0 \in D_0(U)$ and $\varphi, \varphi' \in D_1(U)$.
\end{lemma}
\hrulefill

Now let $\sigma$ be the action of $S_{\chi_1}$ on $W$, so that $\Omega_b(W)$ is $W$ with the action $b\sigma$. Let us write
\[
\rho_1(t \breve{s}) = \chi_1(t) \sigma(s), \qquad \rho_2(t \breve{s}) = \chi_2(t) b(s) \sigma(s)
\]
for the two actions of $G_{\chi_1} = G_{\chi_2}$ on $W$. We apply Lemma~\ref{lem:modelonW} to both, with the same $x$ if $S'_{\chi_1} \neq S_{\chi_1}$. Let us write $D^{(1)}(W)$ and $D^{(2)}(W)$ for the resulting algebras. As $\rho_1(g)$ and $\rho_2(g)$ are proportional for all $g \in G_{\chi_1}$, we have $D^{(1)}_0(W) = D^{(2)}_0(W)$. If $S'_{\chi_1} = S_{\chi_1}$ there is nothing more to prove.
Otherwise, if $S'_{\chi_1} \neq  S_{\chi_1}$, the map
\[
D^{(1)}(W) \rightarrow D^{(2)}(W), \qquad \varphi_0 + \varphi \mapsto \varphi_0 + b(x) \varphi
\]
is an explicit isomorphism of $C_2$-graded algebras, which can be easily verified by a direct calculation.
\end{just}

Let us prove Lemma~\ref{lem:modelonW} now.

\begin{proof}
  We have that $V = \oplus_s \breve{s} U$ over the cosets of $S_{\chi}$ in $S$, and as a $T$-representation each $\breve{s} U$ is the $s \cdot \chi$-isotypical component of $V$. 
  Now for $s \in S$, $s \cdot \chi = \overline{\chi}$ if and only if $(s, \epsilon^{\pi(s) + 1}) \in S_{\chi}'$, and this happens for some $s \in S$ if and only if $S_{\chi}' \neq S_{\chi}$, as all elements of $S_{\chi}' \setminus S_{\chi}$ have this form. In particular, if $S_{\chi}' = S_{\chi}$, the $\overline{\chi}$-isotypic component of $V$ is zero, and otherwise is given by $\breve{x}^{-1} U$.

  Now as $V = \Ind_{G_{\chi}}^G(U)$, restriction is a bijection
  \[
  D(V) = \End_{\bR[G]}(V) \xrightarrow{\sim} \Hom_{\bR[G_{\chi}]}(U,V).
  \]
  For the first claim, note that this restricts to a bijection 
  \[
  D(V)_0 = \End_{\bC_{\star \hspace{-0.1em}}[G]}(V) \xrightarrow{\sim} \Hom_{\bC_{\star \hspace{-0.1em}}[G_{\chi}]}(U,V),
  \]
  and $\Hom_{\bC_{\star \hspace{-0.1em}}[G_{\chi}]}(U,V) = \Hom_{\bC_{\star \hspace{-0.1em}}[G_{\chi}]}(U,U) = D_0(U)$, because as a $T$-representation $U$ is the $\chi$-isotypic component of $V$. For the second claim, note that any $\Phi \in D(V)_1$, being antilinear, has restriction in 
  \[
  \Hom_{\bR[G_{\chi}]}(U,V_{\overline{\chi}}) \subset \Hom_{\bR[G_{\chi}]}(U,V),
  \]
  where $V_{\overline{\chi}}$ is the $\overline{\chi}$-isotypic part of $V$. Thus from above, if $S_{\chi}' = S_{\chi}$ then $\Hom_{\bR[G_{\chi}]}(U, V_{\overline{\chi}}) = 0$ and hence $D(V)_1 = 0$. If instead $S_{\chi}' \neq S_{\chi}$, then the restriction map restricts to a bijection 
  \[
  D(V)_1 \xrightarrow{\sim} \{ \psi \in \Hom_{\bR[G_{\chi}]}(U, \breve{x}^{-1} U) \mid \psi \text{ is $\bC$-antilinear}\}, \qquad \Phi \mapsto \Phi|_U.
  \]

   Now any map $\psi \colon U \rightarrow \breve{x}^{-1} U$ can be written uniquely as $\psi = \breve{x}^{-1} \varphi$ with $\varphi \colon U \rightarrow U$. As $\pi(x) \neq \kappa(x)$, $\psi$ is antilinear if and only if $\varphi(\lambda u) = {}^{x}\lambda \, \varphi(u)$ for all $\lambda \in \bC$ and $u \in U$, and $\psi$ commutes with $g \in G_{\chi}$ if and only if $\varphi \, \rho(g) = \rho(\breve{x} g \breve{x}^{-1}) \, \varphi$. Thus $\psi \mapsto \breve{x} \psi$ is a bijection from this set to $D_1(U)$, and as $\breve{x} \Phi|_U = \Phi \breve{x}|_U$, the map $\Phi \mapsto \Phi \breve{x} |_U$ is the required bijection.
  
  For the product formula, if $\Phi_0 \in D(V)_0$ and $\Phi, \Phi' \in D(V)_1$ with $\Phi_0|_U = \varphi_0 \in D_0(U)$ and $\Phi \breve{x}|_U = \varphi, \Phi' \breve{x}|_U =\varphi' \in D_1(U)$ respectively, then
  \[
  ((\Phi_0\Phi) \breve{x})|_U = \Phi_0|_U \Phi \breve{x}|_U = \varphi_0 \varphi, \qquad ((\Phi \Phi_0) \breve{x})|_U = (\Phi \breve{x}) \Phi_0|_U = \varphi \varphi_0,
  \]
  and
  \[
  (\Phi \Phi')|_U = (\Phi \breve{x} \Phi' \breve{x} \breve{x}^{-2})|_U = \varphi \varphi' \rho(\breve{x}^2)^{-1}
  \]
  as required.
\end{proof}

%% file: appendix-tablesfor24types.tex
\chapter{Properties of the 24 Types}\label{chapter:tablesfor24ways}

\begin{table}[h!] 
     \vspace{1cm}
\centering
\bgroup
\caption{The 24 Possible Types for $V$}
\label{tab:24types}
\def\arraystretch{1.1}
\begin{tabular}{c r c  c  c  c  c  c c c}
    \toprule
    & &  $\bK[G] \subset \bL[G]$ & $\bK[H] \subset \bL[H]$ & $\bL[H] \subset \bL_{\star \hspace{-0.1em}}[G]$ & $\bL[H] \subset \bL[G]$ & $\bK[H] \subset \bK[G]$ & $\bK[G] \subset \bL_{\star \hspace{-0.1em}}[G]$ \\
    \midrule
I & 1.  &  \textbf{A} & \textbf{A} & \textbf{A} & \textbf{1} & \textbf{1} & \textbf{2}  \\ 
III & 2.  &  \textbf{A} & \textbf{A} & \textbf{B} & \textbf{2} & \textbf{2} & \textbf{3(b)}  \\
& 3.  &  \textbf{A} & \textbf{A} & \textbf{A} & \textbf{3(a)} & \textbf{3(a)} & \textbf{3(b)}  \\
II.V & 4.  &  \textbf{B} & \textbf{A} & \textbf{A} & \textbf{1} & \textbf{3(b)} & \textbf{3(a)}  \\
II & 5.  &  \textbf{B} & \textbf{A} & \textbf{C} & \textbf{1} & \textbf{3(b)} & \textbf{3(b)}  \\
& 6.  &  \textbf{C} & \textbf{A} & \textbf{A} & \textbf{3(a)} & \textbf{3(b)} & \textbf{3(a)}  \\
& 7.  &  \textbf{C} & \textbf{A} & \textbf{C} & \textbf{3(a)} & \textbf{3(b)} & \textbf{3(b)}  \\ 
& 8.  &  \textbf{A} & \textbf{A} & \textbf{A} & \textbf{3(b)} & \textbf{3(b)} & \textbf{3(a)}  \\
& 9.  &  \textbf{A} & \textbf{A} & \textbf{C} & \textbf{3(b)} & \textbf{3(b)} & \textbf{3(b)}  \\
IV & 10. &  \textbf{B} & \textbf{B} & \textbf{B} & \textbf{1} & \textbf{1} & \textbf{2}  \\ 
VII & 11. &  \textbf{B} & \textbf{B} & \textbf{B} & \textbf{2} & \textbf{2} & \textbf{1}  \\ 
V & 12. &  \textbf{A} & \textbf{B} & \textbf{A} & \textbf{2} & \textbf{3(a)} & \textbf{1}  \\ 
& 13. &  \textbf{B} & \textbf{B} & \textbf{B} & \textbf{3(a)} & \textbf{3(a)} & \textbf{3(b)}  \\ 
& 14. &  \textbf{B} & \textbf{B} & \textbf{B} & \textbf{3(b)} & \textbf{3(b)} & \textbf{3(a)}  \\ 
VI & 15. &  \textbf{C} & \textbf{B} & \textbf{C} & \textbf{2} & \textbf{3(b)} & \textbf{1}  \\ 
VIII & 16. &  \textbf{C} & \textbf{C} & \textbf{C} & \textbf{1} & \textbf{1} & \textbf{2}  \\ 
X & 17. &  \textbf{C} & \textbf{C} & \textbf{B} & \textbf{2} & \textbf{2} & \textbf{3(a)}  \\ 
IX & 18. &  \textbf{B} & \textbf{C} & \textbf{A} & \textbf{1} & \textbf{3(a)} & \textbf{3(a)}  \\ 
IX.V & 19. &  \textbf{B} & \textbf{C} & \textbf{C} & \textbf{1} & \textbf{3(a)} & \textbf{3(b)}  \\ 
& 20. &  \textbf{C} & \textbf{C} & \textbf{A} & \textbf{3(a)} & \textbf{3(a)} & \textbf{3(a)}  \\ 
& 21. &  \textbf{C} & \textbf{C} & \textbf{C} & \textbf{3(a)} & \textbf{3(a)} & \textbf{3(b)}  \\ 
& 22. &  \textbf{A} & \textbf{C} & \textbf{A} & \textbf{3(b)} & \textbf{3(a)} & \textbf{3(a)}  \\ 
& 23. &  \textbf{A} & \textbf{C} & \textbf{C} & \textbf{3(b)} & \textbf{3(a)} & \textbf{3(b)}  \\
& 24. &  \textbf{C} & \textbf{C} & \textbf{C} & \textbf{3(b)} & \textbf{3(b)} & \textbf{3(a)}  \\
\bottomrule
\end{tabular}
\egroup
\end{table}

\begin{table}[p] 
    \vspace{4cm}
\centering
\bgroup
\caption{The Behaviour of the 24 Types for $V$ (i)}
\label{tab:24types2}
\def\arraystretch{1.1}
\begin{tabular}{c r c c c c c c c c}
    \toprule
    & & $\Ind_{\bK[G]}^{\bL[G]}(U)$ & $\Res_{\bK[G]}^{\bL[G]}(Z)$ & $\Ind_{\bK[H]}^{\bL[H]}(X)$ & $\Res_{\bK[H]}^{\bL[H]}(V)$ & $\Ind_{\bL[H]}^{\bL_{\star \hspace{-0.1em}}[G]}(V)$ & $\Res_{\bL[H]}^{\bL_{\star \hspace{-0.1em}}[G]}(W)$ \\
    \midrule
I & 1.  &  $Z$ & $U \oplus U$ & $V$ & $X \oplus X$ & $W \oplus W$ & $V$ \\ 
III & 2.  &  $Z$ & $U \oplus U$ & $V$ & $X \oplus X$ & $W$ & $V \oplus \epsilon * V$  \\ 
& 3.  &  $Z$ & $U \oplus U$ & $V$ & $X \oplus X$ & $W \oplus W$ & $V$  \\ 
II.V & 4.  &  $Z \oplus \overline{Z}$ & $U$ & $V$ & $X \oplus X$ & $W \oplus W$ & $V$  \\ 
II & 5.  &  $Z \oplus \overline{Z}$ & $U$ & $V$ & $X \oplus X$ & $W$ & $V \oplus V$  \\ 
& 6.  &  $Z \oplus Z$ & $U$ & $V$ & $X \oplus X$ & $W \oplus W$ & $V$  \\ 
& 7.  &  $Z \oplus Z$ & $U$ & $V$ & $X \oplus X$ & $W$ & $V \oplus V$ \\ 
& 8.  &  $Z$ & $U \oplus U$ & $V$ & $X \oplus X$ & $W \oplus W$ & $V$ \\ 
& 9.  &  $Z$ & $U \oplus U$ & $V$ & $X \oplus X$ & $W$ & $V \oplus V$  \\ 
IV & 10. &  $Z \oplus \overline{Z}$ & $U$ & $V \oplus \overline{V}$ & $X$ & $W$ & $V \oplus \epsilon * V$  \\ 
VII & 11. &  $Z \oplus \overline{Z}$ & $U$ & $V \oplus \overline{V}$ & $X$ & $W$ & $V \oplus \epsilon * V$  \\ 
V & 12. &  $Z$ & $U \oplus U$ & $V \oplus \overline{V}$ & $X$ & $W \oplus W$ & $V$ \\ 
& 13. &  $Z \oplus \overline{Z}$ & $U$ & $V \oplus \overline{V}$ & $X$ & $W$ & $V \oplus \epsilon * V$  \\ 
& 14. &  $Z \oplus \overline{Z}$ & $U$ & $V \oplus \overline{V}$ & $X$ & $W$ & $V \oplus \epsilon * V$  \\ 
VI & 15. &  $Z \oplus Z$ & $U$ & $V \oplus \overline{V}$ & $X$ & $W$ & $V \oplus V$ \\ 
VIII & 16. &  $Z \oplus Z$ & $U$ & $V \oplus V$ & $X$ & $W$ & $V \oplus V$  \\ 
X & 17. &  $Z \oplus Z$ & $U$ & $V \oplus V$ & $X$ & $W$ & $V \oplus \epsilon * V$  \\ 
IX & 18. &  $Z \oplus \overline{Z}$ & $U$ & $V \oplus V$ & $X$ & $W \oplus W$ & $V$ \\ 
IX.V & 19. &  $Z \oplus \overline{Z}$ & $U$ & $V \oplus V$ & $X$ & $W$ & $V \oplus V$ \\ 
& 20. &  $Z \oplus Z$ & $U$ & $V \oplus V$ & $X$ & $W \oplus W$ & $V$ \\ 
& 21. &  $Z \oplus Z$ & $U$ & $V \oplus V$ & $X$ & $W$ & $V \oplus V$  \\ 
& 22. &  $Z$ & $U \oplus U$ & $V \oplus V$ & $X$ & $W \oplus W$ & $V$  \\ 
& 23. &  $Z$ & $U \oplus U$ & $V \oplus V$ & $X$ & $W$ & $V \oplus V$  \\ 
& 24. &  $Z \oplus Z$ & $U$ & $V \oplus V$ & $X$ & $W$ & $V \oplus V$ \\ 
\bottomrule
\end{tabular}
\egroup
\end{table}

\begin{table}[p] 
    \vspace{4cm}
\centering
\bgroup
\caption{The Behaviour of the 24 Types for $V$ (ii)}
\label{tab:24types4}
\def\arraystretch{1.1}
\begin{tabular}{c r c c c c c c c c}
    \toprule
    & & $\Ind_{\bL[H]}^{\bL[G]}(V)$ & $\Res_{\bL[H]}^{\bL[G]}(Z)$ & $\Ind_{\bK[H]}^{\bK[G]}(X)$ & $\Res_{\bK[H]}^{\bK[G]}(U)$ & $\Ind_{\bK[G]}^{\bL_{\star \hspace{-0.1em}}[G]}(U)$ & $\Res_{\bK[G]}^{\bL_{\star \hspace{-0.1em}}[G]}(W)$ \\
    \midrule
I & 1.  & $Z \oplus (Z \otimes \pi)$ & $V$ & $U \oplus (U \otimes \pi)$ & $X$ & $W$ & $U \oplus (U \otimes \pi)$ \\
III & 2.  & $Z$ & $V \oplus {}^{\omega} V$ & $U$ & $X \oplus {}^{\omega} X$ & $W$ & $U \oplus U$ \\ 
& 3.  & $Z \oplus Z$ & $V$ & $U \oplus U$ & $X$ & $W$ & $U \oplus U$ \\ 
II.V & 4.  & $Z \oplus (Z \otimes \pi)$ & $V$ & $U$ & $X \oplus X$ & $W \oplus W$ & $U$ \\ 
II & 5.  & $Z \oplus (Z \otimes \pi)$ & $V$ & $U$ & $X \oplus X$ & $W$ & $U \oplus U$ \\ 
& 6.  & $Z \oplus Z$ & $V$ &$U$ & $X \oplus X$ & $W \oplus W$ & $U$ \\ 
& 7.  & $Z \oplus Z$ & $V$ & $U$ & $X \oplus X$ & $W$ & $U \oplus U$ \\ 
& 8.  & $Z$ & $V \oplus V$ & $U$ & $X \oplus X$ & $W \oplus W$ & $U$ \\ 
& 9.  & $Z$ & $V \oplus V$ & $U$ & $X \oplus X$ & $W$ & $U \oplus U$ \\ 
IV & 10. & $Z \oplus (Z \otimes \pi)$ & $V$ & $U \oplus (U \otimes \pi)$ & $X$ & $W$ & $U \oplus (U \otimes \pi)$  \\ 
VII & 11. & $Z$ & $V \oplus {}^{\omega} V$ & $U$ & $X \oplus {}^{\omega} X$ & $W \oplus \overline{W}$ & $U$  \\ 
V & 12. & $Z$ & $V \oplus {}^{\omega} V$ & $U \oplus U$ & $X$ & $W \oplus \overline{W}$ & $U$  \\ 
& 13. & $Z \oplus Z$ & $V$ & $U \oplus U$ & $X$ & $W$ & $U \oplus U$   \\ 
& 14. & $Z$ & $V \oplus V$ & $U$ & $X \oplus X$ & $W \oplus W$ & $U$ \\ 
VI & 15. & $Z$ & $V \oplus {}^{\omega} V$ & $U$ & $X \oplus X$ & $W \oplus \overline{W}$ & $U$  \\ 
VIII & 16. & $Z \oplus (Z \otimes \pi)$ & $V$ & $U \oplus (U \otimes \pi)$ & $X$ & $W$ & $U \oplus (U \otimes \pi)$  \\ 
X & 17. & $Z$ & $V \oplus {}^{\omega} V$ & $U$ & $X \oplus {}^{\omega} X$ & $W \oplus W$ & $U$ \\ 
IX & 18. & $Z \oplus (Z \otimes \pi)$ & $V$ & $U \oplus U$ & $X$ & $W \oplus W$ & $U$ \\ 
IX.V & 19. & $Z \oplus (Z \otimes \pi)$ & $V$ & $U \oplus U$ & $X$ & $W$ & $U \oplus U$   \\ 
& 20. & $Z \oplus Z$ & $V$ & $U \oplus U$ & $X$ & $W \oplus W$ & $U$ \\ 
& 21. & $Z \oplus Z$ & $V$ &$U \oplus U$ & $X$ & $W$ & $U \oplus U$   \\ 
& 22. & $Z$ & $V \oplus V$ & $U \oplus U$ & $X$ & $W \oplus W$ & $U$ \\ 
& 23. & $Z$ & $V \oplus V$ & $U \oplus U$ & $X$ & $W$ & $U \oplus U$ \\ 
& 24. & $Z$ & $V \oplus V$ &$U$ & $X \oplus X$ & $W \oplus W$ & $U$ \\
\bottomrule
\end{tabular}
\egroup
\end{table}

\begin{table}[p] 
    \vspace{4cm}
\centering
\bgroup
\caption{The Behaviour of the 24 Types for $V$ (iii)}
\label{tab:24types5}
\def\arraystretch{1.1}
\begin{tabular}{c r c c c c c c c c}
    \toprule
    & & $\Stab_{K_4}(Z)$ & $\Stab_{K_4}(V)$ & $X \cong {}^{\omega}X$ & $U \cong (U \otimes \pi)$ & $W \cong \overline{W}$ & $\overline{m}_{\bbK}^{\bbL}(Z)$ & $\overline{m}_{\bbK}^{\bbL}(V)$ & $m_{\bbK}^{\bbL}(V)$ \\
    \midrule
I & 1.  & $\langle a \rangle$ & $K_4$ & Yes & No & Yes & $1$ & $1$ & $1$ \\ 
III & 2.  & $K_4$ & $\langle a \rangle$ & No & Yes & Yes & $1$ & $1$ & $1$ \\ 
& 3.  & $K_4$ & $K_4$ & Yes & Yes & Yes & $1$ & $1$ & $1$ \\ 
II.V & 4.  & $\langle c \rangle$ & $K_4$ & Yes & Yes & Yes & $1$ & $1$ & $1$ \\ 
II & 5.  & $\langle c \rangle$ & $K_4$ & Yes & Yes & Yes & $1$ & $1$ & $2$ \\ 
& 6.  & $K_4$ & $K_4$ & Yes & Yes & Yes & $2$ & $1$ & $1$ \\ 
& 7.  & $K_4$ & $K_4$ & Yes & Yes & Yes & $2$ & $1$ & $2$ \\ 
& 8.  & $K_4$ & $K_4$ & Yes & Yes & Yes & $1$ & $1$ & $1$ \\ 
& 9.  & $K_4$ & $K_4$ & Yes & Yes & Yes & $1$ & $1$ & $2$ \\ 
IV & 10. & $1$ & $\langle b \rangle$ & Yes & No & Yes & $1$ & $1$ & $1$ \\ 
VII & 11. & $\langle b \rangle$ & 1 & No & Yes & No & $1$ & $1$ & $1$ \\ 
V & 12. & $K_4$ & $\langle c \rangle$ & Yes & Yes & No & $1$ & $1$ & $1$ \\ 
& 13. & $\langle b \rangle$ & $\langle b \rangle$ & Yes & Yes & Yes & $1$ & $1$ & $1$ \\ 
& 14. & $\langle b \rangle$ & $\langle b \rangle$ & Yes & Yes & Yes & $1$ & $1$ & $1$ \\ 
VI & 15. & $K_4$ & $\langle c \rangle$ & Yes & Yes & No & $2$ & $1$ & $2$ \\ 
VIII & 16. & $\langle a \rangle$ & $K_4$ & Yes & No & Yes & $2$ & $2$ & $2$ \\ 
X & 17. & $K_4$ & $\langle a \rangle$ & No & Yes & Yes & $2$ & $2$ & $1$ \\ 
IX & 18. & $\langle c \rangle$ & $K_4$ & Yes & Yes & Yes & $1$ & $2$ & $1$ \\ 
IX.V & 19. & $\langle c \rangle$ & $K_4$ & Yes & Yes & Yes & $1$ & $2$ & $2$ \\ 
& 20. & $K_4$ & $K_4$ & Yes & Yes & Yes & $2$ & $2$ & $1$ \\ 
& 21. & $K_4$ & $K_4$ & Yes & Yes & Yes & $2$  & $2$ & $2$ \\ 
& 22. & $K_4$ & $K_4$ & Yes & Yes & Yes & $1$ & $2$ & $1$ \\ 
& 23. & $K_4$ & $K_4$ & Yes & Yes & Yes & $1$ & $2$ & $2$ \\ 
& 24. & $K_4$ & $K_4$ & Yes & Yes & Yes & $2$ & $2$ & $2$ \\
\bottomrule
\end{tabular}
\egroup
\end{table}

\begin{table}[p] 
    \vspace{4cm}
\centering
\bgroup
\caption{The Behaviour of the 24 Types for $V$ (iv): Dimension Ratios}
\label{tab:24types5.5}
\def\arraystretch{1.1}
\begin{tabular}{c r c c c c c c c c c}
    \toprule
    & & $\dim_{\bbL}Z$ & $\dim_{\bbL} V$ & $\dim_{\bbK}X$ & $\dim_{\bbK}U$ & $\dim_{\bbL} W$ \\
    \midrule
I & 1.  & 1 & 1 & 1 & 1 & 1 \\  
III & 2.  & 2 & 1 & 1 & 2 & 2 \\ 
& 3.  & 2 & 2 & 2 & 2 & 2 \\ 
II.V & 4.  & 1 & 1 & 1 & 2 & 1 \\ 
II & 5.  & 1 & 1 & 1 & 2 & 2 \\ 
& 6.  & 2 & 2 & 2 & 4 & 2 \\ 
& 7.  & 2 & 2 & 2 & 4 & 4 \\ 
& 8.  & 2 & 1 & 1 & 2 & 1 \\ 
& 9.  & 2 & 1 & 1 & 2 & 2 \\ 
IV & 10. & 1 & 1 & 2 & 2 & 2 \\ 
VII & 11. & 2 & 1 & 2 & 4 & 2 \\ 
V & 12. & 2 & 1 & 2 & 2 & 1 \\ 
& 13. & 2 & 2 & 4 & 4 & 4 \\ 
& 14. & 2 & 1 & 2 & 4 & 2 \\ 
VI & 15. & 2 & 1 & 2 & 4 & 2 \\ 
VIII & 16. & 2 & 2 & 4 & 4 & 4 \\ 
X & 17. & 4 & 2 & 4 & 8 & 4 \\ 
IX & 18. & 2 & 2 & 4 & 4 & 2 \\
IX.V & 19. & 2 & 2 & 4 & 4 & 4 \\ 
& 20. & 4 & 4 & 8 & 8 & 4 \\ 
& 21. & 4 & 4 & 8 & 8 & 8 \\ 
& 22. & 4 & 2 & 4 & 4 & 2 \\ 
& 23. & 4 & 2 & 4 & 4 & 4 \\ 
& 24. & 4 & 2 & 4 & 8 & 4 \\
\bottomrule
\end{tabular}
\egroup
\end{table}

\begin{table}[p] 
    \vspace{4cm}
\centering
\bgroup
\caption{The Behaviour of the 24 Types for $V$ (v): Endomorphism Dimension Ratios}
\label{tab:24types6}
\def\arraystretch{1.1}
\begin{tabular}{c r c c c c c c c c c}
    \toprule
    &  & $\dim_{\bbL}\End_{\bL[G]}(Z)$ & $\dim_{\bbL}\End_{\bL[H]}(V)$ & $\dim_{\bbK}\End_{\bK[H]}(X)$ & $\dim_{\bbK}\End_{\bK[G]}(U)$ & $\dim_{\bbK}\End_{\bL_{\star \hspace{-0.1em}}[G]}(W)$ \\
    \midrule
I & 1.  & 1 & 1 & 1 & 1 & 1 \\ 
III & 2.  & 1 & 1 & 1 & 1 & 2 \\ 
& 3.  & 1 & 2 & 2 & 1 & 2 \\ 
II.V & 4.  & 1 & 1 & 1 & 2 & 1 \\ 
II & 5.  & 1 & 1 & 1 & 2 & 4 \\ 
& 6.  & 1 & 2 & 2 & 4 & 2 \\ 
& 7.  & 1 & 2 & 2 & 4 & 8 \\ 
& 8.  & 2 & 1 & 1 & 2 & 1 \\ 
& 9.  & 2 & 1 & 1 & 2 & 4 \\ 
IV & 10. & 1 & 1 & 2 & 2 & 2 \\ 
VII & 11. & 1 & 1 & 2 & 2 & 2 \\ 
V & 12. & 1 & 1 & 2 & 1 & 1 \\ 
& 13. & 1 & 2 & 4 & 2 & 4 \\ 
& 14. & 2 & 1 & 2 & 4 & 2 \\ 
VI & 15. & 1 & 1 & 2 & 4 & 4 \\ 
VIII & 16. & 1 & 1 & 4 & 4 & 4 \\ 
X & 17. & 1 & 1 & 4 & 4 & 2 \\ 
IX & 18. & 1 & 1 & 4 & 2 & 1 \\ 
IX.V & 19. & 1 & 1 & 4 & 2 & 4 \\ 
& 20. & 1 & 2 & 8 & 4 & 2 \\ 
& 21. & 1 & 2 & 8 & 4 & 8 \\ 
& 22. & 2 & 1 & 4 & 2 & 1 \\ 
& 23. & 2 & 1 & 4 & 2 & 4 \\ 
& 24. & 2 & 1 & 4 & 8 & 4 \\
\bottomrule
\end{tabular}
\egroup
\end{table}

\begin{table}[p] 
    \vspace{7cm}
\centering
\bgroup
\caption{The Behaviour of the 12 Types for $V$}
\label{tab:12types}
\def\arraystretch{1.1}
\begin{tabular}{r c c c c c c c c c}
    \toprule
    & $\End_{\bL[G]}(Z)$ & $\End_{\bL[H]}(V)$ & $\End_{\bK[H]}(X)$ & $\End_{\bK[G]}(U)$ & $\End_{\bL_{\star \hspace{-0.1em}}[G]}(W)$ \\
    \midrule
I  & $\bL$ & $\bL$ & $\bK$ & $\bK$ & $\bK$ \\
II  & $\bL$ & $\bL$ & $\bK$ & $\bL$ & $D$ \\
II.V  & $\bL$ & $\bL$ & $\bK$ & $\bL$ & $\bK$ \\
III  & $\bL$ & $\bL$ & $\bK$ & $\bK$ & $\bL$ \\ 
IV & $\bL$ & $\bL$ & $\bL$ & $\bL$ & $\bL$ \\ 
V & $\bL$ & $\bL$ & $\bL$ & $\bK$ & $\bK$ \\ 
VI & $\bL$ & $\bL$ & $\bL$ & $D$ & $D$ \\ 
VII & $\bL$ & $\bL$ & $\bL$ & $\bL$ & $\bL$ \\ 
VIII & $\bL$ & $\bL$ & $D$ & $D$ & $D$ \\ 
IX & $\bL$ & $\bL$ & $D$ & $\bL$ & $\bK$ \\ 
IX.V & $\bL$ & $\bL$ & $D_1$ & $\bL$ & $D_2$ \\ 
X & $\bL$ & $\bL$ & $D$ & $D$ & $\bL$ \\
\bottomrule
\end{tabular}
\egroup
\end{table}

%% file: appendix-classandcohomC2modules.tex
\chapter{Classification and Cohomology of $C_2$-Modules}\label{appendix:classandcohomC2modules}

In this appendix we first recall some basic facts about the cohomology of $C_2$-modules. We then give a classification theorem for finitely-generated $C_2$-modules, which reduces the classification of arbitrary $C_2$-modules to that of $C_2$-modules which are $2$-power torsion, a complicated but ultimately algorithmically finite problem.

\section{Cohomology of $C_2$-modules}\label{sect:cohomC2modules}
Write $C_2 = \langle \epsilon \rangle$ for the group with two elements, where $\epsilon^2 = 1$. A $C_2$-module is just a module over the group algebra $\bZ [C_2]$, i.e. an additive group $A=(A,+)$ with an order two automorphism. 
To describe the cohomology of a $C_2$-module $A$, we need invariants and skew-invariants
\[
A^0 \coloneqq \{ x\in A \;\mid\; \epsilon x = x \}, \ \
A^1 \coloneqq \{ x\in A \;\mid\; \epsilon x = -x \}.
\]
\begin{prop} \label{prop:C2cohomology}
  The normalised bar resolution of the trivial $\bZ[C_2]$-module $\bZ$ is the free resolution
  \[
  0 \leftarrow \bZ \xleftarrow{\emph{\text{aug}}} \bZ[C_2] \xleftarrow{\epsilon - 1} \bZ[C_2] \xleftarrow{\epsilon + 1}  \bZ[C_2] \xleftarrow{\epsilon - 1} \bZ[C_2] \xleftarrow{\epsilon + 1} \cdots .
  \]
In particular, the cohomology groups of a $C_2$-module $A$ are as follows:
\[
\HH^{n}(A) \coloneqq \HH^{n}(C_2, A) =
\begin{cases}
  A^0, & \mbox{ if } n=0, \\
  A^0/(\epsilon + 1) A, & \mbox{ if } n \mbox{ is even and } n>0, \\
  A^1/(\epsilon - 1) A, & \mbox{ if } n \mbox{ is odd}.
\end{cases}
\]
\end{prop}

\begin{proof}
  This is \cite[Ex.\ 6.5.2(4)]{Weib} and \cite[Thm.\ 6.2.2]{Weib}.
\end{proof}

\begin{remark}\label{rem:connectingmaps}
  Suppose that $0 \rightarrow A \xrightarrow{\iota} B \rightarrow C \rightarrow 0$ is a short exact sequence of $C_2$-modules. From the explicit resolution of Proposition \ref{prop:C2cohomology} we can also explicitly describe the connecting homomorphisms in the associated long exact sequence on cohomology. We have
\[
\left\{
\begin{array}{r@{}l@{\quad}r@{}l@{\qquad}l}
  \delta_0      & \colon C^0 \rightarrow A^1/(\epsilon-1)A,               & \: \delta_0(c)    & = \iota^{-1}(\epsilon(b)-b), & \text{if } n=0, \\
  \delta_n      & \colon C^0/(\epsilon+1)C \rightarrow A^1/(\epsilon-1)A, & \: \delta_n(c)    & = \iota^{-1}(\epsilon(b)-b), & \text{if } n \text{ is even and } n>0, \\
  \delta_n      & \colon C^1/(\epsilon-1)C \rightarrow A^0/(\epsilon+1)A, & \: \delta_n(c)    & = \iota^{-1}(\epsilon(b)+b), & \text{if } n \text{ is odd},
\end{array} \right.
\]
  where in each case $b \in B$ is any lift of $c \in C$ to $B$.
\end{remark}

Recall that a group is called an elementary 2-group if it is a direct sum of copies of $\bZ/2\bZ$ (and so any such group is naturally viewed as a vector space over the field $\bF_2$ of two elements).
\begin{cor}\label{cor:cohomisvectorspace}
  The group $\HH^n(A)$ is an elementary 2-group for all $n\geq 1$.
\end{cor}
\begin{proof}
If $(1 \pm \epsilon)x=0$, then $(1 \mp \epsilon)x = 2x$, so $2x=0$ in $\HH^n(A)$ for all $x$.     
\end{proof}  

Let us now compute explicitly the cohomology of some classes of $C_2$-modules.
Let $B$ be an abelian group. We can create three $C_2$-modules from $B$:
\[
B_+ = B, \ \epsilon x = x; \qquad B_- = B, \ \epsilon x = -x; \qquad   B_d = B\oplus B, \ \epsilon (x,y) = (y,x). 
\]
For a prime $p$, we also consider the torsion subgroups
\[
B_p = \{ x\in B \,\mid\, px =0\}, \ \
B_{(p)} = \{ x\in B \,\mid\, p^nx =0 \text{ for some $n > 0$}\}.
\]
For these modules, the answer immediately follows from Proposition~\ref{prop:C2cohomology}:
\begin{equation}\label{eqn:cohomofspecialmodules}
\HH^{n}(B_+) =
\begin{cases}
  B, \\
  B/2B, \\
  B_2,
\end{cases}
\HH^{n}(B_-) =
\begin{cases}
  B_2, \\
  B_2, \\
  B/2B, 
\end{cases}
\HH^{n}(B_d) =
\begin{cases}
  B, & \mbox{ if } n=0, \\
  0, & \mbox{ if } n=2k >0, \\
  0, & \mbox{ if } n=2k+1 .
\end{cases}
\end{equation}

An abelian group $B$ admits division by 2 if  $x\mapsto 2x$ is an automorphism of $B$. We denote its inverse by $x\mapsto x/2$. 
\begin{lemma} \label{cohom_div_by2}
  Suppose a $C_2$-module $A$ admits division by $2$. Then $A\cong (A^0)_+ \oplus (A^1)_-$ as $C_2$-modules and
  $\HH^n(A)=0$ for all $n\geq 1$.  
\end{lemma}
\begin{proof} Clearly,  $A^0 + A^1 = A^0 \oplus A^1$ is a $C_2$-submodule of $A$. The inverse to the embedding map is given by
  $x \mapsto ((\epsilon(x) + x)/2, (\epsilon(x) - x)/2 )$. The second statement follows directly from the formulas (\ref{eqn:cohomofspecialmodules}), noting that if $B$ admits division by $2$ then both $B / 2 B = 0$ and $B_2 = 0$.
\end{proof}  

In our classification, we will be interested in $C_2$-modules of the following form. Let $B$ be a $\bZ[C_2]$-module, and consider a pair of abelian group homomorphisms 
  \[
  \Theta = (\Theta_+, \Theta_-), \ \ \Theta_+ : \bZ^n \rightarrow B , \ \  \Theta_- : \bZ^m \rightarrow B .
  \]

\begin{lemma} \label{even_part}
  Consider an operator
  \[\epsilon (x,y,z) = (\epsilon x +\Theta_+ (y) + \Theta_- (z), y, -z)\]
  on the abelian group $B_\Theta (\bZ) \coloneqq B \oplus \bZ^n \oplus \bZ^m$. 
It defines a $\bZ [C_2]$-module structure if and only if 
\newline
$\Theta_+ (\bZ^n ) \subseteq B^1 = \HZ^1 (B)$ and  $\Theta_- (\bZ^m) \subseteq B^0=\HZ^2(B)$.
\end{lemma}
\begin{proof} We need to have $\epsilon^2=1$. This means  
  \[(x,y,z) = \epsilon^2 (x,y,z) = (x + (\epsilon + 1) \Theta_+ (y) + (\epsilon - 1)\Theta_- (z), y, z)\]
  for all $x,y,z$. This proves the lemma.
\end{proof}

Whenever $B_\Theta (\bZ)$ is a $C_2$-module, $B_{\Theta}(\bZ)$ can be naturally viewed as an extension
\[
0 \rightarrow B \rightarrow B_\Theta (\bZ) \rightarrow \bZ^n_+ \oplus \bZ^m_- \rightarrow 0
\]
of $\bZ^n_+ \oplus \bZ^m_-$ by $B$. By definition, two $C_2$-modules $B_\Theta(\bZ)$ and $B_{\Psi}(\bZ)$ define equivalent extensions when there is an isomorphism $\gamma \colon B_\Theta(\bZ) \rightarrow B_{\Psi}(\bZ)$ such that the diagram
\begin{equation} \label{diagram_BTheta}
  \begin{tikzcd}
	0 & B & {{B}_\Theta(\bZ)} & {\bZ^n_+\oplus \bZ^m_-} & 0 \\
	0 & B & {{B}_\Psi(\bZ)} & {\bZ^n_+\oplus \bZ^m_-} & 0
	\arrow[from=1-1, to=1-2]
	\arrow[from=1-2, to=1-3]
	\arrow["\id"', from=1-2, to=2-2]
	\arrow[from=1-3, to=1-4]
	\arrow["\gamma"', from=1-3, to=2-3]
	\arrow[from=1-4, to=1-5]
	\arrow["\id"', from=1-4, to=2-4]
	\arrow[from=2-1, to=2-2]
	\arrow[from=2-2, to=2-3]
	\arrow[from=2-3, to=2-4]
	\arrow[from=2-4, to=2-5]
\end{tikzcd}\end{equation}
commutes.

\begin{cor} \label{even_part_ext_new}
Two $C_2$-modules $B_\Theta(\bZ)$ and $B_{\Psi}(\bZ)$ define equivalent extensions if and only if $\Theta_+-\Psi_+$ is valued in $(\epsilon - 1)(B) = B^1(C_2, B)$ and $\Theta_--\Psi_-$ is valued in $(\epsilon + 1)(B) = B^2(C_2, B)$, and any extension of $
\bZ^n_+ \oplus \bZ^m_-$ by $B$ is of the form $B_\Theta(\bZ)$. In particular,
\[\Ext_{\bZ[C_2]}^1 (\bZ^n_+\oplus \bZ^m_-, B) \cong \Hom (\bZ^n,\HH^1 (B)) \oplus \Hom (\bZ^m,\HH^2(B))\, .\] 
\end{cor}  

\begin{proof}
  From the definition, any equivalence has to be of the form $\gamma(x,y,z) = (x + u_+(y) + u_-(z), y,z)$ for some $u_+ \in \Hom(\bZ^n, B)$ and $u_- \in \Hom(\bZ^m, B)$, and the requirement that $\gamma$ commutes with the $C_2$-action is exactly that $\Theta_+-\Psi_+ = (\epsilon - 1)u_+$ and $\Theta_- -\Psi_- = (\epsilon + 1)u_-$. Conversely, suppose that $(\Theta_+-\Psi_+)$ is valued in $(\epsilon - 1)(B)$ and $(\Theta_--\Psi_-)$ is valued in $(\epsilon + 1)(B)$. Letting $e_1, ... , e_n$ denote the standard basis of $\bZ^n$, one can define $u_+ \in \Hom(\bZ^n, B)$ with $\Theta_+-\Psi_+ = (\epsilon - 1)u_+$ by setting $u_+(e_i) = b_i$ for $i = 1,...,n$, where $(\Theta_+-\Psi_+)(e_i) = (\epsilon - 1)(b_i)$. Similarly, one obtains $u_+ \in \Hom(\bZ^n, B)$ with $\Theta_+-\Psi_+ = (\epsilon - 1)u_+$, and therefore $\gamma(x,y,z) = (x + u_+(y) + u_-(z), y,z)$ defines the required isomorphism.
  
  For the final claim, note that any $C_2$-module structure on $B \oplus \bZ^n \oplus \bZ^m$ for which the exact sequence
  \[
  0 \rightarrow B \rightarrow B \oplus \bZ^n \oplus \bZ^m \rightarrow \bZ^n_+ \oplus \bZ^m_- \rightarrow 0
  \]
  is a sequence of $C_2$-modules must be of the form $B_{\Theta}(\bZ)$ for some $\Theta = (\Theta_+, \Theta_-)$.
\end{proof}

We can describe the cohomology of the modules $B_{\Theta}(\bZ)$, when $B$ is finite, in terms of the maps
\begin{equation}\label{BTheta_reduced}
\overline{\Theta_+} : \bF_2^n = (\bZ/2\bZ)^n  \rightarrow \HH^1(B) , \ \  \overline{\Theta_- } : \bF_2^m = (\bZ/2\bZ)^m \rightarrow \HH^2 (B).
\end{equation}

\begin{prop}\label{prop:cohomofextensionexplicit} Suppose that $B$ is a finite $C_2$-module. Then
\[
\HH^{0}(B_\Theta (\bZ)) \cong \HH^0 (B) \oplus \bZ^n,
\]
and for $r \geq 1$,
\[
 \dim_{\bF_2} \HH^r (B_\Theta (\bZ)) = 
\begin{cases}
\dim_{\bF_2} \HH^r(B) + n - d , & \mbox{ if $r$ is even}, \\
\dim_{\bF_2} \HH^r(B) + m - d, & \mbox{ if $r$ is odd},
\end{cases}
\]
where $d \coloneqq \rank_{\bF_2} \overline{\Theta_+} + \rank_{\bF_2} \overline{\Theta_-}$.
\end{prop}

\begin{proof}
  Considering the long exact sequence on cohomology, one can use Remark \ref{rem:connectingmaps} to explicitly compute the connecting maps as 
  \[
  \begin{tikzcd}
  0 \rar & \HH^0(B) \rar & \HH^0(B_{\Theta}(\bZ)) \rar & \bZ^n \ar[out=-20, in=160, "\Theta_+"]{dll} \\
  & \HH^1(B) \rar & \HH^1(B_{\Theta}(\bZ)) \rar & \bF_2^m \ar[out=-20, in=160, "\overline{\Theta_-}"]{dll} \\
  & \HH^2(B) \rar & \HH^2(B_{\Theta}(\bZ)) \rar & \bF_2^n \ar[out=-20, in=160, "\overline{\Theta_+}"]{dll}\\
  & \HH^3(B) \rar & \cdots
\end{tikzcd}
\]
We have that $\ker(\Theta_+ \colon \bZ^n \rightarrow \HH^1(B)) \cong \bZ^n$, as $\HH^1(B)$ is finite, and $\bZ^n$ is free, hence $\HH^{0}(B_\Theta (\bZ)) \cong \HH^0 (B) \oplus \bZ^n$. For the remaining formula, note that all terms from $\HH^1(B)$ onwards in the sequence are $\bF_2$-vector spaces by Corollary \ref{cor:cohomisvectorspace}. Therefore if $r$ is odd, we have that 
\[
\HH^r(B_{\Theta}(\bZ)) \cong \HH^r(B) / \overline{\Theta_+}(\bF_2^n) \oplus \ker(\overline{\Theta_-} \colon \bF_2^m \rightarrow \HH^{r+1}(B)),
\]
and thus $\dim_{\bF_2} \HH^r (B_\Theta (\bZ)) = \HH^r(B) - \dim_{\bF_2}(\overline{\Theta_+}(\bF_2^n)) + m - \rank_{\bF_2}(\overline{\Theta_-})$. Similarly for even $r$.
\end{proof}

\section{Classification of $C_2$-modules}
In this section give a classification theorem for finitely-generated $C_2$-modules.
As part of our classification, we show that any such module decomposes uniquely into a direct sum of indecomposable modules, i.e., the Krull-Schmidt Theorem holds.
This is mildly surprising since the Krull-Schmidt Theorem fails for finitely generated $\bZ [G]$-modules for a finite group $G$ in general \cite{ReinerKS}.

We show that the indecomposable $C_2$-modules fall into three categories: $p$-torsion for odd $p$, torsion-free and mixed.
The latter are the modules $B_\Theta (\bZ)$, where $B$ is non-trivial and has order a power of $2$, which have 2-power-torsion and torsion-free parts.
Their classification is not explicit, reducing to the classification of Levy structures, which we introduce below.

\subsection{Torsion-free modules}

We first consider $C_2$-modules which are finite-free as $\bZ$-modules. Examples are the modules $\bZ_d$, $\bZ_+$ and $\bZ_-$, from Section \ref{sect:cohomC2modules}. The following shows that these suffice to describe all such modules.

\begin{prop}\label{prop:torsionfreeKS}
  The $C_2$-modules $\bZ_d$, $\bZ_+$ and $\bZ_-$ are indecomposable, and any finitely generated $C_2$-module $M$ which is free as a $\bZ$-module is of the form
  \[
    M \cong \bZ_d^{r_d} \oplus \bZ_{+}^{r_+} \oplus \bZ_{-}^{r_-}
  \]
  for unique $r_d, r_+, r_- \in \bZ_{\geq 0}$. 
\end{prop}

More generally, one can classify finitely generated $C_p$-modules which are free over $\bZ$, where $p$ is any prime \cite[Thm.\ 74.3, Ex.\ 74(4)]{CR0}. We give a proof of Proposition \ref{prop:torsionfreeKS}, which is much simpler than the proof in the general case. We also note that Proposition \ref{prop:torsionfreeKS} is equivalent to a classification of self-inverse integer matrices, which can be done in a very elementary, but again longer, way \cite{Han}.

\begin{proof}
It is direct to see that $\bZ_+$ and $\bZ_-$ indecomposable, this being true even as $\bZ$-modules. For $\bZ_d$, note that $\bZ_d \cong \bZ[C_2]$ as $\bZ[C_2]$-modules, and the reduction modulo $2$ is $\bF_2[C_2]$. Therefore any decomposition $M_1 \oplus M_2 = \bZ_d$ as $\bZ[C_2]$-modules induces a decomposition $\overline{M_1} \oplus \overline{M_2} = \bF_2[C_2]$ of the indecomposable $\bF_2[C_2]$-module $\bF_2[C_2]$, and therefore one of $M_1$ and $M_2$ must have $\bZ$-rank $2$. Both $M_1$ and $M_2$ are free over $\bZ$, being submodules of a free module, hence one of $M_1$ and $M_2$ must be $0$.

Now suppose that $M$ is a arbitrary $C_2$-module, free over $\bZ$, and consider $X \coloneqq M / M^1$. We first show that $X^1 = 0$. Suppose that $\epsilon(x) = - x$, and let $m \in M$ be a lift of $x$, so $\epsilon(m) = - m + m'$ for some $m' \in M^1$. Then $\epsilon(m') = -m'$, but also
\[
  \epsilon(m') = \epsilon(\epsilon(m) + m) = \epsilon(m) + m = m',
\]
so $2m' = 0$, hence $m' = 0$, as $M$ is torsion-free. In particular, $m \in M^1$, so $x = 0$ in $X$. This shows that $X^1 = 0$, from which it follows that $X^0 = X$: for any $x \in X$, $\epsilon(x) = x$, because
\[
\epsilon(\epsilon(x) - x) = x - \epsilon(x) = - (\epsilon(x) - x).
\]
In summary, $M^1 \cong \bZ_{-}^{a}$ and $X \cong \bZ_{+}^{b}$ for some $a, b \geq 0$, and $M$ fits into a short exact sequence
\[
0 \rightarrow \bZ_{-}^{a} \rightarrow M \rightarrow \bZ_{+}^{b} \rightarrow 0.
\]

Choosing a splitting of $M \rightarrow \bZ_{+}^{b}$ as a $\bZ$-module provides a basis $e_1, ... , e_a, f_1, ... , f_b$ of $M$ for which
\[
\begin{pmatrix}
  -I_a & U \\ 0 & I_b
\end{pmatrix}
\]
is the matrix of $\epsilon$. If $P \in \GL_a(\bZ)$ and $Q \in \GL_b(\bZ)$, then
\[
\begin{pmatrix}
  P & 0 \\ 0 & Q^{-1}
\end{pmatrix}
\begin{pmatrix}
  -I_a & U \\ 0 & I_b
\end{pmatrix}
\begin{pmatrix}
  P^{-1} & 0 \\ 0 & Q
\end{pmatrix} = \begin{pmatrix}
  -I_a & PUQ \\ 0 & I_b
\end{pmatrix},
\]
and therefore by changing basis we may assume that $U$ is in Smith normal form. We therefore have that $M$ is the direct sum of the $C_2$-submodules $\langle e_i, f_i \rangle$ for $1 \leq i \leq \min(a,b)$ and either modules $\langle e_i \rangle \cong \bZ_-$ for $i > b$, if $a>b$, or modules $\langle f_i \rangle \cong \bZ_+$ for $i>a$, if $a<b$. It therefore remains to understand the submodules $M_i = \langle e_i, f_i \rangle$ for $1 \leq i \leq \min(a,b)$, where $\epsilon(e_i) = -e_i$ and $\epsilon(f_i) = u_{ii} e_i + f_i$.

If $u_{ii} = 2k$ is even, then with respect to the basis $e_i' \coloneqq e_i, f_i' \coloneqq f_i + k e_i$, $\epsilon(e_i') = - e_i'$ and
\[
\epsilon(f_i') = \epsilon(f_i + k e_i) = u_{ii} e_i + f_i - k e_i = f_i + k e_i = f_i',
\]
hence $M_i \cong \bZ_+ \oplus \bZ_-$.

If instead $u_{ii} = 2k + 1$ is odd, then with respect to the basis $e_i' \coloneqq e_i + (f_i + k e_i), f_i' \coloneqq f_i + k e_i$,
\begin{align*}
\epsilon(f_i') &= \epsilon(f_i + k e_i) = u_{ii} e_i + f_i - k e_i = (f_i + k e_i) + e_i = e_i', \\
\epsilon(e_i') &= \epsilon((f_i + k e_i) + e_i) = e_i' - e_i = f_i + k e_i = f_i',
\end{align*}
and so $M_i \cong \bZ_d$.

To see that this description is unique, suppose that 
\[
\bZ_d^{r_d} \oplus \bZ_{+}^{r_+} \oplus \bZ_{-}^{r_-} \cong \bZ_d^{s_d} \oplus \bZ_{+}^{s_+} \oplus \bZ_{-}^{s_-}
\]
for $r_d, r_+, r_-, s_d, s_+, s_- \in \bZ_{\geq 0}$. Applying $-\otimes_{\bZ} \bQ$, $r_+ + r_d = s_+ + s_d$ and $r_- + r_d = s_- + s_d$, and applying $-\otimes_{\bZ} \bF_2$, $r_d = s_d$, as $\bZ_+$ and $\bZ_-$ reduce to the trivial $\bF_2[C_2]$-module, and $\bZ_d \cong \bZ[C_2]$ reduces to the indecomposable module $\bF_2[C_2]$. Therefore $r_+ = s_+$, $r_- = s_-$, and $r_d = s_d$.
\end{proof}

\subsection{Levy Structures}

We now introduce the notion of a Levy structure, which will be used below to understand general $C_2$-modules, which aren't necessarily torsion-free.

Finitely generated $\bZ [C_p]$-modules, where $p$ is prime, were classified by Lawrence Levy \cite{Levy81,Levy85,Levy85a}.
For $p=2$, his classification relies on and is stated in terms of the presentation 
\[ \bZ [C_2] \ \cong \ 
\bZ \times_{\bF_2} \bZ \  = \
\{ (x,y) \in \bZ \times \bZ \,\mid\, x+y \mbox{ is even}\}.
\]
Whilst our classification method is different and more explicit, the core of both methods is to reduce the classification
to some linear algebraic data. To acknowledge Levy's contribution, we make the following definition.

\begin{defn}
A Levy structure is a triple $Q=(B, V_+, V_-, {\theta_+},{\theta_-})$, where
\begin{itemize}
  \item $B$ is a finite $C_2$-module,
  \item $V_+, V_-$ are finite-dimensional $\bF_2$-vector spaces,
  \item $\theta_+ \colon V_+ \rightarrow \HH^1(B)$, $\theta_- \colon V_- \rightarrow \HH^2(B)$ are $\bF_2$-linear maps.
\end{itemize}
\end{defn}

Levy structures form a category $\sL$. The morphisms are triples
\[f=(f_0,f_+,f_-): (B, V_+, V_-, {\theta_+}, {\theta_-} ) \to (B', V_+', V_-', {\eta_+}, {\eta_-} ),\]
where $f_0: B\to B'$ is a homomorphism of $C_2$-modules and $f_+ : V_+ \to V_+'$,
$f_- : V_- \to V_-'$ are $\bF_2$-linear maps, making the following diagrams commutative:
\begin{equation} \label{equality_1}
  \begin{tikzcd}
	V_+ & \HH^1(B)  \\
	V_+' & \HH^1(B') 
	\arrow["\theta_+"', from=1-1, to=1-2]
	\arrow["H^1(f_0)", from=1-2, to=2-2]
	\arrow["\eta_+", from=2-1, to=2-2]
	\arrow["f_+"', from=1-1, to=2-1]
  \end{tikzcd}
\quad \quad
  \begin{tikzcd}
	V_- & \HH^2(B)  \\
	V_-' & \HH^2(B') 
	\arrow["\theta_-"', from=1-1, to=1-2]
	\arrow["H^2(f_0)", from=1-2, to=2-2]
	\arrow["\eta_-", from=2-1, to=2-2]
	\arrow["f_-"', from=1-1, to=2-1]
  \end{tikzcd}
\end{equation}
The category 
$\sL$ is additive and has direct sums
\[(B, V_+, V_-, {\theta_+}, {\theta_-} ) \oplus (B', V_+', V_-', {\eta_+}, {\eta_-} ) \coloneqq (B\oplus B', V_+ \oplus V_+', V_- \oplus V_-', {\theta_+}\oplus{\eta_+}, {\theta_-}\oplus{\eta_-} )\]
so that we can talk about indecomposable objects. Moreover, $\sL$ has split idempotents, i.e., for $Q \in \sL$, every $e=e^2 \in \End_{\sL}(Q)$ yields a direct sum decomposition of $Q$. 
\begin{prop} \label{sLk}
  The category $\sL$ is a Krull-Schmidt category.
\end{prop}

  Recall that a Krull-Schmidt category is an additive category, where every object decomposes into a finite direct sum of objects
  with local endomorphism rings \cite{Krause}.

\begin{proof}
Given the finite nature of  $\sL$, every object there decomposes into a direct sum of indecomposable objects. Suppose $X\in\sL$ is indecomposable.
  Then the ring $\End(Q) = \sL (Q,Q)$ has no non-trivial idempotents.
  On the other hand, $\End(Q)$ is a finite ring. It remains to observe that a finite ring $R$ without non-trivial idempotents is necessarily local.

  Recall that a ring $R$ is local if the sum of two non-units is again non-unit.
  The Jacobson radical $J(R)$ is nilpotent since $R$ is finite.
  Hence, it suffices to observe that $R/J(R)$ is a division ring. Suppose not. Then it has a non-trivial idempotent $e\in R/J(R)$, $e\not\in\{0,1\}$.
  But an idempotent can be lifted modulo a nilpotent ideal. This gives a non-trivial idempotent $e_1 \in R$.
\end{proof}  

\begin{cor} \cite[Theorem 4.2]{Krause}\label{cor:KSforLevyStr}
  Consider $X\in \sL$ and two decompositions
  $X= X_1 \oplus \ldots \oplus X_r= Y_1 \oplus \ldots \oplus Y_s$ into non-zero indecomposable objects.
  Then $r=s$ and there exists a permutation $\sigma$ such that $X_i \cong Y_{\sigma (i)}$ for all $i$.
\end{cor}

We can classify indecomposable Levy structures into three types.

\begin{lemma}\label{lem:indecomlevystr}
  Any indecomposable object $X \in \sL$ is of one of the following types:
  \begin{enumerate}
    \item $X \cong (B,0,0,0,0)$, where $B = (\bZ / p^r \bZ)_+$ or $B = (\bZ / p^r \bZ)_-$, for $r \geq 1$ and $p$ an odd prime.
    \item $X \cong (0,\bF_2,0,0,0)$ or $X = (0, 0, \bF_2,0,0)$.
    \item $X \cong (B, V_+, V_-, {\theta_+},{\theta_-})$, indecomposable, where $B$ is $2$-power-torsion and non-zero.
  \end{enumerate}
\end{lemma}

\begin{proof}
  Let $X = (B, V_+, V_-, {\theta_+},{\theta_-})$ be indecomposable, and first suppose that $B_{(p)} \neq 0$ for some odd prime $p$. Then $\HH^1(B_{(p)}) = 0$ by Lemma \ref{cohom_div_by2}, and so $\theta_+ \colon V_+ \rightarrow \HH^1(B)$ is valued in 
  \[
  \HH^1(\oplus_{l \neq p} B_{(l)}) \subset \HH^1(\oplus_{l \neq p} B_{(l)}) \oplus \HH^1(B_{(p)}) = \HH^1(B).
  \]
  Similarly for $\HH^2(B_{(p)})$ and $\theta_- \colon V_- \rightarrow \HH^2(B)$, and therefore 
  \[
    (B_{(p)}, 0, 0, 0, 0) \oplus (\oplus_{l \neq p} B_{(l)}, V_+, V_-, \theta_+, \theta_-) \cong (B, V_+, V_-, {\theta_+},{\theta_-})
  \]
  is a direct sum decomposition in $\sL$. As $B_{(p)} \neq 0$, this forces $B_{(l)} = 0$ for any prime $l \neq p$, and therefore $B_{(p)} = B$. Now $\HH^1(B) = 0$ and $\HH^2(B) = 0$, again by Lemma \ref{cohom_div_by2}, so $\theta_+ = 0$, $\theta_- = 0$, and
  \[
  (B_{(p)}, 0, 0, 0, 0) \oplus (0, V_+, V_-, 0, 0) \cong (B, V_+, V_-, {\theta_+},{\theta_-})
  \]
  is a direct sum decomposition in $\sL$, and hence $V_+ = 0$ and $V_- = 0$, as $B \neq 0$. Because $B$ admits division by $2$, $B = (B^0)_+ \oplus (B^1)_-$ by Lemma \ref{cohom_div_by2}, and therefore 
  \[
  X = ((B^0)_+,0,0,0,0) \oplus ((B^1)_-,0,0,0,0)
  \]
  which forces $B = A_+$ or $B = A_-$, where $A$ is the underlying abelian group of $B$. In either case, any proper non-trivial decomposition of $A$ as an abelian group induces a proper non-trivial decomposition of $B$ as a $C_2$-module, so $A$ must be indecomposable as an abelian group, hence cyclic of odd prime power order, and we see that $X$ is of type (1).
  
  If instead $B_{(p)} = 0$ for all odd primes, then $B = B_{(2)}$. If $B=0$, then $\theta_+ = 0$ and $\theta_- = 0$, and so any decomposition of either $V_+$ or $V_-$ into proper non-trivial subspaces induces a decomposition of $X$. In particular, $X = (0,\bF_2,0,0,0)$ or $X = (0, 0, \bF_2,0,0)$, by the indecomposability of $X$, and $X$ is of type (2). If instead $B \neq 0$, then $B = B_{(2)}$, and $X$ is of type (3). 
\end{proof}

\subsection{Mixed Modules}

We are now interested in the following class of modules. For a $C_2$-module $M$, we write $M_T$ for the $\bZ$-torsion subgroup, which is always a $C_2$-submodule. By Proposition \ref{prop:torsionfreeKS}, we know that the $C_2$-module $M / M_T$ is isomorphic to a direct sum of copies of $\bZ_+$, $\bZ_-$ and $\bZ_d$.

\begin{lemma}\label{lem:mixedequiv}
  Let $M$ be a finitely generated $C_2$-module. The following are equivalent:
  \begin{enumerate}
    \item $M / M_T \cong \bZ_+^n \oplus \bZ^m_-$ for some $m,n \geq 0$,
    \item $M$ is an extension of $\bZ_+^n \oplus \bZ^m_-$ by a finite $C_2$-module, for some $m,n \geq 0$,
    \item $M \cong B_{\Theta}(\bZ)$ for some $\Theta$ and finite $C_2$-module $B$.
  \end{enumerate}
\end{lemma}

\begin{proof}
  The fact that (1) implies (2) follows from the exact sequence $0 \rightarrow M_T \rightarrow M \rightarrow M / M_T \rightarrow 0$, and the fact that (2) implies (3) from Corollary \ref{even_part_ext_new}. The fact that (3) implies (1) is direct from the definition of the $C_2$-module $B_{\Theta}(\bZ)$.
\end{proof}

\begin{defn}
A \emph{mixed module} is a finitely generated $C_2$-module which satisfies the equivalent conditions of Lemma \ref{lem:mixedequiv}
\end{defn}

Let $\sM$ denote the full subcategory of $C_2$-modules, with objects the mixed modules. We now define a functor
\[
\sL(-) \colon \sM \rightarrow \sL
\]
to the category of Levy structures, which will use to deduce properties of $\sM$ from properties of $\sL$.

Recall the invariants $(-)^0$ and skew-invariants $(-)^1$, as defined at the start of Section \ref{sect:cohomC2modules}.
\begin{defn}
  For $M \in \sM$, we set
  \[
  M(0) \coloneqq (M / M_T)^0 \in \Ab, \qquad M(1) \coloneqq (M / M_T)^1 \in \Ab.
  \]
\end{defn}

\begin{remark}
  We note that for any $M \in \sM$, $M / M_T = M(0) \oplus M(1)$ as abelian groups.
\end{remark}

To define $\sL(-)\colon \sM \rightarrow \sL$ on objects, let $M$ be a mixed module, and let $s_M \colon M / M_T \rightarrow M$ be a splitting of $\pi_M \colon M \rightarrow M / M_T$. We note that this is equivalent to the data of an abelian group isomorphism $\gamma_M \colon M \xrightarrow{\sim} M_T \oplus M / M_T$ such that the diagram
\[\begin{tikzcd}
	0 & {M_T} & M & {M / M_T} & 0 \\
	0 & {M_T} & {M_T \oplus M / M_T} & {M / M_T} & 0
	\arrow[from=1-1, to=1-2]
	\arrow[from=1-2, to=1-3]
	\arrow["\id", from=1-2, to=2-2]
	\arrow[from=1-3, to=1-4]
	\arrow[from=1-4, to=1-5]
	\arrow["\gamma_M", from=1-3, to=2-3]
	\arrow["\id", from=1-4, to=2-4]
	\arrow[from=2-1, to=2-2]
	\arrow[from=2-2, to=2-3]
	\arrow[from=2-3, to=2-4]
	\arrow[from=2-4, to=2-5]
\end{tikzcd}\]
commutes.

\begin{defn}
  For $M \in \sM$, we set
  \[
  \overline{M(0)} \coloneqq M(0) / 2 M(0) \in \Vect_{\bF_2}, \qquad \overline{M(1)} \coloneqq M(1) / 2 M(1) \in~\Vect_{\bF_2}.
  \]
\end{defn}

Using a fixed splitting $s_M$, we can make the following definition.

\begin{defn}
  For $M \in \sM$, we set
  \[
  \sL(M) \coloneqq (M_T, \overline{M(0)}, \overline{M(1)}, \theta_+, \theta_-),
  \]
  where
  \begin{align*}
    \widetilde{\theta}_+ \colon M(0) \rightarrow M_T^1, \qquad \widetilde{\theta}_+ = \pi_{M_T} \circ \gamma_M \circ \epsilon_M \circ s_M, \\
    \widetilde{\theta}_- \colon M(1) \rightarrow M_T^0, \qquad \widetilde{\theta}_- = \pi_{M_T} \circ \gamma_M \circ \epsilon_M \circ s_M,
  \end{align*}
  and
  \begin{align*}
    \theta_+ \colon \overline{M(0)} \rightarrow M_T^1 / (\epsilon - 1)M_T = \HH^1(M_T), \\
    \theta_- \colon \overline{M(1)} \rightarrow M_T^0 / (\epsilon + 1)M_T = \HH^2(M_T),
  \end{align*}
  are the induced maps. If $f \colon M \rightarrow N$ is a morphism in $\sM$, we set
  \[
  \sL(f) = (f_T, f(0), f(1)),
  \]
  where
  \[
    f_T \colon M_T \rightarrow N_T, \qquad
    f(0) \colon \overline{M(0)} \rightarrow \overline{N(0)}, \qquad
    f(1) \colon \overline{M(1)} \rightarrow \overline{N(1)},
  \]
  are the natural maps induced by $f$.
\end{defn}

Let us verify that this is well-defined.

\begin{lemma}
  For $M \in \sM$, the Levy structure $\sL(M)$ is independent of the choice of splitting $s_M$, and $\sL(f) \colon \sL(M) \rightarrow \sL(N)$ is a well-defined morphism in $\sL$.
\end{lemma}

\begin{proof}
  Suppose that $s_M, s_M' \colon M / M_T \rightarrow M$ are two choices of splitting, with corresponding isomorphisms $\gamma_M, \gamma_M' \colon M \xrightarrow{\sim} M_T \oplus M/ M_T$ respectively. An alternative description of the maps $\widetilde{\theta}_+, \widetilde{\theta}_+' \colon M(0) \rightarrow M_T^1$, are that $\widetilde{\theta}_+(x) = \pi_{M_T} (\epsilon_*(0,x))$ and $\widetilde{\theta}_+'(x) = \pi_{M_T} (\epsilon_*'(0,x))$, where $\epsilon_* = \gamma_M \circ \epsilon_M \circ \gamma_M^{-1}$ and $\epsilon_*' = \gamma_M' \circ \epsilon_M \circ (\gamma_M')^{-1}$ are the transported actions of $\epsilon \in \bZ[C_2]$ on $M_T \oplus M / M_T$. Because these are defined to make $\gamma_M$ and $\gamma_M'$ both $\bZ[C_2$]-equivariant, we have that $\phi \circ \epsilon_* = \epsilon_*' \circ \phi$, where $\phi \coloneqq \gamma_M' \circ \gamma_M^{-1}$. Now $\phi$ is an abelian group automorphism of $M_T \oplus M / M_T$ which is the identity on $M_T$ and the identity after projecting to $M / M_T$, are therefore $\phi$ and $\phi^{-1}$ have the form
  \[
  \phi(m, x) = (m + u(x), x), \qquad \phi^{-1}(m,x) = (m - u(x), x)
  \]
  for some $u \in \Hom_{\bZ}(M/M_T, M_T)$. We can therefore compute that for $x \in M(0)$, 
  \begin{align*}
    \widetilde{\theta}_+(x) &= \pi_{M_T} (\epsilon_*(0,x)), \\
    &= \pi_{M_T}((\phi^{-1} \circ \epsilon_*' \circ \phi)(0,x)), \\
    &= \pi_{M_T}(\phi^{-1} (\epsilon_*'(u(x), x))), \\
    &= \pi_{M_T}(\phi^{-1} (\epsilon_M(u(x)) + \widetilde{\theta}_+'(x), x)), \\
    &= \pi_{M_T}(\epsilon_M(u(x)) + \widetilde{\theta}_+'(x) - u(x), x), \\
    &= \widetilde{\theta}_+'(x) + (\epsilon_M - 1)u(x),
  \end{align*}
  and therefore the induced maps $\theta_+$ and $\theta_+'$ are equal. Taking instead $x \in M(1)$, from an analogous computation one obtains that $\theta_- = \theta_-'$. 
  
  For the second claim, given $f \colon M \rightarrow N$ a morphism in $\sM$, to check that $\sL(f) \colon \sL(M) \rightarrow \sL(N)$ is a morphism in $\sL$ one needs to show that $\theta_{+,N} \circ f(0) = \HH^1(f_T) \circ \theta_{+,M}$ and $\theta_{+,N} \circ f(1) = \HH^2(f_T) \circ \theta_{+,M}$. We show the first, the second following mutatis mutandis. The induced map
  \[
  \widetilde{f} \coloneqq \gamma_N \circ f \circ \gamma_M^{-1} \colon M_T \oplus M/M_T \rightarrow N_T \oplus N / N_T
  \]
  is $\bZ[C_2]$-linear when $M_T \oplus M/M_T$ and $N_T \oplus N / N_T$ are given the transported actions $\epsilon_{*,M} = \gamma_M \circ \epsilon_M \circ \gamma_M^{-1}$ and $\epsilon_{*,N} = \gamma_N \circ \epsilon_N \circ \gamma_N^{-1}$ respectively as above. Similarly to above, because $M_T$ is torsion and $N / N_T$ is torsion free, the map $\widetilde{f}$ has the form
  \[
  \widetilde{f}(m,x) = (f_T(m) + u(x), f_{M/M_T}(x))
  \]
  for some $u \in \Hom_{\bZ}(M/M_T, N_T)$. We can then compute that for $x \in M(0)$,
  \[
    (\widetilde{f} \circ \epsilon_{*,M})(0,x) = \widetilde{f}(\widetilde{\theta}_{+, M}(x),x) = (f_T(\widetilde{\theta}_{+, M}(x)) + u(x), f_{M/M_T}(x)),
  \]
  whereas on the other hand
  \[
  (\epsilon_{*,N} \circ \widetilde{f})(0,x) = \epsilon_{*,N}(u(x), f_{M/M_T}(x)) = (\epsilon_N(u(x)) + \widetilde{\theta}_{+,N}(f_{M/M_T}(x)), f_{M/M_T}(x)).
  \]
  Then comparing the first coefficient,
  \[
  f_T(\widetilde{\theta}_{+, M}(x)) - \widetilde{\theta}_{+,N}(f_{M/M_T}(x))  = (\epsilon_N - 1)(u(x))
  \]
  which when taking the induced map on the quotient implies that $\HH^1(f_T) \circ \theta_{+,M} - \theta_{+,N} \circ f(0) = 0$.
\end{proof}

\begin{eg}\label{eg:stdexample}
Suppose that $M = B_{\Theta}(\bZ)$, where $M / M_T$ is canonically identified with $\bZ_+^n \oplus \bZ_-^m$, and take $s_M \colon \bZ_+^n \oplus \bZ_-^m \rightarrow B_{\Theta}(\bZ)$ to be the standard splitting. Then $M(0) = \bZ_+^n$, $M(1) = \bZ_-^m$, which induce identifications $\overline{M(0)} = \bF_2^n$ and $\overline{M(1)} = \bF_2^m$. The maps $\theta_+ \colon \bF_2^n \rightarrow \HH^1(B)$ and $\theta_- \colon \bF_2^m \rightarrow \HH^2(B)$ are simply those induced by $\Theta_+ \colon \bZ^n \rightarrow B^1$ and $\Theta_- \colon \bZ^m \rightarrow B^0$.
\end{eg}

We can now describe how $\sM$ and $\sL$ are related.

\begin{prop}\label{prop:Lfunctorproperties}
  The functor $\sL(-) \colon \sM \rightarrow \sL$ satisfies the following properties.
  \begin{enumerate}
    \item $\sL(-) \colon \sM \rightarrow \sL$ is additive, commutes with direct sums, and is essentially surjective.
    \item If $M, N \in \sM$, then $M \cong N$ in $\sM$ if and only if $\sL(M) \cong \sL(N)$ in $\sL$.
    \item $M$ is indecomposable in $\sM$ if and only if $\sL(M)$ is indecomposable in $\sL$.
  \end{enumerate}
\end{prop}

\begin{remark}\label{rem:notionsofdecomp}
  Note that $\sM$ is closed under direct summands in the category of all $C_2$-modules, and therefore indecomposability in $\sM$ is equivalent to indecomposability as a $C_2$-module. Indeed, if $M \in \sM$ and $M = M_1 \oplus M_2$ is a decomposition as $C_2$-modules, then $M_1 / M_{1, T} \oplus M_2 / M_{2, T} \cong M / M_T$, and therefore both $M_1$ and $M_2$ are mixed modules, by Proposition \ref{prop:torsionfreeKS}.
\end{remark}

\begin{proof}
  It is clear that $\sL(-)$ is additive and commutes with direct sums. For the essential surjectivity, let $(B, V_+, V_-, \theta_+, \theta_-) \in \sL$, and choose a basis for $V_+$ and $V_-$. Any lift of the maps $\theta_+ \colon \bF_2^n \rightarrow \HH^1(B)$ and $\theta_- \colon \bF_2^m \rightarrow \HH^2(B)$ to $\Theta_+ \colon \bZ^n \rightarrow B^1$ and $\Theta_- \colon \bZ^m \rightarrow B^0$ respectively allow us to define $B_{\Theta}(\bZ) \in \sM$, with the property that $\sL(\sM) = (B, V_+, V_-, \theta_+, \theta_-)$, by Example \ref{eg:stdexample}.

  For the second point, note that by functoriality $\sL(M) \cong \sL(N)$ whenever $M \cong N$, so suppose conversely that $\sL(M) \cong \sL(N)$. Then $M_T \cong N_T$, $\dim_{\bF_2}\overline{M(0)} = \dim_{\bF_2}\overline{M(0)}$ and $\dim_{\bF_2}\overline{M(0)} = \dim_{\bF_2}\overline{M(0)}$, and we may therefore pick isomorphisms $M \cong B_{\Theta}(\bZ)$, $N \cong B_{\Psi}(\bZ)$ for the same finite $C_2$-module $B$ and $\Theta_+, \Psi_+ \colon \bZ^n \rightarrow B$, $\Theta_-, \Psi_- \colon \bZ^m \rightarrow B$ with the same domains. Fixing an isomorphism $(f_0,f_+,f_-) \colon \sL(B_{\Theta}(\bZ)) \rightarrow \sL(B_{\Psi}(\bZ))$, we may lift the isomorphisms $f_+ \colon \bF_2^n \rightarrow \bF_2^n$ and $f_- \colon \bF_2^m \rightarrow \bF_2^m$ to isomorphisms $S_+ \colon \bZ^n \rightarrow \bZ^n$ and $S_- \colon \bZ^m \rightarrow \bZ^m$, using that $\GL_n(\bZ) \rightarrow \GL_2(\bF_2)$ is surjective, the latter being generated by transvections. 

  Writing $\sL(B_{\Theta}(\bZ)) = (B, \bF_2^n, \bF_2^m, \theta_+, \theta_-)$, $\sL(B_{\Psi}(\bZ)) = (B, \bF_2^n, \bF_2^m, \psi_+, \psi_-)$, we have that 
  \begin{align*}
    \delta_+ &\coloneqq \Psi_+ \circ S_+ - f_0 \circ \Theta_+ \colon \bZ^n \rightarrow B^1, \\
    \delta_- &\coloneqq \Psi_- \circ S_- - f_0 \circ \Theta_- \colon \bZ^m \rightarrow B^0
  \end{align*}
  reduce to
  
  \begin{align*}
    \psi_+ \circ f_+ - \HH^1(f_0) \circ \theta_+ &\colon \bF_2^n \rightarrow \HH^1(B), \\
    \psi_- \circ f_- - \HH^2(f_0) \circ \theta_+ &\colon \bF_2^n \rightarrow \HH^2(B)
  \end{align*}
  respectively, and therefore $\delta_+$ and $\delta_-$ are valued in $(\epsilon - 1)B$ and $(\epsilon + 1)B$ respectively. Defining
  \begin{align*}
    \Psi_+' &\coloneqq \Psi_+ - \delta_+ \circ S_+^{-1} \colon \bZ^n \rightarrow B^1, \\
    \Psi_-' &\coloneqq \Psi_- - \delta_- \circ S_-^{-1} \colon \bZ^m \rightarrow B^0,
  \end{align*}
  we have that $\Psi_+ - \Psi_+'$ is valued in $(\epsilon - 1)B$ and $\Psi_- - \Psi_-'$ is valued in $(\epsilon + 1)B$, hence $B_{\Psi}(\bZ) \cong B_{\Psi'}(\bZ)$, even as extensions, by Corollary \ref{even_part_ext_new}. It is thus sufficient to give an isomorphism $\phi \colon B_{\Theta}(\bZ) \xrightarrow{\sim} B_{\Psi'}(\bZ)$, which we do by taking $\phi$ to be explicitly defined by
  \[
  \phi(b,x,y) \coloneqq (f_0(b),S_+(x), S_-(y)).
  \]
  which is a $C_2$-module isomorphism, using the fact that $\Psi_+' \circ S_+ = f_0 \circ \Theta_+$ and $\Psi_-' \circ S_- = f_0 \circ \Theta_-$.

  For the third and final claim, note that any proper direct sum decomposition of $M$ induces a proper direct sum decomposition of $\sL(M)$. Conversely, if $\sL(M) \cong \sL_1 \oplus \sL_2$ is a proper direct sum decomposition of $\sL(M)$, then taking $M_1, M_2 \in \sM$ with $\sL(M_1) \cong \sL_1$, $\sL(M_2) \cong \sL_2$, which is possible by the essential surjectivity of $\sL(-)$, we have that $\sL(M) \cong \sL(M_1) \oplus \sL(M_2) \cong \sL(M_1 \oplus M_2)$, and therefore, using part (2), $M \cong M_1 \oplus M_2$ is a proper direct sum decomposition of $M$.
\end{proof}

\begin{cor}\label{cor:KSformixedmodules}
Any non-zero mixed module $M$ can be written as a direct sum of indecomposable $C_2$-modules. This direct sum description is unique up to isomorphism and permutation, and each direct-summand is a mixed module.
\end{cor}

\begin{proof}
  For $M \in \sM$, we may write $\sL(M) = \sL_1 \oplus \cdots \oplus \sL_r$ as a direct sum of indecomposable objects of $\sL$ by Corollary \ref{cor:KSforLevyStr}, and lift them to $M_1, ... , M_r \in \sM$ by Proposition \ref{prop:Lfunctorproperties}(1). Because $\sL(-)$ commutes with direct sums, $\sL(M) \cong \sL(M_1) \oplus \cdots \oplus \sL(M_r)$, and therefore $M \cong M_1 \oplus \cdots \oplus M_r$ by Proposition \ref{prop:Lfunctorproperties}(2), with each $M_i$ indecomposable by Proposition \ref{prop:Lfunctorproperties}(3) and Remark \ref{rem:notionsofdecomp}. 

  If $M_1 \oplus \cdots \oplus M_r \cong M \cong N_1 \oplus \cdots \oplus N_s$ are two indecomposable decompositions of $M$, then $\sL(M_1) \oplus \cdots \oplus \sL(M_r) \cong \sL(M) \cong \sL(N_1) \oplus \cdots \oplus \sL(N_r)$ are indecomposable decompositions in $\sL$, by Proposition \ref{prop:Lfunctorproperties}(3), and therefore the uniqueness follows from the uniqueness of Corollary \ref{cor:KSforLevyStr} and Proposition \ref{prop:Lfunctorproperties}(2). The final claim follows from Remark \ref{rem:notionsofdecomp}.
\end{proof}

\subsection{General classification} We can now state our classification of all $C_2$-modules.

\begin{thm}\label{thm:generalclassification} 
  Any finitely generated $C_2$-module is a direct sum of indecomposable finitely generated $C_2$-modules, and this direct sum decomposition is unique up to isomorphism and permutation of the factors. The indecomposable $C_2$-modules fall into one of the following types.
\begin{enumerate}
\item $p$-Torsion: $B_+$ and $B_-$, where $B$ is a cyclic group of order $p^n$, $n \geq 1$, $p$ odd.
\item Torsion-free: $\bZ_d$, $\bZ_+$ and $\bZ_-$.
\item Mixed: $B_\Theta(\bZ)$ for a $2$-power-torsion $C_2$-module $B \neq 0$, with $\sL(B_\Theta(\bZ))$ indecomposable.
\end{enumerate}
\end{thm}

\begin{remark}
  Note that the isomorphism class of the modules $B_\Theta(\bZ)$ depends only on the isomorphism class of the associated Levy structure $\sL(B_\Theta(\bZ))$ in $\sL$, by Proposition \ref{prop:Lfunctorproperties}(2).
\end{remark}

\begin{proof}
  Suppose that $A$ is a finitely generated $C_2$-module. From Proposition \ref{prop:torsionfreeKS}, $A/A_T \cong \bZ_+^n \oplus \bZ_-^m \oplus \bZ_d^r$ for some $r,m,n \geq 0$. The module $\bZ_d \cong \bZ[C_2]$ is free, and therefore we may split the projection $A \twoheadrightarrow \bZ_d^r$ to write $A \cong M \oplus \bZ_d^r$ for some $C_2$-submodule $M$ of $A$. We have that $M / M_T \cong \bZ_+^n \oplus \bZ_-^m$, and therefore $M$ is a mixed module. In particular, $M$ decomposes into indecomposable $C_2$-modules by Corollary \ref{cor:KSformixedmodules}, and therefore so does $A$. 

  To see that all indecomposable $C_2$-modules are included in our list, and all modules on our list our indecomposable, note that our list consists of exactly the $C_2$-module $\bZ_d$, indecomposable by Proposition \ref{prop:torsionfreeKS}, and the indecomposable $C_2$-modules contained in $\sM$, by Lemma \ref{lem:mixedequiv}, Proposition \ref{prop:Lfunctorproperties}(3) and Remark \ref{rem:notionsofdecomp}, as these are exactly the mixed modules reducing to the indecomposable Levy structures of Lemma \ref{lem:indecomlevystr}. To see that any finitely generated indecomposable $C_2$-module is either $\bZ_d$ or contained in $\sM$, note that if $A$ is a finitely generated indecomposable $C_2$-module, then we can write $A \cong M \oplus \bZ_d^r$ as above. If $M = 0$, then $r = 1$ and $A = \bZ_d$ by Proposition \ref{prop:torsionfreeKS}, and if $M \neq 0$, then $r = 0$ and $A = M$ is a mixed module.

  Finally, to show that the decomposition into indecomposables is unique, suppose that we have two descriptions of $A$ as a direct sum of indecomposable modules on our list. Grouping together the terms other than $\bZ_d$, we therefore have an isomorphism
  \[
    M \oplus \bZ_d^r \cong A \cong N \oplus \bZ_d^s
  \]
for mixed modules $M,N$ and $r,s \geq 0$. From the uniqueness of the decomposition of the torsion-free $C_2$-module $A / A_T$ (Proposition \ref{prop:torsionfreeKS}), we have that $r = s$, and therefore
\[
M^{\oplus n} \cong N^{\oplus n}
\]
for some $n \geq 1$ because $\bZ[C_2]$, the endomorphism ring of $\bZ_d$, is a commutative $\bZ$-algebra which is finite-free over $\bZ$ \cite[Cor.\ 2.2, Thm.\ 3.10]{GOO}. Then writing $M$ and $N$ as their indecomposable decompositions, these are the same up to isomorphism and permutation by Corollary \ref{cor:KSformixedmodules} as $M$ and $N$ are mixed, and therefore the same is true of the two indecomposable decompositions of $A$.
\end{proof}

In light of Theorem \ref{thm:generalclassification}, to understand general $C_2$-modules one is therefore reduced to understanding when the $C_2$-module $B_{\Theta}(\bZ)$ is indecomposable for a finite non-trivial $2$-power-torsion $C_2$-module $B$. From Proposition \ref{prop:Lfunctorproperties}(3) this is equivalent to the following problem.

\begin{prob}
Classify indecomposable Levy structures $(B, V_+, V_-, \theta_+, \theta_-)$ for which $B$ is a finite non-trivial $2$-power-torsion $C_2$-module.
\end{prob}

For example, is it possible to characterise the indecomposability of $(B, V_+, V_-, \theta_+, \theta_-)$ in terms of the indecomposability of $B$ and conditions on the maps $\theta_+$ and $\theta_-$?

%% file: appendix-gradedalgebras.tex
\chapter{$C_2$-Graded Algebras and the Brauer-Wall Group}\label{appendix:gradedalgebras}

In this appendix we collect together some basic results on $C_2$-graded algebras which we make use of in Sections \ref{sect:BrauerWall} and \ref{sect:recoverDyson10}, which we were unable to find in the literature.

We first recall some basic definitions, all of which can be found in \cite[Chap.\ IV]{LAM}. Throughout, $\bK$ is a field with $\charfield(\bK) \neq 2$.

\begin{defn}
    By a \emph{$C_2$-graded algebra} over $\bK$, we mean a finite-dimensional $\bK$-algebra $A$, with $\bK$-subspaces $A_0, A_1 \leq A$, such that $\bK \subset A_0$, $A = A_0 \oplus A_1$, and $A_i \cdot A_j \subset A_{i+j}$, where the subscripts are read in $\bZ / 2 \bZ$.

    A $\bK$-subspace $S \leq A$ of a $C_2$-graded algebra $A$ is called \emph{graded} if $S = S \cap A_0 + S \cap A_1$. 

    A $C_2$-graded algebra $A$ is called \emph{simple} if $A$ has no proper non-zero two-sided graded ideals, and called a \emph{$C_2$-graded division algebra} if every non-zero homogeneous element is invertible.

    The \emph{graded centre} of a $C_2$-graded algebra $A$ is the graded $\bK$-algebra $\widehat{Z}(A)$ such that 
    \begin{align*}
        \widehat{Z}(A)_0 &= \{a \in A_0 \mid ab = ba \text{ for all } b \in A\}, \\
        \widehat{Z}(A)_1 &= \{a \in A_1 \mid ab = ba \text{ for all } b \in A_0, \: ab = - ba \text{ for all } b \in A_1\},
    \end{align*}
and $A$ is called \emph{central} if $\widehat{Z}(A) = \bK$.
\end{defn}

We will use the following basic facts.

\begin{lemma}\label{lem:centreofC2graddiv}
    Suppose that $A$ is a $C_2$-graded simple algebra over $\bK$. Then $\widehat{Z}(A)_1=0$, so
        \[
            \widehat{Z}(A) = A_0 \cap Z(A) \leq Z(A_0).
        \] 
        In particular, if $D$ is a $C_2$-graded division algebra, then $\widehat{Z}(D)$ is a finite field extension of $\bK$.
\end{lemma}

\begin{proof}
    As $\widehat{Z}(A)$ is a graded subalgebra, it suffices to show that $\widehat{Z}(A) \cap A_1 = 0$. Suppose that $z \in \widehat{Z}(A)_1$ and $z \neq 0$. Then $(z)$, the two-sided ideal of $A$ generated by $z$, is homogeneous (as $z$ is homogeneous) and non-trivial, hence $(z) = A$. Furthermore, from the condition that $z$ either commutes or anticommutes with homogeneous elements of $A$, $(z) = Az = \{a z \mid a \in A\}$, so $z$ is invertible. On the other hand, that $z$ is in $\widehat{Z}(A)_1$ means that $z \cdot z = - z \cdot z$, hence $2 z^2 = 0$, so $z^2 = 0$, as $\charfield(\bK) \neq 2$, which contradicts the invertibility of $z$.
\end{proof}

\begin{lemma}\label{lem:invoddeltgradsimple}
    Suppose that $A$ is a $C_2$-graded algebra such that $A_0$ is simple and $A^{\times} \cap A_1 \neq \emptyset$. Then $A$ is simple as a $C_2$-graded algebra.
\end{lemma}

\begin{proof}
    If $I$ is a two-sided graded ideal of $A$, then $I_0$ is a two-sided ideal of $A_0$, hence either $I_0 = 0$ or $I_0 = A_0$. If $u \in A^{\times} \cap A_1$, then left-multiplication by $u$ defines a bijection between $I_0$ and $I_1$, hence either $I = 0$ or $I = A$.
\end{proof}

We now describe simple $C_2$-graded $\bK$-algebras in terms of $C_2$-graded division algebras over $\bK$.

\begin{defn}
    For $p,q \geq 0$ with $p+q \geq 1$, and $A$ a $C_2$-graded $\bK$-algebra, $M_{p,q}(A)$ is $M_{p+q}(A)$ as a $\bK$-algebra, with grading
    \[
    M_{p,q}(A)_0 =
\begin{pmatrix}
M_{p\times p}(A_0) & M_{p\times q}(A_1) \\
M_{q\times p}(A_1) & M_{q\times q}(A_0)
\end{pmatrix},
\qquad
M_{p,q}(A)_1 =
\begin{pmatrix}
M_{p\times p}(A_1) & M_{p\times q}(A_0) \\
M_{q\times p}(A_0) & M_{q\times q}(A_1)
\end{pmatrix}.
\]
\end{defn}

Any central simple $C_2$-graded $\bK$-algebra $A$ is isomorphic to $M_{p,q}(D)$ for some $C_2$-graded division $\bK$-algebra $D$ \cite[Thm.\ 2.10.10]{NO}. To show that this $D$ is unique, we introduce the following.

\begin{defn}
    Suppose that $A$ is a $C_2$-graded $\bK$-algebra. A $C_2$-graded $A$-module $M$ is a $C_2$-graded abelian group $M = M_0 \oplus M_1$, which is also an $A$-module, such that $A_i \cdot M_j \subset M_{i+j}$ (where the subscripts are read in $\bZ / 2 \bZ$). A morphism $f \colon M \rightarrow N$ is $C_2$-graded $A$-modules is an $A$-module homomorphism such that $f(M_i) \subset N_i$ for $i = 0,1$.
    
    We write $\Mod_A^{\wedge}$ for the category of $C_2$-graded $A$-modules.
\end{defn}

The category $\Mod_A^{\wedge}$ has an endofunctor $T \colon \Mod_A^{\wedge} \rightarrow \Mod_A^{\wedge}$ which on objects has $T(M) = M$ as an $A$-module but with the grading $T(M)_0 = M_1$ and $T(M)_1 = M_0$, and is the identity on morphisms.

\begin{defn}
    For $C_2$-graded $\bK$-algebras $A$ and $B$, we say that a functor $F \colon \Mod_A^{\wedge} \rightarrow \Mod_B^{\wedge}$ is \emph{$\bK$-linear} if for any $M,N \in \Mod_A^{\wedge}$, the map
    \[
    F(-) \colon \Hom_A^{\wedge}(M,N) \rightarrow \Hom_B^{\wedge}(F(M), F(N)) 
    \]
    is $\bK$-linear, and \emph{$C_2$-graded} if $F \circ T_A = T_B \circ F$. We say that $\Mod_A^{\wedge}$ and $\Mod_B^{\wedge}$ are equivalent as $\bK$-linear $C_2$-graded categories if there is a $\bK$-linear $C_2$-graded equivalence from $\Mod_A^{\wedge}$ to $\Mod_B^{\wedge}$.
\end{defn}

For a $C_2$-graded $\bK$-algebra $A$, there is an internal $\Hom$ in the category $\Mod_A^{\wedge}$, given by setting 
\[
\underline{\Hom}_A^{\wedge}(M,N) \coloneqq \Hom_A(M,N),
\]
with grading 
\[
\underline{\Hom}_A^{\wedge}(M,N)_0 = \Hom_A^{\wedge}(M,N), \qquad \underline{\Hom}_A^{\wedge}(M,N)_1 = \Hom_A^{\wedge}(M,T(N)).
\]
When $M=N$, this makes $\underline{\End}_A^{\wedge}(M)$ a $C_2$-graded $\bK$-algebra. The following is direct to verify.

\begin{lemma}\label{lem:gradedequivbasicproperties}
    Suppose that $A$ and $B$ are $C_2$-graded $\bK$-algebras. Then 
    \begin{itemize}
        \item For any isomorphism $f \colon A \rightarrow B$ is $C_2$-graded $\bK$-algebras, pushforward
        \[
        f_* \colon \Mod_{B}^{\wedge} \rightarrow \Mod_{A}^{\wedge}
        \] 
        is a $\bK$-linear $C_2$-graded equivalence.
        \item Writing $\! {}_A A$ for $A$ viewed as a left $A$-module and with same grading as $A$, the natural map
        \[
A^{\opp} \rightarrow \underline{\End}_A^{\wedge}({}_A A)
\]
is an isomorphism of $C_2$-graded algebras (where $A^{\opp}$ is the opposite algebra of $A$ (not the graded opposite algebra) and carries the same grading as $A$).
    \end{itemize}
\end{lemma}

We are interested in understanding when two categories $\Mod_{A}^{\wedge}$ and $\Mod_{B}^{\wedge}$ are equivalent. One way to produce functors $\Mod_{A}^{\wedge} \rightarrow \Mod_{B}^{\wedge}$ is through the following construction.

Suppose that $X$ is a $C_2$-graded $(B,A)$-bimodule, meaning that $X$ is a $C_2$-graded left $B$-module and a $C_2$-graded right $A$-module, such that the $C_2$-grading is the same and the actions of $A$ and $B$ commute. For any $M \in \Mod_A^{\wedge}$ we can consider the left $B$-module $X \otimes_A M$, and define
\[
(X \otimes_A M)_i = \langle x \otimes m \mid x \in X, m \in M \text{ homogeneous with } \deg(x) + \deg(m) = i\rangle_{\bZ},
\]
which gives $X \otimes_A M$ the structure of a $C_2$-graded $B$-module \cite[\S 2.4]{NO}. The resulting functor
\[
X \otimes_A (-) \colon \Mod_{A}^{\wedge} \rightarrow \Mod_{B}^{\wedge}
\]
is $C_2$-graded, and $\bK$-linear whenever the $\bK$-linear structures on $X$ from $A$ and $B$ agree.

For the pair $A, M_{p,q}(A)$, with $n \coloneqq p+q \geq 1$, we consider the $C_2$-graded $(A, M_{p,q}(A))$-bimodule
\[
A^{1 \times n}_{p,q}, \: \text{ with } \: (A^{1 \times n}_{p,q})_0 = A_0^{1 \times p} \oplus A_1^{1 \times q}, \: \text{ and } \: (A^{1 \times n}_{p,q})_1 = A_1^{1 \times p} \oplus A_0^{1 \times q}
\]
of row vectors, and the $C_2$-graded $(M_{p,q}(A), A)$-bimodule of column vectors
\[
A^{n \times 1}_{p,q}, \: \text{ with } \: (A^{n \times 1}_{p,q})_0 = A_0^{p \times 1} \oplus A_1^{q \times 1}, \: \text{ and } \: (A^{n \times 1}_{p,q})_1 = A_1^{p \times 1} \oplus A_0^{q \times 1}.
\]

\begin{lemma}\label{lem:gradedmorita}
    Let $A$ be a $C_2$-graded $\bK$-algebra and $p,q \geq 0$ with $n \coloneqq p+q \geq 1$. Then the functors 
        \[
        F \colon \Mod_A^{\wedge} \rightarrow \Mod_{M_{p,q}(A)}^{\wedge}, \qquad F(M) \coloneqq A^{1 \times n}_{p,q} \otimes_{M_{p,q}(A)} M,
        \]
        and
        \[
        G \colon \Mod_{M_{p,q}(A)}^{\wedge} \rightarrow \Mod_A^{\wedge}, \qquad G(N) \coloneqq A^{n \times 1}_{p,q} \otimes_{A} N,
        \]
        are mutually inverse equivalences of $\bK$-linear $C_2$-graded categories.
\end{lemma}

\begin{proof}
    From the above discussion, both functors are $C_2$-graded, and $\bK$-linear because the two possible $\bK$-linear structures on each of the bimodules $A^{n \times 1}_{p,q}$ and $A^{1 \times n}_{p,q}$ coincide. To see that they are mutually equivalent, it suffices to show that there are isomorphisms
    \[
        \phi \colon A^{n \times 1}_{p,q} \otimes_A A^{1 \times n}_{p,q} \xrightarrow{\sim} M_{p,q}(A), \qquad \psi \colon A^{1 \times n}_{p,q} \otimes_A A^{n \times 1}_{p,q} \xrightarrow{\sim} A
    \]
    of $C_2$-graded $(M_{p,q}(A), M_{p,q}(A))$-bimodules and $(A,A)$-bimodules respectively. Let $e_1, ... , e_n$ denote the standard $A$-basis of $A^{n \times 1}_{p,q}$ and $f_1, \cdots, f_n$ the standard $A$-basis of $A^{1 \times n}_{p,q}$. We define $\phi$ and $\psi$ by
    \[
        \phi\left(\sum_{i,j} a_{ij} e_i \otimes f_j\right) = \sum_{i,j} a_{ij} E_{ij}, \qquad  \psi\left(\sum_{i,j} a_{ij} f_i \otimes e_j\right) = \sum_{i} a_{ii}.
    \]
    It is direct to check that these are bimodule isomorphisms, which respect the $C_2$-grading.
    \end{proof}

\begin{lemma}\label{lem:divalgcat}
    Suppose that $D$ is a $C_2$-graded division algebra over $\bK$. Then
    \begin{itemize}
        \item If $D_1 \neq 0$, every object of $\Mod_{D}^{\wedge}$ is isomorphic to a direct sum of copies of ${}_D D$.
        \item If $D_1 = 0$, then object of $\Mod_{D}^{\wedge}$ is isomorphic to a direct sum of copies of ${}_D D$ and $T({}_D D)$.
    \end{itemize}
\end{lemma}

\begin{proof}
    Suppose that $D_1 \neq 0$ and $M$ is an object of $\Mod_D^{\wedge}$. $D_0$ is a $\bK$-division algebra, so there is a basis $\{e_i\}_{i \in I}$ of $M_0$ as a $D_0$-module. It is straightforward to show that
    \[
    \phi \colon \bigoplus_{i \in I} {}_D D \rightarrow M, \qquad \phi\left(\sum_i x_i\right) \coloneqq x_i e_i
    \]
    is an isomorphism in $\Mod_D^{\wedge}$, using that multiplication by any $x \in D_1 \setminus 0$ defines a $\bK$-linear isomorphism from $M_0$ to $M_1$. 

    Now suppose that $D_1 = 0$, and $M$ an object of $\Mod_D^{\wedge}$. We may similarly take a basis $\{e_i\}_{i \in I}$ of $M_0$ and $\{f_j\}_{j \in J}$ of $M_1$ as $D_0$-modules, and define
    \[
    \psi \colon \bigoplus_{i \in I} {}_D D \oplus \bigoplus_{j \in J} T({}_D D) \rightarrow M, \qquad \phi\left(\sum_i x_i, \sum_j y_j \right) \coloneqq \sum_i x_i e_i + \sum_j y_j f_j,
    \]
    which is also easily seen to be an isomorphism in $\Mod_D^{\wedge}$.
\end{proof}

\begin{lemma}\label{lem:recovergradeddivisionalgebra}
    Suppose that $D$ and $D'$ are $C_2$-graded division algebras over $\bK$. If $\Mod_{D}^{\wedge}$ and $\Mod_{D'}^{\wedge}$ are equivalent as $\bK$-linear $C_2$-graded categories, then $D_A \cong D_B$ as $C_2$-graded $\bK$-algebras.
\end{lemma}

\begin{proof}
    Let $F \colon \Mod_{D'}^{\wedge} \rightarrow \Mod_{D}^{\wedge}$ be such an equivalence, and set $X \coloneqq {}_{D'} D'$. The object $Y \coloneqq F(X) \in \Mod_{D}^{\wedge}$ is simple, and therefore isomorphic to ${}_D D$ or $T({}_D D)$ by Lemma \ref{lem:divalgcat}. The functor $T \colon \Mod_{D}^{\wedge} \rightarrow \Mod_{D}^{\wedge}$ is a $\bK$-linear $C_2$-graded equivalence, and therefore potentially replacing $F$ with $T \circ F$ we may assume that $Y \cong {}_D D$. Using that $F(T'(X)) = T(F(X)) = T(Y)$, the functor $F$ induces identifications
\[\begin{tikzcd}
	{\underline{\End}_{D'}^{\wedge}(X)} & {\Hom_{D'}^{\wedge}(X,X)} & {\Hom_{D'}^{\wedge}(X,T'(X))} \\
	{\underline{\End}_D^{\wedge}(Y)} & {\Hom_{D}^{\wedge}(Y,Y)} & {\Hom_{D}^{\wedge}(Y,T(Y))}
	\arrow["{=}"{description}, draw=none, from=1-1, to=1-2]
	\arrow["\bigoplus"{description}, draw=none, from=1-2, to=1-3]
	\arrow["{F(-)}"', from=1-2, to=2-2]
	\arrow["{F(-)}"', from=1-3, to=2-3]
	\arrow["{=}"{description}, draw=none, from=2-1, to=2-2]
	\arrow["\bigoplus"{description}, draw=none, from=2-2, to=2-3]
\end{tikzcd}\]
The equivalence $F$ is between additive categories, hence $F$ is additive, and hence the induced identification $\Phi \colon \underline{\End}_{D'}^{\wedge}(X) \rightarrow \underline{\End}_D^{\wedge}(Y)$ commutes with addition. For multiplication,
\begin{align*}
    \Phi(f_0, f_1) \cdot \Phi(g_0, g_1) &= (F(f_0), F(f_1)) \cdot (F(g_0), F(g_1)), \\
    &= (F(f_0) \circ F(g_0) + F(f_1) \circ F(g_1), F(f_0) \circ F(g_1) + F(f_1) \circ F(g_0)), \\
    &= (F(f_0 \circ g_0) + F(f_1 \circ g_1), F(f_0 \circ g_1) + F(f_1 \circ g_0)), \\
    &= \Phi((f_0, f_1) \cdot (g_0,g_1)),
\end{align*}
again using that $F$ is additive. It is also easy to see that $\Phi$ preserves the identity, and therefore the $\bK$-algebra structure because $F$ is $\bK$-linear. Using Lemma \ref{lem:gradedequivbasicproperties}, we therefore have isomorphisms
\[
(D')^{\opp} \cong \underline{\End}_{D'}^{\wedge}({}_{D'}D') = \underline{\End}_{D'}^{\wedge}(X) \xrightarrow{\Phi} \underline{\End}_D^{\wedge}(Y) \cong \underline{\End}_{D}^{\wedge}({}_D D) \cong D^{\opp}
\]
of $C_2$-graded $\bK$-algebras, and thus $D \cong D'$ as $C_2$-graded $\bK$-algebras.
\end{proof}

\begin{cor}\label{cor:stdformforgradedalg}
If $A$ is a simple $C_2$-graded algebra over $\bK$, then 
    \[
    A \cong M_{p,q}(D_A)
    \]
for a unique $C_2$-graded division algebra $D_A$ over $\bK$.
\end{cor}

\begin{proof}
    The existence is \cite[Thm.\ 2.10.10]{NO}. For uniqueness, if additionally $A \cong M_{p',q'}(D')$ for some other $C_2$-graded division algebra $D'$ over $\bK$, then $\Mod_{D}^{\wedge}$ and $\Mod_{D'}^{\wedge}$ are equivalent as $\bK$-linear $C_2$-graded categories by Lemma \ref{lem:gradedequivbasicproperties} and Lemma \ref{lem:gradedmorita}, and hence $D \cong D'$ as $C_2$-graded $\bK$-algebras by Lemma \ref{lem:recovergradeddivisionalgebra}.
\end{proof}

\begin{cor}\label{cor:morita}
Let $A$ and $B$ be simple $C_2$-graded $\bK$-algebras. Then $\Mod_A^{\wedge}$ and $\Mod_B^{\wedge}$ are equivalent as $\bK$-linear $C_2$-graded categories if and only if $D_A \cong D_B$ as $C_2$-graded $\bK$-algebras.
\end{cor}

\begin{proof}
    If $D_A \cong D_B$ as $C_2$-graded $\bK$-algebras, then 
    \[
    \Mod_A^{\wedge} \cong \Mod_{D_A}^{\wedge} \cong \Mod_{D_B}^{\wedge} \cong \Mod_B^{\wedge}
    \]
    $\bK$-linear $C_2$-graded categories by Lemma \ref{lem:gradedequivbasicproperties} and Lemma \ref{lem:gradedmorita}. Conversely, if $\Mod_A^{\wedge}$ and $\Mod_B^{\wedge}$ are equivalent as $\bK$-linear $C_2$-graded categories, then so are $\Mod_{D_A}^{\wedge}$ and $\Mod_{D_B}^{\wedge}$ by Lemma \ref{lem:gradedmorita}, and thus $D_A \cong D_B$ as $C_2$-graded $\bK$-algebras by Lemma \ref{lem:recovergradeddivisionalgebra}.
\end{proof}

\begin{remark}
    We note that in all the above one could replace $\Mod_A^{\wedge}$ with the full subcategory of \emph{finite-dimensional} $C_2$-graded $A$-modules, and all statements and proofs remain valid.
\end{remark}

We can use the above results to prove a $C_2$-graded double centraliser theorem. We recall \cite[\S IV.2]{LAM} that for a $C_2$-graded $\bK$-algebra $A$ and graded subspace $S \subset A$, the graded centraliser is the graded subalgebra $\widehat{C}_A(S)$ with
\begin{align*}
        \widehat{C}_A(S)_0 &= \{a \in A_0 \mid ab = ba \text{ for all } b \in S\}, \\
        \widehat{C}_A(S)_1 &= \{a \in A_1 \mid ab = ba \text{ for all } b \in S_0, \: ab = - ba \text{ for all } b \in S_1\}.
    \end{align*}
For example, $\widehat{Z}(A) = \widehat{C}_A(A)$.
\begin{prop}\label{prop:graddoublecent}
    Suppose that $A$ is a central simple $C_2$-graded algebra over $\bK$, and let $B$ be a simple $C_2$-graded subalgebra of $A$. Then
    \begin{itemize}
        \item $C \coloneqq \widehat{C}_A(B)$ is simple as a $C_2$-graded algebra over $\bK$,
        \item $B = \widehat{C}_A(C)$,
        \item $\dim_{\bK} A = \dim_{\bK} B \cdot \dim_{\bK} C$.
    \end{itemize}

    Furthermore, $\widehat{Z}(B) = \bK$ and and only if $\widehat{Z}(C) = \bK$, in which case the natural map
    \[
    B \, \hat{\otimes} \, C \rightarrow A
    \]
    is an isomorphism of $C_2$-graded algebras over $\bK$.
\end{prop}

\begin{proof}
    Write $A^*$ for the graded opposite algebra of $A$ \cite[\S IV.4]{LAM}, and $A^{\opp}$ for the usual opposite algebra. From the proof of \cite[Prop.\ IV.4.3]{LAM} there is an isomorphism
    \[
    \Theta \colon A \hat{\otimes} A^* \xrightarrow{\sim} \End_{\bK}(A)
    \]
    of $C_2$-graded $\bK$-algebras which on homogeneous elements is given by
    \[
    \Theta(b \otimes a^*)(x) = (-1)^{|a||x|} bxa.
    \]
    We can restrict this left $A \hat{\otimes} A^*$-module structure on $A$ to $A^*$, and it is straightforward to show that
    \[
        f_{(-)} \colon (A^*)^{\opp} \rightarrow \underline{\End}_{A^*}^{\wedge}(A), \qquad f_a(x) \coloneqq (-1)^{|a||x|}a x
    \]
    is an isomorphism of $C_2$-graded $\bK$-algebras. Setting $E \coloneqq B \hat{\otimes} A^*$, a direct computation shows that this restricts to a $C_2$-graded $\bK$-algebra isomorphism
    \[
        f_{(-)} \colon (\widehat{C}_A(B)^*)^{\opp} \xrightarrow{\sim} \underline{\End}_{E}^{\wedge}(A).
    \]

    By \cite[Thm.\ IV.2.3]{LAM}, $E$ is a simple $C_2$-graded algebra over $\bK$, and therefore $E \cong M_{p,q}(D_E)$ by Corollary \ref{cor:stdformforgradedalg}. From Lemma \ref{lem:gradedmorita}, there is an equivalence $F \colon \Mod_{D_E}^{\wedge} \xrightarrow{\sim} \Mod_E^{\wedge}$, and therefore by Lemma \ref{lem:divalgcat} any object of $\Mod_E^{\wedge}$ is isomorphic to a direct sum of copies of $X \coloneqq F({}_{D_E} D_E)$ and $T(X)$ (which themselves are isomorphic if and only if $(D_E)_1 \neq 0$). In particular, $A \cong X^{\oplus m_1} \oplus T(X)^{\oplus m_2}$ for some $m_1, m_2 \geq 0$. We have that $D_E^{\opp} \xrightarrow{\sim} \underline{\End}_E^{\wedge}(X)$ by as $C_2$-graded $\bK$-algebras by Lemma \ref{lem:gradedequivbasicproperties}, and
    \[
    T(-) \colon \underline{\End}_E^{\wedge}(X) \xrightarrow{\sim} \underline{\End}_E^{\wedge}(T(X))
    \]
    is an isomorphism of $C_2$-graded $\bK$-algebras (following the same reasoning as in the proof of Lemma \ref{lem:recovergradeddivisionalgebra}). Furthermore, the natural map
    \[
    M_{m_1, m_2}(\underline{\End}_E^{\wedge}(X)) \rightarrow \underline{\End}_E^{\wedge}(X^{\oplus m_1} \oplus T(X)^{\oplus m_2})
    \]
    is an isomorphism of $C_2$-graded $\bK$-algebras, and therefore
    \[
        M_{m_1, m_2}(D_E^{\opp}) \cong \underline{\End}_E^{\wedge}(X^{\oplus m_1} \oplus T(X)^{\oplus m_2}) \cong \underline{\End}_E^{\wedge}(A) \cong (\widehat{C}_A(B)^*)^{\opp}
    \]
    as $C_2$-graded $\bK$-algebras. In particular, $C = \widehat{C}_A(B)$ is simple as a $C_2$-graded $\bK$-algebra.

    We have that $B \subset \widehat{C}_A(C)$, and therefore to see that $B \subset \widehat{C}_A(C)$ it suffices to show the third claim that $\dim_{\bK} A = \dim_{\bK} B \cdot \dim_{\bK} C$. Let us denote these quantities as $a,b$ and $c$, and similarly for $d, e, x$. 

    By definition, $e = ab$, and setting $m \coloneqq m_1 + m_2$, $a = m x$ and $c = m^2 d$ from the isomorphism $(\widehat{C}_A(B)^*)^{\opp} \cong M_{m_1, m_2}(D_E^{\opp})$. Setting $n \coloneqq p+q$, then from the isomorphism $E \cong M_{p,q}(D_E)$, $e = n^2 d$ and $x = n d$ from the definition of the equivalence of Lemma \ref{lem:gradedmorita}. Together,
    \[
    ab = e = n^2 d = \left(\frac{x}{d} \right)^2 d = \frac{x^2}{d} = \left(\frac{a}{m} \right)^2 \frac{1}{d}= \frac{a^2}{m^2 d} = \frac{a^2}{c},
    \]
    and so $a = bc$ as required.
    
    Finally, note that if $\widehat{Z}(B) = \bK$, then $B \hat{\otimes} C$ is simple by \cite[Thm.\ IV.2.3]{LAM}, and therefore the natural map $B \, \hat{\otimes} \, C \rightarrow A$ has trivial kernel, and thus is an isomorphism by comparing the dimensions. As this induces an isomorphism $\widehat{Z}(B) \hat{\otimes} \widehat{Z}(C) \xrightarrow{\sim} \widehat{Z}(A)$ \cite[Thm.\ IV.2.3]{LAM}, it follows that $\widehat{Z}(C) = \bK$. The same argument with the roles of $B$ and $C$ reverses shows that $\widehat{Z}(C) = \bK$ implies that $\widehat{Z}(B) = \bK$.
\end{proof}

The above results have consequences for the Brauer-Wall group \cite[\S IV.4]{LAM}.

\begin{defn}\label{def:BW}
    The Brauer-Wall group of $\bK$, $\BrW(\bK)$, is the set of central simple $C_2$-graded algebras over $\bK$ up to equivalence, where $[A] = [B]$ in $\BrW(\bK)$ if and only if $M_{p,q}(A) \cong M_{r,s}(B)$ for some $p,q,r,s \geq 0$ with $p+q, r+s \geq 1$. 
    
    The multiplication is defined by the graded tensor product
    \[
    [A] \cdot [B] \coloneqq [A \hat{\otimes} B].
    \]
\end{defn}

This is related to the classical Brauer group as follows. First, viewing any simple $\bF$-algebra $A$ as a simple $C_2$-graded $\bF$-algebra with the trivial grading defines a homomorphism
\[
\Br(\bF) \rightarrow \BrW(\bF).
\]
This forms the first term of the exact sequence \cite[Thm.\ IV.4.4]{LAM}
\begin{equation}\label{eqn:exactseqbrauer}
0 \rightarrow \Br(\bF) \rightarrow \BrW(\bF) \xrightarrow{j} Q(\bF) \rightarrow 0,
\end{equation}
where 
\[
Q(\bF) \coloneqq \bZ / 2 \bZ \: \dot{\times} \: \bF^{\times} / \bF^{\times 2},
\]
the direct product as a set, but with twisted multiplication 
\[
(e,r) * (e', r') = (e+e', (-1)^{ee'} rr').
\]
The group homomorphism $j \colon \BrW(\bF) \rightarrow Q(\bF)$ is defined in \cite[Thm.\ IV.4.3]{LAM}.

We write $\GradDiv_{\bK}^{C_2}$ for the set of \emph{central} $C_2$-graded division $\bK$-algebras up to isomorphism. The following shows that we can think of $\BrW(\bK)$ is defining a group structure on this set.

\begin{cor}\label{cor:divalgrep}
    The map
    \[
    \GradDiv_{\bK}^{C_2} \rightarrow \BrW(\bK), \qquad D \mapsto [D]
    \]
    is an isomorphism: every class in $\BrW(\bK)$ has a unique representative in $\GradDiv_{\bK}^{C_2}$.
\end{cor}

\begin{proof}
    There is a canonical $\bK$-algebra identification $\widehat{Z}(D) \rightarrow \widehat{Z}(M_{p,q}(D))$ defined by
    \[
    z_0 + z_1 \mapsto \begin{pmatrix}
        (z_0 + z_1) I_p & 0_{p,q} \\ 0_{q,p} & (z_0 - z_1) I_q
    \end{pmatrix},
    \]
    and therefore the claim follows directly from Corollary \ref{cor:stdformforgradedalg}.
\end{proof}

\begin{cor}
    If $A$ and $B$ are central simple $C_2$-graded $\bK$-algebras, then $[A] = [B]$ in $\BrW(\bK)$ if and only if $\Mod_A^{\wedge}$ and $\Mod_B^{\wedge}$ are equivalent as $\bK$-linear $C_2$-graded categories.
\end{cor}

\begin{proof}
    This follows directly from Corollary \ref{cor:morita} and Corollary \ref{cor:divalgrep}.
\end{proof}

%% file: appendix-classicalresults.tex
\chapter{Some Results from Linear Representation Theory}\label{chapter:appendix}

\section{Induction and Inflation Between Cyclic Groups}

Let $n \geq 1$. Let $\sR_{\bZ}(C_n)$ denote the representation ring of $\bC$-linear representations of $C_n$ over $\bZ$, and let $\sR_{\bC}(C_n)$ denote the base change to $\bC$. For any $m \mid n$ there is an inflation ring homomorphism
\[
\Inf_m^n \colon \sR_{\bZ}(C_m) \rightarrow \sR_{\bZ}(C_n),
\]
and a $\bZ$-linear induction map
\[
\Ind_m^n \colon \sR_{\bZ}(C_m) \rightarrow \sR_{\bZ}(C_n).
\]

Let $Q_n \subset \sR_{\bZ}(C_n)$ denote the image of $\Inf_m^n$ over all $m \mid n$ with $m < n$, which consists of the $\bZ$-span of the non-primitive characters of $C_n$. Let $S_n \subset \sR_{\bZ}(C_n)$ denote the image of $\Inf_m^n$ and $\Ind_m^n$ over all $m \mid n$ with $m < n$.

The aim of this section is to prove the following.

\begin{thm}\label{thm:wheniscyclicnotneeded}
    Suppose that $n \geq 1$. Then $S_n = \sR_{\bZ}(C_n)$ if and only if $n \equiv 2 \! \! \mod 4$.
\end{thm}

The equality when $n \equiv 2 \! \! \mod 4$ is relatively straightforward, but for the converse we will need the following. To state it, for $n \geq 1$ let
\[
\Delta(n) \coloneqq \rank_{\bZ} \sR_{\bZ}(C_n) - \rank_{\bZ} S_n \geq 0.
\]
Let $\mu \colon \bN_{\geq 1} \rightarrow \{-1,0,1\}$ denote the M\"{o}bius function, and $\varphi \colon \bN_{\geq 1} \rightarrow \bN_{\geq 1}$ Euler's totient function.

\begin{prop}\label{prop:descofrankdifference}
    For any $n \geq 1$,
    \[
        \Delta(n) = \sum_{d \mid n} \mu(d) \, \varphi \! \left( \frac{n}{d} \right).
    \]
    In particular, $\Delta \colon \bN_{\geq 1} \rightarrow \bN_{\geq 0}$ is multiplicative, and $\Delta(n) = 0$ if and only if $n \equiv 2 \! \! \mod 4$.
\end{prop}

To prove Proposition \ref{prop:descofrankdifference}, we first show that $\Delta$ is multiplicative.

\begin{lemma}\label{lem:deltaismultiplicitive}
    Suppose that $m,n \geq 1$ with $(m,n) = 1$. Then $\Delta(mn) = \Delta(m) \Delta(n)$.
\end{lemma}

\begin{proof}
    Because $(m,n) = 1$, we have a natural isomorphism
    \[
    \Phi_{m,n} \colon \sR_{\bC}(C_m) \otimes_{\bC} \sR_{\bC}(C_n) \xrightarrow{\sim} \sR_{\bC}(C_m \times C_n).
    \]
    We have that $\Delta(mn)$ is the dimension of the quotient of $\sR_{\bC}(C_m \times C_n)$ by $S_{m,n}$, where $S_{m,n}$ is the sum of the images of the induction maps from proper subgroups, and the inflation maps from proper quotients. By the transitivity of induction and inflation, we may restrict to maximal subgroups, and quotients by minimal subgroups. These are all of the form $C_{m/p} \times C_n$ or $C_m \times C_{n/p}$, where $p$ is a prime divisor of $m$ or $n$. For a prime divisor $p$ of $mn$, if $p \mid m$ there are commutative diagrams
\[\begin{tikzcd}
	{\sR_{\bC}(C_m) \otimes_{\bC} \sR_{\bC}(C_n)} & {\sR_{\bC}(C_m \times C_n)} & {\sR_{\bC}(C_m) \otimes_{\bC} \sR_{\bC}(C_n)} & {\sR_{\bC}(C_m \times C_n)} \\
	{\sR_{\bC}(C_{m/p}) \otimes_{\bC} \sR_{\bC}(C_n)} & {\sR_{\bC}(C_{m/p} \times C_n)} & {\sR_{\bC}(C_{m/p}) \otimes_{\bC} \sR_{\bC}(C_n)} & {\sR_{\bC}(C_{m/p} \times C_n)}
	\arrow["{\Phi_{m,n}}"', from=1-1, to=1-2]
	\arrow["\sim"', no head, from=1-2, to=1-1]
	\arrow["{\Phi_{m,n}}"', from=1-3, to=1-4]
	\arrow["\sim"', no head, from=1-4, to=1-3]
	\arrow["{\tilde{\Ind}_{m/p}^m \otimes 1}", from=2-1, to=1-1]
	\arrow["{\Phi_{m/p, n}}"', from=2-1, to=2-2]
	\arrow["{\tilde{\Ind}_{m/p, n}^{m,n}}"', from=2-2, to=1-2]
	\arrow["\sim"', no head, from=2-2, to=2-1]
	\arrow["{\tilde{\Inf}_{m/p}^m \otimes 1}", from=2-3, to=1-3]
	\arrow["{\Phi_{m/p,n}}"', from=2-3, to=2-4]
	\arrow["{\tilde{\Inf}_{m/p, n}^{m,n}}"', from=2-4, to=1-4]
	\arrow["\sim"', no head, from=2-4, to=2-3]
\end{tikzcd}\]
and similarly if instead $p \mid n$. In particular, we have that
\[
\Phi_{m,n}(S_m \otimes_{\bC} \sR_{\bC}(C_n) + \sR_{\bC}(C_m) \otimes_{\bC} S_n) = S_{m,n},
\]
and therefore that $\Phi_{m,n}$ induces an isomorphism
\[
    \tilde{\Phi}_{m,n} \colon (\sR_{\bC}(C_m) / S_m) \otimes_{\bC} (\sR_{\bC}(C_n) / S_n) \xrightarrow{\sim} \sR_{\bC}(C_m \times C_n)/ S_{n,m}.
\]
Taking dimensions we see that $\Delta(mn) = \Delta(m) \Delta(n)$ as required.
\end{proof}

\begin{proof}[Proof of Proposition \ref{prop:descofrankdifference}]
    We want to show that $\Delta = \mu * \varphi$, the Dirichlet convolution of $\mu$ and $\varphi$. Because both $\mu$ and $\varphi$ are multiplicative, $\mu * \varphi$ is too, and therefore it suffices to show that $\Delta = \mu * \varphi$ on prime powers. For $p$ a prime and $r \geq 1$, 
    \[
        (\mu * \varphi)(p^r) = \sum_{k = 0}^r \mu(p^k) \varphi(p^{r-k}) = \varphi(p^r) - \varphi(p^{r-1}).
    \]
    On the other hand, $\Delta(p^r)$ is the dimension of the cokernel of
    \[
        \Ind_{p^{r-1}}^{p^r} \oplus \Inf_{p^{r-1}}^{p^r} \colon \sR_{\bC}(C_{p^{r-1}}) \oplus \sR_{\bC}(C_{p^{r-1}}) \rightarrow \sR_{\bC}(C_{p^r}).
    \]
    The image of the inflation map has dimension $p^r - \varphi(p^r)$, consisting of the imprimitive characters, and the image of the induction map has dimension $p^{r-1}$, being injective. The intersection of the images is the image of the non-primitive characters of $C_{p^{r-1}}$, which has dimension $p^{r-1} - \varphi(p^{r-1})$. We may therefore compute that
    \begin{align*}
        \Delta(p^r) = p^r - (p^r - \varphi(p^r)) - p^{r-1} + (p^{r-1} - \varphi(p^{r-1})) = \varphi(p^r) - \varphi(p^{r-1}),
    \end{align*}
    and so $\Delta = \mu * \varphi$, as required. Using this, we have 
    \[
        \Delta(p^r) = \varphi(p^r) - \varphi(p^{r-1}) = \twopartdef{(p-1) - 1}{r = 1}{p^{r-2}(p-1)^2}{r \geq 2},
    \]
    which equals $0$ if and only if $p = 2$ and $r = 1$. In particular, if $n = 2^e p_1^{e_1} \cdots p_k^{e_k}$, then $\Delta(n) = 0$ if and only if $e = 1$, which is the case exactly when $n \equiv 2 \! \! \mod 4$.
\end{proof}

\begin{remark}
    There is a more conceptual proof of Proposition \ref{prop:descofrankdifference} which explains why the M\"{o}bius function appears, at the cost of being much longer than what we have given above.
\end{remark}

\begin{proof}[Proof of Theorem \ref{thm:wheniscyclicnotneeded}]
    When $S_n = \sR_{\bZ}(C_n)$, $\Delta(n) = 0$, and so $n \equiv 2 \! \! \mod 4$ by Proposition \ref{prop:descofrankdifference}. Conversely, suppose that $n \equiv 2 \! \! \mod 4$, so $n = 2m$ for some odd $m \geq 1$. Let $x$ be a generator of $C_n$, let $\chi$ be an irreducible complex character of $C_n$, and set $\zeta \coloneqq \chi(x)$, and $n$th root of $1$. If $\zeta$ is not primitive, then $\chi \in Q_n \subset S_n$. If instead $\zeta$ is primitive, then we can write
    \[
        \chi = \Ind_{m}^n (\psi_{2}) - \Inf_{m}^{n}(\psi_{1+m}),
    \]
    where for $\langle y \rangle = C_m \leq C_n$, $y = x^2$, $\psi_{2}(y) = \zeta^2$ and $\psi_{1+m} = \zeta^{1+m}$, both of which are well-defined characters of $C_m$ because $2$ and $1+m$ are even. Then as $\chi \in S_n$ for any irreducible $\chi$, $S_n = \sR_{\bZ}(C_n)$.
\end{proof}

%% file: appendix-wallpaper.tex
\chapter{Real Structures on Wallpaper Groups} \label{appendix_wallpaper}

In this appendix, we rederive the classification of magnetic wallpaper groups, using our result on the classification of Real structures from Chapter \ref{chapter:realstructures}. We also give a second proof of the classification, by explicitly counting index two subgroups of wallpaper groups. We use this explicit presentation when we compute the antilinear representation theory of these groups, in Appendix \ref{appendix_repswallpaper}.

\section{Wallpaper Groups}

\subsection{Presentation} Let us first recall the classification of wallpaper groups, that is, $2$-dimensional crystallographic groups, as considered in Chapter \ref{chap:physics}.

We give generators for each wallpaper group inside $\Euc_2 = \HO_2 (\bR)\ltimes \bR^2$. The following matrices are useful. Their index denotes their order.
\[
I =\begin{pmatrix}1&0\\0&1\end{pmatrix},\quad
U_2=\begin{pmatrix}0&1\\1&0\end{pmatrix},\quad
W_2=\begin{pmatrix}-1&0\\0&-1\end{pmatrix},\quad
X_2=\begin{pmatrix}1&0\\0&-1\end{pmatrix},\quad
X_3=\begin{pmatrix}-1&-1\\1&0\end{pmatrix},
\]

\[
X_4=\begin{pmatrix}0&-1\\1&0\end{pmatrix}, \quad
X_6=\begin{pmatrix}0&-1\\1&1\end{pmatrix},\quad
Y_2=\begin{pmatrix}-1&0\\0&1\end{pmatrix}, \quad
Z_2=\begin{pmatrix}0&-1\\-1&0\end{pmatrix}, \quad
X_{\infty}=\begin{pmatrix}1&1\\0&1\end{pmatrix}.
\]

Recall that $[A\mid x]$ denotes the affine transformation $y\mapsto Ay+x$, for $A \in \GL_2(\bR)$ and $x,y \in \bR^2$, column vectors. Any element $v \in \bR^2$ defines a translation $[I_2 \mid v]$. In particular, we will implicitly $\bR^2$ as a subset of $\Euc_2$. Let $e_1$ and $e_2$ denote the standard basis of $\bR^2$. In Table \ref{table:wallpaper1}, we list the 17 isomorphism classes of wallpaper groups, each such $G$ given by generators, so that the lattice it defines $T \coloneqq G \cap \bR^2$ is the subgroup $\langle e_1, e_2 \rangle$. 

\begin{table}[!htbp]
  \centering
\small
\begin{tabular}{lllllll}
\toprule
\textbf{Group} & \textbf{Orbifold} & \textbf{Cell} & \textbf{Centre} & $\boldsymbol{Q}$ & $\boldsymbol{G}=\langle e_1,e_2,$ & $\boldsymbol{Z(G)}$ \\
\midrule
p1   & $o$       & oblique     &  & $C_1=\langle \emptyset \rangle$ &  $\emptyset \rangle$ & $T=\langle e_1,e_2\rangle$ \\
p2   & $2222$    & oblique     &  & $C_2=\langle W_2\rangle$ &  $W_2\rangle$ & $0$ \\
\midrule
pm   & $**$      & rectangular &  & $C_2=\langle X_2\rangle$ & $X_2\rangle$ & $\bZ =\langle e_1\rangle$ \\
pg   &  $xx$      & rectangular &  & $C_2=\langle X_2\rangle$ & $[X_2\mid \tfrac12 e_1]\rangle$ & $\bZ=\langle e_1\rangle$ \\
cm   & $*x$      & rhombic & $\tfrac12 (e_1+e_2) $ & $C_2=\langle U_2\rangle$ & $U_2\rangle$ &  $\bZ =\langle e_1+e_2\rangle$ \\
pmm  & $*2222$   & rectangular &  & $K_4=\langle X_2,Y_2\rangle$ & $X_2,Y_2\rangle$ & $0$\\
pmg  &  $22*$     & rectangular &  & $K_4=\langle X_2,Y_2\rangle$ &$X_2, [Y_2\mid\tfrac12 e_2] \rangle$ &  $0$\\
pgg  & $22x$     & rectangular &  & $K_4=\langle X_2,Y_2\rangle$ & $[X_2\mid\tfrac12 e_1] , [Y_2\mid\tfrac12 e_2] \rangle$ & $0$ \\
cmm  &  $2\hspace{-0.25em}*\hspace{-0.25em}22$    & rhombic & $\tfrac12 (e_1+e_2) $ & $K_4=\langle U_2,W_2\rangle$ &$U_2,W_2\rangle$ & $0$ \\
\midrule
p4   & $442$     & square      &  & $C_4=\langle X_4\rangle$ &  $X_4\rangle$ &   $0$  \\
p4m  & $*442$    &  square      & & $D_4=\langle X_4,X_2\rangle$ &  $X_4,X_2\rangle$ &  $0$   \\
p4g  & $4\hspace{-0.25em}*\hspace{-0.25em}2$     &  square      & & $D_4=\langle X_4,X_2\rangle$ &  $X_4,[X_2\mid\tfrac12 (e_1+e_2)]\rangle$ &  $0$   \\
\midrule
p3   & $333$     &  hexagonal   &  & $C_3=\langle X_3\rangle$ &  $X_3\rangle$ &  $0$  \\
p3m1 & $*333$    & hexagonal   & & $D_3=\langle X_3,U_2\rangle$ &  $X_3,U_2\rangle$ &   $0$  \\
p31m & $3\hspace{-0.25em}*\hspace{-0.25em}3$     & hexagonal   & & $D_3=\langle X_3,Z_2\rangle$ &  $X_3,Z_2 \rangle$ &   $0$   \\
p6   &  $632$     & hexagonal   &  & $C_6=\langle X_6\rangle$ &$X_6\rangle$ &  $0$  \\
p6m  & $*632$    & hexagonal   &  & $D_6=\langle X_6,U_2\rangle$ &  $X_6,U_2\rangle$ &  $0$   \\
\bottomrule
\end{tabular}
\caption{The 17 wallpaper groups with generators, orbifold notation and centres.}
\label{table:wallpaper1}
\end{table}

Let us explain the rest of the entries of the table.

\begin{itemize}
  \item  The cell is the parallelogram spanned by $e_1,e_2$, for which we also record the centre of the cell. The type of the cell corresponds to restriction on the basis.
Namely, oblique means no restriction,
rectangular  -- $e_1,e_2$ are orthogonal,
rhombic  -- $e_1,e_2$ are of the same length,
square --  $e_1,e_2$ is an orthonormal basis,
hexagonal -- $e_1,e_2$ is a pair of vectors of the same length
with the interior angle having 60 degrees. 
\item $Z(G)$ is the centre of $G$, and $Q$ is the point group $G / T$ of $G$.
\end{itemize}

Note that there are four conventional notations for the wallpaper groups. We are not using the spectroscopy (Schoenflies) and algebra (Coxeter) notations. Our main notation is the crystallography (Hermann-Mauguin) notation. It is essentially a recipe for both constructing the group and recognising it from the pattern:
\begin{itemize}
\item The first letter stands for {\bfseries p}rimitive or {\bfseries c}entred cell. 
  A centred cell contains a centre listed in column 4; a primitive cell does not.
\item $n\in \{1,2, 3, 4,6\}$ stand for a rotation symmetry by $\frac{360}{n}$ degrees.
\item {\bfseries m} stands for mirror (a reflection is required to generate $G$).
\item {\bfseries g} stands for glide (a glide reflection is required to generate $G$).
\item Both {\bfseries p3m1} and {\bfseries p31m} are semidirect products of $D_3=S_3$ with $T$ but the representations
  are not equivalent. Interestingly enough, if one embeds $S_3$ into $\GL_2 (\bZ)$, one has two representations on $\bZ^2$, one
  on rows, another on columns. The corresponding semidirect products will be {\bfseries p3m1} and {\bfseries p31m}.
  For another algebraic difference, observe
  $H^1(Q,T)$ in Table~\ref{table:wallpaper2}.
\item There are also two geometric differences between {\bfseries p3m1} and {\bfseries p31m}
that help to distinguish them visually.
One difference is that  in {\bfseries p3m1} the mirrors lie between the vectors of $T$, while
in {\bfseries p31m} the mirrors run parallel to the vectors in $T$.
The second difference is that
 in {\bfseries p3m1} all the centres of 3-fold reflections lie on the mirrors, while
 in {\bfseries p31m} there are two types of 3-fold rotation centres: lying on the mirrors and not lying on the mirrors.
\end{itemize} 

We also use the orbifold (Conway) notation. The notation is simply a list of all singular points on the orbifold $\bR^2/G$.
\begin{itemize}
\item $o$ stands for the empty set $\emptyset$, i.e., there are no singular points.
\item 2, 3, 4 and 6 stand for a cone point of the corresponding order.
\item $*$ stands for a boundary component: a connected set of mirror points.
\item 2, 3, 4 and 6 after $*$ means that the cone point lies on the boundary.  
\item $x$ stands for a boundary component: a connected set of glide points.
\end{itemize}   

\subsection{Outer Automorphism Groups} Let us now consider the outer automorphism groups of the 17 wallpaper groups, these being needed to compute their Real structures via the methods of Chapter \ref{chapter:realstructures}. We record these, and the associated data we need to compute them in Table \ref{table:wallpaper2}. We give explicit generators for these outer automorphism groups in Table \ref{table:wallpaper2_geb}.

\begin{table}[!htbp]
\centering
\small
\begin{tabular}{lllllllll}
\toprule
 $\boldsymbol{G}$ &  $\boldsymbol{Q}$ &  $\boldsymbol{H^1(Q,T)}$  & $\boldsymbol{\Out(Q)}$ & $\boldsymbol{C_{\GL_2 (\bZ) }(Q)}$ & $\boldsymbol{M}$  & $\boldsymbol{N(Q)}$ & $\boldsymbol{N(Q,o)}$ &  $\boldsymbol{\Out(G)}$  \\
\midrule
p1   & $C_1$       & $0$      &  $C_1$ & $\GL_2 (\bZ)$  &   & $\GL_2 (\bZ)$ &  $\GL_2 (\bZ)$ &  $\GL_2 (\bZ)$ \\
p2   & $C_2$    & $(\bZ/2)^2$  & $C_1$  &  $\GL_2 (\bZ)$  &   &  $\PGL_2 (\bZ)$ & $\PGL_2 (\bZ)$ & $\PGL_2 (\bZ)\ltimes K_4$  \\
\midrule
pm   & $C_2$      &  $\bZ/2$ & $C_1$  & $K^a_4$ &   & $C_2$ & $C_2$ &  $K_4$ \\
pg   &  $C_2$      & $\bZ/2$  & $C_1$  & $K^a_4$ &   & $C_2$ & $C_2$ &  $K_4$\\
cm   & $C_2$      & $0$  & $C_1$  & $K^b_4$ &   & $C_2$ & $C_2$ &  $C_2$ \\
pmm  & $K_4$   & $(\bZ/2)^2$  & $\boldsymbol{D_3}$  & $K_4^a$ & $X_4$ & $C_2$ &$C_2$ &  $D_4$\\
pmg  &  $K_4$     & $(\bZ/2)^2$  & $\boldsymbol{D_3}$  & $K_4^a$ & $X_4$ &   $C_2$ & $C_1$ &  $K_4$\\
pgg  & $K_4$     & $(\bZ/2)^2$  & $\boldsymbol{D_3}$  & $K_4^a$ & $X_4$ &    $C_2$ & $C_2$ &  $D_4$\\
cmm  &  $K_4$    & $\bZ/2$ & $\boldsymbol{D_3}$  &  $K^b_4$ & $X_4$ &    $C_2$ & $C_2$ &  $K_4$\\
\midrule
p4   & $C_4$     & $\bZ/2$       & $C_2$ & $C^a_4$  & $X_2$ &  $C_2$  &  $C_2$ &  $K_4$\\
p4m  & $D_4$    &  $\bZ/2$       & $\boldsymbol{C_2}$ & $C^a_2$ &   &  $C_1$ &  $C_1$ &  $C_2$   \\
p4g  & $D_4$     & $\bZ/2$      & $\boldsymbol{C_2}$ & $C^a_2$  &   &  $C_1$ &  $C_1$ &  $C_2$   \\
\midrule
p3   & $C_3$     & $\bZ/3$    & $C_2$ &  $C^a_6$ & $U_2$ &  $K_4$ & $K_4$ &  $D_6$\\
p3m1 & $D_3$    & $\bZ/3$   & $C_1$ & $C^a_2$ &   & $C_2$ &  $C_2$ &  $D_3$   \\
p31m & $D_3$     & $0$    & $C_1$ & $C^a_2$  &   & $C_2$ &  $C_2$ &  $C_2$  \\
p6   &  $C_6$     & $0$       & $C_2$ &  $C^a_6$ & $U_2$       & $C_2$ &  $C_2$ &  $C_2$  \\
p6m  & $D_6$    & $0$   & $\boldsymbol{C_2}$ & $C^a_2$ &   &  $C_1$ &  $C_1$ &  $C_1$   \\
\bottomrule
\end{tabular}
\caption{The outer automorphisms of the 17 wallpaper groups.}
\label{table:wallpaper2}
\end{table}

\begin{table}[!htbp]
\centering
\small
\begin{tabular}{ll|l|l}
\toprule
 $\boldsymbol{G}$ & $\boldsymbol{\Out(G)}$  & $A, \ B, \ C$ &  Relations of $\boldsymbol{\Out (G)}$ \\
\midrule
p1   & $\GL_2 (\bZ)$  & $X_2, \, X_3, \, X_4$ &  $A^2$, $B^3$, $C^4$, $(AC)^2$, $(CAB)^2$, $[C^2,B]$ \\
p2   & $\PGL_2 (\bZ)\ltimes K_4$   & $[I \mid\tfrac12 e_1], \ U_2, \ X_{\infty}$ & $A^2$, $B^2$, $[A,BAB]$,  $[A,C]$, $CBABC^{-1}(AB)^{-2}$, \\
   &   &  & $(C^{-2}BCB)^2$, $(BCBC^{-1})^3$ \\
\midrule
pm   & $K_4$  &  $W_2, \ [I \mid\tfrac12 e_2]$ & $A^2$, $B^2$, $[A,B]$ \\
pg   & $K_4$ &  $Y_2, \ [I \mid\tfrac12 e_2]$ & $A^2$, $B^2$, $[A,B]$ \\
cm   & $C_2$  & $W_2$ & $A^2$ \\
pmm  & $D_4$ &  $U_2, \ [I \mid\tfrac12 e_2]$ & $A^2$, $B^2$, $(AB)^4$ \\
pmg  & $D_4$ &  $[I \mid\tfrac12 e_1], \ [I \mid\tfrac12 e_2]$ & $A^2$, $B^2$, $[A,B]$ \\
pgg  & $D_4$ &  $U_2, \ [I \mid\tfrac12 e_2]$ & $A^2$, $B^2$, $(AB)^4$ \\
cmm  & $K_4$ &  $X_4, \ [I \mid\tfrac12 (e_1+ e_2)]$ &  $A^2$, $B^2$, $[A,B]$ \\
\midrule
p4   &   $K_4$ & $X_2, \ [I \mid\tfrac12 (e_1+ e_2)]$ &  $A^2$, $B^2$, $[A,B]$ \\
p4m  &  $C_2$    & $[I \mid\tfrac12 (e_1+ e_2)]$ & $A^2$ \\
p4g  &  $C_2$    & $[I \mid\tfrac12 (e_1+ e_2)]$ & $A^2$ \\
\midrule
p3   &  $D_6$   & $W_2, \ U_2, \ [I \mid\tfrac13 (e_1+ e_2)]$ &  $A^2$, $B^2$, $C^3$, $[A,B]$, $[B,C]$, $(AC)^2$ \\
p3m1 &  $D_3$    &  $X_6, \ [I \mid\tfrac13 (e_1+ e_2)]$ &  $A^2$, $B^3$, $(AB)^2$ \\
p31m &  $C_2$   &  $X_6$ &  $A^2$ \\
p6   &  $C_2$   & $U_2$ & $A^2$ \\
p6m  &  $C_1$    &  &  none \\
\bottomrule
\end{tabular}
\caption{Generators of $\Aut (G)$ and relations of $\Out (G)$ for the 17 wallpaper groups.}
\label{table:wallpaper2_geb}
\end{table}

As each point group $Q$ is finite, we may immediately write down $\Aut(Q)$ and $\Out (Q)$ for the points groups $Q$ for all the wallpaper groups.
To compute the normaliser $N_{\GL_2 (\bZ) }(Q)$, we list the centraliser $C_{\GL_2 (\bZ) }(Q)$ and an additional matrix $M$ so that
\[ N_{\GL_2 (\bZ) }(Q) = \langle Q, C_{\GL_2 (\bZ) }(Q), M \rangle \, .\]
Note that  $N_{\GL_2 (\bZ) }(Q) / C_{\GL_2 (\bZ) }(Q)$ is a subgroup of $\Out (Q)$:
in Table~\ref{table:wallpaper2} we make $\Out (Q)$ bold in those cases when it is a proper subgroup.
There we use the following subgroups of $\GL_2 (\bZ)$:
\[
K^a_4 = \langle X_2,Y_2\rangle = \begin{pmatrix} \ast & 0 \\ 0 & \ast \\ \end{pmatrix}, \
K^b_4 = \langle U_2,W_2\rangle = \begin{pmatrix} x & y \\ y & x \\ \end{pmatrix}, \
C^a_2=\langle W_2\rangle = \begin{pmatrix} x & 0 \\ 0 & x \\ \end{pmatrix},
\]
\[
C^a_4=\langle X_4\rangle = \begin{pmatrix} x & y \\ -y & x \\ \end{pmatrix}, \
C^a_6 = \langle X_6\rangle = \begin{pmatrix} x & y \\ -y & x+y \\ \end{pmatrix}.
\]
The short exact sequence~\eqref{Physics_OutSequence} always splits in the dimension 2,
 so that $\Out (G)$ is computed for Table~\ref{table:wallpaper2} as
\begin{equation} \label{table:OutG}
\Out (G) = N(Q,o) \ltimes H^1(Q,T) \, .
\end{equation}
We also tabulate the companion group
$N (Q) \coloneqq N_{\GL_2 (\bZ) }(Q)/Q$,
which turns out to be equal to $N(Q,o)$ in all cases except {\bfseries pmg}.

An inquisitive reader may have further questions about~\eqref{table:OutG}.
We know it is split automatically, only if the sequence $1 \to T \to G \to Q \to 1$ is split. 
This is true in the thirteen cases except
{\bfseries pg},
{\bfseries pmg},
{\bfseries pgg}
and
{\bfseries p4g}.
In two of them, 
{\bfseries pg}
and
{\bfseries p4g},
the sequence~\eqref{table:OutG} collapses because one of the groups is trivial
but the other two still present a challenge.
Moreover, in five of the groups, namely,
{\bfseries pmm},
{\bfseries pmg},
{\bfseries pgg},
{\bfseries p3}
and
{\bfseries p31m},
there are two possible actions of
$N_{\GL_2 (\bZ) }(Q)/Q)$ on  $H^1(Q,T)$, trivial and non-trivial.
In all five case, the action is non-trivial.

One method to determine the action of $N_{\GL_2 (\bZ) }(Q)/Q$ on $H^1(Q,T)$ and the structure of $\Out (G)$ is to use~\eqref{formula:AUT1},
which is essentially the second Bieberbach Theorem. We carry this out in Table~\ref{table:wallpaper2_geb}. In the third column we list
$A,B,C \in \Aff_2$ so that
\[ \Aut (G) = \langle \gamma_A,\gamma_B,\gamma_C,\gamma_g \,\mid\, g\in G \rangle \quad \mbox{ where } \
\gamma_x (y) = xyx^{-1} \, .\]
If we need less then three generators, we use only $A$ or $A,B$. The group $G$ is always generated as in Table~\ref{table:wallpaper1}.
The last column are relations in $\Out (G)$. For instance,
the first row tells us that in the type {\bfseries p1} we have isomorphisms
\[
\Out(G) \cong \GL_2 (\bZ) \cong \langle A,B,C \, \mid \, A^2 = B^3 = C^4 = (AC)^2= (CAB)^2 = C^2BC^{-2}B^{-1} = 1 \rangle \,
\]
and the quotient map $\Aut (G) \to \Out (G)$ is given by
\[ \gamma_{X_2} \mapsto A, \
\gamma_{X_3} \mapsto B, \
\gamma_{X_4} \mapsto C. \] 

\section{Index Two Subgroups of Wallpaper Groups}

In Table~\ref{table:wallpaper3} we list the index two subgroups of each wallpaper group, grouped by $\Aut (G)$-orbit.
For instance, 3p1 in the first row means 3 subgroups of type p1, which form a single $\Aut(G)$-orbit.
Of course, we simply write p1 instead of 1p1. This is a complete list of the black-and-white magnetic groups.

A subgroup $H$ of index 2 is necessarily necessarily normal and $H\cap T$ has index at most 2 in $T$.
By {\em the translation index} or simply t-index, we understand this index $|T:T\cap H|$, and we group subgroups by this index into type 1 and type 2 subgroups. Note that any subgroup $H$ of index 2 contains $2T$. Therefore, the subgroups of index 2 in $G$ and $G/2T$ are in one-to-one correspondence. When identifying the group $G / 2T$, we use $[n,k]$ to denote SmallGroup(n,k).
Notice that this group is finite and fits into an exact sequence $1\to C_2^n \to G/2T \to Q \to 1$.

An easy method to work out the subgroups of index 2 and $G/2T$ is to use our explicit presentation of $G$ found in the magma code at \cite{GITHUB}. This also makes the presentation of $G/2T$ immediate. The index 2 subgroups are kernels of homomorphism to $C_2$, which one lists from the presentations.

For example, in type {\bfseries p4m} the presentation is
\[
\langle s,t,x,y \mid
    s^2,
    t^4,
    t^s = t^3,
    x^s = x,
    y^s = y^{-1},
    x^t = y^{-1},
    y^t = x,
    y^x = y
\rangle
\]
where $x,y$ are the basis of $T$, and for example $x^t = t^{-1} x t$. Homomorphisms $\phi: G \to C_2$ correspond to choices of $a,b,c\in C_2$ so that
\[ \phi (r) = a, \ \phi (s) = b, \ \phi (x) = c, \ \phi (y) = c. \] 
This gives 7 non-trivial homomorphisms and thus 7 subgroups of index 2.

\begin{table}[!htbp]
\centering
\small
\begin{tabular}{llll|llll}
\toprule
 $\boldsymbol{G}$ &  $\boldsymbol{G/2T}$ &  type 1  & type 2 & $\boldsymbol{\sN}$ & $\boldsymbol{\Theta(G)}$ & $\boldsymbol{Z(G)_{\theta}}$ & $\boldsymbol{\fR_{\theta}(G)}$ \\
\midrule
p1   & $K_4 = [4,2]$       &      &  $3\mbox{p1}$ & 6  & $I$  & $\bZ^2_{+}$ & 1+3 \\
&       &      &   &   & $U_2$  & $\bZ_{d}$  & 1 \\
&       &      &   &   & $X_2$  & $\bZ_{+}\oplus \bZ_{-}$  &  1+1 \\
&       &      &   &   & $W_2$  & $\bZ^2_{-}$  & 1 \\
\midrule
p2   & $C_2^3 = [8,5]$    & $\mbox{p1}$  & $6\mbox{p2}$  &  7  &  7 & &  1 \\
\midrule
pm   & $C_2^3 = [8,5]$     &  $\mbox{p1}$ & $\mbox{pm},\mbox{pg},2\mbox{pm},2\mbox{cm}$ & 6  & $I$ & $\bZ_{+}$  & 1+1 \\
&       &      &   &   & $W_2$  & $\bZ_{-}$  & 1 \\
&       &      &   &   & $[I \mid\tfrac12 e_2]$  & $\bZ_{+}$  & 1+1 \\
&       &      &   &   & $[W_2 \mid\tfrac12 e_2]$  & $\bZ_{-}$  & 1 \\
\midrule
pg   & $C_2\times C_4 = [8,2]$     & $\mbox{p1}$  & $2\mbox{pg}$ & 6 &  $I$ & $\bZ_{+}$  & 1+1 \\
&       &      &   &   & $Y_2$  & $\bZ_{-}$  & 1 \\
&       &      &   &   & $[I \mid\tfrac12 e_2]$  & $\bZ_{+}$  & 1+1 \\
&       &      &   &   & $[Y_2 \mid\tfrac12 e_2]$  & $\bZ_{-}$  & 1 \\
\midrule
cm   & $D_8 = [8,3]$    & $\mbox{p1}$  & $\mbox{pm}, \mbox{pg}$  & 3  & $I$ & $\bZ_{+}$  & 1+1 \\
&       &      &   &   & $W_2$  & $\bZ_{-}$  & 1 \\
\midrule
pmm  & $C_2^4 =[16,14]$   & $\mbox{p2},2\mbox{pm}$  & $4\mbox{pmm},4\mbox{pmg},4\mbox{cmm}$      & 4  & 4  & & 1 \\
pmg  & $C_2\times D_4 =[16,11]$   & $\mbox{p2},\mbox{pm},\mbox{pg}$ & $2\mbox{pmg},2\mbox{pgg}$  & 4  & 4  & & 1  \\
pgg  & $C_4\ltimes K_4 =[16,3]$   & $\mbox{p2},2\mbox{pg}$          &                            & 4  & 4  & & 1  \\
cmm  &  $C_2\times D_4 =[16,11]$    & $\mbox{p2},2\mbox{cm}$ & $2\mbox{pmg},\mbox{pmm},\mbox{pgg}$ & 4 & 4 & & 1 \\
\midrule
p4   & $C_4\ltimes K_4 =[16,3]$      & $\mbox{p2}$  & $2 \mbox{p4}$    & 4  & 4 &  &  1  \\
p4m  & $K_4 \wr C_2 =[32,27]$     & $\mbox{p4},\mbox{pmm},\mbox{cmm}$  & $2\mbox{p4m},2\mbox{p4g}$       & 2  & 2 &  & 1  \\
p4g  & $C_4 \ltimes C_2^3 =[32,6]$    & $\mbox{p4},\mbox{cmm},\mbox{pgg}$  &       & 2  & 2 &  &  1 \\
\midrule
p3   & $A_4 = [12,3]$     &    &  &  4  & 4 &  &  1 \\
p3m1 & $S_4= [24,12]$     & $\mbox{p3}$   &  & 2 & 2 &  & 1 \\
p31m & $S_4= [24,12]$    & $\mbox{p3}$    & & 2  & 2 & &  1 \\
p6   & $C_2\times A_4 =[24,13]$     & $\mbox{p3}$       &  &  2 & 2       &  &  1  \\
p6m  & $C_2 \times S_4 = [48,48]$    & $\mbox{p6},\mbox{p3m1},\mbox{p31m}$   &  & 1 & 1  & &  1 \\
\bottomrule
\end{tabular}
\caption{The real structures of the 17 wallpaper groups.}
\label{table:wallpaper3}
\end{table}

\subsection{Real Structures}

In the right half of Table~\ref{table:wallpaper3} we collect all the data related to the classification of Real structures. As explained in Section \ref{sect:are_there}, the equivalence relation we are considering on magnetic groups corresponds to the weak-equivalence on Real structures considered in Section \ref{sect:weakequiv}. From Lemma ~\ref{lemma:weak_eq}, we have that 
\begin{equation} \label{Wallpaper:RTheta}
  \fR(G)/\overset{\text{magn}}{\sim}  = \bigsqcup_{\theta \in \Theta(G)} \Out_{\theta}(G) \backslash \fR_{\theta}(G) \, , \quad
  \Theta (G) \coloneqq \Out(G) \backslash \{ \theta \in \Out(G)\,\mid\, \theta^2 =1\}\, .
\end{equation}

The sets  The sets $\fR_{\theta}(G)$ are also easier in the crystallographic case because $Z(G)=T^Q$ is always torsion-free,
hence, classified by Proposition~\ref{prop:torsionfreeKS}.
Furthermore, $Z(G)=0$ in most  of the cases.

The column ${\sN}$ simply lists the size of the left-hand set in~\eqref{Wallpaper:RTheta}.
It can be counted using Fact~\ref{sect:are_there}.(I).
Indeed, ${\sN}$ is simply equal to 1 plus the number of times the corresponding orbit
appears columns 3 and 4.
For instance, we can see that {\bfseries pmm}
appears once in column 3 and two times in column 4. This means that $G$ of type {\bfseries pmm} has 4 Real structure:
one black-and-white of type 1,
two black-and-whites of type 2,
and finally one grey:
\[
(\mathrm{pmm} \leq \mathrm{p4m}), \ 
(\mathrm{pmm} \leq \mathrm{cmm}), \ 
(\mathrm{pmm} \leq \mathrm{pmm}), \
(\mathrm{pmm} \leq \mathrm{pmm}\times C_2).
\]
Simply listing the types of both $H$ and $G$ appears to be a convenient notation for the Real structures in the 2D case.
The only exception is that $(\mathrm{pm} \leq \mathrm{pm})$ may have two different meanings.
Therefore, we use  $(2\mathrm{pm} \leq \mathrm{pm})$ for one of them, where there are two subgroups conjugate by an automorphism.

In the row {\bfseries p2} and from the row {\bfseries pmm} down, 
the centre $Z(G)$ is trivial so that
all $\fR_{\theta}(G)$ are one point sets: these sizes are in the last column.
We also do not list individual members of $\Theta (G)$ -- we just list its size, equal to  ${\sN}$.
%
%
In the type  {\bfseries p2}, the group $\Out (G)$ is infinite. Yet the conjugacy classes in $\Theta (G)$ are easy to describe -- 
here are their representatives:
\[ I_2, \; U_2, \; X_2, \; X_4, \; e_1, \; (X_2,e_1), \; (X_2,e_1+e_2). \]

In the four types with infinite centre, the right hand side of Table~\ref{table:wallpaper3}
is filled slightly differently, with more information. The third column lists representatives in $\Theta(G)$. 
Then we use Theorem~\ref{thm:EMLv1} and Proposition~\ref{prop:C2cohomology} to derive
\[ \fR_{\theta}(G) =
\HH^{2}(C_2, Z(G)_{\theta}) =
\dfrac{Z(G)^0_{\theta}}{(\theta + 1)Z(G)_{\theta}} \, , \]
subject to vanishing of the obstruction, which always vanishes in the 2D case.
For each individual $\theta \in \Theta (G)$, we list $Z(G)_{\theta}$ in the penultimate column, cf. $C_2$-modules in Proposition~\ref{prop:torsionfreeKS}.
The cardinality $|\fR_{\theta}(G)|$ is found by Formula~\eqref{eqn:cohomofspecialmodules} that, in our case, boils down to
  \begin{equation}
  \HH^{2}(\bZ_+) = \bZ/2\bZ, \quad
  \HH^{2}(\bZ_-) =
  \HH^{2}(\bZ_d) = 0 \, .
  \end{equation}
  Finally, we fill the cardinalities of $\Out_\theta (G)$-orbits on $\fR_{\theta}(G)$ into the last column.
  For example, consider the first row for {\bfseries p1}. Here $G=\bZ^2$. Since $I$ acts trivially, we need to add an element $t$ to $G$
  such that $tg=gt$ for all $g\in G$. The four elements of $\fR_{\theta}(G)$ correspond to four possible choices of $t^2$:
  \[ t^2 = 0 \ \mbox{ or } \ t^2 = e_1 \ \mbox{ or } \ t^2 = e_2 \ \mbox{ or } \ t^2 = e_1+e_2 \, . \]
  Clearly, the last three choices are all magnetic equivalent, that is, they comprise a single $\Out_\theta (G)$-orbit. We write $1+3$
  because $\fR_{\theta}(G)$ has  two orbits, with 1 and 3 elements.

%% file: appendix-repswallpaper.tex
\chapter{Representations of Wallpaper Groups} \label{appendix_repswallpaper}

In this Appendix we describe the D-types of the irreducible representations of wallpaper groups, for the 46 non-trivial Real structures listed in Appendix \ref{appendix_wallpaper}.

Recall that for $H \leq G$ a Real structure, and an abelian normal subgroup $T \triangleleft G$ contained in $H$, the antilinear Clifford Theory of Section \ref{sect:cliffordtheory} describes a canonical bijection \eqref{Cliff_2f},
\[
\bigcup_{\chi \in S \backslash T^{\vee}}  {\Irr_{\chi}^{\hspace{0.1em} \star}(S_\chi)} \stackrel{1:1}{\longleftrightarrow}
{\Irr_{\bC}^{\hspace{0.1em} \star}(G)},
\]
where $S \coloneqq G / T$, $S \backslash T^{\vee}$ is the set of orbits of the action of $S$ on the group $T^{\vee}$ of $\bC$-valued characters of $T$, $S_{\chi}$ is the stabiliser of $\chi$ in $S$, and $\Irr_{\chi}^{\hspace{0.1em} \star}(S_\chi)$ is the set of irreducible antilinear $c_{\chi}$-projective representations of $S_{\chi}$.

The D-type of the representation of $\Irr_{\bC}^{\hspace{0.1em} \star}(G)$ corresponding to a given $W \in \Irr_{\chi}^{\hspace{0.1em} \star}(S_\chi)$ depends only on the triple ($S_{\chi}'$, $c_{\chi}'$, $W$) (in an appropriate sense, see Section \ref{sect:extendedstabiliser}), where $S_{\chi}'$ is the extended stabiliser of $\chi$, and $c_{\chi}' \in \HH^2(S_{\chi}', \bC^{\times}_{\kappa})$.

Therefore, determining the possible D-types of $G$ becomes an essentially finite problem, even though the set $T^{\vee}$ is infinite, there being only finitely options for the triple ($S_{\chi}'$, $c_{\chi}'$, $W$). Indeed, there are finitely many options for the subgroup $S_{\chi}'$ because $S$ is finite, $c_{\chi}'$ is an element of $\HH^2(S_{\chi}', \bC^{\times}_{\kappa})$, which is finite, and $W \in \Irr_{\chi}^{\hspace{0.1em} \star}(S_\chi)$, which is finite because $S_{\chi}$ is a finite group.

In this appendix we give, for each wallpaper group $G$ and index two subgroup $H$ of $G$ (up to equivalence), a table that lists the possible D-types for the Real structure $H \leq G$. The tables are structured as follows.

\begin{itemize}
    \item Generators are given describing $H$ as a subgroup of $G$.
    \item Here $G$ is given as a finitely presented group, with presentations found at \cite{GITHUB}.
    \item The final two given generators $f_1, f_2$ of $H$ freely generate the lattice $T \cong \bZ^2$ of $H$.
    \item For each $S_{\chi}$, we list generators for $S_{\chi}$ inside $S$, the ungraded part $S_{\chi,0}$, and $\HH^2(S_{\chi}, \bC^{\times}_{\pi})$.
    \item For each $S_{\chi}'$ extending $S_{\chi}$, we list generators for $S_{\chi}'$ inside $S \times C_2 = S \times \langle c \rangle$, where we use the notation $s \coloneqq (s, 1)$ and $\overline{s} \coloneqq (s, c)$, and $\HH^2(S_{\chi}', \bC^{\times}_{\kappa})$.
    \item $\chi$ is a character of $T$, presented as a pair $[r_1,r_2]$, where 
    \[
    \chi(f_1) = e^{2 \pi i r_1}, \qquad \chi(f_2) = e^{2 \pi i r_2},
    \]
    which realises the pair $S_{\chi}'$ and cocycle $c_{\chi}'$ (and so also the pair $S_{\chi}, c_{\chi}$).
    \item The final column $\fT(V)$ lists the D-types of the representations corresponding to $\Irr_{\chi}^{\hspace{0.1em} \star}(S_\chi)$.
\end{itemize}

We note that not all pairs $(K, c)$, where $K \leq S$ and $c \in \HH^2(S_{\chi}, \bC^{\times}_{\pi})$ are realised by some $\chi$ in these tables, and similarly not all pairs $(K', c')$, where $K' \leq S \times C_2$ and $c \in \HH^2(S_{\chi}', \bC^{\times}_{\kappa})$ are realised by some $\chi$. We only list those quadruples $(K, K', c_{\chi}, c_{\chi}')$ that are realised by some $\chi$. In producing the tables, we searched over all 40,000 characters $\chi$ of $T$ with $r_1,r_2 \in \frac{1}{25} \bZ / \bZ$. We expect this that all possible quadruples are already realised by characters from this set. Indeed, examining the tables, we see that increasing the search from $r_1,r_2 \in \frac{1}{8} \bZ / \bZ$ to $r_1, r_2 \in \frac{1}{25} \bZ / \bZ$ no new quadruple was found.

All data in the tables was computed using magma - the code we used to produce them is available at \cite{GITHUB}. Implementing the theory of Section \ref{sect:cliffordtheory} algorithmically raised some interesting mathematical points about how one computes certain quantities in finite time, and we describe these considerations in the final section of this appendix, Section \ref{sect:algconsiderations}, after giving the tables.

\newgeometry{left=1cm,right=1cm}

\section{Tables for Antilinear Representations}

\begin{table}[!htbp]
\centering
\small
\begin{tabular}{llll|lll|llll}
\toprule
$\boldsymbol{S_{\chi}}$ & & $\boldsymbol{S_{\chi,0}}$ & $\boldsymbol{\HH^2(S_{\chi}, \bC^{\times}_{\pi})}$ & $\boldsymbol{S'_{\chi}}$ & & $\boldsymbol{\HH^2(S'_{\chi}, \bC^{\times}_{\kappa})}$ & $\boldsymbol{\chi}$ & $\boldsymbol{c_{\chi}}$ & $\boldsymbol{c'_{\chi}}$ & $\boldsymbol{\fT(V)}$ \\
\midrule
$1$ & $\langle\emptyset\rangle$ & $1$ & $0$ & $C_2$ & $\langle a\rangle$ & $0$ & $[0,\frac{1}{4}]$ & $(0)$ & $(0)$ & IV \\
\midrule
$C_2$ & $\langle a\rangle$ & $1$ & $\mathbb Z/2\mathbb Z$ & $C_2^2$ & $\langle \overline{1},a\rangle$ & $(\mathbb Z/2\mathbb Z)^{2}$ & $[0,0]$ & $(0)$ & $(0,0)$ & I \\
 &  &  &  &  &  &  & $[0,\frac{1}{2}]$ & $(1)$ & $(0,1)$ & II \\
\bottomrule
\end{tabular}
\caption{$G = \text{p1}$, $H = \langle ab,ab^{-1} \rangle \cong \text{p1}$, so $S \cong C_2$ and $S_0 \cong 1$}
\label{table:repsGp1H1}
\end{table}

\begin{table}[!htbp]
\centering
\small
\begin{tabular}{llll|lll|llll}
\toprule
$\boldsymbol{S_{\chi}}$ & & $\boldsymbol{S_{\chi,0}}$ & $\boldsymbol{\HH^2(S_{\chi}, \bC^{\times}_{\pi})}$ & $\boldsymbol{S'_{\chi}}$ & & $\boldsymbol{\HH^2(S'_{\chi}, \bC^{\times}_{\kappa})}$ & $\boldsymbol{\chi}$ & $\boldsymbol{c_{\chi}}$ & $\boldsymbol{c'_{\chi}}$ & $\boldsymbol{\fT(V)}$ \\
\midrule
$C_2$ & $\langle r\rangle$ & $1$ & $\mathbb Z/2\mathbb Z$ & $C_2$ & $\langle \overline{r}\rangle$ & $\mathbb Z/2\mathbb Z$ & $[0,\frac{1}{4}]$ & $(0)$ & $(0)$ & V \\
\cmidrule(l){5-11}
 &  &  &  & $C_2^2$ & $\langle \overline{1},r\rangle$ & $(\mathbb Z/2\mathbb Z)^{2}$ & $[0,0]$ & $(0)$ & $(0,0)$ & I \\
\bottomrule
\end{tabular}
\caption{$G = \text{p2}$, $H = \langle u,v \rangle \cong \text{p1}$, so $S \cong C_2$ and $S_0 \cong 1$}
\label{table:repsGp2H1}
\end{table}

\begin{table}[!htbp]
\centering
\small
\begin{tabular}{llll|lll|llll}
\toprule
$\boldsymbol{S_{\chi}}$ & & $\boldsymbol{S_{\chi,0}}$ & $\boldsymbol{\HH^2(S_{\chi}, \bC^{\times}_{\pi})}$ & $\boldsymbol{S'_{\chi}}$ & & $\boldsymbol{\HH^2(S'_{\chi}, \bC^{\times}_{\kappa})}$ & $\boldsymbol{\chi}$ & $\boldsymbol{c_{\chi}}$ & $\boldsymbol{c'_{\chi}}$ & $\boldsymbol{\fT(V)}$ \\
\midrule
$C_2$ & $\langle ru\rangle$ & $1$ & $\mathbb Z/2\mathbb Z$ & $C_2^2$ & $\langle u,\overline{r}\rangle$ & $(\mathbb Z/2\mathbb Z)^{2}$ & $[0,\frac{1}{4}]$ & $(0)$ & $(0,0)$ & I \\
\midrule
$C_2^2$ & $\langle r,u\rangle$ & $C_2$ & $(\mathbb Z/2\mathbb Z)^{2}$ & $C_2^3$ & $\langle \overline{1},r,u\rangle$ & $(\mathbb Z/2\mathbb Z)^{4}$ & $[0,0]$ & $(0,0)$ & $(0,0,0,0)$ & I, I \\
 &  &  &  &  &  &  & $[0,\frac{1}{2}]$ & $(1,1)$ & $(0,0,1,1)$ & III \\
\bottomrule
\end{tabular}
\caption{$G = \text{p2}$, $H = \langle r,uv,uv^{-1} \rangle \cong \text{p2}$, so $S \cong C_2^2$ and $S_0 \cong C_2$}
\label{table:repsGp2H2}
\end{table}

\begin{table}[!htbp]
\centering
\small
\begin{tabular}{llll|lll|llll}
\toprule
$\boldsymbol{S_{\chi}}$ & & $\boldsymbol{S_{\chi,0}}$ & $\boldsymbol{\HH^2(S_{\chi}, \bC^{\times}_{\pi})}$ & $\boldsymbol{S'_{\chi}}$ & & $\boldsymbol{\HH^2(S'_{\chi}, \bC^{\times}_{\kappa})}$ & $\boldsymbol{\chi}$ & $\boldsymbol{c_{\chi}}$ & $\boldsymbol{c'_{\chi}}$ & $\boldsymbol{\fT(V)}$ \\
\midrule
$1$ & $\langle\emptyset\rangle$ & $1$ & $0$ & $1$ & $\langle\emptyset\rangle$ & $0$ & $[\frac{1}{4},\frac{1}{4}]$ & $(0)$ & $(0)$ & VII \\
\cmidrule(l){5-11}
 &  &  &  & $C_2$ & $\langle r\rangle$ & $0$ & $[\frac{1}{4},0]$ & $(0)$ & $(0)$ & IV \\
\midrule
$C_2$ & $\langle r\rangle$ & $1$ & $\mathbb Z/2\mathbb Z$ & $C_2$ & $\langle \overline{r}\rangle$ & $\mathbb Z/2\mathbb Z$ & $[0,\frac{1}{4}]$ & $(0)$ & $(0)$ & V \\
\cmidrule(l){5-11}
 &  &  &  & $C_2^2$ & $\langle \overline{1},r\rangle$ & $(\mathbb Z/2\mathbb Z)^{2}$ & $[0,0]$ & $(0)$ & $(0,0)$ & I \\
\bottomrule
\end{tabular}
\caption{$G = \text{pm}$, $H = \langle u,v \rangle \cong \text{p1}$, so $S \cong C_2$ and $S_0 \cong 1$}
\label{table:repsGpmH1}
\end{table}

\begin{table}[!htbp]
\centering
\small
\begin{tabular}{llll|lll|llll}
\toprule
$\boldsymbol{S_{\chi}}$ & & $\boldsymbol{S_{\chi,0}}$ & $\boldsymbol{\HH^2(S_{\chi}, \bC^{\times}_{\pi})}$ & $\boldsymbol{S'_{\chi}}$ & & $\boldsymbol{\HH^2(S'_{\chi}, \bC^{\times}_{\kappa})}$ & $\boldsymbol{\chi}$ & $\boldsymbol{c_{\chi}}$ & $\boldsymbol{c'_{\chi}}$ & $\boldsymbol{\fT(V)}$ \\
\midrule
$1$ & $\langle\emptyset\rangle$ & $1$ & $0$ & $C_2$ & $\langle u\rangle$ & $0$ & $[0,\frac{1}{4}]$ & $(0)$ & $(0)$ & IV \\
\midrule
$C_2$ & $\langle r\rangle$ & $C_2$ & $0$ & $C_2^2$ & $\langle r,u\rangle$ & $\mathbb Z/2\mathbb Z$ & $[\frac{1}{4},\frac{1}{4}]$ & $(0)$ & $(0)$ & IV, IV \\
\midrule
$C_2$ & $\langle ru\rangle$ & $1$ & $\mathbb Z/2\mathbb Z$ & $C_2^2$ & $\langle u,\overline{r}\rangle$ & $(\mathbb Z/2\mathbb Z)^{2}$ & $[\frac{1}{4},\frac{3}{4}]$ & $(0)$ & $(0,0)$ & I \\
\midrule
$C_2$ & $\langle u\rangle$ & $1$ & $\mathbb Z/2\mathbb Z$ & $C_2^2$ & $\langle \overline{1},u\rangle$ & $(\mathbb Z/2\mathbb Z)^{2}$ & $[0,\frac{1}{2}]$ & $(1)$ & $(0,1)$ & II \\
\midrule
$C_2^2$ & $\langle r,u\rangle$ & $C_2$ & $(\mathbb Z/2\mathbb Z)^{2}$ & $C_2^3$ & $\langle \overline{1},r,u\rangle$ & $(\mathbb Z/2\mathbb Z)^{4}$ & $[0,0]$ & $(0,0)$ & $(0,0,0,0)$ & I, I \\
\bottomrule
\end{tabular}
\caption{$G = \text{pm}$, $H = \langle r,uv,uv^{-1} \rangle \cong \text{cm}$, so $S \cong C_2^2$ and $S_0 \cong C_2$}
\label{table:repsGpmH2}
\end{table}

\begin{table}[!htbp]
\centering
\small
\begin{tabular}{llll|lll|llll}
\toprule
$\boldsymbol{S_{\chi}}$ & & $\boldsymbol{S_{\chi,0}}$ & $\boldsymbol{\HH^2(S_{\chi}, \bC^{\times}_{\pi})}$ & $\boldsymbol{S'_{\chi}}$ & & $\boldsymbol{\HH^2(S'_{\chi}, \bC^{\times}_{\kappa})}$ & $\boldsymbol{\chi}$ & $\boldsymbol{c_{\chi}}$ & $\boldsymbol{c'_{\chi}}$ & $\boldsymbol{\fT(V)}$ \\
\midrule
$1$ & $\langle\emptyset\rangle$ & $1$ & $0$ & $C_2$ & $\langle v\rangle$ & $0$ & $[\frac{1}{4},\frac{1}{4}]$ & $(0)$ & $(0)$ & IV \\
\midrule
$C_2$ & $\langle rv\rangle$ & $1$ & $\mathbb Z/2\mathbb Z$ & $C_2^2$ & $\langle v,\overline{r}\rangle$ & $(\mathbb Z/2\mathbb Z)^{2}$ & $[0,\frac{1}{4}]$ & $(0)$ & $(0,0)$ & I \\
\midrule
$C_2$ & $\langle r\rangle$ & $C_2$ & $0$ & $C_2^2$ & $\langle r,v\rangle$ & $\mathbb Z/2\mathbb Z$ & $[\frac{1}{4},0]$ & $(0)$ & $(0)$ & IV, IV \\
 &  &  &  &  &  &  & $[\frac{1}{4},\frac{1}{2}]$ & $(0)$ & $(1)$ & VII, VII \\
\midrule
$C_2^2$ & $\langle r,v\rangle$ & $C_2$ & $(\mathbb Z/2\mathbb Z)^{2}$ & $C_2^3$ & $\langle \overline{1},r,v\rangle$ & $(\mathbb Z/2\mathbb Z)^{4}$ & $[0,0]$ & $(0,0)$ & $(0,0,0,0)$ & I, I \\
 &  &  &  &  &  &  & $[0,\frac{1}{2}]$ & $(1,1)$ & $(0,0,1,1)$ & III \\
\bottomrule
\end{tabular}
\caption{$G = \text{pm}$, $H = \langle r,u,v^{2} \rangle \cong \text{2pm}$, so $S \cong C_2^2$ and $S_0 \cong C_2$}
\label{table:repsGpmH3}
\end{table}

\begin{table}[!htbp]
\centering
\small
\begin{tabular}{llll|lll|llll}
\toprule
$\boldsymbol{S_{\chi}}$ & & $\boldsymbol{S_{\chi,0}}$ & $\boldsymbol{\HH^2(S_{\chi}, \bC^{\times}_{\pi})}$ & $\boldsymbol{S'_{\chi}}$ & & $\boldsymbol{\HH^2(S'_{\chi}, \bC^{\times}_{\kappa})}$ & $\boldsymbol{\chi}$ & $\boldsymbol{c_{\chi}}$ & $\boldsymbol{c'_{\chi}}$ & $\boldsymbol{\fT(V)}$ \\
\midrule
$1$ & $\langle\emptyset\rangle$ & $1$ & $0$ & $C_2$ & $\langle u\rangle$ & $0$ & $[\frac{1}{4},\frac{1}{4}]$ & $(0)$ & $(0)$ & IV \\
\midrule
$C_2$ & $\langle ru\rangle$ & $1$ & $\mathbb Z/2\mathbb Z$ & $C_2^2$ & $\langle u,\overline{r}\rangle$ & $(\mathbb Z/2\mathbb Z)^{2}$ & $[0,\frac{1}{4}]$ & $(0)$ & $(0,0)$ & I \\
 &  &  &  &  &  &  & $[\frac{1}{2},\frac{1}{4}]$ & $(1)$ & $(0,1)$ & II \\
\midrule
$C_2$ & $\langle r\rangle$ & $C_2$ & $0$ & $C_2^2$ & $\langle r,u\rangle$ & $\mathbb Z/2\mathbb Z$ & $[\frac{1}{4},0]$ & $(0)$ & $(0)$ & IV, IV \\
\midrule
$C_2^2$ & $\langle r,u\rangle$ & $C_2$ & $(\mathbb Z/2\mathbb Z)^{2}$ & $C_2^3$ & $\langle \overline{1},r,u\rangle$ & $(\mathbb Z/2\mathbb Z)^{4}$ & $[0,0]$ & $(0,0)$ & $(0,0,0,0)$ & I, I \\
 &  &  &  &  &  &  & $[\frac{1}{2},0]$ & $(0,1)$ & $(0,0,1,0)$ & II, II \\
\bottomrule
\end{tabular}
\caption{$G = \text{pm}$, $H = \langle r,u^{2},v \rangle \cong \text{pm}$, so $S \cong C_2^2$ and $S_0 \cong C_2$}
\label{table:repsGpmH4}
\end{table}

\begin{table}[!htbp]
\centering
\small
\begin{tabular}{llll|lll|llll}
\toprule
$\boldsymbol{S_{\chi}}$ & & $\boldsymbol{S_{\chi,0}}$ & $\boldsymbol{\HH^2(S_{\chi}, \bC^{\times}_{\pi})}$ & $\boldsymbol{S'_{\chi}}$ & & $\boldsymbol{\HH^2(S'_{\chi}, \bC^{\times}_{\kappa})}$ & $\boldsymbol{\chi}$ & $\boldsymbol{c_{\chi}}$ & $\boldsymbol{c'_{\chi}}$ & $\boldsymbol{\fT(V)}$ \\
\midrule
$1$ & $\langle\emptyset\rangle$ & $1$ & $0$ & $C_2$ & $\langle u\rangle$ & $0$ & $[\frac{1}{4},\frac{1}{4}]$ & $(0)$ & $(0)$ & IV \\
\midrule
$C_2$ & $\langle r\rangle$ & $1$ & $\mathbb Z/2\mathbb Z$ & $C_2^2$ & $\langle u,\overline{r}\rangle$ & $(\mathbb Z/2\mathbb Z)^{2}$ & $[0,\frac{1}{4}]$ & $(0)$ & $(0,0)$ & I \\
 &  &  &  &  &  &  & $[\frac{1}{2},\frac{1}{4}]$ & $(0)$ & $(0,1)$ & IX \\
\midrule
$C_2$ & $\langle ru\rangle$ & $C_2$ & $0$ & $C_2^2$ & $\langle r,u\rangle$ & $\mathbb Z/2\mathbb Z$ & $[\frac{1}{4},0]$ & $(0)$ & $(0)$ & IV, IV \\
\midrule
$C_2^2$ & $\langle r,u\rangle$ & $C_2$ & $(\mathbb Z/2\mathbb Z)^{2}$ & $C_2^3$ & $\langle \overline{1},r,u\rangle$ & $(\mathbb Z/2\mathbb Z)^{4}$ & $[0,0]$ & $(0,0)$ & $(0,0,0,0)$ & I, I \\
 &  &  &  &  &  &  & $[\frac{1}{2},0]$ & $(0,1)$ & $(0,0,1,0)$ & IV \\
\bottomrule
\end{tabular}
\caption{$G = \text{pm}$, $H = \langle ru,u^{2},v \rangle \cong \text{pg}$, so $S \cong C_2^2$ and $S_0 \cong C_2$}
\label{table:repsGpmH5}
\end{table}

\begin{table}[!htbp]
\centering
\small
\begin{tabular}{llll|lll|llll}
\toprule
$\boldsymbol{S_{\chi}}$ & & $\boldsymbol{S_{\chi,0}}$ & $\boldsymbol{\HH^2(S_{\chi}, \bC^{\times}_{\pi})}$ & $\boldsymbol{S'_{\chi}}$ & & $\boldsymbol{\HH^2(S'_{\chi}, \bC^{\times}_{\kappa})}$ & $\boldsymbol{\chi}$ & $\boldsymbol{c_{\chi}}$ & $\boldsymbol{c'_{\chi}}$ & $\boldsymbol{\fT(V)}$ \\
\midrule
$1$ & $\langle\emptyset\rangle$ & $1$ & $0$ & $1$ & $\langle\emptyset\rangle$ & $0$ & $[\frac{1}{4},\frac{1}{4}]$ & $(0)$ & $(0)$ & VII \\
\cmidrule(l){5-11}
 &  &  &  & $C_2$ & $\langle r\rangle$ & $0$ & $[\frac{1}{4},0]$ & $(0)$ & $(0)$ & IV \\
\midrule
$C_2$ & $\langle r\rangle$ & $1$ & $\mathbb Z/2\mathbb Z$ & $C_2$ & $\langle \overline{r}\rangle$ & $\mathbb Z/2\mathbb Z$ & $[0,\frac{1}{4}]$ & $(0)$ & $(0)$ & V \\
 &  &  &  &  &  &  & $[\frac{1}{2},\frac{1}{4}]$ & $(1)$ & $(1)$ & VI \\
\cmidrule(l){5-11}
 &  &  &  & $C_2^2$ & $\langle \overline{1},r\rangle$ & $(\mathbb Z/2\mathbb Z)^{2}$ & $[0,0]$ & $(0)$ & $(0,0)$ & I \\
 &  &  &  &  &  &  & $[\frac{1}{2},0]$ & $(1)$ & $(0,1)$ & II \\
\bottomrule
\end{tabular}
\caption{$G = \text{pg}$, $H = \langle u,v \rangle \cong \text{p1}$, so $S \cong C_2$ and $S_0 \cong 1$}
\label{table:repsGpgH1}
\end{table}

\begin{table}[!htbp]
\centering
\small
\begin{tabular}{llll|lll|llll}
\toprule
$\boldsymbol{S_{\chi}}$ & & $\boldsymbol{S_{\chi,0}}$ & $\boldsymbol{\HH^2(S_{\chi}, \bC^{\times}_{\pi})}$ & $\boldsymbol{S'_{\chi}}$ & & $\boldsymbol{\HH^2(S'_{\chi}, \bC^{\times}_{\kappa})}$ & $\boldsymbol{\chi}$ & $\boldsymbol{c_{\chi}}$ & $\boldsymbol{c'_{\chi}}$ & $\boldsymbol{\fT(V)}$ \\
\midrule
$1$ & $\langle\emptyset\rangle$ & $1$ & $0$ & $C_2$ & $\langle v\rangle$ & $0$ & $[\frac{1}{4},\frac{1}{4}]$ & $(0)$ & $(0)$ & IV \\
\midrule
$C_2$ & $\langle rv\rangle$ & $1$ & $\mathbb Z/2\mathbb Z$ & $C_2^2$ & $\langle v,\overline{r}\rangle$ & $(\mathbb Z/2\mathbb Z)^{2}$ & $[0,\frac{1}{4}]$ & $(0)$ & $(0,0)$ & I \\
 &  &  &  &  &  &  & $[\frac{1}{2},\frac{1}{4}]$ & $(1)$ & $(1,0)$ & VIII \\
\midrule
$C_2$ & $\langle r\rangle$ & $C_2$ & $0$ & $C_2^2$ & $\langle r,v\rangle$ & $\mathbb Z/2\mathbb Z$ & $[\frac{1}{4},0]$ & $(0)$ & $(0)$ & IV, IV \\
 &  &  &  &  &  &  & $[\frac{1}{4},\frac{1}{2}]$ & $(0)$ & $(1)$ & VII, VII \\
\midrule
$C_2^2$ & $\langle r,v\rangle$ & $C_2$ & $(\mathbb Z/2\mathbb Z)^{2}$ & $C_2^3$ & $\langle \overline{1},r,v\rangle$ & $(\mathbb Z/2\mathbb Z)^{4}$ & $[0,0]$ & $(0,0)$ & $(0,0,0,0)$ & I, I \\
 &  &  &  &  &  &  & $[0,\frac{1}{2}]$ & $(1,1)$ & $(0,0,1,1)$ & III \\
 &  &  &  &  &  &  & $[\frac{1}{2},0]$ & $(1,0)$ & $(0,1,0,0)$ & IV \\
 &  &  &  &  &  &  & $[\frac{1}{2},\frac{1}{2}]$ & $(0,1)$ & $(0,1,1,1)$ & VI, VI \\
\bottomrule
\end{tabular}
\caption{$G = \text{pg}$, $H = \langle r,u,v^{2} \rangle \cong \text{pg}$, so $S \cong C_2^2$ and $S_0 \cong C_2$}
\label{table:repsGpgH2}
\end{table}

\begin{table}[!htbp]
\centering
\small
\begin{tabular}{llll|lll|llll}
\toprule
$\boldsymbol{S_{\chi}}$ & & $\boldsymbol{S_{\chi,0}}$ & $\boldsymbol{\HH^2(S_{\chi}, \bC^{\times}_{\pi})}$ & $\boldsymbol{S'_{\chi}}$ & & $\boldsymbol{\HH^2(S'_{\chi}, \bC^{\times}_{\kappa})}$ & $\boldsymbol{\chi}$ & $\boldsymbol{c_{\chi}}$ & $\boldsymbol{c'_{\chi}}$ & $\boldsymbol{\fT(V)}$ \\
\midrule
$1$ & $\langle\emptyset\rangle$ & $1$ & $0$ & $1$ & $\langle\emptyset\rangle$ & $0$ & $[0,\frac{1}{4}]$ & $(0)$ & $(0)$ & VII \\
\cmidrule(l){5-11}
 &  &  &  & $C_2$ & $\langle r\rangle$ & $0$ & $[\frac{1}{4},\frac{1}{4}]$ & $(0)$ & $(0)$ & IV \\
\cmidrule(l){5-11}
 &  &  &  & $C_2$ & $\langle \overline{1}\rangle$ & $\mathbb Z/2\mathbb Z$ & $[0,\frac{1}{2}]$ & $(0)$ & $(0)$ & III \\
\midrule
$C_2$ & $\langle r\rangle$ & $1$ & $\mathbb Z/2\mathbb Z$ & $C_2$ & $\langle \overline{r}\rangle$ & $\mathbb Z/2\mathbb Z$ & $[\frac{1}{4},\frac{3}{4}]$ & $(0)$ & $(0)$ & V \\
\cmidrule(l){5-11}
 &  &  &  & $C_2^2$ & $\langle \overline{1},r\rangle$ & $(\mathbb Z/2\mathbb Z)^{2}$ & $[0,0]$ & $(0)$ & $(0,0)$ & I \\
\bottomrule
\end{tabular}
\caption{$G = \text{cm}$, $H = \langle u,v \rangle \cong \text{p1}$, so $S \cong C_2$ and $S_0 \cong 1$}
\label{table:repsGcmH1}
\end{table}

\begin{table}[!htbp]
\centering
\small
\begin{tabular}{llll|lll|llll}
\toprule
$\boldsymbol{S_{\chi}}$ & & $\boldsymbol{S_{\chi,0}}$ & $\boldsymbol{\HH^2(S_{\chi}, \bC^{\times}_{\pi})}$ & $\boldsymbol{S'_{\chi}}$ & & $\boldsymbol{\HH^2(S'_{\chi}, \bC^{\times}_{\kappa})}$ & $\boldsymbol{\chi}$ & $\boldsymbol{c_{\chi}}$ & $\boldsymbol{c'_{\chi}}$ & $\boldsymbol{\fT(V)}$ \\
\midrule
$1$ & $\langle\emptyset\rangle$ & $1$ & $0$ & $C_2$ & $\langle u\rangle$ & $0$ & $[\frac{1}{4},\frac{1}{4}]$ & $(0)$ & $(0)$ & IV \\
\midrule
$C_2$ & $\langle ru\rangle$ & $1$ & $\mathbb Z/2\mathbb Z$ & $C_2^2$ & $\langle u,\overline{r}\rangle$ & $(\mathbb Z/2\mathbb Z)^{2}$ & $[0,\frac{1}{4}]$ & $(0)$ & $(0,0)$ & I \\
 &  &  &  &  &  &  & $[\frac{1}{2},\frac{1}{4}]$ & $(1)$ & $(0,1)$ & II \\
\midrule
$C_2$ & $\langle r\rangle$ & $C_2$ & $0$ & $C_2^2$ & $\langle r,u\rangle$ & $\mathbb Z/2\mathbb Z$ & $[\frac{1}{4},0]$ & $(0)$ & $(0)$ & IV, IV \\
 &  &  &  &  &  &  & $[\frac{1}{4},\frac{1}{2}]$ & $(0)$ & $(1)$ & VII, VII \\
\midrule
$C_2^2$ & $\langle r,u\rangle$ & $C_2$ & $(\mathbb Z/2\mathbb Z)^{2}$ & $C_2^3$ & $\langle \overline{1},r,u\rangle$ & $(\mathbb Z/2\mathbb Z)^{4}$ & $[0,0]$ & $(0,0)$ & $(0,0,0,0)$ & I, I \\
 &  &  &  &  &  &  & $[0,\frac{1}{2}]$ & $(1,1)$ & $(0,0,1,1)$ & III \\
 &  &  &  &  &  &  & $[\frac{1}{2},0]$ & $(0,1)$ & $(0,0,1,0)$ & II, II \\
 &  &  &  &  &  &  & $[\frac{1}{2},\frac{1}{2}]$ & $(1,0)$ & $(0,0,0,1)$ & III \\
\bottomrule
\end{tabular}
\caption{$G = \text{cm}$, $H = \langle r,uv,uv^{-1} \rangle \cong \text{pm}$, so $S \cong C_2^2$ and $S_0 \cong C_2$}
\label{table:repsGcmH2}
\end{table}

\begin{table}[!htbp]
\centering
\small
\begin{tabular}{llll|lll|llll}
\toprule
$\boldsymbol{S_{\chi}}$ & & $\boldsymbol{S_{\chi,0}}$ & $\boldsymbol{\HH^2(S_{\chi}, \bC^{\times}_{\pi})}$ & $\boldsymbol{S'_{\chi}}$ & & $\boldsymbol{\HH^2(S'_{\chi}, \bC^{\times}_{\kappa})}$ & $\boldsymbol{\chi}$ & $\boldsymbol{c_{\chi}}$ & $\boldsymbol{c'_{\chi}}$ & $\boldsymbol{\fT(V)}$ \\
\midrule
$1$ & $\langle\emptyset\rangle$ & $1$ & $0$ & $C_2$ & $\langle u\rangle$ & $0$ & $[\frac{1}{4},\frac{1}{4}]$ & $(0)$ & $(0)$ & IV \\
\midrule
$C_2$ & $\langle r\rangle$ & $1$ & $\mathbb Z/2\mathbb Z$ & $C_2^2$ & $\langle u,\overline{r}\rangle$ & $(\mathbb Z/2\mathbb Z)^{2}$ & $[0,\frac{1}{4}]$ & $(0)$ & $(0,0)$ & I \\
 &  &  &  &  &  &  & $[\frac{1}{2},\frac{1}{4}]$ & $(0)$ & $(0,1)$ & IX \\
\midrule
$C_2$ & $\langle ru\rangle$ & $C_2$ & $0$ & $C_2^2$ & $\langle r,u\rangle$ & $\mathbb Z/2\mathbb Z$ & $[\frac{1}{4},0]$ & $(0)$ & $(0)$ & IV, IV \\
 &  &  &  &  &  &  & $[\frac{1}{4},\frac{1}{2}]$ & $(0)$ & $(1)$ & VII, VII \\
\midrule
$C_2^2$ & $\langle r,u\rangle$ & $C_2$ & $(\mathbb Z/2\mathbb Z)^{2}$ & $C_2^3$ & $\langle \overline{1},r,u\rangle$ & $(\mathbb Z/2\mathbb Z)^{4}$ & $[0,0]$ & $(0,0)$ & $(0,0,0,0)$ & I, I \\
 &  &  &  &  &  &  & $[0,\frac{1}{2}]$ & $(0,1)$ & $(0,0,1,1)$ & III \\
 &  &  &  &  &  &  & $[\frac{1}{2},0]$ & $(0,1)$ & $(0,0,1,0)$ & IV \\
 &  &  &  &  &  &  & $[\frac{1}{2},\frac{1}{2}]$ & $(0,0)$ & $(0,0,0,1)$ & V, V \\
\bottomrule
\end{tabular}
\caption{$G = \text{cm}$, $H = \langle ru,uv,uv^{-1} \rangle \cong \text{pg}$, so $S \cong C_2^2$ and $S_0 \cong C_2$}
\label{table:repsGcmH3}
\end{table}

\begin{table}[!htbp]
\centering
\small
\begin{tabular}{llll|lll|llll}
\toprule
$\boldsymbol{S_{\chi}}$ & & $\boldsymbol{S_{\chi,0}}$ & $\boldsymbol{\HH^2(S_{\chi}, \bC^{\times}_{\pi})}$ & $\boldsymbol{S'_{\chi}}$ & & $\boldsymbol{\HH^2(S'_{\chi}, \bC^{\times}_{\kappa})}$ & $\boldsymbol{\chi}$ & $\boldsymbol{c_{\chi}}$ & $\boldsymbol{c'_{\chi}}$ & $\boldsymbol{\fT(V)}$ \\
\midrule
$1$ & $\langle\emptyset\rangle$ & $1$ & $0$ & $C_2$ & $\langle \overline{st}\rangle$ & $\mathbb Z/2\mathbb Z$ & $[\frac{1}{4},\frac{1}{4}]$ & $(0)$ & $(0)$ & III \\
\midrule
$C_2$ & $\langle t\rangle$ & $1$ & $\mathbb Z/2\mathbb Z$ & $C_2^2$ & $\langle s,\overline{t}\rangle$ & $(\mathbb Z/2\mathbb Z)^{2}$ & $[0,\frac{1}{4}]$ & $(0)$ & $(0,0)$ & I \\
\midrule
$C_2$ & $\langle s\rangle$ & $1$ & $\mathbb Z/2\mathbb Z$ & $C_2^2$ & $\langle t,\overline{s}\rangle$ & $(\mathbb Z/2\mathbb Z)^{2}$ & $[\frac{1}{4},0]$ & $(0)$ & $(0,0)$ & I \\
\midrule
$C_2^2$ & $\langle s,t\rangle$ & $C_2$ & $(\mathbb Z/2\mathbb Z)^{2}$ & $C_2^3$ & $\langle \overline{1},s,t\rangle$ & $(\mathbb Z/2\mathbb Z)^{4}$ & $[0,0]$ & $(0,0)$ & $(0,0,0,0)$ & I, I \\
\bottomrule
\end{tabular}
\caption{$G = \text{pmm}$, $H = \langle st,x,y \rangle \cong \text{p2}$, so $S \cong C_2^2$ and $S_0 \cong C_2$}
\label{table:repsGpmmH1}
\end{table}

\begin{table}[!htbp]
\centering
\small
\begin{tabular}{llll|lll|llll}
\toprule
$\boldsymbol{S_{\chi}}$ & & $\boldsymbol{S_{\chi,0}}$ & $\boldsymbol{\HH^2(S_{\chi}, \bC^{\times}_{\pi})}$ & $\boldsymbol{S'_{\chi}}$ & & $\boldsymbol{\HH^2(S'_{\chi}, \bC^{\times}_{\kappa})}$ & $\boldsymbol{\chi}$ & $\boldsymbol{c_{\chi}}$ & $\boldsymbol{c'_{\chi}}$ & $\boldsymbol{\fT(V)}$ \\
\midrule
$C_2$ & $\langle st\rangle$ & $1$ & $\mathbb Z/2\mathbb Z$ & $C_2$ & $\langle \overline{st}\rangle$ & $\mathbb Z/2\mathbb Z$ & $[\frac{1}{4},\frac{1}{4}]$ & $(0)$ & $(0)$ & V \\
\cmidrule(l){5-11}
 &  &  &  & $C_2^2$ & $\langle t,\overline{s}\rangle$ & $(\mathbb Z/2\mathbb Z)^{2}$ & $[\frac{1}{4},0]$ & $(0)$ & $(0,0)$ & I \\
\midrule
$C_2^2$ & $\langle s,t\rangle$ & $C_2$ & $(\mathbb Z/2\mathbb Z)^{2}$ & $C_2^2$ & $\langle s,\overline{t}\rangle$ & $(\mathbb Z/2\mathbb Z)^{2}$ & $[0,\frac{1}{4}]$ & $(0,0)$ & $(0,0)$ & V, V \\
\cmidrule(l){5-11}
 &  &  &  & $C_2^3$ & $\langle \overline{1},s,t\rangle$ & $(\mathbb Z/2\mathbb Z)^{4}$ & $[0,0]$ & $(0,0)$ & $(0,0,0,0)$ & I, I \\
\bottomrule
\end{tabular}
\caption{$G = \text{pmm}$, $H = \langle s,x,y \rangle \cong \text{pm}$, so $S \cong C_2^2$ and $S_0 \cong C_2$}
\label{table:repsGpmmH2}
\end{table}

\begin{table}[!htbp]
\centering
\small
\begin{tabular}{llll|lll|llll}
\toprule
$\boldsymbol{S_{\chi}}$ & & $\boldsymbol{S_{\chi,0}}$ & $\boldsymbol{\HH^2(S_{\chi}, \bC^{\times}_{\pi})}$ & $\boldsymbol{S'_{\chi}}$ & & $\boldsymbol{\HH^2(S'_{\chi}, \bC^{\times}_{\kappa})}$ & $\boldsymbol{\chi}$ & $\boldsymbol{c_{\chi}}$ & $\boldsymbol{c'_{\chi}}$ & $\boldsymbol{\fT(V)}$ \\
\midrule
$C_2$ & $\langle sty\rangle$ & $1$ & $\mathbb Z/2\mathbb Z$ & $C_2^2$ & $\langle y,\overline{st}\rangle$ & $(\mathbb Z/2\mathbb Z)^{2}$ & $[\frac{1}{4},\frac{1}{4}]$ & $(0)$ & $(0,0)$ & I \\
\midrule
$C_2^2$ & $\langle s,ty\rangle$ & $C_2$ & $(\mathbb Z/2\mathbb Z)^{2}$ & $C_2^3$ & $\langle s,y,\overline{t}\rangle$ & $(\mathbb Z/2\mathbb Z)^{4}$ & $[0,\frac{1}{4}]$ & $(0,0)$ & $(0,0,0,0)$ & I, I \\
\midrule
$C_2^2$ & $\langle t,sy\rangle$ & $C_2$ & $(\mathbb Z/2\mathbb Z)^{2}$ & $C_2^3$ & $\langle t,y,\overline{s}\rangle$ & $(\mathbb Z/2\mathbb Z)^{4}$ & $[\frac{1}{4},0]$ & $(0,0)$ & $(0,0,0,0)$ & I, I \\
 &  &  &  &  &  &  & $[\frac{1}{4},\frac{1}{2}]$ & $(1,1)$ & $(0,0,1,1)$ & III \\
\midrule
$C_2^3$ & $\langle s,t,y\rangle$ & $C_2^2$ & $(\mathbb Z/2\mathbb Z)^{4}$ & $C_2^4$ & $\langle \overline{1},s,t,y\rangle$ & $(\mathbb Z/2\mathbb Z)^{7}$ & $[0,0]$ & $(0,0,0,0)$ & $(0,0,0,0,0,0,0)$ & I, I, I, I \\
 &  &  &  &  &  &  & $[0,\frac{1}{2}]$ & $(0,1,1,0)$ & $(0,0,0,1,0,0,1)$ & III, III \\
\bottomrule
\end{tabular}
\caption{$G = \text{pmm}$, $H = \langle s,t,x,y^{2} \rangle \cong \text{pmm}$, so $S \cong C_2^3$ and $S_0 \cong C_2^2$}
\label{table:repsGpmmH3}
\end{table}

\begin{table}[!htbp]
\centering
\small
\begin{tabular}{llll|lll|llll}
\toprule
$\boldsymbol{S_{\chi}}$ & & $\boldsymbol{S_{\chi,0}}$ & $\boldsymbol{\HH^2(S_{\chi}, \bC^{\times}_{\pi})}$ & $\boldsymbol{S'_{\chi}}$ & & $\boldsymbol{\HH^2(S'_{\chi}, \bC^{\times}_{\kappa})}$ & $\boldsymbol{\chi}$ & $\boldsymbol{c_{\chi}}$ & $\boldsymbol{c'_{\chi}}$ & $\boldsymbol{\fT(V)}$ \\
\midrule
$C_2$ & $\langle st\rangle$ & $1$ & $\mathbb Z/2\mathbb Z$ & $C_2^2$ & $\langle y,\overline{st}\rangle$ & $(\mathbb Z/2\mathbb Z)^{2}$ & $[\frac{1}{4},\frac{1}{4}]$ & $(0)$ & $(0,0)$ & I \\
\midrule
$C_2^2$ & $\langle st,sy\rangle$ & $C_2$ & $(\mathbb Z/2\mathbb Z)^{2}$ & $C_2^3$ & $\langle s,y,\overline{t}\rangle$ & $(\mathbb Z/2\mathbb Z)^{4}$ & $[0,\frac{1}{4}]$ & $(0,0)$ & $(0,0,0,0)$ & I, I \\
\midrule
$C_2^2$ & $\langle s,t\rangle$ & $C_2$ & $(\mathbb Z/2\mathbb Z)^{2}$ & $C_2^3$ & $\langle t,y,\overline{s}\rangle$ & $(\mathbb Z/2\mathbb Z)^{4}$ & $[\frac{1}{4},0]$ & $(0,0)$ & $(0,0,0,0)$ & I, I \\
 &  &  &  &  &  &  & $[\frac{1}{4},\frac{1}{2}]$ & $(0,0)$ & $(0,0,1,1)$ & V, V \\
\midrule
$C_2^3$ & $\langle s,t,y\rangle$ & $C_2^2$ & $(\mathbb Z/2\mathbb Z)^{4}$ & $C_2^4$ & $\langle \overline{1},s,t,y\rangle$ & $(\mathbb Z/2\mathbb Z)^{7}$ & $[0,0]$ & $(0,0,0,0)$ & $(0,0,0,0,0,0,0)$ & I, I, I, I \\
 &  &  &  &  &  &  & $[0,\frac{1}{2}]$ & $(1,1,0,0)$ & $(0,0,0,1,0,0,1)$ & I \\
\bottomrule
\end{tabular}
\caption{$G = \text{pmm}$, $H = \langle sy,t,x,y^{2} \rangle \cong \text{pmg}$, so $S \cong C_2^3$ and $S_0 \cong C_2^2$}
\label{table:repsGpmmH4}
\end{table}

\begin{table}[!htbp]
\centering
\small
\begin{tabular}{llll|lll|llll}
\toprule
$\boldsymbol{S_{\chi}}$ & & $\boldsymbol{S_{\chi,0}}$ & $\boldsymbol{\HH^2(S_{\chi}, \bC^{\times}_{\pi})}$ & $\boldsymbol{S'_{\chi}}$ & & $\boldsymbol{\HH^2(S'_{\chi}, \bC^{\times}_{\kappa})}$ & $\boldsymbol{\chi}$ & $\boldsymbol{c_{\chi}}$ & $\boldsymbol{c'_{\chi}}$ & $\boldsymbol{\fT(V)}$ \\
\midrule
$C_2$ & $\langle stx\rangle$ & $1$ & $\mathbb Z/2\mathbb Z$ & $C_2^2$ & $\langle x,\overline{st}\rangle$ & $(\mathbb Z/2\mathbb Z)^{2}$ & $[0,\frac{1}{4}]$ & $(0)$ & $(0,0)$ & I \\
\midrule
$C_2^2$ & $\langle t,sx\rangle$ & $C_2$ & $(\mathbb Z/2\mathbb Z)^{2}$ & $C_2^3$ & $\langle t,x,\overline{s}\rangle$ & $(\mathbb Z/2\mathbb Z)^{4}$ & $[\frac{1}{4},\frac{1}{4}]$ & $(0,0)$ & $(0,0,0,0)$ & I, I \\
\midrule
$C_2^2$ & $\langle s,tx\rangle$ & $C_2$ & $(\mathbb Z/2\mathbb Z)^{2}$ & $C_2^3$ & $\langle s,x,\overline{t}\rangle$ & $(\mathbb Z/2\mathbb Z)^{4}$ & $[\frac{1}{4},\frac{3}{4}]$ & $(0,0)$ & $(0,0,0,0)$ & I, I \\
\midrule
$C_2^2$ & $\langle x,st\rangle$ & $C_2$ & $(\mathbb Z/2\mathbb Z)^{2}$ & $C_2^3$ & $\langle \overline{1},x,st\rangle$ & $(\mathbb Z/2\mathbb Z)^{4}$ & $[0,\frac{1}{2}]$ & $(1,1)$ & $(0,1,0,1)$ & III \\
\midrule
$C_2^3$ & $\langle s,t,x\rangle$ & $C_2^2$ & $(\mathbb Z/2\mathbb Z)^{4}$ & $C_2^4$ & $\langle \overline{1},s,t,x\rangle$ & $(\mathbb Z/2\mathbb Z)^{7}$ & $[0,0]$ & $(0,0,0,0)$ & $(0,0,0,0,0,0,0)$ & I, I, I, I \\
\bottomrule
\end{tabular}
\caption{$G = \text{pmm}$, $H = \langle s,t,xy,xy^{-1} \rangle \cong \text{cmm}$, so $S \cong C_2^3$ and $S_0 \cong C_2^2$}
\label{table:repsGpmmH5}
\end{table}

\begin{table}[!htbp]
\centering
\small
\begin{tabular}{llll|lll|llll}
\toprule
$\boldsymbol{S_{\chi}}$ & & $\boldsymbol{S_{\chi,0}}$ & $\boldsymbol{\HH^2(S_{\chi}, \bC^{\times}_{\pi})}$ & $\boldsymbol{S'_{\chi}}$ & & $\boldsymbol{\HH^2(S'_{\chi}, \bC^{\times}_{\kappa})}$ & $\boldsymbol{\chi}$ & $\boldsymbol{c_{\chi}}$ & $\boldsymbol{c'_{\chi}}$ & $\boldsymbol{\fT(V)}$ \\
\midrule
$1$ & $\langle\emptyset\rangle$ & $1$ & $0$ & $C_2$ & $\langle \overline{st}\rangle$ & $\mathbb Z/2\mathbb Z$ & $[\frac{1}{4},\frac{1}{4}]$ & $(0)$ & $(0)$ & III \\
\midrule
$C_2$ & $\langle s\rangle$ & $1$ & $\mathbb Z/2\mathbb Z$ & $C_2^2$ & $\langle t,\overline{s}\rangle$ & $(\mathbb Z/2\mathbb Z)^{2}$ & $[0,\frac{1}{4}]$ & $(0)$ & $(0,0)$ & I \\
\midrule
$C_2$ & $\langle t\rangle$ & $1$ & $\mathbb Z/2\mathbb Z$ & $C_2^2$ & $\langle s,\overline{t}\rangle$ & $(\mathbb Z/2\mathbb Z)^{2}$ & $[\frac{1}{4},0]$ & $(0)$ & $(0,0)$ & I \\
 &  &  &  &  &  &  & $[\frac{1}{4},\frac{1}{2}]$ & $(1)$ & $(1,1)$ & II \\
\midrule
$C_2^2$ & $\langle s,t\rangle$ & $C_2$ & $(\mathbb Z/2\mathbb Z)^{2}$ & $C_2^3$ & $\langle \overline{1},s,t\rangle$ & $(\mathbb Z/2\mathbb Z)^{4}$ & $[0,0]$ & $(0,0)$ & $(0,0,0,0)$ & I, I \\
 &  &  &  &  &  &  & $[0,\frac{1}{2}]$ & $(0,1)$ & $(0,0,1,1)$ & III \\
\bottomrule
\end{tabular}
\caption{$G = \text{pmg}$, $H = \langle st,x,y \rangle \cong \text{p2}$, so $S \cong C_2^2$ and $S_0 \cong C_2$}
\label{table:repsGpmgH1}
\end{table}

\begin{table}[!htbp]
\centering
\small
\begin{tabular}{llll|lll|llll}
\toprule
$\boldsymbol{S_{\chi}}$ & & $\boldsymbol{S_{\chi,0}}$ & $\boldsymbol{\HH^2(S_{\chi}, \bC^{\times}_{\pi})}$ & $\boldsymbol{S'_{\chi}}$ & & $\boldsymbol{\HH^2(S'_{\chi}, \bC^{\times}_{\kappa})}$ & $\boldsymbol{\chi}$ & $\boldsymbol{c_{\chi}}$ & $\boldsymbol{c'_{\chi}}$ & $\boldsymbol{\fT(V)}$ \\
\midrule
$C_2$ & $\langle st\rangle$ & $1$ & $\mathbb Z/2\mathbb Z$ & $C_2$ & $\langle \overline{st}\rangle$ & $\mathbb Z/2\mathbb Z$ & $[\frac{1}{4},\frac{1}{4}]$ & $(0)$ & $(0)$ & V \\
\cmidrule(l){5-11}
 &  &  &  & $C_2^2$ & $\langle t,\overline{s}\rangle$ & $(\mathbb Z/2\mathbb Z)^{2}$ & $[0,\frac{1}{4}]$ & $(0)$ & $(0,0)$ & I \\
\midrule
$C_2^2$ & $\langle s,t\rangle$ & $C_2$ & $(\mathbb Z/2\mathbb Z)^{2}$ & $C_2^2$ & $\langle s,\overline{t}\rangle$ & $(\mathbb Z/2\mathbb Z)^{2}$ & $[\frac{1}{4},0]$ & $(0,0)$ & $(0,0)$ & V, V \\
 &  &  &  &  &  &  & $[\frac{1}{4},\frac{1}{2}]$ & $(1,1)$ & $(1,1)$ & VII \\
\cmidrule(l){5-11}
 &  &  &  & $C_2^3$ & $\langle \overline{1},s,t\rangle$ & $(\mathbb Z/2\mathbb Z)^{4}$ & $[0,0]$ & $(0,0)$ & $(0,0,0,0)$ & I, I \\
 &  &  &  &  &  &  & $[0,\frac{1}{2}]$ & $(1,1)$ & $(0,0,1,1)$ & III \\
\bottomrule
\end{tabular}
\caption{$G = \text{pmg}$, $H = \langle s,x,y \rangle \cong \text{pm}$, so $S \cong C_2^2$ and $S_0 \cong C_2$}
\label{table:repsGpmgH2}
\end{table}

\begin{table}[!htbp]
\centering
\small
\begin{tabular}{llll|lll|llll}
\toprule
$\boldsymbol{S_{\chi}}$ & & $\boldsymbol{S_{\chi,0}}$ & $\boldsymbol{\HH^2(S_{\chi}, \bC^{\times}_{\pi})}$ & $\boldsymbol{S'_{\chi}}$ & & $\boldsymbol{\HH^2(S'_{\chi}, \bC^{\times}_{\kappa})}$ & $\boldsymbol{\chi}$ & $\boldsymbol{c_{\chi}}$ & $\boldsymbol{c'_{\chi}}$ & $\boldsymbol{\fT(V)}$ \\
\midrule
$C_2$ & $\langle st\rangle$ & $1$ & $\mathbb Z/2\mathbb Z$ & $C_2$ & $\langle \overline{st}\rangle$ & $\mathbb Z/2\mathbb Z$ & $[\frac{1}{4},\frac{1}{4}]$ & $(0)$ & $(0)$ & V \\
\cmidrule(l){5-11}
 &  &  &  & $C_2^2$ & $\langle s,\overline{t}\rangle$ & $(\mathbb Z/2\mathbb Z)^{2}$ & $[\frac{1}{4},0]$ & $(0)$ & $(0,0)$ & I \\
 &  &  &  &  &  &  & $[\frac{1}{4},\frac{1}{2}]$ & $(0)$ & $(1,1)$ & IX \\
\midrule
$C_2^2$ & $\langle s,t\rangle$ & $C_2$ & $(\mathbb Z/2\mathbb Z)^{2}$ & $C_2^2$ & $\langle t,\overline{s}\rangle$ & $(\mathbb Z/2\mathbb Z)^{2}$ & $[0,\frac{1}{4}]$ & $(0,0)$ & $(0,0)$ & V, V \\
\cmidrule(l){5-11}
 &  &  &  & $C_2^3$ & $\langle \overline{1},s,t\rangle$ & $(\mathbb Z/2\mathbb Z)^{4}$ & $[0,0]$ & $(0,0)$ & $(0,0,0,0)$ & I, I \\
 &  &  &  &  &  &  & $[0,\frac{1}{2}]$ & $(0,0)$ & $(0,0,1,1)$ & V, V \\
\bottomrule
\end{tabular}
\caption{$G = \text{pmg}$, $H = \langle t,x,y \rangle \cong \text{pg}$, so $S \cong C_2^2$ and $S_0 \cong C_2$}
\label{table:repsGpmgH3}
\end{table}

\begin{table}[!htbp]
\centering
\small
\begin{tabular}{llll|lll|llll}
\toprule
$\boldsymbol{S_{\chi}}$ & & $\boldsymbol{S_{\chi,0}}$ & $\boldsymbol{\HH^2(S_{\chi}, \bC^{\times}_{\pi})}$ & $\boldsymbol{S'_{\chi}}$ & & $\boldsymbol{\HH^2(S'_{\chi}, \bC^{\times}_{\kappa})}$ & $\boldsymbol{\chi}$ & $\boldsymbol{c_{\chi}}$ & $\boldsymbol{c'_{\chi}}$ & $\boldsymbol{\fT(V)}$ \\
\midrule
$C_2$ & $\langle stx\rangle$ & $1$ & $\mathbb Z/2\mathbb Z$ & $C_2^2$ & $\langle x,\overline{st}\rangle$ & $(\mathbb Z/2\mathbb Z)^{2}$ & $[\frac{1}{4},\frac{1}{4}]$ & $(0)$ & $(0,0)$ & I \\
\midrule
$C_2^2$ & $\langle t,sx\rangle$ & $C_2$ & $(\mathbb Z/2\mathbb Z)^{2}$ & $C_2^3$ & $\langle t,x,\overline{s}\rangle$ & $(\mathbb Z/2\mathbb Z)^{4}$ & $[0,\frac{1}{4}]$ & $(0,0)$ & $(0,0,0,0)$ & I, I \\
 &  &  &  &  &  &  & $[\frac{1}{2},\frac{1}{4}]$ & $(1,1)$ & $(0,0,1,1)$ & III \\
\midrule
$C_2^2$ & $\langle s,tx\rangle$ & $C_2$ & $(\mathbb Z/2\mathbb Z)^{2}$ & $C_2^3$ & $\langle s,x,\overline{t}\rangle$ & $(\mathbb Z/2\mathbb Z)^{4}$ & $[\frac{1}{4},0]$ & $(0,0)$ & $(0,0,0,0)$ & I, I \\
 &  &  &  &  &  &  & $[\frac{1}{4},\frac{1}{2}]$ & $(1,1)$ & $(1,1,0,0)$ & IV \\
\midrule
$C_2^3$ & $\langle s,t,x\rangle$ & $C_2^2$ & $(\mathbb Z/2\mathbb Z)^{4}$ & $C_2^4$ & $\langle \overline{1},s,t,x\rangle$ & $(\mathbb Z/2\mathbb Z)^{7}$ & $[0,0]$ & $(0,0,0,0)$ & $(0,0,0,0,0,0,0)$ & I, I, I, I \\
 &  &  &  &  &  &  & $[0,\frac{1}{2}]$ & $(0,1,0,1)$ & $(0,1,0,0,0,1,0)$ & I \\
 &  &  &  &  &  &  & $[\frac{1}{2},0]$ & $(0,1,1,0)$ & $(0,0,0,1,0,0,1)$ & III, III \\
 &  &  &  &  &  &  & $[\frac{1}{2},\frac{1}{2}]$ & $(0,0,1,1)$ & $(0,1,0,1,0,1,1)$ & II \\
\bottomrule
\end{tabular}
\caption{$G = \text{pmg}$, $H = \langle s,t,x^{2},y \rangle \cong \text{pmg}$, so $S \cong C_2^3$ and $S_0 \cong C_2^2$}
\label{table:repsGpmgH4}
\end{table}

\begin{table}[!htbp]
\centering
\small
\begin{tabular}{llll|lll|llll}
\toprule
$\boldsymbol{S_{\chi}}$ & & $\boldsymbol{S_{\chi,0}}$ & $\boldsymbol{\HH^2(S_{\chi}, \bC^{\times}_{\pi})}$ & $\boldsymbol{S'_{\chi}}$ & & $\boldsymbol{\HH^2(S'_{\chi}, \bC^{\times}_{\kappa})}$ & $\boldsymbol{\chi}$ & $\boldsymbol{c_{\chi}}$ & $\boldsymbol{c'_{\chi}}$ & $\boldsymbol{\fT(V)}$ \\
\midrule
$C_2$ & $\langle st\rangle$ & $1$ & $\mathbb Z/2\mathbb Z$ & $C_2^2$ & $\langle x,\overline{st}\rangle$ & $(\mathbb Z/2\mathbb Z)^{2}$ & $[\frac{1}{4},\frac{1}{4}]$ & $(0)$ & $(0,0)$ & I \\
\midrule
$C_2^2$ & $\langle s,t\rangle$ & $C_2$ & $(\mathbb Z/2\mathbb Z)^{2}$ & $C_2^3$ & $\langle t,x,\overline{s}\rangle$ & $(\mathbb Z/2\mathbb Z)^{4}$ & $[0,\frac{1}{4}]$ & $(0,0)$ & $(0,0,0,0)$ & I, I \\
 &  &  &  &  &  &  & $[\frac{1}{2},\frac{1}{4}]$ & $(0,0)$ & $(0,0,1,1)$ & V, V \\
\midrule
$C_2^2$ & $\langle st,sx\rangle$ & $C_2$ & $(\mathbb Z/2\mathbb Z)^{2}$ & $C_2^3$ & $\langle s,x,\overline{t}\rangle$ & $(\mathbb Z/2\mathbb Z)^{4}$ & $[\frac{1}{4},0]$ & $(0,0)$ & $(0,0,0,0)$ & I, I \\
 &  &  &  &  &  &  & $[\frac{1}{4},\frac{1}{2}]$ & $(0,1)$ & $(1,1,0,0)$ & IV \\
\midrule
$C_2^3$ & $\langle s,t,x\rangle$ & $C_2^2$ & $(\mathbb Z/2\mathbb Z)^{4}$ & $C_2^4$ & $\langle \overline{1},s,t,x\rangle$ & $(\mathbb Z/2\mathbb Z)^{7}$ & $[0,0]$ & $(0,0,0,0)$ & $(0,0,0,0,0,0,0)$ & I, I, I, I \\
 &  &  &  &  &  &  & $[0,\frac{1}{2}]$ & $(0,1,0,0)$ & $(0,1,0,0,0,1,0)$ & I \\
 &  &  &  &  &  &  & $[\frac{1}{2},0]$ & $(1,1,0,0)$ & $(0,0,0,1,0,0,1)$ & I \\
 &  &  &  &  &  &  & $[\frac{1}{2},\frac{1}{2}]$ & $(1,0,0,0)$ & $(0,1,0,1,0,1,1)$ & VII, VII \\
\bottomrule
\end{tabular}
\caption{$G = \text{pmg}$, $H = \langle sx,t,x^{2},y \rangle \cong \text{pgg}$, so $S \cong C_2^3$ and $S_0 \cong C_2^2$}
\label{table:repsGpmgH5}
\end{table}

\begin{table}[!htbp]
\centering
\small
\begin{tabular}{llll|lll|llll}
\toprule
$\boldsymbol{S_{\chi}}$ & & $\boldsymbol{S_{\chi,0}}$ & $\boldsymbol{\HH^2(S_{\chi}, \bC^{\times}_{\pi})}$ & $\boldsymbol{S'_{\chi}}$ & & $\boldsymbol{\HH^2(S'_{\chi}, \bC^{\times}_{\kappa})}$ & $\boldsymbol{\chi}$ & $\boldsymbol{c_{\chi}}$ & $\boldsymbol{c'_{\chi}}$ & $\boldsymbol{\fT(V)}$ \\
\midrule
$1$ & $\langle\emptyset\rangle$ & $1$ & $0$ & $C_2$ & $\langle \overline{st}\rangle$ & $\mathbb Z/2\mathbb Z$ & $[\frac{1}{4},\frac{1}{4}]$ & $(0)$ & $(0)$ & III \\
\midrule
$C_2$ & $\langle s\rangle$ & $1$ & $\mathbb Z/2\mathbb Z$ & $C_2^2$ & $\langle t,\overline{s}\rangle$ & $(\mathbb Z/2\mathbb Z)^{2}$ & $[0,\frac{1}{4}]$ & $(0)$ & $(0,0)$ & I \\
 &  &  &  &  &  &  & $[\frac{1}{2},\frac{1}{4}]$ & $(1)$ & $(1,1)$ & II \\
\midrule
$C_2$ & $\langle t\rangle$ & $1$ & $\mathbb Z/2\mathbb Z$ & $C_2^2$ & $\langle s,\overline{t}\rangle$ & $(\mathbb Z/2\mathbb Z)^{2}$ & $[\frac{1}{4},0]$ & $(0)$ & $(0,0)$ & I \\
 &  &  &  &  &  &  & $[\frac{1}{4},\frac{1}{2}]$ & $(1)$ & $(1,1)$ & II \\
\midrule
$C_2^2$ & $\langle s,t\rangle$ & $C_2$ & $(\mathbb Z/2\mathbb Z)^{2}$ & $C_2^3$ & $\langle \overline{1},s,t\rangle$ & $(\mathbb Z/2\mathbb Z)^{4}$ & $[0,0]$ & $(0,0)$ & $(0,0,0,0)$ & I, I \\
 &  &  &  &  &  &  & $[0,\frac{1}{2}]$ & $(0,1)$ & $(0,0,1,1)$ & III \\
 &  &  &  &  &  &  & $[\frac{1}{2},0]$ & $(1,0)$ & $(0,1,0,1)$ & III \\
 &  &  &  &  &  &  & $[\frac{1}{2},\frac{1}{2}]$ & $(1,1)$ & $(0,1,1,0)$ & II, II \\
\bottomrule
\end{tabular}
\caption{$G = \text{pgg}$, $H = \langle st,x,y \rangle \cong \text{p2}$, so $S \cong C_2^2$ and $S_0 \cong C_2$}
\label{table:repsGpggH1}
\end{table}

\begin{table}[!htbp]
\centering
\small
\begin{tabular}{llll|lll|llll}
\toprule
$\boldsymbol{S_{\chi}}$ & & $\boldsymbol{S_{\chi,0}}$ & $\boldsymbol{\HH^2(S_{\chi}, \bC^{\times}_{\pi})}$ & $\boldsymbol{S'_{\chi}}$ & & $\boldsymbol{\HH^2(S'_{\chi}, \bC^{\times}_{\kappa})}$ & $\boldsymbol{\chi}$ & $\boldsymbol{c_{\chi}}$ & $\boldsymbol{c'_{\chi}}$ & $\boldsymbol{\fT(V)}$ \\
\midrule
$C_2$ & $\langle st\rangle$ & $1$ & $\mathbb Z/2\mathbb Z$ & $C_2$ & $\langle \overline{st}\rangle$ & $\mathbb Z/2\mathbb Z$ & $[\frac{1}{4},\frac{1}{4}]$ & $(0)$ & $(0)$ & V \\
\cmidrule(l){5-11}
 &  &  &  & $C_2^2$ & $\langle t,\overline{s}\rangle$ & $(\mathbb Z/2\mathbb Z)^{2}$ & $[0,\frac{1}{4}]$ & $(0)$ & $(0,0)$ & I \\
 &  &  &  &  &  &  & $[\frac{1}{2},\frac{1}{4}]$ & $(0)$ & $(1,1)$ & IX \\
\midrule
$C_2^2$ & $\langle s,t\rangle$ & $C_2$ & $(\mathbb Z/2\mathbb Z)^{2}$ & $C_2^2$ & $\langle s,\overline{t}\rangle$ & $(\mathbb Z/2\mathbb Z)^{2}$ & $[\frac{1}{4},0]$ & $(0,0)$ & $(0,0)$ & V, V \\
 &  &  &  &  &  &  & $[\frac{1}{4},\frac{1}{2}]$ & $(1,1)$ & $(1,1)$ & VII \\
\cmidrule(l){5-11}
 &  &  &  & $C_2^3$ & $\langle \overline{1},s,t\rangle$ & $(\mathbb Z/2\mathbb Z)^{4}$ & $[0,0]$ & $(0,0)$ & $(0,0,0,0)$ & I, I \\
 &  &  &  &  &  &  & $[0,\frac{1}{2}]$ & $(1,1)$ & $(0,0,1,1)$ & III \\
 &  &  &  &  &  &  & $[\frac{1}{2},0]$ & $(0,0)$ & $(0,1,0,1)$ & V, V \\
 &  &  &  &  &  &  & $[\frac{1}{2},\frac{1}{2}]$ & $(1,1)$ & $(0,1,1,0)$ & IV \\
\bottomrule
\end{tabular}
\caption{$G = \text{pgg}$, $H = \langle s,x,y \rangle \cong \text{pg}$, so $S \cong C_2^2$ and $S_0 \cong C_2$}
\label{table:repsGpggH2}
\end{table}

\begin{table}[!htbp]
\centering
\small
\begin{tabular}{llll|lll|llll}
\toprule
$\boldsymbol{S_{\chi}}$ & & $\boldsymbol{S_{\chi,0}}$ & $\boldsymbol{\HH^2(S_{\chi}, \bC^{\times}_{\pi})}$ & $\boldsymbol{S'_{\chi}}$ & & $\boldsymbol{\HH^2(S'_{\chi}, \bC^{\times}_{\kappa})}$ & $\boldsymbol{\chi}$ & $\boldsymbol{c_{\chi}}$ & $\boldsymbol{c'_{\chi}}$ & $\boldsymbol{\fT(V)}$ \\
\midrule
$1$ & $\langle\emptyset\rangle$ & $1$ & $0$ & $C_2$ & $\langle \overline{st}\rangle$ & $\mathbb Z/2\mathbb Z$ & $[0,\frac{1}{4}]$ & $(0)$ & $(0)$ & III \\
\midrule
$C_2$ & $\langle t\rangle$ & $1$ & $\mathbb Z/2\mathbb Z$ & $C_2^2$ & $\langle s,\overline{t}\rangle$ & $(\mathbb Z/2\mathbb Z)^{2}$ & $[\frac{1}{4},\frac{1}{4}]$ & $(0)$ & $(0,0)$ & I \\
\midrule
$C_2$ & $\langle s\rangle$ & $1$ & $\mathbb Z/2\mathbb Z$ & $C_2^2$ & $\langle t,\overline{s}\rangle$ & $(\mathbb Z/2\mathbb Z)^{2}$ & $[\frac{1}{4},\frac{3}{4}]$ & $(0)$ & $(0,0)$ & I \\
\midrule
$C_2$ & $\langle st\rangle$ & $C_2$ & $0$ & $C_2^2$ & $\langle \overline{1},st\rangle$ & $(\mathbb Z/2\mathbb Z)^{2}$ & $[0,\frac{1}{2}]$ & $(0)$ & $(0,0)$ & III, III \\
\midrule
$C_2^2$ & $\langle s,t\rangle$ & $C_2$ & $(\mathbb Z/2\mathbb Z)^{2}$ & $C_2^3$ & $\langle \overline{1},s,t\rangle$ & $(\mathbb Z/2\mathbb Z)^{4}$ & $[0,0]$ & $(0,0)$ & $(0,0,0,0)$ & I, I \\
\bottomrule
\end{tabular}
\caption{$G = \text{cmm}$, $H = \langle st,x,y \rangle \cong \text{p2}$, so $S \cong C_2^2$ and $S_0 \cong C_2$}
\label{table:repsGcmmH1}
\end{table}

\begin{table}[!htbp]
\centering
\small
\begin{tabular}{llll|lll|llll}
\toprule
$\boldsymbol{S_{\chi}}$ & & $\boldsymbol{S_{\chi,0}}$ & $\boldsymbol{\HH^2(S_{\chi}, \bC^{\times}_{\pi})}$ & $\boldsymbol{S'_{\chi}}$ & & $\boldsymbol{\HH^2(S'_{\chi}, \bC^{\times}_{\kappa})}$ & $\boldsymbol{\chi}$ & $\boldsymbol{c_{\chi}}$ & $\boldsymbol{c'_{\chi}}$ & $\boldsymbol{\fT(V)}$ \\
\midrule
$C_2$ & $\langle st\rangle$ & $1$ & $\mathbb Z/2\mathbb Z$ & $C_2$ & $\langle \overline{st}\rangle$ & $\mathbb Z/2\mathbb Z$ & $[0,\frac{1}{4}]$ & $(0)$ & $(0)$ & V \\
\cmidrule(l){5-11}
 &  &  &  & $C_2^2$ & $\langle t,\overline{s}\rangle$ & $(\mathbb Z/2\mathbb Z)^{2}$ & $[\frac{1}{4},\frac{3}{4}]$ & $(0)$ & $(0,0)$ & I \\
\cmidrule(l){5-11}
 &  &  &  & $C_2^2$ & $\langle \overline{1},st\rangle$ & $(\mathbb Z/2\mathbb Z)^{2}$ & $[0,\frac{1}{2}]$ & $(0)$ & $(0,0)$ & I \\
\midrule
$C_2^2$ & $\langle s,t\rangle$ & $C_2$ & $(\mathbb Z/2\mathbb Z)^{2}$ & $C_2^2$ & $\langle s,\overline{t}\rangle$ & $(\mathbb Z/2\mathbb Z)^{2}$ & $[\frac{1}{4},\frac{1}{4}]$ & $(0,0)$ & $(0,0)$ & V, V \\
\cmidrule(l){5-11}
 &  &  &  & $C_2^3$ & $\langle \overline{1},s,t\rangle$ & $(\mathbb Z/2\mathbb Z)^{4}$ & $[0,0]$ & $(0,0)$ & $(0,0,0,0)$ & I, I \\
\bottomrule
\end{tabular}
\caption{$G = \text{cmm}$, $H = \langle s,x,y \rangle \cong \text{cm}$, so $S \cong C_2^2$ and $S_0 \cong C_2$}
\label{table:repsGcmmH2}
\end{table}

\begin{table}[!htbp]
\centering
\small
\begin{tabular}{llll|lll|llll}
\toprule
$\boldsymbol{S_{\chi}}$ & & $\boldsymbol{S_{\chi,0}}$ & $\boldsymbol{\HH^2(S_{\chi}, \bC^{\times}_{\pi})}$ & $\boldsymbol{S'_{\chi}}$ & & $\boldsymbol{\HH^2(S'_{\chi}, \bC^{\times}_{\kappa})}$ & $\boldsymbol{\chi}$ & $\boldsymbol{c_{\chi}}$ & $\boldsymbol{c'_{\chi}}$ & $\boldsymbol{\fT(V)}$ \\
\midrule
$C_2$ & $\langle st\rangle$ & $1$ & $\mathbb Z/2\mathbb Z$ & $C_2^2$ & $\langle x,\overline{st}\rangle$ & $(\mathbb Z/2\mathbb Z)^{2}$ & $[\frac{1}{4},\frac{1}{4}]$ & $(0)$ & $(0,0)$ & I \\
\midrule
$C_2^2$ & $\langle st,sx\rangle$ & $C_2$ & $(\mathbb Z/2\mathbb Z)^{2}$ & $C_2^3$ & $\langle t,x,\overline{s}\rangle$ & $(\mathbb Z/2\mathbb Z)^{4}$ & $[0,\frac{1}{4}]$ & $(0,0)$ & $(0,0,0,0)$ & I, I \\
 &  &  &  &  &  &  & $[\frac{1}{2},\frac{1}{4}]$ & $(0,1)$ & $(0,0,1,1)$ & III \\
\midrule
$C_2^2$ & $\langle s,t\rangle$ & $C_2$ & $(\mathbb Z/2\mathbb Z)^{2}$ & $C_2^3$ & $\langle s,x,\overline{t}\rangle$ & $(\mathbb Z/2\mathbb Z)^{4}$ & $[\frac{1}{4},0]$ & $(0,0)$ & $(0,0,0,0)$ & I, I \\
 &  &  &  &  &  &  & $[\frac{1}{4},\frac{1}{2}]$ & $(0,0)$ & $(0,0,1,1)$ & V, V \\
\midrule
$C_2^3$ & $\langle s,t,x\rangle$ & $C_2^2$ & $(\mathbb Z/2\mathbb Z)^{4}$ & $C_2^4$ & $\langle \overline{1},s,t,x\rangle$ & $(\mathbb Z/2\mathbb Z)^{7}$ & $[0,0]$ & $(0,0,0,0)$ & $(0,0,0,0,0,0,0)$ & I, I, I, I \\
 &  &  &  &  &  &  & $[0,\frac{1}{2}]$ & $(1,1,0,0)$ & $(0,0,1,0,0,0,1)$ & I \\
 &  &  &  &  &  &  & $[\frac{1}{2},0]$ & $(1,0,0,0)$ & $(0,0,0,1,0,0,1)$ & III, III \\
 &  &  &  &  &  &  & $[\frac{1}{2},\frac{1}{2}]$ & $(0,1,0,0)$ & $(0,0,1,1,0,0,0)$ & I \\
\bottomrule
\end{tabular}
\caption{$G = \text{cmm}$, $H = \langle s,tx,xy,xy^{-1} \rangle \cong \text{pmg}$, so $S \cong C_2^3$ and $S_0 \cong C_2^2$}
\label{table:repsGcmmH3}
\end{table}

\begin{table}[!htbp]
\centering
\small
\begin{tabular}{llll|lll|llll}
\toprule
$\boldsymbol{S_{\chi}}$ & & $\boldsymbol{S_{\chi,0}}$ & $\boldsymbol{\HH^2(S_{\chi}, \bC^{\times}_{\pi})}$ & $\boldsymbol{S'_{\chi}}$ & & $\boldsymbol{\HH^2(S'_{\chi}, \bC^{\times}_{\kappa})}$ & $\boldsymbol{\chi}$ & $\boldsymbol{c_{\chi}}$ & $\boldsymbol{c'_{\chi}}$ & $\boldsymbol{\fT(V)}$ \\
\midrule
$C_2$ & $\langle stx\rangle$ & $1$ & $\mathbb Z/2\mathbb Z$ & $C_2^2$ & $\langle x,\overline{st}\rangle$ & $(\mathbb Z/2\mathbb Z)^{2}$ & $[\frac{1}{4},\frac{1}{4}]$ & $(0)$ & $(0,0)$ & I \\
\midrule
$C_2^2$ & $\langle t,sx\rangle$ & $C_2$ & $(\mathbb Z/2\mathbb Z)^{2}$ & $C_2^3$ & $\langle t,x,\overline{s}\rangle$ & $(\mathbb Z/2\mathbb Z)^{4}$ & $[0,\frac{1}{4}]$ & $(0,0)$ & $(0,0,0,0)$ & I, I \\
 &  &  &  &  &  &  & $[\frac{1}{2},\frac{1}{4}]$ & $(1,1)$ & $(0,0,1,1)$ & III \\
\midrule
$C_2^2$ & $\langle s,tx\rangle$ & $C_2$ & $(\mathbb Z/2\mathbb Z)^{2}$ & $C_2^3$ & $\langle s,x,\overline{t}\rangle$ & $(\mathbb Z/2\mathbb Z)^{4}$ & $[\frac{1}{4},0]$ & $(0,0)$ & $(0,0,0,0)$ & I, I \\
 &  &  &  &  &  &  & $[\frac{1}{4},\frac{1}{2}]$ & $(1,1)$ & $(0,0,1,1)$ & III \\
\midrule
$C_2^3$ & $\langle s,t,x\rangle$ & $C_2^2$ & $(\mathbb Z/2\mathbb Z)^{4}$ & $C_2^4$ & $\langle \overline{1},s,t,x\rangle$ & $(\mathbb Z/2\mathbb Z)^{7}$ & $[0,0]$ & $(0,0,0,0)$ & $(0,0,0,0,0,0,0)$ & I, I, I, I \\
 &  &  &  &  &  &  & $[0,\frac{1}{2}]$ & $(1,0,1,0)$ & $(0,0,1,0,0,0,1)$ & III, III \\
 &  &  &  &  &  &  & $[\frac{1}{2},0]$ & $(0,1,1,0)$ & $(0,0,0,1,0,0,1)$ & III, III \\
 &  &  &  &  &  &  & $[\frac{1}{2},\frac{1}{2}]$ & $(1,1,0,0)$ & $(0,0,1,1,0,0,0)$ & III, III \\
\bottomrule
\end{tabular}
\caption{$G = \text{cmm}$, $H = \langle s,t,xy,xy^{-1} \rangle \cong \text{pmm}$, so $S \cong C_2^3$ and $S_0 \cong C_2^2$}
\label{table:repsGcmmH4}
\end{table}

\begin{table}[!htbp]
\centering
\small
\begin{tabular}{llll|lll|llll}
\toprule
$\boldsymbol{S_{\chi}}$ & & $\boldsymbol{S_{\chi,0}}$ & $\boldsymbol{\HH^2(S_{\chi}, \bC^{\times}_{\pi})}$ & $\boldsymbol{S'_{\chi}}$ & & $\boldsymbol{\HH^2(S'_{\chi}, \bC^{\times}_{\kappa})}$ & $\boldsymbol{\chi}$ & $\boldsymbol{c_{\chi}}$ & $\boldsymbol{c'_{\chi}}$ & $\boldsymbol{\fT(V)}$ \\
\midrule
$C_2$ & $\langle stx\rangle$ & $1$ & $\mathbb Z/2\mathbb Z$ & $C_2^2$ & $\langle x,\overline{st}\rangle$ & $(\mathbb Z/2\mathbb Z)^{2}$ & $[\frac{1}{4},\frac{1}{4}]$ & $(0)$ & $(0,0)$ & I \\
\midrule
$C_2^2$ & $\langle s,tx\rangle$ & $C_2$ & $(\mathbb Z/2\mathbb Z)^{2}$ & $C_2^3$ & $\langle t,x,\overline{s}\rangle$ & $(\mathbb Z/2\mathbb Z)^{4}$ & $[0,\frac{1}{4}]$ & $(0,0)$ & $(0,0,0,0)$ & I, I \\
 &  &  &  &  &  &  & $[\frac{1}{2},\frac{1}{4}]$ & $(0,0)$ & $(0,0,1,1)$ & V, V \\
\midrule
$C_2^2$ & $\langle t,sx\rangle$ & $C_2$ & $(\mathbb Z/2\mathbb Z)^{2}$ & $C_2^3$ & $\langle s,x,\overline{t}\rangle$ & $(\mathbb Z/2\mathbb Z)^{4}$ & $[\frac{1}{4},0]$ & $(0,0)$ & $(0,0,0,0)$ & I, I \\
 &  &  &  &  &  &  & $[\frac{1}{4},\frac{1}{2}]$ & $(0,0)$ & $(0,0,1,1)$ & V, V \\
\midrule
$C_2^3$ & $\langle s,t,x\rangle$ & $C_2^2$ & $(\mathbb Z/2\mathbb Z)^{4}$ & $C_2^4$ & $\langle \overline{1},s,t,x\rangle$ & $(\mathbb Z/2\mathbb Z)^{7}$ & $[0,0]$ & $(0,0,0,0)$ & $(0,0,0,0,0,0,0)$ & I, I, I, I \\
 &  &  &  &  &  &  & $[0,\frac{1}{2}]$ & $(1,0,0,1)$ & $(0,0,1,0,0,0,1)$ & I \\
 &  &  &  &  &  &  & $[\frac{1}{2},0]$ & $(1,0,0,1)$ & $(0,0,0,1,0,0,1)$ & I \\
 &  &  &  &  &  &  & $[\frac{1}{2},\frac{1}{2}]$ & $(0,0,0,0)$ & $(0,0,1,1,0,0,0)$ & V, V, V, V \\
\bottomrule
\end{tabular}
\caption{$G = \text{cmm}$, $H = \langle sx,tx,xy,xy^{-1} \rangle \cong \text{pgg}$, so $S \cong C_2^3$ and $S_0 \cong C_2^2$}
\label{table:repsGcmmH5}
\end{table}

\begin{table}[!htbp]
\centering
\small
\begin{tabular}{llll|lll|llll}
\toprule
$\boldsymbol{S_{\chi}}$ & & $\boldsymbol{S_{\chi,0}}$ & $\boldsymbol{\HH^2(S_{\chi}, \bC^{\times}_{\pi})}$ & $\boldsymbol{S'_{\chi}}$ & & $\boldsymbol{\HH^2(S'_{\chi}, \bC^{\times}_{\kappa})}$ & $\boldsymbol{\chi}$ & $\boldsymbol{c_{\chi}}$ & $\boldsymbol{c'_{\chi}}$ & $\boldsymbol{\fT(V)}$ \\
\midrule
$1$ & $\langle\emptyset\rangle$ & $1$ & $0$ & $C_2$ & $\langle \overline{r^{2}}\rangle$ & $\mathbb Z/2\mathbb Z$ & $[0,\frac{1}{4}]$ & $(0)$ & $(0)$ & III \\
\midrule
$C_2$ & $\langle r^{2}\rangle$ & $C_2$ & $0$ & $C_2^2$ & $\langle \overline{1},r^{2}\rangle$ & $(\mathbb Z/2\mathbb Z)^{2}$ & $[0,\frac{1}{2}]$ & $(0)$ & $(0,0)$ & III, III \\
\midrule
$C_4$ & $\langle r^{-1}\rangle$ & $C_2$ & $\mathbb Z/2\mathbb Z$ & $C_2 \times C_4$ & $\langle \overline{1},r^{-1}\rangle$ & $(\mathbb Z/2\mathbb Z)^{2}$ & $[0,0]$ & $(0)$ & $(0,0)$ & I, II \\
\bottomrule
\end{tabular}
\caption{$G = \text{p4}$, $H = \langle r^{2},u,v \rangle \cong \text{p2}$, so $S \cong C_4$ and $S_0 \cong C_2$}
\label{table:repsGp4H1}
\end{table}

\begin{table}[!htbp]
\centering
\small
\begin{tabular}{llll|lll|llll}
\toprule
$\boldsymbol{S_{\chi}}$ & & $\boldsymbol{S_{\chi,0}}$ & $\boldsymbol{\HH^2(S_{\chi}, \bC^{\times}_{\pi})}$ & $\boldsymbol{S'_{\chi}}$ & & $\boldsymbol{\HH^2(S'_{\chi}, \bC^{\times}_{\kappa})}$ & $\boldsymbol{\chi}$ & $\boldsymbol{c_{\chi}}$ & $\boldsymbol{c'_{\chi}}$ & $\boldsymbol{\fT(V)}$ \\
\midrule
$C_2$ & $\langle r^{2}u\rangle$ & $1$ & $\mathbb Z/2\mathbb Z$ & $C_2^2$ & $\langle u,\overline{r^{2}}\rangle$ & $(\mathbb Z/2\mathbb Z)^{2}$ & $[0,\frac{1}{4}]$ & $(0)$ & $(0,0)$ & I \\
\midrule
$C_2^2$ & $\langle u,r^{2}\rangle$ & $C_2$ & $(\mathbb Z/2\mathbb Z)^{2}$ & $C_2^3$ & $\langle \overline{1},u,r^{2}\rangle$ & $(\mathbb Z/2\mathbb Z)^{4}$ & $[0,\frac{1}{2}]$ & $(1,1)$ & $(0,1,0,1)$ & III \\
\midrule
$C_2 \times C_4$ & $\langle r^{-1},u\rangle$ & $C_4$ & $(\mathbb Z/2\mathbb Z)^{2}$ & $C_2^2 \times C_4$ & $\langle \overline{1},r^{-1},u\rangle$ & $(\mathbb Z/2\mathbb Z)^{4}$ & $[0,0]$ & $(0,0)$ & $(0,0,0,0)$ & I, I, IV \\
 &  &  &  &  &  &  & $[\frac{1}{2},\frac{1}{2}]$ & $(0,0)$ & $(0,0,0,1)$ & III, V, V \\
\bottomrule
\end{tabular}
\caption{$G = \text{p4}$, $H = \langle r,uv,uv^{-1} \rangle \cong \text{p4}$, so $S \cong C_2 \times C_4$ and $S_0 \cong C_4$}
\label{table:repsGp4H2}
\end{table}

\begin{table}[!htbp]
\centering
\small
\begin{tabular}{llll|lll|llll}
\toprule
$\boldsymbol{S_{\chi}}$ & & $\boldsymbol{S_{\chi,0}}$ & $\boldsymbol{\HH^2(S_{\chi}, \bC^{\times}_{\pi})}$ & $\boldsymbol{S'_{\chi}}$ & & $\boldsymbol{\HH^2(S'_{\chi}, \bC^{\times}_{\kappa})}$ & $\boldsymbol{\chi}$ & $\boldsymbol{c_{\chi}}$ & $\boldsymbol{c'_{\chi}}$ & $\boldsymbol{\fT(V)}$ \\
\midrule
$1$ & $\langle\emptyset\rangle$ & $1$ & $0$ & $C_2$ & $\langle \overline{t^{2}}\rangle$ & $\mathbb Z/2\mathbb Z$ & $[\frac{1}{8},\frac{1}{4}]$ & $(0)$ & $(0)$ & III \\
\midrule
$C_2$ & $\langle st\rangle$ & $1$ & $\mathbb Z/2\mathbb Z$ & $C_2^2$ & $\langle st^{-1},\overline{st}\rangle$ & $(\mathbb Z/2\mathbb Z)^{2}$ & $[\frac{1}{4},\frac{1}{4}]$ & $(0)$ & $(0,0)$ & I \\
\midrule
$C_2$ & $\langle s\rangle$ & $1$ & $\mathbb Z/2\mathbb Z$ & $C_2^2$ & $\langle \overline{s},\overline{t^{2}}\rangle$ & $(\mathbb Z/2\mathbb Z)^{2}$ & $[0,\frac{1}{4}]$ & $(0)$ & $(0,0)$ & I \\
\midrule
$C_2$ & $\langle st^{-1}\rangle$ & $1$ & $\mathbb Z/2\mathbb Z$ & $C_2^2$ & $\langle st,\overline{st^{-1}}\rangle$ & $(\mathbb Z/2\mathbb Z)^{2}$ & $[\frac{1}{4},\frac{3}{4}]$ & $(0)$ & $(0,0)$ & I \\
\midrule
$C_2$ & $\langle st^{2}\rangle$ & $1$ & $\mathbb Z/2\mathbb Z$ & $C_2^2$ & $\langle s,\overline{t^{2}}\rangle$ & $(\mathbb Z/2\mathbb Z)^{2}$ & $[\frac{1}{4},0]$ & $(0)$ & $(0,0)$ & I \\
\midrule
$C_2^2$ & $\langle s,t^{2}\rangle$ & $C_2$ & $(\mathbb Z/2\mathbb Z)^{2}$ & $C_2^3$ & $\langle \overline{1},s,t^{2}\rangle$ & $(\mathbb Z/2\mathbb Z)^{4}$ & $[0,\frac{1}{2}]$ & $(0,0)$ & $(0,0,0,0)$ & I, I \\
\midrule
$D_4$ & $\langle t^{-1},s\rangle$ & $C_4$ & $(\mathbb Z/2\mathbb Z)^{2}$ & $C_2 \times D_4$ & $\langle \overline{1},t^{-1},s\rangle$ & $(\mathbb Z/2\mathbb Z)^{4}$ & $[0,0]$ & $(0,0)$ & $(0,0,0,0)$ & I, I, V, V \\
\bottomrule
\end{tabular}
\caption{$G = \text{p4m}$, $H = \langle t,x,y \rangle \cong \text{p4}$, so $S \cong D_4$ and $S_0 \cong C_4$}
\label{table:repsGp4mH1}
\end{table}

\begin{table}[!htbp]
\centering
\small
\begin{tabular}{llll|lll|llll}
\toprule
$\boldsymbol{S_{\chi}}$ & & $\boldsymbol{S_{\chi,0}}$ & $\boldsymbol{\HH^2(S_{\chi}, \bC^{\times}_{\pi})}$ & $\boldsymbol{S'_{\chi}}$ & & $\boldsymbol{\HH^2(S'_{\chi}, \bC^{\times}_{\kappa})}$ & $\boldsymbol{\chi}$ & $\boldsymbol{c_{\chi}}$ & $\boldsymbol{c'_{\chi}}$ & $\boldsymbol{\fT(V)}$ \\
\midrule
$1$ & $\langle\emptyset\rangle$ & $1$ & $0$ & $C_2$ & $\langle \overline{t^{2}}\rangle$ & $\mathbb Z/2\mathbb Z$ & $[\frac{1}{8},\frac{1}{4}]$ & $(0)$ & $(0)$ & III \\
\midrule
$C_2$ & $\langle st\rangle$ & $1$ & $\mathbb Z/2\mathbb Z$ & $C_2^2$ & $\langle st^{-1},\overline{st}\rangle$ & $(\mathbb Z/2\mathbb Z)^{2}$ & $[\frac{1}{4},\frac{1}{4}]$ & $(0)$ & $(0,0)$ & I \\
\midrule
$C_2$ & $\langle st^{2}\rangle$ & $C_2$ & $0$ & $C_2^2$ & $\langle \overline{s},\overline{t^{2}}\rangle$ & $(\mathbb Z/2\mathbb Z)^{2}$ & $[0,\frac{1}{4}]$ & $(0)$ & $(0,0)$ & III, III \\
\midrule
$C_2$ & $\langle st^{-1}\rangle$ & $1$ & $\mathbb Z/2\mathbb Z$ & $C_2^2$ & $\langle st,\overline{st^{-1}}\rangle$ & $(\mathbb Z/2\mathbb Z)^{2}$ & $[\frac{1}{4},\frac{3}{4}]$ & $(0)$ & $(0,0)$ & I \\
\midrule
$C_2$ & $\langle s\rangle$ & $C_2$ & $0$ & $C_2^2$ & $\langle s,\overline{t^{2}}\rangle$ & $(\mathbb Z/2\mathbb Z)^{2}$ & $[\frac{1}{4},0]$ & $(0)$ & $(0,0)$ & III, III \\
\midrule
$C_2^2$ & $\langle s,t^{2}\rangle$ & $C_2^2$ & $\mathbb Z/2\mathbb Z$ & $C_2^3$ & $\langle \overline{1},s,t^{2}\rangle$ & $(\mathbb Z/2\mathbb Z)^{4}$ & $[0,\frac{1}{2}]$ & $(0)$ & $(0,0,0,0)$ & III, III, III, III \\
\midrule
$D_4$ & $\langle t^{-1},s\rangle$ & $C_2^2$ & $(\mathbb Z/2\mathbb Z)^{2}$ & $C_2 \times D_4$ & $\langle \overline{1},t^{-1},s\rangle$ & $(\mathbb Z/2\mathbb Z)^{4}$ & $[0,0]$ & $(0,0)$ & $(0,0,0,0)$ & I, I, III \\
\bottomrule
\end{tabular}
\caption{$G = \text{p4m}$, $H = \langle s,t^{2},x,y \rangle \cong \text{pmm}$, so $S \cong D_4$ and $S_0 \cong C_2^2$}
\label{table:repsGp4mH2}
\end{table}

\begin{table}[!htbp]
\centering
\small
\begin{tabular}{llll|lll|llll}
\toprule
$\boldsymbol{S_{\chi}}$ & & $\boldsymbol{S_{\chi,0}}$ & $\boldsymbol{\HH^2(S_{\chi}, \bC^{\times}_{\pi})}$ & $\boldsymbol{S'_{\chi}}$ & & $\boldsymbol{\HH^2(S'_{\chi}, \bC^{\times}_{\kappa})}$ & $\boldsymbol{\chi}$ & $\boldsymbol{c_{\chi}}$ & $\boldsymbol{c'_{\chi}}$ & $\boldsymbol{\fT(V)}$ \\
\midrule
$1$ & $\langle\emptyset\rangle$ & $1$ & $0$ & $C_2$ & $\langle \overline{t^{2}}\rangle$ & $\mathbb Z/2\mathbb Z$ & $[\frac{1}{8},\frac{1}{4}]$ & $(0)$ & $(0)$ & III \\
\midrule
$C_2$ & $\langle st^{-1}\rangle$ & $C_2$ & $0$ & $C_2^2$ & $\langle st^{-1},\overline{st}\rangle$ & $(\mathbb Z/2\mathbb Z)^{2}$ & $[\frac{1}{4},\frac{1}{4}]$ & $(0)$ & $(0,0)$ & III, III \\
\midrule
$C_2$ & $\langle s\rangle$ & $1$ & $\mathbb Z/2\mathbb Z$ & $C_2^2$ & $\langle \overline{s},\overline{t^{2}}\rangle$ & $(\mathbb Z/2\mathbb Z)^{2}$ & $[0,\frac{1}{4}]$ & $(0)$ & $(0,0)$ & I \\
\midrule
$C_2$ & $\langle st\rangle$ & $C_2$ & $0$ & $C_2^2$ & $\langle st,\overline{st^{-1}}\rangle$ & $(\mathbb Z/2\mathbb Z)^{2}$ & $[\frac{1}{4},\frac{3}{4}]$ & $(0)$ & $(0,0)$ & III, III \\
\midrule
$C_2$ & $\langle st^{2}\rangle$ & $1$ & $\mathbb Z/2\mathbb Z$ & $C_2^2$ & $\langle s,\overline{t^{2}}\rangle$ & $(\mathbb Z/2\mathbb Z)^{2}$ & $[\frac{1}{4},0]$ & $(0)$ & $(0,0)$ & I \\
\midrule
$C_2^2$ & $\langle s,t^{2}\rangle$ & $C_2$ & $(\mathbb Z/2\mathbb Z)^{2}$ & $C_2^3$ & $\langle \overline{1},s,t^{2}\rangle$ & $(\mathbb Z/2\mathbb Z)^{4}$ & $[0,\frac{1}{2}]$ & $(0,0)$ & $(0,0,0,0)$ & I, I \\
\midrule
$D_4$ & $\langle t^{-1},s\rangle$ & $C_2^2$ & $(\mathbb Z/2\mathbb Z)^{2}$ & $C_2 \times D_4$ & $\langle \overline{1},t^{-1},s\rangle$ & $(\mathbb Z/2\mathbb Z)^{4}$ & $[0,0]$ & $(0,0)$ & $(0,0,0,0)$ & I, I, III \\
\bottomrule
\end{tabular}
\caption{$G = \text{p4m}$, $H = \langle st,st^{-1},x,y \rangle \cong \text{cmm}$, so $S \cong D_4$ and $S_0 \cong C_2^2$}
\label{table:repsGp4mH3}
\end{table}

\begin{table}[!htbp]
\centering
\small
\begin{tabular}{llll|lll|llll}
\toprule
$\boldsymbol{S_{\chi}}$ & & $\boldsymbol{S_{\chi,0}}$ & $\boldsymbol{\HH^2(S_{\chi}, \bC^{\times}_{\pi})}$ & $\boldsymbol{S'_{\chi}}$ & & $\boldsymbol{\HH^2(S'_{\chi}, \bC^{\times}_{\kappa})}$ & $\boldsymbol{\chi}$ & $\boldsymbol{c_{\chi}}$ & $\boldsymbol{c'_{\chi}}$ & $\boldsymbol{\fT(V)}$ \\
\midrule
$C_2$ & $\langle t^{2}x\rangle$ & $1$ & $\mathbb Z/2\mathbb Z$ & $C_2^2$ & $\langle x,\overline{t^{2}}\rangle$ & $(\mathbb Z/2\mathbb Z)^{2}$ & $[\frac{1}{8},\frac{1}{4}]$ & $(0)$ & $(0,0)$ & I \\
\midrule
$C_2^2$ & $\langle s,t^{2}x\rangle$ & $C_2$ & $(\mathbb Z/2\mathbb Z)^{2}$ & $C_2^3$ & $\langle s,x,\overline{t^{2}}\rangle$ & $(\mathbb Z/2\mathbb Z)^{4}$ & $[\frac{1}{4},\frac{1}{4}]$ & $(0,0)$ & $(0,0,0,0)$ & I, I \\
\midrule
$C_2^2$ & $\langle st,sxt^{-1}\rangle$ & $C_2$ & $(\mathbb Z/2\mathbb Z)^{2}$ & $C_2^3$ & $\langle x,st,\overline{st^{-1}}\rangle$ & $(\mathbb Z/2\mathbb Z)^{4}$ & $[0,\frac{1}{4}]$ & $(0,0)$ & $(0,0,0,0)$ & I, I \\
 &  &  &  &  &  &  & $[\frac{1}{2},\frac{1}{4}]$ & $(0,1)$ & $(0,0,0,1)$ & III \\
\midrule
$C_2^2$ & $\langle sx,st^{2}\rangle$ & $C_2$ & $(\mathbb Z/2\mathbb Z)^{2}$ & $C_2^3$ & $\langle x,\overline{s},\overline{t^{2}}\rangle$ & $(\mathbb Z/2\mathbb Z)^{4}$ & $[\frac{1}{4},\frac{3}{4}]$ & $(0,0)$ & $(0,0,0,0)$ & I, I \\
\midrule
$C_2^2$ & $\langle st^{-1},stx\rangle$ & $C_2$ & $(\mathbb Z/2\mathbb Z)^{2}$ & $C_2^3$ & $\langle x,st^{-1},\overline{st}\rangle$ & $(\mathbb Z/2\mathbb Z)^{4}$ & $[\frac{1}{4},0]$ & $(0,0)$ & $(0,0,0,0)$ & I, I \\
 &  &  &  &  &  &  & $[\frac{1}{4},\frac{1}{2}]$ & $(0,1)$ & $(0,1,0,0)$ & III \\
\bottomrule
\end{tabular}
\caption{$G = \text{p4m}$, $H = \langle s,t,xy,xy^{-1} \rangle \cong \text{p4m}$, so $S \cong C_2 \times D_4$ and $S_0 \cong D_4$}
\label{table:repsGp4mH4}
\end{table}

\begin{table}[!htbp]
\centering
\small
\begin{tabular}{llll|lll|llll}
\toprule
$\boldsymbol{S_{\chi}}$ & & $\boldsymbol{S_{\chi,0}}$ & $\boldsymbol{\HH^2(S_{\chi}, \bC^{\times}_{\pi})}$ & $\boldsymbol{S'_{\chi}}$ & & $\boldsymbol{\HH^2(S'_{\chi}, \bC^{\times}_{\kappa})}$ & $\boldsymbol{\chi}$ & $\boldsymbol{c_{\chi}}$ & $\boldsymbol{c'_{\chi}}$ & $\boldsymbol{\fT(V)}$ \\
\midrule
$C_2$ & $\langle t^{2}x\rangle$ & $1$ & $\mathbb Z/2\mathbb Z$ & $C_2^2$ & $\langle x,\overline{t^{2}}\rangle$ & $(\mathbb Z/2\mathbb Z)^{2}$ & $[\frac{1}{8},\frac{1}{4}]$ & $(0)$ & $(0,0)$ & I \\
\midrule
$C_2^2$ & $\langle sx,st^{2}\rangle$ & $C_2$ & $(\mathbb Z/2\mathbb Z)^{2}$ & $C_2^3$ & $\langle s,x,\overline{t^{2}}\rangle$ & $(\mathbb Z/2\mathbb Z)^{4}$ & $[\frac{1}{4},\frac{1}{4}]$ & $(0,0)$ & $(0,0,0,0)$ & I, I \\
\midrule
$C_2^2$ & $\langle st^{-1},stx\rangle$ & $C_2$ & $(\mathbb Z/2\mathbb Z)^{2}$ & $C_2^3$ & $\langle x,st,\overline{st^{-1}}\rangle$ & $(\mathbb Z/2\mathbb Z)^{4}$ & $[0,\frac{1}{4}]$ & $(0,0)$ & $(0,0,0,0)$ & I, I \\
 &  &  &  &  &  &  & $[\frac{1}{2},\frac{1}{4}]$ & $(0,0)$ & $(0,0,0,1)$ & V, V \\
\midrule
$C_2^2$ & $\langle s,t^{2}x\rangle$ & $C_2$ & $(\mathbb Z/2\mathbb Z)^{2}$ & $C_2^3$ & $\langle x,\overline{s},\overline{t^{2}}\rangle$ & $(\mathbb Z/2\mathbb Z)^{4}$ & $[\frac{1}{4},\frac{3}{4}]$ & $(0,0)$ & $(0,0,0,0)$ & I, I \\
\midrule
$C_2^2$ & $\langle st,sxt^{-1}\rangle$ & $C_2$ & $(\mathbb Z/2\mathbb Z)^{2}$ & $C_2^3$ & $\langle x,st^{-1},\overline{st}\rangle$ & $(\mathbb Z/2\mathbb Z)^{4}$ & $[\frac{1}{4},0]$ & $(0,0)$ & $(0,0,0,0)$ & I, I \\
 &  &  &  &  &  &  & $[\frac{1}{4},\frac{1}{2}]$ & $(0,0)$ & $(0,1,0,0)$ & V, V \\
\bottomrule
\end{tabular}
\caption{$G = \text{p4m}$, $H = \langle sx,t,xy,xy^{-1} \rangle \cong \text{p4g}$, so $S \cong C_2 \times D_4$ and $S_0 \cong D_4$}
\label{table:repsGp4mH5}
\end{table}

\begin{table}[!htbp]
\centering
\small
%
\caption{$G = \text{p4g}$, $H = \langle t,x,y \rangle \cong \text{p4}$, so $S \cong D_4$ and $S_0 \cong C_4$}
\label{table:repsGp4gH1}
\end{table}

\begin{table}[!htbp]
\centering
\small
%
\caption{$G = \text{p4g}$, $H = \langle st,t^{2},x,y \rangle \cong \text{cmm}$, so $S \cong D_4$ and $S_0 \cong C_2^2$}
\label{table:repsGp4gH2}
\end{table}

\begin{table}[!htbp]
\centering
\small
%
\caption{$G = \text{p4g}$, $H = \langle s,t^{2},x,y \rangle \cong \text{pgg}$, so $S \cong D_4$ and $S_0 \cong C_2^2$}
\label{table:repsGp4gH3}
\end{table}

\begin{table}[!htbp]
\centering
\small
\begin{tabular}{llll|lll|llll}
\toprule
$\boldsymbol{S_{\chi}}$ & & $\boldsymbol{S_{\chi,0}}$ & $\boldsymbol{\HH^2(S_{\chi}, \bC^{\times}_{\pi})}$ & $\boldsymbol{S'_{\chi}}$ & & $\boldsymbol{\HH^2(S'_{\chi}, \bC^{\times}_{\kappa})}$ & $\boldsymbol{\chi}$ & $\boldsymbol{c_{\chi}}$ & $\boldsymbol{c'_{\chi}}$ & $\boldsymbol{\fT(V)}$ \\
\midrule
$1$ & $\langle\emptyset\rangle$ & $1$ & $0$ & $1$ & $\langle\emptyset\rangle$ & $0$ & $[\frac{1}{8},\frac{3}{8}]$ & $(0)$ & $(0)$ & VII \\
\cmidrule(l){5-11}
 &  &  &  & $C_2$ & $\langle s\rangle$ & $0$ & $[\frac{1}{4},\frac{1}{4}]$ & $(0)$ & $(0)$ & IV \\
\cmidrule(l){5-11}
 &  &  &  & $C_2$ & $\langle st^{-1}\rangle$ & $0$ & $[0,\frac{1}{4}]$ & $(0)$ & $(0)$ & IV \\
\cmidrule(l){5-11}
 &  &  &  & $C_2$ & $\langle st\rangle$ & $0$ & $[\frac{1}{4},0]$ & $(0)$ & $(0)$ & IV \\
\midrule
$C_2$ & $\langle st\rangle$ & $1$ & $\mathbb Z/2\mathbb Z$ & $C_2$ & $\langle \overline{st}\rangle$ & $\mathbb Z/2\mathbb Z$ & $[\frac{1}{4},\frac{1}{2}]$ & $(0)$ & $(0)$ & V \\
\cmidrule(l){5-11}
 &  &  &  & $C_2^2$ & $\langle \overline{1},st\rangle$ & $(\mathbb Z/2\mathbb Z)^{2}$ & $[\frac{1}{2},0]$ & $(0)$ & $(0,0)$ & I \\
\midrule
$C_2$ & $\langle s\rangle$ & $1$ & $\mathbb Z/2\mathbb Z$ & $C_2$ & $\langle \overline{s}\rangle$ & $\mathbb Z/2\mathbb Z$ & $[\frac{1}{4},\frac{3}{4}]$ & $(0)$ & $(0)$ & V \\
\bottomrule
\end{tabular}
\caption{$G = \text{p3m1}$, $H = \langle t,x,y \rangle \cong \text{p3}$, so $S \cong S_3$ and $S_0 \cong C_3$}
\label{table:repsGp3m1H1}
\end{table}

\begin{table}[!htbp]
\centering
\small
%
\caption{$G = \text{p31m}$, $H = \langle t,x,y \rangle \cong \text{p3}$, so $S \cong S_3$ and $S_0 \cong C_3$}
\label{table:repsGp31mH1}
\end{table}

\begin{table}[!htbp]
\centering
\small
%
\caption{$G = \text{p6}$, $H = \langle r^{2},u,v \rangle \cong \text{p3}$, so $S \cong C_6$ and $S_0 \cong C_3$}
\label{table:repsGp6H1}
\end{table}

\begin{table}[!htbp]
\centering
\small
%
\caption{$G = \text{p6m}$, $H = \langle t,x,y \rangle \cong \text{p6}$, so $S \cong D_6$ and $S_0 \cong C_6$}
\label{table:repsGp6mH1}
\end{table}

\begin{table}[!htbp]
\centering
\small
%
\caption{$G = \text{p6m}$, $H = \langle s,t^{2},x,y \rangle \cong \text{p3m1}$, so $S \cong D_6$ and $S_0 \cong S_3$}
\label{table:repsGp6mH2}
\end{table}

\begin{table}[!htbp]
\centering
\small
%
\caption{$G = \text{p6m}$, $H = \langle st,t^{2},x,y \rangle \cong \text{p31m}$, so $S \cong D_6$ and $S_0 \cong S_3$}
\label{table:repsGp6mH3}
\end{table}

\restoregeometry

\section{Mathematical Basis for the Algorithm}\label{sect:algconsiderations}

\subsection{Representing and Computing with the Schur Multiplier of a Finite Group}

Suppose in this section that $G$ is a finite group, and that $H$ is a subgroup of $G$ of index $1$ or $2$. Let $e$ denote the exponent of $\HH^2(G, \mu_{|G|, \pi})$, and suppose that $N$ is a multiple of $e$; for example, $N$ can be any multiple of $|G|$.

Let $\pi \colon G \rightarrow \Gal(\bC / \bR)$ be the homomorphism with kernel $H$, and let $\bC^{\times}_{\pi}$ denote the $G$-module $\bC^{\times}$ with this action of $G$. 

Letting $\mu_{N, \pi}$ denote the $G$-submodule of $N$th roots of $1$, there is a short exact sequence of $G$-modules
\[
1 \rightarrow \mu_{N, \pi} \rightarrow \bC^{\times}_{\pi} \xrightarrow{\times N} \bC^{\times}_{\pi} \rightarrow 0
\]
and thus an exact sequence
\[
\HH^1(G, \bC^{\times}_{\pi}) \xrightarrow {\times N} \HH^1(G, \bC^{\times}_{\pi}) \xrightarrow{\delta} \HH^2(G, \mu_{N, \pi}) \rightarrow \HH^2(G, \bC^{\times}_{\pi}) \xrightarrow{\times N} \HH^2(G, \bC^{\times}_{\pi}).
\]
Here, for $\chi \in \HH^1(G, \bC^{\times}_{\pi})$, $\delta(\chi)$ is defined by choosing some set-theoretic map $a \colon G \rightarrow \bC^{\times}_{\pi}$ with $a(g)^N = \chi(g)$ for all $g \in G$, and setting
\[
    \delta(\chi)(g,h) \coloneqq a(g) \cdot {}^{\pi(g)}a(h) \cdot a(gh)^{-1}.
\]
Because $G$ is finite, both $\HH^1(G, \bC^{\times}_{\pi})$ and $\HH^2(G, \bC^{\times}_{\pi})$ are killed by multiplication by $|G|$ \cite[Thm.\ 6.5.8]{Weib}. Taking $N = |G|$ in the exact sequence above therefore shows that $\HH^1(G, \bC^{\times}_{\pi})$ is a subgroup, and $\HH^2(G, \bC^{\times}_{\pi})$ a quotient, of $\HH^2(G, \mu_{|G|, \pi})$, so both are killed by $e$. Hence for any multiple $N$ of $e$ the first and final maps are zero, and we have an exact sequence
\begin{equation}\label{eqn:schurindexseq}
0 \rightarrow \HH^1(G, \bC^{\times}_{\pi}) \xrightarrow{\delta} \HH^2(G, \mu_{N, \pi}) \rightarrow \HH^2(G, \bC^{\times}_{\pi}) \rightarrow 0.
\end{equation}

We compute $\HH^2(G, \bC^{\times}_{\pi})$ by computing $\HH^1(G, \bC^{\times}_{\pi})$ and $\HH^2(G, \mu_{N, \pi})$, and the connecting map between them.

\subsubsection{Computation of $\HH^1(G, \bC^{\times}_{\pi})$}

The $G$-module $\bC^{\times}_{\pi}$ is not finite, and so we cannot compute $\HH^1(G, \bC^{\times}_{\pi})$ directly in magma. Instead we need a little theory.

There are two cases. If $H = G$, then $\HH^1(G, \bC^{\times}_{\pi}) = \Hom(G, \bC^{\times})$, the character group of $H$, which can be computed in magma as $\Hom(G, C_n)$, working with $C_n$ as a model of $\mu_n \subset \bC^{\times}$.

Suppose instead from now on that $H$ has index two in $G$. Then $\HH^1(G, \bC^{\times}_{\pi})$ is canonically identified with the set of $1$-dimensional antilinear representations of $G$. Here a $1$-cocycle $\mu \colon G \rightarrow \bC^{\times}_{\pi}$ corresponds to the antilinear representation $\bC_{\mu}$, which is $\bC$ as a $\bC$-vector space and action
\[
g \cdot \lambda \coloneqq \overline{\lambda} \mu(g), \qquad \text{ for $g \in G$ and $\lambda \in \bC_{\mu}$.}
\]

In particular, from general theory, we know that the restriction map
\[
\HH^1(G, \bC^{\times}_{\pi}) \rightarrow \Hom(H, \bC^{\times})
\]
is injective, and has image those characters of antilinear type $\bR$. Thus, to compute representatives for $\HH^1(G, \bC^{\times}_{\pi})$, we just need to provide a process which determines which elements of $\Hom(H, \bC^{\times})$ have antilinear type $\bR$, and for each gives an extension to an $1$-cocycle of $G$ valued in $\bC^{\times}_{\pi}$.

The following lemma allows us to do this very efficiently.

\begin{lemma}
  Suppose that $\chi$ is a complex unitary character of a group $H$, $G$ is a Real structure on $H$, and $\chi$ is fixed by the action of $C_2$. The following are equivalent:
  \begin{enumerate}
    \item $\chi$ has antilinear type $\bR$,
    \item $\chi$ is extendable to an element of $Z^1(G, \bC^{\times}_{\pi})$,
    \item $\chi(\omega^2) = 1$ for all $\omega \in G \setminus H$,
    \item $\chi(\omega^2) = 1$ for some $\omega \in G \setminus H$.
  \end{enumerate}
  If these are satisfied, then for each $\omega \in G \setminus H$, an extension $\mu$ of $\chi$ to $Z^1(G, \bC^{\times}_{\pi})$ is defined by setting 
  \[
  \mu(h \omega) \coloneqq \chi(h)
  \]
  for each $h \in H$.
\end{lemma}

\begin{proof}
  The first two points are equivalent from the general theory, and it is clear that the the third implies the fourth. To see the second implies the third, note that if $\chi$ admits an extension $\mu$, then for any $\omega \in G \setminus H$,
  \[
  \chi(\omega^2) = \mu(\omega^2) = \mu(\omega) \overline{\mu(\omega)} \in \bR_{>0}.
  \]
  Thus $\chi(\omega^2) = 1$, as $\chi(\omega^2) \in S^1$ by the assumption that $\chi$ is unitary. Finally we show that the fourth point implies the second, so suppose that $\chi(\omega^2) = 1$ for some $\omega \in G \setminus H$. Setting
  \[
  \mu(h \omega) \coloneqq \chi(h),
  \]
  it is straightforward to check that $\mu(gh) = \mu(g) \cdot {}^{\pi(g)} \mu(h)$ is satisfied for all $g,h \in G$, by considering the four possible cases for $g$ and $h$ being in $H$ or in $G \setminus H$, the properties of $\chi$ needed being that $\chi$ is unitary, ${}^{\omega}\chi = \overline{\chi}$, and $\chi(\omega^2) = 1$. This also shows the final claim.
\end{proof}

\subsubsection{Computation of $\HH^2(G, \mu_{N, \pi})$}

Now, because the underlying abelian group of $\mu_{N, \pi}$ is finite, magma has built in functions for computing $\HH^2(G, \mu_{N, \pi})$. The only care needed is that magma works with right actions and cocycles, and we with left, so we need to convert between these.

Recall that if $G$ is a group and $M$ is a left $G$-module, we have the bar resolution $(C^i(G,M), d_i)$, where $d_i \colon C^i(G,M) \rightarrow C^{i+1}(G,M)$,
\begin{align*}
d_i (f)(g_1,\ldots,g_{n+1}) \coloneqq &g_1 f(g_2, ... , g_{n+1}) + \sum_{i=1}^n (-1)^i
 f(g_1,\ldots,g_ig_{i+1},\ldots,g_{n+1}), \\
 &\qquad + (-1)^{n+1}f(g_1,\ldots,g_n).
\end{align*}

On the other hand, if $N$ is a right $G$-module, the bar resolution used in magma (as far as we can gather) is $(C^i_r(G,M), d_i^r)$, where $d_i^r \colon C^i_r(G,M) \rightarrow C_r^{i+1}(G,M)$ is defined by
\begin{align*}
d_i^r (f)(g_1,\ldots,g_{n+1}) \coloneqq &f(g_2, ... , g_{n+1}) + \sum_{i=1}^n (-1)^i
 f(g_1,\ldots,g_ig_{i+1},\ldots,g_{n+1}), \\
 &\qquad + (-1)^{n+1}f(g_1,\ldots,g_n) \cdot g_{n+1}.
\end{align*}

Now, if $M$ is a left $G$-module, we can obtain a right $G$-module structure $M_r$ on $M$ be setting $M * g \coloneqq g^{-1}(m)$. The identification 
\[
\HH^i_r(G,M_r) \xrightarrow{\sim} \HH^i(G,M)
\]
that we use in the code is that induced by the isomorphism of chain complexes,
\[
\phi_i(f)(g_1, ... , g_i) \coloneqq g_1 \cdots g_i f(g_1, ... ,g_i).
\]

The only real check to make is that this is an isomorphism of chain complexes.

\begin{lemma}
  We have that $\phi_{i+1} \circ d_i^r = d_i \circ \phi_i$, for all $i \geq 0$.
\end{lemma}

\begin{proof}
We may compute:
\[
\begin{aligned}
&d_i\phi_i(f)(g_1,\ldots,g_{n+1})\\
&=g_1\phi_i(f)(g_2,\ldots,g_{n+1})
 +\sum_{i=1}^n(-1)^i
 \phi_i(f)(g_1,\ldots,g_ig_{i+1},\ldots,g_{n+1})\\
&\qquad+(-1)^{n+1}\phi_i(f)(g_1,\ldots,g_n)\\[2mm]
&=g_1(g_2\cdots g_{n+1})\,f(g_2,\ldots,g_{n+1})\\
&\qquad+\sum_{i=1}^n(-1)^i
 (g_1\cdots g_ig_{i+1}\cdots g_{n+1})
 f(g_1,\ldots,g_ig_{i+1},\ldots,g_{n+1})\\
&\qquad+(-1)^{n+1}(g_1\cdots g_n)f(g_1,\ldots,g_n)\\[2mm]
&=(g_1\cdots g_{n+1})
\Bigg[
f(g_2,\ldots,g_{n+1})
+\sum_{i=1}^n(-1)^i
f(g_1,\ldots,g_ig_{i+1},\ldots,g_{n+1})\\
&\hspace{55mm}
+(-1)^{n+1}f(g_1,\ldots,g_n)*g_{n+1}
\Bigg]\\[2mm]
&=(g_1\cdots g_{n+1})(d_rf)(g_1,\ldots,g_{n+1})\\
&=\phi_{i+1}(d_i^rf)(g_1,\ldots,g_{n+1}).
\end{aligned}
\]
\end{proof}

Concretely, from a (right) $2$-cocycle $t$ in magma, we obtain a (left) $2$-cocycle
\[
c(g,h) \coloneqq gh t(g,h),
\]
and from a function $b \colon G \rightarrow M$ with $t_1 - t_2 = d_1^r b$, we have that
\[
c_1 - c_2 = \phi_2(t_1) - \phi_2(t_2) = \phi_2(d_1^rb) = d_1(\phi_1(b)),
\] 
where $\phi_1(b)(g) = g b(g)$.

\subsection{Projective Antilinear Representations}\label{sect:projantireps}

We continue to assume that $G$ is a finite $C_2$-graded group with ungraded subgroup $H$, and $N$ is a multiple of $n \coloneqq |G|$. 

In this section we how we can compute with categories of projective antilinear representations but reducing to ordinary (non-projective) antilinear representations, where the theory is well-developed.

\begin{defn}
  For $\alpha \in \HH^2(G, \mu_{N, \pi})$, an $\alpha$-projective antilinear representation of $G$ is a finite-dimensional complex vector space $V$ together with, for each $g \in G$, an invertible $\pi(g)$-linear operator $\rho(g) \in \End_{\bR}(V)$ such that
  \[
  \alpha(g,h) \rho(gh) = \rho(g)\rho(h)
  \]
  for any $g,h \in G$. Let $\Rep_{\bC, \alpha}^{\hspace{0.1em} \star}(G)$ denote the category of such representations.
\end{defn}

Recall that $\alpha \in  \HH^2(G, \mu_{N, \pi})$ is normalised if $\alpha(g,h) = 1$ whenever $g$ or $h$ is trivial. 

\begin{defn}
  Suppose that $\alpha \in \HH^2(G, \mu_{N, \pi})$ is normalised, and define a group $E_{\alpha}(G)$ by setting
  \[
  E_{\alpha}(G) \coloneqq \mu_N \times G
  \]
  as a set, with multiplication
  \[
  (\lambda, g) \cdot (\mu, h) \coloneqq (\lambda ({}^{\pi(g)} \! \mu) \alpha(g,h), gh).
  \]
  The natural inclusion and projection defines a short exact sequence of groups
  \[
  1 \rightarrow \mu_N \rightarrow E_{\alpha}(G) \rightarrow G \rightarrow 1.
  \]
  There is also a canonical section
  \[
    s \colon G \rightarrow E_{\alpha}(G), \qquad s(g) = (1, g),
  \]
  and subgroup 
  \[
  E_{\alpha}(G)_0 = \mu_N \times H \subset E_{\alpha}(G)
  \]
  which gives  $E_{\alpha}(G)$ the structure of a $C_2$-graded group.
\end{defn}

The fact that $\alpha$ is normalised needed for $E_{\alpha}(G)$ to define a group.

\begin{defn}
  For normalised $\alpha \in \HH^2(G, \mu_{N, \pi})$, let $\Rep_{\bC, \mu_N}^{\hspace{0.1em} \star}(E_{\alpha}(G))$ denote the full subcategory of the category $\Rep_{\bC}^{\hspace{0.1em} \star}(E_{\alpha}(G))$ of finite-dimensional antilinear representations of $E_{\alpha}(G)$, where $\mu_N$ acts as scalar multiplication.

  We define a functor 
  \[
  \Phi \colon \Rep_{\bC, \alpha}^{\hspace{0.1em} \star}(G)
\; \rightarrow \;
\Rep_{\bC, \mu_N}^{\hspace{0.1em} \star}(E_{\alpha}(G))
  \]
  to be the enrichment of the identity functor of $\Vect_{\bC}$, where for $V \in \Rep_{\bC, \alpha}^{\hspace{0.1em} \star}(G)$, $E_{\alpha}(G)$ acts on $\Phi(V)$ by
  \[
  (\lambda, g) \cdot v \coloneqq \lambda g(w).
  \] 

  In the other direction, we define a functor 
  \[
  \Psi \colon \Rep_{\bC, \mu_N}^{\hspace{0.1em} \star}(E_{\alpha}(G))
\; \rightarrow \;
\Rep_{\bC, \alpha}^{\hspace{0.1em} \star}(G)
  \]
  which again is an enrichment of the identity functor of $\Vect_{\bC}$, where now for $W \in \Rep_{\bC, \mu_N}^{\hspace{0.1em} \star}(E_{\alpha}(G))$, $G$ acts on $\Psi(W)$ via the section $s$:
  \[
  g \cdot w \coloneqq s(g)(w).
  \]
\end{defn}

That $\alpha$ is normalised is used to show that $\Psi(W)$ defines a projective representation. The following is direct to verify.

\begin{lemma}
  For $\alpha \in \HH^2(G, \mu_{N, \pi})$, normalised, the functors $\Phi$ and $\Psi$ are well-defined, preserve irreducibility, and define mutually inverse equivalences of categories
    \[
\Phi \colon \Rep_{\bC, \alpha}^{\hspace{0.1em} \star}(G)
\; \leftrightarrow \;
\Rep_{\bC, \mu_N}^{\hspace{0.1em} \star}(E_{\alpha}(G)) \vcentcolon \Psi
\]
and in particular induce an identification of irreducible objects.
\end{lemma}

Suppose now that $\alpha, \beta$ are normalised $2$-cocycles which define the same class in $\HH^2(G, \mu_{N, \pi})$, so there is some set-theoretic map $b \colon G \rightarrow \mu_N$ with 
\[
\beta(g,h)\alpha(g,h)^{-1} = db = b(g) \cdot {}^{\pi(g)}b(h) \cdot b(gh)^{-1}.
\]
We can do two things with $b$. On the one hand, we can use it to define equivalences between the corresponding categories of projective representations. The following is direct to verify.

\begin{lemma}
  For $\alpha, \beta$ and $b$ as above, the functor
  \[
  \Pi_{\alpha, \beta}^b \colon \Rep_{\bC, \alpha}^{\hspace{0.1em} \star}(G) \rightarrow \Rep_{\bC, \beta}^{\hspace{0.1em} \star}(G),
  \]
  where
  \[
  \Pi_{\alpha, \beta}^b(V, \rho) = (V, \rho_b), \qquad \rho_b(g) \coloneqq b(g)\rho(g)
  \]
  and $\Pi_{\alpha, \beta}^b$ is the identity on morphisms, is an equivalence of categories, with inverse $\Pi_{\beta, \alpha}^{b^{-1}}$.
\end{lemma}

The second thing we can do with $b$ is use it do define an isomorphism between $E_{\alpha}(G)$ and $E_{\beta}(G)$.

\begin{lemma}
  The map
  \[
    \phi_{\beta, \alpha}^b \colon E_{\beta}(G) \rightarrow E_{\alpha}(G), \qquad \phi_{\beta, \alpha}^b(\lambda, g) \coloneqq (\lambda b(g), g)
  \]
  is an isomorphism of $C_2$-graded groups, with inverse $\phi_{\alpha, \beta}^{b^{-1}}$.
\end{lemma}

Again, this is immediate to verify. From this isomorphism, we obtain a pull-back functor
\[
(\phi_{\beta, \alpha}^b)^* \colon \Rep_{\bC, \mu_N}^{\hspace{0.1em} \star}(E_{\alpha}(G)) \rightarrow \Rep_{\bC, \mu_N}^{\hspace{0.1em} \star}(E_{\beta}(G)),
\]
where for $W \in \Rep_{\bC, \mu_N}^{\hspace{0.1em} \star}(E_{\alpha}(G))$, $(\phi_{\beta, \alpha}^b)^*(W)$ is the $\bC$-vector space $W$ with action of $E_{\beta}(G)$ via composition with $\phi_{\beta, \alpha}^b$. As $\alpha$ and $\beta$ are normalised, $b(1) = 1$, and therefore $\mu_N \subset E_{\beta}(G)$ acts on $(\phi_{\beta, \alpha}^b)^*(W)$ via the natural multiplication action.

These are compatible in the following sense (which is again direct to verify).

\begin{lemma}\label{lem:projcommdiag}
  The diagram of functors
\[\begin{tikzcd}
	{\Rep_{\bC, \mu_N}^{\hspace{0.1em} \star}(E_{\alpha}(G))} & {\Rep_{\bC, \alpha}^{\hspace{0.1em} \star}(G)} \\
	{\Rep_{\bC, \mu_N}^{\hspace{0.1em} \star}(E_{\beta}(G))} & {\Rep_{\bC, \beta}^{\hspace{0.1em} \star}(G)}
	\arrow["{\Psi_{\alpha}}", from=1-1, to=1-2]
	\arrow["{(\phi_{\beta, \alpha}^b)^*}"', from=1-1, to=2-1]
	\arrow["{\Pi_{\alpha, \beta}^b}", from=1-2, to=2-2]
	\arrow["{\Psi_{\beta}}"', from=2-1, to=2-2]
\end{tikzcd}\]
  commutes.
\end{lemma}

Note that the diagram commutes on the nose, we don't mean up to natural isomorphism here.




\subsection{Exact Computation of Frobenius--Schur Indicators}

Suppose that we are in the situation of Section \ref{sect:cliffordtheory}, so $T \leq H \leq G$, with $[G: H] = 2$, $T$ abelian and $T$ is finite index and normal in $G$. In this section, we additionally assume that $T$ is free-abelian of finite rank $r$.

In this case, we can compute the Frobenius--Schur indicators of a unitary representation $V$ of $H$ as limits of averages over a F{\o}lner sequence - see Proposition \ref{prop:FSamenable}. The aim of this section is to show that when $V$ arises from a unitary character of finite order, this limit can be computed exactly as a finite sum.

Fix a character $\chi \in T^{\vee}$ of finite order $N$, and let $U \coloneqq S_{\chi, 0} = S_{\chi} \cap S_0$, and let $G_U$ be the preimage of $U$ in $G$. Let $\eta \colon G_U \rightarrow \bC$ be a class function of the form
\[
\eta(t \breve{s}) = \chi(t)\sigma(s),
\]
for $t \in T$ and $s \in U$, where $\sigma \colon U \rightarrow \bC$ is any function. Let $\theta \colon H \rightarrow \bC$ be the induced class function
\begin{equation}\label{eqn:thetainduced}
\theta(h) \coloneqq \sum_{\substack{y \in G_U \backslash H \\ y h y^{-1} \in G_U}} \eta(y h y^{-1}), \qquad h \in H.
\end{equation}
Let $T^N \coloneqq \{t^N \mid t \in T\}$ denote the $N$th power subgroup of $T$, which is normal in $G$.

\begin{lemma}\label{lem:basicpropertiesinvariance}
  For all $x \in T^N$,
  \begin{enumerate}
    \item $\theta(xh) = \theta(h)$ for all $h \in H$,
    \item $\theta((xg)^2) = \theta(g^2)$ for all $g \in G$,
    \item The same hold for $\theta$ replaced by $\overline{\theta}$ or ${}^{\omega}\theta$, for any $\omega \in G \setminus H$.
  \end{enumerate}
\end{lemma}
Recall that $({}^{\omega}\theta)(h) = \theta(\omega^{-1} h \omega)$ for any $h \in H$ by definition.

\begin{proof}
  For (1), by (\eqref{eqn:thetainduced}), it suffices to show that for any $y \in G_U \backslash H$, $y h y^{-1} \in G_U$ if and only if $y xh y^{-1} \in G_U$, and in this case that $\eta(yx hy^{-1}) = \eta(yhy^{-1})$. The first is direct, as $yxhy^{-1} = (yxy^{-1})(yhy^{-1})$, and $yxy^{-1} \in T \subset G_U$. For the second,
  \[
  \eta(yx hy^{-1}) = \eta(yxy^{-1} y h y^{-1}) = \chi(yxy^{-1})\eta(yhy^{-1}) = \eta(yhy^{-1})
  \]
  because $\chi(yxy^{-1}) = 1$ as $yxy^{-1} \in T^N$ and $\eta(th) = \chi(t) \eta(h)$ for any $t \in T$ and $h \in H$. 
  
  Point (2) follows from (1), as $(xg)^2 = xg xg = (x gxg^{-1}) g^2$, and $x, gxg^{-1} \in T^N$. For the final point, note that $\overline{\theta}$ is of the same form as $\theta$, for $\chi$ and $\sigma$ replaced by $\overline{\chi}$ and $\overline{\sigma}$, and the claim for ${}^{\omega}\theta$ follows from the fact that $\omega^{-1} T^N \omega = T^N$.
\end{proof}

Using this, we have the following general exact formula for the limits of Proposition \ref{prop:FSamenable}.

\begin{lemma}\label{lem:meanperiodic}
Let $H$ be a discrete countable amenable group with F{\o}lner sequence $(F_n)_n$, let $M \trianglelefteq H$ be a normal subgroup of finite index, and let $f \colon H \rightarrow \bC$ with $f(xh) = f(x)$ for all $x \in M$ and $h \in H$. Then
\[
\lim_{n \rightarrow \infty} \frac{1}{|F_n|} \sum_{h \in F_n} f(h) = \frac{1}{[H : M]} \sum_{h \in H / M} f(x).
\]
\end{lemma}

\begin{proof}
  From the fact that $f(xh) = f(h)$, it suffices to take $f$ to be the indicator function of a coset $h M$. In this case, the claim is exactly that the indicator function $a_n(h) \coloneqq |F_n \cap hM| / |F_n|$ converges to $1 / [H:M]$. From the definition of the F{\o}lner sequence, we have that 
  \[
  |a_n(gh) - a_n(h)| = \frac{\left| |g^{-1} F_n \cap hM| - |F_n \cap hM| \right|}{|F_n|} \leq \frac{|g^{-1} F_n \oplus F_n|}{|F_n|} \longrightarrow 0
  \]
  for any $g \in H$. As $H$ acts transitively on $H/M$, the limits of the $a_n(h)$ are the same for any $h \in H$, and thus as $\sum_{h \in H / M} a_n(h) = 1$ for any $n \geq 1$, this forces the limit of any $a_n(h)$ to be $1 / [H:M]$.
\end{proof}

We can make this more precise in the situation of interest. We now take $\sigma$ to be the character of an irreducible projective antilinear representation $W$ of $S_{\chi}$, and so $\theta$ is the character of the corresponding irreducible antilinear representation of $G$. To fix notation, suppose that $T$ has basis $e_1, \ldots, e_r$, and for $a = (a_1, \ldots, a_r) \in \bZ^r$ we write $e^a \coloneqq e_1^{a_1} \cdots e_r^{a_r} \in \bC^{\times}$.

\begin{cor}
  Let $\chi \in T^{\vee}$ of order $N$, $W \in \Irr_{\chi}^{\hspace{0.1em} \star}(S_\chi)$ with character $\sigma$, let $(F_n)_n$ be a F{\o}lner sequence for $H$, and let $\omega \in G \setminus H$. Then $\theta = \chi_V$, for the associated $V \in \Irr_{\bC}^{\hspace{0.1em} \star}(G)$, and
  \begin{align*}
\sF_{\Sigma}(V) = \lim_{n \rightarrow \infty} \frac{1}{|F_n|} \sum_{x \in F_n} \chi_V(x^2) &= \frac{1}{N^r |S_0|} \sum_{a \in \{0, \ldots, N-1\}^r} \ \sum_{s \in S_0} \chi_V \left( (e^a \breve{s})^2 \right), \\
\sF_{\Sigma}^{\sharp}(V) = \lim_{n \rightarrow \infty} \frac{1}{|F_n|} \sum_{x \in F_n} \chi_V((\omega x)^2) &= \frac{1}{N^r |S_0|} \sum_{a \in \{0, \ldots, N-1\}^r} \ \sum_{s \in S \setminus S_0} \chi_V \left( (e^a \breve{s})^2 \right).
\end{align*}
Moreover, $\chi_V = {}^{\omega} \chi_V$ (respectively $\chi_V = \overline{\chi_V}$) if and only if the two functions agree on the finite set $\{ e^a \breve{s} \mid a \in \{0, \ldots, N-1\}^r, \ s \in S_0 \}$.
\end{cor}

\begin{proof}
  The first formula follows from Lemma \ref{lem:meanperiodic} with $M = T^N$ and $f(x) \coloneqq \chi_V(x^2)$, which is invariant under left-multiplication by Lemma \ref{lem:basicpropertiesinvariance}. The second formula similarly follows from Lemma \ref{lem:meanperiodic} with $M = T^N$ and $f(x) \coloneqq \chi_V((\omega x)^2)$, again using Lemma \ref{lem:basicpropertiesinvariance}. The final claim follows from Lemma \ref{lem:basicpropertiesinvariance}, as two periodic functions on $H$ are equal if and only if they agree on a transversal of $T^N$ in $H$.
\end{proof}